\documentclass{dis}
\usepackage{amsmath}
\usepackage{amsfonts}
\usepackage{amsthm}
\usepackage{amstext}
\usepackage{amsmath}
\usepackage{amscd}
\usepackage{amssymb}
\usepackage[mathscr]{eucal}
\usepackage{diagmac2}
\usepackage{booktabs}
\usepackage{rotating}
\usepackage{breakcites}
\usepackage{graphicx,epstopdf}
\usepackage{floatrow}
\usepackage{caption}
\usepackage{imakeidx}
\makeindex[name=symbol,title={List of symbols},intoc]
\makeindex[intoc]

\usepackage{remreset}
\makeatletter
\@removefromreset{figure}{chapter}
\makeatother

\newcommand\mL{L\kern-0.08cm\char39}
\newcommand{\pr}{\mathop{\rm pr}}

\newcommand{\Id}{\mathop{\rm Id}}
\newcommand{\id}{\mathop{\rm id}}
\newcommand{\card}{\mathop{\rm card}}

\newcommand{\tor}{\mathop{\rm tor}}
\newcommand{\tork}{\mathop{\rm tor_k}}
\newcommand{\Tor}{\mathop{\rm Tor}}
\newcommand{\im}{\mathop{\rm im}}
\newcommand{\ord}{\mathop{\rm ord}}
\newcommand{\Div}{\mathop{\rm Div}}
\newcommand{\Hom}{\mathop{\rm Hom}}
\newcommand{\Mon}{\mathop{\rm Mon}}
\newcommand{\Ext}{\mathop{\rm Ext}}
\newcommand{\Flow}{\mathop{\rm \mathcal F\colon\Gamma\curvearrowright X}}
\newcommand{\Coc}{\mathop{\rm \mathsf{CAGpZ}_{\mathcal F}}}
\newcommand{\Cob}{\mathop{\rm {\bf B}_{\mathcal F}}}
\newcommand{\Cobcl}{\mathop{\rm \overline{\bf{B}}_{\mathcal F}}}
\newcommand{\Coctd}{\mathop{\rm {\bf Z}^{td}_{\mathcal F}}}
\newcommand{\Coccn}{\mathop{\rm {\bf Z}^{cn}_{\mathcal F}}}

\newcommand{\Coch}{\mathop{\rm {\bf H}_{\mathcal F}}}
\newcommand{\Cochtd}{\mathop{\rm {\bf H}^{td}_{\mathcal F}}}
\newcommand{\Cochcn}{\mathop{\rm {\bf H}^{cn}_{\mathcal F}}}

\newcommand{\Cocc}{\mathop{\rm {\bf Z}_{\mathcal F}}}
\newcommand{\Cocm}{\mathop{\rm {\bf MinZ}_{\mathcal F}}}
\newcommand{\Cochm}{\mathop{\rm {\bf MinH}_{\mathcal F}}}
\newcommand{\co}{\mathop{\rm co}}
\newcommand{\BSdim}{\mathop{\rm dim_{\mathbb R BS}}}
\newcommand{\LSdim}{\mathop{\rm dim_{\mathbb R LS}}}
\newcommand{\QLSdim}{\mathop{\rm dim_{\mathbb Q LS}}}

\newcommand{\Gimel}{\mathop{\rm \gimel_{\mathcal F}}}
\newcommand{\Daleth}{\mathop{\rm \daleth_{\mathcal F}}}
\newcommand{\cGimel}{\mathop{\rm \gimel_{\mathcal F}^c}}
\newcommand{\cDaleth}{\mathop{\rm \daleth_{\mathcal F}^c}}
\newcommand{\mGimel}{\mathop{\rm \gimel_{\mathcal F}^m}}
\newcommand{\mDaleth}{\mathop{\rm \daleth_{\mathcal F}^m}}

\newcommand{\rank}{\mathop{\rm rank}}
\newcommand{\w}{\mathop{\rm w}}

\newcommand{\per}{\mathop{\rm per}}
\newcommand{\sbgp}{\mathrel{\lhd}}

\newcommand{\md}{d\kern-0.035cm\char39\kern-0.03cm}
\newcommand{\mt}{t\kern-0.035cm\char39\kern-0.03cm}
\newcommand{\ml}{l\kern-0.035cm\char39\kern-0.03cm}

\newtheorem{theorem}{Theorem}[chapter]
\newtheorem{proposition}[theorem]{Proposition}
\newtheorem{lemma}[theorem]{Lemma}
\newtheorem{corollary}[theorem]{Corollary}

\theoremstyle{definition}
\newtheorem{definition}[theorem]{Definition}
\newtheorem{example}[theorem]{Example}
\newtheorem{remark}[theorem]{Remark}

\chaptersnewpage

\refnewpage

\begin{document}

\mathclass{Primary 37B05, 54H20; Secondary 22C05, 43A40.}
\thanks{This work was supported by the Slovak Research and
Development Agency under the contract APVV-15-0439 and by VEGA~1/0786/15.}

\abbrevauthors{M. Dirb\'ak}
\abbrevtitle{First Cohomology Groups of Minimal Flows}

\title{First Cohomology Groups of Minimal Flows}

\author{Mat\'u\v s Dirb\'ak}
\address{Department of Mathematics, Faculty of Natural Sciences, Matej Bel University, Bansk\'a Bystrica, Slovakia\\
E-mail: Matus.Dirbak@umb.sk}

\maketitledis

\tableofcontents
\begin{abstract}
Our interest in this work is in group extensions of minimal flows with compact abelian groups in the fibres. We study their structure from categorical and algebraic points of view, and describe relations of their dynamics to the one-dimensional algebraic-topological invariants. We determine the first cohomology groups of flows with simply connected acting groups and those of to\-po\-lo\-gi\-cally free flows possessing a free cycle. As an application we show that minimal extensions of these flows not only do exist, but they have a rich algebraic structure.
\end{abstract}
\makeabstract

\chapter{Introduction}

\section{Minimal flows}

The theory of dynamical systems is nowadays a deeply developed field of ma\-the\-ma\-tics, which interacts with many mathematical disciplines. As the field possesses many branches, one encounters several opinions on what the subject of the theory of dynamical systems is (or should be). One point of view is that the theory of dynamical systems is the study of the global properties of groups of transformations. The subject has its roots in celestial and statistical mechanics, where the group is $\mathbb R$. In the abstract part of the theory ``larger'' groups are considered, such as Lie groups or, more generally, locally compact groups and Polish groups.

In this work we are interested in the topological aspect of the theory of dynamical systems, where the main object of study is a flow (or a semiflow). By a \emph{flow} or a \emph{transformation group} or a \emph{(dynamical) system} $\mathcal F$ we mean a representation of a topological (usually locally compact) group $\Gamma$ as a group of homeomorphisms on a topological (usually compact or locally compact) space $X$, which is continuous as a mapping $\Gamma\times X\to X$; we write $\Flow$.

One of the most important classes of dynamical systems in this setting is formed by the minimal ones, which were defined for the first time by Birkhoff in \cite{Bir}. Recall that a flow $\mathcal F$ is \emph{minimal} if all its orbits are dense. The importance of minimal flows in topological dynamics is highlighted by the following facts. First, by a Birkhoff's theorem, every flow on a compact space contains a (closed) minimal subflow. Second, minimal systems are the irreducible dynamical systems in the sense that they are the ones that do not possess any proper (closed) subsystem. Third, minimal systems are often viewed as topological analogues of ergodic systems from ergodic theory, the latter ones being probably the most important systems in the theory of dynamical systems as a whole.

In this work we will be interested mainly in minimal flows with \emph{connected} acting groups, a special attention will be paid to actions of connected Lie groups. (Notice that a minimal flow with a connected acting group has automatically a connected phase space.) In fact, our main results will concern first cohomology groups of minimal flows $\mathcal F$ of the following two types.
\begin{itemize}
\item The acting group $\Gamma$ of $\mathcal F$ is a connected Lie group and the phase space $X$ of $\mathcal F$ is a compact (connected) manifold.
\item The acting group $\Gamma$ of $\mathcal F$ is a simply connected Lie group (more generally, a simply connected topological group) and the phase space $X$ of $\mathcal F$ is a (not necessarily locally connected) continuum.
\end{itemize}
In the first of these two situations we shall be particularly interested in flows which are to\-po\-lo\-gi\-cally free and/or possess a free cycle (see Subsection~\ref{Sub:top.frnss.fr.ccl} for definitions of these two concepts). In the second situation we shall pay a special attention to flows whose phase spaces $X$ have a non-trivial first cohomotopy group $\pi^1(X)$. As a matter of fact, the algebraic properties of the one dimensional invariants of the space $X$ (the first weak homology group in the first case and the first cohomotopy group in the second case) will occur frequently in our discussions and results on the first cohomology groups of the flow $\mathcal F$.

\section{The structure of minimal flows}\label{S:strct.theory}

It is often the case in mathematics that in order to understand a general object, one tries to express it in terms of special objects, special morphisms and concrete operations. Such is the situation also in the theory of minimal flows. Let us briefly recall several classical results in this direction (for a detailed exposition with some interesting applications, we refer to \cite{Gla2}). Throughout this whole section we shall assume that phase spaces of the considered flows are compact.

The first result in this direction was obtained by Furstenberg (see \cite[Theorem~2.4]{Fur2}), who showed that every metrizable minimal distal flow can be obtained from the trivial flow by using isometric extensions and inverse limits along segments of countable ordinals. As important corollaries of this structure theorem he proved that every distal flow possesses an invariant probability measure (\cite[Theorem~12.3]{Fur2}) and that no simply connected space supports a minimal distal flow with an abelian locally compact acting group (\cite[Theorem~11.1]{Fur2}). One of the essential tools in Furstenberg's proof was the enveloping semigroup of a flow, defined by Ellis in \cite{Ell3} (see also \cite{Ell5}), which turns out to be a group if the flow under consideration is distal. (For many interesting results on Ellis enveloping semigroups or Ellis actions we refer to \cite{Ell} and \cite{AkiAusGla}.) A few years later after Furstenberg published his structure theorem, a different proof of a strengthening of his result to the class of quasi-separable minimal distal flows was presented by Ellis (see \cite[Chapters~13,~15]{Ell}), the main ingredient in his new proof being his own discovery of a Galois theory of distal extensions. (For yet another proof of the Furstenberg's result we refer to \cite[Chapter~20]{EllEll}, where minimal flows are treated as invariant closed equivalence relations on the universal minimal flow.)

Soon after the discovery of the structure of minimal distal flows, an analogous structure theorem was obtained also for metrizable point-distal minimal flows. Under the assumption of residuality of the set of distal points, Veech proved that such systems possess almost automorphic extensions, which can be obtained from the trivial flow by using (proper) isometric extensions, almost automorphic extensions and inverse limits (see \cite[Theorem~7.2]{Vee}). As was conjectured by Veech and later proved by Ellis, the assumption of residuality of the set of distal points is redundant (see \cite{Ell4}, in which the Veech structure theorem was extended to the class of quasi-separable point-distal flows). Later Bron\v ste\u\i n discovered a structure theorem for minimal prodal flows. He showed that such flows lift across proximal epimorphisms to flows, which can be constructed from the trivial flow by using isometric and proximal extensions and inverse limits (see \cite{Bro1}).

A structure theorem for general minimal systems was proved by Ellis, Glasner, Sha\-pi\-ro, McMahon and Veech (see \cite{EllGlaSha}, \cite{McM}, \cite{Vee2} and \cite{Gla3}). It involves the notion of PI-towers, relatively incontractible extensions, isometric extensions, proximal extensions, weakly mixing extensions and inverse limits. The picture is more complex than in the structure theorems mentioned above (and we shall not go into details), but it should be noted that some of the earlier structure theorems follow from this general one (this is the case, for instance, with the structure of point-distal flows (see \cite[pp.~197--198]{Gla2}) and the structure of prodal flows (see \cite[pp.~198--199]{Gla2})). 

A structure theorem for minimal normal flows has been proved by Glasner, Mentzen and Siemaszko (see \cite{GlaMenSie}). They showed that every minimal normal flow lifts across a proximal epimorphism to a minimal normal flow, which descends via the composition of a virtually weakly mixing extension and a group extension to a (minimal) proximal flow. The situation simplifies if the acting group of the flow under consideration is abelian, in which case the normal minimal flow descends onto its maximal Kronecker factor via a virtually weakly mixing extension map.

We conclude our short discussion on the structure of minimal flows by mentioning several results on the structure of minimal tame flows. In \cite{Gla6} Glasner proved that minimal tame flows with abelian acting groups are PI and have zero topological entropy. Later, Huang, Kerr, Li and Glasner showed that such flows are almost automorphic and uniquely ergodic (see \cite{Hua}, \cite{KerLi}, \cite{Gla4}). Recently, a structure theorem for minimal tame flows with an \emph{arbitrary} acting group was proved by Glasner (see \cite[Theorem~5.3]{Gla5}). As a corollary of this very general theorem, Glasner proved that the earlier results on minimal tame flows remain true for flows with amenable acting groups (see \cite[Corollary~5.4]{Gla5}). More concretely, every mi\-ni\-mal tame flow with an amenable acting group is almost automorphic and uniquely ergodic.

\section{Cochain complexes associated to a minimal flow}\label{S:coc.compl.min.flw}

With every (minimal) flow and every abelian topological group one can associate in a natural way a cochain complex and the corresponding cohomology groups. Before describing this construction, let us recall a more general concept, namely the cohomology of (right) $\Gamma$-modules for a given group $\Gamma$. (For further information on the results of this section we refer to \cite[Chapter~7]{Ser} and \cite[Section~1.21]{Tao}.)

Let $\Gamma$ be a group and $\Lambda$ be a right $\Gamma$-module; that is, $\Lambda$ is an abelian group, on which $\Gamma$ acts from the right by automorphisms. For $n\geq1$ let $F^n(\Gamma,\Lambda)$ consist of all functions $\Gamma^n\to\Lambda$. With operations defined point-wise, $F^n(\Gamma,\Lambda)$ is an abelian group. Further, let $F^0(\Gamma,\Lambda)=\Lambda$ and $F^n(\Gamma,\Lambda)=0$ for $n\leq-1$. Then there is a cochain complex of abelian groups
\begin{equation}\label{Eq:coc.compl.group}
\dots\stackrel{\delta_{n-2}}{\longrightarrow}F^{n-1}(\Gamma,\Lambda)\stackrel{\delta{n-1}}{\longrightarrow}F^n(\Gamma,\Lambda)\stackrel{\delta_n}{\longrightarrow}F^{n+1}(\Gamma,\Lambda)\stackrel{\delta_{n+1}}{\longrightarrow}\dots,
\end{equation}
where the morphisms $\delta_n$ are defined as follows. First, if $n\geq1$ then
\begin{equation*}
\begin{split}
(\delta_n\varphi)(\gamma_1,\dots,\gamma_{n+1})&=\varphi(\gamma_2,\dots,\gamma_{n+1})\gamma_1\\
&+\sum_{i=1}^n(-1)^i\varphi(\gamma_1,\dots,\gamma_{i-1},\gamma_{i+1}\gamma_i,\gamma_{i+2},\dots,\gamma_{n+1})\\
&+(-1)^{n+1}\varphi(\gamma_1,\dots,\gamma_n)
\end{split}
\end{equation*}
for all $\varphi\in F^n(\Gamma,\Lambda)$ and $\gamma_1,\dots,\gamma_{n+1}\in\Gamma$. Further, if $n=0$ then
\begin{equation*}
(\delta_0\lambda)(\gamma)=\lambda\gamma-\lambda
\end{equation*}
for all $\lambda\in\Lambda$ and $\gamma\in\Gamma$. Finally, if $n\leq-1$ then, of course, $\delta_n=0$. The identity $\delta_{n+1}\delta_n=0$ holds true and (\ref{Eq:coc.compl.group}) is thus a cochain complex of abelian groups. The cohomology groups of this complex are then defined in the usual way.

Now let $\Gamma$ be a topological group, $X$ be a topological space and $\Flow$ be a minimal flow with acting homeomorphisms $T_{\gamma}$ ($\gamma\in\Gamma$). Fix an abelian topological group $G$. With operations defined point-wise, the set $\Lambda=C(X,G)$ of all continuous maps $\lambda\colon X\to G$ is an abelian group, on which the group $\Gamma$ acts from the right via $\Lambda\times\Gamma\ni(\lambda,\gamma)\mapsto\lambda T_{\gamma}\in\Lambda$. The corresponding acting transformations are clearly automorphisms of $\Lambda$ and so $\Lambda$ is a right $\Gamma$-module. The cochain complex (\ref{Eq:coc.compl.group}) and the morphisms $\delta_n$ are then defined as above.

Since we are working in the topological category, we want to restrict ourselves to continuous maps. Hence, for $n\geq1$ we consider only those elements $\varphi\in F^n(\Gamma,\Lambda)$, for which the map $\Gamma^n\times X\ni(\gamma_1,\dots,\gamma_n,x)\mapsto(\varphi(\gamma_1,\dots,\gamma_n))(x)\in G$ is continuous. Such maps $\varphi$ form a subgroup of $F^n(\Gamma,\Lambda)$, which is naturally isomorphic to the group $C^n(\mathcal F,G)$ of all continuous maps $\Gamma^n\times X\to G$. This identification leads to a cochain complex obtained by restricting the one in (\ref{Eq:coc.compl.group})
\begin{equation}\label{Eq:cochain.compl}
\dots\stackrel{\delta_{n-2}}{\longrightarrow}C^{n-1}(\mathcal F,G)\stackrel{\delta_{n-1}}{\longrightarrow}C^n(\mathcal F,G)\stackrel{\delta_n}{\longrightarrow}C^{n+1}(\mathcal F,G)\stackrel{\delta_{n+1}}{\longrightarrow}\dots,
\end{equation}
where $C^0(\mathcal F,G)=C(X,G)$, $C^n(\mathcal F,G)=0$ for $n\leq-1$ and the morphisms $\delta_n$ are defined as follows. First, if $n\geq1$ then
\begin{equation*}
\begin{split}
(\delta_n f)(\gamma_1,\dots,\gamma_{n+1},x)&=f(\gamma_2,\dots,\gamma_{n+1},T_{\gamma_1}x)\\
&+\sum_{i=1}^n(-1)^if(\gamma_1,\dots,\gamma_{i-1},\gamma_{i+1}\gamma_i,\gamma_{i+2},\dots,\gamma_{n+1},x)\\
&+(-1)^{n+1}f(\gamma_1,\dots,\gamma_n,x)
\end{split}
\end{equation*}
for all $f\in C^n(\mathcal F,G)$, $\gamma_1,\dots,\gamma_{n+1}\in\Gamma$ and $x\in X$. Further, if $n=0$ then
\begin{equation*}
(\delta_0f)(\gamma,x)=f(T_{\gamma}x)-f(x)
\end{equation*}
for all $f\in C^0(\mathcal F,G)$, $\gamma\in\Gamma$ and $x\in X$. Finally, if $n\leq-1$ then, of course, $\delta_n=0$.

\section{Cohomology groups of a minimal flow}

Given the cochain complex (\ref{Eq:cochain.compl}) and $n\in\mathbb Z$, one can consider the group of \emph{$n$-cocycles} $Z^n(\mathcal F,G)=\ker(\delta_n)$, the group of \emph{$n$-coboundaries} $B^n(\mathcal F,G)=\im(\delta_{n-1})$ and the corresponding $n^{\text{th}}$ cohomology group
\begin{equation*}
H^n(\mathcal F,G)=Z^n(\mathcal F,G)/B^n(\mathcal F,G).
\end{equation*}
A map $f\in C^1(\mathcal F,G)$ is a $1$-cocycle if and only if it satisfies the identity 
\begin{equation}\label{Eq:1-cocycle}
f(\alpha,T_{\beta}x)-f(\alpha\beta,x)+f(\beta,x)=0.
\end{equation}
Moreover, the $1$-coboundaries are the maps $f\in C^1(\mathcal F,G)$ expressible in the form
\begin{equation}\label{Eq:1-coboundary}
f(\gamma,x)=g(T_{\gamma}x)-g(x)
\end{equation}
for an appropriate $g\in C^0(\mathcal F,G)$, the map $g$ being called a \emph{transfer function} of $f$. (Let us mention that a transfer function of a given $1$-coboundary $f$ is not unique in general, but if the flow $\mathcal F$ is minimal then a transfer function of $f$ is unique up to an additive constant.)

We have, obviously, $H^n(\mathcal F,G)=0$ for every $n\leq-1$. The $0^{\text{th}}$ cohomology group $H^0(\mathcal F,G)$ is also determined easily. Indeed, by minimality of $\mathcal F$, the kernel $\ker(\delta_0)$ of $\delta_0$ consists of the constant maps $X\to G$ and hence is isomorphic to $G$. Thus, since $\im(\delta_{-1})=0$,
\begin{equation*}
H^0(\mathcal F,G)=\ker(\delta_0)/\im(\delta_{-1})\cong\ker(\delta_0)\cong G.
\end{equation*}
Naturally, it is much more difficult to determine the higher-dimensional cohomology groups $H^n(\mathcal F,G)$ ($n\geq1$), already the case $n=1$ being highly non-trivial. In fact, it is one of the main aims of our work to determine the first cohomology group $H^1(\mathcal F,G)$ in the case when $G$ is a \emph{compact} (connected) abelian group. One motivation for the study of this ``compact case'' comes from the theory of topological groups (\cite[p.~369]{HofMor}):
\begin{enumerate}
\item[] \emph{``Compact abelian groups form the most important class of abelian topological\\ groups,as is evidenced by the considerable literature on them and their applications.''}
\end{enumerate}
As a matter of fact, compact abelian groups have their honored place in the category of all compact groups (\cite[p.~447]{HofMor}):
\begin{enumerate}
\item[] \emph{``So it emerges that compact abelian groups are not simply examples of compact groups but basic ingredients of the structure of all compact groups.''} 
\end{enumerate}
To illustrate this statement, recall the following results from the structure theory of compact groups.
\begin{itemize}
\item Every compact connected group is a (topological) semi-direct product $G\ltimes H$ of an abelian compact connected group $G$ and a semisimple compact connected group $H$ (see \cite[Theorem~9.39(i), p.~478]{HofMor}).
\item Every compact connected group is a union of its closed connected abelian subgroups (see \cite[Theorem~9.32(ii), p.~473]{HofMor}).
\item The first Betti number $\beta_1(K)$ of a compact connected Lie group $K$ coincides with the topological dimension $\dim(K/K')$ of the quotient group of $K$ modulo the commutator subgroup $K'$ of $K$ and this quotient group is a torus; this follows from the second structure theorem for connected compact Lie groups (see \cite[Theorem~6.41(i), p.~221]{HofMor}) and from the Weyl's theorem (see \cite[Corollary~4, p.~285]{Bou}).
\end{itemize}

Since we shall deal only with the one-dimensional cohomology in this work, let us adopt the following notation for the corresponding groups:
\begin{equation*}
\Cocc(G)=Z^1(\mathcal F,G),\hspace{3mm}\Cob(G)=B^1(\mathcal F,G),\hspace{3mm}\Coch(G)=H^1(\mathcal F,G).
\end{equation*}
Finally, to simplify our notation even more, we shall write
\begin{equation*}
\Cocc=\Cocc(\mathbb T^1),\hspace{3mm}\Cob=\Cob(\mathbb T^1),\hspace{3mm}\Coch=\Coch(\mathbb T^1),
\end{equation*}
where $\mathbb T^1$ stands for the circle group.

\section{Dynamical interpretation of the first cohomology}\label{Sub:dyn.int.cohom}

Let $\Flow$ be a mi\-ni\-mal flow with acting homeomorphisms $T_{\gamma}$ ($\gamma\in\Gamma$), $G$ be an abelian topological group and $\mathcal C\colon\Gamma\times X\to G$ be a continuous map. Given $\gamma\in\Gamma$, consider the map
\begin{equation*}
\widetilde{T}_{\gamma}\colon X\times G\ni(x,g)\mapsto(T_{\gamma}x,\mathcal C(\gamma,x)+g)\in X\times G.
\end{equation*}
Clearly, all $\widetilde{T}_{\gamma}$ are homeomorphisms on $X\times G$. Moreover, they are acting homeomorphisms of a flow $\mathcal F_{\mathcal C}\colon\Gamma\curvearrowright X\times G$ if and only if
\begin{equation*}
\mathcal C(\alpha,T_{\beta}x)+\mathcal C(\beta,x)=\mathcal C(\alpha\beta,x)
\end{equation*}
for all $\alpha,\beta\in\Gamma$ and $x\in X$, that is, if and only if $\mathcal C\in\Cocc(G)=Z^1(\mathcal F,G)$ is a $1$-cocycle of the cochain complex (\ref{Eq:cochain.compl}). If this is the case then we call the flow $\mathcal F_{\mathcal C}$ (or the $1$-cocycle $\mathcal C$) a \emph{group extension of $\mathcal F$}.

Further, if $\mathcal C\in\Cob(G)=B^1(\mathcal F,G)$ is a $1$-coboundary with a transfer function $\xi\colon X\to G$ then the orbit closures of the flow $\mathcal F_{\mathcal C}$ are the vertical translations of the graph of $\xi$ in $X\times G$ and hence they form a decomposition of $X\times G$ into closed $\mathcal F_{\mathcal C}$-invariant sets, on each of which the corresponding restriction of $\mathcal F_{\mathcal C}$ is isomorphic to $\mathcal F$. Therefore, we may think of the flows $\mathcal F_{\mathcal C}$ with $\mathcal C\in\Cob(G)$ as trivial group extensions of $\mathcal F$ and we are led, also from the dynamical point of view, to focus our attention on the first cohomology group
\begin{equation*}
\Coch(G) =\Cocc(G)/\Cob(G) =Z^1(\mathcal F,G)/B^1(\mathcal F,G) =H^1(\mathcal F,G).
\end{equation*}

There is a standard and commonly used approach to the study of group extensions in the topological setting, namely by means of vertical sections of orbit closures; let us recall it briefly (for more details see Section~\ref{S:funct.appr.sect}). Let $\Flow$ be a minimal flow with acting homeomorphisms $T_{\gamma}$ ($\gamma\in\Gamma$), $G$ be a \emph{compact} abelian group and $\mathcal C\in\Cocc(G)$ be a group extension of $\mathcal F$. Choose a base point $z$ for $X$ and let 
\begin{equation*}
F(\mathcal C)=\overline{\mathcal O}_{\mathcal F_{\mathcal C}}(z,0)|_z
\end{equation*}
be the vertical $z$-section of the orbit closure $\overline{\mathcal O}_{\mathcal F_{\mathcal C}}(z,0)$ of $(z,0)$ in $X\times G$ under the action of $\mathcal F_{\mathcal C}$. Then $g\in G$ is an element of $F(\mathcal C)$ if and only if there is a net $(\gamma_i)$ in $\Gamma$ with $T_{\gamma_i}z\to z$ in $X$ and $\mathcal C(\gamma_i,z)\to g$ in $G$. The set $F(\mathcal C)$ is a closed subgroup of $G$ and it does not depend on the choice of a base point $z$ for $X$. Moreover, the group $F(\mathcal C)$ is a cohomology invariant of $\mathcal C$ and it can be used to reconstruct the whole orbit closure of $(z,0)$. Indeed, we have
\begin{equation*}
\overline{\mathcal O}_{\mathcal F_{\mathcal C}}(z,0)=\overline{\bigcup_{\gamma\in\Gamma}T_{\gamma}z\times\left(\mathcal C(\gamma,z)+F(\mathcal C)\right)}
\end{equation*}
and every vertical section of $\overline{\mathcal O}_{\mathcal F_{\mathcal C}}(z,0)$ is an element of the quotient group $G/F(\mathcal C)$. Interestingly, $F$ can be viewed as a covariant functor from the category $\Coc$ of all group extensions of $\mathcal F$ with compact abelian fibre groups into the category $\mathsf{CAGp}$ of all compact abelian groups. It detects coboundaries as well as minimal extensions and is continuous in the sense that it preserves limits of inverse systems. The functor $F$ will be one of the main tools for our study of group extensions in this work.

Despite a formal analogy of minimal systems with the ergodic ones, there are some essential differences between them. One of the most important distinctions lies in the fact that in topological dynamics there is no analogue of the notion of ergodic decomposition from ergodic theory. Though (compact) topological systems from certain classes do admit a decomposition into minimal subsystems (in other words, are point-wise almost periodic or semi-simple), the phase space of a general topological system does not decompose into its minimal subsystems. The class of point-wise almost periodic systems includes, for instance,
\begin{itemize}
\item distal systems (see \cite[Theorem~1]{Ell3}),
\item products of distal and point-wise almost periodic systems (see \cite[Corollary~16, p.~90]{Aus}),
\item distal extensions of point-wise almost periodic systems (see \cite[p.~93]{Aus}).
\end{itemize}
Independently on whether or not the phase space $X$ of a minimal flow $\Flow$ is compact, group extensions of $\mathcal F$ always decompose into minimal subsystems. Indeed, if $G$ is a compact abelian group and $\mathcal C\in\Cocc(G)$ is a group extension of $\mathcal F$ then all the orbit closures $\overline{\mathcal O}_{\mathcal F_{\mathcal C}}(x,g)$ of $\mathcal F_{\mathcal C}$ ($x\in X$, $g\in G$) are minimal sets for $\mathcal F_{\mathcal C}$ with full projections onto $X$, they form a decomposition of the space $X\times G$ and are permuted by the vertical rotations of $X\times G$. In particular, by our discussion from the preceding paragraph, minimality of $\mathcal F_{\mathcal C}$ reduces to the equality $F(\mathcal C)=G$.

\section{Minimal group extensions}

Following our setting from Section~\ref{S:coc.compl.min.flw}, we are interested in group extensions of a minimal flow $\mathcal F$ with \emph{abelian} fibre groups $G$. In the theory of dynamical systems it is customary (and natural from various viewpoints) to restrict to locally compact groups $G$ (although extensions with other types of fibre groups $G$, Polish in particular, are also of an importance). In topological dynamics there is a considerable difference between the compact and the non-compact case---different phenomena occur in these two situations and different techniques are used in their study. One of these phenomena is the (non-)existence of minimal group extensions; to illustrate this statement, let us recall the following results.
\begin{itemize}
\item If $\mathcal F$ is a minimal $\mathbb Z$-flow on a compact metrizable space then the group extension $\mathcal F_{\mathcal C}$ of $\mathcal F$ is non-minimal for every $\mathcal C\in\Cocc(\mathbb R)$; this follows basically from \cite[Section~1]{Bes} (see also \cite[Section~2]{MenSie}).
\item If $\mathcal F$ is a free minimal $\mathbb Z$-flow (or $\mathbb R$-flow) on a compact metrizable space and $G$ is a compact connected second countable abelian group then there is $\mathcal C\in\Cocc(G)$ such that the group extension $\mathcal F_{\mathcal C}$ of $\mathcal F$ is minimal (see, e.g., \cite{Fur1}, \cite{Ell1}, \cite{Par}, \cite{JonPar}, \cite{KeyNew}).
\end{itemize}

The study of minimal group extensions of a given minimal flow $\mathcal F$ is motivated also by the following problem, which is one of the central problems in the topological theory of dynamical systems.
\begin{enumerate}
\item[($\mathcal P$)] \emph{Given a locally compact group $\Gamma$, describe the class $\mathcal M_{\Gamma}$ of all compact (metrizable) spaces $Y$, on which the group $\Gamma$ acts in a (free) minimal way.}
\end{enumerate}
If our space $Y$ under consideration has the form of a direct product $Y=X\times Z$ with $X\in\mathcal M_{\Gamma}$ then one can try to construct a minimal $\Gamma$-flow on $Y$ in the form of a skew product over a minimal $\Gamma$-flow on $X$. This approach turns out to be useful for two reasons. First, because of their special form, skew products on $Y$ are easier to handle and understand than general flows on $Y$. Second, the space $Y$ often admits a minimal $\Gamma$-flow in the form of a skew product already under mild assumptions on the fibre $Z$; to illustrate this statement, let us recall the following results.
\begin{itemize}
\item If a compact metrizable space $Z$ admits a minimal action of a Polish group with a dense arc-wise connected subgroup then every minimal $\mathbb Z$-flow $\mathcal F$ on an infinite compact metrizable space can be extended to a minimal skew product with the fibre $Z$ (see \cite[Section~1]{GlaWei}).
\item The above mentioned result holds true also for free minimal actions of locally compact non-compact second countable amenable groups (see \cite[Section~4]{Dir3}).
\end{itemize}
These and similar results can help us gain some insight into the structure of the classes $\mathcal M_{\Gamma}$.
\begin{itemize}
\item If $\Gamma$ is a locally compact non-compact second countable amenable group then the class $\mathcal M_{\Gamma}$ is closed with respect to direct products with the following spaces: compact connected manifolds without boundary, compact connected Hilbert cube manifolds, homogeneous spaces of compact connected metrizable groups, compact metrizable spaces admitting a minimal homeomorphism isotopic to the identity, and others (see \cite[Section~4]{Dir3} and references therein).
\item If $\Gamma$ is a locally compact non-compact second countable amenable group possessing a dense arc-wise connected subgroup (say, an amenable connected non-compact Lie group) then the class $\mathcal M_{\Gamma}$ is closed with respect to (at most) countable direct products (see \cite[Section~5]{Dir3}).
\end{itemize}

Now, if our fibre is a compact abelian group $Z=G$ then it is natural to exploit its additional algebraic structure and to try to construct a minimal skew product on $Y=X\times G$ in the form of a group extension of $\mathcal F$. For flows with amenable acting groups $\Gamma$, the existence of minimal group extensions follows from \cite[Theorem~8]{Dir3}. As the amenability assumption is used heavily in the proof of the mentioned result, it is unlikely that the techniques used therein can be modified in such a way that they apply also to flows with non-amenable acting groups. Therefore, in the present work we aim at finding techniques that would yield minimal group extensions of $\mathcal F$ regardless of whether or not the acting group $\Gamma$ of $\mathcal F$ is amenable.

\section{Minimal flows and one-dimensional invariants}
Let $X$ be a (Hausdorff) topological space and $\Gamma$ be an abelian topological group. If $\Flow$ is a minimal flow with acting homeomorphisms $T_{\gamma}$ ($\gamma\in\Gamma$) and $T_{\gamma}(x)=x$ for some $\gamma\in\Gamma$ and $x\in X$ then $T_{\gamma}$ is the identity on $X$, because the set of fixed points of $T_{\gamma}$ is closed and, by commutativity of $\Gamma$, also $\mathcal F$-invariant; in other words, if $\Gamma$ acts on $X$ in an effective and minimal way then the action is free (see, e.g., \cite[p.~23]{Aus}). Thus, an elementary obstruction to the exis\-ten\-ce of a minimal $\Gamma$-flow on $X$ is a fixed point property of $X$. Since the fixed point property of a space is often detected by algebraic-topological methods, one can try to employ techniques of algebraic topology to the study of minimal flows, particularly to their (non-)exis\-tence on given spaces. A famous sufficient condition for the existence of a fixed point is the Lefschetz theorem (see, e.g., \cite[Theorem~2C.3, p.~179--181]{Hat}), which states that a continuous map on a finite polyhedron $X$ with non-zero Euler characteristic $\chi(X)$ possesses a fixed point, provided it is homotopic to the identity. This theorem has the following (well known) immediate implication on minimality: if an abelian arc-wise connected group $\Gamma$ acts in a minimal way on a (non-trivial) finite polyhedron $X$ then $\chi(X)=0$.

It was believed for some time (see, e.g., \cite{ChuGer}, \cite{Haj}, \cite[Problems, p.~515]{AusGot}) that the fundamental group $\pi_1(X)$ of a compact connected manifold $X$ can carry some information about its (non-)minimality. Specifically, it was conjectured that a simply connected compact manifold $X$ can not support a minimal continuous flow, that is, a minimal action $\mathbb R\curvearrowright X$. This conjecture was partially supported by a theorem of Furstenberg (see \cite[Theorem~11.1]{Fur2}), who showed that a simply connected compact metric space can not admit a \emph{distal} minimal action of an abelian locally compact group. However, the conjecture turned out to be false. Indeed, for every smooth compact connected manifold without boundary $Y$, the product $X=\mathbb S^3\times\mathbb S^3\times Y$ supports a (smooth) minimal flow $\mathbb R\curvearrowright X$; this follows from a result of Fathi and Herman (see \cite[Th\'eor\`eme~2]{FatHer}). Now, since a smooth compact connected manifold without boundary $Y$ may have an arbitrary finitely generated group as its fundamental group $\pi_1(Y)$ already in dimension $4$, the fundamental group $\pi_1(X)$ of a compact manifold $X$ admitting a minimal flow $\mathbb R\curvearrowright X$ can be an arbitrary finitely generated group (in particular, it can be trivial) for instance in dimension $10$.

The situation does not improve if we take general continua $X$ and their first cohomotopy (i.e., first \v Cech cohomology) groups $\pi^1(X)$ into consideration: given an arbitrary torsion-free abelian group $A$ with $\rank(A)\leq\mathfrak{c}$, there is a continuum $X$ having $A$ as its first cohomotopy group, $\pi^1(X)\cong A$, and supporting a minimal flow $\mathbb R\curvearrowright X$. Indeed, it is sufficient to take $X=A^*$, the Pontryagin dual of the discrete group $A$. Clearly, the group $X$ has a dense one-parameter subgroup $q\colon\mathbb R\to X$, which gives rise to a minimal flow $\mathcal F\colon\mathbb R\times X\ni(t,x)\mapsto q(t)x\in X$. Moreover, there are isomorphisms of groups $\pi^1(X)\cong X^*\cong A^{**}\cong A$ (see \cite[Theorem~8.57, p.~420]{HofMor} or our Subsection~\ref{Sub:chmtp.cpt.gps}).

Now let $\Flow$ be a minimal flow, $G$ be a compact abelian group and $\mathcal C\in\Cocc(G)$ be a group extension of $\mathcal F$. By commutativity of the group $G$, all its irreducible unitary representations are of dimension one and so the one-dimensional representations of $G$ (that is, the characters $\chi\in G^*$) carry all the information about $G$. Since the flow $\mathcal F_{\mathcal C}$ induced by $\mathcal C$ on $X\times G$ respects the group structure of $G$ to a large extent, it seems reasonable to expect that there might exist deeper links between the dynamics of $\mathcal F_{\mathcal C}$ on one side and the structure of the \emph{one-dimensional} algebraic-topological invariants of the space $X$ on the other side. We are looking for such links in the present work, concentrating on flows with the following two types of phase spaces $X$.
\begin{itemize}
\item The space $X$ is a compact connected manifold. In this case we would like to find a relation of the dynamics of group extensions of $\mathcal F$ to the fundamental group $\pi_1(X)$ (or the first homology group $H_1(X)$) of $X$.
\item The space $X$ is an arbitrary (not necessarily locally connected) compact connected space. In this case we are searching for links between the dynamics of group extensions of $\mathcal F$ and the first cohomotopy (that is, the first \v Cech cohomology) group $\pi^1(X)$ of $X$.
\end{itemize}

Let us describe one situation in support of our conjecture from the previous paragraph. Let $\Gamma$ be a connected (Lie) group, $X$ be a compact manifold and $\Flow$ be a minimal flow. Fix a compact abelian group $G$ and a \emph{finite} subgroup $H$ of $G$. Given the existence of a group extension $\mathcal C\in\Cocc(G)$ of $\mathcal F$ with $F(\mathcal C)=H$, can we find a relation between $H$ and $\pi_1(X)$? To this end, choose a base point $z$ for $X$ and write $M=\overline{\mathcal O}_{\mathcal F_{\mathcal C}}(z,0)$ for the orbit closure of $(z,0)$ in $X\times G$ under the action of $\mathcal F_{\mathcal C}$. By our discussion from Section~\ref{Sub:dyn.int.cohom}, $M$ is a minimal set for $\mathcal F_{\mathcal C}$, hence it is connected by connectedness of $\Gamma$ and compact by compactness of $X$ and $G$. Further, the group $H$ acts freely on $M$ by vertical rotations and the corresponding homogeneous space $M/H$ is homeomorphic to $X$. Thus, it follows that $H$ is a quotient group of the fundamental group $\pi_1(X)$ of $X$ and, by commutativity of $H$, also of the first homology group $H_1(X)$ of $X$ (see, e.g., \cite[Chapter~13, \S~81]{Mun}). In particular, if the group $H_1(X)$ is torsion-free then the number of elementary divisors of $H$ does not exceed the Betti number (that is, the rank) of $H_1(X)$.

\section{Free group extensions}\label{S:Free.gp.ext}

Let $\mathsf{Cat}$ be a category and $K$ be an object of $\mathsf{Cat}$. We shall call $K$ a \emph{free object} for $\mathsf{Cat}$ if the set $\Hom(K,L)$ of morphisms $K\to L$ is a singleton for every object $L$ of $\mathsf{Cat}$. Clearly, if a free object of $\mathsf{Cat}$ exists then it is unique up to an isomorphism.

Now let $\mathsf{Cat}$ be a category, $\mathsf{Tp}$ be a subcategory of the category $\mathsf{Top}$ of topological spaces and $\phi\colon\mathsf{Cat}\to\mathsf{Tp}$ be a covariant functor. Given a topological space $X$ from $\mathsf{Tp}$, consider the category $\mathsf{Map}_{\phi}(X,\mathsf{Cat})$ defined as follows:
\begin{enumerate}
\item[$\centerdot$] the objects of $\mathsf{Map}_{\phi}(X,\mathsf{Cat})$ are the pairs $(B,\psi)$, where $B$ is an object of $\mathsf{Cat}$ and $\psi\colon X\to\phi(B)$ is a morphism in $\mathsf{Tp}$,
\item[$\centerdot$] given objects $(B,\psi)$, $(C,\chi)$ in $\mathsf{Map}_{\phi}(X,\mathsf{Cat})$, a morphism $k\colon(B,\psi)\to(C,\chi)$ is a morphism $k\in\Hom(B,C)$ in $\mathsf{Cat}$ with $\phi(k)\psi=\chi$,
\item[$\centerdot$] the composition of two morphisms in $\mathsf{Map}_{\phi}(X,\mathsf{Cat})$ is defined as their composition in $\mathsf{Cat}$ and the identity on $(B,\psi)$ is the identity on $B$ in $\mathsf{Cat}$.
\end{enumerate}
In this way $\mathsf{Map}_{\phi}(X,\mathsf{Cat})$ becomes a category in its own right. A free object for the category $\mathsf{Map}_{\phi}(X,\mathsf{Cat})$ consists of an object $A$ of $\mathsf{Cat}$ and a morphism $\varphi\colon X\to\phi(A)$ in $\mathsf{Tp}$ such that for every object $B$ of $\mathsf{Cat}$ and every morphism $\psi\colon X\to\phi(B)$ in $\mathsf{Tp}$ there is a unique morphism $h\in\Hom(A,B)$ in $\mathsf{Cat}$ with $\phi(h)\varphi=\psi$; see Figure~\ref{Fig:free.obj.catg}. The same construction works if $\mathsf{Tp}$ is a subcategory of the category $\mathsf{PtTop}$ of pointed topological spaces and base point preserving maps; in this situation we assume that every topological group carries its identity as the base point.
\begin{figure}[ht]
\[\minCDarrowwidth20pt\begin{CD}
X @>\varphi>> \phi(A) @. \hspace{10mm}A\\
@| @VV\phi(h)V \hspace{10mm}@VVhV\\
X @>>\psi> \phi(B) @. \hspace{10mm}B
\end{CD}\]
\caption{Free object $(A,\varphi)$ for the category $\mathsf{Map}_{\phi}(X,\mathsf{Cat})$}
\label{Fig:free.obj.catg}
\end{figure}

Let us recall some well-known special cases of the construction described above.
\begin{itemize}
\item \emph{The Stone-\v Cech compactification of a space $X$.} In this case $\mathsf{Cat}=\mathsf{CpTop}$ is the category of compact topological spaces, $\mathsf{Tp}=\mathsf{Top}$ is the category of topological spaces, $\phi$ is the inclusion functor $\mathsf{CpTop}\to\mathsf{Top}$ and $X$ is a Tychonov space (see, e.g., \cite[pp.~152--154]{Kel}).
\item \emph{The Bohr compactification of a group $X$.} In this case $\mathsf{Cat}=\mathsf{CAGp}$ is the category of compact abelian groups, $\mathsf{Tp}=\mathsf{LCAGp}$ is the category of locally compact abelian groups, $\phi$ is the inclusion functor $\mathsf{CAGp}\to\mathsf{LCAGp}$ and $X$ is a locally compact abelian group (see, e.g., \cite[Section~4.7]{Fol}). (Some properties of the Bohr compactification are summarized also in our Subsection~\ref{Sub:Bohr.compctf}.)
\item \emph{The free compact abelian group over a space $X$.} In this case $\mathsf{Cat}=\mathsf{CAGp}$ is the ca\-te\-go\-ry of compact abelian groups, $\mathsf{Tp}=\mathsf{PtTop}$ is the category of pointed topological spaces, $\phi$ is the functor $\mathsf{CAGp}\to\mathsf{PtTop}$ forgetting the algebraic structure and $X$ is a pointed topological space (see \cite[Chapter~8, p.~415]{HofMor}). (Some properties of free compact abelian groups are summarized also in our Subsection~\ref{Sub:free.cpt.abel.gp}.)
\item \emph{The free compact group over a space $X$.} In this case $\mathsf{Cat}=\mathsf{CGp}$ is the category of compact groups, $\mathsf{Tp}=\mathsf{PtTop}$ is the category of pointed topological spaces, $\phi$ is the functor $\mathsf{CGp}\to\mathsf{PtTop}$ forgetting the algebraic structure and $X$ is a pointed topological space (see \cite[Chapter~11, pp.~577--578]{HofMor}).
\end{itemize}

Now let $\Flow$ be a minimal flow, $\Gimel$ be a compact abelian group and $\Daleth\in\Cocc(\Gimel)$ be a group extension of $\mathcal F$. Following the ideas from the preceding paragraphs, we call $(\Gimel,\Daleth)$ a \emph{free group extension of $\mathcal F$} if for every compact abelian group $G$ and every group extension $\mathcal C\in\Cocc(G)$ of $\mathcal F$ there is a unique topological morphism $h_{\mathcal C}\colon\Gimel\to G$ with $h_{\mathcal C}\Daleth=\mathcal C$; see Figure~\ref{Fig:free.gp.ext}.
\begin{figure}[ht]
\[\minCDarrowwidth20pt\begin{CD}
\Gamma\times X @>\Daleth>> \Gimel @. \hspace{10mm}\Gimel\\
@| @VVh_{\mathcal C}V \hspace{10mm}@VVh_{\mathcal C}V\\
\Gamma\times X @>>\mathcal C> G @. \hspace{10mm}G
\end{CD}\]
\caption{Free group extension $(\Gimel,\Daleth)$ of $\mathcal F$}
\label{Fig:free.gp.ext}
\end{figure}
One of our aims in this work is to construct a free group extension of an arbitrary minimal flow $\mathcal F$. Intuitively, a free group extension of $\mathcal F$ should contain in itself encoded information about all group extensions of $\mathcal F$ and we therefore expect to be able to express the most important properties of group extensions of $\mathcal F$ in terms of the free group extension. In particular, we wish to find an explicit relation between $\Gimel$ and $F(\Daleth)$ on one side and the groups $\Cocc$, $\Cob$ and $\Coch$ on the other side, and express the group $\Gimel$ in terms of familiar compact abelian groups.

Important categories arising in topological dynamics fail to possess free objects, but turn out to possess objects universal in a weaker sense. To give an example, let $\Gamma$ be a topological group and $\mathsf{MinFlow}_{\Gamma}$ be the category of all minimal $\Gamma$-flows with compact phase spaces; notice that by minimality and compactness assumptions, every morphism in $\mathsf{MinFlow}_{\Gamma}$ is an epimorphism. By a theorem of Ellis from \cite{Ell2} (see also \cite[Proposition~7.13, p.~62]{Ell}), the category $\mathsf{MinFlow}_{\Gamma}$ possesses an object $\mathcal F_u$ such that the set of morphisms $\Hom(\mathcal F_u,\mathcal F)$ is non-empty for every object $\mathcal F$ in $\mathsf{MinFlow}_{\Gamma}$. Traditionally, $\mathcal F_u$ is referred to as a \emph{universal minimal $\Gamma$-flow}; it is unique up to an isomorphism (see \cite[Corollary~7.16, p.~62]{Ell}). (For many other examples of categories in topological dynamics that possess a universal object, we refer to \cite{Gla1} and \cite[Chapter~7]{Ell}.)

A universal minimal $\Gamma$-flow $\mathcal F_u$ is not a free object for the category $\mathsf{MinFlow}_{\Gamma}$ as defined at the beginning of this section. Indeed, the set of morphisms $\Hom(\mathcal F_u,\mathcal F_u)$, which coincides with the group $\text{Aut}(\mathcal F_u)$ of automorphisms of $\mathcal F_u$ by coalescence of $\mathcal F_u$ (see \cite[pp.~115--116]{Aus}), is far from being a singleton. On the contrary, if $(x,y)$ is an arbitrary element of a minimal set for the product $\Gamma$-flow $\mathcal F_u\times\mathcal F_u$ then there is $h\in\text{Aut}(\mathcal F_u)$ with $h(x)=y$; the converse statement is also true: if $h\in\text{Aut}(\mathcal F_u)$ then $(x,h(x))$ is necessarily an element of a minimal set for $\mathcal F_u\times\mathcal F_u$ for every $x$ from the base space of $\mathcal F_u$ (indeed, since the map $z\mapsto(z,h(z))$ is an (injective) morphism of flows $\mathcal F_u\to\mathcal F_u\times\mathcal F_u$, its image, that is, the graph of $h$, is a minimal set of $\mathcal F_u\times\mathcal F_u$). Thus, since $\mathcal F_u$ is the only possible candidate for a free object of $\mathsf{MinFlow}_{\Gamma}$, it follows that $\mathsf{MinFlow}_{\Gamma}$ fails to possess a free object in general.

Now, given a minimal flow $\Flow$, we are interested in the existence of a universal object for the category of all \emph{minimal} group extensions of $\mathcal F$. The universality is now understood in the sense that every minimal group extension of $\mathcal F$ is an epimorphic image of the universal one. To be more precise, let $\mGimel$ be a compact abelian group and $\mDaleth\in\Cocc(\mGimel)$ be a minimal group extension of $\mathcal F$. We call $(\mGimel,\mDaleth)$ a \emph{universal minimal group extension of $\mathcal F$} if for every compact abelian group $G$ and every minimal group extension $\mathcal C\in\Cocc(G)$ of $\mathcal F$ there is a topological morphism $q\colon\mGimel\to G$ with $q\mDaleth=\mathcal C$; see Figure~\ref{Fig:free.min.gp.ext}.
\begin{figure}[ht]
\[\minCDarrowwidth20pt\begin{CD}
\Gamma\times X @>\mDaleth>> \mGimel @. \hspace{10mm}\mGimel\\
@| @VVqV \hspace{10mm}@VVqV\\
\Gamma\times X @>>\mathcal C> G @. \hspace{10mm}G
\end{CD}\]
\caption{Universal minimal group extension $(\mGimel,\mDaleth)$ of $\mathcal F$}
\label{Fig:free.min.gp.ext}
\end{figure}
Notice that, by minimality and compactness assumptions, such a morphism $q$ is necessarily an epimorphism of groups. Even more importantly, such a morphism $q$ is unique (this follows, for instance, from Lemma~\ref{L:repr.min.ext.gp}(i) and Remark~\ref{R:repr.min.ext.gp} in Section~\ref{Sub:free.min.ext}); consequently, if a universal minimal group extension $(\mGimel,\mDaleth)$ of $\mathcal F$ does exist then it is automatically a \emph{free} object for the category of all minimal group extensions of $\mathcal F$. One of our aims in this work is to find out whether or not this universal/free minimal group extension of $\mathcal F$ exists. Contrary to the situation of all group extensions of $\mathcal F$, this turns out not to be the case.

Finally, we wish to mention that a notion of a universal minimal group extension of a minimal flow $\Flow$ has been defined and studied by Glasner (see \cite[Section~8.1]{Gla1}). This concept is, however, not to be confused with the one used in this work. In fact, while a universal minimal group extension in the sense of Glasner exists for every minimal flow $\mathcal F$ (and coincides with the universal minimal almost periodic extension of $\mathcal F$), a universal minimal group extension in our sense (which has been defined in analogy to the notion of a free group extension of $\mathcal F$) typically does not exist. On the other hand, if $\mathcal M$ is an arbitrary minimal subflow of $\mathcal F_{\Daleth}$ (that is, if $\mathcal M$ is a restriction of $\mathcal F_{\Daleth}$ onto one of its orbit closures) then every minimal group extension $\mathcal F_{\mathcal C}$ of $\mathcal F$ is an epimorphic image of $\mathcal M$; indeed, if $h_{\mathcal C}$ is the unique morphism $\Daleth\to\mathcal C$ then $\id_X\times h_{\mathcal C}\colon\mathcal M\to\mathcal F_{\mathcal C}$ is an epimorphism of flows. (Let us mention that $\mathcal M$ is not a group extension of $\mathcal F$ in the sense used in this work, though it is a group extension of $\mathcal F$ in the sense used in \cite{Gla1}.)

\section{Torsions in cohomology groups}

Let $A$ be an abelian group and $\tor(A)$ be the set of all elements of $A$ of a finite order. Then $\tor(A)$ forms a pure subgroup of $A$, traditionally called the \emph{torsion subgroup of $A$}, with a torsion-free quotient group $C=A/\tor(A)$. The torsion subgroup $\tor(A)$ of $A$ is one of the most important objects associated to $A$, since (\cite[p.~35]{Fuchs1})
\begin{enumerate}
\item[] ``\emph{... the structure theory of abelian groups splits into the theories of torsion and torsion-free groups, and investigations of how these are glued together to form mixed groups.}''
\end{enumerate}
To illustrate this, recall that there is an isomorphism $A\cong C\times_{\varphi}\tor(A)$, where $\varphi\in Z^2(C,\tor(A))$ is an appropriate symmetric $2$-cocycle (see \cite[Chapter~9]{Fuchs1} or our Subsections~\ref{Sub:2-coc.2-cob} and~\ref{Sub:ext.ab.gps}). The simplest situation occurs when $\varphi$ is a $2$-coboundary, $\varphi\in B^2(C,\tor(A))$, in which case $\tor(A)$ is a direct summand in $A$, that is, $A=\tor(A)\oplus C$. This is the case, for instance, in the following situations.
\begin{itemize}
\item The group $A$ is divisible. In this case $\tor(A)$ is also a divisible group and so it is a direct summand in $A$ by a theorem of Baer (see \cite[Theorem~21.2, p.~100]{Fuchs1}).
\item The group $A$ is finitely generated. In this case $\tor(A)$ is a direct sum of finitely many finite cyclic groups and it is a direct summand in $A$ by the fundamental theorem on finitely generated abelian groups (see, e.g., \cite[Theorem~15.5, p.~79]{Fuchs1}).
\item The torsion subgroup $\tor(A)$ of $A$ is of a bounded order. In this case $\tor(A)$ is a direct summand in $A$ by a theorem of Kulikov, which asserts that a pure subgroup of a bounded order is necessarily a direct summand (see \cite[Theorem~27.5, p.~118]{Fuchs1}).
\end{itemize}

When abelian groups are assigned to other mathematical objects, the torsion subgroups of these groups are usually of more than merely a group-theoretical interest. Moreover, the torsion elements of these groups often correspond to elements in the studied objects with some special features. Let us give some examples of this phenomenon.
\begin{itemize}
\item If $G$ is a compact abelian group then the torsion subgroup $\tor(G^*)$ of the Pontryagin dual $G^*$ of $G$ coincides with the annihilator $(G_0)^{\perp}$ of the identity component $G_0$ of $G$ in $G^*$; in particular, the group $G$ is connected (respectively, totally disconnected) if and only if the group $G^*$ is torsion-free (respectively, a torsion group) (see, e.g., \cite[Theorem~8.4, p.~371]{HofMor}).
\item A compact connected Lie group $G$ is semi-simple if and only if its fundamental group $\pi_1(G)$ (that is, its first homology group $H_1(G)$) is a torsion group (see \cite[Corollary~4, p.~285]{Bou}). Moreover, $\tor(\pi_1(G))$ is isomorphic to a finite central subgroup of the universal cover $\widetilde{G}$ of $G$.
\item A closed connected surface $X$ is orientable if and only if the torsion subgroup $\tor(H_1(X))$ of the first homology group $H_1(X)$ of $X$ vanishes (see, e.g., \cite[Chapter~12]{Mun}).
\end{itemize}

In this context we also wish to mention a recent work of Glasner and Host (see \cite{GlaHos}) in which a remarkable link was found between torsions and group extensions in the theory of dimension groups of minimal $\mathbb Z$-flows on the Cantor space. (For necessary definitions and notation we refer the reader to \cite{GlaHos}. We also mention that the notion of a group extension in \cite{GlaHos} is understood in the sense of \cite{Gla1}; see our discussion towards the end of Section~\ref{S:Free.gp.ext}.) 
\begin{itemize}
\item Let $\mathcal F_1,\mathcal F_2$ be minimal $\mathbb Z$-flows on the Cantor space and let $p\colon\mathcal F_1\to\mathcal F_2$ be an epimorphism. Denote by $p^*\colon K_0(\mathcal F_2)\to K_0(\mathcal F_1)$ the induced morphism between the associated dimension groups of $\mathcal F_2$ and $\mathcal F_1$, respectively. Then, as was shown in \cite{GlaHos}, there is a one-to-one correspondence between elements of the quotient group $K_0(\mathcal F_1)/p^*K_0(\mathcal F_2)$ and pairs of epimorphisms $\mathcal F_1\stackrel{r}{\longrightarrow}\mathcal F\stackrel{q}{\longrightarrow}\mathcal F_2$ such that $p=q\circ r$ and $q$ is a group extension by a finite abelian group. Moreover, if $H=\tor(K_0(\mathcal F_1)/p^*K_0(\mathcal F_2))$ is the (discrete) torsion group of $K_0(\mathcal F_1)/p^*K_0(\mathcal F_2)$ and $K=H^*$ is the Pontryagin dual of $H$ then there is a pair of epimorphisms $\mathcal F_1\stackrel{r_u}{\longrightarrow}\mathcal F_u\stackrel{q_u}{\longrightarrow}\mathcal F_2$ such that $p=q_u\circ r_u$ and $q_u$ is a group extension by $K$. For further properties of $q_u$ and $\mathcal F_u$ and for other related results, we refer to \cite[Theorems~2.1,~2.3]{GlaHos}.
\end{itemize}

Now let $\Flow$ be a minimal flow. In the present work we are searching for a dynamical interpretation of the torsion elements in the cohomology group $\Coch$ of $\mathcal F$, as well as for a relation between the structure of the torsion subgroup $\tor(\Coch)$ of $\Coch$ on one side and the intrinsic properties of $\mathcal F$ on the other side. We are also interested in sufficient conditions, under which the group $\tor(\Coch)$ is a direct summand in $\Coch$. Since this is true, for instance, if the group $\Coch$ is divisible, we are naturally led also to the study of divisibility of $\Coch$. Recall that every quotient group of a divisible group is again divisible and so it suffices to find sufficient conditions for divisibility of the group $\Cocc$.

\section{The main aims of the present work}

Let $\Flow$ be a minimal flow. Being motivated by our discussion from the preceding sections, we set ourselves the following aims in the present work.
\begin{enumerate}
\item[(1)] \emph{Given a compact abelian group $G$, express the inclusion $\Cob(G)\subseteq\Cocc(G)$ in terms of familiar (or elementary) abelian groups and compute the corresponding first cohomology group $\Coch(G)$.}
\end{enumerate}
We accomplish this task in the following two situations.
\begin{itemize}
\item The acting group $\Gamma$ of $\mathcal F$ is a simply connected Lie group and the phase space $X$ of $\mathcal F$ is a compact second countable space with a non-trivial first cohomotopy group $\pi^1(X)$; see Theorem~\ref{T:F.in.Gim.CLAC} in Section~\ref{S:circ.case} and Theorem~\ref{T:Co.in.Cocc.G} in Section~\ref{S:gen.fib} (other results on $\Coch(G)$ in this situation are obtained in Theorems~\ref{T:exact.sequence} and~\ref{P:tor.gen.G} in Section~\ref{S:coh.cls.gp.mrph}).
\item The acting group $\Gamma$ of $\mathcal F$ is a connected Lie group, the phase space $X$ of $\mathcal F$ is a compact manifold, the flow $\mathcal F$ is topologically free and possesses a free cycle (the last two properties are defined in Subsection~\ref{Sub:top.frnss.fr.ccl}); see Theorem~\ref{T:F.in.Gim.free.ccl} in Section~\ref{S:circ.case} and Theorem~\ref{C:GpExt.frccl.gen.G} in Section~\ref{S:gen.fib} (other results on $\Coch(G)$ in this situation are obtained in Theorems~\ref{T:cn.td.LG.mf} and~\ref{P:CochG.closed.Y.N} in Section~\ref{S:fur.struct.res}).
\end{itemize}
In the first of these two situations, the cohomology group $\Coch(G)$ of $\mathcal F$ can be expressed in terms of elementary abelian groups, the Pontryagin dual $G^*$ of $G$ and the first cohomotopy group $\pi^1(X)$ of $X$. In the second situation we express the group $\Coch(G)$ in terms of elementary abelian groups, the Pontryagin dual $G^*$ of $G$, the first weak homology group $H_1^w(X)$ of $X$ and the group $H_1^w(\mathcal F)$ (the latter group is defined in Subsection~\ref{Sub:top.frnss.fr.ccl}).
\medskip
\begin{enumerate}
\item[(2)] \emph{Given a compact abelian group $G$, find sufficient conditions for the exis\-ten\-ce of minimal extensions in $\Cocc(G)$. More generally, determine the existence of large abelian subgroups within the groupoid $\Cocm(G)$.} (The groupoid $\Cocm(G)$ is formed by the minimal extensions from $\Cocc(G)$ together with the identity of $\Cocc(G)$, see Subsection~\ref{Sub:gpd.min.ext}.)
\end{enumerate}
In solving this problem we concentrate on the following three situations.
\begin{itemize}
\item The group $\Gamma$ is locally compact, non-compact, second countable, amenable, the space $X$ is compact second countable and the flow $\mathcal F$ possesses a free point; see Theorem~\ref{T:grp.min.ext.am} in Section~\ref{Sub:gp.min.ext}.
\item The group $\Gamma$ is simply connected and the space $X$ is compact second countable with $\pi^1(X)\neq0$; see Theorem~\ref{T:grp.min.ext.simp} in Section~\ref{Sub:gp.min.ext}.
\item The group $\Gamma$ is a connected Lie group, the space $X$ is a compact manifold and the flow $\mathcal F$ possesses a free cycle; see Theorem~\ref{T:gp.min.free.ccl} in Section~\ref{Sub:gp.min.ext}.
\end{itemize}
In all these situations we show that for every (non-trivial) compact connected abelian group $G$ with weight at most $\mathfrak{c}=2^{\aleph_0}$, the groupoid $\Cocm(G)$ contains an isomorphic copy of the additive group of real numbers $\mathbb R$ (and hence also an isomorphic copy of every torsion-free abelian group with rank at most $\mathfrak{c}$). Moreover, under some mild additional assumptions in the last two situations, we give necessary and sufficient conditions for a torsion-free abelian group $R$ to have an isomorphic copy within the groupoid $\Cocm(G)$ that would form a dense subset of $\Div(\Cocc(G))$ or $\Cocc(G)$; see Theorem~\ref{T:dns.sbgp.min.div} in Section~\ref{S:top.min.ext}.
\medskip
\begin{enumerate}
\item[(3)] \emph{Find relations of dynamical properties of group extensions of $\mathcal F$ to one-dimensional algebraic-topological invariants (the fundamental group $\pi_1$ and the first cohomotopy group $\pi^1$ in particular) and to Pontryagin duality.}
\end{enumerate}
In our study of this problem we focus our attention on the following two situations.
\begin{itemize}
\item The group $\Gamma$ is connected and the space $X$ is a (compact) manifold; see Theorem~\ref{T:F.and.pi} and Corollary~\ref{C:F.and.pi.fprt} in Section~\ref{S:rel.mndrm.act}.
\item The group $\Gamma$ is connected and the space $X$ is compact; see Theorems~\ref{T:F.and.Cech} and~\ref{T:F.and.Cech.td} in Section~\ref{Sub:F.and.Cech}.
\end{itemize}
In the first of these two situations we show that if the section $F(\mathcal C)$ of a group extension $\mathcal C\in\Cocc(G)$ of $\mathcal F$ is totally disconnected (in particular, if it is finite) then it can be computed by using the monodromy action of the fundamental group of an appropriate quotient group of $G$. In the second situation it is shown that a totally disconnected section $F(\mathcal C)$ can be computed by using the first cohomotopy functor $\pi^1$ and the Pontryagin duality. These results enable us to determine, whether or not a given closed subgroup $K$ of a compact abelian group $G$ takes the form of a section $K=F(\mathcal C)$ of some group extension $\mathcal C\in\Cocc(G)$ of $\mathcal F$; in the first situation we find necessary and sufficient conditions in Sections~\ref{S:ex.prescr.sect.1} and~\ref{S:ex.prscr.sect.2}, in the second situation we accomplish the same task in Section~\ref{S:ex.prescr.sect.3}.
\medskip
\begin{enumerate}
\item[(4)] \emph{Construct a free object for the category $\Coc$ of all group extensions of $\mathcal F$ with compact abelian fibre groups, gain some information about its dynamical properties and express its range in terms of familiar compact abelian groups.}
\end{enumerate}
Basic properties of a free object for $\Coc$ are collected in Section~\ref{Sub:free.ext}. In Theorem~\ref{T:exist.free.ext} a free object for $\Coc$ is constructed - it consists of a compact abelian group $\Gimel$ and an extension $\Daleth\in\Cocc(\Gimel)$; the freeness of $\Daleth$ lies in the fact that for every object $\mathcal C$ of $\Coc$ there is a unique morphism $h_{\mathcal C}\colon\Daleth\to\mathcal C$. In the rest of Section~\ref{Sub:free.ext} we
\begin{enumerate}
\item[$\centerdot$] compute the section $F(\Daleth)$ of $\Daleth$,
\item[$\centerdot$] find a relation between $\Daleth$ and $F(\Daleth)$ on one side and a free compact abelian group over the base space $X$ of $\mathcal F$ on the other side,
\item[$\centerdot$] express selected properties of group extensions $\mathcal C$ of $\mathcal F$ in terms of concepts involving $\Daleth$ and $h_{\mathcal C}$.
\end{enumerate}
In Section~\ref{S:div.and.tor-free} we relate topological and algebraic properties of the group $\Gimel$ to algebraic properties of the groups $\Cocc(G)$. In Section~\ref{Sub:free.min.ext} we show that a free object for the category of all \emph{minimal} group extensions of $\mathcal F$ does not exist. Finally, since the inclusion $F(\Daleth)\subseteq\Gimel$ is intimately related to the inclusion $\Cob\subseteq\Cocc$, we may use our results on the problem (1) and express the first inclusion in terms of familiar compact abelian groups. This is done in the following two situations.
\begin{itemize}
\item The group $\Gamma$ is a simply connected Lie group and the space $X$ is compact second countable with $\pi^1(X)\neq0$.
\item The group $\Gamma$ is a connected Lie group, the space $X$ is a compact manifold, the flow $\mathcal F$ is topologically free and possesses a free cycle.
\end{itemize}
These results are obtained in Section~\ref{S:circ.case}; see Theorems~\ref{T:F.in.Gim.CLAC} and~\ref{T:F.in.Gim.free.ccl}, respectively.
\medskip

\begin{enumerate}
\item[(5)] \emph{Find a dynamical interpretation of torsions in $\Coch$. Express the torsion subgroup $\tor(\Coch)$ of $\Coch$ in terms of familiar (or elementary) abelian groups.  Find sufficient conditions for $\tor(\Coch)$ to be a direct summand in $\Coch$ and determine its complementary summand $\Coch/\tor(\Coch)$.}
\end{enumerate}
In Corollary~\ref{C:MinH.rel.torH} in Section~\ref{S:funct.appr.sect} we show that the torsion subgroup $\tor(\Coch)$ of $\Coch$ coincides with the group of cohomology classes of non-minimal extensions from $\Cocc$; consequently, if $\tor(\Coch)$ is a direct summand in $\Coch$ then each of its complementary summands in $\Coch$ is a subgroup of the groupoid $\Cochm$ of the cohomology classes of extensions from $\Cocm$ ({see Subsection~\ref{Sub:gpd.min.ext}). Further, we express the group $\tor(\Coch)$ in terms of familiar abelian groups in the following situations.
\begin{itemize}
\item The group $\Gamma$ is simply connected and the space $X$ is compact; see Theorem~\ref{T:tor.struct} in Section~\ref{Sub:tor.elem}.
\item The group $\Gamma$ is a connected Lie group and the space $X$ is a compact manifold; see Theorem~\ref{T:tor.Lie.mnfld} in Section~\ref{Sub:tor.elem}.
\end{itemize}
Moreover, in the first of these two situations, the group $\tor(\Coch)$ is a direct summand in $\Coch$; see Theorem~\ref{T:min.exist.simp} in Section~\ref{Sub:lift.of.ext}. In the second situation this is also the case, provided the flow $\mathcal F$ is topologically free; see Theorem~\ref{T:div.gp.min.ext} in Section~\ref{S:tor.top.free.flw}. Finally, it follows from our results on the structure of the group $\Coch$ in Section~\ref{S:circ.case} that, under further mild assumptions on $\mathcal F$ or $X$, the complementary summand to $\tor(\Coch)$ in $\Coch$ is isomorphic to the additive group of real numbers $\mathbb R$ in both these situations. Indeed, in the first situation, this is true if the space $X$ is second countable with $\pi^1(X)\neq0$ (see Theorem~\ref{T:F.in.Gim.CLAC} in Section~\ref{S:circ.case}) and, in the second situation, this is true if the flow $\mathcal F$ possesses a free cycle (see Theorem~\ref{T:F.in.Gim.free.ccl} in Section~\ref{S:circ.case}).
\bigskip

Problems (1)--(5) formulated above are studied in Chapters~\ref{S:E.and.alg-top}, \ref{S:alg.top.asp} and~\ref{S:struct.res}. Chapter~\ref{S:def.of.E} is, in a sense, auxiliary - it is devoted to the study of tools and techniques, which are used to solve our five main problems. There are basically five main ideas that occur in Chapter~\ref{S:def.of.E}.
\begin{itemize}
\item[(i)] For every compact abelian group $G$, there is a (topological) isomorphism between the group of extensions $\Cocc(G)$ and the group of morphisms $\Hom(G^*,\Cocc)$, where $G^*$ stands for the Pontryagin dual of $G$; see Theorem~\ref{T:structure.coc} and its corollaries in Section~\ref{Sub:gp.ext}.
\item[(ii)] There is a functorial approach to sections of orbit closures $F(\mathcal C)$ of group extensions $\mathcal C$ of $\mathcal F$; see Theorem~\ref{T:functor.E.def} and its corollaries in Section~\ref{S:funct.appr.sect}. This approach, in which $F$ is viewed as a covariant functor from the category $\Coc$ of all group extensions of $\mathcal F$ to the category $\mathsf{CAGp}$ of all compact abelian groups, is very convenient and makes many proofs throughout the whole present work more transparent.
\item[(iii)] Given a compact abelian group $G$, the functor $F$ gives rise to a map $F\colon\Cocc(G)\to 2^G$ with values in the hypersemigroup $2^G$ over $G$ (see Subsection~\ref{Sub:hypsmgps} for the definition of $2^G$). In Section~\ref{S:sum.prop} it is shown that the map $F$ respects infinite sums in $\Cocc(G)$ and $2^G$ in a certain, rather feeble, way; see Theorem~\ref{T:first.ineq.yes}, Proposition~\ref{P:first.ineq.yes} and Theorem~\ref{T:second.ineq.not}. Although this property of $F$ from Theorem~\ref{T:first.ineq.yes} seems fairly weak at a first glance, we shall find it extremely useful on many occasions.
\item[(iv)] Since the group $\Cob(G)$ of coboundaries is rarely closed in the group $\Cocc(G)$ with respect to the topology of uniform convergence on compact sets, the quotient cohomology group $\Coch(G)$ does not carry a Hausdorff topology. In Section~\ref{S:ext-top.intr} we introduce on $\Coch(G)$ a new topology, which we call the \emph{ext-topology}. This is a Hausdorff group topology induced by a $2^G$-valued translation-invariant metric, which makes the map $F\colon\Coch(G)\to 2^G$ continuous; see Theorem~\ref{P:group.top.cont}. This topology turns out to be more than a mere abstract construction; as a matter of fact, in Chapter~\ref{S:struct.res} we will be able to express the group $\Coch(G)$ with the ext-topology in terms of familiar abelian topological groups for large classes of minimal flows $\mathcal F$.
\item[(v)] A special form of group extensions $\mathcal C\in\Cocc(G)$ as maps $\mathcal C\colon\Gamma\times X\to G$ allows us to lift them to group extensions across covering morphisms also in situations when the phase space $X$ of our flow $\mathcal F$ is not locally connected. More concretely, this is always possible (even for non-compact abelian groups $G$) if the acting group $\Gamma$ of $\mathcal F$ is simply connected; see Theorem~\ref{T:lift.simply.con} in Section~\ref{Sub:lift.of.ext}. An important corollary of this result is Theorem~\ref{T:min.exist.simp} from Section~\ref{Sub:lift.of.ext}, which binds the groups $\pi^1(X)$, $\Coch(\mathbb R)$ and $\Coch$ into a short exact sequence; this theorem will turn out to be an important tool in our study of group extensions of minimal flows with simply connected acting groups $\Gamma$.
\end{itemize}
\bigskip

Chapter~\ref{S:Preliminaries} of this work contains preliminary material from topology, theory of (compact) abelian groups and topological dynamics. The first four sections from this chapter are devoted mainly to recalling those concepts and results, which are used frequently in this work. However, in Section~\ref{S:top.dyn} we also introduce some new concepts and study their basic properties. It is of course not necessary to read the whole Chapter~\ref{S:Preliminaries} before turning to the core of this work beginning with Chapter~\ref{S:def.of.E}. In fact, the reader will find it sufficient to leaf through Subsections~\ref{Sub:Sets}, \ref{Sub:cat.gps}, \ref{Sub:bas.fcts.ab.gp}, \ref{Sub:some.conv.top} and \ref{Sub:Pontr.dual} for basic notation and conventions, and then continue with Subsections~\ref{Sub:min.flws}, \ref{Sub:group.ext}, \ref{Sub:cob.and.chmlg}, \ref{Sub:gpd.min.ext}, \ref{Sub:cat.of.ext.Coc} and~\ref{Sub:ind.morph.ext}. Afterwards, the reader may return to selected parts of Chapter~\ref{S:Preliminaries}, whenever he/she finds it useful to recall a relevant concept or result. We hope that this first chapter will make the reading of our work more comfortable.

\chapter{Preliminaries}\label{S:Preliminaries}

\section{Categories}

\subsection{Sets}\label{Sub:Sets}

We denote by $\mathsf{Sets}$\index[symbol]{$\mathsf{Sets}$} the category of sets. The cardinality of a set $A$ is denoted by $\card(A)$\index[symbol]{$\card(A)$} and a singleton set $\{x\}$ by a shorter symbol $x$. For the composition of maps $g$ and $f$ we write $gf$ instead of $g\circ f$. If $(X_j)_{j\in J}$ is a collection of sets and $i\in J$ then $\pr_i$\index[symbol]{$\pr_i$} stands for the projection $\prod_{j\in J}X_j\to X_i$. We denote the set of all integral, rational, real and complex numbers by $\mathbb Z$\index[symbol]{$\mathbb Z$}, $\mathbb Q$\index[symbol]{$\mathbb Q$}, $\mathbb R$\index[symbol]{$\mathbb R$} and $\mathbb C$\index[symbol]{$\mathbb C$}, respectively, and the symbol $\mathbb N$\index[symbol]{$\mathbb N$} is used to denote the set of all positive integers. Finally, we write $\mathfrak{c}\index[symbol]{$\mathfrak{c}$}=2^{\aleph_0}$ for the cardinality of the continuum.

\subsection{Categories of groups}\label{Sub:cat.gps}
If $\mathsf{Cat}$ is a category and $K$ is an object of $\mathsf{Cat}$ then we write $K\in\mathsf{Cat}$. We will be working with the following categories of groups:
\begin{itemize}
\item $\mathsf{AbGp}$\index[symbol]{$\mathsf{AbGp}$}--the category of abelian groups,
\item $\mathsf{AbTpGp}$\index[symbol]{$\mathsf{AbTpGp}$}--the category of abelian topological groups,
\item $\mathsf{DAGp}$\index[symbol]{$\mathsf{DAGp}$}--the category of discrete abelian groups,
\item $\mathsf{CAGp}$\index[symbol]{$\mathsf{CAGp}$}--the category of compact abelian groups,
\item $\mathsf{LCGp}$\index[symbol]{$\mathsf{LCGp}$}--the category of locally compact groups,
\item $\mathsf{LCAGp}$\index[symbol]{$\mathsf{LCAGp}$}--the category of locally compact abelian groups,
\item $\mathsf{LieGp}$[symbol]{$\mathsf{LieGp}$}--the cateogry of Lie groups.
\end{itemize}
The categories $\mathsf{AbGp}$ and $\mathsf{DAGp}$ are isomorphic in a natural way and we shall therefore identify them when convenient. If $G,H\in\mathsf{AbTpGp}$ then $\Hom(G,H)$\index[symbol]{$\Hom(G,H)$} denotes the set of all topological (that is, continuous) morphisms $G\to H$.

Every group $G\in\mathsf{AbTpGp}$ is a uniform space, the uniformity being derived in a usual way from the topological group structure on $G$. The group $G$ is called complete\index{complete!group} if every Cauchy net in $G$ is convergent in $G$. Recall that each $G\in\mathsf{LCAGp}$ is a complete group.

\subsection{Category $\mathsf{Hom}(K,\mathsf{Cat})$}

Given a category $\mathsf{Cat}$ and an object $K\in\mathsf{Cat}$, we denote by $\mathsf{Hom}(K,\mathsf{Cat})$\index[symbol]{$\mathsf{Hom}(K,\mathsf{Cat})$} the category defined as follows. The objects of $\mathsf{Hom}(K,\mathsf{Cat})$ are the morphisms $l\in\Hom(K,L)$ with $L\in\mathsf{Cat}$. Further, if $l_1\in\Hom(K,L_1)$ and $l_2\in\Hom(K,L_2)$ then $\Hom(l_1,l_2)$ consists of the morphisms $k\in\Hom(L_1,L_2)$ with $kl_1=l_2$. The composition of morphisms in $\mathsf{Hom}(K,\mathsf{Cat})$ is simply their composition in $\mathsf{Cat}$. Finally, if $L\in\mathsf{Cat}$ and $\id_L$ denotes the unit $\id_L\in\Hom(L,L)$ then $\id_L\in\Hom(l,l)$ for every $l\in\Hom(K,L)$.

\subsection{Category $\mathsf{Hom}(\mathsf{Cat},L)$}

Given a category $\mathsf{Cat}$ and an object $L\in\mathsf{Cat}$, we denote by $\mathsf{Hom}(\mathsf{Cat},L)$\index[symbol]{$\mathsf{Hom}(\mathsf{Cat},L)$} the category defined as follows. The objects of $\mathsf{Hom}(\mathsf{Cat},L)$ are the morphisms $l\in\Hom(K,L)$ with $K\in\mathsf{Cat}$. Further, if $l_1\in\Hom(K_1,L)$ and $l_2\in\Hom(K_2,L)$ then $\Hom(l_1,l_2)$ consists of the morphisms $k\in\Hom(K_1,K_2)$ with $l_2k=l_1$. The composition of morphisms in $\mathsf{Hom}(\mathsf{Cat},L)$ is their composition in $\mathsf{Cat}$. Finally, if $K\in\mathsf{Cat}$ and $\id_K$ denotes the unit $\id_K\in\Hom(K,K)$ then $\id_K\in\Hom(l,l)$ for every $l\in\Hom(K,L)$.

\subsection{Free objects and universal elements}

A \emph{free object}\index{free!object} for a category $\mathsf{Cat}$ is an object $K\in\mathsf{Cat}$ such that the set $\Hom(K,L)$ is a singleton for every object $L\in\mathsf{Cat}$. Any two free objects for a given category $\mathsf{Cat}$ are isomorphic within $\mathsf{Cat}$; we shall therefore often speak of \emph{the} free object $K$ for $\mathsf{Cat}$.

Let $\mathsf{Cat}$ be a category and $F\colon\mathsf{Cat}\to\mathsf{Sets}$ be a covariant (respectively, contravariant) functor. A \emph{universal element}\index{universal!element} for $F$ consists of an object $K\in\mathsf{Cat}$ and an element $k\in F(K)$ such that for every $L\in\mathsf{Cat}$ and every $l\in F(L)$ there is a unique $h\in\Hom(K,L)$ (respectively, $h\in\Hom(L,K)$) with $F(h)\colon k\mapsto l$. A universal element for a given functor $F$ is essentially unique: if $K'\in\mathsf{Cat}$ and $k'\in F(K')$ also constitute a universal element for $F$ then there is an isomorphism $h\colon K\to K'$ (respectively, $h\colon K'\to K$) with $F(h)\colon k\mapsto k'$. We shall therefore often speak of \emph{the} universal element for $F$.

\subsection{Direct systems and their limits}

A \emph{direct system}\index{direct!system} in a category $\mathsf{Cat}$ consists of a directed set $J$, a family of objects $K_j\in\mathsf{Cat}$ ($j\in J$) and a family of morphisms $q_{ik}\in\Hom(K_i,K_k)$ ($i\leq k$), such that $q_{jj}=\id_{K_j}$ for every $j\in J$ and $q_{kl}q_{ik}=q_{il}$ for all $i\leq k\leq l$. A \emph{direct limit}\index{direct!limit} of this system consists of an object $K\in\mathsf{Cat}$ and a family of morphisms $q_j\in\Hom(K_j,K)$ ($j\in J$) with $q_kq_{ik}=q_i$ for $i\leq k$, such that the following condition holds:
\begin{itemize}
\item given an object $L\in\mathsf{Cat}$ and a family of morphisms $s_j\in\Hom(K_j,L)$ ($j\in J$) with $s_kq_{ik}=s_i$ for $i\leq k$, there is a unique morphism $r\in\Hom(K,L)$ with $rq_j=s_j$ for every $j\in J$. 
\end{itemize}
Being a universal element for an appropriate covariant functor $\mathsf{Cat}\to\mathsf{Sets}$, a direct limit of a given direct system is unique up to isomorphism. We call a category $\mathsf{Cat}$ \emph{complete}\index{complete!category} if every direct system in $\mathsf{Cat}$ has a limit.

\subsection{Inverse systems and their limits}

An \emph{inverse system}\index{inverse!system} in a category $\mathsf{Cat}$ consists of a directed set $J$, a family of objects $K_j\in\mathsf{Cat}$ ($j\in J$) and a family of morphisms $p_{ik}\in\Hom(K_k,K_i)$ ($i\leq k$), such that $p_{jj}=\id_{K_j}$ for every $j\in J$ and $p_{ik}p_{kl}=p_{il}$ for all $i\leq k\leq l$. An \emph{inverse limit}\index{inverse!limit} of this system consists of an object $K\in\mathsf{Cat}$ and a family of morphisms $p_j\in\Hom(K,K_j)$ ($j\in J$) with $p_{ik}p_k=p_i$ for $i\leq k$, such that the following condition holds:
\begin{itemize}
\item given an object $L\in\mathsf{Cat}$ and a family of morphisms $s_j\in\Hom(L,K_j)$ ($j\in J$) with $p_{ik}s_k=s_i$ for $i\leq k$, there is a unique morphism $r\in\Hom(L,K)$ with $p_jr=s_j$ for every $j\in J$.
\end{itemize}
Being a universal element for an appropriate contravariant functor $\mathsf{Cat}\to\mathsf{Sets}$, an inverse limit of a given inverse system is unique up to isomorphism. We call a category $\mathsf{Cat}$ \emph{co-complete}\index{co-complete category} if every inverse system in $\mathsf{Cat}$ has a limit.

\section{Abelian groups}\label{S:abel.gps}

Our references in this section are monographs \cite{Fuchs1} and \cite{HofMor}.

\subsection{Basic facts}\label{Sub:bas.fcts.ab.gp}

Let $A$ be an abelian group. Depending on what is more convenient, we shall write the operations in $A$ either additively or multiplicatively; in the first case, $0$ stands for the neutral element of $A$, in the second case, the neutral element of $A$ is denoted by $e$.

Given $k\in\mathbb N$, the \emph{$k$-endomorphism of $A$}\index{$k$-endomorphism} is the morphism $\kappa_k\colon A\to A$\index[symbol]{$\kappa_k$}, acting by the rule $\kappa_k(a)=ka$ for every $a\in A$. The group $A$ is called
\begin{itemize}
\item $k$-divisible\index{$k$-divisible group}, if $\im(\kappa_k)=A$,
\item divisible\index{divisible!group}, if it is $k$-divisible for every $k\in\mathbb N$,
\item a torsion group\index{torsion!group}, if $A=\bigcup_{k\in\mathbb N}\ker(\kappa_k)$,
\item torsion-free\index{torsion-free group}, if $\ker(\kappa_k)=0$ for every $k\in\mathbb N$.
\end{itemize}

The subgroup of $A$ generated by a set $S\subseteq A$ is denoted by $\langle S\rangle$\index[symbol]{$\langle S\rangle$}; if $S=\{a\}$ is a singleton then we write $\langle a\rangle$\index[symbol]{$\langle a\rangle$} instead of $\langle\{a\}\rangle$. The notation $B\sbgp A$ is used to express that $B$ is a subgroup of $A$; if the group $A$ is topological then the symbol $B\sbgp A$\index[symbol]{$B\sbgp A$} means that $B$ is a closed subgroup of $A$.

Let $A$ be an abelian group, $B$ be a subgroup of $A$ and $d$ be a positive integer. We say that $B$ is a \emph{$d$-pure subgroup}\index{$d$-pure subgroup} of $A$ (or that $B$ is $d$-pure in $A$) if $B\cap\kappa_d(A)=\kappa_d(B)$; if $A$ is torsion-free then this is equivalent to $\kappa_d^{-1}(B)\subseteq B$. The group $B$ is \emph{pure}\index{pure subgroup} in $A$ if it is $d$-pure in $A$ for every $d\in\mathbb N$.

The direct sum and the direct product of a family of abelian groups $A_i$ ($i\in I$) is denoted by $\bigoplus_{i\in I}A_i$\index[symbol]{$\bigoplus_{i\in I}A_i$} and $\prod_{i\in I}A_i$\index[symbol]{$\prod_{i\in I}A_i$}, respectively. Finally, given a cardinal $\mathfrak{k}$, the direct sum and the direct product of $\mathfrak{k}$ copies of a group $A\in\mathsf{AbGp}$ is denoted by $A^{(\mathfrak{k})}$\index[symbol]{$A^{(\mathfrak{k})}$} and $A^{\mathfrak{k}}$\index[symbol]{$A^{\mathfrak{k}}$}, respectively.

\subsection{Groups of morphisms}

Let $A,B\in\mathsf{AbGp}$. We denote by $\Hom(A,B)$ the set of all morphisms $A\to B$. With the point-wise defined operations, $\Hom(A,B)$ is an abelian group. Further, if $C\sbgp A$ and $D\sbgp B$, we denote by $\Hom(A,C;B,D)$\index[symbol]{$\Hom(A,C;B,D)$} the set of all morphisms $h\colon A\to B$ with $h(C)\subseteq D$. Clearly, $\Hom(A,C;B,D)$ is a subgroup of $\Hom(A,B)$. Finally, given abelian groups $A$, $B$, $A_i$, $B_i$ ($i\in I$), there are isomorphisms
\begin{equation*}
\Hom\left(\bigoplus_{i\in I}A_i,B\right)\ni h\mapsto(h\sigma_i)_{i\in I}\in\prod_{i\in I}\Hom(A_i,B),
\end{equation*}
where $\sigma_j$ is the inclusion morphism $A_j\to\bigoplus_{i\in I}A_i$ for every $j\in I$, and
\begin{equation*}
\Hom\left(A,\prod_{i\in I}B_i\right)\ni h\mapsto(\text{pr}_ih)_{i\in I}\in\prod_{i\in I}\Hom(A,B_i),
\end{equation*}
where $\pr_j$ is the projection morphism $\prod_{i\in I}B_i\to B_j$ for every $j\in I$.

\subsection{Groupoids of monomorphisms}

Given $A,B\in\mathsf{AbGp}$, we use the symbol $\Mon(A,B)$\index[symbol]{$\Mon(A,B)$} to denote the set of all monomorphisms $A\to B$ together with the trivial morphism $0$. Algebraically, $\Mon(A,B)$ is a groupoid; that is, the product in $\Hom(A,B)$ restricts to a partial operation on $\Mon(A,B)$, which is associative and commutative, and for which $0$ is a neutral element. Clearly, every element $h\in\Mon(A,B)$ has an inverse $-h\in\Mon(A,B)$. We shall refer to $\Mon(A,B)$ as the \emph{groupoid of monomorphisms}\index{groupoid!of monomorphisms} $A\to B$.

\subsection{Elementary divisors}

Let $F$ be a free abelian group of a finite rank $r$ and let $E$ be a subgroup of $F$. Then there exist a basis $f_1,\dots,f_r$ for $F$, an integer $0\leq n\leq r$ and positive integers $d_1,\dots,d_n$, such that $d_1f_1,\dots,d_nf_n$ is a basis for $E$ and $d_i$ divides $d_{i+1}$ for $i=1,\dots,n-1$. Such numbers $d_1,\dots,d_n$ are unique and are called the \emph{elementary divisors of $E$ in $F$}\index{elementary divisors}.

\subsection{Torsions and divisibility}

Let $A\in\mathsf{AbGp}$. The union of all divisible subgroups of $A$ is the largest divisible subgroup of $A$, we denote it by $\Div(A)$\index[symbol]{$\Div(A)$}. If $A$ is either torsion-free or admits a topology of a compact group then $\Div(A)=\bigcap_{k\in\mathbb N}\kappa_k(A)$. The group $\Div(A)$ is a direct summand in $A$; each of its complementary summands in $A$ is isomorphic to $A/\Div(A)$ and is a reduced group\index{reduced group} (that is, it contains no non-trivial divisible subgroups).

Given a positive integer $k$, we denote by $\tor_k(A)$\index[symbol]{$\tor_k(A)$} the group of the $k$-torsions of $A$, that is, the group of all $a\in A$ with $ka=\kappa_k(a)=0$. The union $\tor(A)=\bigcup_{k\in\mathbb N}\tor_k(A)$\index[symbol]{$\tor(A)$} is called the torsion subgroup of $A$. If $a\in\tor(A)$ then the order\index{order of a group element} (or the period)\index{period of a group element} of $a$ in $A$, denoted by $\ord(a)$\index[symbol]{$\ord(a)$}, is the smallest of the positive integers $k$ with $ka=0$. For every $k\in\mathbb N$, we denote by $\per_k(A)$\index[symbol]{$\per_k(A)$} the set of all elements $a\in A$ with order $\ord(a)=k$. If $a\in A\setminus\tor(A)$ then we say that $a$ is of an infinite order in $A$ and write $\ord(a)=\infty$.

We shall view the group $\mathbb T^1$ as a subset of $\mathbb C$ with the complex multiplication. Then for every $k\in\mathbb N$, $\tor_k(\mathbb T^1)$ consists of the $k^{\text{th}}$ roots of unity. Since there is an obvious isomorphism between $\tor_k(\mathbb T^1)$ and $\mathbb Z_k$, we shall often identify these two groups.

Given co-prime integers $q,r\geq2$, there is an isomorphism
\begin{equation*}
\text{tor}_{qr}(A)=\text{tor}_q(A)\oplus\text{tor}_r(A)\cong\text{tor}_q(A)\times\text{tor}_r(A),
\end{equation*}
under which the set $\per_{qr}(A)\subseteq\tor_{qr}(A)$ corresponds to $\per_q(A)\times\per_r(A)$. Consequently, if $k=p_1^{m_1}\dots p_n^{m_n}$ is a prime decomposition of an integer $k\geq2$ then there is an isomorphism
\begin{equation*}
\text{tor}_k(A)=\text{tor}_{p_1^{m_1}}(A)\oplus\dots\oplus\text{tor}_{p_n^{m_n}}(A)
\cong\text{tor}_{p_1^{m_1}}(A)\times\dots\times\text{tor}_{p_n^{m_n}}(A),
\end{equation*}
under which the set $\per_k(A)\subseteq\tor_k(A)$ corresponds to $\per_{p_1^{m_1}}(A)\times\dots\times\per_{p_n^{m_n}}(A)$.

If the group $A$ is divisible then so is its torsion subgroup $\tor(A)$; in such a situation, $\tor(A)$ is a direct summand in $A$, each of its complementary summands in $A$ is isomorphic to $A/\tor(A)$ and is a torsion-free group. Finally, for every $A\in\mathsf{AbGp}$,
\begin{equation*}
\tor(\Div(A))=\Div(\tor(A))=\tor(A)\cap\Div(A).
\end{equation*}

\subsection{Betti number}

Let $A$ be a finitely generated abelian group. Then the torsion subgroup $\tor(A)$ of $A$ is a direct summand in $A$ and its complementary summand $F$ in $A$ is a free abelian group of a finite rank. We denote the rank of $F$ by $\beta(A)$\index[symbol]{$\beta(A)$} and call it the Betti number\index{Betti number of a group} of $A$.

\subsection{2-cocycles and 2-coboundaries}\label{Sub:2-coc.2-cob}

Let $A,C\in\mathsf{AbGp}$. A \emph{(symmetric) $2$-cocycle\index{$2$-cocycle} over $C$ with values in $A$} is a map $\varphi\colon C\times C\to A$, satisfying the following identities:
\begin{enumerate}
\item[($C_1$)] $\varphi(c_1,c_2)+\varphi(c_1+c_2,c_3)=\varphi(c_1,c_2+c_3)+\varphi(c_2,c_3)$,
\item[($C_2$)] $\varphi(0,c)=\varphi(c,0)=0$,
\item[($C_3$)] $\varphi(c_1,c_2)=\varphi(c_2,c_1)$.
\end{enumerate}
Such maps $\varphi$ form a group $Z^2(C,A)$\index[symbol]{$Z^2(C,A)$} with the point-wise operations. Every map $f\colon C\to A$ with $f(0)=0$ gives rise to a $2$-cocycle $\varphi\in Z^2(C,A)$, given by
\begin{equation*}
\varphi(c_1,c_2)=f(c_1)+f(c_2)-f(c_1+c_2);
\end{equation*}
in this situation we say that $\varphi$ is a \emph{$2$-coboundary}\index{$2$-coboundary} with a \emph{transfer function}\index{transfer function!for a $2$-coboundary} $f$. The $2$-coboundaries form a subgroup $B^2(C,A)$\index[symbol]{$B^2(C,A)$} of $Z^2(C,A)$; we shall denote the corres\-pon\-ding quotient cohomology group $Z^2(C,A)/B^2(C,A)$ by $\Ext(C,A)$\index[symbol]{$\Ext(C,A)$}.

Given a map $\varphi\colon C\times C\to A$, we denote by $C\times_{\varphi}A$\index[symbol]{$C\times_{\varphi}A$} the product $C\times A$ equipped with the binary operation
\begin{equation*}
(c_1,a_1)+(c_2,a_2)=(c_1+c_2,a_1+a_2+\varphi(c_1,c_2)).
\end{equation*}
In order that $C\times_{\varphi}A$ be a group, it is necessary and sufficient that conditions ($C_1$) and ($C_2$) hold true. Moreover, the group $C\times_{\varphi}A$ is then commutative if and only if $\varphi$ satisfies condition ($C_3$). Consequently, $C\times_{\varphi}A\in\mathsf{AbGp}$ if and only if $\varphi\in Z^2(C,A)$.

\subsection{Extensions of abelian groups}\label{Sub:ext.ab.gps}

A short exact sequence of abelian groups
\begin{equation}\label{Eq:ext.abel.gps}
E\colon0\longrightarrow A\stackrel{\mu}{\longrightarrow}B\stackrel{\nu}{\longrightarrow}C\longrightarrow0
\end{equation}
is called an \emph{extension of $C$ by $A$}\index{extension!of a group}. Two such extensions
\begin{equation*}
E_1\colon 0\longrightarrow A\stackrel{\mu_1}{\longrightarrow}B_1\stackrel{\nu_1}{\longrightarrow}C\longrightarrow0\hspace{4mm}\text{and}\hspace{4mm}
E_2\colon 0\longrightarrow A\stackrel{\mu_2}{\longrightarrow}B_2\stackrel{\nu_2}{\longrightarrow}C\longrightarrow0
\end{equation*}
are called \emph{equivalent}\index{equivalent extensions!of a group} if there is a morphism of groups $h\colon B_1\to B_2$ with $h\mu_1=\mu_2$ and $\nu_2h=\nu_1$; by the five-lemma, such a morphism $h$ is automatically an isomorphism. We write $E_1\cong E_2$ to express that $E_1$ and $E_2$ are equivalent extensions. The extension $E$ in (\ref{Eq:ext.abel.gps}) is said to \emph{split}\index{split extension of a group} if it is equivalent to the extension
\begin{equation*}
E_0(C,A)\index[symbol]{$E_0(C,A)$}=E_0\index[symbol]{$E_0$}\colon 0\longrightarrow A\stackrel{\mu_0}{\longrightarrow}C\oplus A\stackrel{\nu_0}{\longrightarrow}C\longrightarrow0,
\end{equation*}
where $\mu_0(a)=(0,a)$\index[symbol]{$\mu_0$} and $\nu_0(c,a)=c$\index[symbol]{$\nu_0$} for all $a\in A$ and $c\in C$; this occurs if and only if one of the following conditions holds:
\begin{itemize}
\item there is $k\in\Hom(C,B)$ with $\nu k=\Id_C$,
\item there is $l\in\Hom(B,A)$ with $l\mu=\Id_A$.
\end{itemize}

Every $2$-cocycle $\varphi\in Z^2(C,A)$ gives rise to an extension
\begin{equation}\label{Eq:ext.abel.gps.coc}
E(\varphi)\index[symbol]{$E(\varphi)$}=E_{\varphi}\index[symbol]{$E_{\varphi}$}\colon 0\longrightarrow A\stackrel{\mu_{\varphi}}{\longrightarrow}C\times_{\varphi}A\stackrel{\nu_{\varphi}}{\longrightarrow}C\longrightarrow0,
\end{equation}
where $\mu_{\varphi}(a)=(0,a)$\index[symbol]{$\mu_{\varphi}$} and $\nu_{\varphi}(c,a)=c$\index[symbol]{$\nu_{\varphi}$} for all $a\in A$ and $c\in C$. Recall the following facts:
\begin{itemize}
\item every extension $E$ of the form (\ref{Eq:ext.abel.gps}) is equivalent to an extension $E_{\varphi}$ from (\ref{Eq:ext.abel.gps.coc}) for some $\varphi\in Z^2(C,A)$,
\item given $\varphi_1,\varphi_2\in Z^2(C,A)$, the extensions $E_{\varphi_1}$, $E_{\varphi_2}$ are equivalent if and only if the $2$-cocycles $\varphi_1,\varphi_2$ are cohomologous (that is, if and only if $\varphi_1-\varphi_2\in B^2(C,A)$).
\end{itemize}
It follows from these two statements that the set of equivalence classes of the extensions of $C$ by $A$ can be identified with the cohomology group $\Ext(C,A)$.

Recall the following properties of $\Ext$:
\begin{itemize}
\item given abelian groups $A,C,A_i,C_i$ ($i\in I$), there are isomorphisms
\begin{equation*}
\Ext\left(\bigoplus_{i\in I}C_i,A\right)\cong\prod_{i\in I}\Ext(C_i,A)
\end{equation*}
and
\begin{equation*}
\Ext\left(C,\prod_{i\in I}A_i\right)\cong\prod_{i\in I}\Ext(C,A_i),
\end{equation*}
\item if $C\in\mathsf{AbGp}$ then the equality $\Ext(C,A)=0$ holds for every $A\in\mathsf{AbGp}$ if and only if the group $C$ is free abelian (in particular, if $B\in\mathsf{AbGp}$, $A\sbgp B$ and the group $B/A$ is free abelian, then $A$ is a direct summand in $B$),
\item if $A\in\mathsf{AbGp}$ then the equality $\Ext(C,A)=0$ holds for every $C\in\mathsf{AbGp}$ if and only if the group $A$ is divisible (in particular, if $A$ is a divisible subgroup of an abelian group $B$ then $A$ is a direct summand in $B$),
\item $\Ext(\mathbb Z_k,A)\cong A/kA$ for all $A\in\mathsf{AbGp}$ and $k\in\mathbb N$.
\end{itemize}

\subsection{Some concrete extensions}

Let $k\in\mathbb N$. Consider the short exact sequence in $\mathsf{AbGp}$
\begin{equation}\label{Eq:ext.for.Zk}
E(\mathbb Z_k)\index[symbol]{$E(\mathbb Z_k)$}\colon 0\longrightarrow\mathbb Z\stackrel{\kappa_k}{\longrightarrow}\mathbb Z\stackrel{p_k}{\longrightarrow}\mathbb Z_k\longrightarrow0,
\end{equation}
where $\kappa_k(m)=km$ and $p_k(m)=m$~(mod $k$) for every $m\in\mathbb Z$. This sequence is equivalent to the sequence $E_{\varphi_k}$, where $\varphi_k\in Z^2(\mathbb Z_k,\mathbb Z)$ acts by the rule
\begin{equation}\label{Eq:2-coc.for.Zk}
\varphi_k(r,s)=
\begin{cases}
0; & \text{if $r+s<k$},\\
1; & \text{if $r+s\geq k$},
\end{cases}
\end{equation}
for all $0\leq r,s\leq k-1$. Indeed, the map
\begin{equation*}
h\colon\mathbb Z\ni m\mapsto\left(p_k(m),\left[\frac{m}{k}\right]\right)\in\mathbb Z_k\times_{\varphi_k}\mathbb Z
\end{equation*}
is an isomorphism of groups with $h\kappa_k=\mu_{\varphi_k}$ and $\nu_{\varphi_k}h=p_k$, where $[m/k]$ stands for the integer part of $m/k$.

Now consider the short exact sequence in $\mathsf{AbGp}$
\begin{equation}\label{Eq:ext.for.Q/Z}
E(\mathbb Q/\mathbb Z)\index[symbol]{$E(\mathbb Q/\mathbb Z)$}\colon0\longrightarrow\mathbb Z\stackrel{j}{\longrightarrow}\mathbb Q\stackrel{p}{\longrightarrow}\mathbb Q/\mathbb Z\longrightarrow0,
\end{equation}
where $j$ is the inclusion morphism and $p$ is the quotient morphism. This sequence is equivalent to the sequence $E_{\psi}$, where $\psi\in Z^2(\mathbb Q/\mathbb Z,\mathbb Z)$ acts by the rule
\begin{equation}\label{Eq:2-coc.for.Q/Z}
\psi(r,s)=
\begin{cases}
0; & \text{if $r+s<1$},\\
1; & \text{if $r+s\geq 1$},
\end{cases}
\end{equation}
for all rationals $0\leq r,s<1$. Indeed, the map
\begin{equation*}
h\colon\mathbb Q\ni q\mapsto(p(q),[q])\in(\mathbb Q/\mathbb Z)\times_{\psi}\mathbb Z
\end{equation*}
is an isomorphism of groups with $hj=\mu_{\psi}$ and $\nu_{\psi}h=p$, where $[q]$ stands for the integer part of $q$.

\subsection{Direct sum of sequences and $2$-cocycles}\label{Sub:dir.sum.seq.coc}

Given short exact sequences of abelian groups
\begin{equation*}
E_i\colon0\longrightarrow A_i\stackrel{\mu_i}{\longrightarrow}B_i\stackrel{\nu_i}{\longrightarrow}C_i\longrightarrow0,
\end{equation*}
$i=1,\dots,n$, there is an exact sequence of abelian groups
\begin{equation*}
E\colon0\longrightarrow\bigoplus_{i=1}^nA_i\stackrel{\mu}{\longrightarrow}\bigoplus_{i=1}^nB_i\stackrel{\nu}{\longrightarrow}\bigoplus_{i=1}^nC_i\longrightarrow0,
\end{equation*}
where $\mu=\mu_1\times\dots\times\mu_n$ and $\nu=\nu_1\times\dots\times\nu_n$. We call $E$ the direct sum\index{direct sum!of sequences} of the sequences $E_i$ and write $E=\bigoplus_{i=1}^nE_i$.

Write $A=\bigoplus_{i=1}^nA_i$ and $C=\bigoplus_{i=1}^nC_i$. For $i=1,\dots,n$, let $\pi_i$ and $\vartheta_i$ denote the projection morphisms $C\to C_i$ and $A\to A_i$, respectively, and fix $\varphi_i\in Z^2(C_i,A_i)$ with $E_i\cong E_{\varphi_i}$. Define $\varphi\in Z^2(C,A)$ by the rule $\vartheta_i\varphi=\varphi_i(\pi_i\times\pi_i)$ for $i=1,\dots,n$. We call $\varphi$ the direct sum\index{direct sum!of $2$-cocycles} of the $2$-cocycles $\varphi_i$ and write $\varphi=\bigoplus_{i=1}^n\varphi_i$. There are equivalences
\begin{equation*}
\bigoplus_{i=1}^nE_i\cong\bigoplus_{i=1}^nE(\varphi_i)\cong E\left(\bigoplus_{i=1}^n\varphi_i\right).
\end{equation*}

\subsection{Exact Hom-Ext sequences}

Let $A,B,C\in\mathsf{AbGp}$ and $h\in\Hom(A,B)$. The morphism $Z^2(C,A)\ni\varphi\mapsto h\varphi\in Z^2(C,B)$ maps $2$-coboundaries to $2$-co\-boun\-daries and so it gives rise to a morphism $\Ext(C,h)\colon\Ext(C,A)\to\Ext(C,B)$\index[symbol]{$\Ext(C,h)$}. Clearly, the assignment $h\mapsto\Ext(C,h)$ is covariant functorial. Similarly, if $A,C,D\in\mathsf{AbGp}$ and $h\in\Hom(D,C)$ then the morphism $Z^2(C,A)\ni\varphi\mapsto\varphi(h\times h)\in Z^2(D,A)$ induces a morphism $\Ext(h,A)\colon\Ext(C,A)\to\Ext(D,A)$\index[symbol]{$\Ext(h,A)$}. This time, the assignment $h\mapsto\Ext(h,A)$ is contravariant functorial.

Now let $E$ be a short exact sequence of abelian groups as in (\ref{Eq:ext.abel.gps}) and $D\in\mathsf{AbGp}$. By a theorem due to Cartan and Eilenberg, there is an exact sequence
\begin{equation*}
\begin{split}
0&\longrightarrow\Hom(D,A)\xrightarrow{\Hom(D,\mu)}\Hom(D,B)\xrightarrow{\Hom(D,\nu)}\Hom(D,C)\\
&\stackrel{\epsilon}{\longrightarrow}\Ext(D,A)\xrightarrow{\Ext(D,\mu)}\Ext(D,B)\xrightarrow{\Ext(D,\nu)}\Ext(D,C)\longrightarrow0,
\end{split}
\end{equation*}
which we shall refer to as the \emph{covariant Hom-Ext sequence}\index{covariant Hom-Ext sequence} derived from $E$ and associated to $D$. Recall that the \emph{connecting morphism}\index{connecting morphism} $\epsilon$ assigns to every morphism $h\in\Hom(D,C)$ the cohomology class of $\varphi(h\times h)\in Z^2(D,A)$, where $\varphi\in Z^2(C,A)$ is chosen in such a way that $E\cong E_{\varphi}$.

Further, by a theorem due to Cartan and Eilenberg, $E$ and $D$ give rise to an exact sequence
\begin{equation*}
\begin{split}
0&\longrightarrow\Hom(C,D)\xrightarrow{\Hom(\nu,D)}\Hom(B,D)\xrightarrow{\Hom(\mu,D)}\Hom(A,D)\\
&\stackrel{\epsilon}{\longrightarrow}\Ext(C,D)\xrightarrow{\Ext(\nu,D)}\Ext(B,D)\xrightarrow{\Ext(\mu,D)}\Ext(A,D)\longrightarrow0,
\end{split}
\end{equation*}
which we shall refer to as the \emph{contravariant Hom-Ext sequence}\index{contravariant Hom-Ext sequence} derived from $E$ and associated to $D$. This time, the \emph{connecting morphism}\index{connecting morphism} $\epsilon$  assigns to every morphism $k\in\Hom(A,D)$ the cohomology class of $k\varphi\in Z^2(C,D)$, where, as above, $\varphi\in Z^2(C,A)$ is chosen in such a way that $E\cong E_{\varphi}$.

\subsection{Whitehead groups}

An abelian group $A$ is called a \emph{Whitehead group}\index{Whitehead group} if it satisfies the equality $\Ext(A,\mathbb Z)=0$ or, equivalently, if every morphism $A\to\mathbb T^1$ lifts to a morphism $A\to\mathbb R$ across the usual quotient morphism $\mathbb R\to\mathbb T^1$. Clearly, every free abelian group is a Whitehead group. Though the converse implication is undecidable in ZFC by the Shelah's independence theorem, every \emph{countable} Whitehead group is free abelian by the Pontryagin theorem.

\subsection{Tensor and torsion products}

Given groups $A,C\in\mathsf{AbGp}$, we denote their tensor product\index{tensor product} by $A\otimes C$\index[symbol]{$A\otimes C$}. Recall that $A\otimes C\in\mathsf{AbGp}$ by definition and that the monomial tensors $a\otimes c$ ($a\in A$, $c\in C$) form a set of generators of $A\otimes C$. Further, given groups $A,B,C,D\in\mathsf{AbGp}$ and morphisms $h\in\Hom(A,B)$, $k\in\Hom(C,D)$, we denote by $h\otimes k$\index[symbol]{$h\otimes k$} the tensor product of $h$ and $k$. By definition, $h\otimes k\in\Hom(A\otimes C,B\otimes D)$. Recall that on monomial tensors from $A\otimes C$, the morphism $h\otimes k$ acts by the rule $(h\otimes k)(a\otimes c)=h(a)\otimes k(c)$ ($a\in A$, $c\in C$).

Recall the following properties of tensor products:
\begin{itemize}
\item there is an isomorphism $A\otimes C\cong C\otimes A$ for all $A,C\in\mathsf{AbGp}$,
\item $A\otimes C=0$ whenever $A\in\mathsf{AbGp}$ is a torsion group and $C\in\mathsf{AbGp}$ is divisible,
\item there are isomorphisms $\mathbb Z\otimes C\cong C$ and $\mathbb Z_k\otimes C\cong C/kC$ for all $C\in\mathsf{AbGp}$ and $k\in\mathbb N$,
\item given abelian groups $A_i$ ($i\in I$), $C_j$ ($j\in J$), there is an isomorphism
\begin{equation*}
\Bigg(\bigoplus_{i\in I}A_i\Bigg)\otimes\Bigg(\bigoplus_{j\in J}C_j\Bigg)\cong\bigoplus_{\substack{i\in I \\ j\in J}}A_i\otimes C_j,
\end{equation*}
\item given abelian groups $A,C,D$, there is an isomorphism of groups
\begin{equation}\label{Eq:tens.prod.hom.top}
\varphi\colon\Hom(A\otimes C,D)\ni h\mapsto \varphi(h)=\widetilde{h}\in\Hom(A,\Hom(C,D)),
\end{equation}
where $\widetilde{h}(a)\colon C\ni c\mapsto h(a\otimes c)\in D$ for every $a\in A$; if the group $D$ is topological and all three $\Hom$ groups above are equipped with the topology of point-wise convergence, then $\varphi$ is in fact a topological isomorphism.
\end{itemize}

Given groups $A,C\in\mathsf{AbGp}$, we denote their torsion product\index{torsion!product} by $\Tor(A,C)$\index[symbol]{$\Tor(A,C)$}. Recall that $\Tor(A,C)=\Tor(\tor(A),\tor(C))\in\mathsf{AbGp}$ is a torsion group by definition. Further, given groups $A,B,C,D\in\mathsf{AbGp}$ and morphisms $h\in\Hom(A,B)$, $k\in\Hom(C,D)$, we denote by $\Tor(h,k)$\index[symbol]{$\Tor(h,k)$} the torsion product of $h$ and $k$. By definition, $\Tor(h,k)$ is an element of $\Hom(\Tor(A,C),\Tor(B,D))$.

Recall the following properties of $\Tor$:
\begin{itemize}
\item there is an isomorphism $\Tor(A,C)\cong\Tor(C,A)$ for all $A,C\in\mathsf{AbGp}$,
\item there are isomorphisms $\Tor(\mathbb Z_k,C)\cong\tor_k(C)$ and $\Tor(\mathbb Q/\mathbb Z,C)\cong\tor(C)$ for all $C\in\mathsf{AbGp}$ and $k\in\mathbb N$,
\item given abelian groups $A_i$ ($i\in I$), $C_j$ ($j\in J$), there is an isomorphism
\begin{equation*}
\Tor\Bigg(\bigoplus_{i\in I}A_i,\bigoplus_{j\in J}C_j\Bigg)\cong\bigoplus_{\substack{i\in I\\ j\in J}}\Tor(A_i,C_j).
\end{equation*}
\end{itemize}

Now let $E$ be a short exact sequence of abelian groups as in (\ref{Eq:ext.abel.gps}) and $D\in\mathsf{AbGp}$. By a theorem due to Cartan and Eilenberg, there is an exact sequence
\begin{equation*}\label{Eq:ten.tor.seq}
\begin{split}
0&\longrightarrow\Tor(A,D)\xrightarrow{\Tor(\mu,\id)}\Tor(B,D)\xrightarrow{\Tor(\nu,\id)}\Tor(C,D)\\
&\stackrel{\epsilon}{\longrightarrow}A\otimes D\xrightarrow{\mu\otimes\id}B\otimes D\xrightarrow{\nu\otimes\id}C\otimes D\longrightarrow0,
\end{split}
\end{equation*}
which we shall refer to as the \emph{exact sequence for tensor and torsion products}\index{exact sequence for tensor and torsion products} derived from $E$ and associated to $D$. The morphism $\epsilon$ is called the \emph{connecting morphism}\index{connecting morphism}.

\subsection{Divisible hull and injective resolution}

Let $A$ be a torsion-free abelian group. The tensor product $\mathbb Q\otimes A$ is a torsion-free abelian group and each of its elements is expressible in the form of a monomial tensor. The mapping $j\colon A\to\mathbb Q\otimes A$, given by $j(a)=1\otimes a$, is a monomorphism of groups with the image $\im(j)=\mathbb Z\otimes A$. In fact, $\mathbb Q\otimes A$ is a divisible hull\index{divisible!hull} of $A$ with respect to the inclusion morphism $j$.

Under the identification $A\cong\im(j)=\mathbb Z\otimes A$, there is an isomorphism
\begin{equation*}
(\mathbb Q\otimes A)/A\cong(\mathbb Q/\mathbb Z)\otimes A.
\end{equation*}
More precisely, if $q$ stands for the quotient morphism $\mathbb Q\to\mathbb Q/\mathbb Z$ then
\begin{equation*}
p=q\otimes\text{id}_A\colon\mathbb Q\otimes A\to(\mathbb Q/\mathbb Z)\otimes A
\end{equation*}
is a quotient morphism with kernel $\ker(p)=\mathbb Z\otimes A$. We will refer to the corresponding short exact sequence
\begin{equation}\label{Eq:inj.res}
0\longrightarrow A\stackrel{j}{\longrightarrow}\mathbb Q\otimes A\stackrel{p}{\longrightarrow}(\mathbb Q/\mathbb Z)\otimes A\longrightarrow0
\end{equation}
as the \emph{injective resolution}\index{injective resolution} of $A$.

Recall that a subset $S$ of $A$ is called \emph{independent}\index{independent set in a group} in $A$ if the subgroup $\langle S\rangle$ of $A$ generated by $S$ is free abelian with the basis $S$. The set $S$ is independent in $A$ if and only if its image $j(S)$ is a set of rationally independent elements of $\mathbb Q\otimes A$. Every independent set $S\subseteq A$ is contained in a maximal independent set. A set $S\subseteq A$ is maximal independent if and only if its image $j(S)$ is a rational linear basis for $\mathbb Q\otimes A$. Consequently, any two maximal independent sets have the same cardinality $\rank(A)$\index[symbol]{$\rank(A)$} called the \emph{rank}\index{rank of a group} of $A$, which coincides with the rational linear dimension $\QLSdim(\mathbb Q\otimes A)$ of $\mathbb Q\otimes A$. Finally, an independent set $S\subseteq A$ is maximal independent if and only if the quotient group $A/\langle S\rangle$ is a torsion group.

\section{Topology}\label{S:topology}

Our references in this section are monographs \cite{Hat} and \cite[Chapter~8]{HofMor}.

\subsection{Some conventions}\label{Sub:some.conv.top}

By a (topological) space we mean a Hausdorff topological space. A space $X$ is often, explicitly or implicitly, assumed pointed, its base point then being denoted by $z$. The base point of a topological group is always its identity. If $X$ and $Y$ are spaces with base points $z$ and $v$, respectively, then the product space $X\times Y$ is equipped with the base point $(z,v)$. We shall denote the weight of $X$ by $\w(X)$\index[symbol]{$\w(X)$} and the density of $X$ by $d(X)$\index[symbol]{$d(X)$}. Finally, the symbol $\mathsf{CLAC}$\index[symbol]{$\mathsf{CLAC}$} stands for the category of the connected and locally arc-wise connected topological spaces; clearly, each space $X\in\mathsf{CLAC}$ is arc-wise connected.

\subsection{First weak homology group}\label{Sub:frst.wk.hm.gp}

Let $X$ be a pointed space with the base point $z$. We denote the fundamental group of $X$ by $\pi_1(X,z)$\index[symbol]{$\pi_1(X,z)$} or, briefly, by $\pi_1(X)$\index[symbol]{$\pi_1(X)$}. To simplify notation, we shall sometimes identify a loop in $X$ based at $z$ with its path-homotopy class in $\pi_1(X)$. The product of two elements $f,g\in\pi_1(X)$ in $\pi_1(X)$ is denoted by $f*g$\index[symbol]{$f*g$}; the same symbol $f*g$ is used to denote the product of two paths $f,g$ in $X$ with $f(1)=g(0)$. (If $X$ is a topological group with the identity $z$ then the group $\pi_1(X)$ is abelian, the product in $\pi_1(X)$ coincides with the point-wise product in $X$ and the inverse of $f\in\pi_1(X)$ in $\pi_1(X)$ is its point-wise inverse in $X$.) Further, given $z'\in X$ and a path $f$ in $X$ from $z$ to $z'$, we denote by $\widehat{f}$\index[symbol]{$\widehat{f}$} the isomorphism $\pi_1(X,z)\ni g\mapsto\overline{f}*g*f\in\pi_1(X,z')$, where $\overline{f}$\index[symbol]{$\overline{f}$} stands for the path reverse to $f$.

The first homology group\index{first homology group} of $X$ is denoted by $H_1(X)$\index[symbol]{$H_1(X)$} and is interpreted exclusively as the abelianisation of $\pi_1(X)$. The first weak homology group\index{first weak homology group} of $X$ is denoted by $H_1^w(X)$\index[symbol]{$H_1^w(X)$}. Recall that $H_1^w(X)$ is defined as the largest torsion-free quotient group of $H_1(X)$, that is, $H_1^w(X)=H_1(X)/\tor(H_1(X))$; we shall denote by $p_X$\index[symbol]{$p_X$} the quotient morphism $\pi_1(X)\to H_1^w(X)$. Should it be necessary to invoke the base point $z$ of $X$, we write $H_1^w(X,z)$\index[symbol]{$H_1^w(X,z)$} and $p_{X,z}$\index[symbol]{$p_{X,z}$} instead of $H_1^w(X)$ and $p_X$, respectively. Given $z'\in X$ and a path $f$ in $X$ from $z$ to $z'$, the isomorphism $\widehat{f}\colon\pi_1(X,z)\to\pi_1(X,z')$ projects to an isomorphism $\widehat{f}\colon H_1^w(X,z)\to H_1^w(X,z')$\index[symbol]{$\widehat{f}$} with $p_{X,z'}\widehat{f}=\widehat{f}p_{X,z}$.

If $Y$ is a space and $h\colon X\to Y$ is a continuous map then the same symbol $h^{\sharp,z}$\index[symbol]{$h^{\sharp,z}$} (briefly, $h^{\sharp}$\index[symbol]{$h^{\sharp}$}) is used to denote the induced morphisms $\pi_1(X,z)\to\pi_1(Y,h(z))$ and $H_1^w(X,z)\to H_1^w(Y,h(z))$; thus,
\begin{equation*}
h^{\sharp,z}p_{X,z}=p_{Y,h(z)}h^{\sharp,z}\colon\pi_1(X,z)\to H_1^w(Y,h(z)).
\end{equation*}
Finally, given $z'\in X$ and a path $f$ in $X$ from $z$ to $z'$, we have
\begin{equation*}
h^{\sharp,z'}\widehat{f}=\widehat{hf}h^{\sharp,z}\colon H_1^w(X,z)\to H_1^w(Y,h(z')).
\end{equation*}

\subsection{Fundamental groups of manifolds and Lie groups}

Recall that the fundamental group $\pi_1(X)$ of a connected manifold $X$ is countable. If $X$ is also compact then the group $\pi_1(X)$ is finitely generated and hence $H_1^w(X)$ is a free abelian group of a finite rank; the rank of $H_1^w(X)$ is sometimes called the first Betti number\index{first Betti number of a space} of $X$.

We shall make use of the following results on homotopy of Lie groups due to Iwasawa and Ma\ml cev. Let $\Gamma\in\mathsf{LieGp}$ be a connected Lie group. Any two ma\-xi\-mal compact subgroups of $\Gamma$ are connected and conjugate in $\Gamma$. Moreover, if $K$ is a maximal compact subgroup of $\Gamma$ then the inclusion morphism $K\to\Gamma$ is a homotopy equivalence and so it induces isomorphisms of groups $\pi_1(K)\to\pi_1(\Gamma)$ and $H_1^w(K)\to H_1^w(\Gamma)$. It follows that the group $\pi_1(\Gamma)$ is finitely generated and that $H_1^w(\Gamma)$ is a free abelian group of a finite rank.

\subsection{Covering maps}

Let $E,B$ be topological spaces and $p\colon E\to B$ be a continuous surjective map. Recall that $p$ is called a covering map\index{covering map} if every element of $B$ has a neighbourhood $U$ such that the set $p^{-1}(U)$ is a disjoint union of open sets $V_i\subseteq E$ ($i\in I$), for which the restrictions $p\colon V_i\to U$ are homeomorphisms. In such a case we say that $U$ is evenly covered\index{evenly covered set} by $p$ and that $V_i$ is a $p$-slice\index{$p$-slice over a set} over $U$ for every $i\in I$. If the set $U$ is connected then there is a unique family of $p$-slices over $U$.

\subsection{Covering morphisms}

Let $G,G'$ be abelian topological groups and let $p\in\Hom(G',G)$ be an epimorphism. Then the following conditions are equivalent:
\begin{itemize}
\item $p$ is a covering morphism,
\item $p$ is a local homeomorphism,
\item $p$ is open (equivalently, quotient) with a discrete kernel $\ker(p)$.
\end{itemize}
If, in addition, both $G',G$ are locally compact and $\sigma$-compact then every epimorphism $G'\to G$ is open and so the following condition is equivalent to the previous three:
\begin{itemize}
\item $p$ has a discrete kernel $\ker(p)$.
\end{itemize}

\subsection{Lifts across covering maps}

Let $X,E,B$ be pointed spaces and let $f\colon X\to B$, $p\colon E\to B$ be continuous base point preserving maps. If $X\in\mathsf{CLAC}$ and $p$ is a covering map then the following conditions are equivalent:
\begin{itemize}
\item $f$ lifts across $p$ to a continuous base point preserving map $X\to E$,
\item $f^{\sharp}\pi_1(X)\subseteq p^{\sharp}\pi_1(E)$.
\end{itemize}
These conditions are fulfilled if $\pi_1(X)=0$, in which case we call the space $X$ simply connected\index{simply connected space}.

\subsection{Monodromy action}

Let $E,B$ be topological spaces and $p\colon E\to B$ be a covering map. Fix $b_0\in B$ and denote by $K$ the fibre $p^{-1}(b_0)$ of $p$ over $b_0$. The fundamental group $\pi_1(B,b_0)$ of $B$ acts from the right on the set $K$ by the rule
\begin{equation*}
K\times\pi_1(B,b_0)\ni(e,f)\mapsto\widetilde{f}_e(1)\in K,
\end{equation*}
where $\widetilde{f}_e$ stands for the lift of $f$ across $p$ starting at $e$ and $\widetilde{f}_e(1)$ is the endpoint of $\widetilde{f}_e$. This action is known as the \emph{monodromy action}\index{monodromy action} of $\pi_1(B,b_0)$ on $K$. For every $e\in K$, the isotropy group of $e$ with respect to this action is $p^{\sharp}\pi_1(E,e)$. Finally, if the group $\pi_1(B,b_0)$ is abelian (which is the case, for instance, if $B$ is a topological group and $b_0$ is the identity of $B$) then we may view the monodromy action as an action of $\pi_1(B,b_0)$ on $K$ from the left.

\subsection{Topological linear dimension}

If $L$ is a linear space over a field $\mathbb F$ then its linear dimension is denoted by $\dim_{\mathbb{F}\text{LS}}(L)$\index[symbol]{$\dim_{\mathbb{F}\text{LS}}(L)$}. If $B$ is a real Banach space then we use the symbol $\BSdim(B)$\index[symbol]{$\dim_{\mathbb{R}\text{BS}}(B)$} to denote the topological linear dimension\index{topological linear dimension} of $B$, that is, the smallest of the linear dimensions $\LSdim(L)$ of the dense linear subspaces $L$ of $B$. We recall the Kruse-Schmidt-Stone theorem\index{Kruse-Schmidt-Stone theorem}, which asserts that for every Banach space $B$, $\card(B)=d(B)^{\aleph_0}$, where $d(B)$ stands for the density of $B$. If the space $B$ is infinite-dimensional then $d(B)=\BSdim(B)$ and, by the Mackey theorem\index{Mackey theorem}, $\LSdim(B)\geq\mathfrak{c}=2^{\aleph_0}$.

\subsection{The space of real functions}\label{Sub:spc.real.fnctns}

Given a pointed compact space $X$ with the base point $z$, we denote by $C_z(X,\mathbb R)$\index[symbol]{$C_z(X,\mathbb R)$} the set of all continuous maps $f\colon X\to\mathbb R$ with $f(z)=0$. With operations defined point-wise and with the supremum norm, $C_z(X,\mathbb R)$ is a Banach space (in fact, a Banach algebra). Recall that if the space $X$ is infinite then
\begin{equation*}
\LSdim(C_z(X,\mathbb R))=\card(C_z(X,\mathbb R))=\w(X)^{\aleph_0};
\end{equation*}
consequently, there is an isomorphism of abelian groups $C_z(X,\mathbb R)\cong\mathbb R^{\left(\w(X)^{\aleph_0}\right)}$.

\section{Compact abelian groups}\label{S:top.comp.abel.gps}

Our references in this section are monographs \cite{HewRos1}, \cite{HofMor}, \cite{Mor} and \cite{Pon}.

\subsection{Co-completeness of $\mathsf{CAGp}$}\label{Sub:CAGp.co-cmplt}

The category $\mathsf{CAGp}$ of compact abelian groups is co-complete; we recall the usual construction of an inverse limit in $\mathsf{CAGp}$. Let $\mathfrak{G}=(G_i\stackrel{p_{ik}}{\longleftarrow}G_k)$ be an inverse system in $\mathsf{CAGp}$, indexed by a directed set $J$. Denote by $G$ the set of all $(g_j)_{j\in J}\in\prod_{j\in J}G_j$, such that $p_{ik}(g_k)=g_i$ for all $i\leq k$. Then $G$ is a closed subgroup of the product $\prod_{j\in J}G_j$ and hence $G\in\mathsf{CAGp}$. Moreover, the restrictions $p_j\colon G\to G_j$ of the projection morphisms $\pr_j$ satisfy $p_{ik}p_k=p_i$ for all $i\leq k$. As is well known, the group $G$ is an inverse limit of $\mathfrak{G}$ in $\mathsf{CAGp}$ and $p_j$ ($j\in J$) are the limit projections. Notice that the morphisms $p_j$ separate the points of $G$. Since an inverse limit of $\mathfrak{G}$ in $\mathsf{CAGp}$ is unique up to isomorphism, it follows that the limit projections of an arbitrary inverse limit of $\mathfrak{G}$ in $\mathsf{CAGp}$ separate points.

\subsection{Pontryagin duality}\label{Sub:Pontr.dual}

Let $G\in\mathsf{LCAGp}$. The Pontryagin dual\index{Pontryagin dual!of a group} of $G$ will be denoted by $G^*$\index[symbol]{$G^*$}. Recall that $G^*$ consists of the topological morphisms $\chi\colon G\to\mathbb T^1$, called the characters\index{character of a group} of $G$. With the compact-open topology and point-wise defined operations, $G^*$ is a locally compact abelian group. If $G$ is compact then $G^*$ is discrete and if $G$ is discrete then $G^*$ is compact. Moreover, $G^*$ is second countable if and only if $G$ is second countable. The groups $G$ and $G^{**}$ are topologically isomorphic via the isomorphism $\omega_G\colon G\to G^{**}$,\index[symbol]{$\omega_G$} acting by the rule $\omega_G(g)\colon G^*\ni\chi\mapsto\chi(g)\in\mathbb T^1$ for every $g\in G$. We call $\omega_G$ the Pontryagin isomorphism\index{Pontryagin isomorphism} associated to $G$. Recall also that for every $G\in\mathsf{CAGp}$, $\dim(G)=\rank(G^*)$ and $\w(G)=\card(G^*)$.

Given $G\in\mathsf{LCAGp}$ and a subset $A$ of $G$, the annihilator\index{annihilator of a set} of $A$ in $G^*$ is denoted by $[A^{\perp},G^*]$\index[symbol]{$[A^{\perp},G^*]$} or, if no misunderstanding can arise, by a shorter symbol $A^{\perp}$\index[symbol]{$A^{\perp}$}. Recall that $A^{\perp}$ is formed by the characters $\chi\in G^*$ with $A\subseteq\ker(\chi)$. The annihilator of $A$ coincides with the annihilator of the closed subgroup of $G$ generated by $A$. Moreover, $A^{\perp}$ is always a closed subgroup of $G^*$. If $H\sbgp G$ then $\omega_G(H)=H^{\perp\perp}$ and $\omega_{G^*}(H^{\perp})=\omega_G(H)^{\perp}$. Given a system of closed subgroups $H_i$ ($i\in I$) of $G$, there are equalities $\bigcap_{i\in I}H_i^{\perp}=\left(\sum_{i\in I}H_i\right)^{\perp}$ and $\left(\bigcap_{i\in I}H_i\right)^{\perp}=\overline{\sum_{i\in I}H_i^{\perp}}$. If $G\in\mathsf{LCAGp}$ is either compact or discrete, $H,K\sbgp G$ and there is a topological direct sum $G=H\oplus K$, then there is also a topological direct sum $G^*=H^{\perp}\oplus K^{\perp}$.

Let $G_1,G_2\in\mathsf{LCAGp}$ and $h\colon G_1\to G_2$ be a topological morphism. The Pontryagin dual\index{Pontryagin dual!of a morphism} $h^*$\index[symbol]{$h^*$} of $h$ is the topological morphism $G_2^*\to G_1^*$, acting by the rule $h^*(\Upsilon)=\Upsilon h$ for every $\Upsilon\in G_2^*$. One has $\omega_{G_2}h=h^{**}\omega_{G_1}$, $\im(h)^{\perp}=\ker(h^*)$ and $\ker(h)^{\perp}=\overline{\im(h^*)}$. In particular, $h$ is a monomorphism if and only if $h^*$ has a dense image and, conversely, $h$ has a dense image if and only if $h^*$ is a monomorphism. More generally, if $H_1\sbgp G_1$ and $H_2\sbgp G_2$ then $h(H_1)^{\perp}=(h^*)^{-1}(H_1^{\perp})$ and $h^{-1}(H_2)^{\perp}=\overline{h^*(H_2^{\perp})}$. Notice also that $(\omega_G)^*=(\omega_{G^*})^{-1}$ for every $G\in\mathsf{LCAGp}$.

Let $G_i\stackrel{p_{ik}}{\longleftarrow}G_k$ be an inverse system in $\mathsf{CAGp}$, $G$ be its inverse limit and $p_i\colon G\to G_i$ be the limit projections. Then $G_i^*\stackrel{p_{ik}^*}{\longrightarrow}G_k^*$ is a direct system in $\mathsf{DAGp}$, $G^*$ is its direct limit and $p_i^*\colon G_i^*\to G^*$ are the limit morphisms. We also recall that for every family $G_i$ ($i\in I$) of compact abelian groups, there is an isomorphism of groups $(\prod_{i\in I}G_i)^*\cong\bigoplus_{i\in I}G_i^*$. Finally, if $G_1,\dots,G_n\in\mathsf{LCAGp}$ then there is a topological isomorphism $(G_1\times\dots\times G_n)^*\cong G_1^*\times\dots\times G_n^*$.

Let $G\in\mathsf{LCAGp}$, $H\sbgp G$ and consider the corresponding short exact sequence in $\mathsf{LCAGp}$
\begin{equation}\label{Eq:seq.H.G.G/H}
0\longrightarrow H\stackrel{j}{\longrightarrow}G\stackrel{p}{\longrightarrow}G/H\longrightarrow0,
\end{equation}
where $j$ is the inclusion morphism and $p$ is the quotient morphism. Then the dual sequence\index{dual sequence}
\begin{equation}\label{Eq:seq.H.G.G/H.prdl}
0\longrightarrow(G/H)^*\stackrel{p^*}{\longrightarrow}G^*\stackrel{j^*}{\longrightarrow}H^*\longrightarrow0
\end{equation}
is also short exact in $\mathsf{LCAGp}$. Moreover, $p^*$ restricts to a topological isomorphism
\begin{equation*}
p^*\colon(G/H)^*\to \im(p^*)=H^{\perp}.
\end{equation*}
Consequently, there is a topological isomorphism $(G/H)^*\cong H^{\perp}$ and the sequence (\ref{Eq:seq.H.G.G/H.prdl}) takes the form
\begin{equation}\label{Eq:seq.H.G.G/H.dual}
0\longrightarrow H^{\perp}\stackrel{\sigma}{\longrightarrow}G^*\stackrel{\pi}{\longrightarrow}H^*\longrightarrow0,
\end{equation}
where $\sigma$ is the inclusion morphism and $\pi=j^*$. In particular, there is a topological isomorphism $G^*/H^{\perp}\cong H^*$.

Given $G\in\mathsf{LCAGp}$, we use the symbol $\mathcal{L}(G)$\index[symbol]{$\mathcal L(G)$} to denote the \emph{Lie algebra}\index{Lie algebra} of $G$. Recall that $\mathcal{L}(G)=\Hom(\mathbb R,G)$ is formed by the one-parameter subgroups of $G$. With the compact-open topology and point-wise operations, $\mathcal L(G)$ is an abelian topological group. The \emph{exponential function}\index{exponential function} $\exp_G$\index[symbol]{$\exp_G$} of $G$ is a morphism of topological groups $\mathcal L(G)\to G$, acting by the rule $\exp_G(h)=h(1)$. (If no misunderstanding can arise concerning the group $G$, we write $\exp$\index[symbol]{$\exp$} instead of $\exp_G$.) The kernel $\mathcal K(G)$\index[symbol]{$\mathcal K(G)$} of $\exp_G$ is topologically isomorphic to $\Hom(\mathbb T^1,G)$; indeed, if $p\colon\mathbb R\to\mathbb T^1$ is the usual covering morphism then $\Hom(p,G)\colon\Hom(\mathbb T^1,G)\ni k\mapsto kp\in\mathcal K(G)$ is a topological isomorphism.

\subsection{Connectedness and identity arc-component}

Let $G\in\mathsf{CAGp}$. Then for every $k\in\mathbb N$, $\tor_k(G^*)=(\kappa_k(G))^{\perp}$ and $\kappa_k(G^*)=\tor_k(G)^{\perp}$.
Also, $\bigcap_{k\in\mathbb N}\kappa_k(G^*)=\tor(G)^{\perp}$ and $\tor(G^*)=\Div(G)^{\perp}$. (Here, $\Div(G)$ is always a closed subgroup of $G$, but $\tor(G)$ need not be closed in $G$.) Consequently, $G$ is divisible, respectively, torsion-free if and only if $G^*$ is torsion-free, respectively, divisible. Further, there are identities
\begin{equation*}
\bigcap_{k\in\mathbb N}\kappa_k(G)=\Div(G)=G_0=\overline{G_a}=\overline{\exp(\mathcal L(G))},
\end{equation*}
where $G_0$\index[symbol]{$G_0$} is the identity component of $G$ and $G_a$\index[symbol]{$G_a$} is the identity arc-component of $G$. In fact, there is an exact sequence of abelian groups
\begin{equation}\label{Eq:ex.seq.Lie.alg}
0\longrightarrow \mathcal K(G)\stackrel{j}{\longrightarrow}\mathcal{L}(G)\stackrel{\exp}{\longrightarrow}G\stackrel{p_a}{\longrightarrow}G/G_a\longrightarrow0,
\end{equation}
where $j$ is the inclusion morphism and $p_a$ is the quotient morphism. In particular, the identity arc-component $G_a$ of $G$ consists precisely of the points lying on one-parameter subgroups of $G$. The group $G$ is connected, respectively, totally disconnected if and only if its dual group $G^*$ is torsion-free, respectively, a torsion group. Further, the sequence (\ref{Eq:ex.seq.Lie.alg}) is equivalent to the covariant Hom-Ext sequence
\begin{equation*}
0\longrightarrow\Hom(G^*,\mathbb Z)\longrightarrow\Hom(G^*,\mathbb R)\longrightarrow\Hom(G^*,\mathbb T^1)\longrightarrow\Ext(G^*,\mathbb Z)\longrightarrow0,
\end{equation*}
derived from the natural short exact sequence
\begin{equation*}
0\longrightarrow\mathbb Z\longrightarrow\mathbb R\longrightarrow\mathbb T^1\longrightarrow0
\end{equation*}
and associated to the group $G^*$. In particular, there is an isomorphism of groups $G/G_a\cong\Ext(G^*,\mathbb Z)$. Thus, $G$ is arc-wise connected if and only if $G^*$ is a Whitehead group. Finally, recall that every epimorphism $q\colon G\to H$ of compact abelian groups induces an epimorphism $q\colon G_a\to H_a$.

\subsection{Maximal toral quotients}

Let $G\in\mathsf{CAGp}$ be a connected group and let $r=\dim(G)=\rank(G^*)$. Fix a maximal independent set $S\subseteq G^*$ and let $F=\langle S\rangle$ be the subgroup of $G^*$ generated by $S$. Clearly, $F$ is a free abelian group with the basis $S$, $\rank(F)=\card(S)=r$ and $G^*/F$ is a torsion group. By dualising the corresponding short exact sequence in $\mathsf{DAGp}$
\begin{equation*}
0\longrightarrow F\longrightarrow G^*\longrightarrow G^*/F\longrightarrow0,
\end{equation*}
we obtain a short exact sequence in $\mathsf{CAGp}$
\begin{equation}\label{Eq:max.tor.quot.seq}
0\longrightarrow \mathbb K(G)\stackrel{j}{\longrightarrow}G\stackrel{q}{\longrightarrow}\mathbb T(G)\longrightarrow0.
\end{equation}
Here, $\mathbb{K}(G)\cong(G^*/F)^*$\index[symbol]{$\mathbb{K}(G)$} is totally disconnected and $\mathbb{T}(G)\cong F^*\cong\mathbb T^r$\index[symbol]{$\mathbb{T}(G)$}. Clearly, $r$ is the largest of all the cardinals $\mathfrak{k}$, such that $G$ factors onto $\mathbb T^{\mathfrak{k}}$. Therefore, we shall refer to (\ref{Eq:max.tor.quot.seq}) as a \emph{maximal toral quotient sequence}\index{maximal toral quotient sequence} of $G$. Finally, if $G,H\in\mathsf{CAGp}$ are connected groups and there exist epimorphisms $G\to H\to\mathbb T(G)$, then $\dim(G)=\dim(H)$ and hence $\mathbb T(G)\cong\mathbb T(H)$.

\subsection{Maximal toral subgroup}

Let $A\in\mathsf{AbGp}$ be a countable torsion-free group. Then $A$ decomposes into a direct sum $A=\mathfrak{f}(A)\oplus\mathfrak{t}(A)$, where $\mathfrak{f}(A)$\index[symbol]{$\mathfrak{f}(A)$} is a free abelian group and the group $\mathfrak{t}(A)$\index[symbol]{$\mathfrak{t}(A)$} is torsion-less\index{torsion-less group} (that is, $\Hom(\mathfrak{t}(A),\mathbb Z)=0$ or, equivalently, $\mathfrak{t}(A)$ contains no non-trivial free abelian group as a direct summand). These ideas translate to the category of compact abelian groups as follows. Let $G\in\mathsf{CAGp}$ be a connected second countable group. Then there is a topological direct sum $G=\mathfrak{f}(G)\oplus\mathfrak{t}(G)$, where $\mathfrak{f}(G)$\index[symbol]{$\mathfrak{f}(G)$} is a torus and $\mathfrak{t}(G)$\index[symbol]{$\mathfrak{t}(G)$} is a torus-free group\index{torus-free group} (that is, $\Hom(\mathbb T^1,\mathfrak{t}(G))=0$). If, in addition, $X$ is a (pointed) topological space and $\xi\colon X\to G$ is a continuous (base point preserving) map then $\xi=\mathfrak{f}(\xi)\oplus\mathfrak{t}(\xi)$, where $\mathfrak{f}(\xi)$\index[symbol]{$\mathfrak{f}(\xi)$} and $\mathfrak{t}(\xi)$\index[symbol]{$\mathfrak{t}(\xi)$} are continuous (base point preserving) maps, obtained as compositions of $\xi$ with the projection morphisms from $G$ onto $\mathfrak{f}(G)$ and $\mathfrak{t}(G)$, respectively.

\subsection{Bohr compactification}\label{Sub:Bohr.compctf}

Let $G\in\mathsf{LCAGp}$. The Bohr compactification\index{Bohr!compactification} of $G$ consists of a group $bG\in\mathsf{CAGp}$\index[symbol]{$bG$} and a topological morphism $p\in\Hom(G,bG)$, such that for every $K\in\mathsf{CAGp}$ and every $q\in\Hom(G,K)$ there is a unique $r\in\Hom(bG,K)$ with $rp=q$. The Bohr compactification of $G$ can be viewed as a universal element for the covariant functor $\Hom(G,-)\colon\mathsf{CAGp}\to\mathsf{Sets}$ and is therefore unique up to a natural equivalence. We shall use a model of $(bG,p)$ defined as follows: let $bG=(G^*)_d^*$ be the dual group of the discretized dual of $G$ and $p=\sigma^*\omega_G$, where $\omega_G$ is the Pontryagin isomorphism associated to $G$ and $\sigma\colon(G^*)_d\to G^*$ is the identity morphism. The morphism $p$ has a dense image. The group $b\mathbb Z$ is the universal compact monothetic group\index{universal!compact monothetic group} and $b\mathbb R$ is the universal compact solenoidal group\index{universal!compact solenoidal group}. We also have a topological isomorphism $(b\mathbb R)^{\mathfrak{c}}\cong b\mathbb R$.

Further, recall that the initial topology on $G$ derived from $bG$ and $p$ is called the Bohr topology\index{Bohr!topology} of $G$. It has a basis consisting of the sets of the form $\bigcap_{k=1}^n\chi_k^{-1}(V_k)$, where $\chi_k\in G^*$, $V_k\subseteq\mathbb T^1$ are open and $n\in\mathbb N$. If the group $G$ is compact then $p=\omega_G$ is a topological isomorphism and so the Bohr topology of $G$ coincides with the original topology of $G$. In particular, for every identity neighbourhood $W$ in $G$ there exist $n\in\mathbb N$, $\chi_k\in G^*$ and $1\in V_k\subseteq\mathbb T^1$ open ($k=1,\dots,n$), such that $\bigcap_{k=1}^n\chi_k^{-1}(V_k)\subseteq W$.

\subsection{Maps into compact abelian groups}\label{Sub:mps.cpt.ab.gps}

Let $X$ be a pointed topological space with the base point $z$ and $G$ be an abelian topological group with the identity $e$. We denote by $C_z(X,G)$\index[symbol]{$C_z(X,G)$} the set of all continuous maps $f\colon X\to G$ with $f(z)=e$. With operations defined point-wise, $C_z(X,G)$ is an abelian group. If the space $X$ is locally compact then we equip the group $C_z(X,G)$ with the compact-open topology, thus turning it into an abelian topological group. (This topology on $C_z(X,G)$ is in fact that of a uniform convergence on compact sets.) Every morphism of abelian topological groups $h\colon G\to H$ gives rise to a morphism of groups
\begin{equation*}
\widehat{h}\index[symbol]{$\widehat{h}$}\colon C_z(X,G)\ni f\mapsto hf\in C_z(X,H),
\end{equation*}
which is topological if the space $X$ is locally compact.

For a general pointed space $X$ and a compact abelian group $G$, consider the mapping
\begin{equation}\label{Eq:maps.are.homs}
\varphi\colon C_z(X,G)\ni f\mapsto\varphi(f)=f^*\in\Hom(G^*,C_z(X,\mathbb T^1)),
\end{equation}
where $f^*$\index[symbol]{$f^*$} associates with every $\chi\in G^*$ the map $\chi f\in C_z(X,\mathbb T^1)$. The map $\varphi$ is an isomorphism of groups. If the space $X$ is locally compact and the group $\Hom(G^*,C_z(X,\mathbb T^1))$ carries the topology of point-wise convergence then $\varphi$ is an isomorphism of topological groups.

\subsection{First cohomotopy group}\label{Sub:frst.chmtp.gp}

Let $X$ be a pointed locally compact space with the base point $z$ and $G$ be an abelian topological group. Two maps $f_1,f_2\in C_z(X,G)$ are homotopic if and only if they are connected by a path in $C_z(X,G)$, that is, if and only if $f_1-f_2$ is an element of the identity arc-component $C_z(X,G)_a$ of $C_z(X,G)$ (we shall call the maps from $C_z(X,G)_a$ null-homotopic\index{null-homotopic map}). Consequently, the set $[X,G]$\index[symbol]{$[X,G]$} of homotopy classes of continuous base point preserving maps $X\to G$ coincides with the set of cosets $C_z(X,G)/C_z(X,G)_a$ and so it carries the structure of an abelian group. If the space $X$ is compact and the group $G$ is locally contractible then all maps $f\in C_z(X,G)$ sufficiently close to the identity $e$ are null-homotopic and hence are contained in $C_z(X,G)_a$. Consequently, $C_z(X,G)_a$ is an open subgroup of $C_z(X,G)$ and it coincides with the identity component $C_z(X,G)_0$ of $C_z(X,G)$. (This applies to the case when $G$ is a Lie group and, in particular, when $G=\mathbb T^1$.)

Now let the space $X$ be compact connected. The usual covering morphism $p\colon\mathbb R\to\mathbb T^1$ induces a monomorphism $\widehat{p}\colon C_z(X,\mathbb R)\to C_z(X,\mathbb T^1)$. The image $\im(\widehat{p})$ of $\widehat{p}$ coincides with $C_z(X,\mathbb T^1)_a=C_z(X,\mathbb T^1)_0$ and it is topologically isomorphic to $C_z(X,\mathbb R)$ via $\widehat{p}$. Consequently, $C_z(X,\mathbb T^1)_a$ is a divisible group and hence it is a direct summand in $C_z(X,\mathbb T^1)$. Fix a complementary summand $\pi^1(X)$\index[symbol]{$\pi^1(X)$} of $C_z(X,\mathbb T^1)_a$ in $C_z(X,\mathbb T^1)$. Then $\pi^1(X)$ is isomorphic to the first cohomotopy group $[X,\mathbb T^1]$ of $X$, the isomorphism associating with each map $f\in\pi^1(X)$ its homotopy class $[f]\in[X,\mathbb T^1]$\index[symbol]{$[f]$}. The groups $\pi^1(X)$ and $[X,\mathbb T^1]$ will often be identified in this work. Recall that $C_z(X,\mathbb T^1)_a$ is an open subgroup of $C_z(X,\mathbb T^1)$. Therefore, the group $\pi^1(X)$ is discrete and the splitting $C_z(X,\mathbb T^1)=C_z(X,\mathbb T^1)_a\oplus\pi^1(X)$ is topological. Also, since $C_z(X,\mathbb T^1)$ is a torsion-free group, it follows that so is $\pi^1(X)$. Finally, if the space $X$ is second countable then $C_z(X,\mathbb T^1)$ is a Polish group and so the group $\pi^1(X)$ is countable.

If $X,Y$ are compact connected pointed spaces and $\xi\colon X\to Y$ is a con\-ti\-nuous base point preserving map then we denote by $\xi^{\flat}$\index[symbol]{$\xi^{\flat}$} the induced morphism $[Y,\mathbb T^1]\to[X,\mathbb T^1]$, acting by the rule $\xi^{\flat}([f])=[f\xi]$. Clearly, the assignment $\xi\mapsto\xi^{\flat}$ is contravariant functorial.

\subsection{Inducing morphisms of the first weak homology groups}\label{Sub:ind.mrph.hmlg.gp}

Let $X$ be a compact connected manifold and $z\in X$. Under the usual identification $H_1^w(\mathbb T^1)\cong\mathbb Z$, there is a morphism of groups
\begin{equation*}
\phi\colon C_z(X,\mathbb T^1)\ni f\mapsto f^{\sharp}\in\Hom(H_1^w(X),\mathbb Z).
\end{equation*}
Since $X$ has the homotopy type of a CW complex, $\phi$ is an epimorphism. The kernel $\ker(\phi)$ of $\phi$ coincides with the identity arc-component $C_z(X,\mathbb T^1)_a$ of $C_z(X,\mathbb T^1)$. Consequently, $\phi$ restricts to an isomorphism $\phi\colon\pi^1(X)\to\Hom(H_1^w(X),\mathbb Z)$. Since $H_1^w(X)$ is a free abelian group of a finite rank, it follows that 
\begin{equation*}
\pi^1(X)\cong\Hom(H_1^w(X),\mathbb Z)\cong H_1^w(X).
\end{equation*}

\subsection{First (co)homotopy of compact abelian groups}\label{Sub:chmtp.cpt.gps}

Let $X$ be a compact connected pointed space with the base point $z$ and $G$ be a compact abelian group. The topological splitting $C_z(X,\mathbb T^1)=C_z(X,\mathbb T^1)_a\oplus\pi^1(X)$ yields a to\-po\-lo\-gi\-cal splitting
\begin{equation*}
\Hom(G^*,C_z(X,\mathbb T^1))=\Hom(G^*,C_z(X,\mathbb T^1)_a)\oplus\Hom(G^*,\pi^1(X)).
\end{equation*}
Since $C_z(X,\mathbb T^1)_a$ is topologically isomorphic to $C_z(X,\mathbb R)$ and the latter group is the additive topological group of a real Banach space, it follows that $\Hom(G^*,C_z(X,\mathbb T^1)_a)$ is an arc-wise connected group. Also, since $\pi^1(X)$ is a discrete group, it follows that the group $\Hom(G^*,\pi^1(X))$ is totally disconnected. From these two observations we obtain that
\begin{enumerate}
\item[(1)] the group $\Hom(G^*,C_z(X,\mathbb T^1)_a)$ coincides with the identity arc-component of the group $\Hom(G^*,C_z(X,\mathbb T^1))$.
\end{enumerate}
Further, since the quotient morphism $C_z(X,\mathbb T^1)\to[X,\mathbb T^1]$ possesses a cross-section, it follows that
\begin{enumerate}
\item[(2)] the induced morphism $\Hom(G^*,C_z(X,\mathbb T^1))\to\Hom(G^*,[X,\mathbb T^1])$ is an epimorphism.
\end{enumerate}
From statements (1) and (2) it follows that the isomorphism $\varphi$ from (\ref{Eq:maps.are.homs}) gives rise to an isomorphism
\begin{equation}\label{Eq:maps.are.homs.hmtp}
\widetilde{\varphi}\colon[X,G]\ni[f]\mapsto\widetilde{\varphi}([f])\in\Hom(G^*,[X,\mathbb T^1]),
\end{equation}
where $\widetilde{\varphi}([f])$ associates with every $\chi\in G^*$ the homotopy class $[\chi f]\in[X,\mathbb T^1]$.

Given $G\in\mathsf{CAGp}$, there are isomorphisms of groups
\begin{equation}\label{Eq:fund.gp.CAGp}
\pi_1(G)\cong[\mathbb T^1,G]\cong\Hom(G^*,\mathbb Z)\cong\Hom(\mathbb T^1,G);
\end{equation}
indeed, the first isomorphism is clear, the second one follows from (\ref{Eq:maps.are.homs.hmtp}) applied to $X=\mathbb T^1$ and the third one is a consequence of the Pontryagin duality. In fact, every homotopy class from $[\mathbb T^1,G]$ contains one and only one morphism from $\Hom(\mathbb T^1,G)$. It follows, in particular, that the inclusion morphism $\overline{\tor(G)}\to G$ induces an isomorphism of fundamental groups $\pi_1\big(\overline{\tor(G)}\big)\to\pi_1(G)$. Further, from (\ref{Eq:fund.gp.CAGp}) we obtain an isomorphism of groups $\mu\colon\pi_1(G)\to\Hom(G^*,\pi_1(\mathbb T^1))$, acting by the rule $\mu(l)\colon\chi\mapsto\chi^{\sharp}(l)$ for every $l\in\pi_1(G)$; consequently, the morphisms $\chi^{\sharp}\colon\pi_1(G)\to\pi_1(\mathbb T^1)$ ($\chi\in G^*$) separate the points of $\pi_1(G)$. Finally, let $G,H\in\mathsf{CAGp}$ and $h\in\Hom(G,H)$. In Figure~\ref{Fig:h.sharp.chng} below we indicate how the identifications from (\ref{Eq:fund.gp.CAGp}) affect the induced morphism $h^{\sharp}\colon\pi_1(G)\to\pi_1(H)$.
\begin{figure}[ht]
\[\minCDarrowwidth20pt\begin{CD}
\pi_1(G) @= \Hom(G^*,\mathbb Z) @= \Hom(\mathbb T^1,G)\\
@VVh^{\sharp}V @VV\Hom(h^*,\mathbb Z)V @VV\Hom(\mathbb T^1,h)V \\
\pi_1(H) @= \Hom(H^*,\mathbb Z) @= \Hom(\mathbb T^1,H)
\end{CD}\]
\caption{Induced morphism of the fundamental groups}
\label{Fig:h.sharp.chng}
\end{figure}

If $G\in\mathsf{CAGp}$ is a connected group then every homotopy class from $[G,\mathbb T^1]$ contains a unique character of $G$. Consequently, there are isomorphisms of groups
\begin{equation*}
\pi^1(G)\cong[G,\mathbb T^1]\cong G^*.
\end{equation*}
Recall the following important corollary of these isomorphisms: two compact connected abelian groups are topologically isomorphic if and only if they are homeomorphic as topological spaces. Finally, if $G,H\in\mathsf{CAGp}$ are connected groups and $h\in\Hom(G,H)$ then the induced morphism $h^{\flat}\colon\pi^1(H)\to\pi^1(G)$ turns into the Pontryagin dual $h^*\colon H^*\to G^*$ of $h$ under the identifications $\pi^1(G)\cong G^*$ and $\pi^1(H)\cong H^*$.

\subsection{Free compact abelian groups}\label{Sub:free.cpt.abel.gp}

Let $X$ be a pointed topological space with the base point $z$, $G_X$ be a compact abelian group and $\varphi_X\colon X\to G_X$ be a continuous base point preserving map. We say that $G_X$\index[symbol]{$G_X$} is a \emph{free compact abelian group\index{free!compact abelian group} over $X$ with respect to the universal map\index{universal!map} $\varphi_X$\index[symbol]{$\varphi_X$}}, if for every group $H\in\mathsf{CAGp}$ and every continuous base point preserving map $\psi\colon X\to H$ there is a unique $h\in\Hom(G_X,H)$ with $h\varphi_X=\psi$. For every $X$, such a pair $(G_X,\varphi_X)$ is unique up to a natural equivalence. Given a pointed space $X$, let $G_X=C_z(X,\mathbb T^1)_d^*$ be the Pontryagin dual to the discretized group $C_z(X,\mathbb T^1)$ and let
\begin{equation*}
\varphi_X\colon X\ni x\mapsto \left(\varphi_X(x)\colon C_z(X,\mathbb T^1)_d\ni f\mapsto f(x)\in\mathbb T^1\right)\in G_X;
\end{equation*}
then $G_X$ is a free compact abelian group over $X$ with respect to the universal map $\varphi_X$. Recall that if the space $X$ is compact then there is a topological isomorphism
\begin{equation}\label{Eq:free.comp.gp}
G_X\cong\left((\mathbb Q_d)^*\right)^{\w(X)^{\aleph_0}}\times\prod_{\substack{p\,\text{prime}}}\left(\mathbb Z_p\right)^{\w(\widetilde{X})-1}\times[X,\mathbb T^1]^*,
\end{equation}
where $\widetilde{X}$ denotes the space of the connected components of $X$.

\subsection{Projective resolution} 

A group $P\in\mathsf{CAGp}$ is called \emph{projective}\index{projective!compact abelian group} if for all $K,H\in\mathsf{CAGp}$, every epimorphism $q\in\Hom(K,H)$ and every $h\in\Hom(P,H)$, there is $k\in\Hom(P,K)$ with $qk=h$. The projective groups $P\in\mathsf{CAGp}$ coincide with the torsion-free elements of $\mathsf{CAGp}$. Given a connected group $G\in\mathsf{CAGp}$, consider the injective resolution of $G^*$ in $\mathsf{DAGp}$
\begin{equation*}
0\longrightarrow G^*\stackrel{j}{\longrightarrow}\mathbb Q\otimes G^*\stackrel{p}{\longrightarrow}(\mathbb Q/\mathbb Z)\otimes G^*\longrightarrow0.
\end{equation*}
By dualising this sequence, we obtain a short exact sequence in $\mathsf{CAGp}$
\begin{equation}\label{Eq:proj.res.G}
0\longrightarrow\mathfrak{K}(G)\stackrel{\mathfrak{i}}{\longrightarrow}\mathfrak{P}(G)\stackrel{\mathfrak{p}}{\longrightarrow}G\longrightarrow0.\index[symbol]{$\mathfrak{K}(G)$}\index[symbol]{$\mathfrak{P}(G)$}\index[symbol]{$\mathfrak{i}$}\index[symbol]{$\mathfrak{p}$}
\end{equation}
We shall refer to (\ref{Eq:proj.res.G}) as the \emph{projective resolution\index{projective!resolution} of $G$}. The following conditions hold:
\begin{itemize}
\item $\mathfrak{P}(G)$ is a projective group and $\mathfrak{p}\colon\mathfrak{P}(G)\to G$ is an epimorphism,
\item given a projective group $P\in\mathsf{CAGp}$ and an epimorphism $p\in\Hom(P,G)$, there is a unique epimorphism $q\in\Hom(P,\mathfrak{P}(G))$ with $\mathfrak{p}q=p$.
\end{itemize}
These two conditions are usually expressed by saying that $\mathfrak{p}\colon\mathfrak{P}(G)\to G$ is the \emph{projective cover\index{projective!cover} of $G$}. Notice that the group $\mathfrak{K}(G)\cong((\mathbb Q/\mathbb Z)\otimes G^*)^*$ is totally disconnected. Also,
\begin{equation*}
\dim(\mathfrak{P}(G))=\rank(\mathbb Q\otimes G^*)=\rank(G^*)=\dim(G).
\end{equation*}
Finally, if $G,K\in\mathsf{CAGp}$ are connected groups and there exist epimorphisms $\mathfrak{P}(G)\to K\to G$, then $\dim(G)=\dim(K)$ and hence there is a topological isomorphism $\mathfrak{P}(G)\cong\mathfrak{P}(K)$.

\subsection{Solenoids}\label{Sub:slnds}

Let ${\bf p}=(p_n)_{n\in\mathbb N}$ be a sequence of integers greater than $1$ and let $p_0=1$. We denote by $\mathbb Z/{\bf p}=\{k/p_0\dots p_m : k\in\mathbb Z, m\geq0\}$\index[symbol]{$\mathbb Z/{\bf p}$} the set of all $\bf{p}$-adic rationals\index{${\bf p}$-adic!rational}. Clearly, $\mathbb Z/{\bf p}$ is a non-trivial subgroup of $\mathbb Q$ and so it has rank $1$. Further, we denote by $S_{{\bf p}}$\index[symbol]{$S_{{\bf p}}$} the ${\bf p}$-adic solenoid\index{${\bf p}$-adic!solenoid}
\begin{equation*}
S_{{\bf p}}=\lim_{\longleftarrow}\left(\mathbb T^1\stackrel{\kappa_{p_n}}{\longleftarrow}\mathbb T^1\right)_{n\in\mathbb N}.
\end{equation*}
By definition, $S_{{\bf p}}$ is a closed subgroup of $\mathbb T^{\mathbb N}$. For every $n\in\mathbb N$, the symbol $\pr_n$ stands for the projection morphism $S_{{\bf p}}\to\mathbb T^1$ on the $n^{\text{th}}$ coordinate. Recall an isomorphism of groups
\begin{equation}\label{Eq:dual.gp.sole.iso}
\lambda\colon\mathbb Z/{\bf p}\ni\frac{k}{p_0\dots p_m}\mapsto\chi_{k,m+1}\in(S_{{\bf p}})^*,
\end{equation}
where $\chi_{k,m+1}$ acts by the rule $\chi_{k,m+1}(z)=\pr_{m+1}(z)^k$ for all $k\in\mathbb Z$, $m\geq0$ and $z\in S_{\bf p}$. Consider the morphism $h\in\mathcal L(S_{{\bf p}})$, given by $h(t)=(\exp(i2\pi t/p_0\dots p_{n-1}))_{n\in\mathbb N}$ for every $t\in\mathbb R$. Then every element $g$ of $\mathcal L(S_{{\bf p}})$ is of the form $g(t)=h(st)$ for some $s\in\mathbb R$ and so the identity arc-component $(S_{{\bf p}})_a$ of $S_{{\bf p}}$ equals $(S_{{\bf p}})_a=\exp(\mathcal L(S_{{\bf p}}))=\im(h)$. Consequently, $(S_{\bf p})_a\cap\tor(S_{\bf p})=e$.

\subsection{Groups of topological morphisms}

Given groups $G,H\in\mathsf{LCAGp}$, denote by $\Hom(G,H)$ the set of all topological morphisms $G\to H$. With the point-wise operations and the compact-open topology, $\Hom(G,H)$ is an abelian topological group. If $K\sbgp H$ then $\Hom(G,K)\sbgp\Hom(G,H)$. Given $G,H_1,\dots,H_n\in\mathsf{LCAGp}$, there is a topological isomorphism $\Hom(G,\prod_{i=1}^nH_i)\cong\prod_{i=1}^n\Hom(G,H_i)$, which associates with every morphism $G\to\prod_{i=1}^nH_i$ the $n$-tuple of its coordinate morphisms $G\to H_i$. If, in addition, $K_i\sbgp H_i$ for $i=1,\dots,n$ then this isomorphism gives rise to a topological isomorphism
\begin{equation*}
\Hom\left(G,\prod_{i=1}^nH_i\right)\Big/\Hom\left(G,\prod_{i=1}^nK_i\right)\cong\prod_{i=1}^n\Hom(G,H_i)/\Hom(G,K_i).
\end{equation*}

Let $G_1,G_2,H\in\mathsf{LCAGp}$ and $h\in\Hom(G_1,G_2)$. Then there is a topological morphism 
\begin{equation*}
\Hom(h,H)\colon\Hom(G_2,H)\ni\varphi\mapsto\varphi h\in\Hom(G_1,H)\index[symbol]{$\Hom(h,H)$}.
\end{equation*}
Clearly, the assignment $h\mapsto\Hom(h,H)$ is contravariant functorial for every $H\in\mathsf{LCAGp}$. Similarly, given $G,H_1,H_2\in\mathsf{LCAGp}$ and $k\in\Hom(H_1,H_2)$, there is a topological morphism
\begin{equation*}
\Hom(G,k)\colon\Hom(G,H_1)\ni\varphi\mapsto k\varphi\in\Hom(G,H_2)\index[symbol]{$\Hom(G,k)$}.
\end{equation*}
This time, the assignment $k\mapsto\Hom(G,k)$ is covariant functorial for every $G\in\mathsf{LCAGp}$.

Finally, if $G,H\in\mathsf{LCAGp}$ then there is a topological isomorphism $\Hom(G,H)\to\Hom(H^*,G^*)$, which associates with every $h\in\Hom(G,H)$ its Pontryagin dual $h^*\in\Hom(H^*,G^*)$.

\subsection{Fr\'echet and Banach spaces of morphisms}

Our aim in this subsection is to prove the following lemma.

\begin{lemma}\label{L:hom.Ban.Fre}
Let $A$ be a countable abelian group and $B$ be a real separable Banach space. Equip the set $\Hom(A,B)$ with the point-wise linear operations and with the topology of point-wise convergence. Then $\Hom(A,B)$ is a real separable Fr\'echet space. If, in addition, the group $A$ is of a finite rank then $\Hom(A,B)$ is a real (separable) Banach space.
\end{lemma}
\begin{proof}
Since a linear combination of morphisms $A\to B$ is again a morphism $A\to B$, $\Hom(A,B)$ is a real linear space. Fix a maximal independent set $S$ in $A$. For every $a\in S$, consider the map $\|.\|_a\colon\Hom(A,B)\ni h\mapsto\|h(a)\|\in\mathbb R$, where $\|.\|$ is the norm of $B$. By an elementary argument, all the maps $\|.\|_a$ are semi-norms on $\Hom(A,B)$ and, for every $h\in\Hom(A,B)$, $h=0$ if and only if $\|h\|_a=0$ for every $a\in S$. Moreover, by our choice of $S$, for every sequence $(h_n)_{n\in\mathbb N}$ in $\Hom(A,B)$ and every $h\in\Hom(A,B)$, the following conditions are equivalent:
\begin{itemize}
\item $h_n\to h$ point-wise on $A$,
\item $h_n(a)\to h(a)$ in $B$ for every $a\in S$.
\end{itemize}
Consequently, the semi-norms $\|.\|_a$ induce on $\Hom(A,B)$ the topology of point-wise convergence. Finally, the completeness of $\Hom(A,B)$ follows from the completeness of $B$, hence $\Hom(A,B)$ is a real Fr\'echet space. If, in addition, the group $A$ is of a finite rank then the set $S$ is finite and so the semi-norms $\|.\|_a$ provide us with the structure of a real Banach space on $\Hom(A,B)$.

To finish the proof of the lemma, we must show that $\Hom(A,B)$ is a separable space. Since the group $A$ is countable and the Banach space $B$ is separable, it follows that the product space $B^A$ is separable and metrizable. Since $\Hom(A,B)$ is a topological subspace of $B^A$ by definition, the separability of $\Hom(A,B)$ now follows.
\end{proof}

\subsection{Hypersemigroups}\label{Sub:hypsmgps}

Let $G\in\mathsf{CAGp}$ and let $2^G$\index[symbol]{$2^G$} be the set of all non-empty closed subsets of $G$. With the operation of the product of sets in $G$, $2^G$ is an abelian semigroup with an identity $e$. With the Vietoris topology, $2^G$ is a compact space and, by virtue of Lemma~\ref{L:Viet.top.hyp} and Remark~\ref{R:Viet.top.hyp} below, the semigroup operation on $2^G$ is continuous. Thus, $2^G$ is a compact semigroup; we shall call it the \emph{hypersemigroup over}\index{hypersemigroup} $G$. The subset $\mathcal S(G)$\index[symbol]{$\mathcal S(G)$} of $2^G$, formed by the closed subgroups of $G$, is a subsemigroup of $2^G$. By compactness of $G$, $\mathcal S(G)$ coincides with the set of idempotents of $2^G$; that is $H\in 2^G$ is an element of $\mathcal S(G)$ if and only if $HH=H$ (this follows from the Weil's lemma). Consequently, $\mathcal S(G)$ is a closed subset of $2^G$ and it is therefore a compact semigroup on its own.

We close this subsection by constructing a basis for the Vietoris topology on $2^G$, which involves the group structure of $G$; it will simplify our later considerations. Let us begin by fixing some notation. Given a nonempty open set $U\subseteq G$, we write
\begin{equation*}
\mathcal N_u(U)\index[symbol]{$\mathcal N_u(U)$}=\{K\in 2^G : K\subseteq U\}\hspace{3mm}\text{and}\hspace{3mm}\mathcal N_d(U)\index[symbol]{$\mathcal N_d(U)$}=\{K\in 2^G : K\cap U\neq\emptyset\}.
\end{equation*}
Further, given $K\in 2^G$ and an identity neighbourhood $V$ in $G$, we write
\begin{equation*}
\mathcal U_V(K)\index[symbol]{$\mathcal U_V(K)$}=\{H\in 2^G : H\subseteq VK\hspace{2mm}\text{and}\hspace{2mm}K\subseteq VH\}.
\end{equation*}
By definition of the Vietoris topology, the sets $\mathcal N_u(U)$, $\mathcal N_d(U)$ form a subbasis for the topology of $2^G$.

\begin{lemma}\label{L:Viet.top.hyp}
Let $G\in\mathsf{CAGp}$. Given $K\in 2^G$, the sets $\mathcal U_V(K)$ form a local base at $K$ in $2^G$, when $V$ runs through the identity neighbourhoods in $G$.
\end{lemma}
\begin{remark}\label{R:Viet.top.hyp}
It follows at once from the lemma that the semigroup operation on $2^G$ is continuous. Indeed, given $H,K\in 2^G$ and an identity neighbourhood $V$ in $G$, choose an identity neighbourhood $W$ in $G$ with $WW\subseteq V$. Then $EF\in\mathcal U_V(HK)$ for all $E\in\mathcal U_W(H)$ and $F\in\mathcal U_W(K)$.
\end{remark}
\begin{proof}[Proof of Lemma~\ref{L:Viet.top.hyp}]
First notice that $K\in\mathcal U_V(K)$ for all $K$ and $V$ by definition. We verify the following statement:
\begin{itemize}
\item[($*$)] given $K\in 2^G$ and a neighbourhood $\mathcal N$ of $K$ in $2^G$, there is an identity neighbourhood $V$ in $G$ with $\mathcal U_V(K)\subseteq\mathcal N$.
\end{itemize}
Since the inclusion $\mathcal U_V(K)\cap\mathcal U_W(K)\supseteq\mathcal U_{V\cap W}(K)$ holds for all identity neighbourhoods $V,W$ in $G$, it suffices to verify ($*$) in the case when $\mathcal N$ is of the form $\mathcal N=\mathcal N_u(U)$ or $\mathcal N=\mathcal N_d(U)$ for some non-empty open set $U\subseteq G$. Fix such a set $U\subseteq G$. If $K\in\mathcal N_u(U)$ then $VK\subseteq U$ for an appropriate identity neighbourhood $V$ in $G$; in this case, $\mathcal U_V(K)\subseteq\mathcal N_u(U)$. If $K\in\mathcal N_d(U)$ then $gV^{-1}\subseteq U$ for an appropriate pair of an identity neighbourhood $V$ in $G$ and an element $g\in K$; in this case, $\mathcal U_V(K)\subseteq\mathcal N_d(U)$.

We shall continue by verifying the following statement:
\begin{itemize}
\item[($\star$)] given $H,K\in 2^G$ and an identity neighbourhood $V$ in $G$ with $H\in\mathcal U_V(K)$, there is an identity neighbourhood $W$ in $G$ with $\mathcal U_W(H)\subseteq\mathcal U_V(K)$.
\end{itemize}
By compactness of $H$ and $K$, there are identity neighbourhoods $W,N$ in $G$ with $H\subseteq NK$, $K\subseteq NH$ and $WN\subseteq V$. We claim that $W$ satisfies the conclusion of ($\star$). So let $C\in\mathcal U_W(H)$. Then the inclusions $C\subseteq WH$ and $H\subseteq NK$ yield $C\subseteq WNK\subseteq VK$ and, similarly, the inclusions $K\subseteq NH$ and $H\subseteq WC$ yield $K\subseteq NWC\subseteq VC$. Thus, $C\in\mathcal U_V(K)$, as was to be shown.

We finish the proof of the lemma by showing that all the sets $\mathcal U_V(K)$ are open in $2^G$. In view of ($\star$), it suffices to verify the following statement:
\begin{itemize}
\item[($\circ$)] given $K\in 2^G$ and an identity neighbourhood $V$ in $G$, there is a neighbourhood $\mathcal N$ of $K$ in $2^G$ with $\mathcal N\subseteq\mathcal U_V(K)$.
\end{itemize}
Fix an identity neighbourhood $W$ in $G$ with $WW\subseteq V$ and points $g_1,\dots,g_n\in K$ with $K\subseteq\bigcup_{i=1}^ng_iW$. Set $U=VK$ and $U_i=g_iW^{-1}$ for $i=1,\dots,n$. An elementary argument shows that
\begin{equation*}
K\in\mathcal N_u(U)\cap\bigcap_{i=1}^n\mathcal N_d(U_i)\subseteq\mathcal U_V(K).
\end{equation*}
Since the set $\mathcal N=\mathcal N_u(U)\cap\bigcap_{i=1}^n\mathcal N_d(U_i)$ is open in $2^G$, this verifies statement ($\circ$).
\end{proof}

\subsection{Group-disjointness}

Let $H_1,H_2\in\mathsf{CAGp}$. For the purpose of this work we call $H_1,H_2$ \emph{group-disjoint}\index{group-disjoint groups} and write $H_1\perp H_2$\index[symbol]{$H_1\perp H_2$}, if $H_1\times H_2$ is the only closed subgroup of $H_1\times H_2$ with full projections onto both $H_1$ and $H_2$. This relation is symmetric and saturated with respect to topological isomorphisms. (Our terminology is influenced by the concept of disjointness in topological dynamics: two minimal flows with a common acting group are called disjoint if their product contains no proper subsystem with full both projections. In terms of universal algebra, the groups $H_1,H_2$ are group-disjoint if and only if the set of all their topological subdirect products consists of their direct product alone.) We shall use the following characterization of group-disjointness.

\begin{lemma}\label{L:grp.disj.char}
Let $H_1,H_2\in\mathsf{CAGp}$. Then the following conditions are equivalent:
\begin{enumerate}
\item $H_1,H_2$ are group-disjoint,
\item $H_1,H_2$ have no non-trivial common quotient group,
\item $H_1^*,H_2^*$ have no non-trivial common subgroup up to isomorphism,
\item $\ord(\chi_1)\neq\ord(\chi_2)$ for all $1\neq\chi_1\in H_1^*$ and $1\neq\chi_2\in H_2^*$.
\end{enumerate}
\end{lemma}
\begin{remark}\label{R:grp.disj.char}
We wish to mention the following facts.
\begin{itemize}
\item From part (4) of the lemma it follows that every family $\mathcal A$, consisting of pair-wise group-disjoint non-trivial compact abelian groups, is at most countable. To see this, fix $\chi_H\in H^*\setminus1$ for every $H\in\mathcal A$ and consider the function $\mathcal A\ni H\mapsto\ord\left(\chi_H\right)\in\mathbb N\cup\{\infty\}$. By part (4) of the lemma, this function is injective. Thus, indeed, the family $\mathcal A$ is at most countable.
\item By part (2) of the lemma, the following statement holds: if $H_i,G_i\in\mathsf{CAGp}$ ($i=1,2$), $H_1\perp H_2$ and there are epimorphisms $q_i\in\Hom(H_i,G_i)$ ($i=1,2$), then $G_1\perp G_2$.
\item Recall that if a group $H\in\mathsf{CAGp}$ is not totally disconnected then there is $\chi\in H^*$ with $\ord(\chi)=\infty$. Thus, by part (4) of the lemma, if $H_i$ ($i\in I$) is a family of pair-wise group-disjoint groups from $\mathsf{CAGp}$ then there is at most one $i\in I$ such that the group $H_i$ is not totally disconnected.
\end{itemize}
\end{remark}
\begin{proof}[Proof of Lemma~\ref{L:grp.disj.char}]
To verify implication (1)$\Rightarrow$(2), fix a group $K\in\mathsf{CAGp}$ and epimorphisms $p_1\colon H_1\to K$, $p_2\colon H_2\to K$. Consider the set
\begin{equation*}
H=\left\{(h_1,h_2)\in H_1\times H_2 : p_1(h_1)=p_2(h_2)\right\}.
\end{equation*}
Clearly, $H$ is a closed subgroup of $H_1\times H_2$. Since both $p_1$ and $p_2$ are epimorphisms, $H$ has full projections onto both $H_1$ and $H_2$. It follows by group-disjointness of $H_1,H_2$ that $H=H_1\times H_2$ and so the group $K$ is trivial.

To verify implication (2)$\Rightarrow$(3), assume that $H_1,H_2$ have no common quotient group other than the trivial one and that $D$ is a common subgroup of $H_1^*,H_2^*$ up to isomorphism. Fix the corresponding monomorphisms $l_1\colon D\to H_1^*$ and $l_2\colon D\to H_2^*$. Their dual morphisms $l_1^*\colon H_1\to D^*$, $l_2^*\colon H_2\to D^*$ are epimorphisms and hence, by virtue of (2), $D^*$ is trivial. It follows that $D$ is also trivial, which verifies (3).

To show that (3)$\Rightarrow$(4), assume that $H_1^*$, $H_2^*$ have no non-trivial common subgroup up to isomorphism and let $\chi_1\in H_1^*$, $\chi_2\in H_2^*$ be non-trivial characters of $H_1,H_2$, respectively. By our assumptions, the cyclic subgroups $\langle\chi_1\rangle$ and $\langle\chi_2\rangle$ of $H_1^*$ and $H_2^*$ generated by $\chi_1$ and $\chi_2$, respectively, are not isomorphic. This means that $\ord(\chi_1)\neq\ord(\chi_2)$, which verifies (4).

We finish the proof by verifying implication (4)$\Rightarrow$(1). To this end, assume that $H_1$, $H_2$ are not group-disjoint and fix a proper closed subgroup $H$ of $H_1\times H_2$ with full projections onto $H_1$ and $H_2$. Choose $\Upsilon\in[H^{\perp},(H_1\times H_2)^*]\setminus1$; clearly, $\Upsilon$ is of the form $\Upsilon(h_1,h_2)=\chi_1(h_1)\chi_2(h_2)^{-1}$ for appropriate characters $\chi_1\in H_1^*$ and $\chi_2\in H_2^*$. Since $H$ has full projections onto both $H_1,H_2$ and $\Upsilon$ vanishes on $H$, we have $\im(\chi_1)=\im(\chi_2)$. It follows, finally, that both $\chi_1$ and $\chi_2$ are non-trivial and $\ord(\chi_1)=\ord(\chi_2)$.
\end{proof}

\begin{lemma}\label{L:prod.of.disj}
Let $G_1,\dots,G_n,H_1,\dots,H_m\in\mathsf{CAGp}$ be pair-wise group-disjoint. Then $G_1\times\dots\times G_n$ and $H_1\times\dots\times H_m$ are also group-disjoint.
\end{lemma}
\begin{proof}
By Lemma~\ref{L:grp.disj.char}(3), we need to show that $G_1^*\times\dots\times G_n^*$ and $H_1^*\times\dots\times H_m^*$ have no non-trivial common subgroup up to isomorphism. We shall proceed by contradiction; let $D$ be a non-trivial group and let $\mu\colon D\to G_1^*\times\dots\times G_n^*$ and $\nu\colon D\to H_1^*\times\dots\times H_m^*$ be monomorphisms. By the third part of Remark~\ref{R:grp.disj.char}, at most one of the groups $G_1,\dots,G_n,H_1,\dots,H_m$ is not totally disconnected. Hence, $D$ is a torsion group. Fix $d\in D$ and $p\in\mathbb N$ prime with $\ord(d)=p$. Write $\mu(d)=(\chi_1,\dots,\chi_n)$ and $\nu(d)=(\Upsilon_1,\dots,\Upsilon_m)$, and choose $i$ and $j$ so that $\chi_i\neq1$ and $\Upsilon_j\neq1$. Since $\chi_i\neq1$, $\chi_i^p=1$ and $p$ is prime, we have $\ord(\chi_i)=p$. A similar argument yields $\ord(\Upsilon_j)=p$. By Lemma~\ref{L:grp.disj.char}(4), it follows that $G_i$ and $H_j$ are not group-disjoint, a contradiction.
\end{proof}

\section{Minimal flows and group extensions}\label{S:top.dyn}

Our references in this section are monographs \cite{Aus}, \cite{Bro}, \cite{Ell}, \cite{Gla} and \cite{deV}.

\subsection{Minimal flows}\label{Sub:min.flws}

By a \emph{flow}\index{flow} $\mathcal F$ we mean a representation of a topological group $\Gamma$ as a group of homeomorphisms on a topological space $X$. We write $\Flow$\index[symbol]{$\Flow$} and always require that the corresponding action $\Gamma\times X\to X$ be jointly continuous. The acting homeomorphisms\index{acting homeomorphism of a flow} (also called the transition maps)\index{transition map of a flow} of $\mathcal F$ are denoted by $T_{\gamma}$\index[symbol]{$T_{\gamma}$} ($\gamma\in\Gamma$). It is often assumed that $X$ is a pointed space with the base point $z$; in such a situation we may consider the motion map\index{motion map of a flow}
\begin{equation*}
\mathcal F_z=\mathcal F(-,z)\colon\Gamma\ni\gamma\mapsto\mathcal F(\gamma,z)=T_{\gamma}z\in X.\index[symbol]{$\mathcal F_z$}
\end{equation*}
Since $\mathcal F_z(1)=z$, $\mathcal F_z$ is a (continuous) base point preserving map.

Given $x\in X$, we write $\mathcal O_{\mathcal F}(x)=\{T_{\gamma}x : \gamma\in\Gamma\}$\index[symbol]{$\mathcal O_{\mathcal F}(x)$} for the $\mathcal F$-orbit of $x$ and $\overline{\mathcal O}_{\mathcal F}(x)$\index[symbol]{$\overline{\mathcal O}_{\mathcal F}(x)$} for the closure of $\mathcal O_{\mathcal F}(x)$ in $X$. The flow $\mathcal F$ is called \emph{minimal}\index{minimal flow} if all its orbits are dense. A point $x\in X$ is \emph{free}\index{free!point for a flow} for $\mathcal F$ if it has a trivial isotropy group, that is, if $T_{\gamma}x=x$ implies $\gamma=1$. The flow $\mathcal F$ is called \emph{point-free}\index{point-free flow}, respectively, \emph{free}\index{free!flow} if the set of its free points is nonempty, respectively, the whole $X$.

\subsection{Topological freeness and free cycles}\label{Sub:top.frnss.fr.ccl}

Let $\Flow$ be a flow with $\Gamma\in\mathsf{LieGp}$ connected and with $X$ a compact connected manifold. Fix $z\in X$ and consider the morphism $\mathcal F_z^{\sharp}\colon H_1^w(\Gamma)\to H_1^w(X)$ induced by the map $\mathcal F_z$. The image $\im(\mathcal F_z^{\sharp})$ of $\mathcal F_z^{\sharp}$ is denoted by $H_1^w(\mathcal F)$\index[symbol]{$H_1^w(\mathcal F)$}. We say that the flow $\mathcal F$
\begin{itemize}
\item is \emph{topologically free}\index{topologically free flow} if $\mathcal F_z^{\sharp}\colon H_1^w(\Gamma)\to H_1^w(X)$ is a monomorphism (equi\-va\-lent\-ly, if the restriction $\mathcal F_z^{\sharp}\colon H_1^w(\Gamma)\to H_1^w(\mathcal F)$ is an isomorphism),
\item \emph{possesses a free cycle}\index{free!cycle for a flow} if $\rank(H_1^w(\mathcal F))<\rank(H_1^w(X))$.
\end{itemize}
We write $n=\rank(H_1^w(\mathcal F))$ and $n+m=\rank(H_1^w(X))$. Also, $d_1,\dots,d_n$ denote the ele\-men\-tary divisors of $H_1^w(\mathcal F)$ in $H_1^w(X)$ (written, as usual, in the non-decreasing order). The following lemma describes what happens if the base point $z$ of the phase space $X$ of $\mathcal F$ is changed.

\begin{lemma}\label{L:top.free.bspt}
Let $\Flow$ be a flow with $\Gamma\in\mathsf{LieGp}$ connected and with $X$ a compact connected manifold. Fix $z,z'\in X$ and let $f$ be a path in $X$ from $z$ to $z'$. Then the diagram in Figure~\ref{Fig:top.free.bspt} below commutes. In particular, we have
\begin{equation*}
\widehat{f}\mathcal F_z^{\sharp}=\mathcal F_{z'}^{\sharp}\colon H_1^w(\Gamma)\to H_1^w(X,z').
\end{equation*}
\end{lemma}
\begin{figure}[ht]
\[
\ctdiagram{
\ctinnermid
\ctv 0,0: {\pi_1(\Gamma)}
\ctv 100,0: {\pi_1(X,z)}
\ctv 50,-30: {\pi_1(X,z')}
\ctv 0,-70: {H_1^w(\Gamma)}
\ctv 100,-70: {H_1^w(X,z)}
\ctv 50,-100: {H_1^w(X,z')}
\put(14,0){\vector(1,0){66}}
\put(16,-70){\line(1,0){32}}
\put(52,-70){\vector(1,0){26}}
\put(0,-7){\vector(0,-1){57}}
\put(100,-7){\vector(0,-1){57}}
\put(50,-37){\vector(0,-1){57}}
\put(9,-7){\vector(1,-0.60){30}}
\put(9,-77){\vector(1,-0.60){28}}
\put(93,-6.2){\vector(-1,-0.6){28}}
\put(93,-76.2){\vector(-1,-0.6){28}}
\ctv -7,-35: {p_{\Gamma}}
\ctv 113,-35: {p_{X,z}}
\ctv 63,-50: {p_{X,z'}}
\ctv 30,7: {\mathcal F_z^{\sharp}}
\ctv 30,-63: {\mathcal F_z^{\sharp}}
\ctv 17,-22: {\mathcal F_{z'}^{\sharp}}
\ctv 17,-92: {\mathcal F_{z'}^{\sharp}}
\ctv 83,-22: {\widehat{f}}
\ctv 83,-92: {\widehat{f}}
\ctnohead
}
\]
\caption{Changing the base point $z$ of the phase space $X$}
\label{Fig:top.free.bspt}
\end{figure}
\begin{remark}
Recall from Subsection~\ref{Sub:frst.wk.hm.gp} that $\widehat{f}\colon H_1^w(X,z)\to H_1^w(X,z')$ is an isomorphism of groups. It follows, in particular, that the concepts and quantities introduced in the paragraph preceding Lemma~\ref{L:top.free.bspt} do not depend on the choice of a base point of $X$.
\end{remark}
\begin{proof}[Proof of Lemma~\ref{L:top.free.bspt}]
The commutativity of the three vertical squares in Figure~\ref{Fig:top.free.bspt} follows from our discussion in Subsection~\ref{Sub:frst.wk.hm.gp}. Thus, since $p_{\Gamma}$ is an epimorphism, it suffices now to verify the commutativity of the upper horizontal triangle. To this end, fix $g\in\pi_1(\Gamma)$ and consider the constant loops $c_1$, $c_z$, $c_{z'}$ at $1\in\Gamma$, $z,z'\in X$, respectively. Then
\begin{equation*}
\begin{split}
\widehat{f}\mathcal F_z^{\sharp}(g)&=\overline{f}*\mathcal F_zg*f=\mathcal F(c_1,\overline{f})*\mathcal F(g,c_z)*\mathcal F(c_1,f)\\
&=\mathcal F^{\sharp}\left((c_1,\overline{f})*(g,c_z)*(c_1,f)\right)=\mathcal F^{\sharp}\left((c_1*g*c_1,\overline{f}*c_z*f)\right)\\
&=\mathcal F^{\sharp}\left((g,c_{z'})\right)=\mathcal F_{z'}^{\sharp}(g),
\end{split}
\end{equation*}
as was to be shown.
\end{proof}

\subsection{Freeness versus topological freeness}\label{Sub:frns.vs.top.frns}

We shall often assume in this work that the considered minimal flow $\mathcal F$ is (point-)free or topologically free. The freeness will be required when applying results from \cite{Dir3}, while topological freeness will be assumed when using techniques developed in this work. There is a formal analogy between the definitions of these two properties:
\begin{itemize}
\item a flow $\Flow$ is free if $\mathcal F_z\colon\Gamma\to X$ is an injective map for every $z\in X$,
\item a flow $\Flow$ is topologically free if $\mathcal F_z^{\sharp}\colon H_1^w(\Gamma)\to H_1^w(X)$ is a mo\-no\-mor\-phism.
\end{itemize}
Despite this formal analogy, these two properties are independent for minimal flows, as we now show.

To see that topological freeness does not imply freeness, let $\Gamma\in\mathsf{LieGp}$ be a simply connected group and let $\Gamma'$ be a non-trivial proper closed co-compact subgroup of $\Gamma$. Consider the natural action $\mathcal F\colon\Gamma\curvearrowright X$ of $\Gamma$ on $X=\Gamma/\Gamma'$. Then $X$ is a compact connected manifold and $\mathcal F$ is a minimal $\Gamma$-flow on $X$. Since the group $\Gamma'$ is non-trivial, the flow $\mathcal F$ is not free; actually, $\mathcal F$ possesses no free points whatsoever and so it is not even point-free. On the other hand, since the group $\Gamma$ is simply connected, we have $H_1^w(\Gamma)=0$ and so the induced morphism $\mathcal F_z^{\sharp}\colon H_1^w(\Gamma)\to H_1^w(X)$ is a monomorphism. Thus, the flow $\mathcal F$ is topologically free.

Now we want to show that freeness does not imply topological freeness; this is done in Proposition~\ref{P:fr.not.imp.topfr} below. Let us start with some auxiliary observations. We shall consider the usual transitive action of the second unitary group $\text{U}(2)$ on the $3$-sphere $\mathbb S^3$. Recall that the center $Z(\text{U}(2))$ of $\text{U}(2)$ consists of the matrices $A_{\zeta}=\zeta\mathbb I$, where $\mathbb I$ is the unit of $\text{U}(2)$ and $\zeta\in\mathbb T^1$. Thus, $AA_{\zeta}=A_{\zeta}A$ for all $A\in\text{U}(2)$ and $\zeta\in\mathbb T^1$. Notice also that the action of $\text{U}(2)$ on $\mathbb S^3$ is effective and that the corresponding restricted action of $Z(\text{U}(2))$ on $\mathbb S^3$ is free. It follows, in particular, that we may regard $\text{U}(2)$ as a subgroup of the group of all homeomorphisms on $\mathbb S^3$.

\begin{proposition}\label{P:fr.not.imp.topfr}
Let $X$ be a compact connected manifold with dimension at least $2$. If $X$ admits a minimal continuous flow $\phi\colon\mathbb R\curvearrowright X$ then there is a free minimal flow $\mathcal F\colon\mathbb R\times\mathbb T^1\curvearrowright X\times\mathbb S^3$ with $H_1^w(\mathcal F)=0$.
\end{proposition}
\begin{remark}\label{R:fr.not.imp.topfr}
The flow $\mathcal F$ from the proposition is not topologically free. Indeed, since $H_1^w(\mathbb R\times\mathbb T^1)\cong\mathbb Z$ and $H_1^w(\mathcal F)=0$, the induced morphism
\begin{equation*}
\mathcal F_z^{\sharp}\colon H_1^w(\mathbb R\times\mathbb T^1)\to H_1^w(\mathcal F)\subseteq H_1^w(X\times\mathbb S^3)
\end{equation*}
is not a monomorphism.
\end{remark}
\begin{proof}[Proof of Proposition~\ref{P:fr.not.imp.topfr}]
We shall proceed in four steps.

\emph{1st step.} We construct an auxiliary minimal flow $\mathcal F'$.

Let $\varphi_t$ ($t\in\mathbb R$) be the acting homeomorphisms of $\phi$. Since the flow $\phi$ is minimal and the dimension of $X$ is at least $2$, it follows that $\phi$ is a free flow. Further, since $\text{U}(2)$ is a compact connected Lie group, it is Polish and arc-wise connected. Finally, since the action of $\text{U}(2)$ on $\mathbb S^3$ is transitive and hence minimal, we may use \cite[Theorem~8]{Dir3} to find a minimal skew product $\mathcal F'\colon\mathbb R\curvearrowright X\times\mathbb S^3$ over $\phi$, whose fibre maps lie in the group $\text{U}(2)$. That is, for every $t\in\mathbb R$, the $t$-transition $T_t'$ of $\mathcal F'$ has the form $T_t'(x,y)=(\varphi_tx,A_{(t,x)}y)$, where $A_{(t,x)}\in\text{U}(2)$ for every $x\in X$.

\emph{2nd step.} We construct the desired minimal flow $\mathcal F$.

Given $t\in\mathbb R$ and $\zeta\in\mathbb T^1$, consider the map
\begin{equation*}
T_{(t,\zeta)}\colon X\times\mathbb S^3\ni(x,y)\mapsto\left(\varphi_tx,A_{\zeta}A_{(t,x)}y\right)\in X\times\mathbb S^3.
\end{equation*}
Since $\mathcal F'\colon\mathbb R\curvearrowright X\times\mathbb S^3$ is a flow and $A_{\zeta}\in Z(\text{U}(2))$ for every $\zeta\in\mathbb T^1$, it follows that the maps $T_{(t,\zeta)}$ constitute an action $\mathcal F$ of $\mathbb R\times\mathbb T^1$ on $X\times\mathbb S^3$. Clearly, $\mathcal F$ is continuous as a map
\begin{equation*}
\mathcal F\colon(\mathbb R\times\mathbb T^1)\times(X\times\mathbb S^3)\to X\times\mathbb S^3
\end{equation*}
and so it is a flow. The minimality of $\mathcal F$ follows immediately from the minimality of $\mathcal F'$.

\emph{3rd step.} We show that the flow $\mathcal F$ is free.

Fix $t\in\mathbb R$, $\zeta\in\mathbb T^1$, $x\in X$, $y\in\mathbb S^3$ and assume that $T_{(t,\zeta)}(x,y)=(x,y)$; we show that $t=0$ and $\zeta=1$. Since $x=\varphi_tx$ and the flow $\phi$ is free, we get $t=0$ and hence $y=A_{\zeta}A_{(t,x)}y=A_{\zeta}y$. Since $Z(\text{U}(2))$ acts freely on $\mathbb S^3$, it follows that $\zeta=1$, as was to be shown.

\emph{4th step.} We show that $H_1^w(\mathcal F)=0$.

Choose a base point $z=(x_0,y_0)$ for $X\times\mathbb S^3$. Consider the topological morphism $j\colon\mathbb T^1\ni\zeta\mapsto(0,\zeta)\in\mathbb R\times\mathbb T^1$. Since $j$ is a homotopy equivalence, the induced morphism $j^{\sharp}\colon H_1^w(\mathbb T^1)\to H_1^w(\mathbb R\times\mathbb T^1)$ is an isomorphism of groups. Consequently,
\begin{equation}\label{Eq:HF.red.comp}
H_1^w(\mathcal F)=\mathcal F_z^{\sharp}\left(H_1^w(\mathbb R\times\mathbb T^1)\right)=\mathcal F_z^{\sharp}j^{\sharp}\left(H_1^w(\mathbb T^1)\right)=(\mathcal F_zj)^{\sharp}\left(H_1^w(\mathbb T^1)\right).
\end{equation}
Further, by definition of $\mathcal F$, the map $\mathcal F_zj$ takes its values in $x_0\times\mathbb S^3$. Since the latter is a simply connected space, we have $(\mathcal F_zj)^{\sharp}(H_1^w(\mathbb T^1))=0$. In view of (\ref{Eq:HF.red.comp}), this means that $H_1^w(\mathcal F)=0$, as was to be shown.
\end{proof}

\subsection{Remarks on topological freeness and free cycles}

Our definition of topological freeness and of the existence of a free cycle for a minimal flow $\Flow$ involves the induced morphism $\mathcal F_z^{\sharp}$ between the first weak homology groups of $\Gamma$ and $X$. In this subsection we link these two properties of $\mathcal F$ to the first and the second homotopy group of the quotient space $X/K$ of $X$ by the maximal compact subgroup $K$ of $\Gamma$. This is intended as an informal discussion and we shall not use the results of this subsection in the rest of this work.

Let $\Flow$ be a minimal flow with $\Gamma\in\mathsf{LieGp}$ connected and with $X$ a compact connected manifold. We shall assume that $\mathcal F$ is a free flow. Fix a maximal compact subgroup $K$ of $\Gamma$. Recall that $K$ is itself a connected Lie group; consequently, the groups $\pi_2(K)$ and $\pi_0(K)=K/K_0$ both vanish, while the group $\pi_1(K)$ is abelian. Fix a base point $z$ for $X$ and consider the restriction $j=\mathcal F_z|_K\colon K\to X$ of $\mathcal F_z$. Further, let $p\colon X\to X/K$ be the canonical projection. We shall work under the assumption that $X$ is a locally trivial fibre bundle over $X/K$ with the fibre $K$; since $K$ acts freely on $X$ by our assumptions, this is the case if both $X$ and the action of $K$ on $X$ are smooth. Consequently, by \cite[Section~9.8]{FomFuch}, we have the corresponding long exact homotopy sequence of the bundle
\begin{equation}\label{Eq:long.ex.hom.seq.tffc}
0\longrightarrow\pi_2(X)\stackrel{p^{\sharp}_2}{\longrightarrow}\pi_2(X/K)\stackrel{\partial}{\longrightarrow}\pi_1(K)\stackrel{j^{\sharp}_1}{\longrightarrow}\pi_1(X)\stackrel{p^{\sharp}_1}{\longrightarrow}\pi_1(X/K)\longrightarrow0,
\end{equation}
where $p^{\sharp}_2$, $p^{\sharp}_1$ and  $j^{\sharp}_1$ stand for the morphisms induced by $p$ and $j$, respectively, and $\partial$ is the connecting boundary morphism. To simplify the situation, we shall further assume that the space $X$ has abelian fundamental group $\pi_1(X)$. (Though certainly restrictive, this assumption still leaves us a large class of manifolds $X$ at our disposal.) By the exactness of (\ref{Eq:long.ex.hom.seq.tffc}) at $\pi_1(X/K)$, this yields also the commutativity of $\pi_1(X/K)$.

Since $\pi_1(K)$, $\pi_1(X)$ and $\pi_1(X/K)$ are finitely generated abelian groups, each of them contains its torsion subgroup as a direct summand and the corresponding complementary summand is the first weak homology group of $K$, $X$ and $X/K$, respectively. Consequently, the diagram in Figure~\ref{Fig:suf.cond.tffc} below commutes and its upper row is an exact sequence by the exactness of (\ref{Eq:long.ex.hom.seq.tffc}). 
\begin{figure}[ht]
\[\minCDarrowwidth20pt\begin{CD}
H_1^w(K)\oplus\tau_K @>j^{\sharp}_1>> H_1^w(X)\oplus\tau_X @>p^{\sharp}_1>> H_1^w(X/K)\oplus\tau_{X/K} @>>> 0\\
@V\pr_1VV @V\pr_1VV @V\pr_1VV \\
H_1^w(K) @>>j^{\sharp}> H_1^w(X) @>>p^{\sharp}> H_1^w(X/K)
\end{CD}\]
\caption{Remarks on topological freeness and free cycles}
\label{Fig:suf.cond.tffc}
\end{figure}
(For typographical reasons, we write $\tau_K$, $\tau_X$ and $\tau_{X/K}$ for the torsion subgroups of $\pi_1(K)$, $\pi_1(X)$ and $\pi_1(X/K)$, respectively.) We denote by $c$ and $d$ the largest torsion coefficients of the groups $\tau_X$ and $\tau_{X/K}$, respectively. Notice also that $p^{\sharp}$ is an epimorphism, since $p_1^{\sharp}$ is an epimorphism.

Now, the flow $\mathcal F$ is topologically free, by definition, if and only if the induced morphism $\mathcal F_z^{\sharp}\colon H_1^w(\Gamma)\to H_1^w(X)$ is a monomorphism. Since the inclusion morphism $K\to\Gamma$ induces an isomorphism $H_1^w(K)\to H_1^w(\Gamma)$, this occurs if and only if the induced morphism $j^{\sharp}\colon H_1^w(K)\to H_1^w(X)$ is a monomorphism. By commutativity of the diagram in Figure~\ref{Fig:suf.cond.tffc} and by the exactness of (\ref{Eq:long.ex.hom.seq.tffc}) at $\pi_1(K)$, we have $c\ker(j^{\sharp})\subseteq\ker(j_1^{\sharp})=\im(\partial)$. Thus, since $H_1^w(K)\supseteq\ker(j^{\sharp})$ is a torsion-free group,
\begin{enumerate}
\item[($*$)] \emph{the flow $\mathcal F$ is topologically free, provided $\pi_2(X/K)$ is a torsion group}.
\end{enumerate}
Since the morphism $p_2^{\sharp}$ from (\ref{Eq:long.ex.hom.seq.tffc}) is a monomorphism, in order that $\pi_2(X/K)$ be a torsion group it is necessary that $\pi_2(X)$ be a torsion group.

Further, recall that the flow $\mathcal F$ possesses a free cycle, by definition, if and only if $\rank(H_1^w(\mathcal F))<\rank(H_1^w(X))$, where $H_1^w(\mathcal F)=\mathcal F_z^{\sharp}H_1^w(\Gamma)=\im(j^{\sharp})$. An elementary argument involving properties of the diagram in Figure~\ref{Fig:suf.cond.tffc} yields $d\ker(p^{\sharp})\subseteq\im(j^{\sharp})$ and $c\im(j^{\sharp})\subseteq\ker(p^{\sharp})$. Consequently,
\begin{equation*}
\rank(\ker(p^{\sharp}))=\rank(\im(j^{\sharp}))=\rank(H_1^w(\mathcal F))
\end{equation*}
and so the following conditions are equivalent:
\begin{itemize}
\item $\rank(H_1^w(\mathcal F))<\rank(H_1^w(X))$,
\item $\rank(\ker(p^{\sharp}))<\rank(H_1^w(X))$,
\item $\beta(H_1^w(X/K))\neq0$.
\end{itemize}
Thus, since $\beta(H_1^w(X/K))=\beta(\pi_1(X/K))$,
\begin{enumerate}
\item[($*$)] \emph{the flow $\mathcal F$ possesses a free cycle if and only if $\beta(\pi_1(X/K))\neq0$}.
\end{enumerate}

\subsection{Freeness and topological freeness - an implication}

We have seen in Subsection~\ref{Sub:frns.vs.top.frns} that freeness and topological freeness of a minimal flow are in general independent properties. Our aim in this subsection is to show that, after all, there does exist a non-trivial relationship between these two properties for a special class of flows; this is done in Proposition~\ref{P:min.eq.flow.TFFC} below. Before formulating the result, let us describe the considered situation.

Let $\Gamma,G\in\mathsf{LieGp}$ be connected Lie groups, $G$ compact, and let $h\in\Hom(\Gamma,G)$ be a topological morphism. Consider the equicontinuous $\Gamma$-flow $\mathcal F$ on $G$ induced by $h$:
\begin{equation}\label{Eq:min.eq.flow.gen}
\mathcal F\colon\Gamma\times G\ni(\gamma,g)\mapsto h(\gamma)g\in G.
\end{equation}
The following statements hold obviously:
\begin{itemize}
\item $\mathcal F$ is minimal if and only if the image $\im(h)$ of $h$ is dense in $G$,
\item $\mathcal F$ is free if and only if $h$ is a monomorphism.
\end{itemize}
We choose the identity $z=e$ as the base point of $G$; this leads to $\mathcal F_z=h$.

\begin{remark}\label{R:which.Gamma.fthfl.rep}
Given a connected Lie group $\Gamma$, there may not exist a pair of a compact Lie group $G$ and a monomorphism $h\colon\Gamma\to G$ with a dense image. Indeed, since every compact Lie group is topologically isomorphic to a subgroup of the $n^{\text{th}}$ unitary group $\text{U}(n)$ for some $n\in\mathbb N$, a necessary and sufficient condition for the existence of such a pair $G,h$ is that $\Gamma$ possess a faithful finite-dimensional unitary representation; in fact, the group $\Gamma$ must possess a non-trivial one-dimensional unitary representation, that is, a non-trivial character $\Gamma\to\mathbb T^1$. (This last claim can be proved by using results from the structure theory of compact (Lie) groups, but we may also refer to Remark~\ref{R:min.qscob}.) This is not the case, for instance, when $\Gamma=\text{SL}(2,\mathbb R)$ is the second real special linear group, since all its non-trivial unitary representations are infinite-dimensional.
\end{remark}

\begin{proposition}\label{P:min.eq.flow.TFFC}
Let $\Gamma$, $G$, $h$ and $\mathcal F$ be as above. Assume that $h$ is a monomorphism with a dense image, so that the flow $\mathcal F$ is free and minimal. Then $\mathcal F$ is also topologically free. Further, if $\Gamma_t$ and $G_t$ are maximal tori in $\Gamma$ and $G$, respectively, then
\begin{equation}\label{Eq:min.eq.flow.FC}
\rank(H_1^w(G))-\rank(H_1^w(\mathcal F))=\dim(G_t)-\dim(\Gamma_t).
\end{equation}
Finally, if $h$ is not an isomorphism onto $G$ then the flow $\mathcal F$ possesses a free cycle.
\end{proposition}
\begin{remark}
We wish to add the following remarks.
\begin{itemize}
\item Interestingly, the freeness of the flow $\mathcal F$ alone does not secure the topological freeness of $\mathcal F$; without the assumption that $\mathcal F$ is minimal, the proposition is no longer true. To see this, let $\Gamma=\mathbb T^2$, $G=\text{U}(2)$ and consider the monomorphism $h$ defined by
\begin{equation*}
h\colon\mathbb T^2\ni(\zeta,\xi)\mapsto
\left(
\begin{array}{ccc}
\zeta & 0 \\
0 & \xi 
\end{array} \right)
\in\text{U}(2).
\end{equation*}
The corresponding flow $\mathcal F\colon\mathbb T^2\curvearrowright\text{U}(2)$ is free but not minimal. Moreover, since $H_1^w(\mathbb T^2)\cong\mathbb Z^2$ and $H_1^w(\text{U}(2))\cong\mathbb Z$, the induced morphism $\mathcal F_z^{\sharp}=h^{\sharp}\colon H_1^w(\mathbb T^2)\to H_1^w(\text{U}(2))$ is not a monomorphism and so the flow $\mathcal F$ is not topologically free.
\item The dimensions on the right-hand side of (\ref{Eq:min.eq.flow.FC}) are the topological (that is, ma\-ni\-fold) dimensions. Recall that any two maximal tori in a given Lie group are conjugate. Hence, $\dim(\Gamma_t)$ and $\dim(G_t)$ do not depend on the choice of $\Gamma_t$ and $G_t$, respectively.
\end{itemize}
\end{remark}
\begin{proof}[Proof of Proposition~\ref{P:min.eq.flow.TFFC}]
Let $h_1^{\sharp}\colon\pi_1(\Gamma)\to\pi_1(G)$ and $h^{\sharp}\colon H_1^w(\Gamma)\to H_1^w(G)$ be the morphisms induced by $h$. By virtue of \cite[Theorem~1.9]{Pog}, the following statements hold:
\begin{enumerate}
\item[(a)] $h_1^{\sharp}$ is a monomorphism,
\item[(b)] the quotient group $\pi_1(G)/\im(h_1^{\sharp})$ is free abelian with rank
\begin{equation*}
\rank\left(\pi_1(G)/\im(h_1^{\sharp})\right)=\dim(G_t)-\dim(\Gamma_t),
\end{equation*}
\item[(c)] $h$ is an isomorphism onto $G$ if and only if $\dim(G_t)=\dim(\Gamma_t)$.
\end{enumerate}

Let $\tau_{\Gamma}$ and $\tau_G$ be the torsion subgroups of $\pi_1(\Gamma)$ and $\pi_1(G)$, respectively. Since both $\pi_1(\Gamma)$ and $\pi_1(G)$ are finitely generated abelian groups, they contain $\tau_{\Gamma}$ and $\tau_G$ as direct summands with the complementary summands $H_1^w(\Gamma)$ and $H_1^w(G)$, respectively. Moreover, the diagram in Figure~\ref{Fig:min.eq.flow.TFFC} commutes.
\begin{figure}[ht]
\[\minCDarrowwidth20pt\begin{CD}
H_1^w(\Gamma)\oplus\tau_{\Gamma} @>h^{\sharp}_1>> H_1^w(G)\oplus\tau_G \\
@V\pr_1VV @VV\pr_1V \\
H_1^w(\Gamma) @>>h^{\sharp}> H_1^w(G) 
\end{CD}\]
\caption{Freeness and topological freeness - an implication}
\label{Fig:min.eq.flow.TFFC}
\end{figure}

Let $d$ stand for the largest torsion coefficient of $\tau_G$. An elementary argument involving commutativity of the diagram in Figure~\ref{Fig:min.eq.flow.TFFC} shows that $d\ker(h^{\sharp})\subseteq\ker(h_1^{\sharp})$. Since $\ker(h_1^{\sharp})=0$ by virtue of (a) and the group $H_1^w(\Gamma)$ is torsion-free, it follows that $\ker(h^{\sharp})=0$ and so $\mathcal F_z^{\sharp}=h^{\sharp}$ is a monomorphism. Thus, $\mathcal F$ is a topologically free flow.

Now we turn to a proof of (\ref{Eq:min.eq.flow.FC}). First, we claim that 
\begin{equation}\label{Eq:im.h.h1}
\im(h_1^{\sharp})=\im(h^{\sharp})\oplus\tau_G.
\end{equation}
The inclusion ``$\subseteq$'' follows from the commutativity of the diagram in Figure~\ref{Fig:min.eq.flow.TFFC} and from the inclusion $h_1^{\sharp}(\tau_{\Gamma})\subseteq\tau_G$. We verify the inclusion ``$\supseteq$''. First, observe that $\tau_G\subseteq\im(h_1^{\sharp})$, for the quotient group $\pi_1(G)/\im(h_1^{\sharp})$ is torsion-free by virtue of (b). Further, since also $\im(h^{\sharp})\subseteq\im(h_1^{\sharp})+\tau_G$, it follows that $\im(h_1^{\sharp})=\im(h_1^{\sharp})+\tau_G\supseteq\im(h^{\sharp})$. Thus, $\im(h^{\sharp})\oplus\tau_G\subseteq\im(h_1^{\sharp})$, which finishes the proof of (\ref{Eq:im.h.h1}).

From (\ref{Eq:im.h.h1}) we obtain
\begin{equation*}
\begin{split}
H_1^w(G)/H_1^w(\mathcal F)&=H_1^w(G)/\im(h^{\sharp})\cong\left(H_1^w(G)\oplus\tau_G\right)/\left(\im(h^{\sharp})\oplus\tau_G\right)\\
&=\pi_1(G)/\im(h_1^{\sharp}).
\end{split}
\end{equation*}
Consequently, by virtue of (b),
\begin{equation*}
\begin{split}
\rank(H_1^w(G))-\rank(H_1^w(\mathcal F))&=\beta\left(H_1^w(G)/H_1^w(\mathcal F)\right)\\
&=\beta\left(\pi_1(G)/\im(h_1^{\sharp})\right)\\
&=\dim(G_t)-\dim(\Gamma_t),
\end{split}
\end{equation*}
which verifies (\ref{Eq:min.eq.flow.FC}). Finally, the last statement of the proposition follows from (\ref{Eq:min.eq.flow.FC}) and from (c).
\end{proof}

The following corollary generalizes Proposition~\ref{P:min.eq.flow.TFFC} to a larger class of minimal equicontinuous flows.

\begin{corollary}\label{C:min.eq.semis.con}
Let $\Gamma,G\in\mathsf{LieGp}$ be connected, $G$ compact, $h\colon\Gamma\to G$ be a monomorphism with a dense image and $H$ be a closed subgroup of $G$. Consider the minimal equicontinuous flow
\begin{equation*}
\mathcal F'\colon\Gamma\times G/H\ni(\gamma,gH)\mapsto h(\gamma)gH\in G/H.
\end{equation*}
Assume that the group $H$ is connected semi-simple and that $h$ is not an isomorphism onto $G$. Then the flow $\mathcal F'$ is topologically free and possesses a free cycle.
\end{corollary}
\begin{remark}\label{R:min.eq.semis.con}
We wish to add the following remarks.
\begin{itemize}
\item Recall that a compact connected Lie group $H$ is semi-simple, by definition, if its commutator subgroup $H'$ is the whole $H$. By virtue of \cite[Corollary~4,~p.~285]{Bou}, this is equivalent to $H$ having a finite fundamental group $\pi_1(H)$. In particular, if $H$ is simply connected then it is semi-simple.
\item The homogeneous space $G/H$ from the corollary need not be a group, that is, $H$ may not be a normal subgroup of $G$. To give an example, recall that the second special unitary group $\text{SU}(2)$ is simply connected (being homeomorphic to $\mathbb S^3$) and hence semi-simple by the previous remark, and it is contained as a closed but not normal subgroup
\begin{equation*}
\text{SU}(2)\subseteq\text{U}(2)\cong\left(
\begin{array}{ccc}
\text{U}(2) & 0 \\
0 & 1 
\end{array} \right)\subseteq\text{U}(3)
\end{equation*}
of the third unitary group $\text{U}(3)$.
\end{itemize}
\end{remark}
\begin{proof}[Proof of Corollary~\ref{C:min.eq.semis.con}]
Let $i\colon H\to G$ be the inclusion morphism and $i_1^{\sharp}\colon\pi_1(H)\to\pi_1(G)$ be the morphism induced by $i$. Similarly, let $q\colon G\to G/H$ be the canonical projection and $q_1^{\sharp}\colon\pi_1(G)\to\pi_1(G/H)$, $q^{\sharp}\colon H_1^w(G)\to H_1^w(G/H)$ be the morphisms induced by $q$. By choosing $z=H$ as the base point for $G/H$, we obtain $\mathcal F'_z=qh$ and hence
\begin{equation*}
(\mathcal F'_z)^{\sharp}=q^{\sharp}h^{\sharp}\colon H_1^w(\Gamma)\to H_1^w(G/H).
\end{equation*}
Thus, by virtue of Proposition~\ref{P:min.eq.flow.TFFC}, in order to prove the corollary it suffices to show that $q^{\sharp}$ is an isomorphism.

Since $G$ is a locally trivial fibre bundle over $G/H$ with the fibre $H$ and since $G$, $G/H$ are arc-wise connected, the long exact homotopy sequence of the bundle yields exactness of the upper row of the diagram in Figure~\ref{Fig:min.eq.semis.con} below (see \cite[Section~9.8]{FomFuch}).
\begin{figure}[ht]
\[\minCDarrowwidth20pt\begin{CD}
\pi_1(H) @>i_1^{\sharp}>> \pi_1(G) @>q_1^{\sharp}>> \pi_1(G/H) @>\partial>> H/H_0 @>>> 0\\
@. @Vp_GVV @Vp_{G/H}VV  \\
@. H_1^w(G) @>>q^{\sharp}> H_1^w(G/H) 
\end{CD}\]
\caption{The flow $\mathcal F'$ is topologically free and possesses a free cycle}
\label{Fig:min.eq.semis.con}
\end{figure}
Moreover, the diagram is clearly commutative.

Now, since the group $H$ is connected, $H/H_0=0$ and so $\im(q_1^{\sharp})=\ker(\partial)=\pi_1(G/H)$. Consequently, $q_1^{\sharp}$ is an epimorphism and hence so is $q^{\sharp}$. Also, since the group $\pi_1(G)$ is abelian finitely generated, it follows that so is the group $\pi_1(G/H)$. Thus, both $\pi_1(G)$ and $\pi_1(G/H)$ contain their torsion subgroups as direct summands with the complementary summands $H_1^w(G)$ and $H_1^w(G/H)$, respectively. Let $c$ and $d$ be the largest torsion coefficients of $\pi_1(G)$ and $\pi_1(G/H)$, respectively. By an argument analogous to the one from Proposition~\ref{P:min.eq.flow.TFFC}, which involved the commutativity of the diagram in Figure~\ref{Fig:min.eq.flow.TFFC}, we obtain $d\ker(q^{\sharp})\subseteq\ker(q_1^{\sharp})$.

Further, since the group $\pi_1(H)$ is finite by the first part of Remark~\ref{R:min.eq.semis.con}, we have $\ker(q_1^{\sharp})=\im(i_1^{\sharp})\subseteq\tor(\pi_1(G))$ and so $c\ker(q_1^{\sharp})=0$. Consequently, $cd\ker(q^{\sharp})\subseteq c\ker(q_1^{\sharp})=0$, whence $\ker(q^{\sharp})=0$, since the group $H_1^w(G)$ is torsion-free. This shows that $q^{\sharp}$ is also a monomorphism, as was to be shown.
\end{proof}

\subsection{Group extensions}\label{Sub:group.ext}

Let $\Flow$ be a minimal flow with acting homeomorphisms $T_{\gamma}$ ($\gamma\in\Gamma$) and let $G\in\mathsf{AbTpGp}$. A continuous map $\mathcal C\colon\Gamma\times X\to G$ is called a (\emph{group}) \emph{extension of $\mathcal F$}\index{extension!of a flow} (or a \emph{$1$-cocycle over}\index{$1$-cocycle over a flow} $\mathcal F$) \emph{with values in} $G$ if it satisfies the \emph{cocycle identity}\index{cocycle identity}
\begin{equation}\label{Eq:coc.id.def}
\mathcal C(\alpha,T_{\beta}x)\mathcal C(\beta,x)=\mathcal C(\alpha\beta,x).
\end{equation}
A flow $\Flow$ and its extension $\mathcal C\colon\Gamma\times X\to G$ give rise to the induced flow\index{induced flow} $\mathcal F_{\mathcal C}\index[symbol]{$\mathcal F_{\mathcal C}$}\colon\Gamma\curvearrowright X\times G$ with acting homeomorphisms
\begin{equation*}
\widetilde{T}_{\gamma}\index[symbol]{$\widetilde{T}_{\gamma}$}(x,g)=\left(T_{\gamma}x,\mathcal C(\gamma,x)g\right),
\end{equation*}
where $\gamma\in\Gamma$, $x\in X$, $g\in G$. (On occasions we shall attribute to $\mathcal C$ the dynamical properties of $\mathcal F_{\mathcal C}$. For example, we say that $\mathcal C$ is minimal when $\mathcal F_{\mathcal C}$ is minimal.) Further, given $g\in G$, consider the vertical translation\index{vertical translation} of $X\times G$ by $g$:
\begin{equation*}
R_g\index[symbol]{$R_g$}\colon X\times G\ni(x,h)\mapsto(x,gh)\in X\times G.
\end{equation*}
Since $G$ is an abelian group, $R_g$ is an automorphism of $\mathcal F_{\mathcal C}$; that is, $R_g$ is a homeomorphism with $R_g\widetilde{T}_{\gamma}=\widetilde{T}_{\gamma}R_g$ for every $\gamma\in\Gamma$. It follows, in particular, that $\mathcal O_{\mathcal F_{\mathcal C}}(x,gh)=R_g(\mathcal O_{\mathcal F_{\mathcal C}}(x,h))$ and $\overline{\mathcal O}_{\mathcal F_{\mathcal C}}(x,gh)=R_g(\overline{\mathcal O}_{\mathcal F_{\mathcal C}}(x,h))$ for all $x\in X$ and $g,h\in G$.

The set of all extensions of $\mathcal F$ with values in a given group $G$ will be denoted by $\Cocc(G)$\index[symbol]{$\Cocc(G)$}. Clearly, $\Cocc(G)$ is an abelian group with point-wise group operations. If both $\Gamma$ and $X$ are locally compact then we equip $\Cocc(G)$ with the compact-open topology, that is, with the topology of uniform convergence on compact sets. Thus $\Cocc(G)$ becomes an abelian topological group. If a net $(\mathcal C_i)_{i\in J}$ converges to $\mathcal C$ in $\Cocc(G)$ then we write $\mathcal C_i\stackrel{ucs}{\longrightarrow}\mathcal C$\index[symbol]{$\mathcal C_i\stackrel{ucs}{\longrightarrow}\mathcal C$} or $\mathcal C_i\to\mathcal C$ u.c.s.\index{u.c.s. convergence} If the group $G$ is complete (in particular, if it is locally compact) then the group $\Cocc(G)$ is also complete; that is, every Cauchy net in $\Cocc(G)$ converges u.c.s. to an element of $\Cocc(G)$. If both $\Gamma$ and $X$ are locally compact second countable and $G$ is a Polish group then $\Cocc(G)$ is also a Polish group; that is, it is separable and completely metrizable. In such a case we have $\card(\Cocc(G))\leq\mathfrak{c}$.

If the acting group $\Gamma$ of the flow $\mathcal F$ is connected then the phase space $X$ of $\mathcal F$ is also connected by minimality of $\mathcal F$. In such a case every extension $\mathcal C\in\Cocc(G)$ takes its values in the identity component $G_0$ of $G$. Thus, since our main interest in this work is in connected topological transformation groups, we shall often restrict ourselves to extensions of $\mathcal F$ with values in connected groups $G$.

\subsection{Quasi-coboundaries}

Let $\Flow$ be a minimal flow and let $G\in\mathsf{AbTpGp}$. By the cocycle identity (\ref{Eq:coc.id.def}), the extensions $\mathcal C\in\Cocc(G)$ of $\mathcal F$ can intuitively be viewed as morphisms of groups $\Gamma\to G$ twisted by the flow $\mathcal F$. Those extensions $\mathcal C$ that ``untwist'' are called quasi-coboundaries. More precisely, let $h\colon\Gamma\to G$ be a topological morphism and consider the map
\begin{equation*}
\mathcal Q_h\index[symbol]{$\mathcal Q_h$}\colon\Gamma\times X\ni(\gamma,x)\mapsto h(\gamma)\in G.
\end{equation*}
Then $\mathcal Q_h\in\Cocc(G)$; we call $\mathcal Q_h$ the (\emph{pure}) \emph{quasi-coboundary}\index{quasi-coboundary over a flow} (\emph{over $\mathcal F$}) \emph{induced by $h$}. More generally, an extension $\mathcal C\in\Cocc(G)$ is called a \emph{quasi-coboundary}, if it is equivalent to $\mathcal Q_h$ for some topological morphism $h\colon\Gamma\to G$. An elementary argument shows that the pure quasi-coboundaries form a closed subgroup of $\Cocc(G)$ with the topology of u.c.s. convergence.

\subsection{Coboundaries and cohomology}\label{Sub:cob.and.chmlg}

Let $\Flow$ be a minimal flow and $G\in\mathsf{AbTpGp}$. An extension $\mathcal C\in\Cocc(G)$ is called a \emph{$1$-coboundary}\index{coboundary over a flow} (briefly, \emph{coboundary}) \emph{over} $\mathcal F$ if there is a continuous map $\xi\colon X\to G$ with $\mathcal C(\gamma,x)=\xi(T_{\gamma}x)\xi(x)^{-1}$ for all $\gamma\in\Gamma$ and $x\in X$; we write $\mathcal C=\co(\xi)$\index[symbol]{$\co(\xi)$} or $\mathcal C=\mathcal C_{\xi}$\index[symbol]{$\mathcal C_{\xi}$}. Such a map $\xi$, called a \emph{transfer function}\index{transfer function!for a $1$-coboundary} of $\mathcal C$, is not unique. Nevertheless, by minimality of $\mathcal F$, if $X$ is a pointed space then every coboundary over $\mathcal F$ has a unique base point preserving transfer function. The coboundaries over $\mathcal F$ with values in $G$ form a subgroup $\Cob(G)$\index[symbol]{$\Cob(G)$} of $\Cocc(G)$. The quotient group $\Cocc(G)/\Cob(G)$ will be denoted by $\Coch(G)$\index[symbol]{$\Coch(G)$} and the underlying short exact sequence in $\mathsf{AbGp}$ will be written as
\begin{equation}\label{Eq:ex.seq.ext.G}
0\longrightarrow\Cob(G)\stackrel{\mu_G}{\longrightarrow}\Cocc(G)\stackrel{\pi_G}{\longrightarrow}\Coch(G)\longrightarrow0.\index[symbol]{$\mu_G$}\index[symbol]{$\pi_G$}
\end{equation}
The congruence on $\Cocc(G)$ determined by $\Cob(G)$ is called the \emph{cohomology relation}\index{cohomology of extensions} (or the \emph{equivalence of extensions})\index{equivalent extensions!of a flow} and we shall denote it by the symbol $\simeq$.\index[symbol]{$\mathcal C\simeq\mathcal D$}

To simplify notation, we shall write $\Cocc$,\index[symbol]{$\Cocc$} $\Cob$,\index[symbol]{$\Cob$} $\Coch$,\index[symbol]{$\Coch$} $\mu$,\index[symbol]{$\mu$} $\pi$\index[symbol]{$\pi$} instead of $\Cocc(\mathbb T^1)$, $\Cob(\mathbb T^1)$, $\Coch(\mathbb T^1)$, $\mu_{\mathbb T^1}$, $\pi_{\mathbb T^1}$, respectively. Thus, the exact sequence (\ref{Eq:ex.seq.ext.G}) for $G=\mathbb T^1$ takes the form
\begin{equation}\label{Eq:ex.seq.ext.T1}
0\longrightarrow\Cob\stackrel{\mu}{\longrightarrow}\Cocc\stackrel{\pi}{\longrightarrow}\Coch\longrightarrow0.
\end{equation}

\subsection{Groupoids of minimal extensions}\label{Sub:gpd.min.ext}

Given a minimal flow $\Flow$ and a group $G\in\mathsf{CAGp}$, we denote by $\Cocm(G)$\index[symbol]{$\Cocm(G)$} the set of all minimal extensions from $\Cocc(G)$ together with the identity $e$ of $\Cocc(G)$. By definition, $\Cocm(G)$ is a sub-groupoid of the group $\Cocc(G)$; we call it the \emph{groupoid of minimal extensions}\index{groupoid!of minimal extensions} (of $\mathcal F$ with values in $G$). Further, we set $\Cochm(G)=\pi_G(\Cocm(G))$.\index[symbol]{$\Cochm(G)$} Clearly, $\Cochm(G)$ is a sub-groupoid of the group $\Coch(G)$. Similarly to the preceding subsection, we write $\Cocm$\index[symbol]{$\Cocm$} and $\Cochm$\index[symbol]{$\Cochm$} instead of $\Cocm(\mathbb T^1)$ and $\Cochm(\mathbb T^1)$, respectively.

\subsection{Extensions and fundamental group}

Our aim in this subsection is to show how the morphisms induced by group extensions act on fundamental groups and first weak homology groups. This is done in Lemma~\ref{L:tor.Lie.mnfld} below. Let us start with the following motivating example.

\begin{example}\label{E:ind.fund.eqc}
Let $\Gamma,G$ be connected Lie groups, $G$ compact, and let $h\in\Hom(\Gamma,G)$ be a morphism with a dense image. Consider the minimal equi\-con\-ti\-nuous flow $\mathcal F\colon\Gamma\curvearrowright G$ induced by $h$ (see (\ref{Eq:min.eq.flow.gen})). By choosing the identity $z=e$ of $G$ as the base point for $G$, we get $\mathcal F_z=h$. Now, given $\xi\in\Hom(G,\mathbb T^1)$, we have $\co(\xi)(\gamma,g)=\xi(h(\gamma)g)\xi(g)^{-1}=\xi h(\gamma)$ for all $\gamma\in\Gamma$ and $g\in G$. Consequently,
\begin{equation}\label{Eq: ind.fund.eqc1}
\co(\xi)^{\sharp}\pi_1(\Gamma\times G)=\xi^{\sharp}h^{\sharp}\pi_1(\Gamma)=\xi^{\sharp}\mathcal F_z^{\sharp}\pi_1(\Gamma).
\end{equation}
Moreover, it follows from Lemma~\ref{L:morph.ind.morph} and Remark~\ref{R:morph.ind.morph} in Section~\ref{S:ex.prscr.sect.2} that every $\xi\in C_z(G,\mathbb T^1)$ is homotopic to a (unique) morphism from $\Hom(G,\mathbb T^1)$. Consequently, the identity (\ref{Eq: ind.fund.eqc1}) is true for every continuous base point preserving map $\xi\colon G\to\mathbb T^1$. The following lemma shows that it is in fact true over an arbitrary (minimal) flow. 
\end{example}

\begin{lemma}\label{L:tor.Lie.mnfld}
Let $\Flow$ be a minimal flow with $\Gamma\in\mathsf{LieGp}$ connected and with $X$ a compact connected manifold. Choose a base point $z$ for $X$ and fix $\mathcal C\in\Cocc$, $\xi\in C_z(X,\mathbb T^1)$. Then for all $f\in\pi_1(\Gamma)$ and $g\in\pi_1(X)$,
\begin{equation}\label{Eq:tor.Lie.mnf.L2}
\mathcal C^{\sharp}\left((f,g)\right)=\mathcal C_z^{\sharp}(f)\hspace{3mm}\text{and}\hspace{3mm}\co(\xi)^{\sharp}\left((f,g)\right)=\xi^{\sharp}\mathcal F_z^{\sharp}(f),
\end{equation}
where $\mathcal C_z=\mathcal C(-,z)\colon\Gamma\to \mathbb T^1$ and $\mathcal F_z=\mathcal F(-,z)\colon\Gamma\to X$. Consequently,
\begin{equation}\label{Eq:tor.Lie.mnf.L1}
\mathcal C^{\sharp}\pi_1(\Gamma\times X)=\mathcal C_z^{\sharp}\pi_1(\Gamma)\hspace{3mm}\text{and}\hspace{3mm}\co(\xi)^{\sharp}\pi_1(\Gamma\times X)=\xi^{\sharp}\mathcal F_z^{\sharp}\pi_1(\Gamma).
\end{equation}
Moreover, the diagram in Figure~\ref{Fig:tor.Lie.mnf.L1} below commutes, its first three vertical arrows from the left are epimorphisms and the fourth one is an isomorphism. Finally,
\begin{equation}\label{Eq:tor.Lie.mnf.L3.hmlg}
\mathcal C^{\sharp}H_1^w(\Gamma\times X)=\mathcal C_z^{\sharp}H_1^w(\Gamma)\hspace{3mm}\text{and}\hspace{3mm}\co(\xi)^{\sharp}H_1^w(\Gamma\times X)=\xi^{\sharp}H_1^w(\mathcal F).
\end{equation}
\end{lemma}
\begin{figure}[ht]
\[\minCDarrowwidth20pt\begin{CD}
\pi_1(\Gamma) @>\mathcal F_z^{\sharp}>> \im(\mathcal F_z^{\sharp}) @. {}\subseteq{} @. \pi_1(X) @>\xi^{\sharp}>> \pi_1(\mathbb T^1) @. {}\cong{} @. \mathbb Z\\
@Vp_{\Gamma}VV @Vp_XVV @. @Vp_XVV @Vp_{\mathbb T^1}VV @. @|\\
H_1^w(\Gamma) @>>\mathcal F_z^{\sharp}> H_1^w(\mathcal F) @. {}\subseteq{} @. H_1^w(X) @>>\xi^{\sharp}> H_1^w(\mathbb T^1) @. {}\cong{} @. \mathbb Z
\end{CD}\]
\caption{Extensions and one-dimensional invariants}
\label{Fig:tor.Lie.mnf.L1}
\end{figure}
\begin{proof}[Proof of Lemma~\ref{L:tor.Lie.mnfld}]
Before turning to the proof, notice that both $\mathcal F_z$ and $\mathcal C_z$ are (con\-ti\-nuous) base point preserving maps. We denote by $c_1$, $c_z$ and $c_e$ the constant loop based at $1\in\Gamma$, $z\in X$ and $e=1\in\mathbb T^1$, respectively.

Fix $f\in\pi_1(\Gamma)$ and $g\in\pi_1(X)$. Then $\mathcal C^{\sharp}((f,c_z))=\mathcal C_z^{\sharp}(f)$ by definition of $\mathcal C_z$ and $\mathcal C^{\sharp}((c_1,g))=c_e$ by the cocycle identity. Thus,
\begin{equation*}
\mathcal C^{\sharp}\left((f,g)\right)=\mathcal C^{\sharp}\left((f,c_z)*(c_1,g)\right)=\mathcal C^{\sharp}\left((f,c_z)\right)*\mathcal C^{\sharp}\left((c_1,g)\right)=\mathcal C_z^{\sharp}(f)*c_e=\mathcal C_z^{\sharp}(f).
\end{equation*}
Similarly,
\begin{equation*}
\co(\xi)^{\sharp}\left((f,g)\right)=\co(\xi)_z^{\sharp}(f)=\xi\mathcal F_zf=\xi^{\sharp}\mathcal F_z^{\sharp}(f),
\end{equation*}
which verifies (\ref{Eq:tor.Lie.mnf.L2}). The two equalities from (\ref{Eq:tor.Lie.mnf.L1}) follow at once from the respective equalities in (\ref{Eq:tor.Lie.mnf.L2}).

The statement of the lemma concerning the diagram in Figure~\ref{Fig:tor.Lie.mnf.L1} follows from our discussion in Subsection~\ref{Sub:frst.wk.hm.gp} and from our definition of the group $H_1^w(\mathcal F)$. To verify (\ref{Eq:tor.Lie.mnf.L3.hmlg}), we use (\ref{Eq:tor.Lie.mnf.L1}) and the commutativity of the diagram in Figure~\ref{Fig:tor.Lie.mnf.L1}. Thus we obtain
\begin{equation*}
\begin{split}
\mathcal C^{\sharp}H_1^w(\Gamma\times X)&=\mathcal C^{\sharp}p_{\Gamma\times X}\pi_1(\Gamma\times X)=p_{\mathbb T^1}\mathcal C^{\sharp}\pi_1(\Gamma\times X)=p_{\mathbb T^1}\mathcal C_z^{\sharp}\pi_1(\Gamma)\\
&=\mathcal C_z^{\sharp}p_{\Gamma}\pi_1(\Gamma)=\mathcal C_z^{\sharp}H_1^w(\Gamma)
\end{split}
\end{equation*}
and
\begin{equation*}
\begin{split}
\co(\xi)^{\sharp}H_1^w(\Gamma\times X)&=\co(\xi)^{\sharp}p_{\Gamma\times X}\pi_1(\Gamma\times X)=p_{\mathbb T^1}\co(\xi)^{\sharp}\pi_1(\Gamma\times X)\\
&=p_{\mathbb T^1}\xi^{\sharp}\im(\mathcal F_z^{\sharp})=\xi^{\sharp}p_X\im(\mathcal F_z^{\sharp})=\xi^{\sharp}H_1^w(\mathcal F).
\end{split}
\end{equation*}
\end{proof}

\subsection{Banach spaces of extensions}\label{Sub:Ban.sp.of.ext}

Our aim in this section is to find the structure of a Banach space in $\Cocc(\mathbb R)$; this is done in Lemma~\ref{L:Cocc.Ban.sp.Gmm} below. First, let us recall some useful facts. Let $\Gamma\in\mathsf{LCGp}$ be a compactly generated group and $K\subseteq\Gamma$ be a symmetric compact set containing the identity $1$ of $\Gamma$ in its interior. Given $n\in\mathbb N$, we denote by $K^n$ the set of all products $\gamma_1\dots\gamma_n$, where $\gamma_i\in K$ for every $i=1,\dots,n$; clearly, $K^n$ is a compact subset of $\Gamma$. Observe that $K$ generates the group $\Gamma$ if and only if $\Gamma=\bigcup_{n\in\mathbb N}K^n$. If this is the case then for every compact set $F\subseteq\Gamma$ there is $n\in\mathbb N$ with $F\subseteq K^n$.

\begin{lemma}\label{L:Cocc.Ban.sp.Gmm}
Let $\Flow$ be a minimal flow with $\Gamma\in\mathsf{LCGp}$ compactly gene\-ra\-ted and with $X$ compact. Fix a generating symmetric compact set $K\subseteq\Gamma$ for $\Gamma$, which contains the identity $1$ of $\Gamma$ in its interior. Then $\Cocc(\mathbb R)$ is a real Banach space with the point-wise linear operations and with the norm
\begin{equation}\label{Eq:norm.ext}
\|\mathcal C\|=\sum_{n\in\mathbb N}\frac{1}{2^n}\max\big\{|\mathcal C(\gamma,x)| : \gamma\in K^n, x\in X\big\}.
\end{equation}
Moreover, the norm topology on $\Cocc(\mathbb R)$ coincides with the topology of uniform convergence on compact sets.
\end{lemma}
\begin{remark}\label{R:Cocc.Ban.sp.Gmm}
Recall from Subsection~\ref{Sub:group.ext} that if both $\Gamma$ and $X$ are additionally assumed second countable then the Banach space $\Cocc(\mathbb R)$ is separable.
\end{remark}
\begin{proof}[Proof of Lemma~\ref{L:Cocc.Ban.sp.Gmm}]
First, $\Cocc(\mathbb R)$ is clearly closed with respect to point-wise linear com\-bi\-na\-tions and so it is a real linear space. Further, given $\mathcal C\in\Cocc(\mathbb R)$, $n\in\mathbb N$, $\gamma_1,\dots,\gamma_n\in K$ and $x\in X$, we have
\begin{equation*}
\mathcal C(\gamma_1\dots\gamma_n,x)=\sum_{i=1}^{n-1}\mathcal C\left(\gamma_i,T_{\gamma_{i+1}\dots\gamma_n}x\right)+\mathcal C(\gamma_n,x)
\end{equation*}
by the cocycle identity. Consequently, for all $\gamma\in K^n$ and $x\in X$,
\begin{equation*}
|\mathcal C(\gamma,x)|\leq n\max\big\{|\mathcal C(\alpha,y)| : \alpha\in K, y\in X\big\}.
\end{equation*}
It follows that the series in (\ref{Eq:norm.ext}) converges and hence $\|\mathcal C\|$ is well defined for every $\mathcal C\in\Cocc(\mathbb R)$. It is a routine procedure to verify that $\|.\|$ satisfies all the axioms of a norm on $\Cocc(\mathbb R)$; for instance, $\|\mathcal C\|=0$ implies $\mathcal C=0$, since $K$ generates $\Gamma$.

We claim that the norm topology on $\Cocc(\mathbb R)$ coincides with the topology of uniform convergence on compact sets. To show this, let $(\mathcal C_i)$ be a net in $\Cocc(\mathbb R)$ and let $\mathcal C\in\Cocc(\mathbb R)$. Since each compact subset of $\Gamma\times X$ is contained in $K^n\times X$ for some $n\in\mathbb N$, we have an equivalence of the following statements:
\begin{itemize}
\item $\mathcal C_i\to\mathcal C$ uniformly on compact sets,
\item $\mathcal C_i\to\mathcal C$ uniformly on $K^n\times X$ for every $n\in\mathbb N$.
\end{itemize}
By definition of the norm $\|.\|$, the second of these two statements is equivalent to the convergence $\mathcal C_i\to\mathcal C$ in the norm topology. This proves the claim.

Finally, it follows from our discussion in Subsection~\ref{Sub:group.ext} that $\Cocc(\mathbb R)$ is complete in the topology of uniform convergence on compact sets. Thus, with the norm $\|.\|$ defined by (\ref{Eq:norm.ext}), $\Cocc(\mathbb R)$ is indeed a Banach space.
\end{proof}

\subsection{The category $\Coc$}\label{Sub:cat.of.ext.Coc}

Let $\Flow$ be a minimal flow. The extensions of $\mathcal F$ with values in compact abelian groups form a category $\Coc$\index[symbol]{$\mathsf{CAGpZ}_{\mathcal F}$} defined as follows. The objects of $\Coc$ are extensions $\mathcal C\in\Cocc(G)$ with $G\in\mathsf{CAGp}$. Further, given $G,H\in\mathsf{CAGp}$, $\mathcal C\in\Cocc(G)$ and $\mathcal D\in\Cocc(H)$, the set $\Hom(\mathcal C,\mathcal D)$\index[symbol]{$\Hom(\mathcal C,\mathcal D)$} of morphisms $\mathcal C\to\mathcal D$ consists of all (topological) morphisms $q\in\Hom(G,H)$ with $q\mathcal C=\mathcal D$. The composition of morphisms in $\Coc$ is simply their composition in $\mathsf{CAGp}$; this is clearly a correct definition. Finally, for all $G\in\mathsf{CAGp}$ and $\mathcal C\in\Cocc(G)$, the identity $\id_G$ on $G$ is an element of $\Hom(\mathcal C,\mathcal C)$. Thus, $\Coc$ is a category in its own right. The following proposition shows that $\Coc$ is always co-complete.

\begin{proposition}\label{P:CAGpZ.co-comp}
Let $\Flow$ be a minimal flow. Then $\Coc$ is a co-complete category.
\end{proposition}
\begin{proof}
Let $\mathfrak{C}=(\mathcal C_i\stackrel{p_{ik}}{\longleftarrow}\mathcal C_k)$ be an inverse system in $\Coc$, indexed by a directed set $J$; we shall find an inverse limit of $\mathfrak{C}$ in $\Coc$. First, for $j\in J$, let $G_j\in\mathsf{CAGp}$ be the group with $\mathcal C_j\in\Cocc(G_j)$. By definition of morphisms and their composition in $\Coc$, there is an inverse system $\mathfrak{G}=(G_i\stackrel{p_{ik}}{\longleftarrow}G_k)$ in $\mathsf{CAGp}$. Let $G\sbgp\prod_{j\in J}G_j$ be the usual inverse limit of $\mathfrak{G}$ in $\mathsf{CAGp}$ with the limit projections $p_j\colon G\to G_j$ (see Subsection~\ref{Sub:CAGp.co-cmplt} for the definition of $G$ and $p_j$). Consider the map
\begin{equation*}
\mathcal C\colon\Gamma\times X\ni(\gamma,x)\mapsto\left(\mathcal C_j(\gamma,x)\right)_{j\in J}\in G\subseteq\prod_{j\in J}G_j;
\end{equation*}
we show that $\mathcal C$ is an inverse limit of $\mathfrak{C}$ in $\Coc$ with $p_j$ ($j\in J$) as the limit projections.

First, since $p_{ik}\mathcal C_k=\mathcal C_i$ holds for all $i\leq k$, the map $\mathcal C$ does indeed take its values in the group $G$. Moreover, the cocycle identity and continuity of $\mathcal C$ are immediate, hence $\mathcal C\in\Cocc(G)$. Further, by our definition of $\mathcal C$ and of the maps $p_j$, we have $p_j\in\Hom(\mathcal C,\mathcal C_j)$ for every $j\in J$ and $p_{ik}p_k=p_i$ for all $i\leq k$. Thus, it remains to verify the universality of $\mathcal C$ and $(p_j)_{j\in J}$. To this end, let $H\in\mathsf{CAGp}$, $\mathcal D\in\Cocc(H)$ and $q_j\in\Hom(\mathcal D,\mathcal C_j)$ ($j\in J$) be such that $p_{ik}q_k=q_i$ for all $i\leq k$; we need to verify the following claim:
\begin{itemize}
\item there is a unique $r\in\Hom(\mathcal D,\mathcal C)$ with $p_jr=q_j$ for every $j\in J$.
\end{itemize}

We begin by verifying the existence part of the claim. By our assumptions above, we have $q_j\in\Hom(H,G_j)$ for every $j\in J$ and $p_{ik}q_k=q_i$ for all $i\leq k$. Thus, since $G$ is an inverse limit of $\mathfrak{G}$ in $\mathsf{CAGp}$ with the morphisms $p_j$ as the limit projections, there is a unique $r\in\Hom(G,H)$ with $p_jr=q_j$ for every $j\in J$. We need to show that $r\in\Hom(\mathcal D,\mathcal C)$ or, equivalently, that $r\mathcal D=\mathcal C$. Let $j\in J$; then
\begin{equation*}
p_j(r\mathcal D)=(p_jr)\mathcal D=q_j\mathcal D=\mathcal C_j=p_j\mathcal C.
\end{equation*}
Since the morphisms $p_j$ separate the points of $G$, it follows that $r\mathcal D=\mathcal C$. This verifies the existence part of the claim. The uniqueness part follows at once from the inclusion $\Hom(\mathcal D,\mathcal C)\subseteq\Hom(H,G)$.
\end{proof}

We shall find the following lemma useful.

\begin{lemma}\label{L:inv.sys.ind}
Let $\mathfrak G=(G_i\stackrel{p_{ik}}{\longleftarrow}G_k)$ be an inverse system in $\mathsf{CAGp}$ indexed by a directed set $J$, $G$ be an inverse limit of $\mathfrak{G}$ in $\mathsf{CAGp}$ and $p_j\colon G\to G_j$ ($j\in J$) be the limit projections. Given a minimal flow $\Flow$ and an extension $\mathcal C\in\Cocc(G)$, set $\mathcal C_j=p_j\mathcal C\in\Cocc(G_j)$ for every $j\in J$. Then there is an inverse system $\mathfrak C=(\mathcal C_i\stackrel{p_{ik}}{\longleftarrow}\mathcal C_k)$ in $\Coc$, $\mathcal C$ is its inverse limit and $p_j\colon\mathcal C\to\mathcal C_j$ ($j\in J$) are the limit projections.
\end{lemma}
\begin{proof}
Given $j\in J$, we have indeed $\mathcal C_j\in\Cocc(G_j)$, since $p_j\in\Hom(G,G_j)$. Further, given $i\leq k$, we have $p_{ik}p_k=p_i$ by the assumptions, hence $p_{ik}\mathcal C_k=p_{ik}p_k\mathcal C=p_i\mathcal C=\mathcal C_i$ and so $p_{ik}\in\Hom(\mathcal C_k,\mathcal C_i)$. Thus, $\mathfrak{C}$ is an inverse system in $\Coc$. Also, by definition of the extensions $\mathcal C_j$, we have $p_j\in\Hom(\mathcal C,\mathcal C_j)$ for every $j\in J$.

To show that $\mathcal C$ is an inverse limit of $\mathfrak{C}$ in $\Coc$ with the limit projections $p_j$, fix $H\in\mathsf{CAGp}$, $\mathcal D\in\Cocc(H)$ and $q_j\in\Hom(\mathcal D,\mathcal C_j)$ ($j\in J$) with $p_{ik}q_k=q_i$ for $i\leq k$; we need to verify the following claim:
\begin{enumerate}
\item[($*$)] there is a unique $r\in\Hom(\mathcal D,\mathcal C)$ with $p_jr=q_j$ for every $j\in J$.
\end{enumerate}
Since $G$ is an inverse limit of $\mathfrak{G}$ in $\mathsf{CAGp}$ with the limit projections $p_j$, there is a unique $r\in\Hom(H,G)$ with $p_jr=q_j$ for every $j\in J$. As $\Hom(\mathcal D,\mathcal C)\subseteq\Hom(H,G)$, this immediately implies the uniqueness part of ($*$). To verify the existence part of ($*$), it remains to show that $r\in\Hom(\mathcal D,\mathcal C)$. This is done similarly to the proof of Proposition~\ref{P:CAGpZ.co-comp}. Indeed, given $j\in J$, we have
\begin{equation*}
p_j(r\mathcal D)=(p_jr)\mathcal D=q_j\mathcal D=\mathcal C_j=p_j\mathcal C.
\end{equation*}
Since the morphisms $p_j$ ($j\in J$) separate the points of $G$, it follows that $r\mathcal D=\mathcal C$. Thus, $r\in\Hom(\mathcal D,\mathcal C)$, as was to be shown.
\end{proof}

\subsection{Induced morphisms}\label{Sub:ind.morph.ext}
Let $\Flow$ be a minimal flow, $G,H\in\mathsf{CAGp}$ and $q\in\Hom(G,H)$. The morphism $q$ induces a morphism of groups
\begin{equation*}
\widehat{q}\index[symbol]{$\widehat{q}$}\colon\Cocc(G)\ni\mathcal C\mapsto q\mathcal C\in\Cocc(H);
\end{equation*}
if both $\Gamma$ and $X$ are locally compact then $\widehat{q}$ is a topological morphism. Clearly, the assignment $q\mapsto\widehat{q}$ is covariant functorial. Finally, since $\widehat{q}$ maps coboundaries to coboundaries, it gives rise to a morphism of groups $\Coch(G)\to\Coch(H)$, which we shall denote by the same symbol $\widehat{q}$.

Now let $G_j$ ($j\in J$) be a family of abelian topological groups. For $i\in J$ consider the projection morphism $\pr_i\colon\prod_{j\in J}G_j\to G_i$. Let 
\begin{equation*}
\varphi\colon\Cocc\Big(\prod_{j\in J}G_j\Big)\ni\mathcal C\mapsto\left(\text{pr}_j\mathcal C\right)_{j\in J}\in\prod_{j\in J}\Cocc(G_j).
\end{equation*}
Then $\varphi$ is a morphism of groups by the preceding paragraph. In fact, an elementary argument shows that $\varphi$ is an isomorphism. If, in addition, both $\Gamma$ and $X$ are locally compact then $\varphi$ is a topological isomorphism. In particular, for every cardinal $\mathfrak{k}$, there is a (topological) isomorphism $\Cocc(\mathbb T^{\mathfrak{k}})\cong(\Cocc)^{\mathfrak{k}}$.

\chapter{Fundamental Tools}\label{S:def.of.E}

\section{Extensions as group morphisms}\label{Sub:gp.ext}

Let $\Flow$ be a minimal flow. Our aim in this section is to show that the group $\Cocc(G)$ of the extensions of $\mathcal F$ with values in a given group $G\in\mathsf{CAGp}$ is (topologically) isomorphic to the group of morphisms $\Hom(G^*,\Cocc)$. We construct this isomorphism in Theorem~\ref{T:structure.coc} and then derive from it some useful corollaries. In particular, we show that the category $\Coc$ can be viewed as an opposite to a full subcategory of the category $\mathsf{Hom}(\mathsf{DAGp},\Cocc)$.

\begin{theorem}\label{T:structure.coc}
Let $\Flow$ be a minimal flow and $G\in\mathsf{CAGp}$. Define a map
\begin{equation}\label{Eq:strctr.coc.Phi}
\Phi_G\colon\Cocc(G)\to\Hom\left(G^*,\Cocc \right)\index[symbol]{$\Phi_G$}
\end{equation}
by the rule
\begin{equation*}
\Phi_G(\mathcal C)=\mathcal C^*\colon G^*\ni\chi\mapsto\chi\mathcal C\in\Cocc .
\end{equation*}
Then $\Phi_G$ is an isomorphism of groups. If, in addition, both $\Gamma$ and $X$ are locally compact, $\Cocc(G)$ and $\Cocc$ carry the topology of u.c.s. convergence and $\Hom(G^*,\Cocc)$ is equipped with the topology of point-wise convergence, then $\Phi_G$ is also a topological isomorphism.
\end{theorem}
\begin{remark}\label{R:structure.coc}
We wish to add the following remarks.
\begin{itemize}
\item Let us mention how the correspondence from the theorem affects induced morphisms between groups of extensions. So assume that $G,H\in\mathsf{CAGp}$ and $h\in\Hom(G,H)$. Recall that $h$ induces a morphism of groups $\widehat{h}\colon\Cocc(G)\to\Cocc(H)$ acting by the rule $\widehat{h}\colon\mathcal C\mapsto h\mathcal C$. One verifies at once that $\Phi_H\widehat{h}\Phi_G^{-1}(\mathcal C^*)=\mathcal C^*h^*$ for every $\mathcal C\in\Cocc(G)$, where $h^*$ stands for the morphism dual to $h$. Thus, under the correspondence from the theorem, $\widehat{h}$ corresponds to $\Phi_H\widehat{h}\Phi_G^{-1}=\Hom(h^*,\Cocc)$.
\item Because of the isomorphism (\ref{Eq:strctr.coc.Phi}), we shall refer to $\Coch(G)$ as the first cohomology group of $\mathcal F$ with coefficients in $G^*$. In particular, since the dual group $(\mathbb T^1)^*$ of $\mathbb T^1$ is isomorphic to $\mathbb Z$, the group $\Coch$ will be sometimes called the first cohomology group of $\mathcal F$ with integer coefficients.
\end{itemize}
\end{remark}
\begin{proof}[Proof of Theorem~\ref{T:structure.coc}]
It is a standard procedure to show that for every $\mathcal C\in\Cocc(G)$, $\Phi_G(\mathcal C)=\mathcal C^*$ is a morphism of groups $G^*\to\Cocc $. Similarly, one verifies at once that $\Phi_G$ is a morphism of groups. It is a monomorphism, since the characters of $G$ separate the points of $G$. We verify that $\Phi_G$ is also an epimorphism. To this end, define a map
\begin{equation*}
\Lambda_G\colon\Hom(G^*,\Cocc)\to\Cocc(G^{**})
\end{equation*}
by the rule
\begin{equation*}
\Lambda_G(\varphi)\colon \Gamma\times X\ni(\gamma,x)\mapsto\left(G^*\ni\chi\mapsto\varphi(\chi)(\gamma,x)\in\mathbb T^1\right)\in G^{**}
\end{equation*}
for all $\gamma\in\Gamma$ and $x\in X$. Again, one readily verifies that $\Lambda_G$ is a (well defined) monomorphism of groups. Now consider the Pontryagin isomorphism $\omega_G\colon G\to G^{**}$ and the induced isomorphism
\begin{equation*}
\widehat{\omega}_G\colon\Cocc(G)\ni\mathcal C\mapsto\omega_G\mathcal C\in\Cocc(G^{**}).
\end{equation*}
A simple computation shows that $\Lambda_G\Phi_G=\widehat{\omega}_G$. Since $\Lambda_G$ is a monomorphism and $\widehat{\omega}_G$ is an isomorphism, this shows that $\Phi_G$ is an epimorphism. Thus, indeed, $\Phi_G$ is an isomorphism of groups.

Now assume that both $\Gamma$ and $X$ are locally compact; we show that $\Phi_G$ is a topological isomorphism. To this end, fix a net $(\mathcal C_i)_{i\in J}$ in $\Cocc(G)$ and an extension $\mathcal C\in\Cocc(G)$. We show that $\mathcal C_i\to\mathcal C$ u.c.s. on $\Gamma\times X$ if and only if $\mathcal C^*_i\to\mathcal C^*$ point-wise on $G^*$.

Assume, first, that $\mathcal C_i\to\mathcal C$ u.c.s. on $\Gamma\times X$. Then for every $\chi\in G^*$, $\chi\mathcal C_i\to\chi\mathcal C$ u.c.s. on $\Gamma\times X$ or, equivalently, $\mathcal C^*_i(\chi)\to\mathcal C^*(\chi)$ in $\Cocc$. This means that $\mathcal C^*_i\to\mathcal C^*$ point-wise on $G^*$.

Assume, conversely, that $\mathcal C^*_i\to\mathcal C^*$ point-wise on $G^*$, and fix a compact set $K\subseteq\Gamma\times X$ and a neighbourhood $W$ of $e$ in $G$. By compactness of $G$, there exist $n\in\mathbb N$, neighbourhoods $V_1,\dots,V_n$ of $1$ in $\mathbb T^1$ and $\chi_1,\dots,\chi_n\in G^*$ such that $\bigcap_{k=1}^n\chi_k^{-1}(V_k)\subseteq W$. Since $\chi_k\mathcal C_i\to\chi_k\mathcal C$ in $\Cocc$ for every $k=1,\dots,n$, there is $i_0\in J$ with $\chi_k\mathcal C_i(\gamma,x)\left(\chi_k\mathcal C(\gamma,x)\right)^{-1}\in V_k$ for all $i\geq i_0$, $k=1,\dots,n$ and $(\gamma,x)\in K$. That is, for $i\geq i_0$ and $(\gamma,x)\in K$, $\mathcal C_i(\gamma,x)\mathcal C(\gamma,x)^{-1}\in W$. This shows that $\mathcal C_i\to\mathcal C$ u.c.s. on $\Gamma\times X$.
\end{proof}

\begin{corollary}\label{C:structure.coc.gen}
Let $\Flow$ be a minimal flow and $G,H\in\mathsf{CAGp}$. Equip the tensor product $G^*\otimes H^*$ with the discrete topology. Then there is an isomorphism of groups
\begin{equation*}
\Hom(G^*,\Cocc(H))\cong\Cocc((G^*\otimes H^*)^*).
\end{equation*}
The isomorphism is topological, provided both $\Gamma$ and $X$ are locally compact, $\Cocc(H)$ and $\Cocc((G^*\otimes H^*)^*)$ carry the topology of u.c.s. convergence and $\Hom(G^*,\Cocc(H))$ is equipped with the topology of point-wise convergence.
\end{corollary}
\begin{proof}
The first statement of the corollary follows from Theorem~\ref{T:structure.coc} and from basic facts concerning morphisms from tensor products of abelian groups:
\begin{equation}\label{Eq:Cocc.Hom.tens.gen}
\begin{split}
\Hom(G^*,\Cocc(H))&\cong\Hom(G^*,\Hom(H^*,\Cocc ))\cong\Hom(G^*\otimes H^*,\Cocc )\\
&\cong\Cocc((G^*\otimes H^*)^*).
\end{split}
\end{equation}
If $\Gamma$ and $X$ are locally compact then the first and the third isomorphism from (\ref{Eq:Cocc.Hom.tens.gen}) are topological by virtue of Theorem~\ref{T:structure.coc}. The second isomorphism is also topological, since all the groups $G^*$, $H^*$ and $G^*\otimes H^*$ carry the discrete topology and hence all the $\Hom$ groups in (\ref{Eq:Cocc.Hom.tens.gen}) carry the topology of point-wise convergence (see (\ref{Eq:tens.prod.hom.top}) and the discussion following it).
\end{proof}

\begin{corollary}\label{C:ex.seq.qt.gp}
Let $\Flow$ be a minimal flow, $G\in\mathsf{CAGp}$ and $H\sbgp G$. Then there is an exact sequence of abelian groups
\begin{equation}\label{Eq:sb.qt.gp.EH}
\begin{split}
0&\longrightarrow\Cocc(H)\stackrel{\widehat{j}}{\longrightarrow}\Cocc(G)\stackrel{\widehat{p}}{\longrightarrow}\Cocc(G/H)\\
&\longrightarrow\Ext(H^*,\Cocc )\longrightarrow\Ext(G^*,\Cocc )\longrightarrow\Ext(H^{\perp},\Cocc )\longrightarrow0,
\end{split}
\end{equation}
where $j$ is the inclusion morphism $H\to G$, $p$ is the quotient morphism $G\to G/H$ and $H^{\perp}$ stands for the annihilator of $H$ in $G$.
\end{corollary}
\begin{proof}
There is an exact sequence of abelian groups
\begin{equation}\label{Eq:sb.qt.gp.dual}
0\longrightarrow (G/H)^*\stackrel{p^*}{\longrightarrow}G^*\stackrel{j^*}{\longrightarrow}H^*\longrightarrow 0.
\end{equation}
Consider the contravariant Hom-Ext sequence derived from (\ref{Eq:sb.qt.gp.dual}) and corresponding to $\Cocc$:
\begin{equation}\label{Eq:sb.qt.gp.EH2}
\begin{split}
0&\longrightarrow\Hom(H^*,\Cocc)\xrightarrow{\Hom(j^*,\Cocc)}\Hom(G^*,\Cocc)\xrightarrow{\Hom(p^*,\Cocc)}\Hom((G/H)^*,\Cocc)\\
&\longrightarrow\Ext(H^*,\Cocc)\longrightarrow\Ext(G^*,\Cocc)
\longrightarrow\Ext((G/H)^*,\Cocc)\longrightarrow0.
\end{split}
\end{equation}
By virtue of Theorem~\ref{T:structure.coc} and Remark~\ref{R:structure.coc}, the diagram in Figure~\ref{Fig:ex.seq.qt.gp} below commutes,
\begin{figure}[ht]
\[\minCDarrowwidth50pt\begin{CD}
\Cocc(H) @>\widehat{j}>> \Cocc(G) @>\widehat{p}>> \Cocc(G/H) \\
@V\Phi_HVV @V\Phi_GVV @V\Phi_{G/H}VV \\
\Hom(H^*,\Cocc) @>>\Hom(j^*,\Cocc)> \Hom(G^*,\Cocc) @>>\Hom(p^*,\Cocc)> \Hom((G/H)^*,\Cocc)
\end{CD}\]
\caption{Obtaining (\ref{Eq:sb.qt.gp.EH}) from (\ref{Eq:sb.qt.gp.EH2})}
\label{Fig:ex.seq.qt.gp}
\end{figure}
where all $\Phi_H$, $\Phi_G$ and $\Phi_{G/H}$ are isomorphisms. Since, in addition, the groups $(G/H)^*$ and $H^{\perp}$ are isomorphic, the exact sequence (\ref{Eq:sb.qt.gp.EH2}) takes the desired form (\ref{Eq:sb.qt.gp.EH}).
\end{proof}

\begin{corollary}\label{C:coc.as.morph}
Given a minimal flow $\Flow$, there is a contravariant functor
\begin{equation*}
\mathfrak{F}\colon\Coc\to\mathsf{Hom}(\mathsf{DAGp},\Cocc)\index[symbol]{$\mathfrak{F}$}
\end{equation*}
defined as follows.
\begin{itemize}
\item[(i)] If $G\in\mathsf{CAGp}$ and $\mathcal C\in\Cocc(G)$ then $\mathfrak{F}(\mathcal C)=\mathcal C^*\in\Hom(G^*,\Cocc)$.
\item[(ii)] If $G,H\in\mathsf{CAGp}$, $\mathcal C\in\Cocc(G)$, $\mathcal D\in\Cocc(H)$ and $q\in\Hom(\mathcal C,\mathcal D)$ then $\mathfrak{F}(q)=q^*\in\Hom(\mathcal D^*,\mathcal C^*)$.
\end{itemize}
\end{corollary}
\begin{remark}\label{R:coc.as.morph.cat}
Under the action of $\mathfrak{F}$, distinct objects correspond to distinct objects and, likewise, distinct morphisms to distinct morphisms. Consequently, $\mathfrak{F}$ represents an isomorphism between $\Coc$ and the opposite of the subcategory $\mathfrak{F}(\Coc)$ of $\mathsf{Hom}(\mathsf{DAGp},\Cocc)$. Notice that the subcategory $\mathfrak{F}(\Coc)$ of $\mathsf{Hom}(\mathsf{DAGp},\Cocc)$ is full, that is, $\mathfrak{F}(\Hom(\mathcal C,\mathcal D))=\Hom(\mathfrak{F}(\mathcal D),\mathfrak{F}(\mathcal C))$ for all $\mathcal C,\mathcal D\in\Coc$. To see this, fix $G,H\in\mathsf{CAGp}$, $\mathcal C\in\Cocc(G)$, $\mathcal D\in\Cocc(H)$ and $r\in\Hom(\mathcal D^*,\mathcal C^*)$. Then $r\in\Hom(H^*,G^*)$ and hence $r=q^*$ for some $q\in\Hom(G,H)$. Since $(q\mathcal C)^*=\mathcal C^*q^*=\mathcal C^*r=\mathcal D^*$, it follows that $q\mathcal C=\mathcal D$. Thus, $q\in\Hom(\mathcal C,\mathcal D)$ and $\mathfrak{F}(q)=q^*=r$.
\end{remark}
\begin{proof}[Proof of Corollary~\ref{C:coc.as.morph}]
By virtue of Theorem~\ref{T:structure.coc}, $\mathfrak{F}$ assigns objects to objects. Also, $\mathfrak{F}$ assigns morphisms to morphisms, for in the situation from (ii), $\mathcal C^*q^*=(q\mathcal C)^*=\mathcal D^*$. To ve\-ri\-fy the axioms of a contravariant functor for $\mathfrak{F}$, fix groups $G,H,K\in\mathsf{CAGp}$ and extensions $\mathcal C,\mathcal D,\mathcal E\in\Coc$ with their values in $G,H,K$, respectively. The identity $\id_G\in\Hom(\mathcal C,\mathcal C)$ corresponds via $\mathfrak{F}$ to the identity $\id_{G^*}\in\Hom(\mathcal C^*,\mathcal C^*)$, for $(\id_G)^*=\id_{G^*}$. Also, if $q\in\Hom(\mathcal C,\mathcal D)$ and $r\in\Hom(\mathcal D,\mathcal E)$ then $\mathfrak{F}(rq)=\mathfrak{F}(q)\mathfrak{F}(r)$, for $(rq)^*=q^*r^*$.
\end{proof}

\section{Functorial approach to sections}\label{S:funct.appr.sect}

An important tool for our study of the first cohomology groups of a minimal flow $\mathcal F$ are sections of orbit closures in induced skew products. It turns out that this tool can be introduced formally, in the categorical language, as a unique covariant functor $F\colon\Coc\to\mathsf{CAGp}$ with three natural properties; this is shown in Theorem~\ref{T:functor.E.def}. This theorem and its reformulation in Theorem~\ref{T:the.functor.Fhat} are the main results of the present section. This abstract approach to sections will turn out to be very useful to us; we will often manage with the defining conditions of this functor (and with its other three properties listed in Theorem~\ref{T:functor.E.def}) and will not have to invoke its explicit construction. This will make many proofs throughout this work much more transparent.

\begin{theorem}\label{T:functor.E.def}
Given a minimal flow $\Flow$, there exists a unique covariant functor 
\begin{equation*}
F\colon\Coc\to\mathsf{CAGp},\index[symbol]{$F$}
\end{equation*}
which satisfies the following conditions:
\begin{enumerate}
\item[(1)] if $G\in\mathsf{CAGp}$ and $\mathcal C\in\Cocc(G)$ then $F(\mathcal C)\sbgp G$,\index[symbol]{$F(\mathcal C)$}
\item[(2)] for every $G\in\mathsf{CAGp}$ and $\mathcal C\in\Cocc(G)$, $F(\mathcal C)=e$ if and only if $\mathcal C\in\Cob(G)$,
\item[(3)] for all $\mathcal C,\mathcal D\in\Coc$ and every $q\in\Hom(\mathcal C,\mathcal D)$, $F(\mathcal D)=qF(\mathcal C)$ and $F(q)\index[symbol]{$F(q)$}\colon F(\mathcal C)\to F(\mathcal D)$ is the restriction of $q$.
\end{enumerate}
Moreover, the functor $F$ satisfies the following conditions:
\begin{enumerate}
\item[(4)] $F(\mathcal C)=F(\mathcal D)$ whenever $\mathcal C,\mathcal D\in\Coc$ satisfy $\mathcal C\simeq\mathcal D$,
\item[(5)] $F\colon\Coc\to\mathsf{CAGp}$ is continuous in the sense that it preserves limits of inverse systems,
\item[(6)] for every $G\in\mathsf{CAGp}$, $\mathcal C\in\Cocc(G)$ is minimal if and only if $F(\mathcal C)=G$.
\end{enumerate}
\end{theorem}

Before turning to the proof of Theorem~\ref{T:functor.E.def}, let us prove the following two corollaries.

\begin{corollary}\label{C:min.contr.char}
Let $\Flow$ be a minimal flow. Then for every $G\in\mathsf{CAGp}$, the isomorphism $\Phi_G$ from Theorem~\ref{T:structure.coc} restricts to an isomorphism
\begin{equation*}
\Phi_G\colon\Cob(G)\to\Hom(G^*,\Cob).
\end{equation*}
Also, the following conditions are equivalent:
\begin{enumerate}
\item[(a)] the sequence
\begin{equation*}
0\longrightarrow\Cob\stackrel{\mu}{\longrightarrow}\Cocc\stackrel{\pi}{\longrightarrow}\Coch\longrightarrow0
\end{equation*}
splits,
\item[(b)] the sequence
\begin{equation*}
0\longrightarrow\Cob(G)\stackrel{\mu_G}{\longrightarrow}\Cocc(G)\stackrel{\pi_G}{\longrightarrow}\Coch(G)\longrightarrow0
\end{equation*}
splits for every $G\in\mathsf{CAGp}$.
\end{enumerate}
Further, given $G\in\mathsf{CAGp}$ and $\mathcal C\in\Cocc(G)$, the following conditions are equivalent:
\begin{enumerate}
\item[(i)] $\mathcal C$ is minimal,
\item[(ii)] for every $\chi\in G^*\setminus 1$, $\chi\mathcal C\notin\Cob$.
\end{enumerate}
Finally, if $G\in\mathsf{CAGp}$ is connected then the isomorphism $\Phi_G$ from Theorem~\ref{T:structure.coc} restricts to an isomorphism of groupoids
\begin{equation*}
\Phi_G\colon\Cocm(G)\to\Mon(G^*,\Cocm).
\end{equation*}
\end{corollary}
\begin{proof}
Fix $G\in\mathsf{CAGp}$. By Theorem~\ref{T:functor.E.def}(2, 3), the following conditions are equivalent for every $\mathcal C\in\Cocc(G)$:
\begin{itemize}
\item $\mathcal C\in\Cob(G)$,
\item $F(\mathcal C)=e$,
\item $\chi F(\mathcal C)=1$ for every $\chi\in G^*$,
\item $F(\chi\mathcal C)=1$ for every $\chi\in G^*$,
\item $\mathcal C^*(\chi)=\chi\mathcal C\in\Cob$ for every $\chi\in G^*$.
\end{itemize}
This verifies the first statement of the corollary.

We show that conditions (a) and (b) are equivalent. Implication (b)$\Rightarrow$(a) is clear. To verify the converse, assume that $\Cob$ is a direct summand in $\Cocc$ and fix $G\in\mathsf{CAGp}$. Then $\Hom(G^*,\Cob)$ is a direct summand in $\Hom(G^*,\Cocc)$. Consequently, by the first statement of the corollary and by Theorem~\ref{T:structure.coc}, $\Cob(G)=\Phi_G^{-1}(\Hom(G^*,\Cob))$ is a direct summand in $\Cocc(G)=\Phi_G^{-1}(\Hom(G^*,\Cocc))$. This verifies condition (b).

Fix $G\in\mathsf{CAGp}$ and $\mathcal C\in\Cocc(G)$; we show that conditions (i) and (ii) are equi\-va\-lent. By Theorem~\ref{T:functor.E.def}(2, 3, 6), the following conditions are equivalent:
\begin{itemize}
\item $\mathcal C$ is minimal,
\item $F(\mathcal C)=G$,
\item $\chi F(\mathcal C)\neq1$ for every $\chi\in G^*\setminus1$,
\item $F(\chi\mathcal C)\neq1$ for every $\chi\in G^*\setminus1$,
\item $\chi\mathcal C\notin\Cob$ for every $\chi\in G^*\setminus1$.
\end{itemize}
This verifies the equivalence of (i) and (ii).

We verify the last statement of the corollary. If $G\in\mathsf{CAGp}$ is connected then the following conditions are equivalent for every $\mathcal C\in\Cocc(G)$ by Theorem~\ref{T:functor.E.def}(3, 6):
\begin{itemize}
\item $\mathcal C$ is minimal, that is, $F(\mathcal C)=G$,
\item $\chi F(\mathcal C)=\mathbb T^1$ for every $\chi\in G^*\setminus1$,
\item $F(\chi\mathcal C)=\mathbb T^1$ for every $\chi\in G^*\setminus1$,
\item $\mathcal C^*(\chi)=\chi\mathcal C$ is minimal for every $\chi\in G^*\setminus1$.
\end{itemize}
Moreover, the last condition is clearly equivalent to
\begin{itemize}
\item $\mathcal C^*$ is a monomorphism with values in $\Cocm$.
\end{itemize}
This leads to the equality $\Phi_G(\Cocm(G)\setminus e)=\Mon(G^*,\Cocm)\setminus1$. Consequently, $\Phi_G$ maps $\Cocm(G)$ onto $\Mon(G^*,\Cocm)$, as was to be shown.
\end{proof}

\begin{corollary}\label{C:MinH.rel.torH}
Let $\Flow$ be a minimal flow. Then
\begin{equation*}
\Cochm\setminus 1=\Coch\setminus\tor(\Coch).
\end{equation*}
\end{corollary}
\begin{proof}
For $k\in\mathbb N$ let $\kappa_k$ be the $k$-endomorphism of $\mathbb T^1$. Then for every $\mathcal C\in\Coch$, the following conditions are equivalent by virtue of Theorem~\ref{T:functor.E.def}:
\begin{itemize}
\item $\mathcal C\in\tor(\Coch)$,
\item $\mathcal C^k=1$ for some $k\in\mathbb N$,
\item $\kappa_kF(\mathcal C)=F(\mathcal C^k)=1$ for some $k\in\mathbb N$,
\item $F(\mathcal C)\subseteq\mathbb Z_k$ for some $k\in\mathbb N$,
\item $F(\mathcal C)\neq\mathbb T^1$,
\item $\mathcal C\notin\Cochm\setminus1$.
\end{itemize}
Consequently, $\tor(\Coch)$ and $\Cochm\setminus1$ are complementary in $\Coch$, as was to be shown.
\end{proof}

We begin our proof of Theorem~\ref{T:functor.E.def} by verifying the uniqueness of the functor~$F$.

\begin{lemma}\label{L:at.most.1E}
Let $\Flow$ be a minimal flow. Then there exists at most one covariant functor $F\colon\Coc\to\mathsf{CAGp}$ with properties $(1)-(3)$ from Theorem~\ref{T:functor.E.def}.
\end{lemma}
\begin{proof}
Assume that $F_1,F_2$ are two such functors and fix $G\in\mathsf{CAGp}$ and $\mathcal C\in\Cocc(G)$. Denote by $q$ the quotient morphism $G\to G/F_1(\mathcal C)$. By condition (3), $F_1(q\mathcal C)=qF_1(\mathcal C)=e$ and hence, by virtue of (2), $q\mathcal C\in\Cob(G/F_1(\mathcal C))$. Applying (2) and (3) once more yields $e=F_2(q\mathcal C)=qF_2(\mathcal C)$. Thus, $F_2(\mathcal C)\subseteq F_1(\mathcal C)$. The converse inclusion follows by symmetry and so $F_1,F_2$ act in the same way on the objects of $\Coc$. The fact that they act in the same way also on morphisms of $\Coc$ follows immediately from condition (3).
\end{proof}

Now we turn to the construction of the functor $F$.

\begin{lemma}\label{L:prop.sec.el}
Let $\Flow$ be a minimal flow and let $z\in X$. Given $G\in\mathsf{CAGp}$ and $\mathcal C\in\Cocc(G)$, set $F(\mathcal C)\index[symbol]{$F(\mathcal C)$}=\overline{\mathcal O}_{\mathcal F_{\mathcal C}}(z,e)|_z$, the vertical $z$-section of the orbit closure of $(z,e)$ under the action of the induced skew product $\mathcal F_{\mathcal C}$ on $X\times G$. Then the following statements hold for every $G\in\mathsf{CAGp}$ and every $\mathcal C\in\Cocc(G)$:
\begin{enumerate}
\item[(a)] $F(\mathcal C)$ is a closed subgroup of $G$,
\item[(b)] for every $x\in X$, $\overline{\mathcal O}_{\mathcal F_{\mathcal C}}(z,e)|_x\in G/F(\mathcal C)$.
\end{enumerate}
\end{lemma}
\begin{remark}\label{R:prop.sec.el}
Observe that by definition of $F(\mathcal C)$, a point $g\in G$ lies in $F(\mathcal C)$ if and only if for every neighbourhood $U$ of $z$ in $X$ and every neighbourhood $V$ of $g$ in $G$ there is $\gamma\in\Gamma$ with $T_{\gamma}z\in U$ and $\mathcal C(\gamma,z)\in V$.
\end{remark}
\begin{proof}[Proof of Lemma~\ref{L:prop.sec.el}]
We verify statement (a). By definition, $F(\mathcal C)$ is a closed subset of $G$ containing the identity $e$. Thus, due to compactness of $G$, it remains to show that $F(\mathcal C)$ is a sub-semigroup of $G$. So fix $g,h\in F(\mathcal C)$, a neighbourhood $U$ of $z$ in $X$ and a neighbourhood $V$ of $e$ in $G$; we shall find $\gamma\in\Gamma$ with $T_{\gamma}z\in U$ and $\mathcal C(\gamma,z)\in ghV$. Fix an identity neighbourhood $W$ in $G$ with $W^2\subseteq V$. Since $h\in F(\mathcal C)$, there is $\alpha\in\Gamma$ with $T_{\alpha}z\in U$ and $\mathcal C(\alpha,z)\in hW$. Choose a neighbourhood $U'$ of $z$ in $X$ with $T_{\alpha}U'\subseteq U$ and $\mathcal C(\alpha\times U')\subseteq hW$. Since $g\in F(\mathcal C)$, there is $\beta\in\Gamma$ with $T_{\beta}z\in U'$ and $\mathcal C(\beta,z)\in gW$. Then
\begin{equation*}
\begin{split}
\widetilde{T}_{\alpha\beta}(z,e)&=(T_{\alpha\beta}z,\mathcal C(\alpha\beta,z))=(T_{\alpha}T_{\beta}z,\mathcal C(\alpha,T_{\beta}z)\mathcal C(\beta,z))\\
&\in T_{\alpha}U'\times\mathcal C(\alpha\times U')gW\subseteq U\times ghW^2\subseteq U\times ghV.
\end{split}
\end{equation*}
Thus, for $\gamma=\alpha\beta$, $T_{\gamma}z\in U$ and $\mathcal C(\gamma,z)\in ghV$, as was to be shown.

Now we turn to the proof of (b). For every $\gamma\in\Gamma$, $\widetilde{T}_{\gamma}\overline{\mathcal O}_{\mathcal F_{\mathcal C}}(z,e)=\overline{\mathcal O}_{\mathcal F_{\mathcal C}}(z,e)$. By applying the vertical $T_{\gamma}z$-section to this equality, we obtain $\overline{\mathcal O}_{\mathcal F_{\mathcal C}}(z,e)|_{T_{\gamma}z}=\mathcal C(\gamma,z)F(\mathcal C)$. Thus, (b) holds for every $x$ from the $\mathcal F$-orbit of $z$. Now let $x\in X$ be arbitrary. By letting elements of the $\mathcal F$-orbit of $z$ approach $x$ and by using compactness of $G$, we find $g_x\in G$ with $\overline{\mathcal O}_{\mathcal F_{\mathcal C}}(z,e)|_x\supseteq g_x\overline{\mathcal O}_{\mathcal F_{\mathcal C}}(z,e)|_z$. Likewise, by letting elements of the $\mathcal F$-orbit of $x$ approach $z$, we find $h_x$ with $\overline{\mathcal O}_{\mathcal F_{\mathcal C}}(z,e)|_z\supseteq h_x\overline{\mathcal O}_{\mathcal F_{\mathcal C}}(z,e)|_x$. These two inclusions yield
\begin{equation*}
F(\mathcal C)=\overline{\mathcal O}_{\mathcal F_{\mathcal C}}(z,e)|_z\supseteq h_x\overline{\mathcal O}_{\mathcal F_{\mathcal C}}(z,e)|_x\supseteq h_xg_x\overline{\mathcal O}_{\mathcal F_{\mathcal C}}(z,e)|_z=h_xg_xF(\mathcal C)
\end{equation*}
and hence $g_xF(\mathcal C)\subseteq\overline{\mathcal O}_{\mathcal F_{\mathcal C}}(z,e)|_x\subseteq h_x^{-1}F(\mathcal C)$. Consequently, $\overline{\mathcal O}_{\mathcal F_{\mathcal C}}(z,e)|_x=g_xF(\mathcal C)\in G/F(\mathcal C)$. This finishes the proof of (b).
\end{proof}

Till the end of our proof of Theorem~\ref{T:functor.E.def}, we shall assume that the point $z$ from Lemma~\ref{L:prop.sec.el} above is fixed (in other words, we shall assume that $X$ is a pointed space). Of course, the result of the construction from Lemma~\ref{L:prop.sec.el} above does not depend on the choice of $z$, this is a consequence of the uniqueness part of Theorem~\ref{T:functor.E.def}.

\begin{lemma}\label{L:exists.E}
Let $\Flow$ be a minimal flow. Then the functor $F$ defined in Lem\-ma~\ref{L:prop.sec.el} satisfies conditions $(1)-(3)$ from Theorem~\ref{T:functor.E.def}.
\end{lemma}
\begin{proof}
Condition (1) follows at once from Lemma~\ref{L:prop.sec.el}(a). To verify condition (2), fix $G\in\mathsf{CAGp}$ and $\mathcal C\in\Cocc(G)$. If $\mathcal C\in\Cob(G)$, let $\xi$ be the transfer function of $\mathcal C$ with $\xi(z)=e$. Then $\overline{\mathcal O}_{\mathcal F_{\mathcal C}}(z,e)$ is the graph of $\xi$ and so $F(\mathcal C)=\xi(z)=e$. Conversely, if $F(\mathcal C)=e$ then all the vertical sections of $\overline{\mathcal O}_{\mathcal F_{\mathcal C}}(z,e)$ are singletons by Lemma~\ref{L:prop.sec.el}(b). By compactness of $G$, $\overline{\mathcal O}_{\mathcal F_{\mathcal C}}(z,e)$ is the graph of a continuous map $\xi\colon X\to G$. We claim that $\mathcal C\in\Cob(G)$ and that $\xi$ is a transfer function of $\mathcal C$. Indeed, by definition of $\xi$, $\mathcal C(\gamma,z)=\xi(T_{\gamma}z)$ for every $\gamma\in\Gamma$ and hence 
\begin{equation*}
\mathcal C(\alpha,T_{\beta}z)=\mathcal C(\alpha\beta,z)\mathcal C(\beta,z)^{-1}=\xi(T_{\alpha}(T_{\beta}z))\xi(T_{\beta}z)^{-1}=\co(\xi)(\alpha,T_{\beta}z)
\end{equation*}
for all $\alpha,\beta\in\Gamma$. Consequently, $\mathcal C=\co(\xi)$ on $\Gamma\times\mathcal O_{\mathcal F}(z)$. By minimality of $\mathcal F$ it follows that $\mathcal C=\co(\xi)$ on the whole $\Gamma\times X$.

To verify condition (3), assume that $G,H\in\mathsf{CAGp}$, $q\in\Hom(G,H)$ and $\mathcal C\in\Cocc(G)$. The map $\Id_X\times q\colon X\times G\to X\times H$ is a morphism of flows $\mathcal F_{\mathcal C}\to\mathcal F_{q\mathcal C}$ and, by compactness of $G$, it is a closed map. Consequently,
\begin{equation*}
(\text{Id}_X\times q)\left(\overline{\mathcal O}_{\mathcal F_{\mathcal C}}(z,e)\right)=\overline{(\text{Id}_X\times q)\left(\mathcal O_{\mathcal F_{\mathcal C}}(z,e)\right)}=\overline{\mathcal O}_{\mathcal F_{q\mathcal C}}(z,e)
\end{equation*}
and hence
\begin{equation*}
qF(\mathcal C)=(\text{Id}_X\times q)\left(\overline{\mathcal O}_{\mathcal F_{\mathcal C}}(z,e)\right)|_z=\overline{\mathcal O}_{\mathcal F_{q\mathcal C}}(z,e)|_z=F(q\mathcal C).
\end{equation*}
Thus, the restriction $q\colon F(\mathcal C)\to F(q\mathcal C)$ is a well defined morphism of topological groups.
\end{proof}

\begin{lemma}\label{L:E.satisfies.6}
Let $\Flow$ be a minimal flow. Then the functor $F$ defined in Lem\-ma~\ref{L:prop.sec.el} satisfies condition $(4)$ from Theorem~\ref{T:functor.E.def}.
\end{lemma}
\begin{proof}
Fix $G\in\mathsf{CAGp}$ and $\mathcal C,\mathcal D\in\Cocc(G)$ with $\mathcal C\simeq\mathcal D$. By symmetry of the cohomology relation, it is sufficient to show that $F(\mathcal C)\subseteq F(\mathcal D)$. Let $q\colon G\to G/F(\mathcal D)$ be the quotient morphism. Then, by Lemma~\ref{L:exists.E}, $F(q\mathcal D)=qF(\mathcal D)=e$ and hence $q\mathcal D\in\Cob(G/F(\mathcal D))$. Consequently, it follows by cohomology of $\mathcal C$ and $\mathcal D$ that $q\mathcal C\in\Cob(G/F(\mathcal D))$ and so $qF(\mathcal C)=F(q\mathcal C)=e$ by Lemma~\ref{L:exists.E}. This verifies the desired inclusion $F(\mathcal C)\subseteq F(\mathcal D)$.
\end{proof}

Now we want to verify that the functor $F$ defined in Lemma~\ref{L:prop.sec.el} satisfies condition $(5)$ from Theorem~\ref{T:functor.E.def}. As a matter of fact, in the following lemma we show that $F$ satisfies a slightly stronger condition than continuity from Theorem~\ref{T:functor.E.def}.

\begin{lemma}\label{P:cohom.inv}
Let $\Flow$ be a minimal flow, $G_i\stackrel{p_{ik}}{\longleftarrow} G_k$ be an inverse system in $\mathsf{CAGp}$ indexed by a directed set $J$ and, for $j\in J$, let $\mathcal C_j\in\Cocc(G_j)$. Assume that $\mathcal C_i\simeq p_{ik}\mathcal C_k$ for all $i\leq k$. Consider the extension $\mathcal C\in\Cocc(\prod_{j\in J}G_j)$ given by $\mathcal C(\gamma,x)=(\mathcal C_j(\gamma,x))_{j\in J}$ for $\gamma\in\Gamma$ and $x\in X$. Then
\begin{equation}\label{Eq:F.cont.funct}
F(\mathcal C)=\lim_{\longleftarrow}\left(F(\mathcal C_i)\stackrel{p_{ik}}{\longleftarrow}F(\mathcal C_k)\right)=\bigcap_{j\in J}\text{\emph{pr}}_j^{-1}F(\mathcal C_j),
\end{equation}
where $\pr_j$ ($j\in J$) are the limit projections for the inverse system $G_i\stackrel{p_{ik}}{\longleftarrow} G_k$.
\end{lemma}
\begin{remark}\label{R:cohom.inv}
The assumptions of the lemma yield an inverse system in $\Coc$, provided the condition $\mathcal C_i\simeq p_{ik}\mathcal C_k$ is strengthened to $\mathcal C_i=p_{ik}\mathcal C_k$ for all $i\leq k$. In this case the inverse system $\mathcal C_i\stackrel{p_{ik}}{\longleftarrow}\mathcal C_k$ has many models of an inverse limit in $\Coc$, but all of them are isomorphic to (a co-restriction of) the extension $\mathcal C$ from the lemma. Recall that these isomorphisms hold in a strict sense, that is, they include the limit projections as well. In view of Lemma~\ref{L:exists.E} this means that formula (\ref{Eq:F.cont.funct}) holds for every model of an inverse limit of $\mathcal C_i\stackrel{p_{ik}}{\longleftarrow}\mathcal C_k$ in $\Coc$.
\end{remark}
\begin{proof}[Proof of Lemma~\ref{P:cohom.inv}]
Before turning to the proof we fix some notation. For $i\in J$ let $e_i$ be the unit of $G_i$ and $\pi_i$ be the projection morphism $\prod_{j\in J}G_j\to G_i$. The inverse limit $\lim_{\leftarrow}(G_i\stackrel{p_{ik}}{\longleftarrow}G_k)$ will be denoted by $G$ and will be considered in the usual way as the subgroup of $\prod_{j\in J}G_j$ consisting of the points $(g_j)_{j\in J}$ with $g_i=p_{ik}(g_k)$ for all $i\leq k$. Finally, for $i\leq k$ we fix a transfer function $\xi_{ik}\colon X\to G_i$ for the coboundary $\mathcal C_i(p_{ik}\mathcal C_k)^{-1}\in\Cob(G_i)$.

By Lemmas~\ref{L:exists.E} and \ref{L:E.satisfies.6}, $F(\mathcal C_i)=F(p_{ik}\mathcal C_k)=p_{ik}F(\mathcal C_k)$ for all $i\leq k$. It follows that the inverse system $G_i\stackrel{p_{ik}}{\longleftarrow} G_k$ restricts to an inverse system $F(\mathcal C_i)\stackrel{p_{ik}}{\longleftarrow}F(\mathcal C_k)$ and so the middle term in (\ref{Eq:F.cont.funct}) makes sense. Also,
\begin{equation*}
\begin{split}
\lim_{\longleftarrow}\left(F(\mathcal C_i)\stackrel{p_{ik}}{\longleftarrow}F(\mathcal C_k)\right)&=\left\{(g_j)_{j\in J}\in G : g_j\in F(\mathcal C_j)\hspace{1mm}\text{for every $j\in J$}\,\right\}\\
&=\bigcap_{j\in J}\text{pr}_j^{-1}F(\mathcal C_j),
\end{split}
\end{equation*}
which verifies the second equality in (\ref{Eq:F.cont.funct}).

We finish the proof by verifying the equality $F(\mathcal C)=\bigcap_{j\in J}\pr_j^{-1}F(\mathcal C_j)$ from (\ref{Eq:F.cont.funct}). To verify inclusion ``$\subseteq$'', fix $i\leq k$. By Lemma~\ref{L:exists.E},
\begin{equation*}
\left(\pi_i.(p_{ik}\pi_k)^{-1}\right)F(\mathcal C)=F\left(\left(\pi_i.(p_{ik}\pi_k)^{-1}\right)\mathcal C\right)=F\left(\mathcal C_i.(p_{ik}\mathcal C_k)^{-1}\right)=F(\co(\xi_{ik}))=e_i.
\end{equation*}
Thus, for all $i\leq k$ and every $g\in F(\mathcal C)$, $\pi_i(g)=p_{ik}\pi_k(g)$. Consequently, $F(\mathcal C)\subseteq G$ and hence $\pr_jF(\mathcal C)=\pi_jF(\mathcal C)=F(\pi_j\mathcal C)=F(\mathcal C_j)$ for every $j\in J$ by Lemma~\ref{L:exists.E}. This shows that $F(\mathcal C)\subseteq\bigcap_{j\in J}\pr_j^{-1}F(\mathcal C_j)$.

To verify the converse inclusion, fix a character $\chi\in F(\mathcal C)^{\perp}$ from the annihilator of $F(\mathcal C)$ in $\prod_{j\in J}G_j$. There exist a finite set $E\subseteq J$ and characters $\chi_j\in G_j^*$ ($j\in E$) such that $\chi=\prod_{j\in E}\chi_j\pi_j$. Since $\chi\in F(\mathcal C)^{\perp}$, we have $F(\chi\mathcal C)=\chi F(\mathcal C)=1$ and hence, by Lemma~\ref{L:exists.E}, $\chi\mathcal C\in\Cob$. Thus,
\begin{equation}\label{Eq:con.F.chi.chij}
\Cob\ni\chi\mathcal C=\prod_{j\in E}\chi_j\pi_j\mathcal C=\prod_{j\in E}\chi_j\mathcal C_j.
\end{equation}
Fix $k\in J$ with $j\leq k$ for every $j\in E$. Then
\begin{equation*}
\chi_j\mathcal C_j=\chi_j\left((p_{jk}\mathcal C_k)\co(\xi_{jk})\right)=\left((\chi_jp_{jk})\mathcal C_k\right)\co(\chi_j\xi_{jk})
\end{equation*}
for every $j\in E$ and hence, by virtue of (\ref{Eq:con.F.chi.chij}), $\Cob\ni(\prod_{j\in E}\chi_jp_{jk})\mathcal C_k$. Applying Lem\-ma~\ref{L:exists.E} once more yields
\begin{equation*}
\Big(\prod_{j\in E}\chi_jp_{jk}\Big)F(\mathcal C_k)=F\Big(\Big(\prod_{j\in E}\chi_jp_{jk}\Big)\mathcal C_k\Big)=1,
\end{equation*}
whence it follows that
\begin{equation*}
\begin{split}
\chi\left(\text{pr}_k^{-1}F(\mathcal C_k)\right)&=\Big(\prod_{j\in E}\chi_j\pi_j\Big)\left(\text{pr}_k^{-1}F(\mathcal C_k)\right)=\Big(\prod_{j\in E}\chi_j\text{pr}_j\Big)\left(\text{pr}_k^{-1}F(\mathcal C_k)\right)\\
&=\Big(\prod_{j\in E}\chi_jp_{jk}\text{pr}_k\Big)\left(\text{pr}_k^{-1}F(\mathcal C_k)\right)\subseteq\Big(\prod_{j\in E}\chi_jp_{jk}\Big)F(\mathcal C_k)=1.
\end{split}
\end{equation*}
This leads to
\begin{equation*}
\chi\in\left(\text{pr}_k^{-1}F(\mathcal C_k)\right)^{\perp}\subseteq\sum_{j\in J}\left(\text{pr}_j^{-1}F(\mathcal C_j)\right)^{\perp}=\Big(\bigcap_{j\in J}\text{pr}_j^{-1}F(\mathcal C_j)\Big)^{\perp},
\end{equation*}
where all the annihilators are taken in $\prod_{j\in J}G_j$. Thus, $F(\mathcal C)^{\perp}\subseteq(\bigcap_{j\in J}\pr_j^{-1}F(\mathcal C_j))^{\perp}$ and so, finally, $F(\mathcal C)\supseteq\bigcap_{j\in J}\pr_j^{-1}F(\mathcal C_j)$.
\end{proof}

\begin{lemma}\label{L:E.satisfies.5}
Let $\Flow$ be a minimal flow. Then the functor $F$ defined in Lem\-ma~\ref{L:prop.sec.el} satisfies condition $(6)$ from Theorem~\ref{T:functor.E.def}.
\end{lemma}
\begin{proof}
Fix $G\in\mathsf{CAGp}$ and $\mathcal C\in\Cocc(G)$. If $\mathcal F_{\mathcal C}$ is minimal then $\overline{\mathcal O}_{\mathcal F_{\mathcal C}}(z,e)=X\times G$ and so $F(\mathcal C)=\overline{\mathcal O}_{\mathcal F_{\mathcal C}}(z,e)|_z=G$. Conversely, assume that $F(\mathcal C)=G$; we show that every orbit of $\mathcal F_{\mathcal C}$ is dense in $X\times G$. First, we have $\overline{\mathcal O}_{\mathcal F_{\mathcal C}}(z,e)=X\times G$ by Lemma~\ref{L:prop.sec.el}(b). Further, given $g\in G$, the right vertical rotation $R_g$ of $X\times G$ by $g$ is an automorphism of $\mathcal F_{\mathcal C}$ and hence
\begin{equation*}
\overline{\mathcal O}_{\mathcal F_{\mathcal C}}(z,g)=\overline{\mathcal O}_{\mathcal F_{\mathcal C}}(R_g(z,e))=R_g\left(\overline{\mathcal O}_{\mathcal F_{\mathcal C}}(z,e)\right)=R_g(X\times G)=X\times G.
\end{equation*}
Finally, if $x\in X$ and $h\in G$ then $(z,g)\in\overline{\mathcal O}_{\mathcal F_{\mathcal C}}(x,h)$ for some $g\in G$ by compactness of $G$ and minimality of $\mathcal F$. Thus, $\overline{\mathcal O}_{\mathcal F_{\mathcal C}}(x,h)\supseteq\overline{\mathcal O}_{\mathcal F_{\mathcal C}}(z,g)=X\times G$. Since $x\in X$ and $h\in G$ were arbitrary, the minimality of $\mathcal F_{\mathcal C}$ follows.
\end{proof}

\begin{theorem}\label{T:the.functor.Fhat}
Given a minimal flow $\Flow$, there exists a unique covariant functor 
\begin{equation*}
F^*\index[symbol]{$F^*$}\colon\mathsf{Hom}(\mathsf{DAGp},\Cocc )\to\mathsf{DAGp},
\end{equation*}
which satisfies the following conditions:
\begin{enumerate}
\item[(1)] if $h\in\mathsf{Hom}(\mathsf{DAGp},\Cocc )$ has domain $D$ then $F^*(h)\sbgp D$,
\item[(2)] for every $h\in\mathsf{Hom}(\mathsf{DAGp},\Cocc )$ with domain $D$, $F^*(h)=D$ if and only if $\im(h)\subseteq \Cob $,
\item[(3)] if $h,k\in\mathsf{Hom}(\mathsf{DAGp},\Cocc)$ and $r\in\Hom(h,k)$, then $F^*(h)=r^{-1}(F^*(k))$ and $F^*(r)\colon F^*(h)\to F^*(k)$ is the restriction of $r$.
\end{enumerate}
Moreover, the functor $F^*$ satisfies the following conditions:
\begin{enumerate}
\item[(4)] if $h,k\in\mathsf{Hom}(\mathsf{DAGp},\Cocc )$ have a common domain $D$ and they differ by an element of $\Hom(D,\Cob )$ then $F^*(h)=F^*(k)$,
\item[(5)] for every $G\in\mathsf{CAGp}$ and every $\mathcal C\in\Cocc(G)$, $F^*(\mathcal C^*)$ equals the annihilator $F(\mathcal C)^{\perp}$ of $F(\mathcal C)$ in $G^*$.
\end{enumerate}
\end{theorem}
\begin{proof}
We start by showing that conditions (1)-(3) are satisfied by at most one functor $F^*$. More concretely, we show that $F^*(h)=h^{-1}(\Cob)$ is the only possible choice for $F^*$. So fix $D\in\mathsf{DAGp}$ and $h\in\Hom(D,\Cocc )$. Denote by $i$ the inclusion morphism $h^{-1}(\Cob)\to D$. Then $\im(hi)\subseteq\Cob $ and hence, by virtue of (2), $F^*(hi)=h^{-1}(\Cob)$. Thus, by virtue of (3),
\begin{equation*}
h^{-1}(\Cob)=F^*(hi)=i^{-1}(F^*(h))\subseteq F^*(h).
\end{equation*}
Now let $j$ be the inclusion morphism $F^*(h)\to D$. Then condition (3) yields $F^*(hj)=j^{-1}(F^*(h))=F^*(h)$ and hence, by condition (2), $h(F^*(h))=\im(hj)\subseteq\Cob $. Thus, $F^*(h)\subseteq h^{-1}(\Cob)$. It follows that
\begin{equation}\label{Eq:def.of.F*}
F^*(h)=h^{-1}(\Cob)
\end{equation}
is indeed the only possible choice for $F^*$.

We verify that $F^*$ defined by (\ref{Eq:def.of.F*}) satisfies conditions (1)-(4). First, if $D\in\mathsf{DAGp}$ and $h\in\Hom(D,\Cocc )$ then $F^*(h)$ is a subgroup of $D$ by definition and condition (1) thus holds. Condition (2) is satisfied obviously. To verify condition (3), fix $h,k\in\mathsf{Hom}(\mathsf{DAGp},\Cocc)$ and $r\in\Hom(h,k)$. The equality $F^*(h)=r^{-1}(F^*(k))$ follows directly from the equality $kr=h$ and from the definition of $F^*$, and this allows to define $F^*(r)$ as in (3). It is now easy to see that $F^*$ is a covariant functor $\mathsf{Hom}(\mathsf{DAGp},\Cocc )\to\mathsf{DAGp}$. Finally, to verify (4), let $D\in\mathsf{DAGp}$, $h,k\in\Hom(D,\Cocc)$ and $l\in\Hom(D,\Cob)$ be such that $h=l.k$. Since $\Cob $ is a subgroup of $\Cocc $,
\begin{equation*}
F^*(h)=F^*(l.k)=(l.k)^{-1}(\Cob)=k^{-1}(\Cob)=F^*(k).
\end{equation*}
This verifies condition (4).

We finish the proof by verifying condition (5). To this end, fix $G\in\mathsf{CAGp}$ and $\mathcal C\in\Cocc(G)$. Given $\chi\in G^*$, the following conditions are equivalent due to Theorem~\ref{T:functor.E.def}:
\begin{itemize}
\item $\chi\in F^*(\mathcal C^*)=(\mathcal C^*)^{-1}(\Cob)$,
\item $\chi\mathcal C=\mathcal C^*(\chi)\in\Cob$,
\item $\chi F(\mathcal C)=F(\chi\mathcal C)=1$,
\item $\chi\in F(\mathcal C)^{\perp}$.
\end{itemize}
This verifies condition (5).
\end{proof}

\section{Summation properties}\label{S:sum.prop}

In Theorem~\ref{T:functor.E.def} we saw that the functor $F\colon\Coc\to\mathsf{CAGp}$ introduced in Section~\ref{S:funct.appr.sect} is continuous in the sense that it preserves limits of inverse systems. In this section we investigate another continuity property of $F$, namely its relation to infinite sums in the group $\Cocc(G)$ and the hypersemigroup $2^G$. The main question studied here is the following. Assume that a series $\sum_{n=1}^{\infty}\mathcal C_n$ converges u.c.s. in $\Cocc(G)$. Does it follow that $F(\sum_{n=1}^{\infty}\mathcal C_n)=\overline{\sum_{n=1}^{\infty}F(\mathcal C_n)}$? If not, does at least one of the inclusions ``$\supseteq$'' or ``$\subseteq$'' hold?

By considering first the case of a finite sum, we find out in Proposition~\ref{P:essen.disj.fin} that an appropriate additional assumption in this problem is a mutual group-disjointness of the groups $F(\mathcal C_n)$. In Theorem~\ref{T:first.ineq.yes} we then show that, under this group-disjointness assumption, the inclusion ``$\supseteq$'' can be secured by passing to a subsequence $(\mathcal C_{k_n})_{n\in\mathbb N}$ of $(\mathcal C_n)_{n\in\mathbb N}$. Moreover, in Proposition~\ref{P:first.ineq.yes} we demonstrate that in order to secure this inclusion, passing to a subsequence is also a necessary procedure. Finally, we close this section by showing that the inclusion ``$\subseteq$'' fails to hold even when we allow passing to a subsequence of $(\mathcal C_n)_{n\in\mathbb N}$.

\begin{proposition}\label{P:essen.disj.fin}
Let $\mathcal F\colon\Gamma\curvearrowright X$ be a minimal flow and let $G_n\in\mathsf{CAGp}$ and $\mathcal C_n\in\Cocc(G_n)$ for $n=1,\dots,m$. If $F(\mathcal C_1),\dots,F(\mathcal C_m)$ are pair-wise group-disjoint then
\begin{equation}\label{Eq:prod.essen.disj}
F(\mathcal C_1,\dots,\mathcal C_m)=F(\mathcal C_1)\times\dots\times F(\mathcal C_m).
\end{equation}
If, in addition, all the groups $G_n=G$ are the same then
\begin{equation}\label{Eq:sum.essen.disj}
F\left(\sum_{n=1}^m\mathcal C_n\right)=\sum_{n=1}^mF(\mathcal C_n).
\end{equation}
\end{proposition}
\begin{proof}
First we verify (\ref{Eq:prod.essen.disj}), proceeding by induction on $m$. So let $m=2$; we show that $F(\mathcal C_1,\mathcal C_2)=F(\mathcal C_1)\times F(\mathcal C_2)$. For $i=1,2$ denote by $\pr_i$ the projection morphism $G_1\times G_2\to G_i$. Then $\pr_iF(\mathcal C_1,\mathcal C_2)=F(\pr_i(\mathcal C_1,\mathcal C_2))=F(\mathcal C_i)$ for both $i=1,2$ by Theorem~\ref{T:functor.E.def}(3). Consequently, $F(\mathcal C_1,\mathcal C_2)$ is a closed subgroup of $F(\mathcal C_1)\times F(\mathcal C_2)$ with full projections. Since $F(\mathcal C_1)$ and $F(\mathcal C_2)$ are group-disjoint by the assumptions, it follows that $F(\mathcal C_1,\mathcal C_2)=F(\mathcal C_1)\times F(\mathcal C_2)$.

Now fix $k\geq2$; we show that (\ref{Eq:prod.essen.disj}) holds for $m=k+1$, provided it holds for every $m\leq k$. To this end, fix extensions $\mathcal C_1,\dots,\mathcal C_k,\mathcal C_{k+1}$ with $F(\mathcal C_1),\dots,F(\mathcal C_k),F(\mathcal C_{k+1})$ mutually group-disjoint. By the induction hypothesis for $m=k$, $F(\mathcal C_1,\dots,\mathcal C_k)=F(\mathcal C_1)\times\dots\times F(\mathcal C_k)$. Consequently, by Lemma~\ref{L:prod.of.disj}, the groups $F(\mathcal C_1,\dots,\mathcal C_k)$ and $F(\mathcal C_{k+1})$ are group-disjoint. Our induction hypothesis for $m=2$ thus yields
\begin{equation*}
\begin{split}
F(\mathcal C_1,\dots,\mathcal C_{k+1})&=F((\mathcal C_1,\dots,\mathcal C_k),\mathcal C_{k+1})=F(\mathcal C_1,\dots,\mathcal C_k)\times F(\mathcal C_{k+1})\\
&=F(\mathcal C_1)\times\dots\times F(\mathcal C_k)\times F(\mathcal C_{k+1}),
\end{split}
\end{equation*}
as was to be shown.

Now we verify (\ref{Eq:sum.essen.disj}). So assume that all the groups $G_n=G$ are the same and consider the map
\begin{equation*}
\varphi\colon G^m\ni(g_1,\dots,g_m)\mapsto g_1+\dots+g_m\in G.
\end{equation*}
By commutativity of $G$, $\varphi$ is a morphism of topological groups. Thus, by Theorem~\ref{T:functor.E.def}(3) and by (\ref{Eq:prod.essen.disj}),
\begin{equation*}
\begin{split}
F\left(\sum_{n=1}^m\mathcal C_n\right)&=F(\varphi(\mathcal C_1,\dots,\mathcal C_m))=\varphi F(\mathcal C_1,\dots,\mathcal C_m)\\
&=\varphi\left(F(\mathcal C_1)\times\dots\times F(\mathcal C_m)\right)=\sum_{n=1}^mF(\mathcal C_n).
\end{split}
\end{equation*}
This verifies (\ref{Eq:sum.essen.disj}).
\end{proof}

\begin{corollary}\label{C:essen.disj.fin}
Let $\Flow$ be a minimal flow, $G\in\mathsf{CAGp}$ and $\mathcal C,\mathcal D\in\Cocc(G)$. Assume that $F(\mathcal C)$ is connected and $F(\mathcal D)$ is totally disconnected. Then $F(\mathcal C\mathcal D)=F(\mathcal C)F(\mathcal D)$.
\end{corollary}
\begin{remark}\label{R:essen.disj.fin}
Under the assumptions of the corollary, let the group $G$ be in addition connected and the extension $\mathcal C$ be minimal (that is, $F(\mathcal C)=G$). Then the extension $\mathcal C\mathcal D$ is also minimal. Indeed, the corollary yields $F(\mathcal C\mathcal D)=F(\mathcal C)F(\mathcal D)=GF(\mathcal D)=G$, which verifies minimality of $\mathcal C\mathcal D$.
\end{remark}
\begin{proof}[Proof of Corollary~\ref{C:essen.disj.fin}]
Recall that a connected group from $\mathsf{CAGp}$ has only connected quotient groups and a totally disconnected group from $\mathsf{CAGp}$ has only totally disconnected quotient groups. Thus, by the assumptions of the corollary, the groups $F(\mathcal C)$, $F(\mathcal D)$ have no common quotient group other than the trivial one and so they are group-disjoint by Lemma~\ref{L:grp.disj.char}. By virtue of (\ref{Eq:sum.essen.disj}) from Proposition~\ref{P:essen.disj.fin} it follows that, indeed, $F(\mathcal C\mathcal D)=F(\mathcal C)F(\mathcal D)$.
\end{proof}

\begin{theorem}\label{T:disj.inf.prod}
Let $\mathcal F\colon\Gamma\curvearrowright X$ be a minimal flow. For $n\in\mathbb N$ let $G_n\in\mathsf{CAGp}$ and $\mathcal C_n\in\Cocc(G_n)$. Set $\mathcal C=(\mathcal C_n)_{n\in\mathbb N}\in\Cocc\left(\prod_{n\in\mathbb N}G_n\right)$. If the groups $F(\mathcal C_n)$ are pair-wise group-disjoint then
\begin{equation*}
F(\mathcal C)=\prod_{n\in\mathbb N}F(\mathcal C_n).
\end{equation*}
\end{theorem}
\begin{proof}
If we interpret the product $\prod_{n\in\mathbb N}G_n$ in the usual way as an inverse limit of the groups $\prod_{n=1}^mG_n$ ($m\in\mathbb N$) then the statement of the theorem follows from Proposition~\ref{P:essen.disj.fin} and from Theorem~\ref{T:functor.E.def}(5).

To be more precise, for $i\leq k$ let $p_{ik}\colon\prod_{n=1}^kG_n\to\prod_{n=1}^iG_n$ be given by $p_{ik}(g_n)_{n=1}^k=(g_n)_{n=1}^i$ and for $i\in\mathbb N$ let $p_i\colon\prod_{n\in\mathbb N}G_n\to\prod_{n=1}^iG_n$ be given by $p_i(g_n)_{n\in\mathbb N}=(g_n)_{n=1}^i$. Then $\prod_{n=1}^iG_n\stackrel{p_{ik}}{\longleftarrow}\prod_{n=1}^kG_n$ is an inverse system in $\mathsf{CAGp}$, $\prod_{n\in\mathbb N}G_n$ is its inverse limit and $p_i$ ($i\in\mathbb N$) are the corresponding limit projections. By virtue of Lemma~\ref{L:inv.sys.ind}, $p_i\mathcal C\stackrel{p_{ik}}{\longleftarrow}p_k\mathcal C$ is an inverse system in $\Coc$, $\mathcal C$ is its inverse limit and $p_i\colon\mathcal C\to p_i\mathcal C$ ($i\in\mathbb N$) are the corresponding limit projections. Moreover, by Proposition~\ref{P:essen.disj.fin}, 
\begin{equation*}
F(p_i\mathcal C)=F(\mathcal C_1,\dots,\mathcal C_i)=F(\mathcal C_1)\times\dots\times F(\mathcal C_i)
\end{equation*}
for every $i\in\mathbb N$. Thus, the identity (\ref{Eq:F.cont.funct}) from Lemma~\ref{P:cohom.inv} yields
\begin{equation*}
\begin{split}
F(\mathcal C)&=\lim_{\longleftarrow}\left(F(\mathcal C_i)\stackrel{p_{ik}}{\longleftarrow}F(\mathcal C_k)\right)=\bigcap_{n\in\mathbb N}p_n^{-1}F(p_n\mathcal C)\\
&=\bigcap_{n\in\mathbb N}p_n^{-1}\left(F(\mathcal C_1)\times\dots\times F(\mathcal C_n)\right)=\prod_{n\in\mathbb N}F(\mathcal C_n),
\end{split}
\end{equation*}
as was to be shown.
\end{proof}

\begin{theorem}\label{T:first.ineq.yes}
Let $\Gamma\in\mathsf{LCGp}$, $X$ be a locally compact space and $G\in\mathsf{CAGp}$. Assume that all $\Gamma$, $X$ and $G$ are second countable. Let $\mathcal F\colon\Gamma\curvearrowright X$ be a minimal flow. For $n\in\mathbb N$ let $\mathcal C_n\in\Cocc(G)$ and assume that the groups $F(\mathcal C_n)$ are pair-wise group-disjoint. If $\mathcal C_n\stackrel{ucs}{\longrightarrow}e$ then there exists an increasing sequence $(k_n)_{n\in\mathbb N}$ of positive integers such that the series $\sum_{n=1}^{\infty}\mathcal C_{k_n}$ converges u.c.s. in $\Cocc(G)$ and
\begin{equation}\label{Eq:frst.inq.main}
F\left(\sum_{n=1}^{\infty}\mathcal C_{k_n}\right)\supseteq\overline{\sum_{n=1}^{\infty}F(\mathcal C_{k_n})},
\end{equation}
where the right-hand side of the inclusion above is defined as the closed subgroup of $G$ generated by the groups $F(\mathcal C_{k_n})$:
\begin{equation}\label{Eq:frst.inq.aux}
\overline{\sum_{n=1}^{\infty}F(\mathcal C_{k_n})}=\overline{\bigcup_{n\in\mathbb N}\left(\sum_{l=1}^nF(\mathcal C_{k_l})\right)}.
\end{equation}
\end{theorem}
\begin{remark}\label{R:first.ineq.yes}
Observe that the right-hand side of (\ref{Eq:frst.inq.aux}) can be expressed as a limit in the hypersemigroup $2^G$:
\begin{equation*}
\overline{\bigcup_{n\in\mathbb N}\left(\sum_{l=1}^nF(\mathcal C_{k_l})\right)}=\lim_{n\to\infty}\sum_{l=1}^nF(\mathcal C_{k_l}).
\end{equation*}
\end{remark}

Before turning to the proof of Theorem~\ref{T:first.ineq.yes}, we prove three auxiliary lemmas.

\begin{lemma}\label{L:convergent.subseries}
Under the assumptions of Theorem~\ref{T:first.ineq.yes}, there exists an increasing sequence of positive integers $(k_n)_{n\in\mathbb N}$ such that the series $\sum_{n=1}^{\infty}\mathcal C_{l_n}$ converges u.c.s. in $\Cocc(G)$ for every subsequence $(l_n)_{n\in\mathbb N}$ of $(k_n)_{n\in\mathbb N}$.
\end{lemma}
\begin{proof}
Let us begin by fixing the following objects:
\begin{itemize}
\item a local base $(W_n)_{n\in\mathbb N}$ at $e$ in $G$ with $W_{n+1}+W_{n+1}\subseteq W_n$ for every $n\in\mathbb N$,
\item a sequence $(K_n)_{n\in\mathbb N}$ of compact subsets of $\Gamma$, whose interiors cover $\Gamma$ and which satisfy $K_{n+1}\supseteq K_n$ for every $n\in\mathbb N$,
\item a sequence $(C_n)_{n\in\mathbb N}$ of compacts subsets of $X$, whose interiors cover $X$ and which satisfy $C_{n+1}\supseteq C_n$ for every $n\in\mathbb N$.
\end{itemize}
Since $\mathcal C_n\stackrel{ucs}{\longrightarrow} e$ in $\Cocc(G)$, there is an increasing sequence $(k_n)_{n\in\mathbb N}$ of positive integers such that for every $n\in\mathbb N$, $\mathcal C_{k_n}(\gamma,x)\in W_n$ for all $\gamma\in K_n$ and $x\in C_n$. We show that the sequence $(k_n)_{n\in\mathbb N}$ fulfills the condition from the conclusion of the lemma.

Let $(l_n)_{n\in\mathbb N}$ be a subsequence of $(k_n)_{n\in\mathbb N}$. Then $l_n=k_{m(n)}$ for every $n\in\mathbb N$, where $(m(n))_{n\in\mathbb N}$ is an increasing sequence of positive integers. (In particular, $m(n)\geq n$ for every $n\in\mathbb N$.) We show that the series $\sum_{n=1}^{\infty}\mathcal C_{l_n}$ converges u.c.s. in $\Cocc(G)$; we accomplish this by showing that $(\sum_{n=1}^m\mathcal C_{l_n})_{m\in\mathbb N}$ is a Cauchy sequence in $\Cocc(G)$. To this end, fix  a compact set $F\subseteq \Gamma\times X$ and a neighbourhood $W$ of $e$ in $G$. Choose $n_0\in\mathbb N$ large enough so that $F\subseteq K_{n_0}\times C_{n_0}$ and $W_{n_0}\subseteq W$. Then for all $n\geq n_0$, $p\in\mathbb N$ and $(\gamma,x)\in F$,
\begin{equation*}
\begin{split}
\sum_{j=1}^{n+p}\mathcal C_{l_j}(\gamma,x)-\sum_{j=1}^{n}\mathcal C_{l_j}(\gamma,x)&=\sum_{j=n+1}^{n+p}\mathcal C_{l_j}(\gamma,x)=\sum_{j=n+1}^{n+p}\mathcal C_{k_{m(j)}}(\gamma,x)\\
&\in\sum_{j=n+1}^{n+p}W_{m(j)}\subseteq \sum_{j=n+1}^{n+p}W_{j}\subseteq W_n\subseteq W_{n_0}\subseteq W.
\end{split}
\end{equation*}
Thus, $\left(\sum_{n=1}^{m}\mathcal C_{l_n}\right)_{m\in\mathbb N}$ is indeed a Cauchy sequence in $\Cocc(G)$. Since the group $\Cocc(G)$ is complete, it follows that the series $\sum_{n=1}^{\infty}\mathcal C_{l_n}$ converges u.c.s. in $\Cocc(G)$. This proves the lemma.
\end{proof}

\begin{lemma}\label{L:transit.coc.circ}
Let $\Flow$ be a minimal flow with $\Gamma$ and $X$ locally compact second countable. Assume that $(\mathcal C_n)_{n\in\mathbb N}$ is a sequence in $\Cocc$ with $\mathcal C_n\stackrel{ucs}{\longrightarrow}1$. If the groups $F(\mathcal C_n)$ are pair-wise group-disjoint and $F(\mathcal C_n)\neq1$ for every $n\in\mathbb N$ then there exists an increasing sequence $(k_n)_{n\in\mathbb N}$ of positive integers such that for every subsequence $(l_n)_{n\in\mathbb N}$ of $(k_n)_{n\in\mathbb N}$, the series $\sum_{n=1}^{\infty}\mathcal C_{l_n}$ converges u.c.s. in $\Cocc $ and
\begin{equation}\label{Eq:tr.coc.circ.two}
F\left(\sum_{n=1}^{\infty}\mathcal C_{l_n}\right)=\mathbb T^1=\overline{\sum_{n=1}^{\infty}F(\mathcal C_{l_n})}.
\end{equation}
\end{lemma}
\begin{proof}
First we indicate how to secure the second equality from (\ref{Eq:tr.coc.circ.two}). By our assumptions on the groups $F(\mathcal C_n)$ ($n\in\mathbb N$), we have $F(\mathcal C_i)=\mathbb T^1$ for at most one $i\in\mathbb N$. By removing the underlying extension $\mathcal C_i$ from our considerations, we get $1\neq F(\mathcal C_n)\subsetneq\mathbb T^1$ for every $n\in\mathbb N$. Thus, for every $n\in\mathbb N$ there is an integer $p(n)\geq2$ such that $F(\mathcal C_n)=\mathbb Z_{p(n)}$. Since the groups $F(\mathcal C_n)$ are pair-wise group-disjoint, it follows that the integers $p(n)$ are mutually co-prime. Consequently, if $(l_n)_{n\in\mathbb N}$ is an arbitrary increasing sequence of positive integers then
\begin{equation*}
\overline{\sum_{n=1}^{\infty}F(\mathcal C_{l_n})}=\overline{\sum_{n=1}^{\infty}\mathbb Z_{p(l_n)}}=\overline{\bigcup_{n\in\mathbb N}\sum_{j=1}^n\mathbb Z_{p(l_j)}}=\overline{\bigcup_{n\in\mathbb N}\mathbb Z_{p(l_1)\dots p(l_n)}}=\mathbb T^1.
\end{equation*}
This verifies the second equality from (\ref{Eq:tr.coc.circ.two}).

Now we turn to the first equality from (\ref{Eq:tr.coc.circ.two}). By Lemma~\ref{L:convergent.subseries}, we may assume that the series $\sum_{n=1}^{\infty}\mathcal C_{k_n}$ converges u.c.s. in $\Cocc $ for every increasing sequence of positive integers $(k_n)_{n\in\mathbb N}$. We fix the following objects:
\begin{itemize}
\item a base point $z\in X$,
\item a basis $(U_n)_{n\in\mathbb N}$ of neighbourhoods of $z$ in $X$ with $U_{n+1}\subseteq U_n$ for every $n\in\mathbb N$,
\item a basis $(V_n)_{n\in\mathbb N}$ for the topology of $\mathbb T^1$,
\item a sequence $(\widetilde{V}_n)_{n\in\mathbb N}$ of nonempty open subsets of $\mathbb T^1$ and a basis $(W_n)_{n\in\mathbb N}$ of neighbourhoods of $1$ in $\mathbb T^1$ such that $\widetilde{V}_n+W_n\subseteq V_n$ and $W_{n+1}+W_{n+1}\subseteq W_n$ for every $n\in\mathbb N$.
\end{itemize}

Now we turn to the next step of the proof. Proceeding by induction on $n$, we shall construct an increasing sequence $(k_n)_{n\in\mathbb N}$ of positive integers and a non-decreasing sequence $(F_n)_{n\in\mathbb N}$ of finite subsets of $\Gamma$. Properties of these sequences will be summarized after the induction step of the construction and their importance will become clear in the final step of the proof.

Since the integers $p(n)\geq2$ ($n\in\mathbb N$) are mutually co-prime, we have $p(n)\to\infty$ as $n\to\infty$ and so there is $k_1\in\mathbb N$ such that $F(\mathcal C_{k_1})\cap\widetilde{V}_1=\mathbb Z_{p(k_1)}\cap\widetilde{V}_1\neq\emptyset$. By definition of the functor $F$ from Lemma~\ref{L:prop.sec.el}, there is $\gamma_0\in\Gamma$ with $T_{\gamma_0}z\in U_1$ and $\mathcal C_{k_1}(\gamma_0,z)\in\widetilde{V}_1$. Set $F_1=\{\gamma_0\}$. Further, let $n\in\mathbb N$ and assume that we have at our disposal positive integers $k_1<\dots<k_n$ and finite sets $F_1\subseteq\dots\subseteq F_n\subseteq\Gamma$ such that for every $m=1,\dots,n$, the following conditions hold:
\begin{enumerate}
\item given $1\leq i\leq m$ and a (possibly empty) set $B\subseteq\{1,\dots,m-1\}$, there is $\gamma\in F_m$ with
\begin{equation*}
T_{\gamma}z\in U_m\hspace{3mm}\text{and}\hspace{3mm}\left(\sum_{l\in B}\mathcal C_{k_l}+\mathcal C_{k_m}\right)(\gamma,z)\in\widetilde{V}_i,
\end{equation*}
\item if $m\geq2$ then $\mathcal C_{k_m}(\gamma,z)\in W_m$ for every $\gamma\in F_{m-1}$.
\end{enumerate}
Now we define $k_{n+1}$ and $F_{n+1}$. Since $p(n)\to\infty$ and $\mathcal C_n\stackrel{ucs}{\longrightarrow}1$ as $n\to\infty$, there is $k_{n+1}>k_n$ such that 
\begin{itemize}
\item $F(\mathcal C_{k_{n+1}})\cap\widetilde{V}_i=\mathbb Z_{p(k_{n+1})}\cap\widetilde{V}_i\neq\emptyset$ for every $i=1,\dots,n+1$ and
\item $\mathcal C_{k_{n+1}}(\gamma,z)\in W_{n+1}$ for every $\gamma\in F_n$.
\end{itemize}
Let $B\subseteq\{1,\dots,n\}$ and $1\leq i\leq n+1$. By Proposition~\ref{P:essen.disj.fin},
\begin{equation*}
F\left(\sum_{l\in B}\mathcal C_{k_l}+\mathcal C_{k_{n+1}}\right)=\sum_{l\in B}F(\mathcal C_{k_l})+F(\mathcal C_{k_{n+1}})\supseteq F(\mathcal C_{k_{n+1}})\cap \widetilde{V}_i\neq\emptyset.
\end{equation*}
Consequently, there is $\gamma_B^i\in\Gamma$ such that 
\begin{equation*}
T_{\gamma_B^i}z\in U_{n+1}\hspace{3mm}\text{and}\hspace{3mm}\left(\sum_{l\in B}\mathcal C_{k_l}+\mathcal C_{k_{n+1}}\right)(\gamma_B^i,z)\in\widetilde{V}_i.
\end{equation*}
Let $F_{n+1}$ be the union of $F_n$ with all the elements $\gamma_B^i$ ($B\subseteq\{1,\dots,n\}$, $1\leq i\leq n+1$). Then the following conditions hold:
\begin{itemize}
\item given $1\leq i\leq n+1$ and $B\subseteq\{1,\dots,n\}$, there is $\gamma\in F_{n+1}$ with
\begin{equation*}
T_{\gamma}z\in U_{n+1}\hspace{3mm}\text{and}\hspace{3mm}\left(\sum_{l\in B}\mathcal C_{k_l}+\mathcal C_{k_{n+1}}\right)(\gamma,z)\in\widetilde{V}_i,
\end{equation*}
\item for every $\gamma\in F_n$, $\mathcal C_{k_{n+1}}(\gamma,z)\in W_{n+1}$.
\end{itemize}
This finishes the induction step. To summarize, we obtained an increasing sequence $(k_n)_{n\in\mathbb N}$ of positive integers and a non-decreasing sequence $(F_n)_{n\in\mathbb N}$ of finite subsets of $\Gamma$ such that conditions (1) and (2) hold for every $m\in\mathbb N$.

We show that the sequence $(k_n)_{n\in\mathbb N}$ satisfies the required condition. To this end, let $(l_n)_{n\in\mathbb N}$ be a subsequence of $(k_n)_{n\in\mathbb N}$ and write $l_n=k_{t(n)}$ for $n\in\mathbb N$. The  convergence of the series $\sum_{n=1}^{\infty}\mathcal C_{l_n}$ u.c.s. in $\Cocc $ is guaranteed by one of our earlier arguments in this proof; denote the sum of this series by $\mathcal C$. We show that $F(\mathcal C)=\mathbb T^1$. To this end, fix a neighbourhood $U$ of $z$ in $X$ and an arbitrary nonempty open subset $V$ of $G$; we need to find $\gamma\in\Gamma$ so that $T_{\gamma}z\in U$ and $\mathcal C(\gamma,z)\in V$. Choose $n\in\mathbb N$ with $U_n\subseteq U$ and $\overline{V_n}\subseteq V$. Clearly, $t(n)\geq n$. Therefore, by definition of the set $F_{t(n)}$, there is $\gamma\in F_{t(n)}$ such that 
\begin{equation*}
T_{\gamma}z\in U_{t(n)}\hspace{3mm}\text{and}\hspace{3mm}\left(\sum_{i=1}^n\mathcal C_{k_{t(i)}}\right)(\gamma,z)\in\widetilde{V}_n.
\end{equation*}
Moreover, since $\gamma\in F_{t(n)}$, we have $\mathcal C_{k_{t(i)}}(\gamma,z)\in W_{t(i)}$ for every $i>n$. Thus, for every $p\in\mathbb N$,
\begin{equation*}
\begin{split}
\left(\sum_{i=1}^{n+p}\mathcal C_{l_i}\right)(\gamma,z)&=\left(\sum_{i=1}^{n+p}\mathcal C_{k_{t(i)}}\right)(\gamma,z)=\left(\sum_{i=1}^n\mathcal C_{k_{t(i)}}\right)(\gamma,z)+\sum_{i=n+1}^{n+p}\mathcal C_{k_{t(i)}}(\gamma,z)\\
&\in\widetilde{V}_n+W_{t(n+1)}+W_{t(n+2)}\dots+W_{t(n+p)}\\
&\subseteq\widetilde{V}_n+W_{t(n)+1}+W_{t(n)+2}\dots+W_{t(n)+p}\\
&\subseteq\widetilde{V}_n+W_{t(n)}\subseteq \widetilde{V}_n+W_n\subseteq V_n.
\end{split}
\end{equation*}
Consequently, $\mathcal C(\gamma,z)=(\sum_{n=1}^{\infty}\mathcal C_{l_n})(\gamma,z)\in\overline{V_n}\subseteq V$. Since, in addition, $T_{\gamma}z\in U_{t(n)}\subseteq U_n\subseteq U$, the proof is finished.
\end{proof}

\begin{lemma}\label{L:not.forget.past}
Let $\Flow$ be a minimal flow with $\Gamma$ and $X$ locally compact second countable. Assume that $\mathcal A\subseteq\Cocc $ is a finite set and for $n\in\mathbb N$ let $\xi_n\colon X\to\mathbb T^1$ be a continuous map. If $\co(\xi_n)\stackrel{ucs}{\longrightarrow}1$ then there exists an increasing sequence $(k_n)_{n\in\mathbb N}$ of positive integers such that for every subsequence $(l_n)_{n\in\mathbb N}$ of $(k_n)_{n\in\mathbb N}$, the series $\sum_{n=1}^{\infty}\co(\xi_{l_n})$ converges u.c.s. in $\Cocc $ and
\begin{equation}\label{Eq:not.frgt.past}
F\left(\mathcal C+\sum_{n=1}^{\infty}\co(\xi_{l_n})\right)\supseteq F(\mathcal C)
\end{equation}
for every $\mathcal C\in\mathcal A$.
\end{lemma}
\begin{proof}
By Lemma~\ref{L:convergent.subseries} we may assume that the series $\sum_{n=1}^{\infty}\co(\xi_{k_n})$ converges u.c.s. in $\Cocc $ for every increasing sequence of positive integers $(k_n)_{n\in\mathbb N}$. Moreover, it follows from the nature of the statement of the lemma that we may restrict to the case when $\mathcal A$ is a singleton, $\mathcal A=\{\mathcal C\}$.

Fix a base point $z\in X$ along with the following objects:
\begin{itemize}
\item a sequence $(V_n)_{n\in\mathbb N}$ containing those and only those elements of a fixed countable basis of $\mathbb T^1$, which intersect the set $F(\mathcal C)$,
\item a local base $(U_n)_{n\in\mathbb N}$ at $z$ in $X$ with $U_{n+1}\subseteq U_n$ for every $n\in\mathbb N$,
\item a local base $(W_n)_{n\in\mathbb N}$ at $1$ in $\mathbb T^1$ and a sequence $(\widetilde{V}_n)_{n\in\mathbb N}$ of open subsets of $\mathbb T^1$ such that $F(\mathcal C)\cap\widetilde{V}_n\neq\emptyset$, $W_{n+1}+W_{n+1}\subseteq W_n$ and $\widetilde{V}_n+\sum_{i=1}^nW_n\subseteq V_n$ for every $n\in\mathbb N$.
\end{itemize}
Proceeding by induction on $n$, we shall construct three sequences:
\begin{itemize}
\item a sequence $(\widetilde{U}_n)_{n\in\mathbb N}$ of neighbourhoods of $z$ in $X$,
\item a sequence $(\gamma_n)_{n\in\mathbb N}$ of elements of $\Gamma$,
\item an increasing sequence $(k_n)_{n\in\mathbb N}$ of positive integers,
\end{itemize}
such that the following conditions hold for every $n\in\mathbb N$:
\begin{enumerate}
\item $\widetilde{U}_{n+1}\subseteq\widetilde{U}_n\subseteq U_n$,
\item $T_{\gamma_n}z\in\widetilde{U}_n$ and $\mathcal C(\gamma_n,z)\in\widetilde{V}_n$,
\item $\xi_{k_l}(\widetilde{U}_n)-\xi_{k_l}(\widetilde{U}_n)\subseteq W_n$ for $1\leq l\leq n-1$,
\item $\co(\xi_{k_n})(\gamma_l,z)=\xi_{k_n}(T_{\gamma_l}z)-\xi_{k_n}(z)\in W_{n+1}$ for $1\leq l\leq n$.
\end{enumerate}

First we define $\widetilde{U}_1$, $\gamma_1$ and $k_1$. Set $\widetilde{U}_1=U_1$. Since $\widetilde{V}_1$ intersects $F(\mathcal C)$, there is $\gamma_1\in\Gamma$ with $T_{\gamma_1}z\in\widetilde{U}_1$ and $\mathcal C(\gamma_1,z)\in\widetilde{V}_1$. Finally, by our assumptions on the coboundaries $\co(\xi_i)$, there is $k_1\in\mathbb N$ such that $\co(\xi_{k_1})(\gamma_1,z)=\xi_{k_1}(T_{\gamma_1}z)-\xi_{k_1}(z)\in W_2$.

Now assume that $\widetilde{U}_l$, $\gamma_l$ and $k_l$ are given for $l=1,\dots,n$; we shall define $\widetilde{U}_{n+1}$, $\gamma_{n+1}$ and $k_{n+1}$. First, use continuity of the maps $\xi_{k_1},\dots,\xi_{k_n}$ to find a neighbourhood $\widetilde{U}_{n+1}\subseteq\widetilde{U}_n\cap U_{n+1}$ of $z$ in $X$ so that $\xi_{k_l}(\widetilde{U}_{n+1})-\xi_{k_l}(\widetilde{U}_{n+1})\subseteq W_{n+1}$ for every $1\leq l\leq n$. Further, since $F(\mathcal C)\cap\widetilde{V}_{n+1}\neq\emptyset$, there is $\gamma_{n+1}\in\Gamma$ with $T_{\gamma_{n+1}}z\in\widetilde{U}_{n+1}$ and $\mathcal C(\gamma_{n+1},z)\in\widetilde{V}_{n+1}$. Finally, by our assumptions on the coboundaries $\co(\xi_i)$, there is $k_{n+1}>k_n$ such that $\co(\xi_{k_{n+1}})(\gamma_l,z)=\xi_{k_{n+1}}(T_{\gamma_l}z)-\xi_{k_{n+1}}(z)\in W_{n+2}$ for $1\leq l\leq n+1$. This finishes the induction step of our construction.

We finish the proof by showing that the sequence $(k_n)_{n\in\mathbb N}$ constructed above satisfies the conclusion of the lemma. First, $(k_n)_{n\in\mathbb N}$ is increasing by definition. Let $(l_n)_{n\in\mathbb N}$ be a subsequence of $(k_n)_{n\in\mathbb N}$ and write $l_n=k_{t(n)}$ for $n\in\mathbb N$. The series $\sum_{n=1}^{\infty}\co(\xi_{l_n})$ converges u.c.s. in $\Cocc $ by our argument from the beginning of the proof; denote the sum of this series by $\mathcal D$. To verify (\ref{Eq:not.frgt.past}), fix a neighbourhood $U$ of $z$ in $X$ and an open set $V\subseteq\mathbb T^1$ intersecting $F(\mathcal C)$; we need to find $\gamma\in\Gamma$ with $T_{\gamma}z\in U$ and $(\mathcal C+\mathcal D)(\gamma,z)\in V$. Choose $n\in\mathbb N$ with $\widetilde{U}_n\subseteq U$ and $\overline{V_n}\subseteq V$. If $p\geq n$ then, by virtue of (1)--(4) and by our choice of the sets $V_i,W_i,\widetilde{V}_i$,
\begin{equation*}
\begin{split}
&\left(\mathcal C+\sum_{m=1}^p\co(\xi_{l_m})\right)(\gamma_n,z)=\mathcal C(\gamma_n,z)+\sum_{m=1}^p\co\left(\xi_{k_{t(m)}}\right)(\gamma_n,z)\\
&=\mathcal C(\gamma_n,z)+\sum_{\substack{1\leq m\leq p \\
t(m)\leq n-1}}\left(\xi_{k_{t(m)}}(T_{\gamma_n}z)-\xi_{k_{t(m)}}(z)\right)+\sum_{\substack{1\leq m\leq p \\
n\leq t(m)}}\left(\xi_{k_{t(m)}}(T_{\gamma_n}z)-\xi_{k_{t(m)}}(z)\right)\\
&\in \widetilde{V}_n+\sum_{\substack{1\leq m\leq p \\ t(m)\leq n-1}} W_n+\sum_{\substack{1\leq m\leq p \\
n\leq t(m)}}W_{t(m)+1}\subseteq\widetilde{V}_n+\sum_{i=1}^{n-1}W_n+\sum_{i=n+1}^{t(p)+1}W_i\\
&\subseteq\widetilde{V}_n+\sum_{i=1}^{n-1}W_n+W_n=\widetilde{V}_n+\sum_{i=1}^nW_n\subseteq V_n.
\end{split}
\end{equation*}
Consequently,
\begin{equation*}
(\mathcal C+\mathcal D)(\gamma_n,z)=\left(\mathcal C+\sum_{m=1}^{\infty}\co(\xi_{l_m})\right)(\gamma_n,z)\in\overline{V_n}\subseteq V.
\end{equation*}
Since also $T_{\gamma_n}z\in\widetilde{U}_n\subseteq U$, the proof of the lemma is finished.
\end{proof}

\begin{proof}[Proof of Theorem~\ref{T:first.ineq.yes}]
Before turning to the proof, notice the following facts. First, since the groups $F(\mathcal C_n)$ ($n\in\mathbb N$) are pair-wise group-disjoint, it follows from Remark~\ref{R:grp.disj.char} that
\begin{equation*}
F(\chi_1\mathcal C_n)=\chi_1F(\mathcal C_n)\perp\chi_2F(\mathcal C_m)=F(\chi_2\mathcal C_m)
\end{equation*}
for all $\chi_1,\chi_2\in G^*$ and $n\neq m$. Further, by virtue of Lemma~\ref{L:grp.disj.char} and Remark~\ref{R:grp.disj.char}, we may assume that all the groups $F(\mathcal C_n)$ ($n\in\mathbb N$) are totally disconnected and hence that all the groups $F(\chi\mathcal C_n)=\chi F(\mathcal C_n)$ ($\chi\in G^*$, $n\in\mathbb N$) are finite. Finally, by Lemma~\ref{L:convergent.subseries}, we may assume that the series $\sum_{n=1}^{\infty}\mathcal C_{k_n}$ converges u.c.s. in $\Cocc(G)$ (and hence the series $\sum_{n=1}^{\infty}\chi\mathcal C_{k_n}$ converges u.c.s. in $\Cocc $ for every $\chi\in G^*$) for every increasing sequence of positive integers $(k_n)_{n\in\mathbb N}$. We divide the proof into two steps.

\emph{1st step.} Since the group $G$ is second countable, its dual group $G^*$ is countable. Write $G^*=\{\chi_m : m\in\mathbb N\}$. We show that there exists a double sequence of positive integers $(k_{m,n})_{m,n\in\mathbb N}$ such that the following conditions hold for every $m\in\mathbb N$:
\begin{enumerate}
\item[(1)] the sequence $(k_{m,n})_{n\in\mathbb N}$ is increasing and contains $(k_{m+1,n})_{n\in\mathbb N}$ as its subsequence,
\item[(2)] either all the groups $F(\chi_m\mathcal C_{k_{m,n}})$ ($n\in\mathbb N$) are trivial or none of them is trivial,
\item[(3)] if $F(\chi_m\mathcal C_{k_{m,n}})\neq1$ for every $n\in\mathbb N$ then
\begin{equation*}
F\left(\sum_{n=1}^{\infty}\chi_m\mathcal C_{l_n}\right)=\mathbb T^1=\overline{\sum_{n=1}^{\infty}F(\chi_m\mathcal C_{l_n})}
\end{equation*}
for every subsequence $(l_n)_{n\in\mathbb N}$ of the sequence $(k_{m,n})_{n\in\mathbb N}$,
\item[(4)] if $F(\chi_m\mathcal C_{k_{m,n}})=1$ for every $n\in\mathbb N$ then
\begin{equation*}
F\left(\sum_{l\in B}\chi_m\mathcal C_l+\sum_{n=1}^{\infty}\chi_m\mathcal C_{l_n}\right)\supseteq\sum_{l\in B}F(\chi_m\mathcal C_l)
\end{equation*}
for every subsequence $(l_n)_{n\in\mathbb N}$ of the sequence $(k_{m,n})_{n\geq m+1}$ and every subset $B$ of the set $\{1,2,\dots,k_{m,m}-1\}$.
\end{enumerate}

We shall construct the required double sequence $(k_{m,n})_{m,n\in\mathbb N}$ by induction on $m$. So fix $m\in\mathbb N$ and assume that $(k_{j,n})_{n\in\mathbb N}$ ($j=1,\dots,m-1$) are sequences of positive integers satisfying conditions expressed in (1)--(4). In order to construct $(k_{m,n})_{n\in\mathbb N}$, we distinguish two cases.

First assume that $F(\chi_m\mathcal C_{k_{m-1,n}})\neq1$ for infinitely many indexes $n\in\mathbb N$. Let $(t_n)_{n\in\mathbb N}$ be a subsequence of $(k_{m-1,n})_{n\in\mathbb N}$ with $F(\chi_m\mathcal C_{t_n})\neq1$ for every $n\in\mathbb N$. Since the groups $F(\chi_m\mathcal C_{t_n})$ ($n\in\mathbb N$) are pair-wise group-disjoint, Lemma~\ref{L:transit.coc.circ} yields a subsequence $(k_{m,n})_{n\in\mathbb N}$ of $(t_n)_{n\in\mathbb N}$ such that
\begin{equation*}
F\left(\sum_{n=1}^{\infty}\chi_m\mathcal C_{l_n}\right)=\mathbb T^1=\overline{\sum_{n=1}^{\infty}F(\chi_m\mathcal C_{l_n})}
\end{equation*}
for every subsequence $(l_n)_{n\in\mathbb N}$ of $(k_{m,n})_{n\in\mathbb N}$. It follows directly from our definition of $(k_{m,n})_{n\in\mathbb N}$ that conditions (1)--(4) hold true.

Now assume that $F(\chi_m\mathcal C_{k_{m-1,n}})=1$ (that is, $\chi_m\mathcal C_{k_{m-1,n}}\in\Cob$) for all but finitely many $n\in\mathbb N$. Let $(t_n)_{n\in\mathbb N}$ be a subsequence of $(k_{m-1,n})_{n\in\mathbb N}$ with $\chi_m\mathcal C_{t_n}\in\Cob$ for every $n\in\mathbb N$. Set $k_{m,i}=t_i$ for $i=1,\dots,m$ and let $\mathcal A$  be the set of all extensions of the form $\sum_{l\in B}\chi_m\mathcal C_l$, where $B\subseteq\{1,2,\dots,k_{m,m}-1\}$. By virtue of Lemma~\ref{L:not.forget.past}, there is a subsequence $(k_{m,n})_{n\geq m+1}$ of $(t_n)_{n\geq m+1}$ such that
\begin{equation*}
F\left(\mathcal C+\sum_{n=1}^{\infty}\chi_m\mathcal C_{l_n}\right)\supseteq F(\mathcal C)
\end{equation*}
for every subsequence $(l_n)_{n\in\mathbb N}$ of $(k_{m,n})_{n\geq m+1}$ and for every $\mathcal C\in\mathcal A$. It follows directly from our definition of $(k_{m,n})_{n\in\mathbb N}$ that conditions (1)--(3) hold true. Condition (4) follows from the definition of $(k_{m,n})_{n\geq m+1}$ and from the equality $F(\sum_{l\in B}\chi_m\mathcal C_l)=\sum_{l\in B}F(\chi_m\mathcal C_l)$ for every subset $B$ of $\{1,2,\dots,k_{m,m}-1\}$, the latter being secured by Proposition~\ref{P:essen.disj.fin}.

\emph{2nd step.} For $n\in\mathbb N$ set $k_n=k_{n,n}$. Clearly, $(k_n)_{n\in\mathbb N}$ is an increasing sequence of positive integers. We show that it satisfies the condition expressed in Theorem~\ref{T:first.ineq.yes}.

The convergence of the series $\sum_{n=1}^{\infty}\mathcal C_{k_n}$ u.c.s. in $\Cocc(G)$ has already been justified at the beginning of the proof. In order to verify (\ref{Eq:frst.inq.main}), fix a character $\chi$ of $G$ and choose $m\in\mathbb N$ with $\chi=\chi_m$. Then $(k_{n,n})_{n=m}^{\infty}$ is a subsequence of $(k_{m,n})_{n\in\mathbb N}$ and hence, by condition (2), there are two possibilities:
\begin{enumerate}
\item[(i)] $F(\chi_m\mathcal C_{k_{n,n}})=1$ for every $n\geq m$,
\item[(ii)] $F(\chi_m\mathcal C_{k_{n,n}})\neq1$ for every $n\geq m$.
\end{enumerate}
We wish to show that in both these cases,
\begin{equation}\label{Eq:frst.inq.main.char}
F\left(\sum_{n=1}^{\infty}\chi_m\mathcal C_{k_n}\right)\supseteq\overline{\sum_{n=1}^{\infty}F(\chi_m\mathcal C_{k_n})}.
\end{equation}

First consider case (i). Set $B=\{k_{n,n} : 1\leq n\leq m-1\}$ and $(l_n)_{n\in\mathbb N}=(k_{n,n})_{n\geq m+1}$. Clearly, $B$ is a subset of $\{1,2,\dots,k_{m,m}-1\}$. Moreover, since $(l_n)_{n\in\mathbb N}$ is a subsequence of $(k_{m+1,n})_{n\geq m+1}$ by definition and since $(k_{m+1,n})_{n\geq m+1}$ is a subsequence of $(k_{m,n})_{n\geq m+1}$, it follows that $(l_n)_{n\in\mathbb N}$ is a subsequence of $(k_{m,n})_{n\geq m+1}$. Thus, by applying condition (4) and the inclusion $\chi_m\mathcal C_{k_{n,n}}\in\Cob$ for every $n\geq m$, we obtain the desired inclusion (\ref{Eq:frst.inq.main.char}):
\begin{equation*}
\begin{split}
F\left(\sum_{n=1}^{\infty}\chi_m\mathcal C_{k_n}\right)&=F\left(\sum_{n=1}^{\infty}\chi_m\mathcal C_{k_{n,n}}\right)\\
&=F\left(\sum_{n=1}^{m-1}\chi_m\mathcal C_{k_{n,n}}+\chi_m\mathcal C_{k_{m,m}}+\sum_{n=m+1}^{\infty}\chi_m\mathcal C_{k_{n,n}}\right)\\
&=F\left(\sum_{l\in B}\chi_m\mathcal C_l+\sum_{n=1}^{\infty}\chi_m\mathcal C_{l_n}\right)\supseteq\sum_{l\in B}F(\chi_m\mathcal C_l)\\
&=\overline{\sum_{n=1}^{m-1}F(\chi_m\mathcal C_{k_{n,n}})}=\overline{\sum_{n=1}^{\infty}F(\chi_m\mathcal C_{k_{n,n}})}=\overline{\sum_{n=1}^{\infty}F(\chi_m\mathcal C_{k_n})}.
\end{split}
\end{equation*}

Now consider case (ii). Set $(l_n)_{n\in\mathbb N}=(k_{n,n})_{n=m}^{\infty}$. Since $(l_n)_{n\in\mathbb N}$ is a subsequence of $(k_{m,n})_{n\in\mathbb N}$, condition (3) yields
\begin{equation*}
F\left(\sum_{n=m}^{\infty}\chi_m\mathcal C_{k_{n,n}}\right)=F\left(\sum_{n=1}^{\infty}\chi_m\mathcal C_{l_n}\right)=\mathbb T^1=\overline{\sum_{n=1}^{\infty}F(\chi_m\mathcal C_{l_n})}=\overline{\sum_{n=m}^{\infty}F(\chi_m\mathcal C_{k_{n,n}})}.
\end{equation*}
Further, since the groups $F(\chi_m\mathcal C_{k_{n,n}})$ ($1\leq n\leq m-1$) are pair-wise group-disjoint and finite, it follows from Proposition~\ref{P:essen.disj.fin} that the group
\begin{equation*}
F\left(\sum_{n=1}^{m-1}\chi_m\mathcal C_{k_{n,n}}\right)=\sum_{n=1}^{m-1}F(\chi_m\mathcal C_{k_{n,n}})
\end{equation*}
is also finite. Thus, by virtue of Corollary~\ref{C:essen.disj.fin},
\begin{equation*}
\begin{split}
F\left(\sum_{n=1}^{\infty}\chi_m\mathcal C_{k_n}\right)&=F\left(\sum_{n=1}^{\infty}\chi_m\mathcal C_{k_{n,n}}\right)=F\left(\sum_{n=1}^{m-1}\chi_m\mathcal C_{k_{n,n}}+\sum_{n=m}^{\infty}\chi_m\mathcal C_{k_{n,n}}\right)\\
&=F\left(\sum_{n=1}^{m-1}\chi_m\mathcal C_{k_{n,n}}\right)+F\left(\sum_{n=m}^{\infty}\chi_m\mathcal C_{k_{n,n}}\right)\\
&=\mathbb T^1\supseteq\overline{\sum_{n=1}^{\infty}F(\chi_m\mathcal C_{k_{n,n}})}=\overline{\sum_{n=1}^{\infty}F(\chi_m\mathcal C_{k_n})}.
\end{split}
\end{equation*}
This verifies the desired inclusion (\ref{Eq:frst.inq.main.char}) also in case (ii).

To summarize, for every $\chi\in G^*$, inclusion (\ref{Eq:frst.inq.main.char}) yields
\begin{equation*}
\begin{split}
\chi F\left(\sum_{n=1}^{\infty}\mathcal C_{k_n}\right)&=F\left(\chi\sum_{n=1}^{\infty}\mathcal C_{k_n}\right)=F\left(\sum_{n=1}^{\infty}\chi\mathcal C_{k_n}\right)
\supseteq\overline{\sum_{n=1}^{\infty}F(\chi\mathcal C_{k_n})}\\
&=\overline{\sum_{n=1}^{\infty}\chi F(\mathcal C_{k_n})}=\overline{\bigcup_{n\in\mathbb N}\sum_{l=1}^n\chi F(\mathcal C_{k_l})}=\chi\left(\overline{\bigcup_{n\in\mathbb N}\sum_{l=1}^n F(\mathcal C_{k_l})}\right)\\
&=\chi\left(\overline{\sum_{n=1}^{\infty}F(\mathcal C_{k_n})}\right).
\end{split}
\end{equation*}
Thus, since $F\left(\sum_{n=1}^{\infty}\mathcal C_{k_n}\right)$ is a closed subgroup of $G$, inclusion (\ref{Eq:frst.inq.main}) from the theorem follows at last.
\end{proof}

\begin{proposition}\label{P:first.ineq.yes}
Let $\Gamma\in\mathsf{LCGp}\setminus\mathsf{CGp}$ be amenable, $X$ be a compact space and $G\in\mathsf{CAGp}$ be connected. Assume that $\Flow$ is a minimal flow with a free point. If all $\Gamma$, $X$ and $G$ are second countable then there exists a sequence $(\mathcal C_n)_{n\in\mathbb N}$ in $\Cocc(G)$ such that the groups $F(\mathcal C_n)$ are pair-wise group-disjoint, the series $\sum_{n=1}^{\infty}\mathcal C_n$ converges u.c.s. in $\Cocc(G)$ and
\begin{equation}\label{Eq:F.may.impld}
F\left(\sum_{n=1}^{\infty}\mathcal C_n\right)=e\subseteq G=
\overline{\sum_{n=1}^{\infty}F\left(\mathcal C_n\right)}.
\end{equation}
\end{proposition}
\begin{proof}
Denote by $\Cobcl(G)$ the closure of $\Cob(G)$ in $\Cocc(G)$. By our assumptions, $G$ contains a dense one-parameter subgroup and so it has a dense identity arc-component $G_a$. Consequently, by virtue of \cite[Theorem~6]{Dir3}, there is $\mathcal C\in\Cobcl(G)$ such that the induced skew product $\mathcal F_{\mathcal C}$ on $X\times G$ is topologically transitive (that is, possesses a dense orbit). By our definition of the functor $F$, $F(\mathcal C)=G$. Now fix a sequence $(\xi_n)_{n\in\mathbb N}$ of continuous maps $X\to G$ such that $\co(\xi_n)\to\mathcal C$ u.c.s. in $\Cocc(G)$. Set $\mathcal C_1=\mathcal C$, $\mathcal C_2=-\co(\xi_1)=\co\left(-\xi_1\right)$ and $\mathcal C_n=\co(\xi_{n-2})-\co(\xi_{n-1})=\co\left(\xi_{n-2}-\xi_{n-1}\right)$ for $n\geq3$. Then $\sum_{i=1}^n\mathcal C_i=\mathcal C-\co(\xi_{n-1})$ for every $n\geq2$ and so the series $\sum_{n=1}^{\infty}\mathcal C_n$ converges u.c.s. to the identity $e$ of $\Cocc(G)$. The two equalities from (\ref{Eq:F.may.impld}) now hold obviously. Finally, since $F(\mathcal C_1)=G$ and $F(\mathcal C_n)=e$ for every $n\geq 2$, the groups $F(\mathcal C_n)$ ($n\in\mathbb N$) are pair-wise group-disjoint.
\end{proof}

\begin{theorem}\label{T:second.ineq.not}
Let $\Gamma\in\mathsf{LCGp}\setminus\mathsf{CGp}$ be amenable, $X$ be a compact space and $G\in\mathsf{CAGp}$ be connected. Assume that $\Flow$ is a minimal flow with a free point $z$. If all $\Gamma$, $X$ and $G$ are second countable then there exists a sequence $(\xi_n)_{n\in\mathbb N}$ of continuous maps $X\to G$ such that for every increasing sequence of positive integers $(k_n)_{n\in\mathbb N}$, the series $\sum_{n=1}^{\infty}\co(\xi_{k_n})$ converges u.c.s. in $\Cocc(G)$ and
\begin{equation*}
F\left(\sum_{n=1}^{\infty}\co(\xi_{k_n})\right)=G\supseteq e=
\overline{\sum_{n=1}^{\infty}F\left(\co(\xi_{k_n})\right)}.
\end{equation*}
\end{theorem}

Before turning to the proof of Theorem~\ref{T:second.ineq.not}, we prove the following auxiliary lemma.

\begin{lemma}\label{L:approx.by.cob}
Under the assumptions of Theorem~\ref{T:second.ineq.not}, the following statement holds. Given a neighbourhood $U$ of $z$ in $X$, nonempty open sets $V_1,\dots,V_n\subseteq G$, a compact set $K\subseteq \Gamma$ and a neighbourhood $W$ of $e$ in $G$, there exist a continuous map $\xi\colon X\to G$ and points $\gamma^1,\dots,\gamma^n\in\Gamma$ such that
\begin{enumerate}
\item[(i)] $T_{\gamma^j}z\in U$ and $\xi(T_{\gamma^j}z)-\xi(z)\in V_j$ for every $j=1,\dots,n$ and
\item[(ii)] $\xi(T_{\gamma}x)-\xi(x)\in W$ for all $\gamma\in K$ and $x\in X$.
\end{enumerate}
\end{lemma}
\begin{proof}
Let $\Cobcl(G)$ be the closure of $\Cob(G)$ in $\Cocc(G)$. By our assumptions, $G$ has a dense identity arc-component $G_a$. Consequently, by virtue of \cite[Theorem~6]{Dir3}, the extensions $\mathcal C\in\Cobcl(G)$ with $F(\mathcal C)=G$ form a residual subset of $\Cobcl(G)$ (see also the proof of Proposition~\ref{P:first.ineq.yes}). In particular, such an extension $\mathcal C$ exists within $[K\times X;W]\cap\Cobcl(G)$. Thus, there are $\mathcal C\in\Cobcl(G)$ and $\gamma^1,\dots,\gamma^n\in\Gamma$ such that
\begin{itemize}
\item $T_{\gamma^j}z\in U$ and $\mathcal C(\gamma^j,z)\in V_j$ for every $j=1,\dots,n$ and
\item $\mathcal C(\gamma,x)\in W$ for all $\gamma\in K$ and $x\in X$.
\end{itemize}
The two conditions above are clearly open with respect to $\mathcal C\in\Cobcl(G)$ and hence they are satisfied also by an appropriate coboundary. That is, there is a continuous map $\xi\colon X\to G$ such that
\begin{itemize}
\item $T_{\gamma^j}z\in U$ and $\co(\xi)(\gamma^j,z)\in V_j$ for every $j=1,\dots,n$ and
\item $\co(\xi)(\gamma,x)\in W$ for all $\gamma\in K$ and $x\in X$.
\end{itemize}
This verifies both conditions (i) and (ii).
\end{proof}

\begin{proof}[Proof of Theorem~\ref{T:second.ineq.not}]
We start by fixing the following objects:
\begin{itemize}
\item a local base $(U_n)_{n\in\mathbb N}$ at $z$ in $X$ with $U_{n+1}\subseteq U_n$ for every $n\in\mathbb N$,
\item a sequence $(K_n)_{n\in\mathbb N}$ of compact subsets of $\Gamma$, whose interiors cover $\Gamma$ and which satisfy $K_n\subseteq K_{n+1}$ for every $n\in\mathbb N$,
\item a basis $(V_n)_{n\in\mathbb N}$ for the topology of $G$,
\item a local base $(W_n)_{n\in\mathbb N}$ at $e$ in $G$ and a sequence $(\widetilde{V}_n)_{n\in\mathbb N}$ of nonempty open subsets of $G$ with $W_{n+1}+W_{n+1}\subseteq W_n$ and $\widetilde{V}_i+W_n+W_n\subseteq V_i$ for $1\leq i\leq n<\infty$.
\end{itemize}
We claim that there exist
\begin{itemize}
\item a sequence $(\xi_n)_{n\in\mathbb N}$ of continuous maps $X\to G$,
\item a sequence $(\widetilde{U}_n)_{n\in\mathbb N}$ of neighbourhoods of $z$ in $X$ and
\item a family $\gamma_n^j$ ($1\leq j\leq n<\infty$) of elements of $\Gamma$,
\end{itemize}
such that the following conditions hold:
\begin{enumerate}
\item[(a)] $\widetilde{U}_n\subseteq U_n$ for every $n\in\mathbb N$,
\item[(b)] $\xi_i(x)-\xi_i(z)\in W_{2n-1}$ for $x\in\widetilde{U}_n$ and $1\leq i\leq n-1$,
\item[(c)] $T_{\gamma_n^j}z\in\widetilde{U}_n$ and $\xi_n(T_{\gamma_n^j}z)-\xi_n(z)\in\widetilde{V}_j$ for $1\leq j\leq n<\infty$,
\item[(d)] $\xi_n(T_{\gamma}x)-\xi_n(x)\in W_n$ for $\gamma\in K_n\cup\{\gamma_k^j : 1\leq j\leq k\leq n-1\}$ and $x\in X$.
\end{enumerate}

The construction of these objects follows easily by induction on $n$. Indeed, assume that the maps $\xi_1,\dots,\xi_{n-1}$, the sets $\widetilde{U}_1,\dots,\widetilde{U}_{n-1}$ and the points $\gamma_k^j$ ($1\leq j\leq k\leq n-1$) are given. Fix a neighbourhood $\widetilde{U}_n$ of $z$ in $X$ small enough so that $\widetilde{U}_n\subseteq U_n$ and $\xi_i(x)-\xi_i(z)\in W_{2n-1}$ for all $x\in\widetilde{U}_n$ and $i=1,\dots,n-1$. Finally, use Lemma~\ref{L:approx.by.cob} to find a continuous map $\xi_n\colon X\to G$ and points $\gamma_n^1,\dots,\gamma_n^n\in\Gamma$ such that
\begin{itemize}
\item $T_{\gamma_n^j}z\in\widetilde{U}_n$ and $\xi_n(T_{\gamma_n^j}z)-\xi_n(z)\in\widetilde{V}_j$ for $j=1,\dots,n$ and
\item $\xi_n(T_{\gamma}x)-\xi_n(x)\in W_n$ for all $\gamma\in K_n\cup\{\gamma_k^j : 1\leq j\leq k\leq n-1\}$ and $x\in X$.
\end{itemize}
Then, clearly, all the conditions (a)--(d) hold true.

We show that the sequence $(\xi_n)_{n\in\mathbb N}$ satisfies the conclusion of the theorem. To this end, fix an increasing sequence of positive integers $(k_n)_{n\in\mathbb N}$. First, we verify the u.c.s. convergence of $\sum_{n=1}^{\infty}\co(\xi_{k_n})$ in $\Cocc(G)$. So let $K\subseteq\Gamma$ be a compact set and $W$ be a neighbourhood of $e$ in $G$. Choose $n_0\in\mathbb N$ with $K\subseteq K_{n_0}$ and $W_{n_0}\subseteq W$. Let $\gamma\in K$, $x\in X$, $n\geq n_0$ and $p\in \mathbb N$. Since $K_{k_l}\supseteq K_l\supseteq K_{n_0}\supseteq K$ for every $l\geq n_0$, it follows from condition (d) above that
\begin{equation*}
\begin{split}
\sum_{l=1}^{n+p}\co(\xi_{k_l})(\gamma,x)-\sum_{l=1}^n\co(\xi_{k_l})(\gamma,x)&=\sum_{l=n+1}^{n+p}\co(\xi_{k_l})(\gamma,x)\\
&=\sum_{l=n+1}^{n+p}\left(\xi_{k_l}(T_{\gamma}x)-\xi_{k_l}(x)\right)\\
&\in\sum_{l=n+1}^{n+p}W_{k_l}\subseteq \sum_{l=n+1}^{n+p}W_l\subseteq W_n\subseteq W_{n_0}\subseteq W.
\end{split}
\end{equation*}
By completeness of $\Cocc(G)$, the series $\sum_{n=1}^{\infty}\co(\xi_{k_n})$ converges u.c.s. in $\Cocc(G)$. Denote the sum of this series by $\mathcal C$.

We show that $\mathcal C(\gamma_{k_n}^n,z)\in \overline{V_n}$ for every $n\in\mathbb N$. So fix $n\in\mathbb N$ and let $p\in\mathbb N$. We shall make use of the following observations.
\begin{itemize}
\item Given $1\leq l\leq n-1$, we have $T_{\gamma_{k_n}^n}z\in\widetilde{U}_{k_n}$ by virtue of (c). Hence $\xi_{k_l}(T_{\gamma_{k_n}^n}z)-\xi_{k_l}(z)\in W_{2k_n-1}\subseteq W_{2n-1}$ by virtue of (b).
\item We have $\xi_{k_n}(T_{\gamma_{k_n}^n}z)-\xi_{k_n}(z)\in\widetilde{V}_n$ by virtue of (c).
\item Given $n+1\leq l\leq n+p$, we have $\xi_{k_l}(T_{\gamma_{k_n}^n}z)-\xi_{k_l}(z)\in W_{k_l}\subseteq W_l$ by virtue of (d).
\end{itemize}
These observations yield
\begin{equation*}
\begin{split}
\sum_{l=1}^{n+p}\co(\xi_{k_l})(\gamma_{k_n}^n,z)&=
\sum_{l=1}^{n-1}\left(\xi_{k_l}(T_{\gamma_{k_n}^n}z)-\xi_{k_l}(z)\right)+\left(\xi_{k_n}(T_{\gamma_{k_n}^n}z)-\xi_{k_n}(z)\right)+\\
&+\sum_{l=n+1}^{n+p}\left(\xi_{k_l}(T_{\gamma_{k_n}^n}z)-\xi_{k_l}(z)\right)\\
&\in\sum_{l=1}^{n-1}W_{2n-1}+\widetilde{V}_n+\sum_{l=n+1}^{n+p}W_l
\subseteq\sum_{l=1}^{n-1}W_{n+l}+\widetilde{V}_n+\sum_{l=n+1}^{n+p}W_l\\
&\subseteq W_n+\widetilde{V}_n+W_n\subseteq V_n.
\end{split}
\end{equation*}
Since these inclusions hold for every $p\in\mathbb N$, it follows that, indeed, $\mathcal C(\gamma_{k_n}^n,z)\in\overline{V_n}$.

We finish the proof of the theorem by showing that $F(\mathcal C)=G$. To this end, fix a neighbourhood $U$ of $z$ in $X$ and a nonempty open set $V\subseteq G$; we need to find $\gamma\in\Gamma$ with $T_{\gamma}z\in U$ and $\mathcal C(\gamma,z)\in V$. Choose $n\in\mathbb N$ with $U_n\subseteq U$ and $\overline{V_n}\subseteq V$, and set $\gamma=\gamma_{k_n}^n$. Then $T_{\gamma}z\in\widetilde{U}_{k_n}\subseteq U_{k_n}\subseteq U_n\subseteq U$ by virtue of (c) and (a), and $\mathcal C(\gamma,z)\in\overline{V_n}\subseteq V$ by our argument above.
\end{proof}

\section{The ext-topology}\label{S:ext-top.intr}

In the previous section we had an opportunity to recall \cite[Theorem~6]{Dir3}, which states that under mild assumptions on a minimal flow $\mathcal F$ and a group $G\in\mathsf{CAGp}$, a generic extension from the closure of $\Cob(G)$ in $\Cocc(G)$ is minimal; in particular, $\Cob(G)$ is not a closed subgroup of $\Cocc(G)$ in the topology of u.c.s. convergence. Consequently, the quotient topology on $\Coch(G)$ is not Hausdorff and so $\Coch(G)$ with this topology is not a topological group in the sense used in this work.

Our aim in this section is to introduce a natural (Hausdorff) group topology $\tau_{ext}$ on $\Coch(G)$. We are motivated by our results from the preceding sections, where we discussed two continuity properties of the functor $F$. Being inspired by these considerations, we are searching for a group topology on $\Coch(G)$, with respect to which the map $F\colon\Coch(G)\to 2^G$ would be continuous. Our main result in this section, namely Theorem~\ref{P:group.top.cont}, shows that the coarsest topology with these properties does exist and is induced by a natural $2^G$-valued translation-invariant metric on $\Coch(G)$. We call this topology the \emph{ext-topology} and denote it by $\tau_{ext}$.

Before turning to our results recall that for $G\in\mathsf{CAGp}$, the hyper-semigroup $2^G$ over $G$ carries the Vietoris topology. Moreover, by Lemma~\ref{L:Viet.top.hyp}, for every $K\in 2^G$, the sets
\begin{equation*}
\mathcal U_V(K)=\left\{H\in 2^G : H\subseteq VK\hspace{2mm}\text{and}\hspace{2mm}K\subseteq VH\right\}
\end{equation*}
form a local base at $K$ in $2^G$, when $V$ runs through the neighbourhoods of $e$ in $G$. Consequently, given a net $(H_i)_{i\in I}$ of closed subgroups of $G$, the following conditions are equivalent:
\begin{itemize}
\item $H_i\to e$ in $2^G$,
\item for every neighbourhood $V$ of $e$ in $G$, there is $i_0\in I$ such that $H_i\subseteq V$ for every $i\geq i_0$.
\end{itemize}

\begin{lemma}\label{L:recall.to.metric}
Let $\Flow$ be a minimal flow and $G\in\mathsf{CAGp}$. If $\mathcal C,\mathcal D\in\Cocc(G)$ then the following statements hold:
\begin{enumerate}
\item[(a)] $F(\mathcal C)=e$ if and only if $\mathcal C\in\Cob(G)$,
\item[(b)] $F(\mathcal C^{-1})=F(\mathcal C)^{-1}=F(\mathcal C)$,
\item[(c)] $F(\mathcal C\mathcal D)\subseteq F(\mathcal C)F(\mathcal D)$.
\end{enumerate}
\end{lemma}
\begin{proof}
Statement (a) follows from Theorem~\ref{T:functor.E.def}(2). Since the inversion in $G$ is an automorphism of $G$ and $F(\mathcal C)$ is a subgroup of $G$, statement (b) follows from Theorem~\ref{T:functor.E.def}(3). To verify statement (c), let $\varphi\colon G\times G\to G$ be the product in $G$ and for $i=1,2$, let $\pr_i\colon G\times G\to G$ be the projection on the $i$-th coordinate. By Theorem~\ref{T:functor.E.def}(3), $\pr_1 F(\mathcal C,\mathcal D)=F(\pr_1(\mathcal C,\mathcal D))=F(\mathcal C)$ and, similarly, $\pr_2 F(\mathcal C,\mathcal D)=F(\pr_2(\mathcal C,\mathcal D))=F(\mathcal D)$. Consequently, $F(\mathcal C,\mathcal D)\subseteq F(\mathcal C)\times F(\mathcal D)$. By applying Theorem~\ref{T:functor.E.def}(3) once more, we obtain
\begin{equation*}
F(\mathcal C\mathcal D)=F(\varphi(\mathcal C,\mathcal D))=\varphi(F(\mathcal C,\mathcal D))\subseteq\varphi(F(\mathcal C)\times F(\mathcal D))=F(\mathcal C)F(\mathcal D),
\end{equation*}
as was to be shown.
\end{proof}

\begin{proposition}\label{P:intro.metric.e}
Let $\Flow$ be a minimal flow and $G\in\mathsf{CAGp}$. For $\mathcal C,\mathcal D\in\Cocc(G)$ set $\delta(\mathcal C,\mathcal D)\index[symbol]{$\delta(\mathcal C,\mathcal D)$}=F(\mathcal C\mathcal D^{-1})$. If $\mathcal B,\mathcal C,\mathcal D\in\Cocc(G)$ then the following statements hold:
\begin{enumerate}
\item[(i)] $\delta(\mathcal C,\mathcal D)=e$ if and only if $\mathcal C\simeq \mathcal D$,
\item[(ii)] $\delta(\mathcal C,\mathcal D)=\delta(\mathcal D,\mathcal C)$,
\item[(iii)] $\delta(\mathcal C,\mathcal D)\subseteq \delta(\mathcal C,\mathcal B)\delta(\mathcal B,\mathcal D)$,
\item[(iv)] $\delta(\mathcal C\mathcal B,\mathcal D\mathcal B)=\delta(\mathcal C,\mathcal D)$.
\end{enumerate}
\end{proposition}
\begin{remark}\label{R:intro.metric.e}
The proposition states that $\delta$ is a translation-invariant $2^G$-va\-lued pseu\-do-metric on $\Cocc(G)$, which identifies precisely the pairs of equivalent extensions. It follows immediately from the definition of $\delta$ and from conditions (i)--(iv) above that $\delta$ gives rise to a translation-invariant $2^G$-valued metric on $\Coch(G)=\Cocc(G)/\Cob(G)=\Cocc(G)/\simeq$.
\end{remark}
\begin{proof}[Proof of Proposition~\ref{P:intro.metric.e}]
Statement (i) follows from statement (a) in Lemma~\ref{L:recall.to.metric}. Statement (ii) is a direct consequence of statement (b) from Lemma~\ref{L:recall.to.metric}. Statement (iii) follows from statement (c) of Lemma~\ref{L:recall.to.metric}:
\begin{equation*}
\begin{split}
\delta(\mathcal C,\mathcal D)&=F(\mathcal C\mathcal D^{-1})=F\left(\left(\mathcal C\mathcal B^{-1}\right)\left(\mathcal B\mathcal D^{-1}\right)\right)\subseteq F\left(\mathcal C\mathcal B^{-1}\right)F\left(\mathcal B\mathcal D^{-1}\right)\\
&=\delta(\mathcal C,\mathcal B)\delta(\mathcal B,\mathcal D).
\end{split}
\end{equation*}
Finally, statement (iv) follows at once from the definition of $\delta$.
\end{proof}

Now let $\Flow$ be a minimal flow and $G\in\mathsf{CAGp}$. We wish to introduce a topology on $\Coch(G)$ by using the metric $\delta$ defined above. Given a neighbourhood $V$ of $e$ in $G$ and an extension $\mathcal C\in\Coch(G)$, let $B(\mathcal C,V)$ be the set of all $\mathcal D\in\Coch(G)$ such that $\delta(\mathcal D,\mathcal C)\subseteq V$. As follows from Theorem~\ref{P:group.top.cont} below, the sets $B(\mathcal C,V)$ form a basis for a group topology on $\Coch(G)$. We shall call it the \emph{ext-topology}\index{ext-topology} and denote it by $\tau_{ext}$.\index[symbol]{$\tau_{ext}$}

\begin{theorem}\label{P:group.top.cont}
Let $\Flow$ be a minimal flow and $G\in\mathsf{CAGp}$. Then the collection of sets $B(\mathcal C,V)$, where $\mathcal C\in\Coch(G)$ and $e\in V\subseteq G$ is open, forms a basis for a topology $\tau_{ext}$ on $\Coch(G)$. Moreover, $\tau_{ext}$ is the coarsest of all the topologies on $\Coch(G)$, with respect to which $\Coch(G)$ is a (Hausdorff) topological group and the mapping $F\colon\Coch(G)\to 2^G$ is continuous.
\end{theorem}
\begin{remark}\label{R:group.top.cont}
The following facts will follow from our arguments in the proof of the theorem.
\begin{itemize}
\item The topology $\tau_{ext}$ is in fact the coarsest of all the (Hausdorff) topologies on $\Coch(G)$, with respect to which the mapping $F\colon\Coch(G)\to 2^G$ is continuous and $\Coch(G)$ is a semi-topological semigroup (that is, the product in $\Coch(G)$ is separately continuous).
\item For a fixed $\mathcal C\in\Coch(G)$, the sets $B(\mathcal C,V)$ form a local base at $\mathcal C$ in $\Coch(G)$ with $\tau_{ext}$, when $V$ runs through the neighbourhoods of $e$ in $G$.
\end{itemize}
The following observations will also be useful to us.
\begin{itemize}
\item A net $(\mathcal C_i)$ in $\Coch(G)$ converges to $\mathcal C\in\Coch(G)$ with respect to $\tau_{ext}$ if and only if the net $(F(\mathcal C_i\mathcal C^{-1}))$ converges to $e$ in $2^G$; this follows from the previous statement.
\item If $G$ is a Lie group (that is, if $G_0$ is a finite-dimensional torus and is open in $G$) then $G$ has no small subgroups and the topology $\tau_{ext}$ on $\Coch(G)$ is therefore discrete. Thus, in particular, the group $\Coch $ is always discrete.
\end{itemize}
\end{remark}
\begin{proof}[Proof of Theorem~\ref{P:group.top.cont}]
First, we verify that the sets $B(\mathcal C,V)$ form a basis for a Hausdorff topology on $\Coch(G)$. Their union is clearly the whole set $\Coch(G)$, since $B(\mathcal C,V)$ contains $\mathcal C$ for all $\mathcal C$ and $V$.

Now let $\mathcal C,\mathcal D,\mathcal E\in \Coch(G)$ and $V,W$ be neighbourhoods of $e$ in $G$ with $\mathcal E\in B(\mathcal C,V)\cap B(\mathcal D,W)$; we claim that $B(\mathcal E,N)\subseteq B(\mathcal C,V)\cap B(\mathcal D,W)$ for an appropriate neighbourhood $N$ of $e$ in $G$. (In particular, this will verify the second and hence also the third statement from Remark~\ref{R:group.top.cont}.) By definition, $\delta(\mathcal E,\mathcal C)$ and $\delta(\mathcal E,\mathcal D)$ are closed subgroups of $G$ and they are contained in $V$ and $W$, respectively. Consequently, there is a neighbourhood $N$ of $e$ in $G$ with $N \delta(\mathcal E,\mathcal C)\subseteq V$ and $N \delta(\mathcal E,\mathcal D)\subseteq W$. We claim that $B(\mathcal E,N)\subseteq B(\mathcal C,V)\cap B(\mathcal D,W)$. Indeed, if $\mathcal B\in B(\mathcal E,N)$ then, by Proposition~\ref{P:intro.metric.e}(iii), $\delta(\mathcal B,\mathcal C)\subseteq \delta(\mathcal B,\mathcal E)\delta(\mathcal E,\mathcal C)\subseteq N \delta(\mathcal E,\mathcal C)\subseteq V$ and hence $\mathcal B\in B(\mathcal C,V)$. By a similar argument, $B(\mathcal E,N)\subseteq B(\mathcal D,W)$.

We show that $\tau_{ext}$ satisfies the Hausdorff condition. If $\mathcal C,\mathcal D\in\Coch(G)$ are distinct then $\delta(\mathcal C,\mathcal D)\neq e$ by Proposition~\ref{P:intro.metric.e}(i). It follows that $\delta(\mathcal C,\mathcal D)\not\subseteq V^2$ for an appropriate neighbourhood $V$ of $e$ in $G$. We claim that $B(\mathcal C,V)\cap B(\mathcal D,V)=\emptyset$. Indeed, if $\mathcal E\in B(\mathcal C,V)\cap B(\mathcal D,V)$ then, by Proposition~\ref{P:intro.metric.e}(iii, ii),
\begin{equation*}
\delta(\mathcal C,\mathcal D)\subseteq\delta(\mathcal C,\mathcal E)\delta(\mathcal E,\mathcal D)=\delta(\mathcal E,\mathcal C)\delta(\mathcal E,\mathcal D)\subseteq VV=V^2,
\end{equation*}
a contradiction.

Now we verify the second statement of the theorem. We start by showing that $\Coch(G)$ is a topological group when equipped with the topology $\tau_{ext}$ and that the mapping $F$ is continuous. First, we verify the continuity of the product in $\Coch(G)$. To this end, fix $\mathcal C,\mathcal D\in\Coch(G)$ and a neighbourhood $V$ of $e$ in $G$; we shall find a neighbourhood $W$ of $e$ in $G$ with $B(\mathcal C,W)B(\mathcal D,W)\subseteq B(\mathcal C\mathcal D,V)$. Choose a neighbourhood $W$ of $e$ in $G$ with $W^2\subseteq V$. Then for all $\mathcal C'\in B(\mathcal C,W)$ and $\mathcal D'\in B(\mathcal D,W)$,
\begin{equation*}
\delta(\mathcal C'\mathcal D',\mathcal C\mathcal D)\subseteq\delta(\mathcal C'\mathcal D',\mathcal C'\mathcal D)\delta(\mathcal C'\mathcal D,\mathcal C\mathcal D)=\delta(\mathcal D',\mathcal D)\delta(\mathcal C',\mathcal C)\subseteq WW\subseteq V
\end{equation*}
by virtue of Proposition~\ref{P:intro.metric.e}. Consequently, $B(\mathcal C,W)B(\mathcal D,W)\subseteq B(\mathcal C\mathcal D,V)$, which verifies the continuity of the product in $\Coch(G)$. The continuity of the inversion in $\Coch(G)$ follows from the identity $B(\mathcal C,V)^{-1}=B(\mathcal C^{-1},V)$, which holds by Lemma~\ref{L:recall.to.metric}(b).

To show that $F$ is continuous, fix $\mathcal C\in\Coch(G)$ and a neighbourhood $V$ of $e$ in $G$. By virtue of Lemma~\ref{L:Viet.top.hyp}, it suffices to show that $F(\mathcal D)\in\mathcal U_V(F(\mathcal C))$ for every $\mathcal D\in B(\mathcal C,V)$. If $\mathcal D\in B(\mathcal C,V)$ then
\begin{equation*}
F(\mathcal D)=F((\mathcal D\mathcal C^{-1})\mathcal C)\subseteq F(\mathcal D\mathcal C^{-1})F(\mathcal C)=\delta(\mathcal D,\mathcal C)F(\mathcal C)\subseteq VF(\mathcal C)
\end{equation*}
by Lemma~\ref{L:recall.to.metric}(c), and $F(\mathcal C)\subseteq VF(\mathcal D)$ by a similar argument. Thus, $F(\mathcal D)\in\mathcal U_V(F(\mathcal C))$ for every $\mathcal D\in B(\mathcal C,V)$, as was to be shown.

Finally, let $\tau$ be a group topology on $\Coch(G)$, with respect to which the map $F$ is continuous; we show that $\tau_{ext}$ is coarser than $\tau$. To this end, fix a net $(\mathcal C_i)$ in $\Coch(G)$ and assume that $\mathcal C_i\to\mathcal C$ in $\Coch(G)$ with respect to $\tau$. Since $\tau$ is a group topology on $\Coch(G)$, it follows that $\mathcal C_i\mathcal C^{-1}\to e$ in $\Coch(G)$ with respect to $\tau$. Consequently, by continuity of $F$, $F(\mathcal C_i\mathcal C^{-1})\to F(e)=e$ in $2^G$. By the third part of Remark~\ref{R:group.top.cont}, it follows that $\mathcal C_i\to\mathcal C$ in $\Coch(G)$ with respect to $\tau_{ext}$. This shows that $\tau_{ext}$ is coarser than $\tau$, as was to be shown.
\end{proof}

\section{Lifts of extensions}\label{Sub:lift.of.ext}

Let $\Flow$ be a minimal flow. In this section we study the problem of lifting extensions of $\mathcal F$ across covering morphisms of abelian topological groups. In Theorem~\ref{T:lift.simply.con} we show that this is always possible if the acting group $\Gamma$ of $\mathcal F$ is simply connected. (On the other hand, it is no longer possible if the assumption of simple connectedness is moved from $\Gamma$ to the phase space $X$ of $\mathcal F$; this is shown in Example~\ref{E:simp.con.where}.) Then we use this result to find the structure of a real separable Fr\'echet or Banach space in the groups $\Cocc(G)$ for connected second countable groups $G\in\mathsf{CAGp}$, see Corollary~\ref{C:structure.coc.simply}. Finally, in Theorem~\ref{T:min.exist.simp} we show that for flows $\mathcal F$ with simply connected acting groups $\Gamma\in\mathsf{CLAC}$ and with compact second countable phase spaces $X$ satisfying $\pi^1(X)\neq0$, the groupoid $\Cochm$ contains an isomorphic copy of the group $\mathbb R$.

\begin{lemma}\label{P:lifting.cocycle}
Let $\mathcal F\colon\Gamma\curvearrowright X$ be a minimal flow with a connected acting group $\Gamma$,  $G,G'$ be abelian topological groups and $p\in\Hom(G',G)$ be a morphism with a totally disconnected kernel. Given $\mathcal C\in\Cocc(G)$ and a continuous map $\mathcal C'\colon\Gamma\times X\to G'$ with $p\mathcal C'=\mathcal C$, the following conditions are equivalent:
\begin{enumerate}
\item[(a)] $\mathcal C'\in\Cocc(G')$,
\item[(b)] $\mathcal C'(1,x)=e$ for every $x\in X$, 
\item[(c)] $\mathcal C'(1,z)=e$ for some $z\in X$.
\end{enumerate}
\end{lemma}
\begin{remark}\label{R:lifting.cocycle}
Our assumptions on $p$ are satisfied if $p$ is a covering morphism. Indeed, recall that an epimorphism $p\in\Hom(G',G)$ between topological groups $G',G$ is a covering map if and only if it is a quotient (equivalently, open) map with a discrete kernel $\ker(p)\subseteq G'$.
\end{remark}
\begin{proof}[Proof of Lemma~\ref{P:lifting.cocycle}]
Implication (a)$\Rightarrow$(b) follows from the cocycle identity and implication (b)$\Rightarrow$(c) is clear. We show that (c)$\Rightarrow$(b). To this end, consider the map $\varphi\colon\Gamma\to G'$, given by $\varphi(\gamma)=\mathcal C'(1,T_{\gamma}z)$. Then $p\varphi(\gamma)=p\mathcal C'(1,T_{\gamma}z)=\mathcal C(1,T_{\gamma}z)=e$ for every $\gamma\in\Gamma$ and so $\varphi$ takes its values in $\ker(p)$. Since $\varphi$ is continuous, $\Gamma$ is connected and $\ker(p)$ is totally disconnected, it follows that $\varphi$ is constant with the value $\varphi(1)=\mathcal C'(1,z)=e$. That is to say, $\mathcal C'(1,x)=e$ for every $x$ from the $\mathcal F$-orbit of $z$. Condition (b) now follows by minimality of $\mathcal F$.

We finish the proof by showing that (b)$\Rightarrow$(a). To this end, fix $x\in X$ and consider the map $\psi\colon\Gamma\times\Gamma\to G'$, given by $\psi(\alpha,\beta)=\mathcal C'(\alpha,T_{\beta}x)\mathcal C'(\beta,x)\mathcal C'(\alpha\beta,x)^{-1}$. Then $p\psi(\alpha,\beta)=\mathcal C(\alpha,T_{\beta}x)\mathcal C(\beta,x)\mathcal C(\alpha\beta,x)^{-1}=e$ for all $\alpha,\beta\in\Gamma$ and so $\psi$ takes its values in $\ker(p)$. As above, it follows that $\psi$ is constant with the value $\psi(1,1)=\mathcal C'(1,x)=e$ and so $\mathcal C'$ satisfies the cocycle identity at $x$. Since $x\in X$ was arbitrary, this verifies condition (a).
\end{proof}

\begin{theorem}\label{T:lift.simply.con}
Let $\mathcal F\colon\Gamma\curvearrowright X$ be a minimal flow with $\Gamma\in\mathsf{CLAC}$ simply connected, $G,G'$ be abelian topological groups and $p\in\Hom(G',G)$ be a covering morphism. Then for every $\mathcal C\in\Cocc(G)$ there is a unique $\mathcal C'\in\Cocc(G')$ with $p\mathcal C'=\mathcal C$.
\end{theorem}
\begin{remark}\label{R:lift.simply.con}
We wish to add the following observations.
\begin{itemize}
\item The theorem applies to the following particular situations:
\begin{enumerate}
\item[(a)] $G=G'=\mathbb T^1$ and $p=\kappa_d$ ($d\in\mathbb N$); this shows that under the assumptions of the theorem, the group $\Cocc $ is divisible,
\item[(b)] $G=\mathbb T^1$, $G'=\mathbb R$ and $p\colon\mathbb R\to\mathbb T^1$ is the usual covering morphism.
\end{enumerate}
\item If the groups $G,G'$ are compact then $p\in\Hom(G',G)$ is a covering morphism if and only if it is an epimorphism with a finite kernel. In Theorem~\ref{T:coc.divisible} we shall see that for compact groups $G,G'$, the extensions $\mathcal C\in\Cocc(G)$ lift uniquely across $p$ also when $p\colon G'\to G$ is an epimorphism with a totally disconnected kernel.
\item It is a frequent assumption in texts dealing with covering maps that the domain of a lifted map is an element of $\mathsf{CLAC}$. A special form of our lifted map $\mathcal C$ will allow us to abandon this assumption. That is, we will not have to assume that $X\in\mathsf{CLAC}$. As a matter of fact, we shall need to apply Theorem~\ref{T:lift.simply.con} several times to a situation when $X$ is a general, not necessarily locally connected continuum.
\end{itemize}
\end{remark}
\begin{proof}[Proof of Theorem~\ref{T:lift.simply.con}]
By virtue of Lemma~\ref{P:lifting.cocycle}, we need to lift $\mathcal C$ across $p$ to a continuous map $\mathcal C'\colon\Gamma\times X\to G'$ with $\mathcal C'(1,x)=e$ for every $x\in X$. We proceed in four steps.

\emph{1st step.} We define the required lift $\mathcal C'$ of $\mathcal C$.

Fix $\gamma\in\Gamma$, $x\in X$ and choose a path $f$ in $\Gamma$ from $1$ to $\gamma$. Then $\widetilde{f}\colon[0,1]\ni t\mapsto\mathcal C(f(t),x)\in G$ is a path in $G$ from $e$ to $\mathcal C(\gamma,x)$. Since $p$ is a covering morphism, $\widetilde{f}$ lifts across $p$ to a unique path $\widehat{f}$ in $G'$ starting at $e$. Since $\Gamma$ is simply connected, the endpoint of $\widehat{f}$ does not depend on $f$, only on $\gamma$ and $x$; we shall denote this endpoint by $\mathcal C'(\gamma,x)$. Clearly, the map $\mathcal C'\colon\Gamma\times X\to G'$ thus defined projects to $\mathcal C$ via $p$ and $\mathcal C'(1,x)=e$ for every $x\in X$. We thus need to show that $\mathcal C'$ is continuous.

\emph{2nd step.} We show that $\mathcal C'$ is continuous in the second variable.

Fix $\gamma\in\Gamma$ and $x\in X$; we show that $\mathcal C'$ is continuous at $(\gamma,x)$ in the second variable. So let $V'$ be a neighbourhood of $\mathcal C'(\gamma,x)$ in $G'$. Fix a path $f$ in $\Gamma$ from $1$ to $\gamma$. From our construction of $\mathcal C'$ it follows that $\widehat{f}(t)=\mathcal C'(f(t),x)$ for every $t\in[0,1]$. Thus, by continuity of $\widehat{f}$, there exist
\begin{itemize}
\item a partition $0=t_0<t_1<\dots<t_n=1$ of $[0,1]$,
\item open subsets $V_1,\dots,V_n$ of $G$, each of them evenly covered by $p$
\item open subsets $V_1',\dots,V_n'$ of $G'$, $V_i'$ being a $p$-slice over $V_i$ for every $i=1,\dots,n$, and $V_n'\subseteq V'$,
\end{itemize}
such that
\begin{enumerate}
\item[($*$)] $\mathcal C'(f([t_{i-1},t_i])\times x)\subseteq V_i'$ for $i=1,\dots,n$.
\end{enumerate}

By virtue of ($*$),
\begin{equation*}
\mathcal C(f([t_{i-1},t_i])\times x)=p\mathcal C'(f([t_{i-1},t_i])\times x)\subseteq p(V_i')=V_i
\end{equation*}
for $i=1,\dots,n$ and $\mathcal C(f(t_i),x)=p\mathcal C'(f(t_i),x)\in p(V_i'\cap V_{i+1}')$ for $i=1,\dots,n-1$. Consequently, by continuity of $\mathcal C$, there is a neighbourhood $U$ of $x$ in $X$ such that
\begin{enumerate}
\item[(a)] $\mathcal C(f([t_{i-1},t_i])\times U)\subseteq V_i$ for $i=1,\dots,n$, and
\item[(b)] $\mathcal C(f(t_i)\times U)\subseteq p(V_i'\cap V_{i+1}')$ for $i=1,\dots,n-1$.
\end{enumerate}

Now let $y\in U$; we show that $\mathcal C'(\gamma,y)\in V'$. Consider the path $\widetilde{g}\colon[0,1]\ni t\mapsto\mathcal C(f(t),y)\in G$ in $G$ from $e$ to $\mathcal C(\gamma,y)$. For $i=1,\dots,n$, let $\widetilde{g}_i$ be the restriction of $\widetilde{g}$ onto $[t_{i-1},t_i]$, $p_i\colon V_i'\to V_i$ be the homeomorphism obtained by restricting $p$ and $\widehat{g}_i=p_i^{-1}\widetilde{g}_i$. By virtue of (a), all the maps $\widehat{g}_i$ are well defined and continuous. By virtue of (b),
\begin{equation*}
\widehat{g}_i(t_i)=p_i^{-1}\widetilde{g}_i(t_i)=p_i^{-1}\mathcal C(f(t_i),y)=p_{i+1}^{-1}\mathcal C(f(t_i),y)=p_{i+1}^{-1}\widetilde{g}_{i+1}(t_i)=\widehat{g}_{i+1}(t_i)
\end{equation*}
for every $i=1,\dots,n-1$. Consequently, the maps $\widehat{g}_i$ are restrictions of a single path $\widehat{g}\colon[0,1]\to G'$ in $G'$, which projects to $\widetilde{g}$ via $p$. By our definition of $\mathcal C'$, it follows that
\begin{equation*}
\mathcal C'(\gamma,y)=\widehat{g}(1)=\widehat{g}_n(1)=p_n^{-1}\widetilde{g}_n(1)\in V_n'\subseteq V',
\end{equation*}
as was to be shown. This shows that $\mathcal C'$ is continuous in the second variable.

\emph{3rd step.} We show that $\mathcal C'$ is (jointly) continuous.

Let $\gamma$, $x$ and $V'$ be as in the second step of the proof and set $V=p(V')$; we shall find a neighbourhood $W$ of $\gamma$ in $\Gamma$ and a neighbourhood $U$ of $x$ in $X$ with $\mathcal C'(W\times U)\subseteq V'$. Without loss of generality, we may assume that $V$ is evenly covered by $p$ and $V'$ is a $p$-slice over $V$. Since $p$ is an open map by the assumptions, $V$ is a neighbourhood of $\mathcal C(\gamma,x)$ in $G$. Fix an arc-wise connected neighbourhood $W$ of $\gamma$ in $\Gamma$ and a neighbourhood $U$ of $x$ in $X$ with $\mathcal C(W\times U)\subseteq V$. Since $\mathcal C'$ is continuous in the second variable by the second step of the proof, we may assume that $\mathcal C'(\gamma\times U)\subseteq V'$.

Fix $\delta\in W$ and $y\in U$. Choose a path $h$ in $\Gamma$ with
\begin{itemize}
\item $h(0)=1$, $h(1/2)=\gamma$, $h(1)=\delta$ and
\item $h(t)\in W$ for every $t\in[1/2,1]$.
\end{itemize}
Write $\widetilde{h}(t)=\mathcal C(h(t),y)$ for $t\in[0,1]$ and let $\widehat{h}$ be the lift of $\widetilde{h}$ across $p$ starting at $e$. Recall from the second step of the proof that $\widehat{h}(t)=\mathcal C'(h(t),y)$ for every $t\in[0,1]$. Denote by $p_0$ the homeomorphism $V'\to V$, obtained by restricting $p$. Since $\widehat{h}(1/2)=\mathcal C'(h(1/2),y)=\mathcal C'(\gamma,y)\in V'$ and $\widetilde{h}(t)\in V$ for every $t\in[1/2,1]$, it follows that $\widehat{h}|_{[1/2,1]}$ and $p_0^{-1}\widetilde{h}|_{[1/2,1]}$ are both lifts of $\widetilde{h}|_{[1/2,1]}$ across $p$, starting at the same point. Thus, by uniqueness of lifts of paths, $\widehat{h}(t)=p_0^{-1}\widetilde{h}(t)$ for every $t\in[1/2,1]$. In particular,
\begin{equation*}
\mathcal C'(\delta,y)=\mathcal C'(h(1),y)=\widehat{h}(1)=p_0^{-1}\widetilde{h}(1)\in p_0^{-1}(V)=V'.
\end{equation*}
Thus, $\mathcal C'(W\times U)\subseteq V'$, as was to be shown.

\emph{4th step.} We finish the proof by verifying the uniqueness part of the theorem.

So assume that $\mathcal C',\mathcal C''\in\Cocc(G')$ satisfy $p\mathcal C'=\mathcal C=p\mathcal C''$. Then $p(\mathcal C'\left(\mathcal C''\right)^{-1})=e$ and hence $\mathcal C'\left(\mathcal C''\right)^{-1}$ takes its values in the kernel $\ker(p)$ of $p$. Since the domain $\Gamma\times X$ of $\mathcal C'(\mathcal C'')^{-1}$ is connected and $\ker(p)$ is totally disconnected, it follows that $\mathcal C'(\mathcal C'')^{-1}$ is constant with the value $e$. This shows that $\mathcal C'=\mathcal C''$.
\end{proof}

By the first part of Remark~\ref{R:lift.simply.con}, if the acting group $\Gamma$ of a minimal flow $\mathcal F$ is simply connected then the group $\Cocc$ is divisible. In the following example we show that this is no longer the case if the assumption of simple connectedness is moved from $\Gamma$ to $X$. (This is contrary to the well known situation of $\Gamma=\mathbb Z$ and $X$ compact simply connected, in which case the group $\Cocc\cong C_z(X,\mathbb T^1)\cong C_z(X,\mathbb R)$ is divisible.)

\begin{example}\label{E:simp.con.where}
Consider the usual transitive action $\mathcal F\colon \text{U}(2)\curvearrowright\mathbb S^3$ (see Subsection~\ref{Sub:frns.vs.top.frns}). The acting group $\text{U}(2)$ of $\mathcal F$ is not simply connected, having an infinite cyclic fundamental group $\pi_1(\text{U}(2))$, but the phase space $\mathbb S^3$ of $\mathcal F$ is simply connected; we show that the group $\Cocc$ is not divisible. To this end, let $h=\det$ be the determinant function on $\text{U}(2)$. Clearly, $h$ is a topological morphism $\text{U}(2)\to\mathbb T^1$. Moreover, since
\begin{equation*}
h\left(
\begin{array}{ccc}
e^{i2\pi t} & 0 \\
0 & 1 
\end{array} \right)
=\det\left(
\begin{array}{ccc}
e^{i2\pi t} & 0 \\
0 & 1 
\end{array} \right)
=e^{i2\pi t}
\end{equation*}
for every $t\in[0,1]$, the morphism $h^{\sharp}\colon\pi_1(\text{U}(2))\to\pi_1(\mathbb T^1)$ induced by $h$ is an epimorphism (in fact, $h^{\sharp}$ is an isomorphism). Now let $\mathcal C=\mathcal Q_h\in\Cocc$ be the quasi-coboundary induced by $h$; we claim that $\mathcal C$ is a non-divisible element of $\Cocc$. If $z$ is a base point for $\mathbb S^3 $ then, by virtue of (\ref{Eq:tor.Lie.mnf.L1}) from Lemma~\ref{L:tor.Lie.mnfld},
\begin{equation*}
\mathcal C^{\sharp}\pi_1(\text{U}(2)\times\mathbb S^3)=\mathcal C_z^{\sharp}\pi_1(\text{U}(2))=h^{\sharp}\pi_1(\text{U}(2))=\pi_1(\mathbb T^1).
\end{equation*}
It follows that $\mathcal C$ does not lift across any of the endomorphisms $\kappa_d$ of $\mathbb T^1$ with $d\geq2$. (Indeed, recall that $\kappa_d^{\sharp}\pi_1(\mathbb T^1)=d(\pi_1(\mathbb T^1))\subsetneq\pi_1(\mathbb T^1)$ for $d\geq2$.) Thus, $\mathcal C$ is a non-divisible element of $\Cocc$, as was to be shown.
\end{example}

\begin{corollary}\label{C:structure.coc.simply}
Let $\Flow$ be a minimal flow with $\Gamma\in\mathsf{LieGp}$ simply connected and with $X$ compact second countable. Let $G\in\mathsf{CAGp}$ be a connected second countable group. Then $\Cocc(G)$ is the additive topological group of a real separable Fr\'echet space. If the group $G$ is finite-dimensional then $\Cocc(G)$ is the additive topological group of a real separable Banach space.
\end{corollary}
\begin{remark}\label{R:structure.coc.simply}
In the proof of the corollary we shall observe that $\Cocc(\mathbb R)$ is a separable Banach space and find out that it is isomorphic to $\Cocc$ as a topological group. This will yield the structure of a separable Banach space in $\Cocc$. Now assume that the space $X$ is infinite, so that $C_z(X,\mathbb R)$ is an infinite-dimensional real linear space. Since $\Cocc(\mathbb R)$ contains $\Cob(\mathbb R)\cong C_z(X,\mathbb R)$ as a linear subspace, it follows from the Mackey theorem that $\LSdim(\Cocc(\mathbb R))\geq\mathfrak{c}$. On the other hand, by separability of $\Cocc(\mathbb R)$, $\LSdim(\Cocc(\mathbb R))\leq\card(\Cocc(\mathbb R))\leq\mathfrak{c}$. Consequently,
\begin{equation*}
\LSdim(\Cocc)=\LSdim(\Cocc(\mathbb R))=\mathfrak{c}.
\end{equation*}
\end{remark}
\begin{proof}[Proof of Corollary~\ref{C:structure.coc.simply}]
Throughout the whole proof, all the groups $\Cocc(G)$ are assumed to carry the topology of u.c.s. convergence. Let $p\colon\mathbb R\to\mathbb T^1$ be the usual covering morphism and $\widehat{p}\colon\Cocc(\mathbb R)\ni\mathcal D\mapsto p\mathcal D\in\Cocc$ be the induced topological morphism. By virtue of Theorem~\ref{T:lift.simply.con}, $\widehat{p}$ is an epimorphism. Since $\Gamma\times X$ is connected and $\ker(p)\cong\mathbb Z$ is totally disconnected, $\widehat{p}$ is also a monomorphism.

We show that $\widehat{p}$ is in fact a topological isomorphism. To this end, let $(\mathcal D_i)_{i\in J}$ be a net in $\Cocc(\mathbb R)$ with $\widehat{p}(\mathcal D_i)=p\mathcal D_i\stackrel{ucs}{\longrightarrow}1$ in $\Cocc $; we show that $\mathcal D_i\stackrel{ucs}{\longrightarrow}0$ in $\Cocc(\mathbb R)$. So let $F\subseteq \Gamma\times X$ be a compact set and $0<\varepsilon<1/2$. Fix compact sets $1\in K\subseteq\Gamma$ and $C\subseteq X$ with $K$ connected and with $F\subseteq K\times C$. Since $p\mathcal D_i\stackrel{ucs}{\longrightarrow}1$ in $\Cocc$, there is $i_0\in J$ such that $p\mathcal D_i(\gamma,x)\in p(-\varepsilon,\varepsilon)$ for all $i\geq i_0$, $\gamma\in K$ and $x\in C$. Since $1\in K$ is connected and $\varepsilon<1/2$, it follows that $\mathcal D_i(K\times x)\subseteq(-\varepsilon,\varepsilon)$ for all $i\geq i_0$ and $x\in C$. Consequently, $\mathcal D_i(F)\subseteq\mathcal D_i(K\times C)\subseteq (-\varepsilon,\varepsilon)$ for every $i\geq i_0$. Hence, $\mathcal D_i\stackrel{ucs}{\longrightarrow}0$ in $\Cocc(\mathbb R)$, as was to be shown.

By our assumptions on $\Gamma$ and $X$, $\Cocc(\mathbb R)$ is the additive topological group of a real separable Banach space, see Lemma~\ref{L:Cocc.Ban.sp.Gmm} and Remark~\ref{R:Cocc.Ban.sp.Gmm} in Subsection~\ref{Sub:Ban.sp.of.ext}. By the first part of the proof, it follows that so is the group $\Cocc$. Moreover, the groups $\Cocc(G)$ and $\Hom(G^*,\Cocc)$ are topologically isomorphic due to Theorem~\ref{T:structure.coc} for every $G\in\mathsf{CAGp}$. Now, if $G\in\mathsf{CAGp}$ is a connected second countable group then its dual group $G^*$ is torsion-free countable and so $\Cocc(G)\cong\Hom(G^*,\Cocc)$ is the additive topological group of a real separable Fr\'echet space by the first part of Lemma~\ref{L:hom.Ban.Fre}. If, in addition, the group $G$ is finite-dimensional then the group $G^*$ is of a finite rank and so $\Cocc(G)$ is the additive topological group of a real separable Banach space by the second part of Lemma~\ref{L:hom.Ban.Fre}.
\end{proof}

\begin{theorem}\label{T:min.exist.simp}
Let $\Flow$ be a minimal flow with $\Gamma\in\mathsf{CLAC}$ simply connected and with $X$ compact. Then there is a short exact sequence of abelian groups
\begin{equation}\label{Eq:min.exist.simp}
0\longrightarrow\pi^1(X)\longrightarrow\Coch(\mathbb R)\longrightarrow\Coch\longrightarrow0.
\end{equation}
Further, there is a divisible torsion-free subgroup $\mathfrak{D}\cong\Coch/\tor(\Coch)$ of the groupoid $\Cochm$ with $\Coch=\tor(\Coch)\oplus\mathfrak{D}$. Finally, if $1\leq\rank(\pi^1(X))<\mathfrak{c}$ then 
\begin{equation*}
\rank(\mathfrak{D})=\QLSdim(\mathfrak{D})\geq\mathfrak{c}.
\end{equation*}
\end{theorem}
\begin{remark}\label{R:min.exist.simp}
Notice the following facts.
\begin{itemize}
\item The assumption $\rank(\pi^1(X))<\mathfrak{c}$ is satisfied if the group $\pi^1(X)$ is countable. This is the case if the space $X$ is second countable.
\item It follows from the last two statements of the theorem that the groupoid $\Cochm$ contains an isomorphic copy of the group $\mathbb Q^{(\mathfrak{c})}\cong\mathbb R$.
\end{itemize}
\end{remark}
\begin{proof}[Proof of Theorem~\ref{T:min.exist.simp}]
We verify each of the three statements of the theorem in a separate step of the proof.

\emph{1st step.} We show that there is an exact sequence of abelian groups (\ref{Eq:min.exist.simp}).

Let $p\colon\mathbb R\to\mathbb T^1$ be the usual covering morphism and $\widehat{p}\colon\Cocc(\mathbb R)\ni\mathcal D\mapsto p\mathcal D\in\Cocc$ be the induced morphism. By an argument similar to that from Corollary~\ref{C:structure.coc.simply}, $\widehat{p}$ is an isomorphism of groups. Further, fix a base point $z$ for $X$ and consider the map $\varphi\colon C_z(X,\mathbb T^1)\ni\xi\mapsto\co(\xi)\in\Cob$. Clearly, $\varphi$ is an isomorphism of groups. Finally, since the space $X$ is connected and the kernel $\ker(p)=\mathbb Z$ of $p$ is totally disconnected, it follows that $q\colon C_z(X,\mathbb R)\ni\zeta\mapsto p\zeta\in C_z(X,\mathbb T^1)$ is a monomorphism of groups. Moreover, by definition of $\pi^1(X)$, there is a direct sum $C_z(X,\mathbb T^1)=\im(q)\oplus\pi^1(X)$. An elementary argument also shows that $\varphi(\im(q))=\widehat{p}(\Cob(\mathbb R))$.

Now, the ideas formulated in the preceding paragraph yield a direct sum
\begin{equation}\label{Eq:min.ex.simp.aux}
\begin{split}
\widehat{p}^{-1}(\Cob)&=\widehat{p}^{-1}\varphi\left(C_z(X,\mathbb T^1)\right)
=\widehat{p}^{-1}\varphi\left(\im(q)\oplus\pi^1(X)\right)\\
&=\widehat{p}^{-1}\varphi(\im(q))\oplus\widehat{p}^{-1}\varphi(\pi^1(X))=\Cob(\mathbb R)\oplus\widehat{p}^{-1}\varphi(\pi^1(X)).
\end{split}
\end{equation}
Further, since $\Cob(\mathbb R)$ is a linear subspace of $\Cocc(\mathbb R)$, it follows that $\Cob(\mathbb R)$ is a direct summand in $\Cocc(\mathbb R)$ and each of its complementary summands in $\Cocc(\mathbb R)$ is isomorphic to $\Cocc(\mathbb R)/\Cob(\mathbb R)=\Coch(\mathbb R)$. Since $\Cob(\mathbb R)\cap\widehat{p}^{-1}\varphi(\pi^1(X))=0$ by virtue of (\ref{Eq:min.ex.simp.aux}), there is a complementary summand $H$ of $\Cob(\mathbb R)$ in $\Cocc(\mathbb R)$ with $\widehat{p}^{-1}\varphi(\pi^1(X))\subseteq H$. Consequently, there are isomorphisms of groups
\begin{equation*}
\begin{split}
\Coch&=\Cocc/\Cob\cong\widehat{p}^{-1}(\Cocc)/\widehat{p}^{-1}(\Cob)=\Cocc(\mathbb R)/\widehat{p}^{-1}(\Cob)\\
&=\left(\Cob(\mathbb R)\oplus H\right)/\left(\Cob(\mathbb R)\oplus\widehat{p}^{-1}\varphi(\pi^1(X))\right)\cong H/\widehat{p}^{-1}\varphi(\pi^1(X)).
\end{split}
\end{equation*}
Since $H\cong\Coch(\mathbb R)$ and $\widehat{p}^{-1}\varphi(\pi^1(X))\cong\pi^1(X)$, this verifies the existence of the desired exact sequence (\ref{Eq:min.exist.simp}).

\emph{2nd step.} We verify the second statement of the theorem.

By the first step of the proof, there is an isomorphism $\Cocc\cong\Cocc(\mathbb R)$. Since the group $\Cocc(\mathbb R)$ is divisible, it follows that so are the groups $\Cocc$, $\Coch$ and $\tor(\Coch)$. Consequently, $\tor(\Coch)$ is a direct summand in $\Coch$ and each of its complementary summands is divisible, torsion-free and isomorphic to $\Coch/\tor(\Coch)$. Fix such a complementary summand $\mathfrak{D}$; we need to show that $\mathfrak{D}\subseteq\Cochm$. This follows immediately from the identity $\Cochm\setminus1=\Coch\setminus\tor(\Coch)$, see Corollary~\ref{C:MinH.rel.torH}.

\emph{3rd step.} We finish the proof of the theorem by verifying its last statement.

So assume that $1\leq\rank(\pi^1(X))<\mathfrak{c}$ and fix a complementary summand $\mathfrak{D}$ of $\tor(\Coch)$ in $\Coch$. By the first step of the proof, we may view $\pi^1(X)$ as a subgroup of $\Coch(\mathbb R)$ in such a way that the quotient group $\Coch(\mathbb R)/\pi^1(X)$ is isomorphic to $\Coch$. Let $A$ be the divisible hull of $\pi^1(X)$ in $\Coch(\mathbb R)$ and $B$ be a complementary summand of $A$ in $\Coch(\mathbb R)$. Then $A/\pi^1(X)$ is a torsion group and $B$ is a torsion-free group. Since $\Coch\cong\Coch(\mathbb R)/\pi^1(X)\cong(A/\pi^1(X))\oplus B$, it follows that $A/\pi^1(X)\cong\tor(\Coch)$ and $B\cong\Coch/\tor(\Coch)\cong\mathfrak{D}$. We claim that the following statements hold:
\begin{itemize}
\item $\rank(\Coch(\mathbb R))\geq\mathfrak{c}$; indeed, $\Coch(\mathbb R)$ is a real linear space, which contains a non-trivial group $\pi^1(X)$,
\item $\rank(A)<\mathfrak{c}$; indeed, by the assumptions of the theorem,
\begin{equation*}
\rank(A)=\rank(\pi^1(X))<\mathfrak{c}.
\end{equation*}
\end{itemize}
Thus, since $\Coch(\mathbb R)=A\oplus B$ and $B\cong\mathfrak{D}$, we have $\rank(\mathfrak{D})=\rank(B)\geq\mathfrak{c}$, as was to be shown.
\end{proof}

\chapter{Topological-Algebraic Aspects}\label{S:E.and.alg-top}

\section{Free objects for categories of maps}\label{S:free.obj.maps}

Our aim in this section is to prove the existence of free objects for certain categories of continuous base point preserving maps from a given pointed space into compact abelian groups. Our results from this section will be applied in the next section, where we show the existence of a free compact abelian group extension for an arbitrary minimal flow.

Let $X$ be a pointed space with the base point $z$ and $A$ be a subgroup of the group $C_z(X,\mathbb T^1)$ of all continuous base point preserving maps $X\to\mathbb T^1$. Given $G\in\mathsf{CAGp}$, let $A_G$\index[symbol]{$A_G$} denote the set of all $f\in C_z(X,G)$ such that $\chi f\in A$ for every $\chi\in G^*$. Clearly, $A_G$ is a subgroup of $C_z(X,G)$. Moreover, if $G,H\in\mathsf{CAGp}$ and $h\in\Hom(G,H)$, then $hf\in A_H$ for every $f\in A_G$. Now consider the category $\mathsf{Maps}_A$\index[symbol]{$\mathsf{Maps}_A$} defined as follows. The objects of $\mathsf{Maps}_A$ are the elements of $A_G$ with $G\in\mathsf{CAGp}$ arbitrary. For $f_1\in A_{G_1}$ and $f_2\in A_{G_2}$, the set of morphisms $\Hom(f_1,f_2)$ in $\mathsf{Maps}_A$ consists of all $h\in\Hom(G_1,G_2)$ with $hf_1=f_2$. The composition of morphisms is obtained by restricting that from $\mathsf{CAGp}$. Clearly, $\id_G\in\Hom(f,f)$ for all $f\in A_G$ and $G\in\mathsf{CAGp}$.

Our goal now is to show that $\mathsf{Maps}_A$ is a category with a free object. This is done in Theorem~\ref{T:free.obj.MapsA} below and some examples are given in Example~\ref{E:free.obj.MapsA}. Recall that a free object for $\mathsf{Maps}_A$ consists of a group $G_A\in\mathsf{CAGp}$\index[symbol]{$G_A$} and of a map $f_A\in A_{G_A}$\index[symbol]{$f_A$} with the following universal property:
\begin{itemize}
\item for every $G\in\mathsf{CAGp}$ and every $f\in A_G$ there is a unique $h\in\Hom(G_A,G)$ with $hf_A=f$.
\end{itemize}

\begin{theorem}\label{T:free.obj.MapsA}
Let $X$ be a pointed space and $A$ be a subgroup of $C_z(X,\mathbb T^1)$ equipped with the discrete topology. Set $G_A=A^*$ and let $f_A\colon X\to G_A$ be given by the rule $f_A(x)\colon A\ni\varphi\mapsto\varphi(x)\in\mathbb T^1$ for every $x\in X$. Then $G_A$ and $f_A$ constitute a free object for the category $\mathsf{Maps}_A$.
\end{theorem}
\begin{proof}
The group $G_A$ is compact abelian, being a dual group of a discrete abelian group $A$. The map $f_A$ is well defined (that is, $f_A(x)\in A^*$ for every $x\in X$), it is base point preserving and continuous. It therefore remains to verify the universality of $f_A$. That is, we must show that for every $G\in\mathsf{CAGp}$ and every $f\in A_G$ there is a unique $h\in\Hom(G_A,G)$ with $hf_A=f$. Since the latter condition is an invariant of topological isomorphism with respect to $G$, we may assume that $G=D^*$ for an appropriate group $D\in\mathsf{DAGp}$. Thus, we need to show that $k^*f_A=f$ for one and only one $k\in\Hom(D,A)$. The uniqueness of $k$ is clear, for the equality $k^*f_A=f$ translates into $k(d)(x)=f(x)(d)$ for all $d\in D$ and $x\in X$, and this in turn is satisfied for at most one function $k\colon D\to A$. We finish the proof by showing that $k$ thus defined is a (well defined) morphism of groups.

To see that $k$ takes its values in $A$, fix $d\in D$ and consider the character $\chi$ of $D^*$ given by $\chi(\Upsilon)=\Upsilon(d)$ for every $\Upsilon\in D^*$. Then $\chi f=k(d)$. Since $f\in A_{D^*}$, it follows that $\chi f\in A$ and hence $k(d)\in A$. The fact that $k$ is a morphism of groups $D\to A$ follows at once from the fact that $f(x)$ is a morphism of groups $D\to\mathbb T^1$ for every $x\in X$.
\end{proof}

\begin{example}\label{E:free.obj.MapsA}
We present three examples of categories $\mathsf{Maps}_A$ and the ranges $G_A$ of their free objects $f_A$. We shall abuse the terminology here and refer to the groups $G_A$ as the free objects for the categories $\mathsf{Maps}_A$.
\begin{itemize}
\item[($\iota$)] Let $X$ be a pointed space and $A$ be the whole group $C_z(X,\mathbb T^1)$. Then $\mathsf{Maps}_A$ is the category of all continuous base point preserving maps $X\to G$ with $G\in\mathsf{CAGp}$ and the free object $G_A$ for $\mathsf{Maps}_A$ is the free compact abelian group over $X$. If, for instance, the space $X$ is compact connected second countable and non-degenerate then $C_z(X,\mathbb R)$ is a real linear space with dimension $\LSdim(C_z(X,\mathbb R))=\w(X)^{\aleph_0}=(\aleph_0)^{\aleph_0}=\mathfrak{c}$, and hence there are topological isomorphisms
\begin{equation*}
\begin{split}
G_A&=A^*=\left(C_z(X,\mathbb T^1)_d\right)^*\cong\left(C_z(X,\mathbb R)_d\oplus\pi^1(X)\right)^*\\
&\cong\left(\mathbb R_d^{(\mathfrak{c})}\oplus\pi^1(X)\right)^*\cong\left(\mathbb R_d\oplus\pi^1(X)\right)^*\\
&\cong b\mathbb R\times\pi^1(X)^*.
\end{split}
\end{equation*}
\item[($\iota\iota$)] Let $X$ be a non-degenerate compact connected second countable pointed space and $A$ be the subgroup of $C_z(X,\mathbb T^1)$ formed by the null-homotopic maps. Then $A$ is the identity arc-component of $C_z(X,\mathbb T^1)$ and it is naturally isomorphic to $C_z(X,\mathbb R)$. The category $\mathsf{Maps}_A$ consists of the continuous base point preserving null-homotopic maps $X\to G$ with $G\in\mathsf{CAGp}$. Similarly to ($\iota$), there are topological isomorphisms
\begin{equation*}
G_A=A^*\cong\left(C_z(X,\mathbb R)_d\right)^*\cong\left(\mathbb R_d^{(\mathfrak{c})}\right)^*\cong(\mathbb R_d)^*\cong b\mathbb R.
\end{equation*}
\item[($\iota\iota\iota$)] Let $X$ be a compact abelian group and $A=X^*$ be the dual group of $X$. Clearly, $A$ is a subgroup of $C_z(X,\mathbb T^1)$. Since the characters of compact abelian groups separate points, we have $A_G=\Hom(X,G)$ for every $G\in\mathsf{CAGp}$ and $\mathsf{Maps}_A$ is thus the category of the (topological) morphisms $X\to G$ with $G\in\mathsf{CAGp}$. The free object $G_A$ of $\mathsf{Maps}_A$ is given by the topological isomorphism $G_A=A^*=X^{**}\cong X$.
\end{itemize}
\end{example}

\begin{corollary}\label{C:MapsA.iso.Hom}
Under the assumptions of Theorem~\ref{T:free.obj.MapsA}, the category $\mathsf{Maps}_A$ is isomorphic to the category $\mathsf{Hom}(G_A,\mathsf{CAGp})$.
\end{corollary}
\begin{proof}
We shall define a covariant functor $\mathsf{F}\colon\mathsf{Maps}_A\to\mathsf{Hom}(G_A,\mathsf{CAGp})$ and show that it is an isomorphism of categories. Given $G\in\mathsf{CAGp}$ and $f\in A_G$, we let $\mathsf{F}(f)$ be the unique element of $\Hom(G_A,G)$ with $\mathsf{F}(f)f_A=f$. By universality of $f_A$, this defines a bijection between $A_G$ and $\Hom(G_A,G)$. Further, given $G_1,G_2\in\mathsf{CAGp}$, $f_1\in A_{G_1}$, $f_2\in A_{G_2}$ and $l\in\Hom(f_1,f_2)$, we set $\mathsf{F}(l)=l\in\Hom(\mathsf{F}(f_1),\mathsf{F}(f_2))$. This definition is correct, for $l\in\Hom(G_1,G_2)$ and $lF(f_1)=F(f_2)$; indeed, $(lF(f_1))f_A=lf_1=f_2=F(f_2)f_A$ and hence, by universality of $f_A$, $lF(f_1)=F(f_2)$. Moreover, we see at once that $\mathsf{F}$ defines a bijection between $\Hom(f_1,f_2)$ and $\Hom(\mathsf{F}(f_1),\mathsf{F}(f_2))$. The functorial properties $\mathsf{F}(l_2l_1)=\mathsf{F}(l_2)\mathsf{F}(l_1)$ and $\mathsf{F}(\id_f)=\id_{\mathsf{F}(f)}$ are straightforward and so $\mathsf{F}$ is indeed an isomorphism of categories.
\end{proof}

\section{Free group extensions}\label{Sub:free.ext}

Our aim in this section is to show that for every minimal flow $\mathcal F$, $\Coc$ is a category with a free object. This free object consists of a group $\Gimel\in\mathsf{CAGp}$ and of an extension $\Daleth\in\Cocc(\Gimel)$, and is universal in the sense that every object from $\Coc$ is obtained as a composition of $\Daleth$ with a unique topological morphism. The existence of a free object in $\Coc$ is interesting from the categorical point of view, for it shows that $\Coc$ can be viewed as a category of morphisms of compact abelian groups. However, from the point of view of topological dynamics this fact becomes even more important: the dynamics of every extension	 $\mathcal C\in\Coc$ is encoded in the dynamics of $\Daleth$.

After showing that the free extension $(\Gimel,\Daleth)$ of $\mathcal F$ exists, we compute $\Gimel$ in the case when $\Gamma$ is a simply connected Lie group and describe relationships of $(\Gimel,\Daleth)$ to some of the concepts introduced earlier. The existence of a free object for $\Coc$ becomes even more valuable to us in Section~\ref{Sub:free.min.ext}, where we show that a free object for the category of the minimal extensions $\mathcal C\in\Coc$ typically does not exist.

\begin{definition}\label{D:free.ext}
Let $\Flow$ be a minimal flow. The \emph{free compact abelian group extension} (briefly, the \emph{free extension})\index{free!extension of a minimal flow} of $\mathcal F$ is defined as the free object for the category $\Coc$. Recall that it consists of a group $\Gimel \in\mathsf{CAGp}$\index[symbol]{$\gimel_{\mathcal F}$} and an extension $\Daleth \in\Cocc(\Gimel )$\index[symbol]{$\daleth_{\mathcal F}$} with the following universal property:
\begin{itemize}
\item for every $G\in\mathsf{CAGp}$ and every $\mathcal C\in\Cocc(G)$ there exists a unique $h_{\mathcal C}\in\Hom(\Gimel ,G)$\index[symbol]{$h_{\mathcal C}$} with $h_{\mathcal C}\Daleth =\mathcal C$.
\end{itemize}
The free extension of $\mathcal F$ is of course unique up to a natural isomorphism. On occasions we shall take the liberty of abusing the terminology by referring either to $\Gimel $ or to $\Daleth $ as the free extension of $\mathcal F$.
\end{definition}

\begin{theorem}\label{T:exist.free.ext}
Let $\Flow$ be a minimal flow. Equip the group $\Cocc $ with the discrete topology and define $\Gimel $ as its Pontryagin dual
\begin{equation*}
\Gimel =\left(\Cocc\right)_d^*.
\end{equation*}
Further, let $\Daleth \in\Cocc(\Gimel )$ be given by
\begin{equation*}
\Daleth (\gamma,x)\colon\Cocc\ni\mathcal C\mapsto\mathcal C(\gamma,x)\in\mathbb T^1
\end{equation*}
for all $\gamma\in\Gamma$ and $x\in X$. Then $(\Gimel ,\Daleth )$ is the free extension of $\mathcal F$.
\end{theorem}
\begin{proof}
Fix $z\in X$ and set $A=\Cocc$. Then $A$ is a subgroup of $C_{(1,z)}(\Gamma\times X,\mathbb T^1)$. Moreover, under the notation from Section~\ref{S:free.obj.maps}, we have $A_G=\Cocc(G)$ for every $G\in\mathsf{CAGp}$, for the characters of compact abelian groups separate points. Consequently, $\Coc=\mathsf{Maps}_A$ and the statement of the theorem thus follows directly from Theorem~\ref{T:free.obj.MapsA}.
\end{proof}

\begin{example}\label{E:gim.simp.con}
Let $\Flow$ be a minimal flow with $\Gamma\in\mathsf{LieGp}$ simply connected and with $X$ compact second countable. Assume that both $\Gamma$ and $X$ are infinite. We show that there is a topological isomorphism $\Gimel \cong b\mathbb R$.
Indeed, by virtue of Corollary~\ref{C:structure.coc.simply} and Remark~\ref{R:structure.coc.simply}, $\Cocc$ is a real linear space with dimension $\LSdim(\Cocc)=\mathfrak{c}$. Consequently,
\begin{equation*}
\Gimel =\left(\Cocc\right)_d^*\cong\left(\mathbb R^{(\mathfrak{c})}\right)_d^*\cong\left(\mathbb R_d\right)^*=b\mathbb R.
\end{equation*}
\end{example}

\begin{theorem}\label{T:prop.of.gimel}
Let $\Flow$ be a minimal flow and $(\Gimel ,\Daleth )$ be the model of the free extension of $\mathcal F$ constructed in Theorem~\ref{T:exist.free.ext}. Then the following statements hold.
\begin{enumerate}
\item[(1)] The category $\Coc$ is isomorphic to the category $\mathsf{Hom}(\Gimel ,\mathsf{CAGp})$.
\item[(2)] Under the identification $(\Gimel )^*=(\Cocc)_d^{**}\cong\left(\Cocc\right)_d$ it holds
\begin{equation*}
F(\Daleth )^{\perp}=\Cob
\hspace{4mm}\text{and}\hspace{4mm}\left(\Daleth\right)^*=\emph{Id}_{\Cocc}.
\end{equation*}
\item[(3)] Let $\cGimel\index[symbol]{$\gimel_{\mathcal F}^\text{c}$}=\Gimel/F(\Daleth)$, $\kappa\colon\Gimel\to\cGimel$ be the quotient morphism, $\cDaleth\index[symbol]{$\daleth_{\mathcal F}^\text{c}$}=\kappa\Daleth$, $z$ be the base point of $X$ and $\vartheta$ be the base point preserving transfer function for $\cDaleth$. Then $\cGimel$ is the free compact abelian group over $X$ with respect to the universal map $\vartheta$.
\item[(4)] For every $G\in\mathsf{CAGp}$ there is an isomorphism of groups
\begin{equation*}
\widetilde{\Phi}_G\index[symbol]{$\widetilde{\Phi}_G$}\colon\Cocc(G)\ni\mathcal C\mapsto \left(\omega_G\right)^{-1}\Phi_G(\mathcal C)^*=h_{\mathcal C}\in\Hom\left(\Gimel,G\right),
\end{equation*}
where $\omega_G$ is the Pontryagin isomorphism $G\to G^{**}$, $\Phi_G(\mathcal C)^*=\mathcal C^{**}$ stands for the morphism dual to $\Phi_G(\mathcal C)=\mathcal C^*$ and the morphism $h_{\mathcal C}$ is defined by the formula $h_{\mathcal C}\Daleth=\mathcal C$.
\end{enumerate}
\end{theorem}
\begin{remark}\label{R:prop.of.gimel}
We wish to make the following observations.
\begin{itemize}
\item Under the identification from statement (2), $\Cob\subseteq\Cocc$ corresponds to $(\Cob)^{\perp\perp}\subseteq(\Cocc)_d^{**}=(\Gimel)^*$. Consequently, the isomorphism $F(\Daleth)^{\perp}\cong\Cob$ yields an equality $F(\Daleth)=(\Cob)^{\perp}$ and a topological isomorphism $F(\Daleth)\cong(\Cocc/\Cob)_d^*=(\Coch)_d^*$. Under this isomorphism the identity component $F(\Daleth)_0$ of $F(\Daleth)$ corresponds to the group $\tor(\Coch)^{\perp}$. Consequently, the group $F(\Daleth)$ is connected if and only if the group $\Coch$ is torsion-free. Similarly, $F(\Daleth)$ is torsion-free if and only if $\Coch$ is divisible.
\item Recall that $\Gimel/F(\Daleth)=\cGimel$ is topologically isomorphic to the right-hand side of (\ref{Eq:free.comp.gp}) if the base space $X$ of the flow $\mathcal F$ is compact. If $X$ is additionally assumed connected (which is the case if $\Gamma$ is connected) then the middle term of the right-hand side of (\ref{Eq:free.comp.gp}) vanishes.
\end{itemize}
\end{remark}
\begin{proof}[Proof of Theorem~\ref{T:prop.of.gimel}]
Statement (1) follows from Theorem~\ref{T:free.obj.MapsA}, Corollary~\ref{C:MapsA.iso.Hom} and Theorem~\ref{T:exist.free.ext}, for $\Coc=\mathsf{Maps}_A$ with $A=\Cocc$.

In order to verify statement (2), we adopt the following notation. Given $\mathcal C\in\Cocc $, let $\Upsilon_{\mathcal C}$ be the image of $\mathcal C$ under the Pontryagin isomorphism 
\begin{equation*}
\omega_{(\Cocc)_d}\colon(\Cocc)_d\to(\Cocc)_d^{**}=(\Gimel)^*.
\end{equation*}
Then $\Upsilon_{\mathcal C}\Daleth =\mathcal C$. From this observation and from the properties of the functor $F$ we obtain the equivalence of the following conditions for every $\mathcal C\in\Cocc$:
\begin{itemize}
\item $\Upsilon_{\mathcal C}\in F(\Daleth )^{\perp}$,
\item $\Upsilon_{\mathcal C}F(\Daleth)=1$,
\item $F(\Upsilon_{\mathcal C}\Daleth )=1$,
\item $F(\mathcal C)=1$,
\item $\mathcal C\in\Cob $.
\end{itemize}
This verifies the first statement from (2). As for the second statement from (2), for every $\mathcal C\in\Cocc $, $(\Daleth)^*(\Upsilon_{\mathcal C})=\Upsilon_{\mathcal C}\Daleth=\mathcal C$, and so the action of $(\Daleth)^*$ on $\Cocc $ is indeed identical.

We verify statement (3). First observe that $\cDaleth=\kappa\Daleth$ is a coboundary and the map $\vartheta$ is thus well defined; indeed, we have $F(\cDaleth)=F(\kappa\Daleth)=\kappa F(\Daleth)=e$ and hence $\cDaleth\in\Cob(\cGimel)$. Now let $G\in\mathsf{CAGp}$ and $\xi\colon X\to G$ be a continuous base point preserving map; we need to show that there is a unique topological morphism $\sigma\colon\cGimel\to G$ with $\sigma\vartheta=\xi$. By universality of $\Daleth$, there exists a unique $h\in\Hom(\Gimel,G)$ with $h\Daleth=\co(\xi)$. Since $e=F(\co(\xi))=F(h\Daleth)=hF(\Daleth)$, $h$ factors through $\kappa$ onto a morphism $\sigma\colon\cGimel\to G$; that is, $\sigma\kappa=h$. Thus,
\begin{equation*}
\co(\xi)=h\Daleth=\sigma\kappa\Daleth=\sigma\cDaleth
=\sigma\co(\vartheta)=\co(\sigma\vartheta)
\end{equation*}
and, since both $\xi$ and $\sigma\vartheta$ are base point preserving maps, it follows from minimality of $\mathcal F$ that $\xi=\sigma\vartheta$. To verify the uniqueness of $\sigma$ with the latter property assume that $\varrho\in\Hom(\cGimel,G)$ also satisfies $\varrho\vartheta=\xi$. Then 
\begin{equation*}
\sigma\kappa\Daleth=h\Daleth=\co(\xi)=\co(\varrho\vartheta)=\varrho\co(\vartheta)=\varrho\cDaleth
=\varrho\kappa\Daleth
\end{equation*}
and hence, by universality of $\Daleth$, $\sigma\kappa=\varrho\kappa$. Finally, since $\kappa$ is an epimorphism, it follows that $\sigma=\varrho$. This finishes the proof of statement (3).

We verify statement (4). Fix $G\in\mathsf{CAGp}$. Recall that there is an isomorphism of groups
\begin{equation*}
\Hom(G^*,\Cocc)\ni k\mapsto \left(\omega_G\right)^{-1}k^*\in\Hom(\Gimel,G).
\end{equation*}
Since $\Phi_G$ is an isomorphism $\Cocc(G)\to\Hom(G^*,\Cocc)$ by virtue of Theorem~\ref{T:structure.coc}, it follows that $\widetilde{\Phi}_G$ is an isomorphism $\Cocc(G)\to\Hom(\Gimel,G)$. Thus, to finish the proof of (4), it remains to show that $(\omega_G)^{-1}\mathcal C^{**}=h_{\mathcal C}$ or, equivalently, that $\mathcal C^{**}=\omega_Gh_{\mathcal C}$ for every $\mathcal C\in\Cocc(G)$. So fix $\mathcal C\in\Cocc(G)$. By universality of $\Daleth$ and by virtue of the identity $h_{\mathcal C}\Daleth=\mathcal C$, the desired equality $\mathcal C^{**}=\omega_Gh_{\mathcal C}$ is equivalent with $\mathcal C^{**}\Daleth=\omega_G\mathcal C$. It is the latter equality that we now verify. Observe that both $\mathcal C^{**}\Daleth$ and $\omega_G\mathcal C$ are elements of $\Cocc(G^{**})$.

Fix $\gamma\in\Gamma$, $x\in X$ and $\chi\in G^*$. Then
\begin{equation*}
\left(\mathcal C^{**}\Daleth\right)(\gamma,x)\colon\chi\mapsto\Daleth(\gamma,x)(\mathcal C^*(\chi))=\Daleth(\gamma,x)(\chi\mathcal C)=\chi\mathcal C(\gamma,x)
\end{equation*}
by definition of $\mathcal C^*$ and $\Daleth$. Likewise,
\begin{equation*}
\left(\omega_G\mathcal C\right)(\gamma,x)\colon\chi\mapsto\omega_G(\mathcal C(\gamma,x))(\chi)=\chi(\mathcal C(\gamma,x))=\chi\mathcal C(\gamma,x)
\end{equation*}
by definition of $\omega_G$. This verifies the equality $\mathcal C^{**}\Daleth=\omega_G\mathcal C$.
\end{proof}

\begin{proposition}\label{P:further.prop.of.gimel}
Let $\Flow$ be a minimal flow and $(\Gimel ,\Daleth )$ be the model of the free extension of $\mathcal F$ constructed in Theorem~\ref{T:exist.free.ext}. Then the following statements hold.
\begin{enumerate}
\item[($\iota$)] For $G\in\mathsf{CAGp}$ and $\mathcal C\in\Cocc(G)$ let $h_{\mathcal C}$ be the unique morphism $\Gimel\to G$ with $h_{\mathcal C}\Daleth=\mathcal C$ and let $H_{\mathcal C}\index[symbol]{$H_{\mathcal C}$}=\ker(h_{\mathcal C})$, $G_{\mathcal C}\index[symbol]{$G_{\mathcal C}$}=\im(h_{\mathcal C})$. Then
\begin{enumerate}
\item[(a)] $F(\mathcal C)=h_{\mathcal C}F(\Daleth)$,
\item[(b)] $(h_{\mathcal C})^*=\mathcal C^*$ (under the identification $(\Gimel )^*\cong\Cocc$),
\item[(c)] $H_{\mathcal C}^{\perp}=\im((h_{\mathcal C})^*)=\im(\mathcal C^*)$,
\item[(d)] $G_{\mathcal C}$ is the closed subgroup of $G$ generated by the values of $\mathcal C$.
\end{enumerate}
\item[($\iota\iota$)] For all $G\in\mathsf{CAGp}$ and $\mathcal C\in \Cocc(G)$, the following conditions are equivalent:
\begin{enumerate}
\item[(i)] $\mathcal C\in\Cob(G)$,
\item[(ii)] $F(\Daleth)\subseteq H_{\mathcal C}$.
\end{enumerate}
\item[($\iota\iota\iota$)] For all $G\in\mathsf{CAGp}$ and $\mathcal C\in \Cocc(G)$, the following conditions are equivalent:
\begin{enumerate}
\item[($\alpha$)] $\mathcal C$ is minimal,
\item[($\beta$)] $G_{\mathcal C}=G$ and $H_{\mathcal C}^{\perp}\cap F(\Daleth)^{\perp}=0$,
\item[($\gamma$)] $G_{\mathcal C}=G$ and $H_{\mathcal C}+F(\Daleth)=\Gimel$.
\end{enumerate}
\end{enumerate}
\end{proposition}
\begin{proof}
We begin by verifying ($\iota$). Statement (a) follows immediately from the definition of $h_{\mathcal C}$ and from the properties of the functor $F$. To verify statement (b), observe that $(\Daleth)^*(h_{\mathcal C})^*=\mathcal C^*$ by definition of $h_{\mathcal C}$. Moreover, under the identification $(\Gimel )^*\cong\Cocc$, $(\Daleth)^*$ equals the identity on $\Cocc$ by statement (2) of Theorem~\ref{T:prop.of.gimel}. Thus, it follows that $(h_{\mathcal C})^*=\mathcal C^*$. The first equality in (c) follows by a simple computation $H_{\mathcal C}^{\perp}=\ker(h_{\mathcal C})^{\perp}=\im((h_{\mathcal C})^*)$ and the second equality follows from (b). Let $V_{\mathcal C}$ be the closed subgroup of $G$ generated by the values of $\mathcal C$. Clearly, $\chi\in G^*$ belongs to $V_{\mathcal C}^{\perp}$ if and only if $\mathcal C^*(\chi)=\chi\mathcal C=1$. Thus, by virtue of (b), $V_{\mathcal C}^{\perp}=\ker(\mathcal C^*)=\ker((h_{\mathcal C})^*)=\im(h_{\mathcal C})^{\perp}=G_{\mathcal C}^{\perp}$. Since $V_{\mathcal C}$ and $G_{\mathcal C}$ are closed subgroups of $G$ by definition, it follows that $G_{\mathcal C}=V_{\mathcal C}$. This verifies statement (d).

Statement ($\iota\iota$) follows at once from part (a) of statement ($\iota$) and from Theorem~\ref{T:functor.E.def}(2).

We finish the proof of the proposition by verifying statement ($\iota\iota\iota$). By part (a) of statement ($\iota$), we have $F(\mathcal C)^{\perp}=(h_{\mathcal C}F(\Daleth))^{\perp}=((h_{\mathcal C})^*)^{-1}(F(\Daleth)^{\perp})$. Thus, by part (c) of statement ($\iota$), the following conditions are equivalent:
\begin{itemize}
\item $\mathcal C$ is minimal,
\item $F(\mathcal C)=G$,
\item $((h_{\mathcal C})^*)^{-1}(F(\Daleth)^{\perp})=0$,
\item $(h_{\mathcal C})^*$ is a monomorphism and $\im((h_{\mathcal C})^*)\cap F(\Daleth)^{\perp}=0$,
\item $h_{\mathcal C}$ is an epimorphism and $H_{\mathcal C}^{\perp}\cap F(\Daleth)^{\perp}=0$,
\item $G_{\mathcal C}=G$ and $H_{\mathcal C}^{\perp}\cap F(\Daleth)^{\perp}=0$.
\end{itemize}
This verifies the equivalence of ($\alpha$) and ($\beta$). The equivalence of ($\beta$) and ($\gamma$) follows from basic properties of the annihilator mechanism.
\end{proof}

\section{Divisibility and torsion-freeness}\label{S:div.and.tor-free}

Let $\Flow$ be a minimal flow. In this section we study two important properties of the groups $\Cocc(G)$ ($G\in\mathsf{CAGp}$), namely their torsion-freeness and divisibility. First, we relate the torsion-freeness of $\Cocc(G)$ to the topological(-algebraic) properties of $X$ and $\Gamma$. Then we relate the divisibility of the groups $\Cocc(G)$ to the divisibility of the single group $\Cocc$ and show how it is reflected in the existence of lifts of extensions across various types of epimorphisms within $\Coc$. We also gain from these results some useful information about the free extension $\Gimel$ of $\mathcal F$.

\begin{theorem}\label{T:gimel.conn}
Let $\Flow$ be a minimal flow. Then the following statements are equivalent:
\begin{enumerate}
\item[(i)] $\Gimel $ is connected,
\item[(ii)] $\Cocc $ is torsion-free,
\item[(iii)] $\Cocc(G)$ is torsion-free for every $G\in\mathsf{CAGp}$,
\item[(iv)] the space $X$ is connected and the group $\Gamma$ has no non-trivial finite abelian quotient groups.
\end{enumerate}
\end{theorem}
\begin{remark}\label{R:gimel.conn}
Let us mention the following facts.
\begin{itemize}
\item In connection with (iv) observe that a group $\Gamma$ may posses a non-trivial finite quotient group and still not possess a non-trivial finite \emph{abelian} quotient group. For instance, a finite non-abelian simple group (say, the alternating group~$A_5$) has no proper normal subgroups and hence it possesses no abelian quotient groups other than the trivial one. Of course, our interest is in infinite topological transformation groups, but our remark now applies for instance to $\Gamma=\mathbb R\times A_5$.
\item A discrete abelian group $\Gamma$ possesses no non-trivial finite abelian quotient groups if and only if it is divisible. This occurs if and only if its dual group $\Gamma^*$ is torsion-free.
\item The two conditions from (iv) are automatically satisfied if the group $\Gamma$ is connected, which is the situation of our main interest.
\end{itemize}
\end{remark}
\begin{proof}[Proof of Theorem~\ref{T:gimel.conn}]
The equivalence of (i) and (ii) is an immediate consequence of the definition of $\Gimel $. Implication (iii)$\Rightarrow$(ii) is clear. The converse implication follows from Theorem~\ref{T:structure.coc}, which yields
\begin{equation*}
\begin{split}
\tor(\Cocc(G))&\cong\tor(\Hom(G^*,\Cocc ))=\bigcup_{k\in\mathbb N}\text{tor}_k(\Hom(G^*,\Cocc))\\
&=\bigcup_{k\in\mathbb N}\Hom(G^*,\text{tor}_k(\Cocc ))=0
\end{split}
\end{equation*}
for every $G\in\mathsf{CAGp}$.

To show that (iv)$\Rightarrow$(ii), assume that $X$ is connected and $\Gamma$ has no non-trivial finite abelian quotient groups. Fix $\mathcal C\in\Cocc $ and $k\in\mathbb N$ with $\mathcal C^k=1$. Then $\mathcal C$ takes its values in $\mathbb Z_k$. By connectedness of $X$, the mapping $l\colon\Gamma\ni\gamma\mapsto\mathcal C(\gamma,x)\in\mathbb Z_k$ does not depend on the choice of $x\in X$ and from the cocycle identity for $\mathcal C$ it follows that $l$ is a morphism of topological groups. It is clearly a quotient morphism onto its image. Since the group $\Gamma$ has no non-trivial finite abelian quotient groups by our assumptions, it follows that $l=1$ and hence $\mathcal C=1$. Thus, indeed, $\Cocc$ is a torsion-free group.

We finish the proof by showing that (ii)$\Rightarrow$(iv). Assume that condition (iv) is not satisfied and consider first the case of $X$ not connected. Fix a continuous map $\xi\colon X\to\mathbb T^1$ with $\xi(X)=\{\pm1\}$. Then, by minimality of $\mathcal F$, $\co(\xi)$ is a non-trivial element of $\Cocc $ with $\co(\xi)^2=1$ and hence $\tor(\Cocc )\neq\emptyset$. Now consider the case when $\Gamma$ has a non-trivial finite abelian quotient group. Then $\Gamma$ factors onto $\mathbb Z_k$ for some $k\geq 2$ and the quasi-coboundary induced by the underlying epimorphism $\Gamma\to\mathbb Z_k$ is a non-trivial element of $\tor_k(\Cocc )$. This verifies the desired implication.
\end{proof}

\begin{theorem}\label{T:lifts.and.gimel}
Let $\Flow$ be a minimal flow. Then the following conditions are equivalent:
\begin{enumerate}
\item[(1)] $\Gimel$ is torsion-free,
\item[(2)] $\Cocc $ is divisible,
\item[(3)] $\Cocc(G)$ is divisible for every $G\in\mathsf{CAGp}$ connected,
\item[(4)] the objects of $\Coc$ lift across the finite-to-one epimorphisms from $\mathsf{CAGp}$,
\item[(5)] the objects of $\Coc$ lift across all epimorphisms from $\mathsf{CAGp}$.
\item[(6)] for every every $G\in\mathsf{CAGp}$ and every $H\sbgp G$, the canonical sequence
\begin{equation*}
0\longrightarrow H\stackrel{j}{\longrightarrow}G\stackrel{p}{\longrightarrow}G/H\longrightarrow0
\end{equation*}
gives rise to the short exact sequence of abelian groups
\begin{equation*}
0\longrightarrow\Cocc(H)\stackrel{\widehat{j}}{\longrightarrow}\Cocc(G)\stackrel{\widehat{p}}{\longrightarrow}\Cocc(G/H)\longrightarrow0.
\end{equation*}
\end{enumerate}
\end{theorem}
\begin{remark}\label{R:lifts.and.gimel}
The assumption of connectedness of $G$ in condition (3) is in general not redundant, for otherwise the morphisms $\kappa_k\colon G\to G$ ($k\in\mathbb N$) would not be surjective. However, if we restrict ourselves to the connected topological transformation groups in the base (which is the situation that interests us most) then condition (3) holds for every $G\in\mathsf{CAGp}$. Indeed, in such case all extensions $\mathcal C\in\Cocc(G)$ take their values in $G_0$ and $\kappa_k(G_0)=G_0$ for every $k\in\mathbb N$.
\end{remark}
\begin{proof}[Proof of Theorem~\ref{T:lifts.and.gimel}]
The equivalence of (1) and (2) follows immediately from the definition of $\Gimel$. Further, implications (5)$\Rightarrow$(4), (3)$\Rightarrow$(2) and (4)$\Rightarrow$(2) hold obviously. Also, (5) implies (3), since the endomorphisms $\kappa_k$ of $G$ ($k\in\mathbb N$) are epimorphisms for every connected group $G\in\mathsf{CAGp}$. Thus, in order to prove the equivalence of (1)--(5), it suffices now to show that (2) implies (5). So assume that $\Cocc $ is divisible. Fix $G,H\in\mathsf{CAGp}$, an epimorphism $p\colon H\to G$ and $\mathcal C\in\Cocc(G)$. Then $\mathcal C^*\in\Hom(G^*,\Cocc )$ and $p^*\colon G^*\to H^*$ is a monomorphism. By divisibility of $\Cocc $, $\mathcal C^*$ extends through $p^*$ to a morphism $h\colon H^*\to\Cocc $; that is, $hp^*=\mathcal C^*$. By virtue of Theorem~\ref{T:structure.coc}, $h=\mathcal D^*$ for an appropriate $\mathcal D\in\Cocc(H)$. Then $\mathcal D^*p^*=\mathcal C^*$, which translates into $p\mathcal D=\mathcal C$. This verifies condition (5).

We finish the proof by verifying the equivalence (2)$\Leftrightarrow$(6). First, assume that the group $\Cocc$ is divisible and fix $G,H$ as in (6). Then $\Ext(H^*,\Cocc)=0$ and statement (6) thus follows from the exactness of (\ref{Eq:sb.qt.gp.EH}) in Corollary~\ref{C:ex.seq.qt.gp}. Conversely, assume that condition (6) holds. We show that the group $\Cocc$ is divisible by verifying the identity $\Ext(A,\Cocc)=0$ for every $A\in\mathsf{DAGp}$. So fix $A\in\mathsf{DAGp}$ and set $H=A^*$. Clearly, $H\sbgp G$ for an appropriate torus $G$. Since the group $G^*$ is free abelian, it follows that $\Ext(G^*,\Cocc)=0$ and so the exact sequence (\ref{Eq:sb.qt.gp.EH}) takes the form
\begin{equation}\label{Eq:ex.ext.sq.Ext0}
0\longrightarrow\Cocc(H)\stackrel{\widehat{j}}{\longrightarrow}\Cocc(G)\stackrel{\widehat{p}}{\longrightarrow}\Cocc(G/H)\stackrel{\epsilon}{\longrightarrow}\Ext(H^*,\Cocc)\longrightarrow0.
\end{equation}
By virtue of condition (6), $\widehat{p}$ is an epimorphism and hence, by exactness of the sequence (\ref{Eq:ex.ext.sq.Ext0}), $\Ext(H^*,\Cocc)=0$. Thus, $\Ext(A,\Cocc)=0$, as was to be shown.
\end{proof}

\begin{theorem}\label{T:coc.divisible}
Let $\Flow$ be a minimal flow. Then the following conditions are equivalent:
\begin{enumerate}
\item[(a)] $\Gimel$ is torsion-free and connected,
\item[(b)] $\Cocc $ is divisible and torsion-free,
\item[(c)] $\Cocc(G)$ is divisible for every $G\in\mathsf{CAGp}$,
\item[(d)] $\Cocc(G)$ is divisible and torsion-free for every $G\in\mathsf{CAGp}$,
\item[(e)] the objects of $\Coc$ lift uniquely across the epimorphisms from $\mathsf{CAGp}$ with finite kernels,
\item[(f)] the objects of $\Coc$ lift uniquely across the epimorphisms from $\mathsf{CAGp}$ with totally disconnected kernels.
\end{enumerate}
\end{theorem}
\begin{remark}\label{R:coc.divisible}
In connection with condition (b), recall the following facts. If $\Flow$ is a minimal flow whose acting group $\Gamma$ has no non-trivial finite abelian quotient groups and whose base space $X$ is connected then the group $\Cocc $ is torsion-free, see Theorem~\ref{T:gimel.conn}. If $\Gamma\in\mathsf{CLAC}$ is simply connected and $X$ is compact then $\Cocc $ is divisible, see Remark~\ref{R:lift.simply.con}.
\end{remark}
\begin{proof}[Proof of Theorem~\ref{T:coc.divisible}]
First, the equivalence of (a) and (b) follows immediately from the definition of $\Gimel$. Further, implications (f)$\Rightarrow$(e) and (d)$\Rightarrow$(c) are obvious. Also, (b) follows from (e) at once, for the endomorphisms $\kappa_k$ ($k\in\mathbb N$) of $\mathbb T^1$ have finite kernels. To finish the proof, it is sufficient to verify implications (b)$\Rightarrow$(d), (c)$\Rightarrow$(b) and (b)$\Rightarrow$(f).

We show that (d) follows from (b). So assume that (b) holds and fix $G\in\mathsf{CAGp}$. Since the group $\Cocc$ is torsion-free, it follows that so is $\Hom(G^*,\Cocc)$. Further, since $\Cocc$ is torsion-free and divisible, all the morphisms $\kappa_k\colon\Cocc\to\Cocc$ ($k\in\mathbb N$) are isomorphisms and hence the group $\Hom(G^*,\Cocc)$ is divisible. Thus, $\Hom(G^*,\Cocc)$ is both divisible and torsion-free and so, by Theorem~\ref{T:structure.coc}, it follows that so is $\Cocc(G)$. This verifies condition (d).

We show that (b) follows from (c). So assume that all the groups $\Cocc(G)$ ($G\in\mathsf{CAGp}$) are divisible; we show that the group $\Cocc$ is divisible and torsion-free. The divisibility of $\Cocc$ is immediate. To see that $\Cocc$ is also torsion-free, fix an integer $k\geq2$. Then, by virtue of Theorem~\ref{T:structure.coc}, $\Cocc(\mathbb Z_k)\cong\Hom(\mathbb Z_k,\Cocc)\cong\tor_k(\Cocc)$. Since $\Cocc(\mathbb Z_k)$ is divisible by the assumptions, it follows that $\tor_k(\Cocc)=0$. This verifies condition (b).

We finish the proof of the theorem by showing that (f) follows from (b). So assume that $\Cocc$ is both divisible and torsion-free. Fix an epimorphism $p\colon H\to G$ in $\mathsf{CAGp}$ with a totally disconnected kernel and an extension $\mathcal C\in\Cocc(G)$; we wish to show that $p\mathcal D=\mathcal C$ holds for a unique extension $\mathcal D\in\Cocc(H)$. The existence of such $\mathcal D\in\Cocc(H)$ follows from the implication (2)$\Rightarrow$(5) from Theorem~\ref{T:lifts.and.gimel}. To verify the uniqueness part, assume that $\mathcal D_1,\mathcal D_2\in\Cocc(H)$ are such that $p\mathcal D_i=\mathcal C$ for $i=1,2$. Then $\mathcal D_1\mathcal D_2^{-1}\in\Cocc(\ker(p))$. Since $\ker(p)\in\mathsf{CAGp}$ is totally disconnected, its dual group $\ker(p)^*$ is a torsion group and hence $\chi(\mathcal D_1\mathcal D_2^{-1})\in\tor(\Cocc)$ for every $\chi\in H^*$. By virtue of (b), the group $\Cocc$ is torsion-free and hence $\chi(\mathcal D_1\mathcal D_2^{-1})=1$ for every $\chi\in H^*$. Thus, $\mathcal D_1=\mathcal D_2$, as was to be shown.
\end{proof}

\begin{proposition}\label{P:coc.div.consq}
Let $\Flow$ be a minimal flow and assume that the group $\Cocc$ is torsion-free and divisible. Then the following statements hold.
\begin{enumerate}
\item[(i)] For every $G\in\mathsf{CAGp}$ and every $H\sbgp G$ totally disconnected, the quotient morphism $p\colon G\to G/H$ induces an isomorphism of groups
\begin{equation*}
\widehat{p}\colon\Cocc(G)\to\Cocc(G/H).
\end{equation*}
The isomorphism is topological, provided $\Gamma\in\mathsf{LieGp}$ is connected and $X$ is compact.
\item[(ii)] For every $G\in\mathsf{CAGp}$ connected, the projective resolution \emph{(\ref{Eq:proj.res.G})} of $G$ induces isomorphisms of groups
\begin{equation*}
\left(\Cocc\right)^{\dim(G)}\cong\Cocc(\mathfrak{P}(G))\stackrel{\widehat{\mathfrak{p}}}{\longrightarrow}\Cocc(G).
\end{equation*}
The isomorphisms are topological, provided $\Gamma\in\mathsf{LieGp}$ is connected and $X$ is compact.
\item[(iii)] For every $G\in\mathsf{CAGp}$ connected, the maximal toral quotient sequence \emph{(\ref{Eq:max.tor.quot.seq})} of $G$ induces isomorphisms of groups
\begin{equation*}
\Cocc(G)\stackrel{\widehat{q}}{\longrightarrow}\Cocc(\mathbb T(G))\cong(\Cocc)^{\dim(G)}.
\end{equation*}
The isomorphisms are topological, provided $\Gamma\in\mathsf{LieGp}$ is connected and $X$ is compact.
\item[(iv)] Assume that both $\Gamma$ and $X$ are locally compact second countable and that $X$ is non-degenerate. Then there are an isomorphism of abelian groups $\Cocc\cong\mathbb R$ and a topological isomorphism $\Gimel\cong b\mathbb R$.
\end{enumerate}
\end{proposition}
\begin{proof}
We begin by verifying statement (i). First, the morphism $\widehat{p}$ is an isomorphism by virtue of statement (f) from Theorem~\ref{T:coc.divisible}. Further, let $\Gamma\in\mathsf{LieGp}$ be connected and $X$ be compact. Then both $\Cocc(G)$, $\Cocc(G/H)$ are topological groups with the topology of u.c.s. convergence and $\widehat{p}$ is a bijective morphism of topological groups. We show that $\widehat{p}$ is in fact a topological isomorphism by showing that it is open at the identity. So let $\mathcal U$ be an identity neighbourhood in $\Cocc(G)$. By our assumptions on $\Gamma$ and $X$, there exist a compact connected set $1\in C\subseteq\Gamma$ and an identity neighbourhood $V$ in $G$ such that $\mathcal U\supseteq[C\times X;V]$. Let $K\subseteq V$ be an open subgroup of $H$. The quotient group $F=H/K$ is then finite. Further, for every $f\in F$, fix $h_f\in f$. Since the cosets $K+h_f$ ($f\in F$) are compact, there is an identity neighbourhood $V'$ in $G$ with $V'+K\subseteq V$ and with the sets $V'+K+h_f$ ($f\in F$) mutually disjoint.

We claim that $\widehat{p}(\mathcal U)\supseteq[C\times X;p(V')]$; since $[C\times X;p(V')]$ is clearly an identity neighbourhood in $\Cocc(G/H)$, this will finish the proof of (i). So let $\mathcal D\in\Cocc(G/H)$ be an extension with $\mathcal D(C\times X)\subseteq p(V')$ and let $\mathcal C\in\Cocc(G)$ be the lift of $\mathcal D$ across $p$; we show that $\mathcal C\in\mathcal U$. Then 
\begin{equation*}
\mathcal C(C\times X)\subseteq p^{-1}p(V')=V'+H=\bigcup_{f\in F}(V'+K+h_f).
\end{equation*}
Since the set $C\times X$ is connected, the sets $V'+K+h_f$ ($f\in F$) are open and mutually disjoint, and $e\in\mathcal C(C\times X)\cap(V'+K+h_0)\neq\emptyset$, it follows that
\begin{equation*}
\mathcal C(C\times X)\subseteq V'+K+h_0=V'+K\subseteq V.
\end{equation*}
Thus, $\mathcal C\in[C\times X;V]\subseteq\mathcal U$, as was to be shown.

We verify statement (ii). The fact that $\widehat{\mathfrak{p}}$ is a (topological) isomorphism follows from statement (i), for the kernel $\mathfrak{K}(G)$ of $\mathfrak{p}$ is totally disconnected. We verify the other isomorphism from (ii). By definition of $\mathfrak{P}(G)$, there are topological isomorphisms
\begin{equation*}
\mathfrak{P}(G)=(\mathbb Q\otimes G^*)^*\cong\left(\mathbb Q^{\left(\rank(G^*)\right)}\right)^*=\left(\mathbb Q^{\left(\dim(G)\right)}\right)^*\cong\left(\mathbb Q^*\right)^{\dim(G)}.
\end{equation*}
Moreover, since $\mathfrak{P}(\mathbb T^1)\cong\mathbb Q^*$, the part of statement (ii) already proved yields an isomorphism $\Cocc(\mathbb Q^*)\cong\Cocc(\mathbb T^1)=\Cocc$ (the isomorphism is in fact topological, provided $\Gamma\in\mathsf{LieGp}$ is connected and $X$ is compact). It follows that there are isomorphisms of abelian groups
\begin{equation*}
\Cocc(\mathfrak{P}(G))\cong\Cocc\left((\mathbb Q^*)^{\dim(G)}\right)\cong\left(\Cocc(\mathbb Q^*)\right)^{\dim(G)}\cong\left(\Cocc\right)^{\dim(G)}.
\end{equation*}
If $\Gamma\in\mathsf{LieGp}$ is connected and $X$ is compact then all the three isomorphisms above are topological. This finishes the proof of (ii).

We verify statement (iii). The fact that $\widehat{q}$ is a (topological) isomorphism follows from statement (i), for the kernel $\mathbb K(G)$ of $q$ is totally disconnected. The other (topological) isomorphism from (iii) follows from (topological) isomorphisms
\begin{equation*}
\Cocc(\mathbb T(G))\cong\Cocc\left(\mathbb T^{\rank(G^*)}\right)=\Cocc\left(\mathbb T^{\dim(G)}\right)\cong\left(\Cocc\right)^{\dim(G)}.
\end{equation*}

We finish the proof of the proposition by verifying statement (iv). The group $\Cocc$, being divisible and torsion-free, is a rational linear space. Write $\mathfrak{k}=\rank(\Cocc)=\QLSdim(\Cocc)$. By our assumptions on $\Gamma$ and $X$, $\Cocc$ is a Polish group with the topology of u.c.s. convergence and hence its cardinality is at most $\mathfrak{c}$. Consequently, 
\begin{equation}\label{Eq:k.eqls.c.1}
\mathfrak{k}=\rank(\Cocc)\leq\card(\Cocc)\leq\mathfrak{c}.
\end{equation}
Further, by our assumptions on $X$, $C_z(X,\mathbb R)\neq0$. Therefore, $\LSdim(C_z(X,\mathbb R))\geq1$ and hence $\QLSdim(C_z(X,\mathbb R))\geq\mathfrak{c}$. Since the space $X$ is connected by virtue of Theorem~\ref{T:gimel.conn}, $C_z(X,\mathbb T^1)\cong\Cob\subseteq\Cocc$ contains an isomorphic copy of $C_z(X,\mathbb R)$. Consequently,
\begin{equation}\label{Eq:k.eqls.c.2}
\mathfrak{k}=\rank(\Cocc)\geq\rank(C_z(X,\mathbb R))=\QLSdim(C_z(X,\mathbb R))\geq\mathfrak{c}.
\end{equation}
From (\ref{Eq:k.eqls.c.1}) and (\ref{Eq:k.eqls.c.2}) we infer that $\mathfrak{k}=\mathfrak{c}$. This yields isomorphisms of groups $\Cocc\cong\mathbb Q^{(\mathfrak{k})}=\mathbb Q^{(\mathfrak{c})}\cong\mathbb R$ and topological isomorphisms $\Gimel\cong(\Cocc)_d^*\cong(\mathbb R_d)^*\cong b\mathbb R$.
\end{proof}

\section{The non-existence of free minimal extensions}\label{Sub:free.min.ext}

Let $\Flow$ be a minimal flow. In Section~\ref{Sub:free.ext} we have seen that the category $\Coc$ always possesses a free object. Our aim in this section is to show that this is not the case for the subcategory of $\Coc$ formed by the minimal extensions of $\mathcal F$. The main results of this section are Theorems~\ref{T:free.min.not} and~\ref{T:max.min.ext}. The first of them states that if both $\Cocm$ and $\Cob$ are non-trivial then a free object for the category of the minimal extensions from $\Coc$ does not exist. (This applies, for instance, to flows satisfying one of the three assumptions listed at the beginning of Section~\ref{Sub:gp.min.ext} in Chapter~\ref{S:struct.res}.) The second theorem states that every minimal extension $\mathcal C$ from $\Coc$ is an epimorphic image of a minimal extension $\beth$, which is maximal in the sense that it can not be expressed as a non-isomorphic epimorphic image of another minimal extension $\beth'$.

Let $\mathcal F$ be a minimal flow. Throughout this whole section the group $\Cocc$ is assumed to carry the discrete topology. We use notation $(\Gimel,\Daleth)$ for the model of the free extension of $\mathcal F$ constructed in Section~\ref{Sub:free.ext}. Given $G\in\mathsf{CAGp}$ and $\mathcal C\in\Cocc(G)$, we let $h_{\mathcal C}$ stand for the unique morphism $\Gimel\to G$ with $h_{\mathcal C}\Daleth=\mathcal C$. We also write $H_{\mathcal C}=\ker(h_{\mathcal C})$ and $G_{\mathcal C}=\im(h_{\mathcal C})$. Recall from Proposition~\ref{P:further.prop.of.gimel}($\iota$) that $G_{\mathcal C}$ is the closed subgroup of $G$ generated by the values of $\mathcal C$. Consequently, if $\mathcal C,\mathcal D\in\Coc$ and $h$ is a morphism $\mathcal C\to\mathcal D$ then $G_{\mathcal D}=h(G_{\mathcal C})$. Finally, given a subgroup $\Lambda$ of $\Cocc$, we use the symbol $p_{\Lambda}$\index[symbol]{$p_{\Lambda}$} to denote the quotient morphism $\Gimel\to\Gimel/\Lambda^{\perp}$.

\begin{lemma}\label{L:Dal.dense.val}
Let $\Flow$ be a minimal flow and $(\Gimel,\Daleth)$ be its free extension. Then the following statements hold.
\begin{enumerate}
\item[(1)] For every subgroup $\Lambda$ of $\Cocc$, the values of $p_{\Lambda}\Daleth$ generate a dense subgroup of $\Gimel/\Lambda^{\perp}$; that is, $G_{p_{\Lambda\Daleth}}=\Gimel/\Lambda^{\perp}$. In particular, $G_{\Daleth}=\Gimel$.
\item[(2)] Given a group $G\in\mathsf{CAGp}$ and an extension $\mathcal C\in\Cocc(G)$ with $G_{\mathcal C}=G$, there is a unique subgroup $\Lambda$ of $\Cocc$ such that $\mathcal C$ is isomorphic to $p_{\Lambda}\Daleth$ in $\Coc$.
\end{enumerate}
\end{lemma}
\begin{proof}
We begin by verifying the identity $G_{\Daleth}=\Gimel$. To this end, we need to show that for every $\Upsilon\in(\Gimel)^*$, $\Upsilon\Daleth=1$ implies $\Upsilon=1$. First, denote by $\omega$ the Pontryagin isomorphism $\Cocc\to(\Cocc)^{**}=(\Gimel)^*$ corresponding to $\Cocc$. Now, fix $\Upsilon\in(\Gimel)^*$ with $\Upsilon\Daleth=1$ and write $\mathcal C=\omega^{-1}(\Upsilon)\in\Cocc$. Then for all $\gamma\in\Gamma$ and $x\in X$,
\begin{equation*}
\mathcal C(\gamma,x)=\Daleth(\gamma,x)(\mathcal C)=\Upsilon(\Daleth(\gamma,x))=\Upsilon\Daleth(\gamma,x)=1.
\end{equation*}
Thus, $\mathcal C=1$ and hence $\Upsilon=1$. Now we turn to the proof of the main statement from (1). So fix a subgroup $\Lambda$ of $\Cocc$. Then
\begin{equation*}
G_{p_{\Lambda}\Daleth}=p_{\Lambda}(G_{\Daleth})=p_{\Lambda}(\Gimel)=\Gimel/\Lambda^{\perp},
\end{equation*}
as was to be shown.

Now fix $G$ and $\mathcal C$ as in (2). First we check the existence part of statement (2). Let $h_{\mathcal C}$ be the morphism $\Gimel\to G$ with $h_{\mathcal C}\Daleth=\mathcal C$. We claim that $h_{\mathcal C}$ is an epimorphism; indeed, by virtue of (1), 
\begin{equation*}
G=G_{\mathcal C}=G_{h_{\mathcal C}\Daleth}=h_{\mathcal C}(G_{\Daleth})=h_{\mathcal C}(\Gimel).
\end{equation*}
Consequently, there exist a subgroup $\Lambda$ of the group $\Cocc$ and a topological isomorphism $\varphi\colon\Gimel/\Lambda^{\perp}\to G$ such that $\varphi p_{\Lambda}=h_{\mathcal C}$. Then $\varphi(p_{\Lambda}\Daleth)=h_{\mathcal C}\Daleth=\mathcal C$, whence it follows that $\varphi$ is an isomorphism $p_{\Lambda}\Daleth\to\mathcal C$.

Now we check the uniqueness part of statement (2). So let $\Lambda$ satisfy the conclusion of (2) and fix an underlying isomorphism $\varphi\colon p_{\Lambda}\Daleth\to\mathcal C$. Then $(\varphi p_{\Lambda})\Daleth=\mathcal C=h_{\mathcal C}\Daleth$, and so $\varphi p_{\Lambda}$ and $h_{\mathcal C}$ agree on the set of all values of $\Daleth$. By virtue of (1), this means that $\varphi p_{\Lambda}=h_{\mathcal C}$ on the whole $\Gimel$. Since $\varphi$ is an isomorphism by the assumptions, we obtain $\Lambda^{\perp}=\ker(p_{\Lambda})=\ker(h_{\mathcal C})$. Thus, $\Lambda=\omega^{-1}(\Lambda^{\perp\perp})=\omega^{-1}(\ker(h_{\mathcal C})^{\perp})$. This verifies the uniqueness part of (2).
\end{proof}

\begin{lemma}\label{L:repr.min.ext.gp}
Let $\Flow$ be a minimal flow. Then the following statements hold.
\begin{enumerate}
\item[(i)] Given groups $H,K\in\mathsf{CAGp}$ and extensions $\mathcal C\in\Cocc(H)$, $\mathcal D\in\Cocc(K)$ with $G_{\mathcal C}=H$ and $G_{\mathcal D}=K$, the set $\Hom(\mathcal C,\mathcal D)$ is either a singleton or empty. In the first case the corresponding morphism from $\Hom(\mathcal C,\mathcal D)$ is an epimorphism.
\item[(ii)] Given subgroups $\Lambda$ and $\Xi$ of $\Cocc$, the set $\Hom(p_{\Lambda}\Daleth,p_{\Xi}\Daleth)$ is non-empty if and only if $\Xi\subseteq\Lambda$.
\end{enumerate}
\end{lemma}
\begin{remark}\label{R:repr.min.ext.gp}
Recall that if $G\in\mathsf{CAGp}$ and $\mathcal C\in\Cocc(G)$ is minimal then $G_{\mathcal C}=G$. Consequently, statement (i) from Lemma~\ref{L:repr.min.ext.gp}, as well as statement (2) from Lemma~\ref{L:Dal.dense.val}, applies to all minimal extensions from $\Coc$.
\end{remark}
\begin{proof}[Proof of Lemma~\ref{L:repr.min.ext.gp}]
Fix $H$, $K$, $\mathcal C$ and $\mathcal D$ as in (i) and assume that $\varphi$ is a morphism $\mathcal C\to\mathcal D$. Consider the morphisms $h_{\mathcal C}\colon\Gimel\to H$ and $h_{\mathcal D}\colon\Gimel\to K$ with $h_{\mathcal C}\Daleth=\mathcal C$ and $h_{\mathcal D}\Daleth=\mathcal D$, respectively. Then $(\varphi h_{\mathcal C})\Daleth=\varphi\mathcal C=\mathcal D=h_{\mathcal D}\Daleth$ and hence, by Lemma~\ref{L:Dal.dense.val}(1), $\varphi h_{\mathcal C}=h_{\mathcal D}$. Since both $h_{\mathcal C}$ and $h_{\mathcal D}$ are epimorphisms by our assumptions on $\mathcal C$ and $\mathcal D$, it follows that such a morphism $\varphi$ is at most one and, if it exists, it is in fact an epimorphism. This verifies statement (i).

Now fix subgroups $\Lambda$ and $\Xi$ of $\Cocc$. First we verify the ``if'' part of statement (ii). So assume that $\Xi\subseteq\Lambda$. Then $\Lambda^{\perp}\subseteq\Xi^{\perp}$ and so there is a morphism $\varphi\colon\Gimel/\Lambda^{\perp}\to\Gimel/\Xi^{\perp}$ with $\varphi p_{\Lambda}=p_{\Xi}$. Consequently, $\varphi(p_{\Lambda}\Daleth)=p_{\Xi}\Daleth$ and hence $\varphi$ is a desired morphism $p_{\Lambda}\Daleth\to p_{\Xi}\Daleth$. Now we turn to proving the ``only if'' part of statement (ii). So let $\varphi$ be a morphism of extensions $p_{\Lambda}\Daleth\to p_{\Xi}\Daleth$. Then $\varphi$ is a morphism of topological groups $\Gimel/\Lambda^{\perp}\to\Gimel/\Xi^{\perp}$ with $(\varphi p_{\Lambda})\Daleth=p_{\Xi}\Daleth$. Thus, by Lemma~\ref{L:Dal.dense.val}(1), $\varphi p_{\Lambda}=p_{\Xi}$. This shows that $\Lambda^{\perp}\subseteq\Xi^{\perp}$ and hence $\Xi\subseteq\Lambda$.
\end{proof}

\begin{lemma}\label{L:dns.val.whn.min}
Let $\Flow$ be a minimal flow, $(\Gimel,\Daleth)$ be the free extension of $\mathcal F$ and $\Lambda$ be a subgroup of $\Cocc$. Then $p_{\Lambda}\Daleth\in\Cocc(\Gimel/\Lambda^{\perp})$ is minimal if and only if $\Lambda\cap\Cob=1$.
\end{lemma}
\begin{proof}
By Lemma~\ref{L:Dal.dense.val}(1), $G_{p_{\Lambda}\Daleth}=\Gimel/\Lambda^{\perp}$. Moreover, $H_{p_{\Lambda}\Daleth}=\ker(p_{\Lambda})=\Lambda^{\perp}$ and, by the first part of Remark~\ref{R:prop.of.gimel}, $F(\Daleth)=(\Cob)^{\perp}$. Consequently, it follows from Proposition~\ref{P:further.prop.of.gimel}($\iota\iota\iota$) that the following conditions are equivalent:
\begin{itemize}
\item $p_{\Lambda}\Daleth$ is minimal,
\item $(H_{p_{\Lambda}\Daleth})^{\perp}\cap F(\Daleth)^{\perp}=1$,
\item $\Lambda^{\perp\perp}\cap(\Cob)^{\perp\perp}=1$
\item $\Lambda\cap\Cob=1$.
\end{itemize}
\end{proof}

Let $\mathcal F$ be a minimal flow, $\mGimel\in\mathsf{CAGp}$\index[symbol]{$\gimel_{\mathcal F}^\text{m}$} and $\mDaleth\in\Cocc(\mGimel)$.\index[symbol]{$\daleth_{\mathcal F}^\text{m}$} We say that $(\mGimel,\mDaleth)$ is a \emph{free minimal extension of $\mathcal F$}\index{free!minimal extension of a minimal flow} if the following condition holds:
\begin{enumerate}
\item[($m$)] for every $G\in\mathsf{CAGp}$ and every minimal extension $\mathcal C\in\Cocc(G)$ there exists $q\in\Hom(\mGimel,G)$ with $q\mDaleth=\mathcal C$.
\end{enumerate}
Observe that the morphism $q$ from ($m$) is unique by Lemma~\ref{L:repr.min.ext.gp}(i). Consequently, if a free minimal extension of $\mathcal F$ does exist then it is unique up to an isomorphism. We shall therefore speak of \emph{the} free minimal extension of $\mathcal F$.

\begin{lemma}\label{L:ex.free.min?}
Let $\Flow$ be a minimal flow. Then the following conditions are equi\-va\-lent:
\begin{enumerate}
\item[(a)] the free minimal extension of $\mathcal F$ exists,
\item[(b)] there is a largest subgroup $\Omega$ of $\Cocc$ with $\Omega\cap\Cob=1$.
\end{enumerate}
\end{lemma}
\begin{remark}\label{R:ex.free.min?}
Assume that the group $\Cocc$ is torsion-free; that is, $X$ is connected and $\Gamma$ has no finite abelian quotient groups. Then a subgroup $\Omega$ of $\Cocc$ satisfies $\Omega\cap\Cob=1$ if and only if it is a subgroup of the groupoid $\Cocm$.
\end{remark}
\begin{proof}[Proof of Lemma~\ref{L:ex.free.min?}]
We show that (b) follows from (a). So assume that $(\mGimel,\mDaleth)$ is the free minimal extension of $\mathcal F$ and denote by $h$ the unique morphism $\Gimel\to\mGimel$ with $h\Daleth=\mDaleth$. By Lemma~\ref{L:Dal.dense.val}(2), we may assume, without loss of generality, that $\mGimel=\Gimel/\Omega^{\perp}$ and $\mDaleth=p_{\Omega}\Daleth$ for an appropriate subgroup $\Omega$ of $\Cocc$. We show that $\Omega$ fulfills condition (b). First, by virtue of Lemma~\ref{L:dns.val.whn.min}, the minimality of $\mDaleth=p_{\Omega}\Daleth$ yields $\Omega\cap\Cob=1$. Further, fix $\Lambda\sbgp\Cocc$ with $\Lambda\cap\Cob=1$; we show that $\Lambda\subseteq\Omega$. By virtue of Lemma~\ref{L:dns.val.whn.min}, $p_{\Lambda}\Daleth\in\Cocc(\Gimel/\Lambda^{\perp})$ is a minimal extension. Since $p_{\Omega}\Daleth$ is the free minimal extension of $\mathcal F$, we have $\Hom(p_{\Omega}\Daleth,p_{\Lambda}\Daleth)\neq\emptyset$. Thus, by Lemma~\ref{L:repr.min.ext.gp}(ii), $\Lambda\subseteq\Omega$.

We show that (a) follows from (b). So let $\Omega$ be as in (b). Set $\mGimel=\Gimel/\Omega^{\perp}$ and $\mDaleth=p_{\Omega}\Daleth$; we show that $(\mGimel,\mDaleth)$ is the free minimal extension of $\mathcal F$. First, the extension $\mDaleth$ is minimal by Lemma~\ref{L:dns.val.whn.min}. Further, fix $G\in\mathsf{CAGp}$ and a minimal extension $\mathcal C\in\Cocc(G)$; we need to find a morphism $q\colon\mGimel\to G$ with $q\mDaleth=\mathcal C$. By Lemma~\ref{L:Dal.dense.val}(2), we may assume, without loss of generality, that $G=\Gimel/\Lambda^{\perp}$ and $\mathcal C=p_{\Lambda}\Daleth$ for an appropriate subgroup $\Lambda$ of $\Cocc$. By virtue of Lemma~\ref{L:dns.val.whn.min}, $\Lambda\cap\Cob=1$ and hence $\Lambda\subseteq\Omega$ by our choice of $\Omega$. Lemma~\ref{L:repr.min.ext.gp}(ii) then yields a desired morphism $\mDaleth\to p_{\Lambda}\Daleth$.
\end{proof}

\begin{theorem}\label{T:free.min.not}
Let $\Flow$ be a minimal flow with both $\Cob$ and $\Cocm$ non-trivial. Then the free minimal extension of $\mathcal F$ does not exist.
\end{theorem}
\begin{remark}\label{R:free.min.not}
The assumptions of the theorem are satisfied, provided the flow $\mathcal F$ fulfills one of the three assumptions listed at the beginning of Section~\ref{Sub:gp.min.ext} in Chapter~\ref{S:struct.res}. First, in all the three situations $X$ is a non-trivial compact space and so the group $\Cob\cong C_z(X,\mathbb T^1)$ is non-trivial. Second, it follows from Theorems~\ref{T:grp.min.ext.am}, \ref{T:grp.min.ext.simp} and~\ref{T:gp.min.free.ccl} that the groupoid $\Cocm$ is also non-trivial.
\end{remark}
\begin{proof}[Proof of Theorem~\ref{T:free.min.not}]
We proceed by contradiction, assuming that $\mathcal F$ does possess the free minimal extension. By virtue of Lemma~\ref{L:ex.free.min?}, there is a largest subgroup $\Omega$ of $\Cocc$ with the property $\Omega\cap\Cob=1$. Fix $\mathcal C\in\Cocm\setminus1$ and $\mathcal D\in\Cob\setminus1$. Then the subgroups $\langle\mathcal C\rangle$ and $\langle\mathcal C\mathcal D\rangle$ of $\Cocc$ generated by $\mathcal C$ and $\mathcal C\mathcal D$, respectively, are contained in the groupoid $\Cocm$, and hence they intersect the group $\Cob$ only at the identity. Thus, by definition of $\Omega$, both $\langle\mathcal C\rangle$ and $\langle\mathcal C\mathcal D\rangle$ are contained in $\Omega$. Consequently, $\mathcal D=(\mathcal C\mathcal D)\mathcal C^{-1}\in\Omega\cap(\Cob\setminus1)$, which contradicts the definition of $\Omega$.
\end{proof}

\begin{lemma}\label{L:Lmb.Bth.max}
Let $\Flow$ be a minimal flow, $\beth$ be a minimal extension from $\Coc$ and $\Lambda$ be the subgroup of $\Cocc$ with $p_{\Lambda}\Daleth$ isomorphic to $\beth$. Then the following conditions are equivalent:
\begin{enumerate}
\item[($\alpha$)] every minimal extension $\beth'$ from $\Coc$ with $\Hom(\beth',\beth)\neq\emptyset$ is isomorphic to $\beth$,
\item[($\beta$)] $\Lambda$ is maximal (with respect to inclusion) among the subgroups of $\Cocc$ intersecting $\Cob$ only at the identity.
\end{enumerate}
\end{lemma}
\begin{remark}\label{R:Lmb.Bth.max}
We wish to mention the following facts.
\begin{itemize}
\item Condition ($\alpha$) can be reformulated by saying that $\beth$ is a maximal object in the category of the minimal extensions from $\Coc$.
\item By Lemma~\ref{L:Dal.dense.val}(2), the isomorphism classes of the extensions $\beth$ satisfying ($\alpha$) are in a one-to-one correspondence with the subgroups $\Lambda$ of $\Cocc$ satisfying ($\beta$). A group $\Lambda$ satisfying ($\beta$) is not unique in general and so two extensions $\beth$ satisfying ($\alpha$) need not be isomorphic; see Example~\ref{Ex:many.max.min.ext} below.
\item Let $\Lambda_1,\Lambda_2$ be subgroups of $\Cocc$ satisfying ($\beta$) and let $\beth_i=p_{\Lambda_i}\Daleth\in\Cocc(\Gimel/\Lambda_i^{\perp})$ for $i=1,2$. As mentioned above, $\beth_1,\beth_2$ need not be isomorphic in $\Coc$. As a matter of fact, it is possible that the groups $\Gimel/\Lambda_1^{\perp}$, $\Gimel/\Lambda_2^{\perp}$ are not isomorphic in $\mathsf{CAGp}$; see Example~\ref{Ex:max.min.ext.dist} below.
\end{itemize}
\end{remark}
\begin{proof}[Proof of Lemma~\ref{L:Lmb.Bth.max}]
Before turning to the proof oberve that by Lemma~\ref{L:Dal.dense.val}(2), there is indeed a unique subgroup $\Lambda$ of $\Cocc$ such that $\beth$ is isomorphic to $p_{\Lambda}\Daleth$. Also, since $\beth$ is minimal, $\Lambda\cap\Cob=1$ by virtue of Lemma~\ref{L:dns.val.whn.min}.

We show that ($\beta$) follows from ($\alpha$). So assume that ($\alpha$) holds and fix a subgroup $\Xi$ of $\Cocc$ with $\Xi\cap\Cob=1$ and $\Lambda\subseteq\Xi$. Then the extension $\beth'=p_{\Xi}\Daleth$ is minimal by Lemma~\ref{L:dns.val.whn.min}. Furthermore, Lemma~\ref{L:repr.min.ext.gp}(ii) yields $\Hom(\beth',\beth)\neq\emptyset$. Consequently, $\beth'$ and $\beth$ are isomorphic in $\Coc$ by virtue of ($\alpha$), and Lemma~\ref{L:repr.min.ext.gp}(ii) thus yields $\Lambda=\Xi$. This verifies ($\beta$).

We show that ($\alpha$) follows from ($\beta$). So assume that ($\beta$) holds and fix a minimal extension $\beth'$ from $\Coc$ with $\Hom(\beth',\beth)\neq\emptyset$.  By Lemma~\ref{L:Dal.dense.val}(2), $\beth'$ is isomorphic to $p_{\Xi}\Daleth$ for an appropriate subgroup $\Xi$ of $\Cocc$. Moreover, Lemmas~\ref{L:dns.val.whn.min} and~\ref{L:repr.min.ext.gp} yield $\Xi\cap\Cob=1$ and $\Lambda\subseteq\Xi$. Thus, by applying ($\beta$), we obtain $\Xi=\Lambda$. Consequently, we have isomorphisms of extensions $\beth'\cong p_{\Xi}\Daleth=p_{\Lambda}\Daleth\cong\beth$. This verifies ($\alpha$).
\end{proof}

\begin{theorem}\label{T:max.min.ext}
Let $\Flow$ be a minimal flow and $\mathcal C$ be a minimal extension from $\Coc$. Then there is a minimal extension $\beth$ in $\Coc$ with the following properties:
\begin{enumerate}
\item[(a)] there is an epimorphism $\beth\to\mathcal C$,
\item[(b)] every minimal extension $\beth'$ from $\Coc$ with $\Hom(\beth',\beth)\neq\emptyset$ is isomorphic to $\beth$.
\end{enumerate}
\end{theorem}
\begin{remark}\label{R:max.min.ext}
In view of our discussion from Remark~\ref{R:Lmb.Bth.max} we may reformulate the statement of the theorem by saying that in the category of the minimal extensions from $\Coc$, every object $\mathcal C$ is an epimorphic image of a maximal object $\beth$.
\end{remark}
\begin{proof}[Proof of Theorem~\ref{T:max.min.ext}]
By Lemma~\ref{L:Dal.dense.val}(2), $\mathcal C$ is isomorphic in $\Coc$ to $p_{\Xi}\Daleth$ for an appropriate subgroup $\Xi$ of $\Cocc$. Moreover, Lemma~\ref{L:dns.val.whn.min} yields $\Xi\cap\Cob=1$. By an elementary argument involving Zorn's lemma one finds a subgroup $\Lambda\supseteq\Xi$ of $\Cocc$, which is maximal (with respect to inclusion) among the subgroups of $\Cocc$ intersecting $\Cob$ only at the identity. Now set $\beth=p_{\Lambda}\Daleth$. Then $\beth$ is a minimal extension of $\mathcal F$ due to Lemma~\ref{L:dns.val.whn.min}. Further, since $\Xi\subseteq\Lambda$, condition (a) holds by virtue of Lemma~\ref{L:repr.min.ext.gp}. Finally, condition (b) follows from Lemma~\ref{L:Lmb.Bth.max} by our choice of $\Lambda$.
\end{proof}

\begin{example}\label{Ex:many.max.min.ext}
Let $\Flow$ be a minimal flow with $X$ non-degenerate locally compact and assume that $\Cocm\neq1$. Let $\mathfrak{A}$ be the family of all subgroups $\Lambda$ of $\Cocc$, which are maximal with respect to the condition $\Lambda\cap\Cob=1$. By the second part of Remark~\ref{R:Lmb.Bth.max}, $\mathfrak{A}$ can be identified with the set of all isomorphism classes of the maximal objects in the category of the minimal extensions from $\Coc$. To demonstrate that there are many such isomorphism classes in general, we show that $\card(\mathfrak{A})\geq\mathfrak{c}$ under our current assumptions on $\mathcal F$.

Fix $\mathcal C\in\Cocm\setminus1$. For $\mathcal D\in\Cob$ let $\Lambda_{\mathcal D}$ be an element of $\mathfrak{A}$ with $\mathcal C\mathcal D\in\Lambda_{\mathcal D}$. (Such a group $\Lambda_{\mathcal D}$ of course exists by Zorn's lemma, since $\langle\mathcal C\mathcal D\rangle\cap\Cob=1$.) We claim that distinct coboundaries $\mathcal D$ correspond to distinct groups $\Lambda_{\mathcal D}$. To see this, fix $\mathcal D_1,\mathcal D_2\in\Cob$ and assume that $\Lambda_{\mathcal D_1}=\Lambda_{\mathcal D_2}$. Then 
\begin{equation*}
\Cob\ni\mathcal D_1\mathcal D_2^{-1}=(\mathcal C\mathcal D_1)(\mathcal C\mathcal D_2)^{-1}\in\Lambda_{\mathcal D_1}\Lambda_{\mathcal D_2}^{-1}=\Lambda_{\mathcal D_1}.
\end{equation*}
Since $\Lambda_{\mathcal D_1}\cap\Cob=1$, it follows that $\mathcal D_1=\mathcal D_2$. Consequently, 
\begin{equation*}
\card(\mathfrak{A})\geq\card(\Cob)=\card(C_z(X,\mathbb T^1))\geq\mathfrak{c},
\end{equation*}
as was to be shown.
\end{example}

\begin{example}\label{Ex:max.min.ext.dist}
In the previous example we have seen that two extensions $\beth$ satisfying condition ($\alpha$) from Lemma~\ref{L:Lmb.Bth.max} need not be isomorphic in $\Coc$. In this example we show that it is in fact possible for their sections $F(\beth)=G_{\beth}$ to be non-isomorphic in $\mathsf{CAGp}$. To this end, we shall use a result from Chapter~\ref{S:struct.res}.

Let $\mathcal F\colon\mathbb R\times\mathbb T^1\curvearrowright\mathbb T^2$ be defined as the product action of a minimal continuous flow on $\mathbb T^1$ with the natural action of $\mathbb T^1$ on itself. Clearly, $\mathcal F$ thus defined is a minimal flow. Choose $z=(1,1)$ as the base point for $\mathbb T^2$ and consider the corresponding motion map $\mathcal F_z\colon\mathbb R\times\mathbb T^1\to\mathbb T^2$. The induced morphism $\mathcal F_z^{\sharp}\colon H_1^w(\mathbb R\times\mathbb T^1)\to H_1^w(\mathbb T^2)$ is a monomorphism and it takes the form of an inclusion of groups $0\oplus\mathbb Z\subseteq\mathbb Z\oplus\mathbb Z$. It follows that the flow $\mathcal F$ is topologically free and possesses a free cycle. By virtue of (\ref{Eq:Cob.in.Cocc}) from Theorem~\ref{T:F.in.Gim.free.ccl}, the inclusion $\Cob\subseteq\Cocc$ takes the form of an inclusion of abelian groups
\begin{equation}\label{Eq:Cb.Cc.Ex.bt}
\begin{array}{ccccccccccc}
\Cob & \cong & 0 & \oplus & \mathbb R & \oplus & \mathbb Z & \oplus & \mathbb Z & \subseteq & {}\\
{} & \subseteq & \mathbb R & \oplus & \mathbb R & \oplus & \mathbb Q & \oplus & \mathbb Z & \cong & \Cocc.
\end{array}
\end{equation}
By writing the inclusion $0\subseteq\mathbb R$ in the form $0\oplus0\subseteq\mathbb R\oplus\mathbb Q$ and by reorganizing the direct summands from the isomorphism above, we put (\ref{Eq:Cb.Cc.Ex.bt}) into the following form:
\begin{equation*}
\begin{array}{ccccccccccc}
\Cob & \cong & 0 & \oplus & \mathbb R & \oplus & \mathbb Z & \oplus & (0\oplus\mathbb Z) & \subseteq & {}\\
{} & \subseteq & \mathbb R & \oplus & \mathbb R & \oplus & \mathbb Q & \oplus & (\mathbb Q\oplus\mathbb Z) & \cong & \Cocc.
\end{array}
\end{equation*}
Set $\Lambda_1=\mathbb R\oplus0\oplus0\oplus A_1$ and $\Lambda_2=\mathbb R\oplus0\oplus0\oplus A_2$, where $A_1=\mathbb Q\oplus0\cong\mathbb Q$ and $A_2=\{(k,k) : k\in\mathbb Z\}\cong\mathbb Z$. For $i=1,2$ set $\beth_i=p_{\Lambda_i}\Daleth$. Elementary arguments show that both $\Lambda_1$ and $\Lambda_2$ are maximal among the subgroups of $\Cocc$ intersecting $\Cob$ only at the identity, and hence both $\beth_1$ and $\beth_2$ are minimal and satisfy condition ($\alpha$) from Lemma~\ref{L:Lmb.Bth.max}. We claim that the groups $G_{\beth_1}$ and $G_{\beth_2}$ are not isomorphic in $\mathsf{CAGp}$. Indeed, there are topological isomorphisms
\begin{equation*}
G_{\beth_1}=F(\beth_1)=\Gimel/\Lambda_1^{\perp}\cong\Lambda_1^*\cong b\mathbb R\times\mathbb Q^*
\end{equation*}
and
\begin{equation*}
G_{\beth_2}=F(\beth_2)=\Gimel/\Lambda_2^{\perp}\cong\Lambda_2^*\cong b\mathbb R\times\mathbb T^1.
\end{equation*}
These two groups are not topologically isomorphic, for the first one is torsion-free, while the other one is not.
\end{example}

\begin{remark}\label{R:dim.G.beth.cnst}
Let $\Flow$ be a minimal flow with $\Cocc$ torsion-free. As has been shown in the preceding example, two distinct extensions $\beth$ satisfying condition ($\alpha$) from Lemma~\ref{L:Lmb.Bth.max} may have non-isomorphic sections $F(\beth)=G_{\beth}$. We would like to show now that, after all, there does exist an invariant for such extensions $\beth$, namely the topological dimension of the section $\dim(G_{\beth})$. So fix $\beth$ satisfying ($\alpha$) from Lemma~\ref{L:Lmb.Bth.max} and let $\Lambda$ be the subgroup of $\Cocc$ with $\beth\cong p_{\Lambda}\Daleth$ in $\Coc$. By Lemma~\ref{L:Lmb.Bth.max}, $\Lambda$ is maximal among the subgroups of $\Cocc$ intersecting $\Cob$ only at the identity. Consequently, the divisible hull $\mathbb Q\otimes\Cocc$ of $\Cocc$ takes the form of a direct sum of abelian groups (and hence also of rational linear spaces) $\mathbb Q\otimes\Cocc=\left(\mathbb Q\otimes\Cob\right)\oplus\left(\mathbb Q\otimes\Lambda\right)$. Thus,
\begin{equation*}
\begin{split}
\dim\left(G_{\beth}\right)&=\dim\left(\Gimel/\Lambda^{\perp}\right)=\dim\left(\Lambda^*\right)=\rank(\Lambda)=\QLSdim(\mathbb Q\otimes\Lambda)\\
&=\QLSdim\left((\mathbb Q\otimes\Cocc)/(\mathbb Q\otimes\Cob)\right)
\end{split}
\end{equation*}
does not depend on the choice of $\beth$.
\end{remark}

\section{Cohomology classes as group morphisms}\label{S:coh.cls.gp.mrph}

Let $\Flow$ be a minimal flow and assume that the group $\Cocc$ is divisible. In Theorem~\ref{T:structure.coc} from Section~\ref{Sub:gp.ext} we have seen that for every $G\in\mathsf{CAGp}$, the group $\Cocc(G)$ can be viewed as the group of morphisms $G^*\to\Cocc$. For cohomology groups such an isomorphism fails to hold. Though the group $\Coch(G)$ can be viewed as a subgroup of $\Hom(G^*,\Coch)$, the extent to which these two groups differ is represented by the group of extensions $\Ext(G^*,\pi^1(X))$; this is proved in Theorems~\ref{T:intro.embed.G(F)} and~\ref{T:exact.sequence}. In fact, as we show in Theorem~\ref{P:tor.gen.G}, the obstruction to an isomorphism $\Coch(G)\cong\Hom(G^*,\Coch)$ lies in the existence of extensions $\mathcal C\in\Cocc$ of $\mathcal F$ with totally disconnected sections $F(\mathcal C)$. This also suggests that there might be a connection between the existence of such extensions $\mathcal C$ on one side and the first cohomotopy group $\pi^1(X)$ of $X$ on the other side. This is indeed the case, as we shall demonstrate later in Chapter~\ref{S:alg.top.asp}.

\begin{theorem}\label{T:intro.embed.G(F)}
Let $\Flow$ be a minimal flow, $G\in\mathsf{CAGp}$ and $\pi$ be the canonical quotient morphism $\Cocc\to\Coch$. Then the map
\begin{equation}\label{Eq:ext.hom.G*}
\Psi_G\index[symbol]{$\Psi_G$}\colon\Coch(G)\ni\mathcal C\mapsto\pi\mathcal C^*\in\Hom(G^*,\Coch )
\end{equation}
is a topological isomorphism onto its image, provided the groups $\Coch$, $\Coch(G)$ are equipped with the ext-topology and the group $\Hom(G^*,\Coch)$ carries the topology of point-wise convergence.
\end{theorem}
\begin{remark}\label{R:intro.embed.G(F)}
We wish to add the following remarks.
\begin{itemize}
\item In order to simplify notation, we are identifying an extension $\mathcal C\in\Cocc(G)$ with its cohomology class in the definition of $\Psi_G$. We shall do so also in the sequel in situations when no misunderstanding can arise.
\item We shall see in the proof of the theorem that the morphism $\Psi_G$ can also be defined formally as a unique morphism $\Coch(G)\to\Hom(G^*,\Coch)$ with $\Psi_G\pi_G=\Hom(G^*,\pi)\Phi_G$.
\item The theorem provides us with the following inclusions of topological groups for every $G\in\mathsf{CAGp}$:
\begin{equation*}
\Coch(G)\subseteq\Hom(G^*,\Coch)\subseteq\left(\Coch\right)^{G^*}.
\end{equation*}
Since $\Coch$ is a discrete group by Remark~\ref{R:group.top.cont}, it follows that $\Coch(G)$ is always totally disconnected.
\item The image $\im(\Psi_G)$ of $\Psi_G$ may or may not be closed in $\Hom(G^*,\Coch)$, depending on a particular situation. See, for instance, Theorem~\ref{P:CochG.closed.Y.N} in Chapter~\ref{S:struct.res}.
\end{itemize}
\end{remark}
\begin{proof}[Proof of Theorem~\ref{T:intro.embed.G(F)}]
Consider the isomorphism $\Phi_G\colon\Cocc(G)\to\Hom(G^*,\Cocc )$ defined in Theorem~\ref{T:structure.coc} and the morphism
\begin{equation*}
\Hom(G^*,\pi)\colon\Hom(G^*,\Cocc )\to\Hom(G^*,\Coch )
\end{equation*}
derived from $\pi$. We claim that
\begin{equation*}
\ker(\pi_G)=\Cob(G)=\ker(\Hom(G^*,\pi)\Phi_G).
\end{equation*}
Indeed, the first equality is an immediate consequence of the definition of $\pi_G$, whereas the second equality follows from the equivalence of the following statements for every $\mathcal C\in\Cocc(G)$:
\begin{itemize}
\item $\mathcal C\in\ker(\Hom(G^*,\pi)\Phi_G)$,
\item $\pi\mathcal C^*=1$,
\item $\chi\mathcal C\in\Cob $ for every $\chi\in G^*$,
\item $\mathcal C\in\Cob(G)$.
\end{itemize}
(For the equivalence of the third and the fourth statement, see Corollary~\ref{C:min.contr.char}.) It follows that there exists a unique morphism of groups 
\begin{equation*}
\Psi_G\colon\Coch(G)\to\Hom(G^*,\Coch )
\end{equation*}
with $\Psi_G\pi_G=\Hom(G^*,\pi_G)\Phi_G$ and that $\Psi_G$ is in fact a monomorphism. One verifies at once that $\Psi_G(\mathcal C)=\pi\mathcal C^*$ for every $\mathcal C\in\Coch(G)$.

Now we show that $\Psi_G$ is a topological isomorphism onto its image. First we verify the continuity of $\Psi_G$. So assume that a convergent net $\mathcal C_i\to e$ is given in $\Coch(G)$. Then $F(\mathcal C_i)\to e$ in $2^G$ and hence $F(\mathcal C^*_i(\chi))=F(\chi\mathcal C_i)=\chi F(\mathcal C_i)\to1$ in $2^{\mathbb T^1}$ for every $\chi\in G^*$. This yields convergence $\pi\mathcal C^*_i\to e$ in $\Hom(G^*,\Coch )$ and verifies the continuity of $\Psi_G$.

Conversely, assume that for a given net $(\mathcal C_i)_{i\in J}$ in $\Coch(G)$, we have convergence $\pi\mathcal C^*_i\to e$ in $\Hom(G^*,\Coch )$. Fix a neighbourhood $V$ of $e$ in $G$. There exist $n\in\mathbb N$, $\chi_1,\dots,\chi_n\in G^*$ and neighbourhoods $W_1\dots,W_n$ of $1$ in $\mathbb T^1$ such that $V\supseteq\bigcap_{l=1}^n\chi_l^{-1}(W_l)$. From the convergence $\pi\mathcal C^*_i\to e$ in $\Hom(G^*,\Coch)$ it follows that there is $i_0\in J$ with $\chi_lF(\mathcal C_i)=F(\chi_l\mathcal C_i)=F(\mathcal C^*_i(\chi_l))\subseteq W_l$ for all $i\geq i_0$ and $l=1,\dots,n$. Consequently, $F(\mathcal C_i)\subseteq \bigcap_{l=1}^n\chi_l^{-1}(W_l)\subseteq V$ for every $i\geq i_0$. This shows that $F(\mathcal C_i)\to e$ in $2^G$ and hence $\mathcal C_i\to e$ in $\Coch(G)$.
\end{proof}

\begin{corollary}\label{C:intro.embed.F(D)}
Let $\Flow$ be a minimal flow and $G\in\mathsf{CAGp}$. Given $\mathcal C\in\Cocc(G)$, let $h_{\mathcal C}$ be the morphism $\Gimel\to G$ with $h_{\mathcal C}\Daleth=\mathcal C$. Then the map
\begin{equation*}
\widetilde{\Psi}_G\index[symbol]{$\widetilde{\Psi}_G$}\colon\Coch(G)\ni\mathcal C\mapsto\left(\omega_G\right)^{-1}\Psi_G(\mathcal C)^*(\pi^*)^{-1}=h_{\mathcal C}\in\Hom(F(\Daleth),G)
\end{equation*}
is a topological isomorphism onto its image, provided $\Coch(G)$ is equipped with the ext-topology and $\Hom(F(\Daleth),G)$ carries the topology of uniform convergence.
\end{corollary}
\begin{remark}\label{R:intro.embed.F(D)}
Similarly to Theorem~\ref{T:intro.embed.G(F)}, in our definition of $\widetilde{\Psi}_G$ we are identifying every extension $\mathcal C\in\Cocc(G)$ with its cohomology class. This time, the abuse of notation is excused by the fact that for equivalent extensions $\mathcal C,\mathcal D\in\Cocc(G)$, the corresponding morphisms $h_{\mathcal C}$ and $h_{\mathcal D}$ agree on $F(\Daleth)$. (Indeed, if $\mathcal C=\mathcal D\co(\xi)$ then $h_{\mathcal C}(h_{\mathcal D})^{-1}=h_{\co(\xi)}$ and $h_{\co(\xi)}F(\Daleth)=F(h_{\co(\xi)}\Daleth)=F(\co(\xi))=e$.) We shall use this identification also in the proof of Corollary~\ref{C:intro.embed.F(D)} below.
\end{remark}
\begin{proof}[Proof of Corollary~\ref{C:intro.embed.F(D)}]
Since both $G^*$ and $\Coch$ carry the discrete topology, there is a topological isomorphism $\Hom(G^*,\Coch)\to\Hom((\Coch)^*,G^{**})$, which associates with every morphism $k\colon G^*\to\Coch$ its Pontryagin dual $k^*\colon G^{**}\to(\Coch)^*$. Moreover, the Pontryagin isomorphism $\omega_G$ corresponding to $G$ is a topological isomorphism $G\to G^{**}$ and the dual morphism $\pi^*$ of $\pi$ is a topological isomorphism $(\Coch)^*\to(\Cob)^{\perp}$. Finally, since $(\Cob)^{\perp}=F(\Daleth)$ by the first part of Remark~\ref{R:prop.of.gimel}, we have a topological isomorphism $\Hom(G^*,\Coch)\to\Hom(F(\Daleth),G)$, which associates with every morphism $k\colon G^*\to\Coch$ the morphism $(\omega_G)^{-1}k^*(\pi^*)^{-1}\colon F(\Daleth)\to G$. Thus, since $\Psi_G\colon\Coch(G)\to\Hom(G^*,\Coch)$ is a topological isomorphism onto its image by virtue of Theorem~\ref{T:intro.embed.G(F)}, it follows that $\widetilde{\Psi}_G$ is also a topological isomorphism onto its image. To finish the proof, we need to show that $(\omega_G)^{-1}\Psi_G(\mathcal C)^*(\pi^*)^{-1}=h_{\mathcal C}$ for every $\mathcal C\in\Cocc(G)$.

Fix $\mathcal C\in\Cocc(G)$. Given $\chi\in G^*$, consider the set
\begin{equation*}
\mathcal A_{\chi,\mathcal C}=\left\{\Lambda\in\Gimel=(\Cocc)^* : \Lambda(\chi\mathcal C)=(\chi h_{\mathcal C})(\Lambda)\right\}.
\end{equation*}
Clearly, $\mathcal A_{\chi,\mathcal C}$ is a closed subgroup of $\Gimel$. Further, $\mathcal A_{\chi,\mathcal C}$ contains all the values of $\Daleth$; indeed, for all $\gamma\in\Gamma$ and $x\in X$,
\begin{equation*}
\Daleth(\gamma,x)(\chi\mathcal C)=(\chi\mathcal C)(\gamma,x)=\chi(\mathcal C(\gamma,x))=\chi(h_{\mathcal C}\Daleth(\gamma,x))=(\chi h_{\mathcal C})(\Daleth(\gamma,x)).
\end{equation*}
These two observations lead to the inclusion $F(\Daleth)\subseteq\mathcal A_{\chi,\mathcal C}$.

We verify the desired identity $(\omega_G)^{-1}\Psi_G(\mathcal C)^*(\pi^*)^{-1}=h_{\mathcal C}$. We accomplish this by showing that $\Psi_G(\mathcal C)^*=\omega_Gh_{\mathcal C}\pi^*$. For all $\Upsilon\in(\Coch)^*$ and $\chi\in G^*$ we have
\begin{equation}\label{Eq:C:int.emb.F(D)1}
\Psi_G(\mathcal C)^*(\Upsilon)=\Upsilon\Psi_G(\mathcal C)\colon\chi\mapsto \Upsilon\left(\Psi_G(\mathcal C)(\chi)\right)=\Upsilon(\pi(\chi\mathcal C))=(\Upsilon\pi)(\chi\mathcal C),
\end{equation}
and
\begin{equation}\label{Eq:C:int.emb.F(D)2}
\left(\omega_G h_{\mathcal C}\pi^*\right)(\Upsilon)=\omega_G\left(h_{\mathcal C}(\Upsilon\pi)\right)\colon\chi\mapsto\chi\left(h_{\mathcal C}(\Upsilon\pi)\right)=(\chi h_{\mathcal C})(\Upsilon\pi).
\end{equation}
Since $\Upsilon\pi=\pi^*(\Upsilon)\in\im(\pi^*)=F(\Daleth)\subseteq\mathcal A_{\chi,\mathcal C}$ for all $\Upsilon$ and $\chi$ by our discussion above, the right-hand sides of (\ref{Eq:C:int.emb.F(D)1}) and (\ref{Eq:C:int.emb.F(D)2}) are equal. This verifies the identity $\Psi_G(\mathcal C)^*=\omega_Gh_{\mathcal C}\pi^*$.
\end{proof}

\begin{theorem}\label{T:exact.sequence}
Let $\Flow$ be a minimal flow on a compact connected space~$X$ and let $G\in\mathsf{CAGp}$. Then there exist 
\begin{itemize}
\item a morphism $\Omega\index[symbol]{$\Omega$}\colon\Hom(G^*,\Coch )\to\Ext(G^*,\pi^1(X))$, and
\item a monomorphism $\nu\index[symbol]{$\nu$}\colon\pi^1(X)\to\Cocc $, 
\end{itemize}
such that the following sequence is exact:
\begin{equation}\label{Eq:exact.sequence}
\begin{split}
0&\longrightarrow\Cob(G)\stackrel{\mu_G}{\longrightarrow}\Cocc(G)\stackrel{\kappa_G}{\longrightarrow}\Hom(G^*,\Coch )
\stackrel{\Omega}{\longrightarrow}\\
&\stackrel{\Omega}{\longrightarrow}\Ext(G^*,\pi^1(X))\stackrel{\nu_{e}}{\longrightarrow}\Ext(G^*,\Cocc )\stackrel{\pi_{e}}{\longrightarrow}\Ext(G^*,\Coch )\longrightarrow 0,
\end{split}
\end{equation}
where $\kappa_G\index[symbol]{$\kappa_G$}=\Psi_G\pi_G$, $\nu_e\index[symbol]{$\nu_e$}=\Ext(G^*,\nu)$ and $\pi_e\index[symbol]{$\pi_e$}=\Ext(G^*,\pi)$. If, in addition, the group $\Cocc$ is divisible then the following sequence is exact:
\begin{equation}\label{Eq:ex.seq.ext.Gamma.CLAC}
\begin{split}
0&\longrightarrow\Cob(G)\stackrel{\mu_G}{\longrightarrow}\Cocc(G)\stackrel{\kappa_G}{\longrightarrow}\Hom(G^*,\Coch )
\stackrel{\Omega}{\longrightarrow}\Ext(G^*,\pi^1(X))\longrightarrow0.
\end{split}
\end{equation}
Finally, if $\Cocc$ is divisible and $G$ is connected then, under the identification $\im(\Psi_G)\cong\Coch(G)$, the latter group is a direct summand in $\Hom(G^*,\Coch)$ with the complementary summand isomorphic to $\Ext(G^*,\pi^1(X))$. Thus, there is a direct sum
\begin{equation}\label{Eq:quot.gp.Gamma.CLAC}
\Hom\left(G^*,\Coch\right)=\Coch(G)\oplus\Ext\left(G^*,\pi^1(X)\right).
\end{equation}
\end{theorem}
\begin{remark}\label{R:exact.sequence}
It follows from the exactness of (\ref{Eq:exact.sequence}) and from the definition of $\kappa_G$ that $\Psi_G$ is an isomorphism, provided $\Ext(G^*,\pi^1(X))=0$. This is the case in the following situations:
\begin{itemize}
\item the group $\pi^1(X)$ is divisible,
\item the group $G^*$ is free abelian, that is, $G$ is a torus,
\item $\pi^1(X)$ is finitely generated (in particular, $X$ is a manifold) and $G^*$ is a Whitehead group (that is, $G$ is arc-wise connected).
\end{itemize}
\end{remark}
\begin{proof}[Proof of Theorem~\ref{T:exact.sequence}]
Consider the covariant Hom-Ext sequence derived from the sequence (\ref{Eq:ex.seq.ext.T1}) and corresponding to the group $G^*$:
\begin{equation}\label{Eq:aux.hom.ext}
\begin{split}
0&\longrightarrow\Hom(G^*,\Cob )\stackrel{\mu_h}{\longrightarrow}\Hom(G^*,\Cocc )\stackrel{\pi_h}{\longrightarrow}\Hom(G^*,\Coch )\\
&\stackrel{\epsilon}{\longrightarrow}\Ext(G^*,\Cob )\stackrel{\mu_e}{\longrightarrow}\Ext(G^*,\Cocc )\stackrel{\pi_e}{\longrightarrow}\Ext(G^*,\Coch )\longrightarrow 0.
\end{split}
\end{equation}
Recall that $\mu_h=\Hom(G^*,\mu)$, $\pi_h=\Hom(G^*,\pi)$, $\mu_e=\Ext(G^*,\mu)$, $\pi_e=\Ext(G^*,\pi)$ and $\epsilon$ is the connecting morphism. Our aim now is to show that the sequence (\ref{Eq:aux.hom.ext}) is equivalent to the sequence (\ref{Eq:exact.sequence}). We shall proceed by replacing all the arrows in (\ref{Eq:aux.hom.ext}) step by step. The result of this construction will be summarized in the diagrams in Figures~\ref{Fig:ex.seq.ext1} and~\ref{Fig:ex.seq.ext2} below.

By virtue of Corollary~\ref{C:min.contr.char}, the isomorphism $\Phi_G$ from Theorem~\ref{T:structure.coc} restricts to an isomorphism $\Phi_G'\colon\Cob(G)\to\Hom(G^*,\Cob )$. Moreover, $\Phi_G\mu_G=\Hom(G^*,\mu)\Phi_G'=\mu_h\Phi_G'$. It follows that the arrow
\begin{equation*}
\mu_h=\Hom(G^*,\mu)\colon\Hom(G^*,\Cob )\longrightarrow\Hom(G^*,\Cocc )
\end{equation*}
in (\ref{Eq:aux.hom.ext}) can be replaced by the arrow $\mu_G\colon\Cob(G)\longrightarrow\Cocc(G)$.

Further, the arrow 
\begin{equation*}
\pi_h=\Hom(G^*,\pi)\colon\Hom(G^*,\Cocc)\to\Hom(G^*,\Coch)
\end{equation*}
in (\ref{Eq:aux.hom.ext}) can be replaced by the arrow $\pi_h\Phi_G\colon\Cocc(G)\to\Hom(G^*,\Coch)$, this is clear. Since $\pi_h\Phi_G=\Psi_G\pi_G=\kappa_G$ by definition of $\Psi_G$ and $\kappa_G$, this yields the third arrow in (\ref{Eq:exact.sequence}).

Fix a base point $z$ for $X$. There are isomorphisms of groups $\Cob\cong C_z(X,\mathbb T^1)\cong C_z(X,\mathbb R)\oplus\pi^1(X)$; denote by $j$ the underlying inclusion morphism $\pi^1(X)\to\Cob$. We claim that $j$ induces an isomorphism of groups
\begin{equation*}
j_e=\Ext(G^*,j)\colon\Ext(G^*,\pi^1(X))\to\Ext(G^*,\Cob ).
\end{equation*}
To see this, let $\pr_1$, $\pr_2$ be the projection morphisms from $\Cob$ onto $C_z(X,\mathbb R)$, $\pi^1(X)$, respectively, and write $q_1=\Ext(G^*,\pr_1)$, $q_2=\Ext(G^*,\pr_2)$. Since $\pr_2 j=\id_{\pi^1(X)}$, it follows that $q_2 j_e=\id_{\Ext(G^*,\pi^1(X))}$ and hence $j_e$ is a monomorphism. Further, from the exactness of the sequence
\begin{equation*}
0\longrightarrow\pi^1(X)\stackrel{j}{\longrightarrow}\Cob\stackrel{\pr_1}{\longrightarrow}C_z(X,\mathbb R)\longrightarrow0
\end{equation*}
we obtain an exact sequence
\begin{equation}\label{Eq:je.epi.HE.part}
\Ext(G^*,\pi^1(X))\stackrel{j_e}{\longrightarrow}\Ext(G^*,\Cob)\stackrel{q_1}{\longrightarrow}\Ext(G^*,C_z(X,\mathbb R))\longrightarrow0.
\end{equation}
By divisibility of $C_z(X,\mathbb R)$, we have $\Ext(G^*,C_z(X,\mathbb R))=0$ and so $j_e$ is an epimorphism due to the exactness of (\ref{Eq:je.epi.HE.part}). Thus, indeed, $j_e$ is an isomorphism of groups. Consequently, the arrow $\epsilon\colon\Hom(G^*,\Coch)\to\Ext(G^*,\Cob)$ in (\ref{Eq:aux.hom.ext}) can be replaced by the arrow $\Omega\colon\Hom(G^*,\Coch)\to\Ext(G^*,\pi^1(X))$, where $\Omega=j_e^{-1}\epsilon$.

Finally, set $\nu=\mu j$ and $\nu_e=\Ext(G^*,\nu)$. Then $\nu_e=\mu_ej_e$ and the arrow
\begin{equation*}
\mu_e=\Ext(G^*,\mu)\colon\Ext(G^*,\Cob)\to\Ext(G^*,\Cocc)
\end{equation*}
in (\ref{Eq:aux.hom.ext}) thus turns into
\begin{equation*}
\nu_e=\Ext(G^*,\nu)\colon \Ext(G^*,\pi^1(X))\to\Ext(G^*,\Cocc),
\end{equation*}
which is precisely the fifth arrow in (\ref{Eq:exact.sequence}). The observations made so far can be summarized in the form of the commutative diagrams from Figures~\ref{Fig:ex.seq.ext1} and~\ref{Fig:ex.seq.ext2}.

\begin{figure}[ht]
\[\minCDarrowwidth20pt\begin{CD}
\Hom(G^*,\Cob) @>\mu_h>> \Hom(G^*,\Cocc) @>\pi_h>> \Hom(G^*,\Coch) \\
@A\Phi_G'AA @A\Phi_GAA @| \\
\Cob(G) @>>\mu_G> \Cocc(G) @>>\kappa_G> \Hom(G^*,\Coch) 
\end{CD}\]
\caption{Obtaining (\ref{Eq:exact.sequence}) from (\ref{Eq:aux.hom.ext}), part 1.}
\label{Fig:ex.seq.ext1}
\end{figure}

\begin{figure}[ht]
\[\minCDarrowwidth20pt\begin{CD}
\Hom(G^*,\Coch) @>\epsilon>> \Ext(G^*,\Cob) @>\mu_e>> \Ext(G^*,\Cocc) @>\pi_e>> \Ext(G^*,\Coch)\\
@| @Aj_eAA @| @| \\
\Hom(G^*,\Coch) @>>\Omega> \Ext(G^*,\pi^1(X)) @>>\nu_e> \Ext(G^*,\Cocc) @>>\pi_e> \Ext(G^*,\Coch)
\end{CD}\]
\caption{Obtaining (\ref{Eq:exact.sequence}) from (\ref{Eq:aux.hom.ext}), part 2.}
\label{Fig:ex.seq.ext2}
\end{figure}

\noindent Since all the vertical arrows in Figures~\ref{Fig:ex.seq.ext1} and~\ref{Fig:ex.seq.ext2} are isomorphisms, the exactness of (\ref{Eq:exact.sequence}) follows from the exactness of (\ref{Eq:aux.hom.ext}). The exactness of (\ref{Eq:ex.seq.ext.Gamma.CLAC}) follows at once from the exactness of (\ref{Eq:exact.sequence}), for $\Ext(G^*,\Cocc )=0$ as soon as the group $\Cocc $ is divisible.

We finish the proof of the theorem by verifying (\ref{Eq:quot.gp.Gamma.CLAC}). First, since $\pi_G\colon\Cocc(G)\to\Coch(G)$ is an epimorphism and $\kappa_G=\Psi_G\pi_G$ by definition, we have $\im(\Psi_G)=\im(\kappa_G)$. Further, by the exactness of (\ref{Eq:ex.seq.ext.Gamma.CLAC}), $\im(\kappa_G)=\ker(\Omega)$ and
\begin{equation*}
\Hom(G^*,\Coch)/\ker(\Omega)\cong\Ext(G^*,\pi^1(X)).
\end{equation*}
Consequently,
\begin{equation}\label{Eq:aux.compl.sum.ext}
\begin{split}
\Hom(G^*,\Coch)/\Coch(G)&=\Hom(G^*,\Coch)/\im(\Psi_G)\\
&=\Hom(G^*,\Coch)/\im(\kappa_G)\\
&=\Hom(G^*,\Coch)/\ker(\Omega)\\
&\cong\Ext\left(G^*,\pi^1(X)\right).
\end{split}
\end{equation}
By virtue of Theorem~\ref{T:lifts.and.gimel}, the group $\Cocc(G)$ is divisible due to divisibility of $\Cocc$ and connectedness of $G$, and hence its quotient group $\Coch(G)$ is also divisible. It follows that $\Coch(G)$ is a direct summand in $\Hom(G^*,\Coch)$ and, by virtue of (\ref{Eq:aux.compl.sum.ext}), its complementary summand is isomorphic to $\Ext\left(G^*,\pi^1(X)\right)$. This verifies (\ref{Eq:quot.gp.Gamma.CLAC}).
\end{proof}

In order to formulate our next result we fix some notation. Given a minimal flow $\mathcal F$ and a group $G\in\mathsf{CAGp}$, we use the symbol $\Coctd(G)$\index[symbol]{$\Coctd(G)$} to denote the set of all extensions $\mathcal C\in\Cocc(G)$ with a totally disconnected section $F(\mathcal C)\subseteq G$. Also, we write $\Coctd$\index[symbol]{$\Coctd$} instead of $\Coctd(\mathbb T^1)$.

\begin{lemma}\label{L:tor.gen.G}
Let $\Flow$ be a minimal flow. Then for every $G\in\mathsf{CAGp}$, $\Coctd(G)$ is a subgroup of $\Cocc(G)$. Moreover, if $\Cocc$ is divisible and $G\in\mathsf{CAGp}$ is connected then $\Coctd(G)$ is also divisible and hence it is a direct summand in $\Cocc(G)$.
\end{lemma}
\begin{proof}
We divide the proof into four steps.

\emph{1st step.} Let $G\in\mathsf{CAGp}$ and $H\sbgp G$. We recall the equivalence of the following statements:
\begin{enumerate}
\item[($\alpha$)] $H$ is totally disconnected,
\item[($\beta$)] $\bigcap_{k\in\mathbb N}\kappa_k(H)=e$,
\item[($\gamma$)] for every identity neighbourhood $U$ in $G$ there is $k\in\mathbb N$ such that $\kappa_l(H)\subseteq U$ for every $l\in\mathbb N$ divisible by $k$.
\end{enumerate}

First, ($\alpha$) follows from ($\beta$), for if $H_0\neq e$ then, for every $k\in\mathbb N$, $\kappa_k(H)\supseteq\kappa_k(H_0)=H_0\neq e$. To show that ($\beta$) follows from ($\alpha$), assume that $H$ is totally disconnected. Then $H$ is isomorphic to a subgroup of the product of a family of finite abelian groups. Since every finite abelian group is annihilated by $\kappa_k$ for some $k\in\mathbb N$, we conclude that for every identity neighbourhood $V$ in $H$ there is $k$ with $\kappa_k(H)\subseteq V$. This verifies ($\beta$). The equivalence of ($\beta$) and ($\gamma$) follows from the fact that $\kappa_l(H)\subseteq\kappa_k(H)$ as soon as $k$ divides $l$.

\emph{2nd step.} We show that if $H,K$ are closed totally disconnected subgroups of $G$ and $h\in\mathbb N$, then both $HK$ and $\kappa_h^{-1}(H)$ are (closed and) totally disconnected. (Here, $\kappa_h$ is the $h$-endomorphism of the whole group $G$.)

Fix an identity neighbourhood $U$ in $G$. Further, choose an identity neighbourhood $V$ in $G$ with $VV\subseteq U$. By virtue of the equivalence ($\alpha$)$\Leftrightarrow$($\gamma$) from the first step of the proof, there are $k,l\in\mathbb N$ with $\kappa_k(H)\subseteq V$ and $\kappa_l(K)\subseteq V$. Then 
\begin{equation*}
\kappa_{kl}(HK)=\kappa_{kl}(H)\kappa_{kl}(K)\subseteq\kappa_k(H)\kappa_l(K)\subseteq VV\subseteq U.
\end{equation*}
Since $U$ was an arbitrary identity neighbourhood in $G$, it follows by the first step of the proof that $HK$ is indeed totally disconnected. Finally, the group $\kappa_h^{-1}(H)$ is totally disconnected by virtue of the equivalence ($\alpha$)$\Leftrightarrow$($\beta$), for
\begin{equation*}
\begin{split}
\bigcap_{k\in\mathbb N}\kappa_k\left(\kappa_h^{-1}(H)\right)&\subseteq\bigcap_{k\in\mathbb N}\kappa_{kh}\left(\kappa_h^{-1}(H)\right)=\bigcap_{k\in\mathbb N}\kappa_k\left(\kappa_h\left(\kappa_h^{-1}(H)\right)\right)\\
&\subseteq\bigcap_{k\in\mathbb N}\kappa_k(H)=e.
\end{split}
\end{equation*}

\emph{3rd step.} Let $G\in\mathsf{CAGp}$. We show that $\Coctd(G)$ is a subgroup of $\Cocc(G)$.

Fix $\mathcal C,\mathcal D\in\Coctd(G)$. Then $F(\mathcal C\mathcal D)\subseteq F(\mathcal C)F(\mathcal D)$ by Lemma~\ref{L:recall.to.metric}(c). Since $F(\mathcal C)F(\mathcal D)$ is a totally disconnected subgroup of $G$ by the second step of the proof, it follows that the group $F(\mathcal C\mathcal D)$ is also totally disconnected and hence $\mathcal C\mathcal D\in\Coctd(G)$. Also, $F(\mathcal C^{-1})=F(\mathcal C)$ by Lemma~\ref{L:recall.to.metric}(b) and so $\mathcal C^{-1}\in\Coctd(G)$. This shows that $\Coctd(G)$ is indeed a subgroup of $\Cocc(G)$. 

\emph{4th step.} Assume that the group $\Cocc$ is divisible and fix $G\in\mathsf{CAGp}$ connected. We show that the group $\Coctd(G)$ is divisible.

Fix $\mathcal C\in\Coctd(G)$ and $h\in\mathbb N$. By the assumptions of the theorem and by virtue of Theorem~\ref{T:lifts.and.gimel}, the group $\Cocc(G)$ is divisible. Consequently, there is $\mathcal D\in\Cocc(G)$ with $\mathcal D^h=\mathcal C$; we show that $\mathcal D\in\Coctd(G)$. We have $\kappa_h(F(\mathcal D))=F(\kappa_h\mathcal D)=F(\mathcal D^h)=F(\mathcal C)$ and so $F(\mathcal D)\subseteq\kappa_h^{-1}(F(\mathcal C))$. By the second step of the proof, $\kappa_h^{-1}(F(\mathcal C))$ is totally disconnected and hence so is $F(\mathcal D)$. Thus, $\mathcal D\in\Coctd(G)$, as was to be shown.
\end{proof}

\begin{theorem}\label{P:tor.gen.G}
Let $\Flow$ be a minimal flow on a compact connected space~$X$. Assume that the group $\Cocc$ is divisible and fix a complementary direct summand $\Coccn$\index[symbol]{$\Coccn$} to $\Coctd$ in $\Cocc$. Given $G\in\mathsf{CAGp}$ connected, set
\begin{itemize}
\item $\Cochtd(G)=\pi_G(\Coctd(G))=\Coctd(G)/\Cob(G)$,\index[symbol]{$\Cochtd(G)$}
\item $\Coccn(G)=\Phi_G^{-1}\left(\Hom(G^*,\Coccn)\right)$,\index[symbol]{$\Coccn(G)$}
\item $\Cochcn(G)=\pi_G\left(\Coccn(G)\right)$,\index[symbol]{$\Cochcn(G)$}
\item $\Cochtd=\Cochtd(\mathbb T^1)$\index[symbol]{$\Cochtd$} and $\Cochcn=\Cochcn(\mathbb T^1)$.\index[symbol]{$\Cochcn$}
\end{itemize}
Then the following statements hold:
\begin{enumerate}
\item[(1)] $F(\mathcal C)\subseteq G$ is connected for every $\mathcal C\in\Coccn(G)$,
\item[(2)] there is a direct sum
\begin{equation*}
\Cocc(G)=\Coctd(G)\oplus\Coccn(G),
\end{equation*}
\item[(3)] there is a direct sum
\begin{equation*}
\Coch(G)=\Cochtd(G)\oplus\Cochcn(G),
\end{equation*}
which is topological, provided $\Coch(G)$ carries the ext-topology,
\item[(4)] the restriction $\Psi_G\colon\Cochcn(G)\to\Hom(G^*,\Cochcn)$ is an isomorphism of groups,
\item[(5)] under the identification $\Psi_G(\Cochtd(G))\cong\Cochtd(G)$, the latter group is a direct summand in $\Hom(G^*,\Cochtd)$ and its complementary summand is isomorphic to the group $\Ext(G^*,\pi^1(X))$; consequently, there is a direct sum
\begin{equation}\label{Eq:tor.gen.G}
\Hom(G^*,\Cochtd)=\Cochtd(G)\oplus\Ext(G^*,\pi^1(X)).
\end{equation}
\end{enumerate}
\end{theorem}
\begin{remark}\label{R:tor.gen.G}
Notice the following facts.
\begin{itemize}
\item Assume that $G$ has no proper non-trivial closed connected subgroups. Then, since $\Cob(G)\subseteq\Coctd(G)$ by definition of $\Coctd(G)$, it follows from statement (1) that $\Coccn(G)$ is a subgroup of the groupoid $\Cocm(G)$. Recall also that for a non-trivial connected group $G\in\mathsf{CAGp}$ the following conditions are equivalent:
\begin{enumerate}
\item[$\circ$] $G$ has no proper non-trivial closed connected subgroups,
\item[$\circ$] $\rank(G^*)=1$.
\end{enumerate}
These conditions are satisfied, for instance, by the circle group $G=\mathbb T^1$ and by the solenoids.
\item By virtue of (3) applied to the group $\mathbb T^1$, there is a direct sum
\begin{equation*}
\Hom(G^*,\Coch)=\Hom(G^*,\Cochtd)\oplus\Hom(G^*,\Cochcn).
\end{equation*}
Consequently, by virtue of (4) and (5), it is the groups $\Cochtd$ and $\Cochtd(G)$ that are ``responsible'' for the non-surjectivity of the monomorphism $\Psi_G$. That is, $\Psi_G$ is an isomorphism if and only if its restriction $\Psi_G\colon\Cochtd(G)\to\Hom(G^*,\Cochtd)$ is an isomorphism.
\item The splitting in statement (2), unlike the one in statement (3), need not be topological if $\Cocc(G)$ carries the topology of u.c.s. convergence. For examples of such a situation we refer to Proposition~\ref{P:CLAC.min.lim.fin} and Remark~\ref{R:CLAC.min.lim.fin} in Chapter~\ref{S:alg.top.asp}.
\end{itemize}
\end{remark}
\begin{proof}[Proof of Theorem~\ref{P:tor.gen.G}]
Before turning to the proof of the theorem notice the following facts. First, by virtue of Lemma~\ref{L:tor.gen.G}, $\Coctd$ is indeed a direct summand in $\Cocc$. Second, though unique only up to an isomorphism, $\Coccn$ is always a subgroup of the groupoid $\Cocm$. Now we turn to the proof of the theorem, proceeding in five steps.

\emph{1st step.} We verify statement (1).

Fix $\mathcal C\in\Coccn(G)$. By definition of $\Coccn(G)$, $\chi\mathcal C=\mathcal C^*(\chi)\in\Coccn$ for every $\chi\in G^*$. It follows that for every $\chi\in G^*$, $\chi F(\mathcal C)=F(\chi\mathcal C)$ is either $1$ or $\mathbb T^1$. Consequently, $F(\mathcal C)$ is a connected subgroup of $G$ and statement (1) thus follows.

\emph{2nd step.} We verify statement (2).

By definition of $\Coccn$, there is a direct sum
\begin{equation}\label{Eq:dec.td.cn.hom}
\Hom(G^*,\Cocc)=\Hom(G^*,\Coctd)\oplus\Hom(G^*,\Coccn).
\end{equation}
Under the isomorphism $\Phi_G$ from Theorem~\ref{T:structure.coc}, the following correspondences hold:
\begin{enumerate}
\item[(i)] $\Coctd(G)$ corresponds to $\Hom(G^*,\Coctd)$; indeed, given $\mathcal C\in\Cocc(G)$, the following conditions are equivalent:
\begin{itemize}
\item $\mathcal C\in\Coctd(G)$,
\item $F(\mathcal C)\subseteq G$ is totally disconnected,
\item $F(\chi\mathcal C)=\chi F(\mathcal C)\subseteq\mathbb T^1$ is totally disconnected for every $\chi\in G^*$,
\item $\mathcal C^*(\chi)=\chi\mathcal C\in\Coctd$ for every $\chi\in G^*$,
\item $\Phi_G(\mathcal C)=\mathcal C^*\in\Hom(G^*,\Coctd)$,
\end{itemize}
\item[(ii)] $\Coccn(G)$ corresponds to $\Hom(G^*,\Coccn)$; this holds by definition of $\Coccn(G)$.
\end{enumerate}
Thus, statement (2) follows from (\ref{Eq:dec.td.cn.hom}), (i) and (ii).

\emph{3rd step.} We verify statement (3).

First we check that the direct sum in (3) holds algebraically. The fact that $\Coch(G)$ is a sum of $\Cochtd(G)$ and $\Cochcn(G)$ follows by applying $\pi_G$ to the direct sum in (2). We need only to show that $\Cochtd(G)\cap\Cochcn(G)=e$, but this is immediate. For if $\mathcal C\in\Cocc(G)$ and $F(\mathcal C)$ is both connected and totally disconnected then $F(\mathcal C)=e$ and hence $\mathcal C\in\Cob$.

Now we show that the splitting in (3) is topological. To this end, fix a directed set $J$ and nets $(\mathcal C_i)_{i\in J}$ in $\Coctd(G)$ and $(\mathcal D_i)_{i\in J}$ in $\Coccn(G)$. We need to show that $F(\mathcal C_i\mathcal D_i)\to e$ in $2^G$ if and only if $F(\mathcal C_i),F(\mathcal D_i)\to e$ in $2^G$. Since all $F(\mathcal C_i)$ ($i\in J$) are totally disconnected and all $F(\mathcal D_i)$ ($i\in J$) are connected, Corollary~\ref{C:essen.disj.fin} yields $F(\mathcal C_i\mathcal D_i)=F(\mathcal C_i)F(\mathcal D_i)$ for every $i\in J$. From this identity the desired equivalence follows at once.

\emph{4th step.} We verify statement (4).

Notice the following facts:
\begin{itemize}
\item the restriction $\Phi_G\colon\Coccn(G)\to\Hom(G^*,\Coccn)$ is an isomorphism; this follows from the definition of $\Coccn(G)$ and from the fact that $\Phi_G$ is an isomorphism,
\item the restriction $\pi_G\colon\Coccn(G)\to\Cochcn(G)$ is an isomorphism; this follows from the definition of $\Cochcn(G)$ and from the identity $\Coccn(G)\cap\Cob(G)=e$ (the latter holds by virtue of (2)),
\item the restriction $\Hom(G^*,\pi)\colon\Hom(G^*,\Coccn)\to\Hom(G^*,\Cochcn)$ is an isomorphism; this follows from the fact that (by the preceding statement applied to the group $\mathbb T^1$) the restriction $\pi\colon\Coccn\to\Cochcn$ is an isomorphism,
\item the morphism $\Psi_G\colon\Coch(G)\to\Hom(G^*,\Coch)$ is defined by the identity $\Psi_G\pi_G=\Hom(G^*,\pi)\Phi_G$; this follows from Remark~\ref{R:intro.embed.G(F)}.
\end{itemize}
The four facts listed above now yield statement (4).

\emph{5th step.} We verify statement (5).

First notice that $\Psi_G$ maps $\Cochtd(G)$ into $\Hom(G^*,\Cochtd)$, since $\Phi_G$ maps $\Coctd(G)$ into $\Hom(G^*,\Coctd)$ by virtue of (i) from the second step of the proof. We shall need the following facts:
\begin{itemize}
\item under the isomorphism $\Psi_G(\Coch(G))\cong\Coch(G)$, the latter group is a direct summand in $\Hom(G^*,\Coch)$ and its complementary summand is isomorphic to the group $\Ext\left(G^*,\pi^1(X)\right)$; this follows from (\ref{Eq:quot.gp.Gamma.CLAC}) in Theorem~\ref{T:exact.sequence},
\item under the isomorphism $\Psi_G(\Cochtd(G))\cong\Cochtd(G)$, the latter group is a direct summand in $\Hom(G^*,\Cochtd)$; indeed, by virtue of Lemma~\ref{L:tor.gen.G}, $\Coctd(G)$ is a divisible group and hence so is its quotient group $\Cochtd(G)$,
\item there is a direct sum
\begin{equation*}
\Coch(G)\cong\Cochtd(G)\oplus\Cochcn(G);
\end{equation*}
this follows from statement (3),
\item there is a direct sum
\begin{equation*}
\Hom(G^*,\Coch)=\Hom(G^*,\Cochtd)\oplus\Hom(G^*,\Cochcn);
\end{equation*}
this follows from statement (3) applied to the group $\mathbb T^1$,
\item $\Psi_G$ maps $\Cochcn(G)$ isomorphically onto $\Hom(G^*,\Cochcn)$; this follows from statement (4).
\end{itemize}
These five facts lead, on one hand, to a diagram in Figure~\ref{Fig:verif.stat.5} and, on the other hand, to the isomorphisms
\begin{equation*}
\begin{split}
\Ext\left(G^*,\pi^1(X)\right)&\cong\frac{\Hom\left(G^*,\Coch\right)}{\Coch(G)}\\
&\cong\frac{\Hom(G^*,\Cochtd)}{\Cochtd(G)}\oplus\frac{\Hom\left(G^*,\Cochcn\right)}{\Cochcn(G)}\\
&\cong\frac{\Hom(G^*,\Cochtd)}{\Cochtd(G)}.
\end{split}
\end{equation*}
\begin{figure}[ht]
\[\minCDarrowwidth20pt\begin{CD}
\Hom(G^*,\Coch) @.{} = {}@. \Hom(G^*,\Cochtd) @. {}\oplus{} @. \Hom(G^*,\Cochcn) \\
@A\Psi_GA\text{mon}A @. @A\Psi_GA\text{mon}A @. @A\Psi_GA\text{iso}A\\
\Coch(G) @. = @. \Cochtd(G) @. \oplus @. \Cochcn(G).
\end{CD}\]
\caption{Monomorphic and isomorphic restrictions of $\Psi_G$.}
\label{Fig:verif.stat.5}
\end{figure}
This verifies statement (5).
\end{proof}

\section{Torsions as non-minimal extensions}\label{Sub:tor.elem}

Let $\Flow$ be a minimal flow. In this section we study the torsion subgroup $\tor(\Coch)$ of $\Coch$ and its subgroups $\tork(\Coch)$ ($k\in\mathbb N$). We are motivated by the fact that the torsion elements of $\Coch$ correspond to the non-minimal extensions from $\Cocc$, see Corollary~\ref{C:MinH.rel.torH}. We begin our investigation with the situation when $\Cocc$ is a divisible group, which includes the case when the acting group $\Gamma$ of $\mathcal F$ is simply connected; see Theorem~\ref{T:tor.struct}. Then we embark on the case when $\Gamma$ is a connected Lie group and $X$ is a compact (connected) manifold; see Theorem~\ref{T:tor.Lie.mnfld}.

\begin{theorem}\label{T:tor.struct}
Let $\Flow$ be a minimal flow on a compact connected space $X$ and assume that the group $\Cocc $ is divisible. Given $k\in\mathbb N$, there are short exact sequences
\begin{equation}\label{Eq:tor.str.k}
0\longrightarrow\Hom(\Gamma,\mathbb Z_k)\longrightarrow\tork(\Coch )\longrightarrow\pi^1(X)/k\pi^1(X)\longrightarrow0
\end{equation}
and
\begin{equation}\label{Eq:tor.str.inf}
0\longrightarrow\bigcup_{l\in\mathbb N}\Hom(\Gamma,\mathbb Z_l)\longrightarrow\tor(\Coch )\longrightarrow(\mathbb Q/\mathbb Z)\otimes\pi^1(X)\longrightarrow0.
\end{equation}
Moreover, if $\Gamma$ has no non-trivial finite abelian quotient groups then there are isomorphisms
\begin{equation}\label{Eq:tork.simp.iso}
\tork(\Coch )\cong\pi^1(X)/k\pi^1(X)
\end{equation}
and
\begin{equation}\label{Eq:tor.simp.iso}
\tor(\Coch )\cong(\mathbb Q/\mathbb Z)\otimes\pi^1(X),
\end{equation}
and, finally, there exists a divisible subgroup $\mathfrak{D}$ of the groupoid $\Cochm$ such that
\begin{equation}\label{Eq:DEM.in.simp}
\Coch=\mathfrak{D}\oplus\tor(\Coch)\cong\mathfrak{D}\oplus\Big((\mathbb Q/\mathbb Z)\otimes\pi^1(X)\Big).
\end{equation}
\end{theorem}
\begin{remark}\label{R:tor.struct}
Notice the following facts.
\begin{itemize}
\item The union $\bigcup_{l\in\mathbb N}\Hom(\Gamma,\mathbb Z_l)$ from (\ref{Eq:tor.str.inf}) is to be interpreted as a subgroup of $\Hom(\Gamma,(\mathbb Q/\mathbb Z)_d)\subseteq\Hom(\Gamma,(\mathbb T^1)_d)$, under the usual identifications $\mathbb Z_l\cong\tor_l(\mathbb T^1)\subseteq\mathbb Q/\mathbb Z\subseteq\mathbb T^1$. An elementary argument shows that if the quotient group $\Gamma/\Gamma_0$ is finitely generated then this union coincides with $\Hom(\Gamma,(\mathbb Q/\mathbb Z)_d)$.
\item An extension $\mathcal C\in\Cocc$ corresponds to an element of $\tork(\Coch)$ if and only if $\mathcal C^k\in\Cob$. This is equivalent to the inclusion $F(\mathcal C)\subseteq\mathbb Z_k$.
\end{itemize}
\end{remark}
\begin{proof}[Proof of Theorem~\ref{T:tor.struct}]
Consider the exact sequence for tensor and torsion products derived from (\ref{Eq:ex.seq.ext.T1}) and corresponding to the group $\mathbb Z_k$:
\begin{equation}\label{Eq:tor.Ext.k.aux}
\begin{split}
0&\longrightarrow\Tor(\mathbb Z_k,\Cob)\longrightarrow\Tor(\mathbb Z_k,\Cocc)\longrightarrow\Tor(\mathbb Z_k,\Coch)\\
&\longrightarrow\mathbb Z_k\otimes\Cob \longrightarrow\mathbb Z_k\otimes\Cocc \longrightarrow\mathbb Z_k\otimes\Coch \longrightarrow0.
\end{split}
\end{equation}
Then (\ref{Eq:tor.Ext.k.aux}) yields (\ref{Eq:tor.str.k}) with the help of the following observations:
\begin{itemize}
\item $\Tor(\mathbb Z_k,\Cob)=0$; this follows from the fact that for an arbitrary $z\in X$, $\Cob \cong C_z(X,\mathbb T^1)$ is a torsion-free group by con\-nec\-ted\-ness of $X$,
\item $\Tor(\mathbb Z_k,\Cocc)\cong\tork(\Cocc )=\Cocc(\mathbb Z_k)\cong\Hom(\Gamma,\mathbb Z_k)$; the first isomorphism and the equality are clear, the second isomorphism follows from the fact that $X$ is connected and $\mathbb Z_k$ is totally disconnected,
\item $\Tor(\mathbb Z_k,\Coch)\cong\tork(\Coch )$; this is clear,
\item $\mathbb Z_k\otimes\Cob\cong\pi^1(X)/k\pi^1(X)$; this follows from the isomorphisms $\Cob\cong C_z(X,\mathbb T^1)\cong C_z(X,\mathbb R)\oplus\pi^1(X)$ and from the divisibility of $C_z(X,\mathbb R)$:
\begin{equation*}
\begin{split}
\mathbb Z_k\otimes\Cob&\cong\mathbb Z_k\otimes C_z(X,\mathbb T^1)\cong\mathbb Z_k\otimes\big(C_z(X,\mathbb R)\oplus\pi^1(X)\big)\\
&\cong\big(\mathbb Z_k\otimes C_z(X,\mathbb R)\big)\oplus\big(\mathbb Z_k\otimes\pi^1(X)\big)\\
&\cong\mathbb Z_k\otimes\pi^1(X)\cong\pi^1(X)/k\pi^1(X),
\end{split}
\end{equation*}
\item $\mathbb Z_k\otimes\Cocc=0$; this follows from the fact that the group $\Cocc$ is divisible by the assumptions of the theorem.
\end{itemize}

Now consider the exact sequence for tensor and torsion products derived from (\ref{Eq:ex.seq.ext.T1}) and corresponding to the group $\mathbb Q/\mathbb Z$:
\begin{equation}\label{Eq:tor.Ext.k.auxx}
\begin{split}
0&\longrightarrow\Tor(\mathbb Q/\mathbb Z,\Cob)\longrightarrow\Tor(\mathbb Q/\mathbb Z,\Cocc)\longrightarrow\Tor(\mathbb Q/\mathbb Z,\Coch)\\
&\longrightarrow(\mathbb Q/\mathbb Z)\otimes\Cob \longrightarrow (\mathbb Q/\mathbb Z)\otimes\Cocc\longrightarrow(\mathbb Q/\mathbb Z)\otimes\Coch \longrightarrow0.
\end{split}
\end{equation}
Then, similarly to the preceding step of the proof, (\ref{Eq:tor.Ext.k.auxx}) yields (\ref{Eq:tor.str.inf}) with the help of the following observations:
\begin{itemize}
\item $\Tor(\mathbb Q/\mathbb Z,\Cob )=0$; this follows from the fact that the group $\Cob $ is torsion-free by connectedness of $X$,
\item $\Tor(\mathbb Q/\mathbb Z,\Cocc )\cong\bigcup_{l\in\mathbb N}\Hom(\Gamma,\mathbb Z_l)$; this follows from the following isomorphisms by connectedness of $X$:
\begin{equation*}
\begin{split}
\Tor(\mathbb Q/\mathbb Z,\Cocc)&\cong\tor(\Cocc)=\bigcup_{l\in\mathbb N}\text{tor}_l(\Cocc)=\bigcup_{l\in\mathbb N}\Cocc(\mathbb Z_l)\\
&\cong\bigcup_{l\in\mathbb N}\Hom(\Gamma,\mathbb Z_l),
\end{split}
\end{equation*}
\item $\Tor(\mathbb Q/\mathbb Z,\Coch )\cong\tor(\Coch )$; this is clear,
\item $(\mathbb Q/\mathbb Z)\otimes\Cob \cong(\mathbb Q/\mathbb Z)\otimes\pi^1(X)$; this follows from the divisibility of $C_z(X,\mathbb R)$:
\begin{equation*}
\begin{split}
(\mathbb Q/\mathbb Z)\otimes\Cob&\cong(\mathbb Q/\mathbb Z)\otimes C_z(X,\mathbb T^1)\\
&\cong(\mathbb Q/\mathbb Z)\otimes\big(C_z(X,\mathbb R)\oplus\pi^1(X)\big)\\
&\cong\big((\mathbb Q/\mathbb Z)\otimes C_z(X,\mathbb R)\big)\oplus\big((\mathbb Q/\mathbb Z)\otimes\pi^1(X)\big)\\
&\cong(\mathbb Q/\mathbb Z)\otimes\pi^1(X),
\end{split}
\end{equation*}
\item $(\mathbb Q/\mathbb Z)\otimes\Cocc =0$; this follows from the fact that $\Cocc $ is a divisible group by the assumptions of the theorem and $\mathbb Q/\mathbb Z$ is a torsion group.
\end{itemize}

Now assume that $\Gamma$ possesses no non-trivial finite abelian quotient groups. Then $\Hom(\Gamma,\mathbb Z_l)=0$ for every $l\in\mathbb N$ and so the isomorphisms (\ref{Eq:tork.simp.iso}) and (\ref{Eq:tor.simp.iso}) follow from the exactness of the sequences (\ref{Eq:tor.str.k}) and (\ref{Eq:tor.str.inf}), respectively. To finish the proof, we need to verify (\ref{Eq:DEM.in.simp}). The second isomorphism from (\ref{Eq:DEM.in.simp}) is a consequence of (\ref{Eq:tor.simp.iso}). The first isomorphism from (\ref{Eq:DEM.in.simp}) follows from the following facts.
\begin{itemize}
\item The torsion subgroup $\tor(\Coch)$ of $\Coch$ is divisible. Indeed, by the assumptions of the theorem, $\Cocc$ is a divisible group and hence so is its quotient group $\Coch$. Thus, $\tor(\Coch)$ is divisible, being the torsion subgroup of a divisible group.
\item The group $\tor(\Coch)$ is contained as a direct summand in $\Coch$ by divisibility.
\item Every complementary summand $\mathfrak{D}$ of $\tor(\Coch)$ in $\Coch$ is divisible. This follows from the divisibility of both $\Coch$ and $\tor(\Coch)$.
\item Given $\mathfrak{D}$ as above, it is always a subgroup of the groupoid $\Cochm$. Indeed, we have $\mathfrak{D}\setminus1\subseteq\Coch\setminus\tor(\Coch)=\Cochm\setminus 1$.
\end{itemize}
This verifies (\ref{Eq:DEM.in.simp}).
\end{proof}

\begin{corollary}\label{C:F.for.simp.con}
Let $\Flow$ be a minimal flow on a compact connected space $X$. Assume that $\Gamma$ has no non-trivial finite abelian quotient groups and that the group $\Cocc$ is divisible. Then the identity component $F(\Daleth)_0$ of $F(\Daleth)$ is a topological direct summand in $F(\Daleth)$. In fact, under the notation from Theorem~\ref{T:tor.struct}, there are topological isomorphisms
\begin{equation*}
\begin{array}{ccccc}
F(\Daleth) & \cong & F(\Daleth)_0 & \times & F(\Daleth)/F(\Daleth)_0 \\
{} & \cong & \mathfrak{D}^* & \times & \tor\left(\Coch\right)^* \\
{} & \cong & \mathfrak{D}^* & \times & \left((\mathbb Q/\mathbb Z)\otimes\pi^1(X)\right)^*,
\end{array}
\end{equation*}
where all $\mathfrak{D}$, $\tor(\Coch)$ and $(\mathbb Q/\mathbb Z)\otimes\pi^1(X)$ are assumed to carry the discrete topology.
\end{corollary}
\begin{proof}
By virtue of Remark~\ref{R:prop.of.gimel} and by virtue of (\ref{Eq:DEM.in.simp}) from Theorem~\ref{T:tor.struct}, the following topological isomorphisms hold 
for an appropriate divisible subgroup $\mathfrak{D}$ of the groupoid $\Cochm$:
\begin{equation}\label{Eq:F.for.simp.con}
\begin{split}
F(\Daleth)=\left(\Cob\right)^{\perp}\cong\left(\Coch\right)^*&\cong(\mathfrak{D}\oplus\tor(\Coch))^*\\
&\cong\mathfrak{D}^*\times\tor\left(\Coch\right)^*\\
&\cong\mathfrak{D}^*\times\left((\mathbb Q/\mathbb Z)\otimes\pi^1(X)\right)^*.
\end{split}
\end{equation}
Since $\tor(\Coch)$ is a torsion group, its dual group $\tor(\Coch)^*$ is totally disconnected. Further, since $\mathfrak{D}$ is a subgroup of $\Cochm$, it is torsion-free and its dual group $\mathfrak{D}^*$ is therefore connected. Consequently, under the isomorphisms from (\ref{Eq:F.for.simp.con}), the group $\mathfrak{D}^*$ corresponds to $F(\Daleth)_0$ and the group $\tor(\Coch)^*$ corresponds to $F(\Daleth)/F(\Daleth)_0$. This proves the corollary.
\end{proof}

\begin{theorem}\label{T:tor.Lie.mnfld}
Let $\Flow$ be a minimal flow with $\Gamma\in\mathsf{LieGp}$ connected and with $X$ a compact manifold. Set $n=\rank(H_1^w(\mathcal F))$, $n+m=\rank(H_1^w(X))$ and denote by $d_1,\dots,d_n$ the elementary divisors of $H_1^w(\mathcal F)$ in $H_1^w(X)$. If $k\geq2$ is an integer and $\delta_i=\gcd(d_i,k)$ ($i=1,\dots,n$) then
\begin{equation}\label{Eq:ktor.Lie.mfld}
\text{$\tor$}_k\left(\Coch\right)\cong\text{$\tor$}_k\left(H_1^w(X)/H_1^w(\mathcal F)\right)\oplus\left(\mathbb Z_k\right)^m\cong\mathbb Z_{\delta_1}\oplus\dots\oplus\mathbb Z_{\delta_n}\oplus\left(\mathbb Z_k\right)^m
\end{equation}
and
\begin{equation}\label{Eq:tor.Lie.mfld}
\tor\left(\Coch\right)\cong\tor\left(H_1^w(X)/H_1^w(\mathcal F)\right)\oplus\left(\mathbb Q/\mathbb Z\right)^m\cong\mathbb Z_{d_1}\oplus\dots\oplus\mathbb Z_{d_n}\oplus\left(\mathbb Q/\mathbb Z\right)^m.
\end{equation}
\end{theorem}
\begin{remark}\label{R:tor.Lie.mnfld}
We wish to add the following remarks.
\begin{itemize}
\item The isomorphisms (\ref{Eq:ktor.Lie.mfld}) and (\ref{Eq:tor.Lie.mfld}) are consistent with the isomorphisms (\ref{Eq:tork.simp.iso}) and (\ref{Eq:tor.simp.iso}) from Theorem~\ref{T:tor.struct}, respectively. For if $\Gamma$ is simply connected then $\Cocc$ is a divisible group and we have $n=\rank(H_1^w(\mathcal F))=0$, $m=\rank(H_1^w(X))$ and $\pi^1(X)\cong H_1^w(X)$, which leads to isomorphisms
\begin{equation*}
\begin{split}
\tork(\Coch)&\cong\pi^1(X)/k\pi^1(X)\cong H_1^w(X)/kH_1^w(X)\\
&\cong\mathbb Z^m/k\mathbb Z^m\cong\left(\mathbb Z_k\right)^m
\end{split}
\end{equation*}
and
\begin{equation*}
\begin{split}
\tor(\Coch)&\cong(\mathbb Q/\mathbb Z)\otimes\pi^1(X)\cong(\mathbb Q/\mathbb Z)\otimes H_1^w(X)\\
&\cong(\mathbb Q/\mathbb Z)\otimes\mathbb Z^m\cong(\mathbb Q/\mathbb Z)^m.
\end{split}
\end{equation*}
\item The isomorphism (\ref{Eq:ktor.Lie.mfld}) can help us understand the structure of the set $\per_k(\Coch)$ of all cohomology classes of extensions $\mathcal C\in\Cocc$ with $F(\mathcal C)=\mathbb Z_k$. Indeed, an element $(z_1,\dots,z_n,\zeta_1,\dots,\zeta_m)$ of $\mathbb Z_{\delta_1}\oplus\dots\oplus\mathbb Z_{\delta_n}\oplus(\mathbb Z_k)^m$ corresponds to an element of $\per_k(\Coch)$ under the isomorphism (\ref{Eq:ktor.Lie.mfld}) if and only if one of the following conditions is satisfied:
\begin{enumerate}
\item[(i)] there is $1\leq i\leq n$ such that $k$ divides $d_i$ (that is, $\delta_i=k$) and $z_i$ generates $\mathbb Z_{\delta_i}=\mathbb Z_k$,
\item[(ii)] there is $1\leq j\leq m$ such that $\zeta_j$ generates $\mathbb Z_k$.
\end{enumerate}
\end{itemize}
\end{remark}
\begin{proof}[Proof of Theorem~\ref{T:tor.Lie.mnfld}]
Throughout the whole proof, both $X$ and $\Gamma$ will be assumed pointed. The base point of $X$ will be denoted by $z$ and the base point of $\Gamma$ will be its identity $1$. We also recall that for abelian groups $A,C$ and their subgroups $B\sbgp A$, $D\sbgp C$, $\Hom(A,B;C,D)$ stands for the group of all morphisms $h\colon A\to C$ with $h(B)\subseteq D$. Now we turn to the proof, proceeding in three steps.

\emph{1st step.} Fix a map $\xi\in C_z(X,\mathbb T^1)$ and an integer $k\geq2$. We show that the following conditions are equivalent:
\begin{enumerate}
\item[(a)] $\co(\xi)$ lifts across $\kappa_k\colon\mathbb T^1\to\mathbb T^1$ to an extension $\mathcal C\in\Cocc$,
\item[(b)] $\xi^{\sharp}H_1^w(\mathcal F)\subseteq k\mathbb Z$, where $\xi^{\sharp}$ is the induced morphism $H_1^w(X)\to H_1^w(\mathbb T^1)$ and $H_1^w(\mathbb T^1)$ is identified with $\mathbb Z$.
\end{enumerate}

We claim that the following conditions are equivalent:
\begin{itemize}
\item $\co(\xi)$ lifts across $\kappa_k$ to an extension $\mathcal C\in\Cocc$,
\item $\co(\xi)$ lifts across $\kappa_k$ to a continuous base point preserving map $\Gamma\times X\to\mathbb T^1$,
\item $\co(\xi)^{\sharp}\pi_1(\Gamma\times X)\subseteq \kappa_k^{\sharp}\pi_1(\mathbb T^1)=k\mathbb Z$,
\item $\xi^{\sharp}\mathcal F_z^{\sharp}\pi_1(\Gamma)\subseteq k\mathbb Z$,
\item $\xi^{\sharp}H_1^w(\mathcal F)\subseteq k\mathbb Z$.
\end{itemize}
Indeed, the equivalence of the first four of these conditions follows from Lemmas~\ref{P:lifting.cocycle} and~\ref{L:tor.Lie.mnfld}. The equivalence of the last two conditions follows from the following simple computation based on the commutativity of the diagram in Lemma~\ref{L:tor.Lie.mnfld}
\begin{equation*}
p_{\mathbb T^1}\left(\xi^{\sharp}\mathcal F_z^{\sharp}\pi_1(\Gamma)\right)=
\xi^{\sharp}\mathcal F_z^{\sharp}p_{\Gamma}(\pi_1(\Gamma))=\xi^{\sharp}
\mathcal F_z^{\sharp}H_1^w(\Gamma)=\xi^{\sharp}H_1^w(\mathcal F)
\end{equation*}
and from the equality $p_{\mathbb T^1}=\id_{\mathbb Z}$, which holds under the standard identifications $\pi_1(\mathbb T^1)\cong H_1^w(\mathbb T^1)\cong\mathbb Z$.

\emph{2nd step.} We verify (\ref{Eq:ktor.Lie.mfld}).

Given a propositional form $\mathcal P$ defined on the set $C_z(X,\mathbb T^1)$, we use the symbol $[\mathcal P(\xi)]$ to denote the set of all $\xi\in C_z(X,\mathbb T^1)$ such that the proposition $\mathcal P(\xi)$ holds true. We shall arrive at (\ref{Eq:ktor.Lie.mfld}) through a series of isomorphisms. Each of these isomorphisms will be dealt with in a separate item.
\begin{enumerate}
\item[(i)] By virtue of the equivalence (a)$\Leftrightarrow$(b) from the first step of the proof,
\begin{equation}\label{Eq:arriv.at.tork.1}
\begin{split}
\tork(\Coch)&=\left\{\mathcal C\in\Cocc : \mathcal C^k\in\Cob\right\}/\Cob\\
&=\left\{\kappa_k^{-1}(\co(\xi)) : \xi^{\sharp}H_1^w(\mathcal F)\subseteq k\mathbb Z\right\}/\Cob,
\end{split}
\end{equation} 
where $\kappa_k$ stands for the $k$-endomorphism of $\Cocc$. Observe that $\kappa_k$ is a monomorphism by connectedness of $\Gamma$ and $X$. Consequently, by virtue of (\ref{Eq:arriv.at.tork.1}),
\begin{equation}\label{Eq:arriv.at.tork.2}
\tork(\Coch)\cong\left\{\co(\xi) : \xi^{\sharp}H_1^w(\mathcal F)\subseteq k\mathbb Z\right\}/\kappa_k(\Cob).
\end{equation}
\item[(ii)] Under the usual isomorphism $\sigma\colon C_z(X,\mathbb T^1)\ni\xi\mapsto\co(\xi)\in\Cob$, the subgroup $\kappa_k(\Cob)$ of $\Cob$ corresponds to the subgroup $\kappa_k(C_z(X,\mathbb T^1))$ of $C_z(X,\mathbb T^1)$. The latter subgroup consists of all those maps $\xi\in C_z(X,\mathbb T^1)$, which lift across $\kappa_k\colon\mathbb T^1\to\mathbb T^1$ to an element of $C_z(X,\mathbb T^1)$; this occurs if and only if $\xi^{\sharp}H_1^w(X)\subseteq k\mathbb Z$. Thus, the isomorphism (\ref{Eq:arriv.at.tork.2}) yields
\begin{equation}\label{Eq:arriv.at.tork.3}
\tork(\Coch)\cong\frac{[\xi^{\sharp}H_1^w(\mathcal F)\subseteq k\mathbb Z]}{[\xi^{\sharp}H_1^w(X)\subseteq k\mathbb Z]}\cong\frac{[\xi^{\sharp}H_1^w(\mathcal F)\subseteq k\mathbb Z]/[\xi^{\sharp}H_1^w(X)=0]}{[\xi^{\sharp}H_1^w(X)\subseteq k\mathbb Z]/[\xi^{\sharp}H_1^w(X)=0]}.
\end{equation}
\item[(iii)] Consider the map $\varphi\colon C_z(X,\mathbb T^1)\ni\xi\mapsto\xi^{\sharp}\in\Hom(H_1^w(X),\mathbb Z)$. By our assumptions on $X$, $\varphi$ is an epimorphism and $\ker(\varphi)=[\xi^{\sharp}H_1^w(X)=0]$. Consequently, it follows from (\ref{Eq:arriv.at.tork.3}) by elementary group-theoretic arguments that
\begin{equation}\label{Eq:arriv.at.tork.4}
\begin{split}
\tork(\Coch)&\cong\frac{\Hom(H_1^w(X),H_1^w(\mathcal F);\mathbb Z,k\mathbb Z)}{\Hom(H_1^w(X),H_1^w(X);\mathbb Z,k\mathbb Z)}\\
&\cong\frac{\Hom(\mathbb Z^n\oplus\mathbb Z^m,d_1\mathbb Z\oplus\dots\oplus d_n\mathbb Z;\mathbb Z,k\mathbb Z)}{\Hom(\mathbb Z^n\oplus\mathbb Z^m,\mathbb Z^n\oplus\mathbb Z^m;\mathbb Z,k\mathbb Z)},
\end{split}
\end{equation}
where $d_1\mathbb Z\oplus\dots\oplus d_n\mathbb Z\subseteq\mathbb Z^n\cong\mathbb Z^n\oplus0\subseteq\mathbb Z^n\oplus\mathbb Z^m$.
\item[(iv)] Consider the isomorphism 
\begin{equation*}
\psi\colon\Hom(\mathbb Z^n\oplus\mathbb Z^m,\mathbb Z)\ni h\mapsto(h(e_1),\dots,h(e_{n+m}))\in\mathbb Z^n\oplus\mathbb Z^m,  
\end{equation*}
where $e_1,\dots,e_{n+m}$ is the standard basis for $\mathbb Z^n\oplus\mathbb Z^m$. Then
\begin{equation*}
\psi(\Hom(\mathbb Z^n\oplus\mathbb Z^m,d_1\mathbb Z\oplus\dots\oplus d_n\mathbb Z;\mathbb Z,k\mathbb Z))=\frac{k}{\delta_1}\mathbb Z\oplus\dots\oplus\frac{k}{\delta_n}\mathbb Z\oplus\mathbb Z^m
\end{equation*}
and
\begin{equation*}
\psi(\Hom(\mathbb Z^n\oplus\mathbb Z^m,\mathbb Z^n\oplus\mathbb Z^m;\mathbb Z,k\mathbb Z))=k\mathbb Z^n\oplus k\mathbb Z^m.
\end{equation*}
Consequently, by virtue of (\ref{Eq:arriv.at.tork.4}),
\begin{equation}\label{Eq:arriv.at.tork.5}
\begin{split}
\tork(\Coch)&\cong\left(\frac{k}{\delta_1}\mathbb Z\oplus\dots\oplus\frac{k}{\delta_n}\mathbb Z\oplus\mathbb Z^m\right)/\left(k\mathbb Z^n\oplus k\mathbb Z^m\right)\\
&\cong\mathbb Z_{\delta_1}\oplus\dots\oplus\mathbb Z_{\delta_n}\oplus\left(\mathbb Z_k\right)^m\\
&\cong\tork(\mathbb Z_{d_1}\oplus\dots\oplus \mathbb Z_{d_n})\oplus\left(\mathbb Z_k\right)^m\\
&\cong\tork(H_1^w(X)/H_1^w(\mathcal F))\oplus\left(\mathbb Z_k\right)^m.
\end{split}
\end{equation}
This verifies the two isomorphisms from (\ref{Eq:ktor.Lie.mfld}) and finishes the second step of the proof.
\end{enumerate}
\begin{figure}[ht]
\[\minCDarrowwidth20pt\begin{CD}
\kappa_k^{-1}(\Cob) @>\sigma^{-1}\kappa_k>> [\xi^{\sharp}H_1^w(\mathcal F)\subseteq k\mathbb Z] @>\psi\varphi>> \frac{k}{\delta_1}\mathbb Z\oplus\dots\oplus\frac{k}{\delta_n}\mathbb Z\oplus\mathbb Z^m\\
{\rotatebox{90}{$\subseteq$}} @. {\rotatebox{90}{$\subseteq$}} @. {\rotatebox{90}{$\subseteq$}}\\
\Cob @>>\sigma^{-1}\kappa_k> [\xi^{\sharp}H_1^w(X)\subseteq k\mathbb Z] @>>\psi\varphi> k\mathbb Z^n\oplus k\mathbb Z^m
\end{CD}\]
\caption{Verifying (\ref{Eq:ktor.Lie.mfld})}
\label{Fig:tor.Lie.mnf.L2}
\end{figure}
Since we shall have to refer to the procedure just finished and also modify it in a certain way, we summarize its steps in the form of a diagram in Figure~\ref{Fig:tor.Lie.mnf.L2}. (We emphasize that all the arrows in the diagram are epimorphisms.)

\emph{3rd step.} We verify (\ref{Eq:tor.Lie.mfld}).

Let $k,l\geq2$ be integers, $l$ divisible by $k$ and $k$ divisible by $d=d_n$. Then $\gcd(d_i,k)=\gcd(d_i,l)=d_i$ for $i=1,\dots,n$ and hence, by the second step of the proof, there are isomorphisms
\begin{equation}\label{Eq:tork.strt}
\text{tor}_k(\Coch)\cong\mathbb Z_{d_1}\oplus\dots\oplus\mathbb Z_{d_n}\oplus\left(\mathbb Z_k\right)^m
\end{equation}
and
\begin{equation}\label{Eq:torl.strt}
\text{tor}_l(\Coch)\cong\mathbb Z_{d_1}\oplus\dots\oplus\mathbb Z_{d_n}\oplus\left(\mathbb Z_l\right)^m.
\end{equation}
We need to determine what form the inclusion $\tor_k(\Coch)\subseteq\tor_l(\Coch)$ takes after passing to the right hand sides of (\ref{Eq:tork.strt}) and (\ref{Eq:torl.strt}). To this end, we invoke the procedure from the second step of the proof, which is summarized in the diagram in Figure~\ref{Fig:tor.Lie.mnf.L2}, but we apply the monomorphism $\kappa_l$ instead of $\kappa_k$. This leads to the inclusions depicted in the diagram in Figure~\ref{Fig:tor.Lie.mnf.L3}, as we now show. (Similarly to Figure~\ref{Fig:tor.Lie.mnf.L2}, all the arrows in Figure~\ref{Fig:tor.Lie.mnf.L3} are claimed to be epimorphisms.) 
\begin{figure}[ht]
\[\minCDarrowwidth20pt\begin{CD}
\kappa_k^{-1}(\Cob) @>\sigma^{-1}\kappa_l>> \kappa_{l/k}[\xi^{\sharp}H_1^w(\mathcal F)\subseteq k\mathbb Z] @>\psi\varphi>> \frac{l}{d_1}\mathbb Z\oplus\dots\oplus\frac{l}{d_n}\mathbb Z\oplus\frac{l}{k}\mathbb Z^m\\
{\rotatebox{90}{$\subseteq$}} @. {\rotatebox{90}{$\subseteq$}} @. {\rotatebox{90}{$\subseteq$}}\\
\Cob @>\sigma^{-1}\kappa_l>> [\xi^{\sharp}H_1^w(X)\subseteq l\mathbb Z] @>\psi\varphi>> l\mathbb Z^n\oplus l\mathbb Z^m\\
{\rotatebox{90}{$\supseteq$}} @. {\rotatebox{90}{$\supseteq$}} @. {\rotatebox{90}{$\supseteq$}}\\
\kappa_l^{-1}(\Cob) @>\sigma^{-1}\kappa_l>> [\xi^{\sharp}H_1^w(\mathcal F)\subseteq l\mathbb Z] @>\psi\varphi>> \frac{l}{d_1}\mathbb Z\oplus\dots\oplus\frac{l}{d_n}\mathbb Z\oplus\mathbb Z^m
\end{CD}\]
\caption{Towards the inclusion $\tor_k(\Coch)\subseteq\tor_l(\Coch)$}
\label{Fig:tor.Lie.mnf.L3}
\end{figure}
The second and the third line of the diagram follow immediately from the second step of the proof. We verify the first line, proceeding in two steps.
\begin{itemize}
\item Recall from the first step of the proof that
\begin{equation*}
\kappa_k^{-1}(\Cob)=\left\{\kappa_k^{-1}(\co(\xi)) : \xi^{\sharp}H_1^w(\mathcal F)\subseteq k\mathbb Z\right\}.
\end{equation*}
By applying $\sigma^{-1}\kappa_l$ to this equality, we obtain
\begin{equation*}
\begin{split}
\sigma^{-1}\kappa_l\left(\kappa_k^{-1}(\Cob)\right)&=\sigma^{-1}\left(\left\{\co\left(\xi^{l/k}\right) : \xi^{\sharp}H_1^w(\mathcal F)\subseteq k\mathbb Z\right\}\right)\\
&=\left\{\xi^{l/k} : \xi^{\sharp}H_1^w(\mathcal F)\subseteq k\mathbb Z\right\}\\
&=\kappa_{l/k}[\xi^{\sharp}H_1^w(\mathcal F)\subseteq k\mathbb Z].
\end{split}
\end{equation*}
\item By the second step of the proof,
\begin{equation*}
\begin{split}
\psi\varphi[\xi^{\sharp}H_1^w(\mathcal F)\subseteq k\mathbb Z]&=\frac{k}{\delta_1}\mathbb Z\oplus\dots\oplus\frac{k}{\delta_n}\mathbb Z\oplus\mathbb Z^m\\
&=\frac{k}{d_1}\mathbb Z\oplus\dots\oplus\frac{k}{d_n}\mathbb Z\oplus\mathbb Z^m.
\end{split}
\end{equation*}
Consequently,
\begin{equation*}
\begin{split}
\psi\varphi\left(\kappa_{l/k}[\xi^{\sharp}H_1^w(\mathcal F)\subseteq k\mathbb Z]\right)&=\kappa_{l/k}\left(\psi\varphi[\xi^{\sharp}H_1^w(\mathcal F)\subseteq k\mathbb Z]\right)\\
&=\kappa_{l/k}\left(\frac{k}{d_1}\mathbb Z\oplus\dots\oplus\frac{k}{d_n}\mathbb Z\oplus\mathbb Z^m\right)\\
&=\frac{l}{d_1}\mathbb Z\oplus\dots\oplus\frac{l}{d_n}\mathbb Z\oplus\frac{l}{k}\mathbb Z^m.
\end{split}
\end{equation*}
\end{itemize}
These two observations verify the first line of the diagram in Figure~\ref{Fig:tor.Lie.mnf.L3}.

Set $\varrho=\psi\varphi\sigma^{-1}\kappa_l\colon\kappa_l^{-1}(\Cob)\to\mathbb Z^n\oplus\mathbb Z^m$. Following the results summarized in Fi\-gure~\ref{Fig:tor.Lie.mnf.L3}, we claim that $\varrho$ induces isomorphisms
\begin{equation}\label{Eq:Cob.ro.kl.1}
\begin{split}
\text{tor}_k(\Coch)&=\frac{\kappa_k^{-1}(\Cob)}{\Cob}\cong\frac{\varrho(\kappa_k^{-1}(\Cob))}{\varrho(\Cob)}\\
&=\left(\frac{l}{d_1}\mathbb Z\oplus\dots\oplus\frac{l}{d_n}\mathbb Z\oplus\frac{l}{k}\mathbb Z^m\right)/(l\mathbb Z^n\oplus l\mathbb Z^m)
\end{split}
\end{equation}
and
\begin{equation}\label{Eq:Cob.ro.kl.2}
\begin{split}
\text{tor}_l(\Coch)&=\frac{\kappa_l^{-1}(\Cob)}{\Cob}\cong\frac{\varrho(\kappa_l^{-1}(\Cob))}{\varrho(\Cob)}\\
&=\left(\frac{l}{d_1}\mathbb Z\oplus\dots\oplus\frac{l}{d_n}\mathbb Z\oplus\mathbb Z^m\right)/(l\mathbb Z^n\oplus l\mathbb Z^m).
\end{split}
\end{equation}
In order to prove this claim, it suffices to show that $\ker(\varrho)\subseteq\Cob$. So let $\mathcal C\in\kappa_l^{-1}(\Cob)$, $\xi\in C_z(X,\mathbb T^1)$ be the transfer function of $\kappa_l(\mathcal C)$ and assume that $\varrho(\mathcal C)=0$. Then $\xi^{\sharp}=0$ and hence $\xi\in\Div(C_z(X,\mathbb T^1))$. Since $\sigma\colon C_z(X,\mathbb T^1)\to\Cob$ is an isomorphism, it follows that $\co(\xi)=\sigma(\xi)\in\Div(\Cob)$. Consequently, $\mathcal C=\kappa_l^{-1}(\co(\xi))\in\Cob$, as was to be shown.

Now, under the natural isomorphisms $\mathbb Z/l\mathbb Z\cong\tor_l(\mathbb T^1)$ and $((l/k)\mathbb Z)/l\mathbb Z\cong\mathbb Z/k\mathbb Z\cong\tor_k(\mathbb T^1)$, the inclusion $((l/k)\mathbb Z)/l\mathbb Z\subseteq\mathbb Z/l\mathbb Z$ takes the form of the inclusion $\tor_k(\mathbb T^1)\subseteq\tor_l(\mathbb T^1)$. Thus, by virtue of (\ref{Eq:Cob.ro.kl.1}) and (\ref{Eq:Cob.ro.kl.2}), the inclusion $\tor_k(\Coch)\subseteq\tor_l(\Coch)$ takes the form of the right hand side of Figure~\ref{Fig:tor.Lie.mnf.L4}.
\begin{figure}[ht]
\[\minCDarrowwidth20pt\begin{CD}
\tor_k(\Coch) @. {}={} @. \kappa_k^{-1}(\Cob)/\Cob @. {}\cong{} @. \mathbb Z_{d_1}\oplus\dots\oplus\mathbb Z_{d_n}\oplus\tor_k(\mathbb T^1)^m\\
{\rotatebox{90}{$\supseteq$}} @. @. {\rotatebox{90}{$\supseteq$}} @. @. {\rotatebox{90}{$\supseteq$}}\\
\tor_l(\Coch) @. {}={} @. \kappa_l^{-1}(\Cob)/\Cob @. {}\cong{} @. \mathbb Z_{d_1}\oplus\dots\oplus\mathbb Z_{d_n}\oplus\tor_l(\mathbb T^1)^m
\end{CD}\]
\caption{Another form of the inclusion $\tor_k(\Coch)\subseteq\tor_l(\Coch)$}
\label{Fig:tor.Lie.mnf.L4}
\end{figure}

Let $(k_j)_{j\in\mathbb N}$ be a sequence of positive integers with the following properties:
\begin{itemize}
\item $d$ divides $k_1$ and $k_j$ divides $k_{j+1}$ for every $j\in\mathbb N$,
\item for every $k\in\mathbb N$ there is $j\in\mathbb N$ such that $k$ divides $k_j$.
\end{itemize}
Then the results of this step of the proof yield isomorphisms
\begin{equation*}
\begin{split}
\tor(\Coch)&=\bigcup_{k\in\mathbb N}\text{$\tor$}_k(\Coch)=\bigcup_{j\in\mathbb N}\text{$\tor$}_{k_j}(\Coch)\cong\bigcup_{j\in\mathbb N}\left(\mathbb Z_{d_1}\oplus\dots\oplus\mathbb Z_{d_n}\oplus\text{tor}_{k_j}(\mathbb T^1)^m\right)\\
&=\mathbb Z_{d_1}\oplus\dots\oplus\mathbb Z_{d_n}\oplus\tor(\mathbb T^1)^m=\mathbb Z_{d_1}\oplus\dots\oplus\mathbb Z_{d_n}\oplus\left(\mathbb Q/\mathbb Z\right)^m\\
&\cong\tor\left(H_1^w(X)/H_1^w(\mathcal F)\right)\oplus(\mathbb Q/\mathbb Z)^m.
\end{split}
\end{equation*}
This verifies (\ref{Eq:tor.Lie.mfld}).
\end{proof}

\section{Torsions over topologically free flows}\label{S:tor.top.free.flw}

In this section we continue our study of the cohomology group $\Coch$, focusing now on the situation when $\Flow$ is a topologically free minimal flow with $\Gamma\in\mathsf{LieGp}$ connected and with $X$ a compact (connected) manifold. We start by determining the complementary summand of $\Div(\Cocc)$ in $\Cocc$ and that of $\Div(\Coch)$ in $\Coch$. Then we express $\Coch$ as a direct sum of its torsion subgroup $\tor(\Coch)$ and a divisible subgroup $\mathfrak{D}$ of the groupoid $\Cochm$.

\begin{theorem}\label{T:ex.seq.Lie.mfld}
Let $\Flow$ be a minimal topologically free flow with $\Gamma\in\mathsf{LieGp}$ connected and with $X$ a compact manifold. Set $n=\rank(H_1^w(\mathcal F))$, $n+m=\rank(H_1^w(X))$ and denote by $d_1,\dots,d_n$ the elementary divisors of $H_1^w(\mathcal F)$ in $H_1^w(X)$. Then there are direct sums
\begin{equation}\label{Eq:Cocc.div.sum}
\Cocc=\Div(\Cocc)\oplus H_1^w(\mathcal F)\cong\Div(\Cocc)\oplus\mathbb Z^n
\end{equation}
and
\begin{equation}\label{Eq:Coch.div.sum}
\Coch=\Div(\Coch)\oplus\tor\Big(H_1^w(X)/H_1^w(\mathcal F)\Big)\cong\Div(\Coch)\oplus\mathbb Z_{d_1}\oplus\dots\oplus\mathbb Z_{d_n}.
\end{equation}
Moreover, for every $G\in\mathsf{CAGp}$, the exact sequence $(\ref{Eq:exact.sequence})$ from Theorem~\ref{T:exact.sequence} takes the form
\begin{equation}\label{Eq:Hom.Ext.Lie.mfld.2}
\begin{split}
0&\longrightarrow\Cob(G)\stackrel{\mu_G}{\longrightarrow}\Cocc(G)\stackrel{\kappa_G}{\longrightarrow}\Hom(G^*,\Coch)\longrightarrow\\
&\longrightarrow\Ext(G^*,\mathbb Z)^{n+m}\longrightarrow\Ext(G^*,\mathbb Z)^n\longrightarrow\bigoplus_{i=1}^n\Ext(G^*,\mathbb Z_{d_i})\longrightarrow0.
\end{split}
\end{equation}
\end{theorem}
\begin{proof}
Before turning to the proof we fix some notation. By definition of $n$, $m$ and $d_1,\dots,d_n$, the inclusion $H_1^w(\mathcal F)\subseteq H_1^w(X)$ can be expressed in the form $d_1\mathbb Z\oplus\dots\oplus d_n\mathbb Z\subseteq\mathbb Z^n\cong\mathbb Z^n\oplus0\subseteq\mathbb Z^n\oplus\mathbb Z^m$. The standard basis of $\mathbb Z^n\oplus\mathbb Z^m\cong\mathbb Z^{n+m}$ will be denoted by $e_1,\dots,e_{n+m}$ and we shall write $\eta=d_1\dots d_n$. We assume throughout the whole proof that the space $X$ is pointed and we denote its base point by $z$. Recall the notation $\mathcal F_z=\mathcal F(-,z)\colon\Gamma\to X$ and $\mathcal C_z=\mathcal C(-,z)\colon\Gamma\to\mathbb T^1$ for every $\mathcal C\in\Cocc$. We shall use freely the usual identifications $\pi_1(\mathbb T^1)\cong H_1^w(\mathbb T^1)\cong\mathbb Z$. Finally, observe that by our assumptions on $X$, for every $h\in\Hom(H_1^w(X),\mathbb Z)$ there is $\xi\in C_z(X,\mathbb T^1)$ such that $\xi^{\sharp}=h$.

\emph{1st step.} We collect some useful facts about lifts of extensions from $\Cocc$ and about the divisible subgroup $\Div(\Cocc)$ of $\Cocc$.

Let $\mathcal C\in\Cocc$, $\xi\in C_z(X,\mathbb T^1)$ and $l\in\mathbb N$. Then the following statements hold:
\begin{enumerate}
\item[(a)] we have
\begin{equation*}
\mathcal C^{\sharp}\pi_1(\Gamma\times X)=\mathcal C_z^{\sharp}H_1^w(\Gamma)\subseteq\mathbb Z;
\end{equation*}
this follows from Lemma~\ref{L:tor.Lie.mnfld},
\item[(b)] we have
\begin{equation*}
\co(\xi)^{\sharp}\pi_1(\Gamma\times X)=\xi^{\sharp}\mathcal F_z^{\sharp}H_1^w(\Gamma)=\xi^{\sharp}H_1^w(\mathcal F)\subseteq\mathbb Z;
\end{equation*}
this follows from Lemma~\ref{L:tor.Lie.mnfld},
\item[(c)] $\mathcal C$ lifts across $\kappa_l\colon\mathbb T^1\to\mathbb T^1$ to an element of $\Cocc$ if and only if $\mathcal C_z^{\sharp}H_1^w(\Gamma)\subseteq l\mathbb Z$ (in particular, $\mathcal C\in\Div(\Cocc)$ if and only if $\mathcal C_z^{\sharp}H_1^w(\Gamma)=0$); this follows from Lemma~\ref{P:lifting.cocycle} and from (a),
\item[(d)] $\co(\xi)$ lifts across $\kappa_l\colon\mathbb T^1\to\mathbb T^1$ to an element of $\Cocc$ if and only if $\xi^{\sharp}H_1^w(\mathcal F)\subseteq l\mathbb Z$ (in particular, $\co(\xi)\in\Div(\Cocc)$ if and only if $\xi^{\sharp}H_1^w(\mathcal F)=0$); this follows from Lemma~\ref{P:lifting.cocycle} and from (b).
\end{enumerate}

\emph{2nd step.} For $i=1,\dots,n$ fix $\xi_i\in C_z(X,\mathbb T^1)$ so that $(\xi_i)^{\sharp}(e_j)=\delta_{ij}\eta/d_j$ for $j=1,\dots,n$ and $(\xi_i)^{\sharp}(e_j)=0$ for $j=n+1,\dots,n+m$, and let $\mathcal C_i\in\Cocc$ be the lift of $\co(\xi_i)$ across $\kappa_{\eta}\colon\mathbb T^1\to\mathbb T^1$. We claim that $\mathcal C_1,\dots,\mathcal C_n$ form a basis for a free abelian subgroup $\mathcal S$ of $\Cocc$ and that there is a direct sum $\Cocc=\Div(\Cocc)\oplus\mathcal S$. This will verify (\ref{Eq:Cocc.div.sum}).

First notice that by virtue of (d) from the first step of the proof, all $\mathcal C_i$ are well defined. Now, in order to prove the two claims, it suffices to verify the following statements.
\begin{enumerate}
\item[($\alpha$)] For all $k_1,\dots,k_n\in\mathbb Z$, $\sum_{i=1}^nk_i\mathcal C_i\in\Div(\Cocc)$ implies $k_1=\dots=k_n=0$.
\item[($\beta$)] The set $\Div(\Cocc)\cup\mathcal S$ generates the group $\Cocc$.
\end{enumerate}

We begin be verifying ($\alpha$). So fix $k_1,\dots,k_n\in\mathbb Z$ with $\sum_{i=1}^nk_i\mathcal C_i\in\Div(\Cocc)$. Then 
\begin{equation*}
\co\left(\sum_{i=1}^nk_i\xi_i\right)=\sum_{i=1}^nk_i\co(\xi_i)=\eta\left(\sum_{i=1}^nk_i\mathcal C_i\right)\in\Div(\Cocc)
\end{equation*}
and hence, by virtue of (d), $\sum_{i=1}^nk_i(\xi_i)^{\sharp}=(\sum_{i=1}^nk_i\xi_i)^{\sharp}$ vanishes on $H_1^w(\mathcal F)$. By our choice of the maps $\xi_i$, this leads to $k_1=\dots=k_n=0$.

To verify ($\beta$), fix $\mathcal C\in\Cocc$. Since the flow $\mathcal F$ is topologically free by the assumptions of the theorem, the morphism $\mathcal F_z^{\sharp}\colon H_1^w(\Gamma)\to H_1^w(\mathcal F)$ is an isomorphism. Thus, since $\eta\mathcal C_z^{\sharp}=(\eta\mathcal C)_z^{\sharp}\colon H_1^w(\Gamma)\to\mathbb Z$ takes its values in $\eta\mathbb Z$, it follows that so does $(\eta\mathcal C)_z^{\sharp}(\mathcal F_z^{\sharp})^{-1}\colon H_1^w(\mathcal F)\to\mathbb Z$. Consequently, $(\eta\mathcal C)_z^{\sharp}(\mathcal F_z^{\sharp})^{-1}$ extends to a morphism $\varrho\colon H_1^w(X)\to\mathbb Z$ with $\varrho(e_j)\in(\eta/d_j)\mathbb Z$ for $j=1,\dots,n$ and $\varrho(e_j)=0$ for $j=n+1,\dots,n+m$. By our choice of the maps $\xi_i$, there are integers $k_1,\dots,k_n$ with $\varrho=\sum_{i=1}^nk_i(\xi_i)^{\sharp}=(\sum_{i=1}^nk_i\xi_i)^{\sharp}$. Thus, with the help of Lemma~\ref{L:tor.Lie.mnfld}, we obtain
\begin{equation*}
\eta\mathcal C_z^{\sharp}=(\eta\mathcal C)_z^{\sharp}=\varrho\mathcal F_z^{\sharp}=\Big(\sum_{i=1}^nk_i\xi_i\Big)^{\sharp}\mathcal F_z^{\sharp}=\co\Big(\sum_{i=1}^nk_i\xi_i\Big)_z^{\sharp}=\eta\Big(\sum_{i=1}^nk_i\mathcal C_i\Big)_z^{\sharp},
\end{equation*}
whence it follows that
\begin{equation*}
\Big(\mathcal C-\sum_{i=1}^nk_i\mathcal C_i\Big)_z^{\sharp}=\mathcal C_z^{\sharp}-\Big(\sum_{i=1}^nk_i\mathcal C_i\Big)_z^{\sharp}=0.
\end{equation*}
In view of (c) from the first step of the proof, this means that $\mathcal C-\sum_{i=1}^nk_i\mathcal C_i\in\Div(\Cocc)$, which verifies ($\beta$).

\emph{3rd step.} We show that $\mathcal S\cap\Cob=\langle d_1\mathcal C_1,\dots,d_n\mathcal C_n\rangle$.

First we verify inclusion ``$\subseteq$''. So fix $k_1,\dots,k_n\in\mathbb Z$ and assume that $\sum_{i=1}^nk_i\mathcal C_i\in\Cob$; we show that $d_i$ divides $k_i$ for every $i=1,\dots,n$. Then 
\begin{equation*}
\eta\left(\sum_{i=1}^nk_i\mathcal C_i\right)=\sum_{i=1}^nk_i\co(\xi_i)=\co\left(\sum_{i=1}^nk_i\xi_i\right)
\end{equation*}
and hence $\co(\sum_{i=1}^nk_i\xi_i)$ lifts across $\kappa_{\eta}\colon\mathbb T^1\to\mathbb T^1$ to an element of $\Cob$. Consequently, $\sum_{i=1}^nk_i\xi_i$ lifts across $\kappa_{\eta}$ to an element of $C_z(X,\mathbb T^1)$, whence it follows that $(\sum_{i=1}^nk_i\xi_i)^{\sharp}=\sum_{i=1}^nk_i(\xi_i)^{\sharp}$ takes its values in $\eta\mathbb Z$. By our choice of the maps $\xi_i$ it follows that $d_i$ divides $k_i$ for every $i=1,\dots,n$.

To verify inclusion ``$\supseteq$'', it suffices to show that $d_i\mathcal C_i\in\Cob$ for every $i=1,\dots,n$. By our choice of the maps $\xi_i$, $(d_i\xi_i)^{\sharp}(H_1^w(X))=d_i(\xi_i)^{\sharp}(H_1^w(X))\subseteq\eta\mathbb Z$ for every $i=1,\dots,n$. Consequently, $d_i\xi_i$ lifts across $\kappa_{\eta}$ to an element of $C_z(X,\mathbb T^1)$ and hence $\co(d_i\xi_i)$ lifts across $\kappa_{\eta}$ to an element of $\Cob$. Since $\kappa_{\eta}(d_i\mathcal C_i)=d_i(\kappa_{\eta}\mathcal C_i)=d_i\co(\xi_i)=\co(d_i\xi_i)$, it follows by uniqueness of a lift that, indeed, $d_i\mathcal C_i\in\Cob$.

\emph{4th step.} Let $\mathcal M$ consist of the maps $\xi\in C_z(X,\mathbb T^1)$ with $\xi^{\sharp}H_1^w(\mathcal F)=0$. Clearly, $\mathcal M$ is a subgroup of $C_z(X,\mathbb T^1)$ and hence $\sigma(\mathcal M)$ is a subgroup of $\Cob$, where $\sigma$ is the usual isomorphism $C_z(X,\mathbb T^1)\to\Cob$. We show that $\Cob=\sigma(\mathcal M)\oplus(\mathcal S\cap\Cob)=\sigma(\mathcal M)\oplus\langle d_1\mathcal C_1,\dots,d_n\mathcal C_n\rangle$.

To verify the desired direct sum, we need to check the following statements:
\begin{enumerate}
\item[(i)] $\sigma(\mathcal M)\cap(\mathcal S\cap\Cob)=0$,
\item[(ii)] $\Cob=\sigma(\mathcal M)+(\mathcal S\cap\Cob)$.
\end{enumerate}

We begin with statement (i). So fix $\xi\in\mathcal M$ and $k_1,\dots,k_n\in\mathbb Z$ with $\co(\xi)=\sum_{i=1}^nk_i(d_i\mathcal C_i)$; we show that $k_1=\dots=k_n=0$ and $\xi=0$. From the last equality we obtain
\begin{equation*}
\co(\eta\xi)=\eta\co(\xi)=\sum_{i=1}^nk_id_i(\eta\mathcal C_i)=\sum_{i=1}^nk_id_i\co(\xi_i)=\co\left(\sum_{i=1}^nk_id_i\xi_i\right),
\end{equation*}
whence $\eta\xi=\sum_{i=1}^nk_id_i\xi_i$ and so $\eta\xi^{\sharp}=\sum_{i=1}^nk_id_i(\xi_i)^{\sharp}$. By applying the last equality to the elements $e_1,\dots,e_n$ and by using the definition of the maps $\xi_1,\dots\xi_n$, we get $k_1=\dots=k_n=0$. This yields $\co(\xi)=0$ and hence also $\xi=0$.

To verify statement (ii), fix a map $\zeta\in C_z(X,\mathbb T^1)$; we show that $\co(\zeta)\in\sigma(\mathcal M)+(\mathcal S\cap\Cob)$. Set $k_i=\zeta^{\sharp}(e_i)$ for $i=1,\dots,n$. By our choice of the maps $\xi_i$, the morphism $(\sum_{i=1}^nk_id_i\xi_i)^{\sharp}=\sum_{i=1}^nk_id_i(\xi_i)^{\sharp}$ takes its values in $\eta\mathbb Z$ and so the map $\sum_{i=1}^nk_id_i\xi_i$ lifts across $\kappa_{\eta}\colon\mathbb T^1\to\mathbb T^1$ to an element $\vartheta$ of $C_z(X,\mathbb T^1)$.

Set $\xi=\zeta-\vartheta$; we show that $\xi\in\mathcal M$. For every $j=1,\dots,n$,
\begin{equation*}
\eta\xi^{\sharp}(e_j)=\eta\zeta^{\sharp}(e_j)-(\eta\vartheta)^{\sharp}(e_j)=\eta\zeta^{\sharp}(e_j)-\sum_{i=1}^nk_id_i(\xi_i)^{\sharp}(e_j)=\eta k_j-k_jd_j\frac{\eta}{d_j}=0
\end{equation*}
and hence $\xi^{\sharp}(e_j)=0$. Consequently, $\xi^{\sharp}H_1^w(\mathcal F)=0$ and, indeed, $\xi\in\mathcal M$.

Now, since $\co(\zeta)=\co(\xi)+\co(\vartheta)$ and $\co(\xi)\in\sigma(\mathcal M)$ by our discussion above, it remains to show that $\co(\vartheta)\in\mathcal S\cap\Cob$. By definition of $\vartheta$,
\begin{equation*}
\eta\co(\vartheta)=\co(\eta\vartheta)=\co\left(\sum_{i=1}^nk_id_i\xi_i\right)=\sum_{i=1}^nk_id_i\co(\xi_i)=
\eta\left(\sum_{i=1}^nk_i(d_i\mathcal C_i)\right),
\end{equation*}
which yields $\co(\vartheta)=\sum_{i=1}^nk_i(d_i\mathcal C_i)\in\langle d_1\mathcal C_1,\dots,d_n\mathcal C_n\rangle=\mathcal S\cap\Cob$.

\emph{5th step.} We verify the isomorphism (\ref{Eq:Coch.div.sum}).

By virtue of (d) from the first step of the proof, $\sigma(\mathcal M)$ is a subgroup of $\Div(\Cocc)$. Further, by the second and the third step of the proof, there is an isomorphism $\mathcal S/(\mathcal S\cap\Cob)\cong\mathbb Z_{d_1}\oplus\dots\oplus\mathbb Z_{d_n}$. Consequently, there are isomorphisms of groups
\begin{equation*}
\begin{split}
\Coch&=\Cocc/\Cob=\frac{\Div(\Cocc)\oplus\mathcal S}{\sigma(\mathcal M)\oplus(\mathcal S\cap\Cob)}\cong\frac{\Div(\Cocc)}{\sigma(\mathcal M)}\oplus\frac{\mathcal S}{\mathcal S\cap\Cob}\\
&\cong\frac{\Div(\Cocc)}{\sigma(\mathcal M)}\oplus\left(\mathbb Z_{d_1}\oplus\dots\oplus\mathbb Z_{d_n}\right).
\end{split}
\end{equation*}
By divisibility of $\Div(\Cocc)$, the group $\Div(\Cocc)/\sigma(\mathcal M)$ corresponds to the divisible subgroup $\Div(\Coch)$ of $\Coch$ under the isomorphism above. This verifies (\ref{Eq:Coch.div.sum}).

\emph{6th step.} We finish the proof of the theorem by verifying (\ref{Eq:Hom.Ext.Lie.mfld.2}).

The existence of an exact sequence (\ref{Eq:Hom.Ext.Lie.mfld.2}) follows from the exactness of (\ref{Eq:exact.sequence}) with the help of the following three observations. First, since $\pi^1(X)\cong H_1^w(X)\cong\mathbb Z^{n+m}$,
\begin{equation*}
\Ext(G^*,\pi^1(X))\cong\Ext(G^*,H_1^w(X))\cong\Ext(G^*,\mathbb Z^{n+m})\cong\Ext(G^*,\mathbb Z)^{n+m}.
\end{equation*}
Further, since $\Cocc\cong\Div(\Cocc)\oplus\mathbb Z^n$ by the second step of the proof,
\begin{equation*}
\Ext(G^*,\Cocc)\cong\Ext(G^*,\Div(\Cocc))\oplus\Ext(G^*,\mathbb Z^n)\cong\Ext(G^*,\mathbb Z)^n.
\end{equation*}
Finally, since $\Coch\cong\Div(\Coch)\oplus(\mathbb Z_{d_1}\oplus\dots\oplus\mathbb Z_{d_n})$ by the fifth step of the proof,
\begin{equation*}
\Ext(G^*,\Coch)\cong\Ext(G^*,\Div(\Coch))\oplus\Ext\left(G^*,\bigoplus_{i=1}^n\mathbb Z_{d_i}\right)\cong\bigoplus_{i=1}^n\Ext(G^*,\mathbb Z_{d_i}).
\end{equation*}
This verifies the exactness of (\ref{Eq:Hom.Ext.Lie.mfld.2}).
\end{proof}

\begin{theorem}\label{T:div.gp.min.ext}
Under the assumptions and notation from Theorem~\ref{T:ex.seq.Lie.mfld}, there exists a divisible subgroup $\mathfrak{D}$ of the groupoid $\Cochm$ such that $\Coch=\mathfrak{D}\oplus\tor(\Coch)$. Moreover, there are isomorphisms of groups
\begin{equation}\label{Eq:div.gp.min.ext}
\begin{split}
\tor(\Coch)&\cong(\mathbb Q/\mathbb Z)^m\oplus\tor\left(H_1^w(X)/H_1^w(\mathcal F)\right)\\
&\cong(\mathbb Q/\mathbb Z)^m\oplus\Big(\mathbb Z_{d_1}\oplus\dots\oplus\mathbb Z_{d_n}\Big).
\end{split}
\end{equation}
Finally, there are topological isomorphisms
\begin{equation}\label{Eq:div.gp.min.ext.dual}
\begin{split}
F(\Daleth)&\cong\mathfrak{D}^*\times\tor\left(\Coch\right)^*\cong\mathfrak{D}^*\times\left((\mathbb Q/\mathbb Z)^*\right)^m\times\mathbb Z_{d_1}\times\dots\times\mathbb Z_{d_n},
\end{split}
\end{equation}
where all $\mathfrak{D}$, $\tor(\Coch)$, $\mathbb Q/\mathbb Z$ and $\mathbb Z_{d_i}$ ($i=1,\dots,n$) are assumed to carry the discrete topology.
\end{theorem}
\begin{remark}\label{R:div.gp.min.ext}
We wish to add the following observations.
\begin{itemize}
\item By the conclusion of the theorem, $\mathfrak{D}$ is a divisible group. Moreover, being a complementary summand to $\tor(\Coch)$ in $\Coch$, $\mathfrak{D}$ is also torsion-free. Thus, algebraically, $\mathfrak{D}$ is a rational linear space.
\item Since the group $\mathfrak{D}$ is divisible and torsion-free, its dual group $\mathfrak{D}^*$ is torsion-free and connected. Further, since $\tor(\Coch)$ is a torsion group, its dual group $\tor(\Coch)^*$ is totally disconnected. Consequently, under the first of the two isomorphisms in (\ref{Eq:div.gp.min.ext.dual}), $\mathfrak{D}^*$ corresponds to $F(\Daleth)_0$. In particular, $F(\Daleth)_0$ is a topological direct summand in $F(\Daleth)$ and $F(\Daleth)/F(\Daleth)_0$ is topologically isomorphic to $\tor(\Coch)^*$.
\end{itemize}
\end{remark}
\begin{proof}[Proof of Theorem~\ref{T:div.gp.min.ext}]
Since $\Div(\Coch)$ is divisible by definition, its torsion subgroup $\tor(\Div(\Coch))$ is a direct summand in $\Div(\Coch)$ with a divisible torsion-free complementary summand $\mathfrak{D}$. Moreover,
\begin{equation*}
\begin{split}
\mathfrak{D}\setminus1&\subseteq\Div(\Coch)\setminus\tor(\Div(\Coch))=\Div(\Coch)\setminus\left(\tor(\Coch)\cap\Div(\Coch)\right)\\
&=\Div(\Coch)\setminus\tor(\Coch)\subseteq\Coch\setminus\tor(\Coch)=\Cochm\setminus1
\end{split}
\end{equation*}
and hence $\mathfrak{D}$ is a subgroup of the groupoid $\Cochm$.

Now, by virtue of (\ref{Eq:Hom.Ext.Lie.mfld.2}) from Theorem~\ref{T:ex.seq.Lie.mfld}, $\Coch$ contains a subgroup $\mathfrak{A}$ isomorphic to $\tor(H_1^w(X)/H_1^w(\mathcal F))\cong\mathbb Z_{d_1}\oplus\dots\oplus\mathbb Z_{d_n}$ such that $\Coch=\Div(\Coch)\oplus\mathfrak{A}$. Consequently,
\begin{equation*}
\Coch=\Div(\Coch)\oplus\mathfrak{A}=\mathfrak{D}\oplus\tor(\Div(\Coch))\oplus\mathfrak{A}.
\end{equation*}
Thus, $\tor(\Coch)=\tor(\Div(\Coch))\oplus\mathfrak{A}$ is a direct summand in $\Coch$ and $\mathfrak{D}$ is its complementary direct summand. Finally, by virtue of (\ref{Eq:tor.Lie.mfld}) from Theorem~\ref{T:tor.Lie.mnfld},
\begin{equation*}
\tor(\Div(\Coch))=\Div(\tor(\Coch))\cong\Div\left((\mathbb Q/\mathbb Z)^m\oplus\mathbb Z_{d_1}\oplus\dots\oplus\mathbb Z_{d_n}\right)\cong(\mathbb Q/\mathbb Z)^m,
\end{equation*}
which finishes the proof of (\ref{Eq:div.gp.min.ext})

The two isomorphisms from (\ref{Eq:div.gp.min.ext.dual}) follow immediately from the topological isomorphism $F(\Daleth)\cong(\Coch)_d^*$ obtained in Remark~\ref{R:prop.of.gimel} and from the first part of Theorem~\ref{T:div.gp.min.ext}.
\end{proof}

\chapter{Algebraic-Topological Aspects}\label{S:alg.top.asp}

\section{Lifts of transfer functions}\label{S:lifts.transfer}

Let $\Flow$ be a minimal flow. For an extension $\mathcal C\in\Coc$ of $\mathcal F$ with a totally disconnected section $F(\mathcal C)$, we investigate the relation between $F(\mathcal C)$ on one side and the existence of lifts of transfer functions for certain quotients of $\mathcal C$ on the other side. This relation is described in Proposition~\ref{P:char.F.as.smallest}.

In order to formulate the next result we fix some notation. Given $G\in\mathsf{CAGp}$ and $E\sbgp K\sbgp G$, we denote by $p_K$ and $p_E$\index[symbol]{$p_E$} the canonical quotient morphisms $G\to G/K$ and $G\to G/E$, respectively. Also, $q_E$\index[symbol]{$q_E$} stands for the (unique) morphism $G/E\to G/K$ with $q_Ep_E=p_K$. Notice that $q_E$ is a quotient morphism with kernel $\ker(q_E)=p_E(K)\cong K/E$ and it is thus equivalent to the canonical quotient morphism $G/E\to(G/E)/(K/E)$.

\begin{lemma}\label{L:cob.vs.lift}
Let $\Flow$ be a minimal flow, $G\in\mathsf{CAGp}$ and $\mathcal C\in\Cocc(G)$. Fix $z\in X$. Assume that $K\sbgp G$ is totally disconnected with $p_K\mathcal C\in\Cob(G/K)$ (that is, with $F(\mathcal C)\subseteq K$), and denote by $\xi$ the transfer function of $p_K\mathcal C$ with $\xi(z)=e$. If $X$ is connected and $\Gamma$ has no non-trivial finite abelian quotient groups then the following conditions are equivalent for every $E\sbgp K$:
\begin{enumerate}
\item[(a)] $p_E\mathcal C\in\Cob(G/E)$ (that is, $F(\mathcal C)\subseteq E$),
\item[(b)] $\xi$ lifts across $q_E$ to a continuous map $\eta\colon X\to G/E$.
\end{enumerate}
The situation is depicted in Figure~\ref{Fig:essen.fund}.
\end{lemma}
\begin{remark}\label{R:cob.vs.lift}
Notice the following facts.
\begin{itemize}
\item The assumption on $\Gamma$ from Lemma~\ref{L:cob.vs.lift} is satisfied for instance if $\Gamma$ is connected, which is the situation of our main interest. Of course, there are also disconnected groups without non-trivial finite abelian quotient groups, say, the discrete group of rationals $\mathbb Q_d$. More generally, a locally compact abelian group $\Gamma$ has no non-trivial finite (abelian) quotient groups if and only if its dual group $\Gamma^*$ is torsion-free. In case when $\Gamma$ is discrete abelian the latter condition is equivalent with the divisibility of $\Gamma$.
\item In connection with condition (b) observe that $\xi$ lifts across $q_E$ to a continuous map $X\to G/E$ if and only if it lifts to a continuous base point preserving map $X\to G/E$. For if $\eta\colon X\to G/E$ is continuous and satisfies $q_E\eta=\xi$ then $\eta':=\eta\eta(z)^{-1}\colon X\to G/E$ is continuous, preserves the base points and, since $\eta(z)\in q_E^{-1}(\xi(z))=\ker(q_E)$, $\eta'$ is also a lift of $\xi$ across~$q_E$.
\end{itemize}
\end{remark}
\begin{figure}[ht]
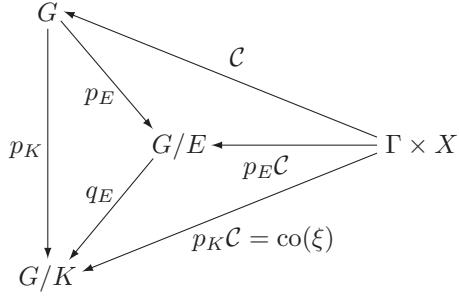

\[
\ctdiagram{
\ctinnermid
\ctv 0,0: {G/E}
\ctv -50,-50: {G/K}
\ctv -50,50: {G}
\ctv 92,0: {\Gamma\times X}
\ctv -58,0: {p_K}
\ctv -30,20: {p_E}
\ctv -30,-17: {q_E}
\ctv 33,-8: {p_E\mathcal C}
\ctv 32,-35: {p_K\mathcal C=\co(\xi)}
\ctv 22,33: {\mathcal C}
\put(75,0){\vector(-1,0){63}}
\put(75,3){\vector(-1,0.4){119}}
\put(75,-3){\vector(-1,-0.4){112}}
\put(-50,44){\vector(0,-1){87}}
\put(-45,47){\vector(1,-1.17){35}}
\put(-10,-5){\vector(-1,-1.2){32}}
\ctnohead
}
\]
\caption{Lifts of transfer functions of coboundaries.}
\label{Fig:essen.fund}
\end{figure}
\begin{proof}[Proof of Lemma~\ref{L:cob.vs.lift}]
Fix a closed subgroup $E$ of $K$. As has been observed at the beginning of this section, the morphism $q_E\colon G/E\to G/K$ is equivalent to the quotient morphism $G/E\to (G/E)/(K/E)$ and its kernel $\ker(q_E)=p_E(K)\cong K/E$ is therefore totally disconnected. 

First assume that $p_E\mathcal C\in\Cob(G/E)$ and denote by $\eta$ the transfer function of $p_E\mathcal C$ with $\eta(z)=e$. Then $\co(\xi)=p_K\mathcal C=q_Ep_E\mathcal C=\co(q_E\eta)$. Since the base flow $\mathcal F$ is minimal and $q_E\eta(z)=e=\xi(z)$, it follows that $q_E\eta=\xi$ on the whole $X$. This shows that (a) implies (b).

Now assume that $\eta\colon X\to G/E$ is continuous and satisfies $q_E\eta=\xi$. Then $q_E\co(\eta)=\co(\xi)=p_K\mathcal C=q_E(p_E\mathcal C)$ and hence $(p_E\mathcal C)^{-1}\co(\eta)\in\Cocc(G/E)$ takes its values in a totally disconnected subgroup $\ker(q_E)$ of $G/E$. Now let $\chi\in(\ker(q_E))^*$. Since $\ker(q_E)$ is (compact and) totally disconnected, $\chi$ is a torsion element of $(\ker(q_E))^*$ and hence $\chi((p_E\mathcal C)^{-1}\co(\eta))$ is a torsion element of $\Cocc$. By Theorem~\ref{T:gimel.conn} it follows that $\chi((p_E\mathcal C)^{-1}\co(\eta))=1$. As the latter identity holds for every character $\chi$ of $\ker(q_E)$, we conclude with $p_E\mathcal C=\co(\eta)\in\Cob(G/E)$. This shows that (b) implies (a).
\end{proof}

\begin{proposition}\label{P:char.F.as.smallest}
Let $\Flow$ be a minimal flow, $G\in\mathsf{CAGp}$ and $\mathcal C\in\Cocc(G)$. Fix $z\in X$. Assume that $K\sbgp G$ is totally disconnected with $p_K\mathcal C\in\Cob(G/K)$ (that is, with $F(\mathcal C)\subseteq K$), and denote by $\xi$ the transfer function of $p_K\mathcal C$ with $\xi(z)=e$. If $X$ is connected and $\Gamma$ has no non-trivial finite abelian quotient groups then the following statements hold:
\begin{enumerate}
\item[(i)] $F(\mathcal C)$ equals the smallest of all the closed subgroups $E$ of $K$, for which $\xi$ lifts across $q_E$ to a continuous map $\eta\colon X\to G/E$,
\item[(ii)] $F(\mathcal C)$ equals the smallest of all the closed subgroups $E$ of $K$, for which the morphism $\xi^*\colon (G/K)^*\to C_z(X,\mathbb T^1)$ extends through $(q_E)^*\colon (G/K)^*\to (G/E)^*$ to a morphism $\sigma\colon (G/E)^*\to C_z(X,\mathbb T^1)$.
\end{enumerate}
\end{proposition}
\begin{remark}\label{R:char.F.as.smallest}
We shall have an occasion to use the proposition also in a special case when $K$ is finite. In such case every subgroup $E$ of $K$ is finite and hence closed, and so the word ``closed'' becomes redundant in the reading of both (i) and (ii).
\end{remark}
\begin{proof}[Proof of Proposition~\ref{P:char.F.as.smallest}]
We verify statement (i) by showing that the following statements are equivalent for every closed subgroup $E$ of $K$:
\begin{enumerate}
\item[($\alpha$)] $F(\mathcal C)\subseteq E$,
\item[($\beta$)] $\xi$ lifts across $q_E$ to a continuous map $\eta\colon X\to G/E$.
\end{enumerate}
So fix a closed subgroup $E$ of $K$ and assume that $F(\mathcal C)\subseteq E$. Then, since $p_{F(\mathcal C)}\mathcal C\in\Cob(G/F(\mathcal C))$, it follows that $p_E\mathcal C\in\Cob(G/E)$ and the existence of a map $\eta$ as in ($\beta$) follows from Lemma~\ref{L:cob.vs.lift}. This verifies that ($\beta$) follows from ($\alpha$). Conversely, assume that $\xi$ lifts across $q_E$ to a continuous map $\eta\colon X\to G/E$. Then, by Lemma~\ref{L:cob.vs.lift}, $p_E\mathcal C\in\Cob(G/E)$ and so $e=F(p_E\mathcal C)=p_EF(\mathcal C)$. Thus, $F(\mathcal C)\subseteq E$, which shows that ($\alpha$) follows from ($\beta$).

Statement (ii) follows from statement (i) by virtue of the following observations:
\begin{itemize}
\item if $\xi$ lifts across $q_E$ to a continuous map $X\to G/E$ then it lifts also to a continuous base point preserving map $X\to G/E$; this has already been observed in Remark~\ref{R:cob.vs.lift},
\item the continuous base point preserving maps $X\to G/E$ are in a one-to-one corres\-pon\-dence with the morphisms $(G/E)^*\to C_z(X,\mathbb T^1)$, the correspondence assigning the morphism $\eta^*$ to every $\eta\in C_z(X,G/E)$; this follows from our discussion in Subsection~\ref{Sub:mps.cpt.ab.gps}, see formula (\ref{Eq:maps.are.homs}),
\item a map $\eta\in C_z(X,G/E)$ satisfies $q_E\eta=\xi$ if and only if the morphism $\eta^*$ satisfies $\eta^*(q_E)^*=\xi^*$; this follows from the fact that the characters of $G/K$ separate the points of $G/K$.
\end{itemize}
\end{proof}

\section{Relation to the monodromy action}\label{S:rel.mndrm.act}

Our aim in this section is to describe a relationship between the functor $F$ and the functor $\pi_1$. The situation that we shall consider is as follows. Assume that $\Flow$ is a minimal flow and $\mathcal C\in\Coc$ is an extension of $\mathcal F$ taking its values in a group $G\in\mathsf{CAGp}$. Suppose that $K$ is a closed totally disconnected subgroup of $G$ with $F(\mathcal C)\subseteq K$ or, equivalently, with $p_K\mathcal C\in\Cob(G/K)$. Let $\xi$ be the base point preserving transfer function for $p_K\mathcal C$. Is it possible to express the group $F(\mathcal C)$ in terms of the map $\xi$ or in terms of the morphism $\xi^{\sharp}\colon\pi_1(X)\to\pi_1(G/K)$ induced by $\xi$? We show that the answer to the proposed question is yes, provided $\Gamma$ is a connected Lie group and $X$ is a compact (connected) manifold; in Theorem~\ref{T:F.and.pi} we show that the group $F(\mathcal C)$ equals the orbit closure of the identity $e$ under the monodromy action of $\xi^{\sharp}\pi_1(X)$ on $K$. Before proving this theorem, we study epimorphisms in $\mathsf{CAGp}$ with totally disconnected kernels. We show that they have properties analogous to those of the covering morphisms; in particular, a monodromy action of $\pi_1(G/K)$ on $K$ is defined in the situation described above. After gaining more information about the monodromy action in this setting, we prove the mentioned Theorem~\ref{T:F.and.pi}, which is our main result in this section.

\begin{lemma}\label{L:lift.tot.disc.ker}
Let $G\in\mathsf{CAGp}$, $K$ be a closed totally disconnected subgroup of $G$ and $p_K\colon G\to G/K$ be the quotient morphism. Let $X\in\mathsf{CLAC}$ be a pointed space with the base point $z$ and $\xi\colon X\to G/K$ be a continuous base point preserving map. Then the following conditions are equivalent:
\begin{enumerate}
\item[(a)] $\xi$ lifts across $p_K$ to a continuous base point preserving map $\zeta\colon X\to G$,
\item[(b)] $\xi^{\sharp}\pi_1(X)\subseteq p_K^{\sharp}\pi_1(G)$.
\end{enumerate}
\end{lemma}
\begin{remark}\label{R:lift.tot.disc.ker}
Notice the following facts.
\begin{itemize}
\item The map $\zeta$ from statement (a) is obviously unique, for the space $X$ is connected and the kernel $K$ of $p_K$ is totally disconnected by the assumptions of the lemma.
\item The statement of the lemma in case of a finite group $K$ follows at once from the well known and standard results from algebraic topology, for in such case the morphism $p_K$ is a covering map (see, e.g., \cite[Proposition~1.33,~p.~61]{Hat}). In the situation when $K$ is infinite the morphism $p_K$ is no longer a covering map and so the usual arguments from algebraic topology do not apply.
\item From our proof of the lemma it will follow that the paths in $G/K$ do lift across $p_K$ as well as do the homotopies of paths. Thus, if two loops $f_1,f_2$ in $G/K$ based at $e$ are path-homotopic, then their lifts $\widetilde{f}_1,\widetilde{f}_2$ across $p_K$ starting at $e$ satisfy $\widetilde{f}_1(1)=\widetilde{f}_2(1)$. Consequently, there is a monodromy action of the group $\pi_1(G/K)$ on the $p_K$-fibre $K$ of $e$, just as in the case when $p_K$ is a covering map. The action is given by the rule $\pi_1(G/K)\times K\ni(f,k)\mapsto k\widetilde{f}(1)\in K$, where $\widetilde{f}$ is the lift of $f$ across $p_K$ starting at $e$.
\end{itemize}
\end{remark}
\begin{proof}[Proof of Lemma~\ref{L:lift.tot.disc.ker}]
We divide the proof into four steps. Before turning into it we remind the reader that the base point for the fundamental group of a topological group $G$ is always the identity $e$ of $G$. We also recall that the operation of composition of paths is denoted by $*$, the reverse of a path $f$ by $\overline{f}$ and the relation of homotopy of paths by $\simeq_p$.

\emph{1st step.} We fix some notation that shall be used in the rest of the proof.

Since $K$ is a totally disconnected compact abelian group, it is a projective limit of finite abelian groups. Hence there is a filter base $(H_j)_{j\in J}$ in $K$ formed by closed subgroups of $K$ with $\bigcap_{j\in J}H_j=e$ and with $K/H_j$ finite for every $j$. For $j,l\in J$ write $j\leq l$ when (and only when) $H_l\subseteq H_j$. The quotient morphism $G\to G/H_j$ will be denoted by $p_j$ for every $j\in J$. Given $j,l\in J$ with $j\leq l$, the inclusion $H_l\subseteq H_j$ yields the existence of a (unique) quotient morphism $p_{jl}\colon G/H_l\to G/H_j$ with $p_{jl}p_l=p_j$. The epimorphisms $p_{jl}$ ($j\leq l$) constitute an inverse system in $\mathsf{CAGp}$ indexed by the directed set $J$. As follows from the equality $\bigcap_{j\in J}H_j=e$, the inverse limit of this system is the group $G$ along with the morphisms $p_j$ ($j\in J$) as the limit projections. Finally, for every $j\in J$, the inclusion $H_j\subseteq K$ yields the existence of a (unique) quotient morphism $r_j\colon G/H_j\to G/K$ with $r_jp_j=p_K$. An elementary argument shows that $r_jp_{jl}=r_l$ for all $j\leq l$.

\emph{2nd step.} Let $f\colon[0,1]\to G/K$ be a path starting at the identity $e$ of $G/K$. We show that $f$ lifts uniquely across $p_K$ to a path $\widetilde{f}\colon[0,1]\to G$ starting at the identity $e$ of $G$.

The uniqueness part follows at once from the connectedness of $[0,1]$ and from the fact that the kernel $K$ of $p_K$ is totally disconnected. To prove the existence part we use an inverse limit argument. Given $j\in J$, the morphism $r_j\colon G/H_j\to G/K$ is a covering morphism, for its kernel $\ker(r_j)=p_j(K)\cong K/H_j$ is finite. Consequently, $f$ lifts across $r_j$ to a path $f_j$ in $G/H_j$ starting at $e$. If $j\leq l$ then $p_{jl}f_l$ is a path in $G/H_j$ starting at $e$ and projecting by $r_j$ to $f$, and hence $p_{jl}f_l=f_j$ by uniqueness of $f_j$. Thus, the paths $f_j$ ($j\in J$) form a system of maps compatible with the inverse system of groups $G/H_j\stackrel{p_{jl}}{\longleftarrow}G/H_l$ and hence they give rise to a continuous map $\widetilde{f}\colon[0,1]\to G$ with $p_j\widetilde{f}=f_j$ for every $j\in J$. Clearly, $\widetilde{f}(0)=e$ and so it remains to show that $\widetilde{f}$ is a lift of $f$ across $p_K$. This is immediate, for if $j\in J$ then $p_K\widetilde{f}=r_jp_j\widetilde{f}=r_jf_j=f$.

\emph{3rd step.} Let $F\colon[0,1]^2\to G/K$ be a homotopy of paths which starts at the identity $e$ of $G/K$. We claim that $F$ lifts uniquely across $p_K$ to a homotopy of paths $\widetilde{F}\colon[0,1]^2\to G$ starting at the identity $e$ of $G$.

The proof follows the line of the argument used in the second step of the proof and so we omit the details. Let us only mention that in proving the claim one uses an inverse limit argument along with the following two facts:
\begin{itemize}
\item homotopies of paths do lift to homotopies of paths across the quotient morphisms in $\mathsf{CAGp}$ with finite kernels,
\item the inverse limit of a system of homotopies of paths is again a homotopy of paths.
\end{itemize}

\emph{4th step.} We finish the proof of the lemma by verifying the equivalence of (a) and (b). Implication (a)$\Rightarrow$(b) is clear and so we verify only the converse implication (b)$\Rightarrow$(a). Our construction of the map $\zeta$ as in (a) will proceed as in the case of $K$ finite, by using the method of lifting appropriate paths across $p_K$. To verify the continuity of $\zeta$ we shall use an inverse limit argument.

We start by verifying the following auxiliary claim:
\begin{enumerate}
\item[($*$)] Let $f_1,f_2$ be paths in $X$ with $f_1(0)=f_2(0)=z$ and $f_1(1)=f_2(1)$, and let $\widetilde{f}_1,\widetilde{f}_2$ be the lifts of the paths $\xi f_1,\xi f_2$ across $p_K$ with $\widetilde{f}_1(0)=\widetilde{f}_2(0)=e$. Then $\widetilde{f}_1(1)=\widetilde{f}_2(1)$.
\end{enumerate}
By virtue of (b), there is $\widetilde{l}\in\pi_1(G)$ with $\xi(f_1*\overline{f_2})\simeq_p p_K\widetilde{l}$. Then
\begin{equation*}
\xi f_1\simeq_p\xi((f_1*\overline{f_2})*f_2)=\xi(f_1*\overline{f_2})*\xi f_2\simeq_p p_K\widetilde{l}*\xi f_2
\end{equation*}
and hence, by the third step of the proof, the lifts $\widetilde{f}_1$ and $\widetilde{l}*\widetilde{f}_2$ of $\xi f_1$ and $p_K\widetilde{l}*\xi f_2$, respectively, are homotopic paths. Thus, $\widetilde{f}_1(1)=(\widetilde{l}*\widetilde{f}_2)(1)=\widetilde{f}_2(1)$, as was to be shown.

Now, given $x\in X$, let $f$ be a path in $X$ from $z$ to $x$. Then $\xi f$ is a path in $G/K$ from $e$ to $\xi(x)$ and hence, by the second step of the proof, it lifts across $p_K$ to a (unique) path $\widetilde{f}$ in $G$ with $\widetilde{f}(0)=e$. We set $\zeta(x)=\widetilde{f}(1)$. By virtue of ($*$), $\zeta$ is a well defined map $X\to G$ with $p_K\zeta=\xi$ and $\zeta(z)=e$. Thus, to finish the proof it remains to show that $\zeta$ is continuous. Given $j\in J$, we have $\xi^{\sharp}\pi_1(X)\subseteq p_K^{\sharp}\pi_1(G)=r_j^{\sharp}p_j^{\sharp}\pi_1(G)\subseteq r_j^{\sharp}\pi_1(G/H_j)$. Since, in addition, $r_j\colon G/H_j\to G/K$ is a covering morphism, it follows that the map $\xi$ lifts uniquely across $r_j$ to a continuous base point preserving map $\zeta_j\colon X\to G/H_j$. If $x\in X$, $f$ is a path in $X$ from $z$ to $x$ and $\widetilde{f}$ is the continuous lift of $\xi f$ across $p_K$ starting at $e$, then both $p_j\widetilde{f}$ and $\zeta_jf$ are continuous lifts of $\xi f$ across $r_j$ starting at $e$. Thus $p_j\widetilde{f}=\zeta_jf$ and hence $p_j\zeta(x)=p_j\widetilde{f}(1)=\zeta_jf(1)=\zeta_j(x)$. This shows that $p_j\zeta=\zeta_j$ is continuous for every $j\in J$ and hence the map $\zeta$ itself is continuous.
\end{proof}

\begin{corollary}\label{C:lift.cob.Lie.gp}
Let $\Flow$ be a minimal flow with $\Gamma\in\mathsf{LieGp}$ connected and with $X$ a compact connected manifold. Let $G\in\mathsf{CAGp}$, $K\sbgp G$ be totally disconnected and $p_K\colon G\to G/K$ be the quotient morphism. Given a continuous base point preserving map $\xi\colon X\to G/K$, the following conditions are equivalent:
\begin{enumerate}
\item[(a)] $\co(\xi)$ lifts across $p_K$ to an extension $\mathcal C\in\Cocc(G)$,
\item[(b)] $\xi^{\sharp}\mathcal F_z^{\sharp}\pi_1(\Gamma)\subseteq p_K^{\sharp}\pi_1(G)$,
\item[(c)] $\xi^{\sharp}H_1^w(\mathcal F)\subseteq p_K^{\sharp}H_1^w(G)$.
\end{enumerate}
\end{corollary}
\begin{proof}
By virtue of (\ref{Eq:tor.Lie.mnf.L1}) from Lemma~\ref{L:tor.Lie.mnfld}, $\co(\xi)^{\sharp}\pi_1(\Gamma\times X)=\xi^{\sharp}\mathcal F_z^{\sharp}\pi_1(\Gamma)$. By applying also Lemma~\ref{L:lift.tot.disc.ker} and Lemma~\ref{P:lifting.cocycle}, we obtain the equivalence of the following conditions:
\begin{itemize}
\item $\co(\xi)$ lifts across $p_K$ to an extension $\mathcal C\in\Cocc(G)$,
\item $\co(\xi)$ lifts across $p_K$ to a continuous base point preserving map $\Gamma\times X\to G$,
\item $\co(\xi)^{\sharp}\pi_1(\Gamma\times X)\subseteq p_K^{\sharp}\pi_1(G)$,
\item $\xi^{\sharp}\mathcal F_z^{\sharp}\pi_1(\Gamma)\subseteq p_K^{\sharp}\pi_1(G)$.
\end{itemize}
This verifies the equivalence of (a) and (b).

Before turning to the proof of the equivalence of (b) and (c) let us add some comments on the notation. The symbols $p_{\Gamma}$ and $p_{G/K}$ stand for the canonical quotient morphisms $\pi_1(\Gamma)\to H_1^w(\Gamma)$ and $\pi_1(G/K)\to H_1^w(G/K)$, respectively, whereas the symbol $p_K$ is used to denote the quotient morphism $G\to G/K$ according to the assumptions of the corollary. Recall also that by definition of $H_1^w(\mathcal F)$, $H_1^w(\mathcal F)=\mathcal F_z^{\sharp}H_1^w(\Gamma)=p_X\mathcal F_z^{\sharp}\pi_1(\Gamma)$, where $p_X$ stands for the canonical quotient morphism $\pi_1(X)\to H_1^w(X)$.

Now, since the group $\pi_1(G/K)$ is abelian and torsion-free, the morphism $p_{G/K}$ is an isomorphism. Thus, in order to verify the equivalence of (b) and (c), it suffices to check the following two identities:
\begin{enumerate}
\item[($\alpha$)] $\xi^{\sharp}H_1^w(\mathcal F)=p_{G/K}\xi^{\sharp}\mathcal F_z^{\sharp}\pi_1(\Gamma)$,
\item[($\beta$)] $p_K^{\sharp}H_1^w(G)=p_{G/K}p_K^{\sharp}\pi_1(G)$.
\end{enumerate}
First, we have
\begin{equation*}
\xi^{\sharp}H_1^w(\mathcal F)=\xi^{\sharp}\mathcal F_z^{\sharp}H_1^w(\Gamma)=\xi^{\sharp}\mathcal F_z^{\sharp}p_{\Gamma}\pi_1(\Gamma)=p_{G/K}\xi^{\sharp}\mathcal F_z^{\sharp}\pi_1(\Gamma),
\end{equation*}
which verifies ($\alpha$). A similar computation yields
\begin{equation*}
p_K^{\sharp}H_1^w(G)=p_K^{\sharp}p_G\pi_1(G)=p_{G/K}p_K^{\sharp}\pi_1(G),
\end{equation*}
where $p_G$ denotes the quotient morphism $\pi_1(G)\to H_1^w(G)$. This verifies ($\beta$).
\end{proof}

In order to formulate the next lemma we fix some notation. Let $G\in\mathsf{CAGp}$ and $K\sbgp G$ be totally disconnected. Consider the quotient morphism $p_K\colon G\to G/K$ and the induced morphism of fundamental groups $p_K^{\sharp}\colon\pi_1(G)\to \pi_1(G/K)$, where, as usual, the identities of the groups are taken as the base points for their fundamental groups. Recall that the groups $\pi_1(G)$ and $\pi_1(G/K)$ are abelian.

Given a subgroup $E$ of $K$, let $\mathcal Q_K(E)$\index[symbol]{$\mathcal Q_K(E)$} consist of those elements of $\pi_1(G/K)$, which lift across $p_K$ to a path both starting and ending in $E$ (or, equivalently, whose lift starting at $e$ ends at an element of $E$). (Recall from Remark~\ref{R:lift.tot.disc.ker} that $\mathcal Q_K(E)$ is a well defined  object since we have a monodromy action of $\pi_1(G/K)$ on $K$ just as in the case of a finite kernel $K$. Thus, the endpoint $\widetilde{f}(1)$ of the lift $\widetilde{f}$ of a loop $f$ in $G/K$ across $p_K$ with $\widetilde{f}(0)=e$ depends indeed only on the path-homotopy class of $f$.) Observe also that the identity $\mathcal Q_K(E)=\mathcal Q_K(E\cap G_a)$ holds by definition of $\mathcal Q_K$. Despite the latter identity we shall not restrict ourselves to the subgroups of $K$ contained in $G_a$. We have two reasons for doing so. First, $K\cap G_a$ need not be a closed subgroup of $K$. Second, some of the subgroups $E$ of $K$, which we will have to deal with later in our considerations, will not be contained in $G_a$. (See Remarks~\ref{R:F.and.pi} and~\ref{R:F.and.pi.E.not.clsd} below.)

Conversely, given a subgroup $Q$ of $\pi_1(G/K)$, let $\mathcal E_K(Q)$\index[symbol]{$\mathcal E_K(Q)$} be the set of all $k\in K$, such that there exists a path in $G$ starting at $e$, ending at $k$ and projecting by $p_K$ to an element of $Q$. Clearly, $\mathcal E_K(Q)=\mathcal E_K(Q+p_K^{\sharp}\pi_1(G))$ and we shall therefore often restrict ourselves to the subgroups $Q$ of $\pi_1(G/K)$ containing $p_K^{\sharp}\pi_1(G)$. Observe that $\mathcal E_K(Q)$ is the orbit of the identity $e$ under the monodromy action of $Q$ on $K$. The following lemma summarizes some of the properties of $\mathcal Q_K$ and $\mathcal E_K$.

\begin{lemma}\label{L:QK.EK.tot.disc}
Let $G\in\mathsf{CAGp}$, $K\sbgp G$ be totally disconnected and $p_K\colon G\to G/K$ be the quotient morphism. Fix a subgroup $E$ of $K$ and a subgroup $Q$ of $\pi_1(G/K)$ containing $p_K^{\sharp}\pi_1(G)$. Then the following statements hold:
\begin{enumerate}
\item[(i)] $\mathcal Q_K(E)$ is a subgroup of $\pi_1(G/K)$ containing $p_K^{\sharp}\pi_1(G)$,
\item[(ii)] $\mathcal E_K(Q)$ is a subgroup of $K$ contained in $G_a$,
\item[(iii)] there are identities
\begin{equation*}
\mathcal E_K(\mathcal Q_K(E))=E\cap G_a\hspace{4mm}\text{and}\hspace{4mm}\mathcal Q_K(\mathcal E_K(Q))=Q,
\end{equation*}
\item[(iv)] there is an isomorphism
\begin{equation*}
\mathcal E_K(Q)\cong Q/p_K^{\sharp}\pi_1(G).
\end{equation*}
\end{enumerate}
\end{lemma}
\begin{remark}\label{R:QK.EK.tot.disc}
We wish to mention the following facts.
\begin{itemize}
\item Under the assumptions of the lemma there are identities $\mathcal E_K(p_K^{\sharp}\pi_1(G))=e$ and $\mathcal E_K(\pi_1(G/K))=K\cap G_a$. In particular, the latter identity and statement (iv) together yield an isomorphism $K\cap G_a\cong\pi_1(G/K)/p_K^{\sharp}\pi_1(G)$.
\item In Examples~\ref{E:ess.discr.torus} and~\ref{E:ess.discr.solen} below we describe the correspondence $\mathcal E_K$ in the situ\-ations when $G=\mathbb T^n$ is a finite-dimensional torus and when $G=S_{{\bf p}}$ is a solenoid.
\end{itemize}
\end{remark}
\begin{proof}[Proof of Lemma~\ref{L:QK.EK.tot.disc}]
We verify statement (i). Inclusion $p_K^{\sharp}\pi_1(G)\subseteq\mathcal Q_K(E)$ is an immediate consequence of the definition of $\mathcal Q_K$. To see that $\mathcal Q_K(E)$ is a subgroup of $\pi_1(G/K)$, fix $f_1,f_2\in\mathcal Q_K(E)$ and let $\widetilde{f}_1,\widetilde{f}_2$ be the lifts of $f_1,f_2$ across $p_K$ starting at $e$. By definition of $\mathcal Q_K(E)$, $k_1=\widetilde{f}_1(1)$ and $k_2=\widetilde{f}_2(1)$ are both elements of $E$. Now, the path $\widetilde{f}_1*\overline{k_1k_2^{-1}\widetilde{f}_2}$ is the lift of $f_1*\overline{f_2}$ across $p_K$ starting at $e$ and its endpoint is $(\widetilde{f}_1*\overline{k_1k_2^{-1}\widetilde{f}_2})(1)=k_1k_2^{-1}\widetilde{f}_2(0)=k_1k_2^{-1}\in E$. Thus, $f_1*\overline{f_2}\in\mathcal Q_K(E)$, as was to be shown.

We verify statement (ii). Inclusion $\mathcal E_K(Q)\subseteq K\cap G_a$ follows directly from the definition of $\mathcal E_K$. To see that $\mathcal E_K(Q)$ is a subgroup of $K$, fix $k_1,k_2\in\mathcal E_K(Q)$ and let $\widetilde{f}_1,\widetilde{f}_2$ be paths in $G$ from $e$ to $k_1,k_2$, projecting by $p_K$ to elements $f_1,f_2$ of $Q$, respectively. Since $Q$ is a subgroup of $\pi_1(G/K)$, $f_1*\overline{f_2}$ is an element of $Q$ and its lift $\widetilde{f}_1*\overline{k_1k_2^{-1}\widetilde{f}_2}$ across $p_K$ starting at $e$ thus ends at an element of $\mathcal E_K(Q)$. That is, $k_1k_2^{-1}=(\widetilde{f}_1*\overline{k_1k_2^{-1}\widetilde{f}_2})(1)\in\mathcal E_K(Q)$, as was to be shown.

The first identity from statement (iii) is an immediate consequence of the definitions of $\mathcal Q_K$ and $\mathcal E_K$, and of the existence and uniqueness of the lifts $\widetilde{f}$ across $p_K$ with $\widetilde{f}(0)=e$ for the loops $f\in\pi_1(G/K)$, see Lemma~\ref{L:lift.tot.disc.ker}. The second identity from (iii) follows also from the definitions of $\mathcal Q_K$ and $\mathcal E_K$, and from the assumption $p_K^{\sharp}\pi_1(G)\subseteq Q$.

Consider the map $\varphi\colon Q\to\mathcal E_K(Q)$, associating with every loop $f\in Q$ the endpoint of its lift $\widetilde{f}$ across $p_K$ with $\widetilde{f}(0)=e$. To verify statement (iv) it suffices to show that $\varphi$ is an epimorphism of groups with kernel $p_K^{\sharp}\pi_1(G)$. First, $\varphi$ is well defined due to the monodromy action of $\pi_1(G/K)$ on $K$, and it is surjective by definition of $\mathcal E_K(Q)$. To see that $\varphi$ is a morphism of groups, fix $f_1,f_2\in Q$ and consider their lifts $\widetilde{f}_1,\widetilde{f}_2$ across $p_K$ starting at $e$. Then $\varphi(f_1*f_2)=(\widetilde{f}_1*\widetilde{f}_1(1)\widetilde{f}_2)(1)=\widetilde{f}_1(1)\widetilde{f}_2(1)=\varphi(f_1)\varphi(f_2)$, as was to be shown. Finally, the identity $\ker(\varphi)=p_K^{\sharp}\pi_1(G)$ is clear and the proof is thus finished.
\end{proof}

\begin{example}\label{E:ess.discr.torus}
We describe the operation $\mathcal E_K$ in the case when the group $K$ is finite and either $G$ or $G/K$ is a finite-dimensional torus. To justify the setting used in this example recall the following well known facts.
\begin{itemize}
\item Every quotient group of $\mathbb T^n$ modulo a finite group is topologically isomorphic to $\mathbb T^n$ itself. Hence, every covering morphism from $\mathbb T^n$ is equivalent to a surjective endomorphism $p$ of $\mathbb T^n$.
\item If $G\in\mathsf{CAGp}$ is connected and $\mathbb T^n$ is a quotient group of $G$ modulo a finite subgroup of $G$ then $G$ is topologically isomorphic to $\mathbb T^n$.
\end{itemize}
Now let $p$ be a surjective endomorphism of $\mathbb T^n$ induced by an integer matrix $A=(a_{jl})$. That is, the matrix $A$ is regular and $p(z_m)_{m=1}^n=\left(z_1^{a_{m1}}\dots z_n^{a_{mn}}\right)_{m=1}^n$ for every $(z_m)_{m=1}^n\in\mathbb T^n$. Recall the isomorphism $\pi_1(\mathbb T^n)\cong\mathbb Z^n$,
\begin{equation*}
\mathbb Z^n\ni(k_1,\dots,k_n)\mapsto f_{k_1,\dots,k_n}\colon [0,1]\ni t\mapsto \left(e^{i2\pi k_1t},\dots,e^{i2\pi k_nt}\right)\in\mathbb T^n,
\end{equation*}
under which the induced morphism $p^{\sharp}\colon\pi_1(\mathbb T^n)\to\pi_1(\mathbb T^n)$ takes the form of the endomorphism of $\mathbb Z^n$ induced by the matrix $A$. That is, $p^{\sharp}\colon\mathbb Z^n\to\mathbb Z^n$ acts by the rule $p^{\sharp}(k_m)_{m=1}^n=(a_{m1}k_1+\dots+a_{mn}k_n)_{m=1}^n$. Let $B=(b_{jl})$ be the inverse of $A$ in $GL_n(\mathbb Q)$. Then for all $k_1,\dots,k_n\in\mathbb Z$, the loop $f_{k_1,\dots,k_n}$ lifts across $p$ to the path
\begin{equation*}
\widetilde{f}_{k_1,\dots,k_n}\colon[0,1]\ni t\mapsto\left(e^{i2\pi\left(\sum_{l=1}^n b_{1l}k_l\right)t},\dots,e^{i2\pi\left(\sum_{l=1}^n b_{nl}k_l\right)t}\right)\in\mathbb T^n.
\end{equation*}
The path $\widetilde{f}_{k_1,\dots,k_n}$ starts at the identity $1$ of $\mathbb T^n$ and ends at
\begin{equation*}
\begin{split}
\widetilde{f}_{k_1,\dots,k_n}(1)&=\left(e^{i2\pi\sum_{l=1}^n b_{1l}k_l},\dots,e^{i2\pi\sum_{l=1}^n b_{nl}k_l}\right)=q\left(\sum_{l=1}^n b_{1l}k_l,\dots,\sum_{l=1}^n b_{nl}k_l\right)\\
&=q\left(\left(B(k_1,\dots,k_n)^{T}\right)^{T}\right),
\end{split}
\end{equation*}
where $(-)^T$ stands for the transposing operator and $q\colon\mathbb R^n\to\mathbb T^n=\mathbb R^n/\mathbb Z^n$ for the usual covering morphism. Consequently, $p^{\sharp}\pi_1(\mathbb T^n)=A\mathbb Z^n$ and the operation $\mathcal E_K$ with $K=\ker(p)$ can be described as follows:
\begin{itemize}
\item if $Q$ is a subgroup of $\mathbb Z^n$ containing $A\mathbb Z^n$ then $\mathcal E_K(Q)=q(BQ)$.
\end{itemize}

With some additional information on $p$, the operation $\mathcal E_K$ can be simplified as follows. Let $e_1,\dots,e_n$ be the standard basis for $\mathbb Z^n$. Then for every monomorphism $h\colon\mathbb Z^n\to\mathbb Z^n$ there are automorphisms $\varphi,\psi$ of $\mathbb Z^n$ and positive integers $d_1,\dots,d_n$ with $d_i$ dividing $d_{i+1}$ for $i=1,\dots,n-1$, such that the composition $k=\varphi h\psi$ acts by the rule $k(e_i)=d_ie_i$ ($i=1,\dots,n$). Consequently, every surjective endomorphism of $\mathbb T^n$ is equivalent to a morphism of the form $p(z_1,\dots,z_n)=(z_1^{d_1},\dots,z_n^{d_n})$ with $d_1,\dots,d_n$ as above. In such case the matrix $A$ of $p$ is the diagonal matrix with $d_1,\dots,d_n$ as its diagonal elements. Hence, the inverse $B=A^{-1}$ of $A$ in $GL_n(\mathbb Q)$ is also a diagonal matrix with the diagonal elements $1/d_1,\dots,1/d_n$. Having expressed $p$ in such form we see that
\begin{equation*}
p^{\sharp}\pi_1(\mathbb T^n)=A\mathbb Z^n=d_1\mathbb Z\oplus\dots\oplus d_n\mathbb Z,
\end{equation*}
and for every subgroup $Q$ of $\mathbb Z^n$ containing $d_1\mathbb Z\oplus\dots\oplus d_n\mathbb Z$,
\begin{equation*}
\mathcal E_K(Q)=q(BQ)=\left\{\left(e^{i2\pi\frac{q_1}{d_1}},\dots,e^{i2\pi\frac{q_n}{d_n}}\right) : (q_1,\dots,q_n)\in Q\right\}
\end{equation*}
with $K=\ker(p)$.
\end{example}

\begin{example}\label{E:tow.ess.discr.solen}
Let $G\in\mathsf{CAGp}$, $K\sbgp G$ be totally disconnected and $p_K\colon G\to G/K$ be the quotient morphism. We show that for every positive integer $n$,
\begin{equation*}
\mathcal E_K\left(n\pi_1(G/K)\right)=n(K\cap G_a).
\end{equation*}
Fix $k\in\mathcal E_K(n\pi_1(G/K))$. There is a loop $f\in n\pi_1(G/K)$, whose lift $\widetilde{f}$ across $p_K$ starting at $e$ ends at $k$. Choose $g\in\pi_1(G/K)$ with $f=ng$ and denote by $\widetilde{g}$ the lift of $g$ across $p_K$ starting at $e$. The endpoint $l$ of $\widetilde{g}$ is clearly an element of $K\cap G_a$. Moreover, the uniqueness of lifts yields $\widetilde{f}=\widetilde{g}*l\widetilde{g}*\dots*l^{n-1}\widetilde{g}$ and hence $k=\widetilde{f}(1)=(\widetilde{g}*l\widetilde{g}*\dots *l^{n-1}\widetilde{g})(1)=l^n\in n(K\cap G_a)$. This verifies the inclusion $\mathcal E_K(n\pi_1(G/K))\subseteq n(K\cap G_a)$. To verify the converse inclusion let $k\in n(K\cap G_a)$ and write $k=l^n$ with $l\in K\cap G_a$. Choose a path $\widetilde{g}$ in $G$ from $e$ to $l$ and set $g=p_K\widetilde{g}\in\pi_1(G/K)$. Then $\widetilde{g}*l\widetilde{g}*\dots *l^{n-1}\widetilde{g}$ is the lift of $ng$ across $p_K$ starting at $e$ and its endpoint is $l^n=k$. Thus, $k\in\mathcal E_K(n\pi_1(G/K))$, as was to be shown.
\end{example}

\begin{example}\label{E:ess.discr.solen}
We describe the operation $\mathcal E_K$ in the situation when $G$ is a solenoid and $K$ is an arbitrary closed totally disconnected subgroup of $G$. So let $S_{{\bf p}}$ be a solenoid and let $K\sbgp S_{{\bf p}}$ be totally disconnected. The dual group $(S_{{\bf p}}/K)^*$ of $S_{{\bf p}}/K$ is isomorphic to $K^{\perp}\sbgp(S_{{\bf p}})^*\cong\mathbb Z/{\bf p}$ and so it has rank $1$. Consequently, $S_{{\bf p}}/K$ is topologically isomorphic either to a solenoid or to the circle group $\mathbb T^1$. In the first case we have $\pi_1(S_{{\bf p}}/K)=0$ and hence $\mathcal E_K=e$ is trivial. In the second case every subgroup $Q$ of $\pi_1(S_{{\bf p}}/K)$ has the form $Q=n\pi_1(S_{{\bf p}}/K)$ for some integer $n$, in which case $\mathcal E_K(Q)=\mathcal E_K(n\pi_1(S_{{\bf p}}/K))=n(K\cap(S_{{\bf p}})_a)$ by Example~\ref{E:tow.ess.discr.solen}.
\end{example}

Before formulating our main result of this section, namely Theorem~\ref{T:F.and.pi}, we fix some more necessary notation. Let $G\in\mathsf{CAGp}$, $K\sbgp G$ be totally disconnected and $p_K\colon G\to G/K$ be the quotient morphism. Given a subgroup $Q$ of $\pi_1(G/K)$, we let $\overline{\mathcal E_K}(Q)$\index[symbol]{$\overline{\mathcal E_K}(Q)$} denote the closure in $K$ of the group $\mathcal E_K(Q)$. Since $\mathcal E_K(Q)=\mathcal E_K(Q+p_K^{\sharp}\pi_1(G))$ for every $Q$, the same identity holds also with $\mathcal E_K$ replaced by $\overline{\mathcal E_K}$. Observe that $\overline{\mathcal E_K}(Q)$ is the orbit closure of the identity $e$ under the monodromy action of $Q$ on $K$. Moreover, by Lemma~\ref{L:QK.EK.tot.disc}(iii), $\overline{\mathcal E_K}(\mathcal Q_K(E))=\overline{E\cap G_a}$ for every subgroup $E$ of $K$ and $\mathcal Q_K(\overline{\mathcal E_K}(Q))\supseteq Q$ for every subgroup $Q$ of $\pi_1(G/K)$ containing $p_K^{\sharp}\pi_1(G)$.

\begin{remark}\label{R:cl.of.EK.QK}
Let $G\in\mathsf{CAGp}$ and $K\sbgp G$ be totally disconnected. Assume that $K\subseteq G_a$. (This is true, for instance, if $G$ is a (finite- or infinite-dimensional) torus, in which case $G_a=G$.) Then, by Lemma~\ref{L:QK.EK.tot.disc}(iii), the following statements hold:
\begin{enumerate}
\item[(a)] $\mathcal E_K(\mathcal Q_K(E))=E$ for every subgroup $E$ of $K$,
\item[(b)] $\mathcal Q_K(\mathcal E_K(Q))=Q$ for every subgroup $Q$ of $\pi_1(G/K)$ containing $p_K^{\sharp}\pi_1(G)$.
\end{enumerate}
Consequently, $\mathcal E_K$ and $\mathcal Q_K$ consitute one-to-one correspondences between the subgroups $E$ of $K$ on one side and the subgroups $Q$ of $\pi_1(G/K)$ containing $p_K^{\sharp}\pi_1(G)$ on the other side. As a consequence of these facts, the following observations emerge:
\begin{itemize}
\item $\mathcal E_K(Q)$ need not be a closed subgroup of $K$; to see this, take $Q=\mathcal Q_K(E)$ with $E$ a non-closed subgroup of $K$,
\item the equality $\mathcal Q_K(\overline{\mathcal E_K}(Q))=Q$ may fail; to see this, take, as above, $Q=\mathcal Q_K(E)$ with $E$ a non-closed subgroup of $K$, then
\begin{equation*}
\mathcal Q_K(\overline{\mathcal E_K}(Q))=\mathcal Q_K(\overline{\mathcal E_K}(\mathcal Q_K(E)))=\mathcal Q_K(\overline{E})\supsetneq\mathcal Q_K(E)=Q.
\end{equation*}
\end{itemize}
\end{remark}

\begin{theorem}\label{T:F.and.pi}
Let $\Flow$ be a minimal flow, let $G\in\mathsf{CAGp}$ and $\mathcal C\in\Cocc(G)$, and fix $z\in X$. Assume that $K\sbgp G$ is totally disconnected with $p_K\mathcal C\in\Cob(G/K)$, and denote by $\xi$ the transfer function of $p_K\mathcal C$ with $\xi(z)=e$. If $X\in\mathsf{CLAC}$ and $\Gamma$ has no non-trivial finite abelian quotient groups then 
\begin{equation}\label{Eq:F.when.K.tot.disc}
F(\mathcal C)=\overline{\mathcal E_K}\left(\xi^{\sharp}\pi_1(X)+p_K^{\sharp}\pi_1(G) \right)=\overline{\mathcal E_K}\left(\xi^{\sharp}\pi_1(X)\right).
\end{equation}
If, in addition, the group $K$ is finite then the following statements hold:
\begin{enumerate}
\item[(1)] $F(\mathcal C)$ is given by the formula
\begin{equation*}
F(\mathcal C)=\mathcal E_K\left(\xi^{\sharp}\pi_1(X)+p_K^{\sharp}\pi_1(G) \right)=\mathcal E_K\left(\xi^{\sharp}\pi_1(X)\right),
\end{equation*}
\item[(2)] $F(\mathcal C)$ satisfies
\begin{equation*}
\mathcal Q_K\left(F(\mathcal C)\right)=\xi^{\sharp}\pi_1(X)+p_K^{\sharp}\pi_1(G) ,
\end{equation*}
\item[(3)] $F(\mathcal C)\subseteq G_a$ and
\begin{equation*}
F(\mathcal C)\cong\left(\xi^{\sharp}\pi_1(X)+p_K^{\sharp}\pi_1(G)\right) \Big/p_K^{\sharp}\pi_1(G)\cong \pi_1(G/F(\mathcal C))\Big/p_{F(\mathcal C)}^{\sharp}\pi_1(G).
\end{equation*}
\item[(4)] the following statements are equivalent:
\begin{enumerate}
\item[(i)] $F(\mathcal C)=K$,
\item[(ii)] $K\subseteq G_a$ and $\xi^{\sharp}\pi_1(X)+p_K^{\sharp}\pi_1(G)=\pi_1(G/K)$.
\end{enumerate}
\end{enumerate}
\end{theorem}
\begin{remark}\label{R:F.and.pi}
We wish to mention the following facts.
\begin{itemize}
\item In the isomorphism (\ref{Eq:F.when.K.tot.disc}) the symbol $\overline{\mathcal E_K}$ can not be replaced by the symbol $\mathcal E_K$, for the group $\xi^{\sharp}\pi_1(X)+p_K^{\sharp}\pi_1(G)$ need not correspond to a closed subgroup of $K$ via the operation $\mathcal E_K$, see Remark~\ref{R:F.and.pi.E.not.clsd} below.
\item If the group $K$ is not finite then $F(\mathcal C)$ may fail to be a subgroup of $G_a$ and so statement (3) is no longer true. Again, we leave an example of such a situation to Remark~\ref{R:F.and.pi.E.not.clsd}.
\item An immediate consequence of (\ref{Eq:F.when.K.tot.disc}) is that the section $F(p_K^{-1}\co(\xi))$ of the lift $p_K^{-1}\co(\xi)$ of $\co(\xi)$ across $p_K$ depends only on the homotopy class of $\xi\in C_z(X,G/K)$. Consequently, if $ C_z(X,G/K)_a$ is an open subgroup of $C_z(X,G/K)$, which is the case if $X$ is a compact space and $G/K$ is a finite-dimensional torus (see Subsection~\ref{Sub:frst.chmtp.gp}), then $F(p_K^{-1}\co(\xi))$ is a locally constant (and hence continuous) function of the variable $\xi\in C_z(X,G/K)$.
\end{itemize}
\end{remark}
\begin{proof}[Proof of Theorem~\ref{T:F.and.pi}]
Let $E$ be a closed subgroup of $K$. Recall the notation $p_K$ and $p_E$ for the quotient morphisms $G\to G/K$ and $G\to G/E$, respectively, and $q_E$ for the (unique quotient) morphism $G/E\to G/K$ with $q_Ep_E=p_K$. Since $E$ and $K$ are closed totally disconnected subgroups of $G$, Lemma~\ref{L:lift.tot.disc.ker} applies to the quotient morphisms $p_E$ and $p_K$. Moreover, it applies also to $q_E$, for $q_E$ is equivalent to the quotient morphism $G/E\to (G/E)/(K/E)$, whose kernel $K/E\cong p_E(K)$ is a closed totally disconnected subgroup of $G/E$.

We claim that $q_E^{\sharp}\pi_1(G/E)=\mathcal Q_K(G_a\cap E)$. To see this let $f$ be a loop in $G/K$ based at $e$ and denote by $\widetilde{f}$ the lift of $f$ across $p_K$ starting at $e$. Then, clearly, $p_E\widetilde{f}$ is the lift of $f$ across $q_E$ starting at $e$. Thus, the equivalence of the following statements follows:
\begin{itemize}
\item $f\in q_E^{\sharp}\pi_1(G/E)$,
\item $p_E\widetilde{f}$ is a loop,
\item $\widetilde{f}$ ends at an element of $E$,
\item $f\in\mathcal Q_K(E)$,
\item $f\in\mathcal Q_K(E\cap G_a)$,
\end{itemize}
and that verifies the claim.

Further, we claim that the following statements are equivalent:
\begin{enumerate}
\item[(a)] $\xi$ lifts across $q_E$ to a continuous base point preserving map $\eta\colon X\to G/E$,
\item[(b)] $\xi^{\sharp}\pi_1(X)\subseteq q_E^{\sharp}\pi_1(G/E)$,
\item[(c)] $\xi^{\sharp}\pi_1(X)\subseteq Q_K(G_a\cap E)$.
\end{enumerate}
Indeed, (a) and (b) are equivalent by virtue of Lemma~\ref{L:lift.tot.disc.ker} and the equivalence of (b) and (c) follows from the equality $q_E^{\sharp}\pi_1(G/E)=\mathcal Q_K(G_a\cap E)$ verified above. 

Now, from statement (i) in Proposition~\ref{P:char.F.as.smallest} and from the equivalence of (a) and (c) it follows that $F(\mathcal C)$ equals the smallest of all the closed subgroups $E$ of $K$ with $\xi^{\sharp}\pi_1(X)\subseteq Q_K(G_a\cap E)$. Since $\mathcal Q_K(G_a\cap E)=q_E^{\sharp}\pi_1(G/E)\supseteq p_K^{\sharp}\pi_1(G)$, the inclusion $\xi^{\sharp}\pi_1(X)\subseteq Q_K(G_a\cap E)$ is equivalent with $\xi^{\sharp}\pi_1(X)+p_K^{\sharp}\pi_1(G)\subseteq\mathcal Q_K(G_a\cap E)$, and the latter inclusion is equivalent with
\begin{equation}\label{Eq:mndrm.pt.1}
\mathcal E_K\left(\xi^{\sharp}\pi_1(X)+p_K^{\sharp}\pi_1(G) \right)\subseteq G_a\cap E
\end{equation}
by virtue of Lemma~\ref{L:QK.EK.tot.disc}(iii). Since the left-hand side of (\ref{Eq:mndrm.pt.1}) is contained in $G_a$ by Lemma~\ref{L:QK.EK.tot.disc}(ii), the smallest closed subgroup $E$ of $K$ for which (\ref{Eq:mndrm.pt.1}) holds is the group
\begin{equation*}
F(\mathcal C)=\overline{\mathcal E_K}\left(\xi^{\sharp}\pi_1(X)+p_K^{\sharp}\pi_1(G) \right)=\overline{\mathcal E_K}\left(\xi^{\sharp}\pi_1(X)\right).
\end{equation*}
This verifies (\ref{Eq:F.when.K.tot.disc}). 

Assume, from now on, that the group $K$ is finite. Then every subgroup $E$ of $K$ is closed and $\mathcal E_K$ therefore coincides with $\overline{\mathcal E_K}$. Thus, statement (1) follows at once from (\ref{Eq:F.when.K.tot.disc}). Statement (2) follows from statement (1) by applying $\mathcal Q_K$ and by using Lemma~\ref{L:QK.EK.tot.disc}.

The inclusion $F(\mathcal C)\subseteq G_a$ in statement (3) follows from statement (1) and from Lemma~\ref{L:QK.EK.tot.disc}(ii). The first isomorphism from (3) follows by applying Lemma~\ref{L:QK.EK.tot.disc}(iv) to $Q=\xi^{\sharp}\pi_1(X)+p_K^{\sharp}\pi_1(G)$. The second isomorphism from (3) is proved similarly, this time one applies Lemma~\ref{L:QK.EK.tot.disc}(iv) to $Q=\pi_1(G/F(\mathcal C))$, with $K$ replaced by $F(\mathcal C)$, to obtain
\begin{equation*}
F(\mathcal C)=F(\mathcal C)\cap G_a=\mathcal E_{F(\mathcal C)}(\pi_1(G/F(\mathcal C)))\cong\pi_1(G/F(\mathcal C))/p_{F(\mathcal C)}^{\sharp}\pi_1(G).
\end{equation*}

We finish the proof of the theorem by verifying the equivalence of (i) and (ii) from statement (4). Assume, first, that $F(\mathcal C)=K$. Then $K\subseteq G_a$ by virtue of (3). Moreover, from statement (1) and from Lemma~\ref{L:QK.EK.tot.disc}(iii) it follows that
\begin{equation*}
\pi_1(G/K)=\mathcal Q_K(K)=\mathcal Q_K\mathcal E_K(\xi^{\sharp}\pi_1(X)+p_K^{\sharp}\pi_1(G))=\xi^{\sharp}\pi_1(X)+
p_K^{\sharp}\pi_1(G),
\end{equation*}
which shows that (i) implies (ii). To verify that (ii) implies (i) we invoke statement (1) once more to obtain
\begin{equation*}
F(\mathcal C)=\mathcal E_K(\xi^{\sharp}\pi_1(X)+p_K^{\sharp}\pi_1(G))=\mathcal E_K(\pi_1(G/K))=K\cap G_a=K.
\end{equation*}
This verifies statement (4).
\end{proof}

Before formulating our next result let us recall some facts and notation. Given $G\in\mathsf{CAGp}$ connected second countable, there is a topological direct sum $G=\mathfrak{f}(G)\oplus\mathfrak{t}(G)$, where $\mathfrak{f}(G)$ is a torus and $\mathfrak{t}(G)$ is a torus-free group (that is, $\Hom(\mathbb T^1,\mathfrak{t}(G))=0$). The group $\mathfrak{f}(G)$ is unique and is called the maximal torus of $G$. The group $\mathfrak{t}(G)$ is unique up to a topological isomorphism. Given, in addition, a topological space $X$ and a continuous map $\xi\colon X\to G$, we write $\mathfrak{f}(\xi)$ and $\mathfrak{t}(\xi)$ for the projections of $\xi$ to $\mathfrak{f}(G)$ and $\mathfrak{t}(G)$, respectively. The maps $\mathfrak{f}(\xi)$, $\mathfrak{t}(\xi)$ are continuous and satisfy $\xi=\mathfrak{f}(\xi)\oplus\mathfrak{t}(\xi)$. If $X$ is a pointed space and the map $\xi$ is base point preserving, then so are the maps $\mathfrak{f}(\xi)$ and $\mathfrak{t}(\xi)$.

\begin{corollary}\label{C:F.and.pi.fprt}
Under the assumptions of Theorem~\ref{T:F.and.pi}, if $X$ is a compact connected manifold and $G\in\mathsf{CAGp}$ is connected second countable then
\begin{equation}\label{Eq:F.and.pi.fprt1}
F(\mathcal C)=\overline{\mathcal E_K}\left(\mathfrak{f}(\xi)^{\sharp}\pi_1(X)\right)\subseteq\left(p_K^{-1}\mathfrak{f}(G/K)\right)_0.
\end{equation}
If, in addition, the group $K$ is finite then
\begin{equation}\label{Eq:F.and.pi.fprt2}
F(\mathcal C)\subseteq\mathfrak{f}(G)\cap\left(p_K^{-1}\mathfrak{f}(G/K)\right)_a.
\end{equation}
\end{corollary}
\begin{remark}\label{R:F.and.pi.fprt}
It follows from the corollary that if $\mathfrak{f}(G/K)=e$ (that is, if $G/K$ contains no torus) then $F(\mathcal C)=e$ and hence $\mathcal C\in\Cob(G)$. Indeed, by virtue of (\ref{Eq:F.and.pi.fprt1}),
\begin{equation*}
F(\mathcal C)\subseteq\left(p_K^{-1}\mathfrak{f}(G/K)\right)_0=\left(p_K^{-1}e\right)_0=K_0=e.
\end{equation*}
We shall have an occasion to use this observation in the situation when $G$ is a solenoid. In such a case every quotient group $G/K$ of $G$ modulo a totally disconnected subgroup $K\sbgp G$ is either a solenoid or else is isomorphic to $\mathbb T^1$. In the first case $\mathfrak{f}(G/K)=e$ and hence $F(\mathcal C)=e$. Thus, in order to obtain a non-trivial totally disconnected section $F(\mathcal C)$, it is necessary to have $G/K\cong\mathbb T^1$ or, equivalently, $K^{\perp}\cong\mathbb Z$, where $K^{\perp}$ stands for the annihilator of $K$ in $G^*$.
\end{remark}
\begin{proof}[Proof of Corollary~\ref{C:F.and.pi.fprt}]
We begin the proof by showing that the maps $\xi,\mathfrak{f}(\xi)\colon X\to G/K$ are homotopic. To this end, it is sufficient to verify that the map $\mathfrak{t}(\xi)\colon X\to\mathfrak{t}(G/K)$ is null-homotopic. Since the group $\pi^1(\mathfrak{t}(G/K))\cong\mathfrak{t}(G/K)^*\cong\mathfrak{t}((G/K)^*)$ is torsion-less (that is, it has no non-trivial free abelian quotient groups) and the group $\pi^1(X)\cong H_1^w(X)$ is free abelian by our assumptions on $X$ (see Lemma~\ref{L:(co)hmtp.rel} for more details), it follows that the induced morphism $\mathfrak{t}(\xi)^{\flat}\colon\pi^1(\mathfrak{t}(G/K))\to\pi^1(X)$ vanishes. Consequently, the isomorphism (\ref{Eq:maps.are.homs.hmtp}) from Subsection~\ref{Sub:chmtp.cpt.gps} yields that the map $\mathfrak{t}(\xi)$ is null-homotopic, as was to be shown. This shows that $\xi$ and $\mathfrak{f}(\xi)$ are indeed homotopic and so the induced morphisms $\xi^{\sharp},\mathfrak{f}(\xi)^{\sharp}\colon\pi_1(X)\to\pi_1(G/K)$ coincide. Consequently, by virtue of (\ref{Eq:F.when.K.tot.disc}) from Theorem~\ref{T:F.and.pi}, 
\begin{equation}\label{Eq:aux.FC.FD}
F(\mathcal C)=\overline{\mathcal E_K}\left(\xi^{\sharp}\pi_1(X)\right)=\overline{\mathcal E_K}\left(\mathfrak{f}(\xi)^{\sharp}\pi_1(X)\right),
\end{equation}
which verifies the first part of (\ref{Eq:F.and.pi.fprt1}). Further, since $\mathfrak{f}(\xi)$ takes its values in $\mathfrak{f}(G/K)$, the group $\mathcal E_K(\mathfrak{f}(\xi)^{\sharp}\pi_1(X))$ is contained in the identity arc-component $(p_K^{-1}\mathfrak{f}(G/K))_a$ of $p_K^{-1}\mathfrak{f}(G/K)$ by definition of $\mathcal E_K$. Thus, by virtue of (\ref{Eq:aux.FC.FD}),
\begin{equation*}
F(\mathcal C)=\overline{\mathcal E_K\left(\mathfrak{f}(\xi)^{\sharp}\pi_1(X)\right)}\subseteq\overline{\left(p_K^{-1}\mathfrak{f}(G/K)\right)_a}=\left(p_K^{-1}\mathfrak{f}(G/K)\right)_0,
\end{equation*}
which verifies the second part of (\ref{Eq:F.and.pi.fprt1}). 

Now assume that the group $K$ is finite. Then, by Theorem~\ref{T:F.and.pi}(1),
\begin{equation*}
F(\mathcal C)=\mathcal E_K\left(\xi^{\sharp}\pi_1(X)\right)=\mathcal E_K\left(\mathfrak{f}(\xi)^{\sharp}\pi_1(X)\right)\subseteq\left(p_K^{-1}\mathfrak{f}(G/K)\right)_a.
\end{equation*}
Finally, to see that $F(\mathcal C)\subseteq\mathfrak{f}(G)$, fix $g\in F(\mathcal C)$. Since $F(\mathcal C)\subseteq G_a=\exp(\mathcal L(G))$ by Theorem~\ref{T:F.and.pi}(3), there exist $h\in\Hom(\mathbb R,G)$ and $t\in\mathbb R\setminus0$ with $g=h(t)$. Let $d$ be the largest torsion coefficient of $K$. Then $e=g^d=h(dt)$ and hence $dt\mathbb Z\subseteq\ker(h)$. Consequently, there is $r\in\Hom(\mathbb T^1,G)$ with $\im(r)=\im(h)\ni g$. Since $\im(r)\subseteq\mathfrak{f}(G)$, it follows that $g\in\mathfrak{f}(G)$. This verifies (\ref{Eq:F.and.pi.fprt2}).
\end{proof}

As promised earlier, we conclude this section by justifying the first two statements from Remark~\ref{R:F.and.pi}.

\begin{proposition}\label{P:F.and.pi.E.not.clsd}
Let $G\in\mathsf{CAGp}$ and $K\sbgp G$ be a totally disconnected subgroup of $G$. Given a (not necessarily closed) subgroup $E$ of $K$, the following conditions are equivalent:
\begin{enumerate}
\item there exist a minimal flow $\mathcal F\colon\mathbb R\curvearrowright X$ on a compact (pointed) manifold $X$ and an extension $\mathcal C\in\Cocc(G)$ with $p_K\mathcal C\in\Cob(G/K)$ such that
\begin{equation*}
E=\mathcal E_K\left(\xi^{\sharp}\pi_1(X)+p_K^{\sharp}\pi_1(G)\right),
\end{equation*}
where $\xi$ is the transfer function of $p_K\mathcal C$ with $\xi(z)=e$,
\item $E$ is finitely generated in the algebraic sense and is contained in $G_a$.
\end{enumerate}
\end{proposition}
\begin{remark}\label{R:F.and.pi.E.not.clsd}
In Remark~\ref{R:F.and.pi} we made two claims, which we now justify using Proposition~\ref{P:F.and.pi.E.not.clsd}. Our fist claim concerned the fact that in the situation of Theorem~\ref{T:F.and.pi}, the group $\xi^{\sharp}\pi_1(X)+p_K^{\sharp}\pi_1(G)$ need not correspond to a closed subgroup of $K$ under $\mathcal E_K$, and hence the operations $\mathcal E_K$ and $\overline{\mathcal E_K}$ are distinct. Our second claim concerned the fact that the inclusion $F(\mathcal C)\subseteq G_a$ from statement (3) of Theorem~\ref{T:F.and.pi} may fail if the group $K$ is not finite. We present an example justifying both these statements. Let $G\sbgp\mathbb T^{\mathbb N}$ be the dyadic solenoid and $K$ be the set of all $(z_n)_{n\in\mathbb N}\in G$ with $z_1=1$. Clearly, $K\sbgp G$ is totally disconnected. Recall that the identity arc-component $G_a$ of $G$ consists of the sequences $(z_n)_{n\in\mathbb N}\in G$ satisfying $z_n\to1$ in $\mathbb T^1$ as $n\to\infty$. Consequently, the following statements hold:
\begin{itemize}
\item $K$ is not contained in $G_a$,
\item $K\cap G_a$ contains a topological generator $k$ of $K$.
\end{itemize}
Let $E=\langle k\rangle$ be the subgroup of $K$ generated by $k$ in the algebraic sense. Then $E$ satisfies both assumptions from statement (2) in Proposition~\ref{P:F.and.pi.E.not.clsd}, and so there exist $\mathcal F$, $\mathcal C$ and $\xi$ as in statement (1) of Proposition~\ref{P:F.and.pi.E.not.clsd} such that $E=\mathcal E_K(\xi^{\sharp}\pi_1(X)+p_K^{\sharp}\pi_1(G))$. Further, by virtue of Theorem~\ref{T:F.and.pi}, $F(\mathcal C)=\overline{\mathcal E_K}(\xi^{\sharp}\pi_1(X)+p_K^{\sharp}\pi_1(G))=\overline{E}=K$. We thus have the following conclusions:
\begin{enumerate}
\item[(i)] the group $\mathcal E_K(\xi^{\sharp}\pi_1(X)+p_K^{\sharp}\pi_1(G))$ is not closed in $K$; indeed, we have $E\subsetneq\overline{E}$, for $E$ is contained in $G_a$, while $\overline{E}=K$ is not,
\item[(ii)] the group $F(\mathcal C)=K$ is not contained in $G_a$.
\end{enumerate}
\end{remark}
\begin{proof}[Proof of Proposition~\ref{P:F.and.pi.E.not.clsd}]
Assume that condition (1) holds. Since $X$ is a compact connected manifold, its fundamental group $\pi_1(X)$ is finitely generated and hence so is the group $\xi^{\sharp}\pi_1(X)$. Let $f_1,\dots,f_m$ be a family of generators for $\xi^{\sharp}\pi_1(X)$. Further, for $1\leq i\leq m$, let $k_i$ be the endpoint of the lift of $f_i$ across $p_K$ starting at $e$. Then for every $f\in\xi^{\sharp}\pi_1(X)$, the endpoint of the lift of $f$ across $p_K$ starting at $e$ is expressible as a product of the integer powers of the elements $k_1,\dots,k_m$. This shows that $E\subseteq\langle k_1,\dots,k_m\rangle$ is finitely generated. Also, the inclusion $E\subseteq G_a$ follows from Lemma~\ref{L:QK.EK.tot.disc}(ii) and so (1) implies (2).

Now assume that condition (2) holds and let $k_1,\dots,k_m$ be a family of algebraic generators of $E$. Since $E\subseteq G_a$ by the assumptions, it follows that for every $i=1,\dots,m$ there is $f_i\in\pi_1(G/K)$ such that the lift of $f_i$ across $p_K$ starting at $e$ ends at $k_i$. By \cite[Theorem~8.57,~p.~420]{HofMor} we may assume that $f_i$ is of the form $f_i=h_ip$, where $p\colon[0,1]\to\mathbb T^1$ is the restriction of the standard covering morphism $\mathbb R\to\mathbb T^1$ and $h_i\in\Hom(\mathbb T^1,G/K)$. Clearly, $\mathcal E_K(\langle f_1,\dots,f_m\rangle)=\langle k_1,\dots,k_m\rangle=E$.

Let $X=\mathbb T^m$ have as a base point its identity $e$ and consider the base point preserving map $\xi\colon\mathbb T^m\to G/K$ be defined by $\xi(z_1,\dots,z_m)=h_1(z_1)\dots h_m(z_m)$. We claim that $\mathcal E_K(\xi^{\sharp}\pi_1(X)+p_K^{\sharp}\pi_1(G))=E$. Given $i\in\{1,\dots,m\}$, let $p_i=(1,\dots,1,p,1,\dots,1)\in\pi_1(\mathbb T^m)$ be the loop having $p$ as its $i^{\text{th}}$ coordinate. Then, clearly, $\pi_1(\mathbb T^m)=\langle p_1,\dots,p_m\rangle$. Moreover, $\xi^{\sharp}(p_i)=\xi p_i=h_ip=f_i$ for every $i$ and hence $\xi^{\sharp}\pi_1(\mathbb T^m)=\langle f_1,\dots,f_m\rangle$. Consequently,
\begin{equation*}
\mathcal E_K\left(\xi^{\sharp}\pi_1(\mathbb T^m)+p_K^{\sharp}\pi_1(G)\right)=\mathcal E_K\left(\xi^{\sharp}\pi_1(\mathbb T^m)\right)=\mathcal E_K(\langle f_1,\dots,f_m\rangle)=E.
\end{equation*}

To finish the proof of (1), let $\mathcal F$ be, say, an irrational continuous flow on $\mathbb T^m$ and consider the coboundary $\co(\xi)\in\Cob(G/K)$ over $\mathcal F$ with the transfer function $\xi$. Since the acting group $\mathbb R$ of $\mathcal F$ is simply connected, we have $\xi^{\sharp}\mathcal F_z^{\sharp}\pi_1(\mathbb R)=0$. Thus, by Corollary~\ref{C:lift.cob.Lie.gp}, $\co(\xi)$ lifts across $p_K$ to an extension $\mathcal C\in\Cocc(G)$. Hence $\mathcal F$ and $\mathcal C$ are as in (1), which shows that (2) implies (1).
\end{proof}

\section{Existence of totally disconnected sections}\label{S:ex.prescr.sect.1}

Given a minimal flow $\Flow$, a group $G\in\mathsf{CAGp}$ and a closed subgroup $K$ of $G$, does there exist an extension $\mathcal C\in\Cocc(G)$ with $F(\mathcal C)=K$? In this section we use our results from Section~\ref{S:rel.mndrm.act} to give a necessary and sufficient condition for the existence of such an extension $\mathcal C$ in the case when the group $K$ is totally disconnected (see Thoerem~\ref{T:F.and.Betti}) and, in particular, when $K$ is finite (see Corollary~\ref{C:F.and.Betti}). The condition involves ranks of the groups $H_1^w(\mathcal F)$ and $H_1^w(X)$, the elementary divisors of $H_1^w(\mathcal F)$ in $H_1^w(X)$, and the number and orders of topological generators of $K$ lying on the arc-component $G_a$ of $G$.

\begin{lemma}\label{L:(co)hmtp.rel}
Let $X$ be a pointed compact connected manifold with the base point $z$. Consider the map $\lambda\colon\pi_1(X)\ni l\mapsto \lambda(l)\in\Hom(\pi^1(X),\pi_1(\mathbb T^1))$ with $\lambda(l)\colon\pi^1(X)\ni f\mapsto f^{\sharp}(l)=fl\in\pi_1(\mathbb T^1)$. Then $\lambda$ is an epimorphism of groups with the kernel $\ker(\lambda)=\ker(p_X)$, where $p_X\colon\pi_1(X)\to H_1^w(X)$ is the canonical quotient morphism. Thus, $\lambda$ induces an isomorphism $\nu\colon H_1^w(X)\to\Hom(\pi^1(X),\pi_1(\mathbb T^1))$ with $\nu p_X=\lambda$.
\end{lemma}
\begin{proof}
Before turning to the proof let us mention that we shall use the standard identifications $\pi_1(\mathbb T^1)\cong H_1^w(\mathbb T^1)\cong\mathbb Z$ and that the group operation of $\pi^1(X)$ will be written additively. Moreover, given $f\in C_z(X,\mathbb T^1)$, we shall use the same symbol $f^{\sharp}$ for the induced morphisms $\pi_1(X)\to\mathbb Z$ and $H_1^w(X)\to\mathbb Z$. Finally, a map $f\in C_z(X,\mathbb T^1)$ will be sometimes identified with its homotopy class in $\pi^1(X)$. We divide the proof of the lemma into three steps.

\emph{1st step.} We show that the mapping
\begin{equation}\label{Eq:mrphsm.R}
\sigma\colon\pi^1(X)\ni f\mapsto f^{\sharp}\in\Hom(\pi_1(X),\mathbb Z)
\end{equation}
is an isomorphism of groups.

Recall from Subsection~\ref{Sub:ind.mrph.hmlg.gp} that there is a morphism of groups
\begin{equation*}
\varrho\colon C_z(X,\mathbb T^1)\ni f\mapsto f^{\sharp}\in\Hom(\pi_1(X),\mathbb Z).
\end{equation*}
Since the induced morphism $f^{\sharp}\colon\pi_1(X)\to\mathbb Z$ depends only on the homotopy class of $f\in C_z(X,\mathbb T^1)$, $\varrho$ gives rise to a morphism $\sigma$ as in (\ref{Eq:mrphsm.R}). By our assumptions on $X$, $\varrho$ is an epimorphism and hence so is $\sigma$. To see that $\sigma$ is in fact an isomorphism observe that if $f\in C_z(X,\mathbb T^1)$ satisfies $\varrho(f)=f^{\sharp}=0$ then $f$ lifts across the usual covering morphism $p\colon\mathbb R\to\mathbb T^1$ to an element of $C_z(X,\mathbb R)$ and hence $f=0$ as an element of $\pi^1(X)$.

\emph{2nd step.} By our assumptions on $X$, $\pi_1(X)$ is a finitely generated group and $H_1^w(X)$ is thus a free abelian group with a finite rank $n$. Assume that $n\neq0$ and choose $l_1,\dots,l_n\in\pi_1(X)$ in such a way that $p_X(l_1),\dots,p_X(l_n)$ form a basis for $H_1^w(X)$. Use the first step of the proof to find $f_1,\dots,f_n\in\pi^1(X)$ with $f_i^{\sharp}(l_j)=\delta_{ij}$ for $i,j=1,\dots,n$. We show that $f_1,\dots,f_n$ form a basis for $\pi^1(X)$. (Observe that, by the first step of the proof,
\begin{equation*}
\pi^1(X)\cong\Hom(\pi_1(X),\mathbb Z)\cong\Hom(H_1^w(X),\mathbb Z)\cong\Hom(\mathbb Z^n,\mathbb Z)\cong\mathbb Z^n
\end{equation*}
is a free abelian group with rank $n$.)

To see that $f_1,\dots,f_n$ are independent in $\pi^1(X)$, let $k_1,\dots,k_n$ be integers with $\sum_{i=1}^nk_if_i=0$. By applying the isomorphism $\sigma$ we obtain $\sum_{i=1}^nk_if_i^{\sharp}=0$. By our choice of the maps $f_i$ it follows at once that $k_1=\dots=k_n=0$.

To see that $f_1,\dots,f_n$ generate $\pi^1(X)$, fix $f\in\pi^1(X)$. For $i=1,\dots,n$ write $k_i=f^{\sharp}(l_i)$ and set $g=\sum_{i=1}^nk_if_i\in\pi^1(X)$. Then $g^{\sharp}(l_j)=k_j=f^{\sharp}(l_j)$ for every $j=1,\dots,n$ and hence, by our choice of the loops $l_1,\dots,l_n$, $\sigma(g)=g^{\sharp}=f^{\sharp}=\sigma(f)$. Since $\sigma$ is an isomorphism by the first step of the proof, it follows that $f=g=\sum_{i=1}^nk_if_i$.

\emph{3rd step.} We finish the proof of the lemma by showing that $\lambda$ is an epimorphism with the kernel $\ker(\lambda)=\ker(p_X)$.

First, it is clear that $\lambda(l)\in\Hom(\pi^1(X),\mathbb Z)$ for every $l\in\pi_1(X)$ and so $\lambda$ is well defined. An elementary argument also shows that $\lambda$ is a morphism of groups. To see that $\lambda$ is an epimorphism assume, first, that $H_1^w(X)=0$. Then $f^{\sharp}=0$ for every $f\in C_z(X,\mathbb T^1)$ and hence $\pi^1(X)=0$. Thus, in this case, $\lambda=0$ is indeed an epimorphism. So assume that $\rank(H_1^w(X))=n\neq0$ and fix $h\in\Hom(\pi^1(X),\mathbb Z)$. We invoke the notation introduced in the second step of the proof. Write $m_i=h(f_i)$ for $i=1,\dots,n$ and set $l=l_1^{m_1}*\dots *l_n^{m_n}\in\pi_1(X)$. Then 
\begin{equation*}
\lambda(l)(f_i)=f_i^{\sharp}(l)=f_i^{\sharp}(l_1^{m_1}*\dots *l_n^{m_n})=\sum_{j=1}^nm_jf_i^{\sharp}(l_j)=m_i=h(f_i)
\end{equation*}
for every $i=1,\dots,n$ and, since $f_1,\dots,f_n$ generate $\pi^1(X)$ by the second step of the proof, it follows that $h=\lambda(l)$. Thus, $\lambda$ is indeed an epimorphism.

We show that $\ker(\lambda)=\ker(p_X)$. To this end, fix $l\in\pi_1(X)$. By our choice of $l_1,\dots,l_n$, there are integers $k_1,\dots,k_n$ with $p_X(l)=\sum_{i=1}^nk_ip_X(l_i)$ in $H_1^w(X)$. Then the following conditions are equivalent by definition of $f_1,\dots,f_n$ and $l_1,\dots,l_n$:
\begin{itemize}
\item $l\in\ker(\lambda)$,
\item $f_i^{\sharp}(l)=0$ for every $i$,
\item $f_i^{\sharp}(p_X(l))=0$ for every $i$,
\item $k_i=0$ for every $i$,
\item $p_X(l)=0$,
\item $l\in\ker(p_X)$.
\end{itemize}
Thus, $\ker(\lambda)=\ker(p_X)$, as was to be shown.
\end{proof}

\begin{lemma}\label{L:real.of.morph.pi1}
Let $X$ be a pointed compact connected manifold with the base point $z$ and let $G\in\mathsf{CAGp}$ be connected. Given $h\in\Hom(\pi_1(X),\pi_1(G))$, there exists $\xi\in C_z(X,G)$ with $h=\xi^{\sharp}$. Such a map $\xi$ is unique up to a homotopy.
\end{lemma}
\begin{proof}
In the proof of the lemma we shall use the usual identification $\pi_1(\mathbb T^1)\cong\mathbb Z$. Further, recall from Subsection~\ref{Sub:chmtp.cpt.gps} that by connectedness of $G$, every element of $\pi^1(G)$ (considered as a homotopy class) contains one and only one character of $G$. Hence $G^*\cong\pi^1(G)$. Finally, we invoke the isomorphism $\mu\colon\pi_1(G)\to\Hom(G^*,\mathbb Z)$, acting by the rule $\mu(l)\colon\chi\mapsto\chi^{\sharp}(l)=\chi l$ for $l\in\pi_1(G)$ and $\chi\in G^*$.

Now we turn to the proof of the lemma, proceeding in four steps. Since the group $\pi_1(G)\cong\Hom(G^*,\mathbb Z)$ is abelian and torsion-free, there is an isomorphism of groups $\Hom(\pi_1(X),\pi_1(G))\cong\Hom(H_1^w(X),H_1^w(G))$, and the statement of the lemma is thus true if $H_1^w(X)=0$. We shall therefore assume throughout the whole proof that
\begin{equation*}
\rank(H_1^w(X))=n\neq0.
\end{equation*}

\emph{1st step.} By Lemma~\ref{L:(co)hmtp.rel}, the map $\lambda$ defined therein is an epimorphism of groups with the kernel $\ker(\lambda)=\ker(p_X)$. Since $\Hom(G^*,\mathbb Z)$ is an abelian torsion-free group, we have $\ker(\mu h)\supseteq\ker(p_X)=\ker(\lambda)$ and so there is a morphism $\kappa\colon\Hom(\pi^1(X),\mathbb Z)\to\Hom(G^*,\mathbb Z)$ with $\kappa\lambda=\mu h$. We show that $\kappa=\Hom(\varrho,\mathbb Z)$ for an appropriate $\varrho\in\Hom(G^*,\pi^1(X))$.

It follows from our proof of Lemma~\ref{L:(co)hmtp.rel} that $\pi^1(X)\cong H_1^w(X)$ is a free abelian group with rank $n\neq0$. Let $f_1,\dots,f_n$ be a basis for $\pi^1(X)$ and $\varphi_1,\dots,\varphi_n$ be the corresponding dual basis for $\Hom(\pi^1(X),\mathbb Z)$. (That is, $\varphi_i(f_j)=\delta_{ij}$ for $i,j=1,\dots,n$.) Write $h_i=\kappa(\varphi_i)\in\Hom(G^*,\mathbb Z)$ for $i=1,\dots,n$ and set $\varrho\colon G^*\ni\chi\mapsto\sum_{i=1}^nh_i(\chi)f_i\in\pi^1(X)$. Then $\varrho$ is clearly a morphism of groups and we claim that $\kappa=\Hom(\varrho,\mathbb Z)$. Given $j\in\{1,\dots,n\}$ and $\chi\in G^*$, we have
\begin{equation*}
\Hom(\varrho,\mathbb Z)(\varphi_j)\colon \chi\mapsto\varphi_j(\varrho(\chi))=\varphi_j\left(\sum_{i=1}^nh_i(\chi)f_i\right)=h_j(\chi)=\kappa(\varphi_j)(\chi).
\end{equation*}
Hence $\Hom(\varrho,\mathbb Z)(\varphi_j)=\kappa(\varphi_j)$ for every $j=1,\dots,n$ and this shows that $\Hom(\varrho,\mathbb Z)=\kappa$.

\emph{2nd step.} We show that there is $\xi\in C_z(X,G)$ with $\xi^{\flat}=\varrho$, where $\pi^1(G)$ is identified with $G^*$ in the usual way.

Recall from Subsection~\ref{Sub:mps.cpt.ab.gps} that there is an isomorphism of groups
\begin{equation}
C_z(X,G)\ni\xi\mapsto\xi^*\in\Hom(G^*,C_z(X,\mathbb T^1)),
\end{equation}
where $\xi^*$ acts by the rule $\xi^*(\chi)=\chi\xi\in C_z(X,\mathbb T^1)$ for $\chi\in G^*$. Since $\pi^1(X)$ can be viewed (by one of its equivalent definitions) as a subgroup of $C_z(X,\mathbb T^1)$, we have $\varrho=\xi^*$ for an appropriate $\xi\in C_z(X,G)$. Then $\varrho=\xi^*=\xi^{\flat}$, as was to be shown.

\emph{3rd step.} We show that the map $\xi$ from the second step of the proof satisfies the conclusion of the lemma.

By our choice of $\xi$ and by definition of $\kappa=\Hom(\varrho,\mathbb Z)$, we have $\Hom(\xi^{\flat},\mathbb Z)\lambda=\mu h$. Further, for every $l\in\pi_1(X)$,
\begin{equation*}
\Hom(\xi^{\flat},\mathbb Z)\lambda\colon l\mapsto\lambda(l)\xi^{\flat}\colon G^*\ni\chi\mapsto\chi^{\sharp}\xi^{\sharp}(l)\in\mathbb Z\cong\pi_1(\mathbb T^1)
\end{equation*}
and
\begin{equation*}
\mu h\colon l\mapsto\mu(h(l))\colon G^*\ni\chi\mapsto\chi^{\sharp}(h(l))\in\mathbb Z\cong\pi_1(\mathbb T^1).
\end{equation*}
To summarize, if $l\in\pi_1(X)$ then $\chi^{\sharp}\xi^{\sharp}(l)=\chi^{\sharp}h(l)$ holds for every $\chi\in G^*$, and hence, since the morphisms $\chi^{\sharp}$ ($\chi\in G^*$) separate the points of $\pi_1(G)$ (see Subsection~\ref{Sub:chmtp.cpt.gps}), it follows that $\xi^{\sharp}(l)=h(l)$. Thus, $\xi^{\sharp}=h$, as was to be shown.

\emph{4th step.} We finish the proof of the lemma by verifying the uniqueness statement.

First recall the following facts:
\begin{enumerate}
\item[(i)] there is an epimorphism $C_z(X,\mathbb T^1)\ni f\mapsto f^{\sharp}\in\Hom(\pi_1(X),\pi_1(\mathbb T^1))$, whose kernel consists of the group $C_z(X,\mathbb T^1)_0$ of the null-homotopic maps; consequently, if two maps $f_1,f_2\in C_z(X,\mathbb T^1)$ satisfy $f_1^{\sharp}=f_2^{\sharp}$ then they are homotopic,
\item[(ii)] there is an isomorphism $[X,G]\ni[\xi]\mapsto\xi^{\flat}\in\Hom(G^*,\pi^1(X))$; consequently, if $\xi_1,\xi_2\in C_z(X,G)$ are such that $\chi\xi_1,\chi\xi_2$ are homotopic for every $\chi\in G^*$, then $\xi_1,\xi_2$ are homotopic.
\end{enumerate}

Now let $\xi_1,\xi_2\in C_z(X,G)$ satisfy $\xi_1^{\sharp}=\xi_2^{\sharp}$. Then, for every $\chi\in G^*$, $(\chi\xi_1)^{\sharp}=\chi^{\sharp}\xi_1^{\sharp}=\chi^{\sharp}\xi_2^{\sharp}
=(\chi\xi_2)^{\sharp}$ and hence $\chi\xi_1$, $\chi\xi_2$ are homotopic by virtue of (i). Thus, by virtue of (ii), the maps $\xi_1,\xi_2$ are homotopic. This verifies the uniqueness statement of the lemma.
\end{proof}

\begin{theorem}\label{T:F.and.Betti}
Let $\Flow$ be a minimal flow with $\Gamma\in\mathsf{LieGp}$ connected and with $X$ a compact connected manifold. Write $n=\rank(H_1^w(\mathcal F))$, $n+m=\rank(H_1^w(X))$ and let $d_1,\dots,d_n$ be the elementary divisors of $H_1^w(\mathcal F)$ in $H_1^w(X)$. Assume that $G\in\mathsf{CAGp}$ is connected and $K\sbgp G$ is totally disconnected. Then the following conditions are equivalent:
\begin{enumerate}
\item[(1)] there is $\mathcal C\in\Cocc(G)$ with $F(\mathcal C)=K$,
\item[(2)] there exist $k_1,\dots,k_{n+m}\in K\cap G_a$ with $K=\overline{\langle k_1,\dots,k_{n+m}\rangle}$ and $k_i^{d_i}=e$ for $i=1,\dots,n$.
\end{enumerate}
\end{theorem}
\begin{remark}\label{R:F.and.Betti}
We wish to add the following observations.
\begin{itemize}
\item The points $k_1,\dots,k_{n+m}$ from condition (2) are not required to be independent in any sense (in particular, some of them may equal the identity $e$ of $K$). The condition requires that the group $K$ possess a subset $S$ of $G_a$ of cardinality at most $n+m$, which generates $K$ topologically. There are no restrictions on $m$ of these generators, but the remaining ones (if any) must be torsion elements of $K$ with orders dividing the elementary divisors of $H_1^w(\mathcal F)$ in $H_1^w(X)$.
\item It follows from the theorem that if the flow $\mathcal F$ does not possess a free cycle (that is, if $m=0$) then every totally disconnected section $F(\mathcal C)$ of an extension $\mathcal C\in\Cocc(G)$ is necessarily a finite subgroup of $G_a$. Moreover, in such a situation the group $F(\mathcal C)$ is a quotient group of $\mathbb Z_{d_1}\oplus\dots\oplus\mathbb Z_{d_n}$, see Corollary~\ref{C:F.and.Betti}.
\item If, on the other hand, we have $n=0$ (that is, $H_1^w(\mathcal F)=0$), then $K=F(\mathcal C)$ for some $\mathcal C\in\Cocc(G)$ if and only if $K\cap G_a$ has a subset with cardinality at most $m$, which generates $K$ topologically.
\item Assume that $n+m=0$, that is, $H_1^w(X)=0$. Given $\xi\in C_z(X,G/K)$, we have $\xi^{\sharp}\pi_1(X)=0$ and hence, by Lemma~\ref{L:lift.tot.disc.ker}, $\xi$ lifts across $p_K$ to an element of $C_z(X,G)$. Consequently, if $\mathcal C\in\Cocc(G)$ satisfies $F(\mathcal C)\subseteq K$ then $\mathcal C\in\Cob(G)$ by virtue of Lemma~\ref{L:cob.vs.lift}. Thus, in this case, condition (2) from Theorem~\ref{T:F.and.Betti} takes the form of the identity $K=e$. In the proof of the theorem we shall restrict to the (non-trivial) case of $H_1^w(X)\neq0$.
\end{itemize}
\end{remark}
\begin{proof}[Proof of Theorem~\ref{T:F.and.Betti}]
We divide the proof into three steps. To simplify notation we shall occasionally identify every loop $f\in\pi_1(X)$ with its underlying cycle $p_X(f)\in H_1^w(X)$. Also, by virtue of part four of Remark~\ref{R:F.and.Betti}, we shall assume that $H_1^w(X)\neq0$, that is, $n+m\neq0$.

\emph{1st step.} Let $p_K\colon G\to G/K$ denote the canonical quotient morphism. Fix loops $f_1,\dots,f_{n+m}\in\pi_1(X)$ in such a way that $p_X(f_1),\dots,p_X(f_{n+m})$ form a basis for $H_1^w(X)$ and $d_1p_X(f_1),\dots,d_np_X(f_n)$ form a basis for $H_1^w(\mathcal F)$. Given a continuous base point preserving map $\xi\colon X\to G/K$, let $k_i$ ($i=1,\dots,n+m$) be the endpoint of the lift of $\xi f_i$ across $p_K$ starting at $e$. (Recall from Lemma~\ref{L:lift.tot.disc.ker} and Remark~\ref{R:lift.tot.disc.ker} that the paths in $G/K$ do lift across the morphism $p_K$.) We claim that the following conditions are equivalent:
\begin{enumerate}
\item[($\alpha$)] $\co(\xi)$ lifts across $p_K$ to an extension $\mathcal C\in\Cocc(G)$,
\item[($\beta$)] $k_i^{d_i}=e$ for $i=1,\dots,n$.
\end{enumerate}

Given $i\in\{1,\dots,n\}$, let $\widetilde{f}_i$ be the lift of $\xi f_i$ across $p_K$ starting at $e$. Then $\widetilde{g}_i:=\widetilde{f}_i*k_i\widetilde{f}_i*\dots*k_i^{d_i-1}
\widetilde{f}_i$ is the lift of $d_i(\xi f_i)$ across $p_K$ starting at $e$ and its endpoint is $k_i^{d_i}$. Thus, by virtue of Corollary~\ref{C:lift.cob.Lie.gp}, Lemma~\ref{L:lift.tot.disc.ker} and Lemma~\ref{L:tor.Lie.mnfld}, the following conditions are equivalent:
\begin{itemize}
\item $\co(\xi)$ lifts across $p_K$ to an element of $\Cocc(G)$,
\item $\xi^{\sharp}H_1^w(\mathcal F)\subseteq p_K^{\sharp}H_1^w(G)$,
\item $\xi^{\sharp}(d_ip_X(f_i))\in p_K^{\sharp}H_1^w(G)$ for $i=1,\dots,n$,
\item $d_i\xi^{\sharp}(f_i)=\xi^{\sharp}(d_if_i)\in p_K^{\sharp}\pi_1(G)$ for $i=1,\dots,n$,
\item $d_i(\xi f_i)$ lifts across $p_K$ to a loop based at $e$ for $i=1,\dots,n$,
\item $\widetilde{g}_i(1)=e$ for $i=1,\dots,n$,
\item $k_i^{d_i}=e$ for $i=1,\dots,n$.
\end{itemize}
This verifies the equivalence of ($\alpha$) and ($\beta$).

\emph{2nd step.} We show that (2) follows from (1).

Let $\mathcal C\in\Cocc(G)$ satisfy $F(\mathcal C)=K$ and denote by $\xi$ the continuous base point preserving map $X\to G/K$ with $p_K\mathcal C=\co(\xi)$. For $i=1,\dots,n+m$ let $k_i$ be the endpoint of the lift of $\xi f_i$ across $p_K$ starting at $e$. Then $k_1,\dots,k_{n+m}\in K\cap G_a$ and, by the equivalence of ($\alpha$) and ($\beta$) from the first step of the proof, $k_i^{d_i}=e$ for $i=1,\dots,n$. Finally, since the group $\pi_1(G/K)$ is abelian torsion-free, we have $\xi^{\sharp}\pi_1(X)=\langle\xi^{\sharp}(f_1),\dots,\xi^{\sharp}(f_{n+m})\rangle$ and hence, by virtue of Theorem~\ref{T:F.and.pi},
\begin{equation*}
K=F(\mathcal C)=\overline{\mathcal E_K}(\xi^{\sharp}\pi_1(X))=\overline{\mathcal E_K}(\langle \xi^{\sharp}(f_1),\dots,\xi^{\sharp}(f_{n+m})\rangle)=\overline{\langle k_1,\dots,k_{n+m}\rangle},
\end{equation*}
which verifies condition (2).

\emph{3rd step.} We show that (1) follows from (2).

Fix elements $k_1,\dots,k_{n+m}$ of $K$ as in (2). For $i=1,\dots,n+m$ let $\widetilde{f}_i$ be a path in $G$ from $e$ to $k_i$. Since $H_1^w(X)$ is a free abelian group with basis $f_1,\dots,f_{n+m}$ and the group $\pi_1(G/K)\cong H_1^w(G/K)$ is abelian, there is a morphism $h\colon\pi_1(X)\to\pi_1(G/K)$ with $h(f_i)=p_K\widetilde{f}_i$ for $i=1,\dots,n+m$. By Lemma~\ref{L:real.of.morph.pi1}, $h=\xi^{\sharp}$ for an appropriate continuous base point preserving map $\xi\colon X\to G/K$. Clearly, $k_i$ ($i=1,\dots,n+m$) is the endpoint of the lift of $\xi f_i$ across $p_K$ starting at $e$. Moreover, given $i\in\{1,\dots,n\}$, we have $k_i^{d_i}=e$ by the assumptions. Consequently, it follows from the equivalence of ($\alpha$) and ($\beta$) from the first step of the proof that $\co(\xi)$ lifts across $p_K$ to an extension $\mathcal C\in\Cocc(G)$. By virtue of Theorem~\ref{T:F.and.pi},
\begin{equation*}
F(\mathcal C)=\overline{\mathcal E_K}\left(\xi^{\sharp}\pi_1(X)\right)=\overline{\mathcal E_K}(\langle p_K\widetilde{f}_1,\dots,p_K\widetilde{f}_{n+m}\rangle)=\overline{\langle k_1,\dots,k_{n+m}\rangle}=K.
\end{equation*}
This verifies condition (1).
\end{proof}

\begin{corollary}\label{C:F.and.Betti}
Let $\Flow$ be a minimal flow with $\Gamma\in\mathsf{LieGp}$ connected and with $X$ a compact connected manifold. Write $n=\rank(H_1^w(\mathcal F))$, $n+m=\rank(H_1^w(X))$ and let $d_1,\dots,d_n$ be the elementary divisors of $H_1^w(\mathcal F)$ in $H_1^w(X)$. Assume that $G\in\mathsf{CAGp}$ is connected and $K\sbgp G$ is finite. Then the following conditions are equivalent:
\begin{enumerate}
\item[(i)] there is $\mathcal C\in\Cocc(G)$ with $F(\mathcal C)=K$,
\item[(ii)] $K\subseteq G_a$ is a (finite) quotient group of $\mathbb Z_{d_1}\oplus\dots\oplus\mathbb Z_{d_n}\oplus\mathbb Z^m$.
\end{enumerate}
\end{corollary}
\begin{proof}
Similarly to the proof of Theorem~\ref{T:F.and.Betti}, we may restrict to the case of $H_1^w(X)\neq0$. Assume that $\mathcal C\in\Cocc(G)$ satisfies $F(\mathcal C)=K$. Then, by Theorem~\ref{T:F.and.pi}(3), $K\subseteq G_a$. Moreover, by Theorem~\ref{T:F.and.Betti}, $K$ possesses a family of (algebraic) generators $k_1,\dots,k_{n+m}$ with $k_i^{d_i}=e$ for $i=1,\dots,n$. This means that $K$ is a quotient group of $\mathbb Z_{d_1}\oplus\dots\oplus\mathbb Z_{d_n}\oplus\mathbb Z^m$. Thus, (ii) follows from (i).

Now assume that $K\subseteq G_a$ is a quotient group of $\mathbb Z_{d_1}\oplus\dots\oplus\mathbb Z_{d_n}\oplus\mathbb Z^m$. Then there is a family of generators $k_1,\dots,k_{n+m}\in G_a$ for $K$ with $k_i^{d_i}=e$ for $i=1,\dots,n$. By applying Theorem~\ref{T:F.and.Betti} once more we see that $F(\mathcal C)=K$ holds for an appropriate extension $\mathcal C\in\Cocc(G)$. Thus, (i) follows from (ii).
\end{proof}

\section{Existence of prescribed sections, part 1}\label{S:ex.prscr.sect.2}

In this section we continue to pursue the problem formulated at the beginning of Section~\ref{S:ex.prescr.sect.1}. First, we relate the existence of finite (or arbitrary) sections for extensions with values in a torus $G=\mathbb T^k$ to the difference of ranks of the weak homology groups $\rank(H_1^w(X))-\rank(H_1^w(\mathcal F))$; this is done in Theorem~\ref{T:ex.sec.tor.k}. Then, in Theorem~\ref{T:prscr.sct.solen}, we concentrate on the situation when $G=S_{{\bf p}}$ is the ${\bf p}$-adic solenoid with an arbitrary sequence ${\bf p}=(p_n)_{n\in\mathbb N}$. First we show that the flow $\mathcal F$ must possess a free cycle in order to possess an extension $\mathcal C\in\Cocc(S_{\bf p})$ with a non-trivial totally disconnected section $F(\mathcal C)$. Then, under the assumption that $\mathcal F$ does possess a free cycle, we show that $\Cocc(S_{\bf p})$ contains an extension with a prescribed non-trivial totally disconnected section $K\sbgp S_{\bf p}$ if and only if $K^{\perp}\cong\mathbb Z$. Finally, we show that for a flow $\mathcal F$ with a free cycle, the group $\Cocc(S_{\bf p})$ contains also minimal extensions.

\begin{theorem}\label{T:ex.sec.tor.k}
Let $\Flow$ be a minimal flow with $\Gamma\in\mathsf{LieGp}$ connected and with $X$ a compact connected manifold, and let $k\in\mathbb N$. Then the following conditions are equivalent:
\begin{enumerate}
\item[(1)] for every $K\sbgp\mathbb T^k$ there is $\mathcal C\in\Cocc(\mathbb T^k)$ with $F(\mathcal C)=K$,
\item[(2)] for every $K\sbgp\mathbb T^k$ finite there is $\mathcal C\in\Cocc(\mathbb T^k)$ with $F(\mathcal C)=K$,
\item[(3)] $\rank(H_1^w(X))-\rank(H_1^w(\mathcal F))\geq k$.
\end{enumerate}
\end{theorem}
\begin{proof}
Write $n=\rank(H_1^w(\mathcal F))$, $n+m=\rank(H_1^w(X))$ and let $d_1,\dots,d_n=d$ be the elementary divisors of $H_1^w(\mathcal F)$ in $H_1^w(X)$ written, as usual, in the increasing order.  The implication (1)$\Rightarrow$(2) is clear. To verify that (2) implies (3), fix an integer $l\geq2$ co-prime with $d$ and set $K=(\mathbb Z_l)^k\sbgp\mathbb T^k$. By virtue of (2), there is $\mathcal C\in\Cocc(\mathbb T^k)$ with $F(\mathcal C)=K$. Set $A=\mathbb Z_{d_1}\oplus\dots\oplus\mathbb Z_{d_n}\oplus\mathbb Z^m$. By Corollary~\ref{C:F.and.Betti}, $K$ is a quotient group of $A$. Denote the underlying quotient morphism $A\to K$ by $q$. Then, by our choice of $l$, $K=dK=dq(A)=q(dA)=q(d\mathbb Z^m)$. It follows that $m\geq k$, which verifies (3).

We show that (3) implies (2). So assume that $m\geq k$ and fix a finite subgroup $K$ of $\mathbb T^k$. The Pontryagin dual to the inclusion morphism $K\to\mathbb T^k$ is (equivalent to) a quotient morphism $\mathbb Z^k\to K$. Since $m\geq k$, it follows that $K$ is a quotient group of $\mathbb Z^m$ and hence also of $\mathbb Z_{d_1}\oplus\dots\oplus\mathbb Z_{d_n}\oplus\mathbb Z^m$. Thus, by Corollary~\ref{C:F.and.Betti}, there is an extension $\mathcal C\in\Cocc(\mathbb T^k)$ with $F(\mathcal C)=K$. This verifies (2).

We finish the proof by showing that (3) implies (1). So assume that (3) holds and fix a closed subgroup $K$ of $\mathbb T^k$. We distinguish two cases.

\emph{Case 1.} The group $K\sbgp\mathbb T^k$ is connected.

In this case $K$ is a torus with dimension at most $k$ and it is therefore a quotient group of $\mathbb T^k$. Let $r\colon\mathbb T^k\to K$ be a quotient morphism and assume that $\mathcal C\in\Cocc(\mathbb T^k)$ is minimal (that is, it satisfies $F(\mathcal C)=\mathbb T^k$). Then $r\mathcal C\in\Cocc(K)\subseteq\Cocc(\mathbb T^k)$ and $F(r\mathcal C)=rF(\mathcal C)=r(\mathbb T^k)=K$. Thus, in order to verify condition (1) in Case~1, it suffices to construct a minimal extension $\mathcal C\in\Cocc(\mathbb T^k)$.

Condition (3) yields an epimorphism $h\colon H_1^w(X)\to H_1^w(\mathbb T^k)$ with $h(H_1^w(\mathcal F))=0$. By virtue of Lemma~\ref{L:real.of.morph.pi1}, $h=\xi^{\sharp}$ for an appropriate continuous base point preserving map $\xi\colon X\to\mathbb T^k$. Let $p\colon\mathbb R^k\to\mathbb T^k$ be the usual covering morphism with kernel $\mathbb Z^k$. Since $\xi^{\sharp}\mathcal F_z^{\sharp}\pi_1(\Gamma)=0$ by our choice of $\xi$, statement (\ref{Eq:tor.Lie.mnf.L1}) from Lemma~\ref{L:tor.Lie.mnfld} yields $\co(\xi)^{\sharp}\pi_1(\Gamma\times X)=0$. Consequently, $\co(\xi)$ lifts across $p$ to an extension $\mathcal D\in\Cocc(\mathbb R^k)$. For every $l\in\mathbb N$ set $\mathcal C_l=p((1/l)\mathcal D)\in\Cocc(\mathbb T^k)$. Then $\kappa_l\mathcal C_l=\kappa_lp((1/l)\mathcal D)=p\mathcal D=\mathcal C$, where $\kappa_l$ stands for the $l$-endomorphism of $\mathbb T^k$.

We claim that the following statements hold:
\begin{enumerate}
\item[(a)] $F(\mathcal C_l)=(\mathbb Z_l)^k$ for every $l\in\mathbb N$; indeed, by virtue of Theorem~\ref{T:F.and.pi}(1),
\begin{equation*}
F(\mathcal C_l)=\mathcal E_{\ker(\kappa_l)}(\xi^{\sharp}\pi_1(X))=\mathcal E_{\ker(\kappa_l)}(\pi_1(\mathbb T^k))=\ker(\kappa_l)=(\mathbb Z_l)^k,
\end{equation*}
\item[(b)] $\mathcal C_l\stackrel{ucs}{\longrightarrow}1$ as $l\to\infty$; this follows from our definition of the extensions $\mathcal C_l$ ($l\in\mathbb N$) and from an obvious convergence $(1/l)\mathcal D\stackrel{ucs}{\longrightarrow}0$ as $l\to\infty$.
\end{enumerate}
Further, since the groups $F(\mathcal C_l)=(\mathbb Z_l)^k$ ($l\in\mathbb N$) become mutually group-disjoint if we restrict to the prime numbers $l$, we may apply Theorem~\ref{T:first.ineq.yes} to obtain an increasing sequence of prime numbers $(p_l)_{l\in\mathbb N}$ and an extension $\mathcal C\in\Cocc(\mathbb T^k)$ with
\begin{itemize}
\item $\mathcal C=\sum_{l=1}^{\infty}\mathcal C_{p_l}$ u.c.s. in $\Cocc(\mathbb T^k)$, and
\item $F(\mathcal C)\supseteq\overline{\sum_{l\in\mathbb N}F(\mathcal C_{p_l})}$.
\end{itemize}
Since $\mathbb T^k\supseteq F(\mathcal C)\supseteq\overline{\sum_{l\in\mathbb N}F(\mathcal C_{p_l})}=\overline{\sum_{l\in\mathbb N}(\mathbb Z_{p_l})^k}=\mathbb T^k$, it follows that $F(\mathcal C)=\mathbb T^k$ and so the extension $\mathcal C\in\Cocc(\mathbb T^k)$ is minimal. This verifies condition (1) in Case~1.

\emph{Case 2.} The group $K\sbgp\mathbb T^k$ is arbitrary.

The group $K$ can be expressed as a topological direct sum $K=K_1\oplus K_2$, where $K_1\sbgp\mathbb T^k$ is connected and $K_2\sbgp\mathbb T^k$ is finite. From Case~1 we know that $K_1=F(\mathcal C_1)$ for an appropriate $\mathcal C_1\in\Cocc(\mathbb T^k)$. Moreover, since the equivalence of (2) and (3) has already been established in an earlier part of the proof, we may apply (2) to obtain $\mathcal C_2\in\Cocc(\mathbb T^k)$ with $F(\mathcal C_2)=K_2$. Now, by virtue of Corollary~\ref{C:essen.disj.fin}, $F(\mathcal C_1+\mathcal C_2)=F(\mathcal C_1)+F(\mathcal C_2)=K_1\oplus K_2=K$. This verifies statement (1) also in Case~2.
\end{proof}

By applying Theorem~\ref{T:ex.sec.tor.k} to minimal equicontinuous flows $\mathcal F$, we will be able to enlarge ranges of dense Lie group morphisms in Corollary~\ref{T:min.qscob.sec}. Before doing so, we prove an auxiliary lemma.

\begin{lemma}\label{L:morph.ind.morph}
Let $G$ be a compact connected Lie group. Then the map
\begin{equation}
\varphi\colon\Hom(G,\mathbb T^1)\ni\xi\mapsto\xi^{\sharp}\in\Hom(\pi_1(G),\pi_1(\mathbb T^1))
\end{equation}
is an isomorphism of groups. In particular, given $\varrho\in\Hom(\pi_1(G),\pi_1(\mathbb T^1))$, there is $\xi\in\Hom(G,\mathbb T^1)$ with $\varrho=\xi^{\sharp}$.
\end{lemma}
\begin{remark}\label{R:morph.ind.morph}
We wish to add the following remarks.
\begin{itemize}
\item Recall from Subsection~\ref{Sub:ind.mrph.hmlg.gp} that there is an epimorphism
\begin{equation*}
C_e(G,\mathbb T^1)\ni f\mapsto f^{\sharp}\in\Hom(\pi_1(G),\pi_1(\mathbb T^1)),
\end{equation*}
which induces an isomorphism between the group $\pi^1(G)$ of the homotopy classes of continuous base point preserving maps $f\colon G\to\mathbb T^1$ and the group of morphisms $\varrho\colon\pi_1(G)\to\pi_1(\mathbb T^1)$. By combining this fact with the statement of Lem\-ma~\ref{L:morph.ind.morph}, we infer that every continuous base point preserving map $f\colon G\to\mathbb T^1$ is homotopic to a unique topological morphism $\xi\colon G\to\mathbb T^1$.

\item We were not able to find the statement of this lemma explicitly stated in the li\-te\-ra\-tu\-re. Therefore, we present its detailed proof based on results from the structure theory of compact (Lie) groups.
\end{itemize}
\end{remark}
\begin{proof}[Proof of Lemma~\ref{L:morph.ind.morph}]
A rough idea of the proof is to reduce the problem to the case when $G$ is a torus, in which case the statement of the lemma clearly holds (see, for instance, our discussion in Subsections~\ref{Sub:ind.mrph.hmlg.gp} and~\ref{Sub:chmtp.cpt.gps}). We divide the proof into three steps.

\emph{1st step.} We start by collecting some preliminary observations.

Let $G'$ be the commutator subgroup of $G$. By virtue of \cite[Theorem~6.18(i), p.~206]{HofMor}, $G'$ is a compact connected semi-simple Lie group and its fundamental group $\pi_1(G')$ is thus finite by the Weyl's theorem (see \cite[Corollary~4, p.~285]{Bou}). Further, by the second structure theorem for connected compact Lie groups (see \cite[Theorem~6.41(i), p.~221]{HofMor}), there is a (finite-dimensional) torus $T$ such that $G$ is the topological semi-direct product $T\ltimes G'$ of $T$ and $G'$; we denote by $p\colon G\to T$ the associated projection morphism. Under the usual identification $\pi_1(G)=\pi_1(T)\oplus\pi_1(G')$, the induced morphism $p^{\sharp}$ takes the form of the projection of $\pi_1(G)$ onto $\pi_1(T)$ and so it is an epimorphism with the kernel $\ker(p^{\sharp})=\pi_1(G')$. Finally, given $\xi\in\Hom(G,\mathbb T^1)$, we have $G'\subseteq\ker(\xi)$ by commutativity of $\mathbb T^1$ and so $\xi=\xi'p$ for a unique $\xi'\in\Hom(T,\mathbb T^1)$.

\emph{2nd step.} We show that $\varphi$ is a monomorphism.

Fix $\xi\in\Hom(G,\mathbb T^1)$ and assume that $\varphi(\xi)=\xi^{\sharp}=0$. Then $0=\xi^{\sharp}=(\xi'p)^{\sharp}=(\xi')^{\sharp}p^{\sharp}$ and hence $(\xi')^{\sharp}=0$ by surjectivity of $p^{\sharp}$. Since the domain of $\xi'$ is a torus, it follows that $\xi'=1$ and so $\xi=\xi'p=1$. Thus, $\varphi$ is indeed a monomorphism.

\emph{3rd step.} We finish the proof by showing that $\varphi$ is an epimorphism.

So let $\varrho\in\Hom(\pi_1(G),\pi_1(\mathbb T^1))$; we shall find $\xi\in\Hom(G,\mathbb T^1)$ with $\varrho=\xi^{\sharp}$. Since the group $\pi_1(G')$ is finite and the group $\pi_1(\mathbb T^1)$ is torsion-free, $\ker(\varrho)\supseteq\pi_1(G')=\ker(p^{\sharp})$. Consequently, there is $\vartheta\in\Hom(\pi_1(T),\pi_1(\mathbb T^1))$ with $\varrho=\vartheta p^{\sharp}$. Let $\zeta\in\Hom(T,\mathbb T^1)$ be such that $\zeta^{\sharp}=\vartheta$ and set $\xi=\zeta p$. Then $\xi\in\Hom(G,\mathbb T^1)$ and $\varrho=\vartheta p^{\sharp}=\zeta^{\sharp}p^{\sharp}=(\zeta p)^{\sharp}=\xi^{\sharp}$. Thus, $\varrho\in\im(\varphi)$, as was to be shown.
\end{proof}

\begin{corollary}\label{T:min.qscob.sec}
Let $\Gamma,G$ be connected Lie groups, $\Gamma$ non-compact, $G$ compact, and let $h\in\Hom(\Gamma,G)$ be a monomorphism with a dense image. Then for every $k\in\mathbb N$ there is $r\in\Hom(\Gamma,\mathbb T^k)$ such that the morphism $h\oplus r\colon\Gamma\to G\oplus\mathbb T^k$ has a dense image, that is, $\overline{\im(h\oplus r)}=G\oplus\mathbb T^k$.
\end{corollary}
\begin{remark}\label{R:min.qscob}
We wish to add the following observations.
\begin{itemize}
\item In connection with Remark~\ref{R:which.Gamma.fthfl.rep} notice the following immediate consequence of the corollary: if for a given connected group $\Gamma\in\mathsf{LieGp}$ there exist a compact group $G\in\mathsf{LieGp}$ and a monomorphism $h\colon\Gamma\to G$ with a dense image (in other words, if $\Gamma$ possesses a faithful finite-dimensional unitary representation) then $\Gamma$ possesses also a non-trivial one-dimensional unitary representation, that is, a non-trivial character $\Gamma\to\mathbb T^1$.
\item The statement of the corollary can be reformulated by saying that the minimal equicontinuous flow $\mathcal F\colon\Gamma\curvearrowright G$ induced by $h$ possesses a minimal pure quasi-coboundary $\mathcal Q_r\in\Cocc(\mathbb T^k)$ for every $k\in\mathbb N$.
\end{itemize}
\end{remark}
\begin{proof}[Proof of Corollary~\ref{T:min.qscob.sec}]
It is clearly sufficient to verify the statement of the corollary in case when $k=1$; the general case then follows by a repeated application of this particular one. Our proof will follow the line of that of Case~1 in the proof of Theorem~\ref{T:ex.sec.tor.k}, but we shall present it with all the details.

Let $\mathcal F\colon\Gamma\curvearrowright G$ be the minimal equicontinuous flow induced by $h$ (see (\ref{Eq:min.eq.flow.gen})). We choose the identity $z=e$ as the base point of $G$; this leads to $\mathcal F_z=h$. Since $h\colon\Gamma\to G$ is not an isomorphism, the flow $\mathcal F$ possesses a free cycle by virtue of Proposition~\ref{P:min.eq.flow.TFFC}; that is, $\rank(H_1^w(\mathcal F))<\rank(H_1^w(G))$. Consequently, there is an epimorphism $\varrho\colon H_1^w(G)\to H_1^w(\mathbb T^1)$ with $H_1^w(\mathcal F)\subseteq\ker(\varrho)$. By virtue of Lemma~\ref{L:morph.ind.morph}, $\varrho=\xi^{\sharp}$ for some $\xi\in\Hom(G,\mathbb T^1)$.

Now, under the usual identification $\pi_1(\mathbb T^1)=H_1^w(\mathbb T^1)$, Lemma~\ref{L:tor.Lie.mnfld} yields
\begin{equation*}
\co(\xi)^{\sharp}\pi_1(\Gamma\times G)=\xi^{\sharp}\mathcal F_z^{\sharp}\pi_1(\Gamma)=\xi^{\sharp}\mathcal F_z^{\sharp}H_1^w(\Gamma)=\xi^{\sharp}H_1^w(\mathcal F)=\varrho H_1^w(\mathcal F)=0.
\end{equation*}
Thus, $\co(\xi)$ lifts across the usual covering morphism $p\colon\mathbb R\to\mathbb T^1$ to an extension $\mathcal D\in\Cocc(\mathbb R)$. For $l\in\mathbb N$ set $\mathcal C_l=p((1/l)\mathcal D)$; clearly, $\mathcal C_l\in\Cocc$ is the lift of $\co(\xi)$ across the $l$-endomorphism $\kappa_l$ of $\mathbb T^1$. Also, since $(1/l)\mathcal D\stackrel{ucs}{\longrightarrow}0$ as $l\to\infty$, it follows that $\mathcal C_l\stackrel{ucs}{\longrightarrow}1$ as $l\to\infty$.

Further, by our choice of $\xi$ and by Theorem~\ref{T:F.and.pi}(1),
\begin{equation*}
F(\mathcal C_l)=\mathcal E_{\ker(\kappa_l)}\left(\xi^{\sharp}\pi_1(G)\right) =\mathcal E_{\ker(\kappa_l)}(\pi_1(\mathbb T^1))=\ker(\kappa_l)=\mathbb Z_l
\end{equation*}
for every $l\in\mathbb N$. Consequently, Theorem~\ref{T:first.ineq.yes} yields an increasing sequence of prime numbers $(p_l)_{l\in\mathbb N}$ and an extension $\mathcal C\in\Cocc$, such that
\begin{itemize}
\item $\mathcal C=\sum_{l=1}^{\infty}\mathcal C_{p_l}$ and
\item $F(\mathcal C)\supseteq\overline{\sum_{l\in\mathbb N}F(\mathcal C_{p_l})}$.
\end{itemize}
Moreover, since $\overline{\sum_{l\in\mathbb N}F(\mathcal C_{p_l})}=\overline{\sum_{l\in\mathbb N}\mathbb Z_{p_l}}=\mathbb T^1$, it follows that the extension $\mathcal C$ of $\mathcal F$ is minimal.

We claim that $\mathcal C$ is a pure quasi-coboundary over $\mathcal F$. First, given $\gamma\in\Gamma$ and $g\in G$, we have
\begin{equation*}
p\mathcal D(\gamma,g)=\co(\xi)(\gamma,g)=\xi(h(\gamma)g)\xi(g)^{-1}=\xi h(\gamma),
\end{equation*}
since $\xi\colon G\to\mathbb T^1$ is a morphism of groups. Consequently, since $G$ is connected and $\ker(p)=\mathbb Z$ is totally disconnected, it follows that $\mathcal D(\gamma,-)$ is constant on $G$ for every $\gamma\in\Gamma$. This means that $\mathcal D$ is a pure quasi-coboundary over $\mathcal F$ and hence so is $\mathcal C_l$ for every $l\in\mathbb N$. Since the pure quasi-coboundaries form a closed subgroup of $\Cocc$, the extension $\mathcal C=\sum_{l=1}^{\infty}\mathcal C_{p_l}$ is also a pure quasi-coboundary over $\mathcal F$, which verifies the claim.

Let $r\in\Hom(\Gamma,\mathbb T^1)$ be the morphism with $\mathcal C=\mathcal Q_r$; we show that $r$ satisfies the conclusion of the corollary. Indeed, by minimality of $\mathcal C$,
\begin{equation*}
G\oplus\mathbb T^1=\overline{\mathcal O}_{\mathcal F_{\mathcal C}}(e,1)=\overline{\{(h(\gamma),r(\gamma)) : \gamma\in\Gamma\}}=\overline{\im(h\oplus r)},
\end{equation*}
as was to be shown.
\end{proof}

\begin{theorem}\label{T:prscr.sct.solen}
Let $\Flow$ be a minimal flow with $\Gamma\in\mathsf{LieGp}$ connected and with $X$ a compact connected manifold, and let $d_1,\dots,d_n$ be the elementary divisors of $H_1^w(\mathcal F)$ in $H_1^w(X)$. Fix a solenoid $S_{\bf{p}}$ and a closed non-trivial proper subgroup $K$ of $S_{\bf{p}}$. Then the following statements hold:
\begin{enumerate}
\item[(A)] If the flow $\mathcal F$ does not possess a free cycle then there is no extension $\mathcal C\in\Cocc(S_{{\bf p}})$ with $F(\mathcal C)=K$.
\item[(B)] If the flow $\mathcal F$ does possess a free cycle then the following conditions are equivalent:
\begin{enumerate}
\item[(i)] $F(\mathcal C)=K$ for an appropriate extension $\mathcal C\in\Cocc(S_{\bf{p}})$,
\item[(ii)] $K^{\perp}$ is an infinite cyclic group, that is, $K^{\perp}\cong\mathbb Z$,
\item[(iii)] $K=\pr_m^{-1}(\mathbb Z_k)$ for an appropriate pair of integers $m\geq1$ and $k\geq2$.
\end{enumerate}
\item[(C)] If the flow $\mathcal F$ possesses a free cycle then there is $\mathcal C\in\Cocc(S_{{\bf p}})$ with $F(\mathcal C)=S_{{\bf p}}$, that is, $\Cocc(S_{{\bf p}})$ contains a minimal extension.
\end{enumerate}
\end{theorem}
\begin{remark}\label{R:prscr.sct.solen}
We wish to add the following remarks.
\begin{itemize}
\item In connection with statements (A) and (C) observe that if $\mathcal F$ does not possess a free cycle then there may or may not exist $\mathcal C\in\Cocc(S_{{\bf p}})$ with $F(\mathcal C)=S_{{\bf p}}$. Firstly, let $\mathcal F\colon\mathbb R\curvearrowright X$ be a minimal continuous flow on a simply connected manifold $X$, say, on $X=\mathbb S^3\times\mathbb S^3$. (The existence of a minimal continuous flow on $\mathbb S^3\times\mathbb S^3$ follows from \cite{FatHer}.) Then $H_1^w(X)=0$ (hence $\mathcal F$ does not possess a free cycle) and the group $\Cocc(G)$ contains a minimal extension for every connected group $G\in\mathsf{CAGp}$ with $\card(G^*)\leq\mathfrak{c}$. (This follows from our Theorem~\ref{T:grp.min.ext.am} in Section~\ref{Sub:gp.min.ext}.) Secondly, let $\mathcal F\colon\mathbb T^1\curvearrowright X$ be the natural action of $\mathbb T^1$ on $X=\mathbb T^1$ by rotations. Then $H_1^w(\mathcal F)=H_1^w(X)\cong\mathbb Z$ (hence $\mathcal F$ does not possess a free cycle) and for every non-trivial group $G\in\mathsf{CAGp}$, the group $\Cocc(G)$ contains no minimal extensions whatsoever. (Indeed, recall that, up to a homeomorphism, there is a unique non-degenerate space - namely $\mathbb T^1$, on which the group $\mathbb T^1$ acts in a minimal way.)
\item Condition (ii) from statement (B) is clearly equivalent to the isomorphism $S_{{\bf p}}/K\cong\mathbb T^1$. This occurs if and only if $p_K$ is equivalent to a non-trivial character of $S_{{\bf p}}$.
\item From statement (A) and from condition (ii) in statement (B) it follows that a closed non-trivial proper subgroup $K$ of $S_{{\bf p}}$ must be infinite in order to be the section $F(\mathcal C)$ of an extension $\mathcal C$ from $\Cocc(S_{{\bf p}})$. Consequently, finite non-trivial sections for extensions in $\Cocc(S_{{\bf p}})$ are not allowed. This observation follows also from Theorem~\ref{T:F.and.pi}(3), for every finite subgroup of a solenoid $S_{{\bf p}}$ intersects the identity arc-component $(S_{{\bf p}})_a$ of $S_{\bf p}$ only at the identity. Alternatively, one could also use Corollary~\ref{C:F.and.pi.fprt} and Remark~\ref{R:F.and.pi.fprt}.
\end{itemize}
\end{remark}
\begin{proof}[Proof of Theorem~\ref{T:prscr.sct.solen}]
We divide the proof into four steps.

\emph{1st step.} We show that the following conditions are equivalent for a non-trivial proper closed subgroup $K$ of $S_{{\bf p}}$:
\begin{enumerate}
\item[(a)] $K^{\perp}$ is an infinite cyclic group, that is, $K^{\perp}\cong\mathbb Z$,
\item[(b)] $K=\pr_m^{-1}(\mathbb Z_k)$ for an appropriate pair of integers $m,k\geq1$,
\item[(c)] $K=\pr_m^{-1}(\mathbb Z_k)$ for an appropriate pair of integers $m\geq1$ and $k\geq2$.
\end{enumerate}

Given integers $m,k\ge1$, the following statements are equivalent:
\begin{itemize}
\item $K=\pr_m^{-1}(\mathbb Z_k)$,
\item $K=\ker(\chi_{k,m})$,
\item $K^{\perp}=\langle\chi_{k,m}\rangle$,
\item $\lambda^{-1}(K^{\perp})=k\mathbb Z/p_0\dots p_{m-1}$,
\end{itemize}
where $\lambda$ is the isomorphism $\mathbb Z/{\bf p}\to(S_{\bf p})^*$ defined by (\ref{Eq:dual.gp.sole.iso}). Since the infinite cyclic subgroups of $\mathbb Z/{\bf p}$ are precisely the groups of the form $k\mathbb Z/p_0\dots p_{m-1}$ with $m,k\in\mathbb N$, the equivalence of (a) and (b) follows. Implication (c)$\Rightarrow$(b) is clear and the converse follows from the identity $\pr_m^{-1}(\mathbb Z_k)=\pr_{m+1}^{-1}(\mathbb Z_{kp_m})$, which holds for all $m,k\in\mathbb N$.

\emph{2nd step.} We verify statement (A).

Assume, contrary to (A), that $K=F(\mathcal C)$ for some extension $\mathcal C\in\Cocc(S_{\bf p})$. By Theorem~\ref{T:F.and.Betti}, there exist $k_1,\dots,k_n\in K\cap (S_{\bf p})_a$ with $K=\overline{\langle k_1,\dots,k_n\rangle}$ and such that $k_i^{d_i}=e$ for $i=1,\dots,n$. It follows that the group $K$ is finite and $K\subseteq(S_{\bf p})_a$. Consequently, $K\subseteq(S_{\bf p})_a\cap\tor(S_{{\bf p}})=e$, which contradicts the assumptions of the theorem.

\emph{3rd step.} We verify statement (B).

First, the equivalence of (ii) and (iii) has been established in the first part of the proof. Further, implication (i)$\Rightarrow$(ii) follows immediately from Remark~\ref{R:F.and.pi.fprt}. Finally, to show that (i) follows from (iii), let $K$ have the form $K=\pr_m^{-1}(\mathbb Z_k)$ with $m\geq1$ and $k\geq2$. Consider the morphism $h\in\mathcal L(S_{{\bf p}})$ defined by $h(t)=(\exp(i2\pi t/p_0\dots p_{n-1}))_{n\in\mathbb N}$. Then $h(p_0\dots p_{m-1}/k)$ is a topological generator of $K$, which lies on the identity arc-component $(S_{{\bf p}})_a$ of $S_{{\bf p}}$. Since the flow $\mathcal F$ possesses a free cycle by the assumptions, we may apply Theorem~\ref{T:F.and.Betti} to obtain $\mathcal C\in\Cocc(S_{{\bf p}})$ with $F(\mathcal C)=K$. This verifies statement (i) and finishes the proof of (B).

\emph{4th step.} We verify statement (C). The argument will be similar to the one used to verify Case~1 in the proof of Theorem~\ref{T:ex.sec.tor.k}.

Since the flow $\mathcal F$ possesses a free cycle by the assumptions, Lemma~\ref{L:real.of.morph.pi1} yields a continuous base point preserving map $\xi\colon X\to\mathbb T^1$ with $\xi^{\sharp}\mathcal F_z^{\sharp}\pi_1(\Gamma)=0$ and $\xi^{\sharp}\pi_1(X)=\pi_1(\mathbb T^1)$. Let $p\colon\mathbb R\to\mathbb T^1$ be the standard covering morphism with the kernel $\mathbb Z$. By virtue of (\ref{Eq:tor.Lie.mnf.L1}) from Lemma~\ref{L:tor.Lie.mnfld}, $\co(\xi)$ lifts across $p$ to an extension $\mathcal D\in\Cocc(\mathbb R)$. Consider the morphism $h\colon\mathbb R\to S_{{\bf p}}$ from the third step of the proof and set $\mathcal C_l=h((1/l)\mathcal D)\in\Cocc(S_{{\bf p}})$ for $l\in\mathbb N$. Since $(1/l)\mathcal D\stackrel{ucs}{\longrightarrow}0$ as $l\to\infty$, it follows that $\mathcal C_l\stackrel{ucs}{\longrightarrow}e$ as $l\to\infty$. Consequently, there exists an increasing sequence of prime numbers $(r_l)_{l\in\mathbb N}$ such that the series $\sum_{l=1}^{\infty}\mathcal C_{q_l}$ converges u.c.s. in $\Cocc(S_{{\bf p}})$ for every subsequence $(q_l)_{l\in\mathbb N}$ of $(r_l)_{l\in\mathbb N}$.

Given $l\in\mathbb N$, we show that $F(\mathcal C_l)=\pr_1^{-1}(\mathbb Z_l)$. First, from the identity $\pr_1 h=p$ it follows that
\begin{equation*}
(\kappa_l\text{pr}_1)\mathcal C_l=\kappa_l\text{pr}_1 h((1/l)\mathcal D)=\kappa_lp((1/l)\mathcal D)=p\mathcal D=\co(\xi).
\end{equation*}
Moreover, since $K_l:=\ker(\kappa_l\pr_1)=\pr_1^{-1}(\mathbb Z_l)$ is a (closed) totally disconnected subgroup of $S_{{\bf p}}$, we may apply (\ref{Eq:F.when.K.tot.disc}) from Theorem~\ref{T:F.and.pi} to obtain
\begin{equation*}
K_l\supseteq F(\mathcal C_l)=\overline{\mathcal E_{K_l}}(\xi^{\sharp}\pi_1(X))=\overline{\mathcal E_{K_l}}(\pi_1(\mathbb T^1))\supseteq\overline{\langle h(1/l)\rangle}=K_l;
\end{equation*}
where we have used the following facts:
\begin{itemize}
\item $\mathcal E_{K_l}(\pi_1(\mathbb T^1))\ni h(1/l)$; this follows from the fact that $h(t/l)$ ($0\leq t\leq 1$) is a path in $S_{\bf p}$ from $e$ to $h(1/l)$, which projects to the element $p(t)$ ($0\leq t\leq 1$) of $\pi_1(\mathbb T^1)$ via $\kappa_l\pr_1$,
\item $h(1/l)$ is a topological generator of $K_l$; this is clear (see also the third step of the proof).
\end{itemize}
Thus it follows that $F(\mathcal C_l)=K_l=\pr_1^{-1}(\mathbb Z_l)$, as was to be shown. As an immediate corollary of this fact we get $F(\pr_1\mathcal C_l)=\pr_1F(\mathcal C_l)=\mathbb Z_l$.

From our discussion so far the following facts emerge:
\begin{itemize}
\item $\pr_1\mathcal C_{r_l}\stackrel{ucs}{\longrightarrow}1$ as $l\to\infty$,
\item $F(\pr_1\mathcal C_{r_l})=\mathbb Z_{r_l}$ ($l\in\mathbb N$) are mutually group-disjoint.
\end{itemize}
Consequently, by Theorem~\ref{T:first.ineq.yes}, there is a subsequence $(q_l)_{l\in\mathbb N}$ of $(r_l)_{l\in\mathbb N}$ such that the series $\sum_{l=1}^{\infty}\pr_1\mathcal C_{q_l}$ converges u.c.s. in $\Cocc$ and
\begin{equation}\label{Eq:aux.eq.T1.pr1}
\mathbb T^1\supseteq F\left(\sum_{l=1}^{\infty}\text{pr}_1\mathcal C_{q_l}\right)\supseteq\overline{\sum_{l=1}^{\infty}F(\text{pr}_1\mathcal C_{q_l})}=\overline{\sum_{l=1}^{\infty}\mathbb Z_{q_l}}=\mathbb T^1.
\end{equation}
Now set $\mathcal C=\sum_{l=1}^{\infty}\mathcal C_{q_l}\in\Cocc(S_{{\bf p}})$. (Recall that the series converges u.c.s. in $\Cocc(S_{\bf p})$ by our choice of the sequence $(r_l)_{l\in\mathbb N}$ and so the extension $\mathcal C$ is well defined.) Moreover, by virtue of (\ref{Eq:aux.eq.T1.pr1}),
\begin{equation*}
\text{pr}_1F(\mathcal C)=F(\text{pr}_1\mathcal C)=F\left(\sum_{l=1}^{\infty}\text{pr}_1\mathcal C_{q_l}\right)=\mathbb T^1,
\end{equation*}
which shows that $F(\mathcal C)=S_{{\bf p}}$. This verifies statement (C).
\end{proof}

\section{Relation to the first cohomotopy group}\label{Sub:F.and.Cech}

Let $\Flow$ be a minimal flow. In this section we continue our study of the problem of determining the section $F(\mathcal C)$ of an extension $\mathcal C\in\Coc$ in terms of morphisms induced by transfer functions of certain quotients of $\mathcal C$ on one-dimensional algebraic-topological invariants. Now we concentrate on the situation when $X$ is an arbitrary compact connected space. When dealing with general, not necessarily locally connected continua, the first cohomotopy group $\pi^1$ is a more appropriate tool than the first homology group $H_1$. In Theorems~\ref{T:F.and.Cech} and~\ref{T:F.and.Cech.td} we explain how the functors $F$ and $\pi^1$ are related, and show that $\pi^1$ can be used to obtain the values of $F$ in some important situations.

Before formulating and proving our results, let us mention that in the whole present section we shall use notation introduced at the beginning of Section~\ref{S:lifts.transfer}.

\begin{theorem}\label{T:F.and.Cech}
Let $\Flow$ be a minimal flow, $G\in\mathsf{CAGp}$ be connected and $\mathcal C\in\Cocc(G)$. Fix $z\in X$. Assume that $K\sbgp G$ is finite with $p_K\mathcal C\in\Cob(G/K)$ (that is, with $F(\mathcal C)\subseteq K$), and denote by $\xi$ the transfer function of $p_K\mathcal C$ with $\xi(z)=e$. If the space $X$ is compact connected and the group $\Gamma$ has no non-trivial finite abelian quotient groups then the following statements hold under the usual identifications $\pi^1(G/K)\cong(G/K)^*\cong K^{\perp}$:
\begin{enumerate}
\item[(1)] the group $F(\mathcal C)$ is the smallest of all the subgroups $E$ of $K$, such that the morphism $\xi^{\flat}\colon (G/K)^*\to \pi^1(X)$ extends through $(q_E)^*\colon (G/K)^*\to (G/E)^*$ to a morphism $\varrho\colon (G/E)^*\to \pi^1(X)$,
\item[(2)] if $d$ stands for the largest torsion coefficient of $K$ and $\eta=|K|$ denotes the cardinality of $K$ then
\begin{equation}\label{Eq:F(C).cohmtp.d}
F(\mathcal C)^{\perp}=\frac{1}{d}\left(\xi^{\flat}\right)^{-1}\left(d \pi^1(X)\right)=\frac{1}{\eta}\left(\xi^{\flat}\right)^{-1}\left(\eta\pi^1(X)\right),
\end{equation}
see also \emph{(\ref{Eq:F.pi.precise})} below,
\item[(3)] the following statements are equivalent:
\begin{enumerate}
\item[($\iota$)] $F(\mathcal C)=K$,
\item[($\iota\iota$)] $dG^*\cap\left(\xi^{\flat}\right)^{-1}(d\pi^1(X))=dK^{\perp}$,
\item[($\iota\iota\iota$)] $\xi^{\flat}\left(dG^*\setminus d K^{\perp}\right)\cap d\pi^1(X)=\emptyset$.
\end{enumerate}
\end{enumerate}
\end{theorem}
\begin{remark}\label{R:F.and.Cech}
We wish to mention the following facts.
\begin{itemize}
\item The connectedness assumption of $G$ is not at all restrictive in the situation of the theorem, for under our assumptions on $\Gamma$ and $X$, every extension $\mathcal C\in\Cocc(G)$ with $G\in\mathsf{CAGp}$ takes its values in the identity component $G_0$ of $G$. Indeed, since the group $G/G_0$ is totally disconnected, its Pontryagin dual $(G/G_0)^*$ is a torsion group. Further, by Theorem~\ref{T:gimel.conn}, $\Cocc$ is a torsion-free group. Hence, by Theorem~\ref{T:structure.coc}, 
\begin{equation*}
\Cocc(G/G_0)\cong\Hom((G/G_0)^*,\Cocc)=0.
\end{equation*}
Thus, if $p_0\colon G\to G/G_0$ denotes the canonical quotient morphism then $p_0\mathcal C\in\Cocc(G/G_0)$ is trivial and hence $\mathcal C$ takes its values in $G_0$.
\item We would like to clarify formula (\ref{Eq:F(C).cohmtp.d}) from statement (2) of the theorem. By definition, $\xi^{\flat}$ is a morphism $\pi^1(G/K)\to\pi^1(X)$. Since the group $G/K\in\mathsf{CAGp}$ is connected, there are isomorphisms $\pi^1(G/K)\cong(G/K)^*\cong K^{\perp}$, see Subsection~\ref{Sub:chmtp.cpt.gps}. Thus, $(\xi^{\flat})^{-1}(d\pi^1(X))$ can be interpreted as a subgroup of $K^{\perp}\subseteq G^*$ and the expression $(1/d)(\xi^{\flat})^{-1}(d\pi^1(X))$ then stands for the group of all $\chi\in G^*$ with $\chi^d\in(\xi^{\flat})^{-1}(d\pi^1(X))$. If one prefers not to use the identification $(G/K)^*\cong K^{\perp}$ and to regard $\xi^{\flat}$ as a morphism $(G/K)^*\to\pi^1(X)$ then the formula (\ref{Eq:F(C).cohmtp.d}) takes the form
\begin{equation}\label{Eq:F.pi.precise}
F(\mathcal C)^{\perp}=\frac{1}{d}\,p_K^*\left(\xi^{\flat}\right)^{-1}\left(d \pi^1(X)\right)=\frac{1}{\eta}\,p_K^*\left(\xi^{\flat}\right)^{-1}\left(\eta\pi^1(X)\right).
\end{equation}
More generally, if the morphism $p_K$ is replaced by a general quotient morphism $q\colon G\to H$ with the kernel $K$ (and the map $\xi$ thus takes its values in $H$) then
\begin{equation}\label{Eq:F.pi.precise.gen}
F(\mathcal C)^{\perp}=\frac{1}{d}\,q^*\left(\xi^{\flat}\right)^{-1}\left(d \pi^1(X)\right)=\frac{1}{\eta}\,q^*\left(\xi^{\flat}\right)^{-1}\left(\eta\pi^1(X)\right),
\end{equation}
where $\xi^{\flat}$ is interpreted as a morphism $H^*\to\pi^1(X)$. This follows from (\ref{Eq:F.pi.precise}) and from the equivalence of $q$ and $p_K$.
\item If the group $\pi^1(X)$ is $d$-divisible (that is, if $d\pi^1(X)=\pi^1(X)$) then formula (\ref{Eq:F(C).cohmtp.d}) yields $F(\mathcal C)^{\perp}=(1/d)K^{\perp}=G^*$, where the last equality holds by definition of $d$. Thus, in this situation, $F(\mathcal C)=e$ and hence $\mathcal C\in\Cob(G)$.
\item Recall that under the assumptions of Theorem~\ref{T:F.and.pi}, $F(\mathcal C)\subseteq G_a$ as soon as $F(\mathcal C)$ is finite. Under the assumptions of Theorem~\ref{T:F.and.Cech} this is no longer the case, not even when the acting group $\Gamma$ of $\mathcal F$ is a (simply) connected Lie group. For an example of such a situation we refer to Section~\ref{S:ex.prescr.sect.3}, where finite subgroups of solenoids $S_{\bf p}$ appear as sections $F(\mathcal C)$ of appropriate extensions $\mathcal C\in\Cocc(S_{\bf p})$. (Recall that a finite subgroup of a solenoid $S_{\bf p}$ intersects the identity arc-component $(S_{\bf p})_a$ of $S_{\bf p}$ only at the identity $e$.)
\end{itemize}
\end{remark}
\begin{proof}[Proof of Theorem~\ref{T:F.and.Cech}]
We begin the proof by verifying statement (1). By virtue of statement (ii) from Proposition~\ref{P:char.F.as.smallest}, it suffices to show that for every subgroup $E$ of $K$, the following conditions are equivalent:
\begin{enumerate}
\item[($\alpha$)] there is $\sigma\in\Hom((G/E)^*,C_z(X,\mathbb T^1))$ with $\sigma (q_E)^*=\xi^*$,
\item[($\beta$)] there is $\varrho\in\Hom((G/E)^*,\pi^1(X))$ with $\varrho (q_E)^*=\xi^{\flat}$.
\end{enumerate}
First recall that there is an isomorphism $\varphi\colon C_z(X,\mathbb T^1)\to C_z(X,\mathbb R)\oplus\pi^1(X)$ with $\pr_2\varphi=\pi$, where $\pi$ stands for the canonical quotient morphism $C_z(X,\mathbb T^1)\to\pi^1(X)$. Set $\nu=\pr_1\varphi\xi^*$. Then $\nu\in\Hom((G/K)^*,C_z(X,\mathbb R))$ and, since $\pi\xi^*=\xi^{\flat}$, $\varphi\xi^*=\nu\oplus\xi^{\flat}$. Consequently, we may use identifications $C_z(X,\mathbb T^1)\cong C_z(X,\mathbb R)\oplus\pi^1(X)$ and $\xi^*\cong\nu\oplus\xi^{\flat}$.

Now let $E$ be a subgroup of $K$. We show that ($\alpha$) and ($\beta$) are equivalent. First, ($\beta$) follows from ($\alpha$) at once, for if $\sigma(q_E)^*=\xi^*$ then $(\pi\sigma)(q_E)^*=\pi(\sigma(q_E)^*)=\pi\xi^*=\xi^{\flat}$ and so we may take $\varrho=\pi\sigma$ to see that ($\beta$) holds. Conversely, assume that a morphism $\varrho\colon (G/E)^*\to\pi^1(X)$ satisfies $\varrho (q_E)^*=\xi^{\flat}$. Since the group $C_z(X,\mathbb R)$ is divisible, $\nu$ extends through $(q_E)^*$ to a morphism $\mu\colon (G/E)^*\to C_z(X,\mathbb R)$. That is, $\mu (q_E)^*=\nu$. Set $\sigma=\mu\oplus\varrho\colon (G/E)^*\to C_z(X,\mathbb T^1)$. Then $\sigma (q_E)^*=\mu (q_E)^*\oplus\varrho (q_E)^*=\nu\oplus\xi^{\flat}=\xi^*$ and this verifies statement ($\alpha$).

We verify formula (\ref{Eq:F(C).cohmtp.d}) from statement (2). In order to simplify notation, write $D=(1/d)\left(\xi^{\flat}\right)^{-1}(d\pi^1(X))$. Given a subgroup $E$ of $K$, we shall use the standard identification of the group $(G/E)^*$ with the annihilator $E^{\perp}$ of $E$ in $G^*$, the isomorphism between them being the Pontryagin dual $(p_E)^*\colon(G/E)^*\to E^{\perp}\subseteq G^*$ to the quotient morphism $p_E\colon G\to G/E$. Under this identification the morphism $(q_E)^*$ takes the form of the inclusion morphism $K^{\perp}\to E^{\perp}$. Further, via the annihilator mechanism in $G^*$, the subgroups $E$ of $K$ are in a one-to-one correspondence with the subgroups $A$ of $G^*$ containing $K^{\perp}$. Thus, by virtue of statement (1) of the theorem, $F(\mathcal C)^{\perp}$ is the largest (with respect to inclusion) of all the subgroups $A$ of $G^*$ containing $K^{\perp}$, such that $\xi^{\flat}$ extends to a morphism $\varrho\colon A\to\pi^1(X)$. Consequently, in order to prove the first equality in (\ref{Eq:F(C).cohmtp.d}), we need to show that the following statements hold:
\begin{enumerate}
\item[(i)] $K^{\perp}\subseteq D$,
\item[(ii)] for every $\chi\in G^*$, $\chi$ is an element of $D$ if and only if $\xi^{\flat}$ extends to a morphism $K^{\perp}+\langle\chi\rangle\to\pi^1(X)$.
\end{enumerate}
Before proving the two statements let us describe when a given character $\chi\in G^*$ belongs to $D$. Given $\chi\in G^*$, we have $\chi^d\in K^{\perp}$ by definition of $d$ and so $\chi^d=\Upsilon p_K$ for a unique $\Upsilon\in (G/K)^*$. Thus, under our identifications $\pi^1(G/K)\cong(G/K)^*\cong K^{\perp}$, $\chi$ belongs to $D$ if and only if $\xi^{\flat}(\Upsilon)\in d\pi^1(X)$.

Now, given $\chi\in K^{\perp}$, fix $\Upsilon\in(G/K)^*$ with $\chi=\Upsilon p_K$. Then $\chi^d=\Upsilon^dp_K$ and, since $\xi^{\flat}(\Upsilon^d)=d\xi^{\flat}(\Upsilon)\in d\pi^1(X)$, statement (i) follows. To verify statement (ii), fix $\chi\in G^*$ and let $\Upsilon\in(G/K)^*$ be such that $\chi^d=\Upsilon p_K$. We need to show that the following conditions are equivalent:
\begin{enumerate}
\item[(a)] $\xi^{\flat}$ extends to a morphism $\varrho\colon K^{\perp}+\langle\chi\rangle\to\pi^1(X)$,
\item[(b)] there is $f\in\pi^1(X)$ with $\xi^{\flat}(\Upsilon)=df$.
\end{enumerate}
Assume that condition (a) is satisfied and set $f=\varrho(\chi)$. Then $f\in\pi^1(X)$ and $df=d\varrho(\chi)=\varrho(\chi^d)=\xi^{\flat}(\chi^d)=\xi^{\flat}(\Upsilon)$, which verifies (b). Conversely, let $f\in\pi^1(X)$ be as in (b) and denote by $p$ the smallest positive integer with $\chi^p\in K^{\perp}$. Since $\chi^d\in K^{\perp}$, it follows that $p$ divides $d$ and so $d=pq$ for an appropriate positive integer $q$. We claim that $\xi^{\flat}(\chi^p)=pf$. Indeed, by our choice of $f$, $q\xi^{\flat}(\chi^p)=\xi^{\flat}(\chi^{pq})=\xi^{\flat}(\Upsilon)=df=q(pf)$ and hence, since $\pi^1(X)$ is a torsion-free group, $\xi^{\flat}(\chi^p)=pf$. Using now the additive notation in $G^*$, define a map $\varrho\colon K^{\perp}+\langle\chi\rangle\to\pi^1(X)$ by the rule $\varrho(\chi'+k\chi)=\xi^{\flat}(\chi')+kf$, where $\chi'\in K^{\perp}$ and $k\in\mathbb Z$. We show that $\varrho$ is well defined and it will immediately follow that $\varrho$ is in fact a morphism extending $\xi^{\flat}$. So let $\chi_1',\chi_2'\in K^{\perp}$ and $k,l\in\mathbb Z$, and assume that $\chi_1'+k\chi=\chi_2'+l\chi$. Then $(l-k)\chi=\chi_1'-\chi_2'\in K^{\perp}$ and hence, by our definition of $p$, $l-k=pr$ for an appropriate integer $r$. Thus,
\begin{equation*}
\begin{split}
\xi^{\flat}(\chi_1')+kf&=\xi^{\flat}(\chi_2'+r(p\chi))+kf=\xi^{\flat}(\chi_2')+r\xi^{\flat}(p\chi)+kf\\
&=\xi^{\flat}(\chi_2')+rpf+kf=\xi^{\flat}(\chi_2')+lf,
\end{split}
\end{equation*}
as was to be shown. This finishes the proof of the first equality from (\ref{Eq:F(C).cohmtp.d}).

We verify the second equality from (\ref{Eq:F(C).cohmtp.d}). Write $\eta=ld$ with $l\in\mathbb N$. Since the group $\pi^1(X)$ is torsion-free, the equivalence of the following conditions follows immediately for every $\chi\in G^*$:
\begin{itemize}
\item $\chi\in(1/\eta)\left(\xi^{\flat}\right)^{-1}(\eta\pi^1(X))$,
\item $\xi^{\flat}(\chi^{\eta})\in\eta\pi^1(X)$,
\item $l\xi^{\flat}(\chi^d)\in l(d\pi^1(X))$,
\item $\xi^{\flat}(\chi^d)\in d\pi^1(X)$,
\item $\chi\in(1/d)\left(\xi^{\flat}\right)^{-1}(d\pi^1(X))$.
\end{itemize}
This verifies the second equality from (\ref{Eq:F(C).cohmtp.d}).

We finish the proof of the theorem be verifying statement (3). By virtue of statement (2), condition ($\iota$) is equivalent with $K^{\perp}=(1/d)\left(\xi^{\flat}\right)^{-1}(d\pi^1(X))$. Moreover, the inclusion $K^{\perp}\subseteq(1/d)\left(\xi^{\flat}\right)^{-1}(d\pi^1(X))$ has been verified above (see condition (i)) and the inclusion $dG^*\cap\left(\xi^{\flat}\right)^{-1}(d\pi^1(X))\supseteq dK^{\perp}$ from ($\iota\iota$) is clear. Consequently, it remains to verify the equivalence of the following statements:
\begin{itemize}
\item $(1/d)\left(\xi^{\flat}\right)^{-1}(d\pi^1(X))\subseteq K^{\perp}$,
\item $dG^*\cap\left(\xi^{\flat}\right)^{-1}(d\pi^1(X))\subseteq d K^{\perp}$,
\item $\xi^{\flat}(dG^*\setminus dK^{\perp})\cap d\pi^1(X)=\emptyset$.
\end{itemize}
The equivalence of the above statements follows by elementary arguments, with the help of the fact that the group $G^*$ is torsion-free by the assumptions of the theorem.
\end{proof}

\begin{corollary}\label{C:when.F.conn}
Let $\Flow$ be a minimal flow on a compact connected space $X$. Consider the following statements:
\begin{enumerate}
\item[(i)] for every $\mathcal C\in\Coc$, the group $F(\mathcal C)$ is connected,
\item[(ii)] the group $\pi^1(X)$ is divisible.
\end{enumerate}
If $\Gamma$ has no non-trivial finite abelian quotient groups then \emph{(i)} follows from \emph{(ii)}. If $\Gamma\in\mathsf{CLAC}$ is simply connected then \emph{(ii)} follows from \emph{(i)}.
\end{corollary}
\begin{remark}\label{R:when.F.conn}
Recall from Remark~\ref{R:tor.gen.G} that the following conditions are equivalent for every non-trivial connected group $G\in\mathsf{CAGp}$:
\begin{itemize}
\item $G$ has no proper non-trivial closed connected subgroup,
\item $\rank(G^*)=1$.
\end{itemize}
These conditions are satisfied by the circle group $G=\mathbb T^1$ and by the solenoids $G=S_{{\bf p}}$. Consequently, if the acting group $\Gamma$ of $\mathcal F$ has no non-trivial finite abelian quotient groups, the space $X$ has a divisible first cohomotopy group $\pi^1(X)$ and $G\in\mathsf{CAGp}$ is a connected group with $\rank(G^*)=1$, then every extension $\mathcal C\in\Cocc(G)$ is either minimal or a coboundary.
\end{remark}
\begin{proof}[Proof of Corollary~\ref{C:when.F.conn}]
We show that (i) follows from (ii). To this end, assume that $\Gamma$ has no non-trivial finite abelian quotient groups and that the group $\pi^1(X)$ is divisible. Fix $G\in\mathsf{CAGp}$ and $\mathcal C\in\Cocc(G)$. We need to show that $F(\mathcal C)$ is connected or, equivalently, that $F(\mathcal C)^*$ is torsion-free. So let $\chi$ be a torsion element of $F(\mathcal C)^*$. Then $\chi F(\mathcal C)\subseteq\mathbb Z_d$ for some positive integer $d$. The morphism $\chi\colon F(\mathcal C)\to\mathbb T^1$ extends to a morphism of topological groups $G\to\mathbb T^1$, denote the latter one by the same symbol $\chi$. Since $F(\chi\mathcal C)=\chi F(\mathcal C)\subseteq\mathbb Z_d$, it follows that $F(\kappa_d(\chi\mathcal C))=\kappa_dF(\chi\mathcal C)=1$ and hence $\kappa_d(\chi\mathcal C)\in\Cob$, where $\kappa_d$ is the $d$-endomorphism of $\mathbb T^1$. Let $\xi\in C_z(X,\mathbb T^1)$ be the transfer function for $\kappa_d(\chi\mathcal C)$. Then by virtue of statement (2) from Theorem~\ref{T:F.and.Cech} and the second part of Remark~\ref{R:F.and.Cech}, we obtain
\begin{equation*}
\begin{split}
F(\chi\mathcal C)^{\perp}&=\frac{1}{d}(\kappa_d)^*\left(\xi^{\flat}\right)^{-1}(d\pi^1(X))=\frac{1}{d}(\kappa_d)^*\left(\xi^{\flat}\right)^{-1}(\pi^1(X))\\
&=\frac{1}{d}(\kappa_d)^*(\mathbb Z)=\frac{1}{d}(d\mathbb Z)=\mathbb Z
\end{split}
\end{equation*}
with the usual identifications $\pi^1(\mathbb T^1)\cong(\mathbb T^1)^*\cong\mathbb Z$. This shows that $\chi F(\mathcal C)=F(\chi\mathcal C)=1$ and so $\chi=1$ as an element of $F(\mathcal C)^*$. Thus, the group $F(\mathcal C)^*$ is torsion-free, as was to be shown.

We show that (ii) follows from (i). So assume that $\Gamma\in\mathsf{CLAC}$ is simply connected. Since $C_z(X,\mathbb T^1)\cong C_z(X,\mathbb R)\oplus\pi^1(X)$ and the group $C_z(X,\mathbb R)$ is divisible, the divisibility of $\pi^1(X)$ is equivalent to the divisibility of $C_z(X,\mathbb T^1)$. Fix $f\in C_z(X,\mathbb T^1)$ and $d\in\mathbb N$. By Theorem~\ref{T:lift.simply.con}, $\co(f)\in\Cob$ lifts across $\kappa_d\colon\mathbb T^1\to\mathbb T^1$ to an extension $\mathcal C\in\Cocc$ and $F(\mathcal C)\subseteq\ker(\kappa_d)=\mathbb Z_d$. Since the group $F(\mathcal C)$ is connected by the assumptions, $F(\mathcal C)=1$ and hence $\mathcal C\in\Cob$. Let $g\in C_z(X,\mathbb T^1)$ be the transfer function for $\mathcal C$. Then $f$ lifts across $\kappa_d$ to $g$, which means that $f=g^d$ holds in $C_z(X,\mathbb T^1)$. Since $f\in C_z(X,\mathbb T^1)$ and $d\in\mathbb N$ were arbitrary, this shows that the group $C_z(X,\mathbb T^1)$ is divisible and hence so is the group $\pi^1(X)$.
\end{proof}

\begin{corollary}\label{C:ex.fin.chmtp.intr}
Let $\Flow$ be a minimal flow with $\Gamma\in\mathsf{CLAC}$ simply connected and with $X$ compact. Given $G\in\mathsf{CAGp}$ connected and $K\sbgp G$ finite with the largest torsion coefficient $d$, the following conditions are equivalent:
\begin{enumerate}
\item[(a)] there is $\mathcal C\in\Cocc(G)$ with $F(\mathcal C)=K$,
\item[(b)] there is $\varrho\in\Hom(K^{\perp},\pi^1(X))$ with $dG^*\cap\varrho^{-1}(d\pi^1(X))=dK^{\perp}$.
\end{enumerate}
\end{corollary}
\begin{proof}
To see that (b) follows from (a), let $\xi$ be the base point preserving transfer function for $p_K\mathcal C$. Then by Theorem~\ref{T:F.and.Cech}(3), condition (b) holds with $\varrho=\xi^{\flat}$. To see that (a) follows from (b), fix $\varrho$ as in (b). By the isomorphism (\ref{Eq:maps.are.homs.hmtp}) from Subsection~\ref{Sub:chmtp.cpt.gps}, $\varrho=\xi^{\flat}$ for an appropriate continuous base point preserving map $\xi\colon X\to G/K$. By the second part of Remark~\ref{R:lift.simply.con}, Theorem~\ref{T:gimel.conn} and Theorem~\ref{T:coc.divisible}(f), $\co(\xi)$ lifts across $p_K$ to an extension $\mathcal C\in\Cocc(G)$. Using Theorem~\ref{T:F.and.Cech}(3) once more now yields $F(\mathcal C)=K$, which verifies condition (a).
\end{proof}

\begin{example}\label{E:mndrm.chmtp.unif}
Suppose that the assumptions of Theorem~\ref{T:F.and.pi} and Theorem~\ref{T:F.and.Cech} are satisfied. Then one can use both (\ref{Eq:F.when.K.tot.disc}) and (\ref{Eq:F(C).cohmtp.d}) to determine the section $F(\mathcal C)$ of a given extension $\mathcal C$. As a matter of fact, the two formulas follow from each other. To show this, we verify the following statement which does not involve any dynamics on $X$ at all:
\begin{enumerate}
\item[($*$)] Let $X$ be a compact connected manifold, $G\in\mathsf{CAGp}$ be connected and $K\sbgp G$ be finite with the largest torsion coefficient $d$. Given a continuous base point preserving map $\xi\colon X\to G/K$, the following identity holds:
\begin{equation}\label{Eq:mndrm.chmtp.unif}
\mathcal E_K\left(\xi^{\sharp}\pi_1(X)\right)^{\perp}=\frac{1}{d}\left(\xi^{\flat}\right)^{-1}(d\pi^1(X)),
\end{equation}
where the annihilator on the left is taken in $G^*$.
\end{enumerate}
\begin{proof}[Proof of statement ($*$)]
We shall prove ($*$) by verifying the equivalence of several pairs of auxiliary statements. This might not be the fastest way to arrive at (\ref{Eq:mndrm.chmtp.unif}) but we believe it is (one of) the most transparent.

Given $f\in\pi_1(G/K)$, let $\widetilde{f}$ be the lift of $f$ across $p_K$ starting at $e$. Fix $\chi\in G^*$. It follows immediately from the definition of $\mathcal E_K$ that the following statements are equivalent:
\begin{enumerate}
\item[($\alpha$)] $\chi\in\mathcal E_K(\xi^{\sharp}\pi_1(X))^{\perp}$,
\item[($\beta$)] $\chi\widetilde{f}$ is a loop for every $f\in\xi^{\sharp}\pi_1(X)$.
\end{enumerate}

Further, let $\kappa_d$ be the $d$-endomorphism of $\mathbb T^1$. By definition of $d$, there is a unique character $\Upsilon\in(G/K)^*$ of $G/K$ with $\kappa_d\chi=\Upsilon p_K$. We show that ($\beta$) is equivalent with
\begin{enumerate}
\item[($\gamma$)] $\Upsilon^{\sharp}\xi^{\sharp}\pi_1(X)\subseteq\kappa_d^{\sharp}\pi_1(\mathbb T^1)$.
\end{enumerate}
Assume, first, that ($\beta$) holds and fix $f\in\xi^{\sharp}\pi_1(X)$. Then $\Upsilon^{\sharp}(f)=\Upsilon f=\Upsilon p_K\widetilde{f}=\kappa_d\chi\widetilde{f}\in\kappa_d^{\sharp}\pi_1(\mathbb T^1)$, since $\chi\widetilde{f}$ is a loop. Conversely, let condition ($\gamma$) be fulfilled and fix $f\in\xi^{\sharp}\pi_1(X)$. By virtue of ($\gamma$), $\Upsilon f$ lifts across $\kappa_d$ to a loop $l$ in $\mathbb T^1$ based at $1$. That is, $\kappa_dl=\Upsilon f=\Upsilon p_K\widetilde{f}=\kappa_d\chi\widetilde{f}$ and hence $l=\chi\widetilde{f}$. This shows that $\chi\widetilde{f}$ is a loop and so statement ($\beta$) holds.

We continue the proof by showing that ($\gamma$) is equivalent with
\begin{enumerate}
\item[($\delta$)] $\Upsilon\xi\in d C_z(X,\mathbb T^1)$.
\end{enumerate}
Assume, first, that ($\delta$) holds and let $\eta\in C_z(X,\mathbb T^1)$ satisfy $\Upsilon\xi=\eta^d$. Then $(\Upsilon\xi)^{\sharp}=(\kappa_d\eta)^{\sharp}=\kappa_d^{\sharp}\eta^{\sharp}$ and so $\Upsilon^{\sharp}\xi^{\sharp}\pi_1(X)=(\Upsilon\xi)^{\sharp}\pi_1(X)=
\kappa_d^{\sharp}\eta^{\sharp}\pi_1(X)\subseteq \kappa_d^{\sharp}\pi_1(\mathbb T^1)$, which verifies ($\gamma$). Conversely, let condition ($\gamma$) be fulfilled and define a morphism $\varrho\colon\pi_1(X)\to\pi_1(\mathbb T^1)$ by the formula $\varrho=(1/d)(\Upsilon\xi)^{\sharp}$. By our assumptions on $X$, $\varrho=\eta^{\sharp}$ for an appropriate $\eta\in C_z(X,\mathbb T^1)$. Then $((\Upsilon\xi).(\eta^d)^{-1})^{\sharp}=0$ and so there is $\zeta\in C_z(X,\mathbb T^1)$ with $(\Upsilon\xi).(\eta^d)^{-1}=\zeta^d$. Consequently, $\Upsilon\xi=(\eta.\zeta)^d\in d C_z(X,\mathbb T^1)$ and statement ($\delta$) follows.

We finish the proof by showing that ($\delta$) is equivalent with
\begin{enumerate}
\item[($\epsilon$)] $\chi\in(1/d)(\xi^{\flat})^{-1}(d\pi^1(X))$;
\end{enumerate}
the equivalence of ($\alpha$) and ($\epsilon$) will then yield (\ref{Eq:mndrm.chmtp.unif}). Let $\pi$ be the canonical projection $C_z(X,\mathbb T^1)\to\pi^1(X)$. Since the complementary summand of $\pi^1(X)$ in $C_z(X,\mathbb T^1)$ is di\-vi\-si\-ble, the following statements are equivalent:
\begin{itemize}
\item $\Upsilon\xi\in dC_z(X,\mathbb T^1)$,
\item $\xi^{\flat}(\Upsilon)=\pi(\Upsilon\xi)\in d\pi^1(X)$,
\item $\chi\in(1/d)(\xi^{\flat})^{-1}(d\pi^1(X))$.
\end{itemize}
This verifies the equivalence of ($\delta$) and ($\epsilon$).
\end{proof}
\end{example}

In order to formulate our next result we fix some notation. Let $K\in\mathsf{CAGp}$ be a totally disconnected group and let $(H_j)_{j\in J}$ be a family of open subgroups of $K$ indexed by a directed set $J$, such that $H_j\subseteq H_i$ for $i\leq j$ and $\bigcap_{j\in J}H_j=e$. For $j\in J$ set $K_j=K/H_j$ and let $q_j\colon K\to K_j$ be the quotient morphism. Clearly, $K_j$ is a finite abelian group; denote by $d_j$ its largest torsion coefficient. Given $i\leq j$, there is a unique (quotient) morphism $q_{ij}\colon K_j\to K_i$ with $q_{ij}q_j=q_i$ and it follows that $d_i$ divides $d_j$. We call $(d_j)_{j\in J}$ a \emph{generating net of torsions\index{generating net of torsions} for $K$}. If there exist $d\in\mathbb N$ and $i\in J$ such that $d_j=d$ for all $j\geq i$ then we say that $(d_j)_{j\in J}$ \emph{terminates at $d$}\index{terminating net of torsions}. Finally, since $\bigcap_{j\in J}H_j=e$, the inverse limit of the inverse system $(K_i\stackrel{q_{ij}}{\longleftarrow}K_j)_{i\leq j}$ in $\mathsf{CAGp}$ is the group $K$ along with the morphisms $q_j$ as the limit projections.

\begin{remark}\label{R:gen.net.tor}
Notice the following facts.
\begin{itemize}
\item It is possible for two non-isomorphic groups to have a common generating net of torsions. For example, if $d\geq2$ is an integer then both $\mathbb Z_d$ and $\mathbb Z_d\oplus\mathbb Z_d$ have as their generating net of torsions the singleton consisting of $d$ alone.
\item A given totally disconnected group $K\in\mathsf{CAGp}$ has in general many distinct ge\-ne\-rating nets of torsions. Nevertheless, if $(d_j)_{j\in J}$ and $(d_l')_{l\in L}$ are two such nets then they are related by the following (symmetric) rules:
\begin{itemize}
\item[$\circ$] for every $j\in J$ there is $l\in L$ such that $d_j$ divides $d_l'$,
\item[$\circ$] for every $l\in L$ there is $j\in J$ such that $d_l'$ divides $d_j$.
\end{itemize}
In particular, if $K$ possesses a generating net of torsions terminating at $d$ then each of its generating nets of torsions terminates at $d$.
\item Assume that $(d_j)_{j\in J}$ is a net in $\mathbb N$ with the property that $d_i$ divides $d_j$ for $i\leq j$. Then there exists a totally disconnected group $K\in\mathsf{CAGp}$, of which $(d_j)_{j\in J}$ is a generating net of torsions. Indeed, for $i\leq j$ let $q_{ij}$ be the natural quotient morphism $\mathbb Z_{d_j}\to\mathbb Z_{d_i}$ (that is, under the usual identifications $\mathbb Z_{d_i}\subseteq\mathbb Z_{d_j}\subseteq\mathbb T^1$, $q_{ij}$ acts by the rule $q_{ij}(z)=z^{d_j/d_i}$). Then the inverse limit $K=\lim_{\leftarrow}(\mathbb Z_{d_i}\stackrel{q_{ij}}{\longleftarrow}\mathbb Z_{d_j})$ in $\mathsf{CAGp}$ is totally disconnected and it has $(d_j)_{j\in J}$ as its generating net of torsions.
\end{itemize}
\end{remark}

We would like to make now one more remark in order to make some of our later arguments more transparent.

\begin{remark}\label{R:Kprp.cnsms.Gstr}
Let $G\in\mathsf{CAGp}$ and $K\sbgp G$ be totally disconnected. Fix a generating net of torsions $(d_j)_{j\in J}$ for $K$. Then
\begin{itemize}
\item[($*$)] for every $\chi\in G^*$ there is $j\in J$ with $\chi^{d_j}\in K^{\perp}$.
\end{itemize}
To verify statement ($*$), fix $\chi\in G^*$. Recall that $V=\{\exp(i2\pi t) : -1/3<t<1/3\}$ is a neighbourhood of $1$ in $\mathbb T^1$ which contains no non-trivial subgroup of $\mathbb T^1$. Fix a neighbourhood $U$ of $e$ in $G$ with $\chi(U)\subseteq V$. By definition of $(d_j)_{j\in J}$, there exist $j\in J$ and an open subgroup $H_j$ of $K$ such that $H_j\subseteq U$ and $K/H_j$ has $d_j$ as its largest torsion coefficient. Since $\chi(H_j)\subseteq\chi(U)\subseteq V$ is a subgroup of $\mathbb T^1$, it follows that $\chi\in H_j^{\perp}$. Consequently, the restriction $\chi|_K$ of $\chi$ to $K$ equals $\chi|_K=\chi_jp_j$ for some $\chi_j\in (K/H_j)^*$, where $p_j$ denotes the quotient morphism $K\to K/H_j$. Then $\chi_j^{d_j}=1$ and hence $\chi^{d_j}\in K^{\perp}$, as was to be shown.
\end{remark}

\begin{theorem}\label{T:F.and.Cech.td}
Let $\Flow$ be a minimal flow and $G\in\mathsf{CAGp}$ be connected. Let $\mathcal C\in\Cocc(G)$ and fix $z\in X$. Assume that $K\sbgp G$ is totally disconnected with $p_K\mathcal C\in\Cob(G/K)$ and denote by $\xi$ the transfer function of $p_K\mathcal C$ with $\xi(z)=e$. If $X$ is compact connected, $\Gamma$ has no non-trivial finite abelian quotient groups and $(d_j)_{j\in J}$ is a generating net of torsions for $K$ then
\begin{equation}\label{Eq:F.and.Cech.td}
\begin{split}
F(\mathcal C)^{\perp}&=\bigcup_{j\in J}\frac{1}{d_j}\left(\xi^{\flat}\right)^{-1}\left(d_j\pi^1(X)\right)=\sum_{j\in J}\frac{1}{d_j}\left(\xi^{\flat}\right)^{-1}\left(d_j\pi^1(X)\right)\\
&=\lim_{\longrightarrow}\frac{1}{d_j}\left(\xi^{\flat}\right)^{-1}\left(d_j\pi^1(X)\right),
\end{split}
\end{equation}
where $\xi^{\flat}$ stands for the morphism $G^*\supseteq K^{\perp}\to\pi^1(X)$ under the usual identifications $\pi^1(G/K)\cong (G/K)^*\cong K^{\perp}$.
\end{theorem}
\begin{remark}\label{R:F.and.Cech.td}
We wish to add the following remarks.
\begin{itemize}
\item In the third step of the proof of Theorem~\ref{T:F.and.Cech.td} we shall observe that the inclusion $(1/d_i)(\xi^{\flat})^{-1}(d_i\pi^1(X))\subseteq(1/d_j)(\xi^{\flat})^{-1}(d_j\pi^1(X))$ holds for all $i,j\in J$ with $i\leq j$. Consequently, if the group $K$ possesses a generating net of torsions terminating at $d\in\mathbb N$ then the formula (\ref{Eq:F.and.Cech.td}) takes the form
\begin{equation*}
F(\mathcal C)^{\perp}=\frac{1}{d}\left(\xi^{\flat}\right)^{-1}(d\pi^1(X)),
\end{equation*}
just as in Theorem~\ref{T:F.and.Cech}(2), where the group $K$ was assumed to be finite with the largest torsion coefficient $d$.
\item It follows from the identity (\ref{Eq:F.and.Cech.td}) that the inclusion
\begin{equation*}
K^{\perp}\subseteq\bigcup_{j\in J}\frac{1}{d_j}\left(\xi^{\flat}\right)^{-1}(d_j\pi^1(X))
\end{equation*}
holds under the assumptions of Theorem~\ref{T:F.and.Cech.td}. Consequently, the equality $F(\mathcal C)=K$ holds if and only if $(1/d_j)(\xi^{\flat})^{-1}(d_j\pi^1(X))\subseteq K^{\perp}$ for every $j\in J$.
\end{itemize}
\end{remark}
\begin{proof}[Proof of Theorem~\ref{T:F.and.Cech.td}]
We divide the proof of the theorem into three steps. A rough idea of the proof is to use Theorem~\ref{T:F.and.Cech} together with the continuity of the functor $F\colon\Coc\to\mathsf{CAGp}$ from Theorem~\ref{T:functor.E.def}.

\emph{1st step.} We begin by fixing some notation and by collecting some useful observations.

Let $H_j,K_j,q_j$ ($j\in J$) and $q_{ij}$ ($i\leq j$) be as in the paragraph preceding Remark~\ref{R:gen.net.tor}. Further, for $j\in J$ let $p_j\colon G\to G/H_j$ be the quotient morphism and for $i\leq j$ let $p_{ij}$ be the uniqe quotient morphism $G/H_j\to G/H_i$ with $p_{ij}p_j=p_i$. Since $\bigcap_{j\in J}H_j=e$, the inverse system $(G/H_i\stackrel{p_{ij}}{\longleftarrow}G/H_j)_{i\leq j}$ in $\mathsf{CAGp}$ has as its inverse limit the group $G$ along with the morphisms $p_j$ as the limit projections. For $j\in J$ let $r_j$ be the (unique) morphism $G/H_j\to G/K$ with $r_jp_j=p_K$. Clearly, $r_j$ is a quotient morphism with a finite kernel $\ker(r_j)=p_j(K)\cong K/H_j=K_j$. Given $j\in J$, write $\mathcal C_j=p_j\mathcal C\in\Cocc(G/H_j)$. Then $r_j\mathcal C_j=r_jp_j\mathcal C=p_K\mathcal C=\co(\xi)$. Finally, since $G=\lim_{\leftarrow}(G/H_i\stackrel{p_{ij}}{\longleftarrow}G/H_j)$ in $\mathsf{CAGp}$, it follows from Lemma~\ref{L:inv.sys.ind} that $\mathcal C=\lim_{\leftarrow}(\mathcal C_i\stackrel{p_{ij}}{\longleftarrow}\mathcal C_j)$ in $\Coc$ and hence, by Theorem~\ref{T:functor.E.def}(5), $F(\mathcal C)=\lim_{\leftarrow}(F(\mathcal C_i)\stackrel{p_{ij}}{\longleftarrow}F(\mathcal C_j))$ in $\mathsf{CAGp}$.

\emph{2nd step.} The map $\xi\colon X\to G/K$ induces a morphism between the first cohomotopy groups $\xi^{\flat}\colon\pi^1(G/K)\to\pi^1(X)$, acting by the rule $\xi^{\flat}(f)=f\xi$. Under the identifications $\pi^1(G/K)\cong(G/K)^*\cong[K^{\perp},G^*]$, the morphism $\xi^{\flat}$ has the form
\begin{equation*}
\xi^{\flat}=\kappa\colon[K^{\perp},G^*]\ni\chi\mapsto\widetilde{\chi}
\xi\in\pi^1(X),
\end{equation*}
where, for a given $\chi\in[K^{\perp},G^*]$, $\widetilde{\chi}$ is the unique character of $G/K$ with $\widetilde{\chi}p_K=\chi$. Now fix $j\in J$. Then there is a topological isomorphism $G/K\cong(G/H_j)/(K/H_j)$, which induces isomorphisms of groups 
\begin{equation*}
\pi^1(G/K)\cong((G/H_j)/(K/H_j))^*\cong[(K/H_j)^{\perp},(G/H_j)^*].
\end{equation*}
By using these isomorphisms the morphism $\xi^{\flat}$ takes the form
\begin{equation*}
\xi^{\flat}=\lambda\colon[(K/H_j)^{\perp},(G/H_j)^*]\ni\Upsilon\mapsto\bar{\Upsilon}\xi\in\pi^1(X),
\end{equation*}
where, for a given $\Upsilon\in[(K/H_j)^{\perp},(G/H_j)^*]$, $\bar{\Upsilon}$ is the unique character of $G/K$ with $\bar{\Upsilon}r_j=\Upsilon$. Our aim in this step of the proof is to verify the formula
\begin{equation}\label{Eq:preim.Fj.annih}
\left[\left(p_j^{-1}F(\mathcal C_j)\right)^{\perp},G^*\right]=\left(\frac{1}{d_j}\kappa^{-1}(d_j\pi^1(X))\right)\cap\left[H_j^{\perp},G^*\right].
\end{equation}
First observe that the identity $\lambda r_j^*=\kappa p_K^*$ holds directly by definitions of $\lambda$ and $\kappa$.

Let $k\colon (K/H_j)^{\perp}\to K^{\perp}$ be the restriction of the morphism $p_j^*\colon(G/H_j)^*\to G^*$ dual to $p_j$ (notice that $p_j$ maps $K$ onto $K/H_j$ and hence $k$ is well defined). We claim that $\lambda=\kappa k$. To verify the claim, fix $\Upsilon\in(K/H_j)^{\perp}$. Then, by definition of the morphisms $\lambda$ and $\kappa$,
\begin{itemize}
\item $\lambda(\Upsilon)=\bar{\Upsilon}\xi$, where $\bar{\Upsilon}\in(G/K)^*$ satisfies $\bar{\Upsilon}r_j=\Upsilon$,
\item $\kappa k(\Upsilon)=\kappa(\Upsilon p_j)=\widetilde{\chi}\xi$, where $\widetilde{\chi}\in(G/K)^*$ satisfies $\widetilde{\chi}p_K=\Upsilon p_j$.
\end{itemize}
Consequently, it suffices to show that $\bar{\Upsilon}=\widetilde{\chi}$ or, equivalently, that $\bar{\Upsilon}p_K=\Upsilon p_j$. The latter identity follows by a simple computation
\begin{equation*}
\bar{\Upsilon}p_K=\bar{\Upsilon}(r_jp_j)=(\bar{\Upsilon}r_j)p_j=\Upsilon p_j.
\end{equation*}
This show that, indeed, $\lambda=\kappa k$ and hence, by definition of $k$,
\begin{equation}\label{Eq:dj.pi.lmbd.prmg}
\lambda^{-1}(d_j\pi^1(X))=(p_j^*)^{-1}\Big(\kappa^{-1}(d_j\pi^1(X))\Big)\cap\left[(K/H_j)^{\perp},(G/H_j)^*\right].
\end{equation}

Further, the morphism $r_j\colon G/H_j\to G/K$ has a finite kernel $\ker(r_j)\cong K/H_j=K_j$ with the largest torsion coefficient $d_j$. We may therefore use (\ref{Eq:F(C).cohmtp.d}) from Theorem~\ref{T:F.and.Cech} (see also the second part of Remark~\ref{R:F.and.Cech}) to obtain
\begin{equation}\label{Eq:preim.Fj.annih.lmbd}
\left[F(\mathcal C_j)^{\perp},(G/H_j)^*\right]=\frac{1}{d_j}r_j^*\left(\xi^{\flat}\right)^{-1}(d_j\pi^1(X))=\frac{1}{d_j}\lambda^{-1}(d_j\pi^1(X)),
\end{equation}
where the symbol $\xi^{\flat}=\kappa p_K^*=\lambda r_j^*$ now stands for the morphism $(G/K)^*\to\pi^1(X)$ under the usual identification $\pi^1(G/K)\cong(G/K)^*$. Now, by using (\ref{Eq:dj.pi.lmbd.prmg}) and (\ref{Eq:preim.Fj.annih.lmbd}), we obtain
\begin{equation*}
\begin{split}
\left[\left(p_j^{-1}F(\mathcal C_j)\right)^{\perp},G^*\right]&=p_j^*\left[F(\mathcal C_j)^{\perp},(G/H_j)^*\right]=p_j^*\left(\frac{1}{d_j}\lambda^{-1}(d_j\pi^1(X))\right)\\
&=p_j^*\left(\frac{1}{d_j}\Bigg((p_j^*)^{-1}\Big(\kappa^{-1}(d_j\pi^1(X))\Big)\cap\left[(K/H_j)^{\perp},(G/H_j)^*\right]\Bigg)\right)\\
&\stackrel{(\text{i})}{=}p_j^*\left(\frac{1}{d_j}\Bigg((p_j^*)^{-1}\Big(\kappa^{-1}(d_j\pi^1(X))\Big)\Bigg)\right)\\
&\stackrel{(\text{ii})}{=}p_j^*\left((p_j^*)^{-1}\left(\frac{1}{d_j}\kappa^{-1}(d_j\pi^1(X))\right)\right)\\
&=\left(\frac{1}{d_j}\kappa^{-1}(d_j\pi^1(X))\right)\cap\im(p_j^*)\\
&=\left(\frac{1}{d_j}\kappa^{-1}(d_j\pi^1(X))\right)\cap\left[H_j^{\perp},G^*\right],
\end{split}
\end{equation*}
where the equalities (i) and (ii) follow, respectively, from the following facts:
\begin{enumerate}
\item[(i)] $\Upsilon^{d_j}\in[(K/H_j)^{\perp},(G/H_j)^*]$ for every $\Upsilon\in(G/H_j)^*$; this follows from the definition of $d_j$ as the largest torsion coefficient of $K_j=K/H_j$,
\item[(ii)] $p_j^*(\Upsilon^{d_j})=(p_j^*(\Upsilon))^{d_j}$ for every $\Upsilon\in(G/H_j)^*$; this is clear.
\end{enumerate}
This verifies the desired identity (\ref{Eq:preim.Fj.annih}) and finishes the second step of the proof.

\emph{3rd step.} We finish the proof of the theorem. First, let us mention that from now on, all the annihilators are taken in $G^*$ and $\xi^{\flat}$ stands for the morphism $\kappa$ as described at the beginning of the second step of the proof.

By Lemma~\ref{P:cohom.inv} and Remark~\ref{R:cohom.inv}, the identity $F(\mathcal C)=\lim_{\leftarrow}(F(\mathcal C_i)\stackrel{p_{ij}}{\longleftarrow}F(\mathcal C_j))$ from the first step of the proof can be reformulated by writing $F(\mathcal C)=\bigcap_{j\in J}p_j^{-1}F(\mathcal C_j)$, where $p_i^{-1}F(\mathcal C_i)\supseteq p_j^{-1}F(\mathcal C_j)$ for all $i\leq j$. Consequently,
\begin{equation}\label{Eq:F.C.td.aux1}
F(\mathcal C)^{\perp}=\Bigg(\bigcap_{j\in J}p_j^{-1}F(\mathcal C_j)\Bigg)^{\perp}=\sum_{j\in J}\left(p_j^{-1}F(\mathcal C_j)\right)^{\perp}=\bigcup_{j\in J}\left(p_j^{-1}F(\mathcal C_j)\right)^{\perp}.
\end{equation}
By definitions of the objects $H_j$, $\mathcal C_j$ and $d_j$ ($j\in J$), the following inclusions hold for all $j\geq i$:
\begin{enumerate}
\item[($\alpha$)] $H_i^{\perp}\subseteq H_j^{\perp}$,
\item[($\beta$)] $(p_i^{-1}F(\mathcal C_i))^{\perp}\subseteq(p_j^{-1}F(\mathcal C_j))^{\perp}$,
\item[($\gamma$)] $(1/d_i)(\xi^{\flat})^{-1}(d_i\pi^1(X))\subseteq(1/d_j)(\xi^{\flat})^{-1}(d_j\pi^1(X))$.
\end{enumerate}
By using these inclusions along with (\ref{Eq:F.C.td.aux1}) and (\ref{Eq:preim.Fj.annih}), we obtain for every $i\in J$,
\begin{equation*}
\begin{split}
F(\mathcal C)^{\perp}\cap H_i^{\perp}&=\Bigg(\bigcup_{\substack{j\in J \\ j\geq i}}\left(p_j^{-1}F(\mathcal C_j)\right)^{\perp}\Bigg)\cap H_i^{\perp}\\
&=\Bigg(\bigcup_{\substack{j\in J \\ j\geq i}}\left(\left(\frac{1}{d_j}\left(\xi^{\flat}\right)^{-1}(d_j\pi^1(X))\right)\cap H_j^{\perp}\right)\Bigg)\cap H_i^{\perp}\\
&=\Bigg(\bigcup_{\substack{j\in J \\ j\geq i}}\frac{1}{d_j}\left(\xi^{\flat}\right)^{-1}(d_j\pi^1(X))\Bigg)\cap H_i^{\perp}\\
&=\Bigg(\bigcup_{j\in J}\frac{1}{d_j}\left(\xi^{\flat}\right)^{-1}(d_j\pi^1(X))\Bigg)\cap H_i^{\perp}.
\end{split}
\end{equation*}
Finally, since $\bigcup_{i\in J}H_i^{\perp}=G^*$, it follows that
\begin{equation*}
\begin{split}
F(\mathcal C)^{\perp}&=\bigcup_{i\in J}\left(F(\mathcal C)^{\perp}\cap H_i^{\perp}\right)=\bigcup_{i\in J}\Bigg(\Bigg(\bigcup_{j\in J}\frac{1}{d_j}\left(\xi^{\flat}\right)^{-1}(d_j\pi^1(X))\Bigg)\cap H_i^{\perp}\Bigg)\\
&=\bigcup_{j\in J}\frac{1}{d_j}\left(\xi^{\flat}\right)^{-1}(d_j\pi^1(X)),
\end{split}
\end{equation*}
which verifies the first equality from the theorem. The second and the third equality follow from the first one by virtue of ($\gamma$).
\end{proof}

\section{Existence of prescribed sections, part 2}\label{S:ex.prescr.sect.3}

In this section we continue our discussion from Sections~\ref{S:ex.prescr.sect.1} and~\ref{S:ex.prscr.sect.2}. Given a minimal flow $\Flow$, a group $G\in\mathsf{CAGp}$ and its subgroup $H\sbgp G$, we ask whether there is an extension $\mathcal C\in\Cocc(G)$ with $F(\mathcal C)=H$. Similarly to Section~\ref{S:ex.prscr.sect.2}, we concentrate mainly on the tori $G=\mathbb T^n$  and the solenoids $G=S_{\bf p}$. Our interest now is in extensions of flows with simply connected acting groups $\Gamma$ and arbitrary continua $X$ as phase spaces. For such flows $\mathcal F$, the results of Section~\ref{Sub:F.and.Cech} will serve us as the main tool.

We begin this section by considering first the situation when $X$ is a compact manifold, $G\in\mathsf{CAGp}$ is an arbitrary connected group and $H$ is totally disconnected; this is done in Theorem~\ref{T:gen.CLAC.mnfl.td}. (This case has already been covered by our Theorem~\ref{T:F.and.Betti} in Section~\ref{S:ex.prescr.sect.1}, but our results from Section~\ref{Sub:F.and.Cech} will enable us to see this case from a different point of view.) Then we turn to the case of a general continuum $X$ in the base and the fibre groups $G=\mathbb T^n$ and $G=S_{\bf p}$; in Theorems~\ref{T:ex.fin.chmtp.mod} and~\ref{T:sbgp.sol.inf.CLAC} we show how the existence of prescribed totally disconnected sections $F(\mathcal C)$ for extensions $\mathcal C\in\Cocc(G)$ is related to the algebraic properties of $\pi^1(X)$. At the end of this section we also give sufficient conditions for the existence of minimal extensions in the groups $\Cocc(G)$.

\begin{theorem}\label{T:gen.CLAC.mnfl.td}
Let $\Flow$ be a minimal flow with $\Gamma\in\mathsf{CLAC}$ simply connected and with $X$ a compact (connected) manifold, and let $G\in\mathsf{CAGp}$ be connected. Assume that $K$ is a closed non-trivial totally disconnected subgroup of $G$ and $(d_j)_{j\in J}$ is a generating net of torsions for $K$. Then the following conditions are equivalent:
\begin{enumerate}
\item[(1)] there exists $\mathcal C\in\Cocc(G)$ with $F(\mathcal C)=K$,
\item[(2)] the group $K^{\perp}$ can be expressed as a direct sum $K^{\perp}=F\oplus Q$, where $F$ is a free abelian group with $\rank(F)\leq\rank(\pi^1(X))$ and $Q$ is a $d_j$-pure subgroup of $G^*$ for every $j\in J$,
\item[(3)] the group $G/K$ contains a topological direct summand $T$, which is a torus with $\dim(T)\leq\rank(\pi^1(X))$ and whose pre-image $p_K^{-1}(T)$ is connected.
\end{enumerate}
\end{theorem}
\begin{remark}\label{R:gen.CLAC.mnfl.td}
We wish to add the following remarks.
\begin{itemize}
\item From our proof of Lemma~\ref{L:(co)hmtp.rel} it follows that $\pi^1(X)$ is isomorphic to $H_1^w(X)$. Since the latter group is finitely generated by our assumptions on $X$, it follows that the rank $\rank(\pi^1(X))=\rank(H_1^w(X))$ is finite.
\item In connection with condition (3) observe that for every closed subgroup $T$ of $G/K$, the group $p_K^{-1}(T)$ is connected if and only if $K\subseteq p_K^{-1}(T)_0$. Indeed, the implication from left to right follows from an obvious inclusion $K\subseteq p_K^{-1}(T)$ and the converse implication follows from the equality $p_K(p_K^{-1}(T)_0)=T$. We shall use this observation also in our proof of the theorem.
\item To clarify the nature of condition (3), let us give an example of a morphism $p=p_K\colon G\to G/K$ which the condition excludes. Let $T$ be an arbitrary torus, $r\colon S_{\bf p}\to S_{\bf q}$ be an epimorphism between two solenoids with a finite kernel $F$ and set $p=\id_T\oplus r\colon T\oplus S_{\bf p}\to T\oplus S_{\bf q}$. Then $p^{-1}(T)=T\oplus F$ is not connected and so condition (3) fails. (Also, in connection with the previous remark observe that the kernel $0\oplus F$ of $p$ is not contained in $p^{-1}(T)_0=T\oplus 0$.)
\end{itemize}
\end{remark}
\begin{proof}[Proof of Theorem~\ref{T:gen.CLAC.mnfl.td}]
We divide the proof into three steps.

\emph{1st step.} We show that (2) follows from (1).

So let $\mathcal C\in\Cocc(G)$ be an extension with $F(\mathcal C)=K$. Denote by $\xi$ the base point preserving transfer function for $p_K\mathcal C$. Let $Q$ be the kernel of the induced morphism $\xi^{\flat}\colon K^{\perp}\to\pi^1(X)$ with the usual identifications $\pi^1(G/K)\cong(G/K)^*\cong K^{\perp}$. By our assumptions on $X$,  $\pi^1(X)$ is a free abelian group with a finite rank, and so the group $Q$ is a direct summand in $K^{\perp}$. Let $F$ be a complementary summand to $Q$ in $K^{\perp}$. Since the restriction $\xi^{\flat}\colon F\to\pi^1(X)$ is a monomorphism, the group $F$ is free abelian with $\rank(F)\leq\rank(\pi^1(X))$. Now, in order to show that the group $Q$ is $d_j$-pure in $G^*$ for every $j\in J$, notice that, by virtue of Theorem~\ref{T:F.and.Cech.td},
\begin{equation*}
K^{\perp}=F(\mathcal C)^{\perp}=\bigcup_{j\in J}\frac{1}{d_j}\left(\xi^{\flat}\right)^{-1}(d_j\pi^1(X)),
\end{equation*}
which yields the inclusion
\begin{equation}\label{Eq:K.incl.dj.xi}
\frac{1}{d_j}\left(\xi^{\flat}\right)^{-1}(d_j\pi^1(X))\subseteq K^{\perp}
\end{equation}
for every $j\in J$. Fix $j\in J$ and $\chi\in G^*$ with $d_j\chi\in Q$; we show that $\chi\in Q$. By definition of $Q$, $\xi^{\flat}(d_j\chi)=0\in d_j\pi^1(X)$, which means that $\chi\in(1/d_j)(\xi^{\flat})^{-1}(d_j\pi^1(X))$. Thus, by virtue of (\ref{Eq:K.incl.dj.xi}), $\chi\in K^{\perp}$. Write $\chi=f+q$ with $f\in F$ and $q\in Q$. Then $d_jf+d_jq=d_j\chi\in Q$ and hence $d_jf=0$. Since the group $F$ is torsion-free, it follows that $f=0$. Consequently, $\chi=q\in Q$, as was to be shown.

\emph{2nd step.} We show that (1) follows from (2).

So let $F$ and $Q$ be as in (2). Fix a morphism $\varrho\colon K^{\perp}\to\pi^1(X)$ with the following properties:
\begin{itemize}
\item $\varrho$ restricted to $F$ is a monomorphism and $\varrho(F)$ is a direct summand in $\pi^1(X)$,
\item $Q$ equals the kernel $\ker(\varrho)$ of $\varrho$.
\end{itemize}
By virtue of the isomorphism (\ref{Eq:maps.are.homs.hmtp}) from Subsection~\ref{Sub:chmtp.cpt.gps}, there is a continuous base point preserving map $\xi\colon X\to G/K$ with $\xi^{\flat}=\varrho$. Further, by the second part of Remark~\ref{R:lift.simply.con}, Theorem~\ref{T:gimel.conn} and Theorem~\ref{T:coc.divisible}(f), $\co(\xi)$ lifts across $p_K$ to an extension $\mathcal C\in\Cocc(G)$. We show that $F(\mathcal C)=K$. By Theorem~\ref{T:F.and.Cech.td} and Remark~\ref{R:F.and.Cech.td}, it suffices to show that $(1/d_j)\varrho^{-1}(d_j\pi^1(X))\subseteq K^{\perp}$ for every $j\in J$. So fix $j\in J$ and let $\chi\in(1/d_j)\varrho^{-1}(d_j\pi^1(X))$. Then, in particular, $d_j\chi\in K^{\perp}$. Write $d_j\chi=f+q$ with $f\in F$ and $q\in Q$. Then $d_j\pi^1(X)\ni\varrho(d_j\chi)=\varrho(f)$ and hence, by definition of $\varrho$, $f=d_jf'$ for some $f'\in F$. Consequently, $d_j\chi=d_jf'+q$ and hence $d_j(\chi-f')=q\in Q$. Since the group $Q$ is $d_j$-pure in $G^*$ by the assumptions, it follows that $\chi-f'\in Q$. Thus, $\chi\in f'+Q\subseteq F\oplus Q=K^{\perp}$, which verifies the desired inclusion $(1/d_j)\varrho^{-1}(d_j\pi^1(X))\subseteq K^{\perp}$.

\emph{3rd step.} We show that (3) follows from (2).

So fix $F,Q$ as in (2). There exist $T,R\sbgp G/K$ such that $F=R^{\perp}$ and $Q=T^{\perp}$. We claim that there is a topological splitting $G/K=T\oplus R$. First, we have $e=F\cap Q=R^{\perp}\cap T^{\perp}=(R+T)^{\perp}$, which shows that $R+T=G/K$. Second, $(G/K)^*=F+Q=R^{\perp}+T^{\perp}=(R\cap T)^{\perp}$ and so $R\cap T=e$. It follows that the splitting $G/K=T\oplus R$ is algebraic and hence, by compactness of $T$ and $R$, it is also topological. Moreover, $T^*\cong(G/K)^*/T^{\perp}=(F\oplus Q)/Q\cong F$, which shows that $T$ is a torus with dimension $\dim(T)=\rank(F)\leq\rank(\pi^1(X))$.

Write $\widetilde{T}=p_K^{-1}(T)$. We need to show that the group $\widetilde{T}$ is connected. Given $j\in J$, the $d_j$-purity of $Q$ in $G^*$ translates into the following, mutually equivalent statements:
\begin{itemize}
\item for every $\chi\in G^*$, $\chi^{d_j}\in\widetilde{T}^{\perp}$ implies $\chi\in\widetilde{T}^{\perp}$,
\item $\kappa_{d_j}(\widetilde{T})^{\perp}=\widetilde{T}^{\perp}$,
\item $\kappa_{d_j}\colon\widetilde{T}\to\widetilde{T}$ is an epimorphism.
\end{itemize}
It follows that $\kappa_{d_j}\colon\widetilde{T}/\widetilde{T}_0\to\widetilde{T}/\widetilde{T}_0$ is an epimorphism for every $j\in J$.

Set $L=K\cap\widetilde{T}_0$. Then $L$ is clearly a closed subgroup of $K$. The inclusion morphism $K\to\widetilde{T}$ induces a morphism $\varphi\colon K/L\to\widetilde{T}/\widetilde{T}_0$, acting by the rule $\varphi(k+L)=k+\widetilde{T}_0$ for every $k\in K$. We claim that $\varphi$ is an isomorphism. First, $\varphi$ is a monomorphism by definition of $L$. Second, since the group $T$ is connected and $p_K\colon\widetilde{T}\to T$ is an epimorphism, the restriction $p_K\colon\widetilde{T}_0\to T$ is also an epimorphism and so every element of $\widetilde{T}/\widetilde{T}_0$ intersects $K$. This shows that $\varphi$ is an epimorphism. Thus, the groups $\widetilde{T}/\widetilde{T}_0$ and $K/L$ are (topologically) isomorphic. Since the morphisms $\kappa_{d_j}\colon\widetilde{T}/\widetilde{T}_0\to\widetilde{T}/\widetilde{T}_0$ are all epic by the previous paragraph, it follows that so are the morphisms $\kappa_{d_j}\colon K/L\to K/L$ for all $j\in J$.

Let $(c_i)_{i\in I}$ be a generating net of torsions for $K/L$. Given $i\in I$, the (topological) group $K/L$ factors onto $\mathbb Z_{c_i}$ and hence so does the group $K$. Consequently, for every $i\in I$ there is $j\in J$ such that $c_i$ divides $d_j$. On the other hand, since $\kappa_{d_j}\colon K/L\to K/L$ is an epimorphism for every $j\in J$, it follows that all $\kappa_{d_j}\colon \mathbb Z_{c_i}\to\mathbb Z_{c_i}$ are epimorphisms. Hence, $c_i$ is co-prime with $d_j$ for all $i\in I$ and $j\in J$. Summarizing, we have $c_i=1$ for every $i\in I$ and so $L=K$. This shows that $K\subseteq\widetilde{T}_0$ and hence, by part one of Remark~\ref{R:gen.CLAC.mnfl.td}, the group $\widetilde{T}=p_K^{-1}(T)$ is connected.

\emph{4th step.} We show that (2) follows from (3).

So fix $T,R$ as in (3) and write $Q=T^{\perp}$, $F=R^{\perp}$. As above, one verifies that $(G/K)^*=F\oplus Q$ and $F^*\cong T$, which shows that $F$ is a free abelian group with $\rank(F)=\dim(T)\leq\rank(\pi^1(X))$. It remains to show that the group $Q$ is $d_j$-pure in $G^*$ for every $j\in J$. That is, we must verify that the following condition holds for every $j\in J$:
\begin{enumerate}
\item[($*$)] given $\chi\in G^*$ and $\Upsilon\in Q$ with $\chi^{d_j}=\Upsilon p_K$, there is $\Lambda\in Q$ with $\chi=\Lambda p_K$.
\end{enumerate}
So let $\chi$ and $\Upsilon$ be as in ($*$). Since $Q=T^{\perp}$ and $\kappa_{d_j}(\widetilde{T})=\widetilde{T}$ by connectedness of $\widetilde{T}$, it follows that
\begin{equation*}
\chi(\widetilde{T})=\chi(\kappa_{d_j}(\widetilde{T}))=\chi^{d_j}(\widetilde{T})=\Upsilon p_K(\widetilde{T})=\Upsilon(T)=1.
\end{equation*}
This means, in particular, that $\chi\in K^{\perp}$, and hence $\chi=\Lambda p_K$ for an appropriate character $\Lambda\in(G/K)^*$ of $G/K$. Moreover, $\Lambda(T)=\Lambda p_K(\widetilde{T})=\chi(\widetilde{T})=1$, which shows that $\Lambda\in T^{\perp}=Q$. This verifies ($*$).
\end{proof}

\begin{theorem}\label{T:ex.fin.chmtp.mod}
Let $\Flow$ be a minimal flow with $\Gamma\in\mathsf{CLAC}$ simply connected and with $X$ compact. Fix positive integers $n$ and $d\geq2$. Then the following conditions are equivalent:
\begin{enumerate}
\item[(1)] for every finite subgroup $K$ of $\mathbb T^n$, whose largest torsion coefficient divides $d$, there exists $\mathcal C\in\Cocc(\mathbb T^n)$ with $F(\mathcal C)=K$,
\item[(2)] there are $f_1,\dots,f_n\in\pi^1(X)$ such that if $k_1,\dots,k_n\in\mathbb Z$ satisfy $\sum_{i=1}^n k_if_i\in d\pi^1(X)$ then $d$ divides $k_i$ for every $i=1,\dots,n$,
\item[(3)] $\pi^1(X)/d\pi^1(X)$ has rank at least $n$, when considered as a $\mathbb Z_d$-module.
\end{enumerate}
\end{theorem}
\begin{proof}
Condition (3) is only a reformulation of condition (2) and so it suffices to verify the equivalence of (1) and (2). We do this in two steps.

\emph{1st step.} We show that (2) follows from (1).

Set $K=(\mathbb Z_d)^n\sbgp\mathbb T^n$ and use (1) to find $\mathcal C\in\Cocc(\mathbb T^n)$ with $F(\mathcal C)=K$. Consider the morphism $\kappa_d\colon\mathbb T^n\to\mathbb T^n$. Since $\ker(\kappa_d)=K$, $\kappa_d\mathcal C\in\Cob(\mathbb T^n)$. Let $\xi=(\xi_1,\dots,\xi_n)\in C_z(X,\mathbb T^n)$ be the base point preserving transfer function for $\kappa_d\mathcal C$. Denote by $\pi$ the canonical projection morphism $C_z(X,\mathbb T^1)\to\pi^1(X)$ and set $f_i=\pi(\xi_i)\in\pi^1(X)$ for $i=1,\dots,n$. Then the induced morphism $\xi^{\flat}\colon\pi^1(\mathbb T^n)\to\pi^1(X)$ acts by the rule $\xi^{\flat}(k_1,\dots,k_n)=\sum_{i=1}^nk_if_i$ under the usual natural identifications $\pi^1(\mathbb T^n)\cong(\mathbb T^n)^*\cong\mathbb Z^n$.

Now, by virtue of Theorem~\ref{T:F.and.Cech} and part two of Remark~\ref{R:F.and.Cech},
\begin{equation*}
K^{\perp}=F(\mathcal C)^{\perp}=\frac{1}{d}\kappa_d^*\left(\xi^{\flat}\right)^{-1}(d\pi^1(X)).
\end{equation*}
Observe also that under the identification $(\mathbb T^n)^*\cong\mathbb Z^n$ used above, $\kappa_d^*\colon\mathbb Z^n\to\mathbb Z^n$ acts by the rule $\kappa_d^*(k_1,\dots,k_n)=(dk_1,\dots,dk_n)$ and $K^{\perp}=\im(\kappa_d^*)=d\mathbb Z^n$. Consequently, the following conditions are equivalent for every $(k_1,\dots,k_n)\in\mathbb Z^n$:
\begin{itemize}
\item $d$ divides $k_i$ for every $i=1,\dots,n$,
\item $(k_1,\dots,k_n)\in d\mathbb Z^n=K^{\perp}$,
\item $(k_1,\dots,k_n)\in(1/d)\kappa_d^*(\xi^{\flat})^{-1}(d\pi^1(X))$,
\item $\xi^{\flat}(k_1,\dots,k_n)\in d\pi^1(X)$,
\item $\sum_{i=1}^nk_if_i\in d\pi^1(X)$.
\end{itemize}
This shows that $f_1,\dots,f_n\in\pi^1(X)$ satisfy the condition from statement (2).

\emph{2nd step.} We show that (1) follows from (2).

Fix $f_1,\dots,f_n\in\pi^1(X)$ as in (2) and set $\xi=(f_1,\dots,f_n)\in C_z(X,\mathbb T^n)$. (Here we regard $\pi^1(X)$ as a subgroup of $C_z(X,\mathbb T^1)$ in a usual way.) By virtue of Theorem~\ref{T:lift.simply.con}, $\co(\xi)$ lifts across $\kappa_d\colon\mathbb T^n\to\mathbb T^n$ to an extension $\mathcal D\in\Cocc(\mathbb T^n)$. From Theorem~\ref{T:F.and.Cech} and part two of Remark~\ref{R:F.and.Cech} it follows that
\begin{equation*}
F(\mathcal D)^{\perp}=\frac{1}{d}\kappa_d^*\left(\xi^{\flat}\right)^{-1}(d\pi^1(X)).
\end{equation*}
Consequently, by our choice of the maps $f_1,\dots,f_n$, the following conditions are equivalent for every $(k_1,\dots,k_n)\in\mathbb Z^n$:
\begin{itemize}
\item $(k_1,\dots,k_n)\in d\mathbb Z^n$,
\item $\sum_{i=1}^nk_if_i\in d\pi^1(X)$,
\item $\xi^{\flat}(k_1,\dots,k_n)\in d\pi^1(X)$,
\item $(k_1,\dots,k_n)\in(1/d)\kappa_d^*(\xi^{\flat})^{-1}(d\pi^1(X))$,
\item $(k_1,\dots,k_n)\in F(\mathcal D)^{\perp}$.
\end{itemize}
This shows that $F(\mathcal D)^{\perp}=d\mathbb Z^n=((\mathbb Z_d)^n)^{\perp}$ and so $F(\mathcal D)=(\mathbb Z_d)^n$.

Now we show that condition (1) holds. So let $K\sbgp\mathbb T^n$ be a finite group with the largest torsion coefficient dividing $d$. There is an epimorphism $p\colon\mathbb T^n\to\mathbb T^n$ with $K=p((\mathbb Z_d)^n)$. Set $\mathcal C=p\mathcal D$. Then $\mathcal C\in\Cocc(\mathbb T^n)$ and $F(\mathcal C)=pF(\mathcal D)=p((\mathbb Z_d)^n)=K$, as was to be shown.
\end{proof}

\begin{proposition}\label{P:CLAC.min.lim.fin}
Let $\Flow$ be a minimal flow with $\Gamma\in\mathsf{LieGp}$ simply connected and with $X$ compact second countable. Assume that $\pi^1(X)\setminus p\pi^1(X)\neq\emptyset$ (that is, the group $\pi^1(X)$ is not $p$-divisible) for infinitely many prime numbers $p$. Then for every $m\in\mathbb N$ there is a sequence $(\mathcal C_n)_{n\in\mathbb N}$ in $\Coctd(\mathbb T^m)$, which converges u.c.s. to a minimal extension $\mathcal C\in\Cocc(\mathbb T^m)$.
\end{proposition}
\begin{remark}\label{R:CLAC.min.lim.fin}
Notice the following facts.
\begin{itemize}
\item In Theorem~\ref{P:tor.gen.G}(2) in Chapter~\ref{S:E.and.alg-top} we obtained a direct sum
\begin{equation*}
\Cocc(G)=\Coctd(G)\oplus\Coccn(G)
\end{equation*}
for every connected group $G\in\mathsf{CAGp}$ and in Remark~\ref{R:tor.gen.G} we claimed that the splitting is not topological in general. Proposition~\ref{P:CLAC.min.lim.fin} justifies this claim for many flows $\mathcal F$ in the base and for the fibre groups $G=\mathbb T^m$ ($m\in\mathbb N$).
\item In the proof of Proposition~\ref{P:CLAC.min.lim.fin} we shall use some of our earlier results from Section~\ref{S:sum.prop}. In order for these results to be applicable, we shall need a sequence $(\mathcal D_n)_{n\in\mathbb N}$ of extensions in $\Cocc$ with the groups $F(\mathcal D_n)$ ($n\in\mathbb N$) non-trivial and mutually group-disjoint. Such a sequence does not exist if only finitely many prime numbers $p$ satisfy $\pi^1(X)\setminus p\pi^1(X)\neq\emptyset$; this follows from Proposition~\ref{P:ex.prscr.F.T1} below.
\end{itemize}
\end{remark}
\begin{proof}[Proof of Proposition~\ref{P:CLAC.min.lim.fin}]
In the whole proof we shall regard $\pi^1(X)$ as a subgroup of $C_z(X,\mathbb T^1)$. Let $p$ be a prime number and $f\in\pi^1(X)\setminus p\pi^1(X)$. By virtue of Theorem~\ref{T:lift.simply.con}, $\co(f)\in\Cob$ lifts across $\kappa_{p^l}\colon\mathbb T^1\to\mathbb T^1$ to an extension $\mathcal D^{(p)}_l\in\Cocc$ for every $l\in\mathbb N$.

\emph{1st step.} We show that $F(\mathcal D^{(p)}_l)=\mathbb Z_{p^l}$ for all $p$ and $l$ as above.

Under the usual identifications $\pi^1(\mathbb T^1)\cong(\mathbb T^1)^*\cong\mathbb Z$, the morphism $f^{\flat}\colon\mathbb Z\to\pi^1(X)$ acts by the rule $f^{\flat}(k)=kf$ and the morphism $\kappa_{p^l}^*\colon\mathbb Z\to\mathbb Z$ acts by the rule $\kappa_{p^l}^*(k)=p^lk$. We claim that $(f^{\flat})^{-1}(p^l\pi^1(X))=p^l\mathbb Z$ or, equivalently, that $kf\in p^l\pi^1(X)$ if and only if $p^l$ divides $k$, where $k\in\mathbb Z$. The implication from right to left is clear. To verify the converse, let $k\in\mathbb Z$ satisfy $kf\in p^l\pi^1(X)$ and let $0\leq t\leq l$ be such that $\gcd(k,p^l)=p^t$. There exist integers $\alpha,\beta$ with $p^t=\alpha k+\beta p^l$. Then $p^tf=\alpha kf+\beta p^lf\in p^l\pi^1(X)$ and hence $f\in p^{l-t}\pi^1(X)$. Since $f\in\pi^1(X)\setminus p\pi^1(X)$ by the assumptions, it follows that $t=l$ and hence $p^l$ divides $k$. This verifies the claim.

Now, by virtue of Theorem~\ref{T:F.and.Cech}, part two of Remark~\ref{R:F.and.Cech} and our discussion above,
\begin{equation*}
F(\mathcal D^{(p)}_l)^{\perp}=\frac{1}{p^l}\kappa_{p^l}^*\left(f^{\flat}\right)^{-1}(p^l\pi^1(X))=\frac{1}{p^l}p^l(p^l\mathbb Z)=p^l\mathbb Z=\left(\mathbb Z_{p^l}\right)^{\perp},
\end{equation*}
which shows that, indeed, $F(\mathcal D^{(p)}_l)=\mathbb Z_{p^l}$.

\emph{2nd step.} Let $p$ be as above. We show that $\mathcal D^{(p)}_l\stackrel{ucs}{\longrightarrow}1$ as $l\to\infty$.

For every $0<s\leq 1/2$ we consider the set $W_s=\{\exp(i2\pi t) : -s<t<s\}\subseteq\mathbb T^1$. The sets $W_s$ form a local base at $1$ in $\mathbb T^1$. Moreover, for every $s\in(0,1/2]$ and every $k\in\mathbb N$, $\kappa_k^{-1}(W_s)_0=W_{s/k}$, where $\kappa_k^{-1}(W_s)_0$ stands for the identity component of the set $\kappa_k^{-1}(W_s)$ in $\mathbb T^1$.

Now fix $s\in(0,1/2]$ and a compact set $F\subseteq\Gamma$. By our assumptions on $\Gamma$ and by compactness of $X$, there is a connected identity neighbourhood $V$ in $\Gamma$ with $\co(f)(V\times X)\subseteq W_s$. Then for every $l\in\mathbb N$, $\kappa_{p^l}\mathcal D^{(p)}_l(V\times X)\subseteq W_s$ and hence, by connectedness of $V$ and $X$, $\mathcal D^{(p)}_l(V\times X)\subseteq\kappa_{p^l}^{-1}(W_s)_0=W_{s/p^l}$. Now let $l_0\in\mathbb N$ be large enough so that $V^{p^{l_0}}\supseteq F$. Then for every $l\geq l_0$, $V^{p^l}\supseteq F$ and so
\begin{equation*}
\mathcal D^{(p)}_l(F\times X)\subseteq\mathcal D^{(p)}_l(V^{p^l}\times X)\stackrel{(*)}{\subseteq}\left(\mathcal D^{(p)}_l(V\times X)\right)^{p^l}\subseteq\left(W_{s/p^l}\right)^{p^l}=W_s,
\end{equation*}
where the inclusion ($*$) follows from the cocycle identity for $\mathcal D^{(p)}_l$. This shows that $\mathcal D^{(p)}_l\stackrel{ucs}{\longrightarrow}1$ as $l\to\infty$.

\emph{3rd step.} We prove the proposition.

Fix $m\in\mathbb N$. From the first and the second step of the proof it follows that there is a sequence $(\mathcal D_n)_{n\in\mathbb N}$ of extensions in $\Cocc$, which satisfies the following conditions:
\begin{itemize}
\item the groups $F(\mathcal D_n)$ ($n\in\mathbb N$) are non-trivial, finite and mutually group-disjoint,
\item $F(\mathcal D_n)\to\mathbb T^1$ in $2^{\mathbb T^1}$ as $n\to\infty$,
\item $\mathcal D_n\stackrel{ucs}{\longrightarrow}1$ in $\Cocc$ as $n\to\infty$.
\end{itemize}
For $n\in\mathbb N$ set
\begin{equation*}
\mathcal E_n=(\mathcal D_{(n-1)m+1},\mathcal D_{(n-1)m+2},\dots,\mathcal D_{nm})\in\Cocc(\mathbb T^m).
\end{equation*}
By virtue of Proposition~\ref{P:essen.disj.fin},
\begin{equation*}
F(\mathcal E_n)=\prod_{i=1}^mF(\mathcal D_{(n-1)m+i})\sbgp\mathbb T^m
\end{equation*}
for every $n\in\mathbb N$, and so the sequence $(\mathcal E_n)_{n\in\mathbb N}$ satisfies the following conditions:
\begin{itemize}
\item the groups $F(\mathcal E_n)$ ($n\in\mathbb N$) are mutually group-disjoint,
\item $F(\mathcal E_n)\to\mathbb T^m$ in $2^{\mathbb T^m}$ as $n\to\infty$,
\item $\mathcal E_n\stackrel{ucs}{\longrightarrow}1$ in $\Cocc(\mathbb T^m)$ as $n\to\infty$.
\end{itemize}
By virtue of Theorem~\ref{T:first.ineq.yes}, we may assume that the series $\sum_{n=1}^{\infty}\mathcal E_n$ converges u.c.s. in $\Cocc(\mathbb T^m)$ and that its sum $\mathcal C\in\Cocc(\mathbb T^m)$ is minimal. Now set $\mathcal C_n=\sum_{i=1}^n\mathcal E_i$ for every $n\in\mathbb N$. Then the groups $F(\mathcal C_n)=\sum_{i=1}^nF(\mathcal E_i)$ are finite by Proposition~\ref{P:essen.disj.fin} and we have $\mathcal C_n\stackrel{ucs}{\longrightarrow}\mathcal C$ as $n\to\infty$. Since the extension $\mathcal C$ is minimal, this finishes the proof.
\end{proof}

\begin{proposition}\label{P:ex.prscr.F.T1}
Let $\Flow$ be a minimal flow with $\Gamma\in\mathsf{CLAC}$ simply connected and with $X$ compact. Let $d\geq 2$ be an integer and $d=p_1^{k_1}\dots p_n^{k_n}$ be the prime decomposition of $d$. Then the following conditions are equivalent:
\begin{enumerate}
\item[(1)] there is $\mathcal C\in\Cocc$ with $F(\mathcal C)=\mathbb Z_d$,
\item[(2)] $\pi^1(X)\setminus p_i\pi^1(X)\neq\emptyset$ for every $i=1,\dots,n$,
\item[(3)] $\bigcap_{i=1}^n\pi^1(X)\setminus p_i\pi^1(X)\neq\emptyset$.
\end{enumerate}
\end{proposition}
\begin{proof}
We shall follow the notation from the first step of our proof of Proposition~\ref{P:CLAC.min.lim.fin}. To show that (1) implies (3), fix $\mathcal C\in\Cocc$ with $F(\mathcal C)=\mathbb Z_d$ and denote by $\xi$ the base point preserving transfer function for $\kappa_d\mathcal C$. Set $f=\xi^{\flat}(1)$. Then $f\in\pi^1(X)$ and $\xi^{\flat}=f^{\flat}\colon\mathbb Z\ni k\mapsto kf\in\pi^1(X)$. By virtue of Theorem~\ref{T:F.and.Cech} and part two of Remark~\ref{R:F.and.Cech},
\begin{equation*}
d\mathbb Z=(\mathbb Z_d)^{\perp}=\frac{1}{d}\kappa_d^*(f^{\flat})^{-1}(d\pi^1(X))=\frac{1}{d}d(f^{\flat})^{-1}(d\pi^1(X))=(f^{\flat})^{-1}(d\pi^1(X)).
\end{equation*}
Thus, for every integer $k$, $kf\in d\pi^1(X)$ occurs if and only if $d$ divides $k$. This shows that $f\in\pi^1(X)\setminus l\pi^1(X)$ for every integer $l\geq2$ dividing $d$ and, in particular, $f\in\bigcap_{i=1}^n\pi^1(X)\setminus p_i\pi^1(X)$. Thus, (1) implies (3).

Implication (3)$\Rightarrow$(2) is clear and so it suffices to show that (2) implies (1). So assume that $\pi^1(X)\setminus p_i\pi^1(X)\neq\emptyset$ for every $i=1,\dots,n$ and write $d_i=p_i^{k_i}$. By the same argument as in the first step of our proof of Proposition~\ref{P:CLAC.min.lim.fin}, we find $\mathcal C_1,\dots,\mathcal C_n\in\Cocc$ with $F(\mathcal C_i)=\mathbb Z_{d_i}$ for every $i=1,\dots,n$. Set $\mathcal C=\sum_{i=1}^n\mathcal C_i\in\Cocc$. Then by virtue of Proposition~\ref{P:essen.disj.fin},
\begin{equation*}
F(\mathcal C)=F\left(\sum_{i=1}^n\mathcal C_i\right)=\sum_{i=1}^nF(\mathcal C_i)=\sum_{i=1}^n\mathbb Z_{d_i}=\mathbb Z_d.
\end{equation*}
This verifies condition (1).
\end{proof}

\begin{theorem}\label{T:sbgp.sol.inf.CLAC}
Let $\Flow$ be a minimal flow with $\Gamma\in\mathsf{CLAC}$ simply connected and with $X$ compact. Assume that $K$ is a non-trivial closed totally disconnected subgroup of a solenoid $S_{{\bf p}}$ and let $(c_j)_{j\in J}$ be a generating net of torsions for $K$. Then the following conditions are equivalent:
\begin{enumerate}
\item[(1)] there is an extension $\mathcal C\in\Cocc(S_{{\bf p}})$ with $F(\mathcal C)=K$,
\item[(2)] $\pi^1(X)$ contains a subgroup $A\cong K^{\perp}$, which is $c_j$-pure in $\pi^1(X)$ for every $j\in J$.
\end{enumerate}
\end{theorem}
\begin{proof}
We divide the proof into three steps.

\emph{1st step.} We collect some preliminary observations that will be of use to us in the proof of the theorem.

Write ${\bf p}=(p_n)_{n\in\mathbb N}$ and $p_0=1$. Given $n\in\mathbb N$, let $d_n$ be the positive integer with $\mathbb Z_{d_n}=\pr_n(K)$. Then $\kappa_{p_n}(\mathbb Z_{d_{n+1}})=\mathbb Z_{d_n}$ for every $n\in\mathbb N$. Consequently, for every $n\in\mathbb N$ there is $q_n\in\mathbb N$ co-prime with $d_n$ such that $d_{n+1}q_n=d_np_n$. Moreover, since $K$ is closed in $S_{\bf p}$,
\begin{equation*}
K=\lim_{\longleftarrow}\big(\mathbb Z_{d_n}\stackrel{\kappa_{p_n}}{\longleftarrow}\mathbb Z_{d_{n+1}}\big)_{n\in\mathbb N}=\left\{(z_n)_{n\in\mathbb N}\in S_{\bf p} : z_n\in\mathbb Z_{d_n}\,\text{for every}\,n\in\mathbb N\right\}.
\end{equation*}
Clearly, $(d_n)_{n\in\mathbb N}$ is a generating net of torsions for $K$. Under the isomorphism $\lambda\colon\mathbb Z/{\bf p}\to(S_{\bf p})^*$ defined by (\ref{Eq:dual.gp.sole.iso}) in Subsection~\ref{Sub:slnds}, the annihilator $K^{\perp}$ of $K$ in $(S_{\bf p})^*$ corresponds to
\begin{equation}\label{Eq:incr.un.annih}
\lambda^{-1}(K^{\perp})=\left\{\frac{ld_n}{p_0\dots p_{n-1}} : l\in\mathbb Z,n\in\mathbb N\right\}=\bigcup_{n\in\mathbb N}\frac{d_n}{p_0\dots p_{n-1}}\mathbb Z.
\end{equation}
We claim that
\begin{equation}\label{Eq:Kperp.bigcap.dn}
K^{\perp}\subseteq\bigcap_{n\in\mathbb N}d_n(S_{\bf p})^*.
\end{equation}
Indeed, the inclusion (\ref{Eq:Kperp.bigcap.dn}) is a consequence of (\ref{Eq:incr.un.annih}) and the following facts:
\begin{itemize}
\item the union on the right hand side of (\ref{Eq:incr.un.annih}) is increasing in $n\in\mathbb N$; this follows from the identity $d_{n+1}q_n=d_np_n$,
\item the intersection on the right hand side of (\ref{Eq:Kperp.bigcap.dn}) is decreasing in $n\in\mathbb N$; this follows from the fact that $d_n$ divides $d_{n+1}$ for every $n\in\mathbb N$,
\item for every $n\in\mathbb N$,
\begin{equation*}
\lambda\left(\frac{d_n}{p_0\dots p_{n-1}}\mathbb Z\right)\subseteq d_n(S_{\bf p})^*.
\end{equation*}
\end{itemize}

\emph{2nd step.} Let $\varrho\colon K^{\perp}\to\pi^1(X)$ be a monomorphism of groups. We show that the following conditions are equivalent for every $n\in\mathbb N$:
\begin{enumerate}
\item[(a)] $\im(\varrho)$ is $d_n$-pure in $\pi^1(X)$,
\item[(b)] $\im(\varrho)\cap d_n\pi^1(X)\subseteq d_n\im(\varrho)$,
\item[(c)] $\varrho^{-1}(d_n\pi^1(X))\subseteq d_nK^{\perp}$,
\item[(d)] $(1/d_n)\varrho^{-1}(d_n\pi^1(X))\subseteq K^{\perp}$.
\end{enumerate}

The equivalence of (a) and (b) follows at once from the definition of $d_n$-purity. Further, the implication (c)$\Rightarrow$(b) is clear and the converse holds since $\varrho$ is a mo\-no\-mor\-phism. Finally, the implication (c)$\Rightarrow$(d) holds since the group $(S_{\bf p})^*$ is torsion-free and the converse follows from the inclusion (\ref{Eq:Kperp.bigcap.dn}) from the first step of the proof.

\emph{3rd step.} We prove the equivalence of conditions (1) and (2) from the theorem.

First, let us mention the following facts.
\begin{itemize}
\item Condition (2) from the statement of the theorem is independent on the choice of a generating net of torsions for $K$; this follows from the second part of Remark~\ref{R:gen.net.tor} and from the fact that $d$-purity implies $d'$-purity for all $d,d'\in\mathbb N$ with $d'$ dividing $d$. As a result of this observation we shall restrict ourselves to the generating net of torsions $(d_n)_{n\in\mathbb N}$ for $K$ that has been constructed in the first step of the proof.
\item Let $\varrho\colon K^{\perp}\to\pi^1(X)$ be a morphism of groups. Then either $\varrho$ is a monomorphism or else it is trivial. This follows from the fact that the group $\pi^1(X)$ is torsion-free and the group $K^{\perp}$ has rank $1$.
\end{itemize}

To show that (2) follows from (1), fix $\mathcal C\in\Cocc(S_{\bf p})$ with $F(\mathcal C)=K$. Consider the quotient morphism $p_K\colon S_{\bf p}\to S_{\bf p}/K$ and denote by $\xi$ the base point preserving transfer function $\xi\colon X\to S_{\bf p}/K$ for $p_K\mathcal C$. By virtue of Theorem~\ref{T:F.and.Cech.td},
\begin{equation}\label{Eq:sbgp.sol.inf.CLAC.1}
K^{\perp}=F(\mathcal C)^{\perp}=\bigcup_{n\in\mathbb N}\frac{1}{d_n}\left(\xi^{\flat}\right)^{-1}(d_n\pi^1(X)).
\end{equation}
Since the group $K$ is non-trivial by the assumptions of the theorem, (\ref{Eq:sbgp.sol.inf.CLAC.1}) and Remark~\ref{R:Kprp.cnsms.Gstr} show that $\xi^{\flat}\colon K^{\perp}\to\pi^1(X)$ is non-trivial and hence it is a monomorphism. By applying (\ref{Eq:sbgp.sol.inf.CLAC.1}) once more and by using the equivalence of (a) and (d) from the second step of the proof, we conclude that the group $A:=\im(\xi^{\flat})\cong K^{\perp}$ is a $d_n$-pure subgroup of $\pi^1(X)$ for every $n\in\mathbb N$. This verifies condition (2).

To show that (1) follows from (2), assume that $\pi^1(X)$ contains a subgroup $A\cong K^{\perp}$, which is $d_n$-pure in $\pi^1(X)$ for every $n\in\mathbb N$. Fix a monomorphism $\varrho\colon K^{\perp}\to\pi^1(X)$ with $\im(\varrho)=A$. The equivalence of (a) and (d) from the second step of the proof yields the inclusion
\begin{equation}\label{Eq:Kprp.cntn.unn}
\frac{1}{d_n}\varrho^{-1}(d_n\pi^1(X))\subseteq K^{\perp}
\end{equation}
for every $n\in\mathbb N$. By virtue of (\ref{Eq:maps.are.homs.hmtp}) from Subsection~\ref{Sub:chmtp.cpt.gps}, $\varrho=\xi^{\flat}$ for an appropriate continuous base point preserving map $\xi\colon X\to S_{\bf p}/K$. Since the group $\Cocc$ is divisible by Remark~\ref{R:lift.simply.con} and torsion-free by connectedness of both $\Gamma$ and $X$, we may apply Theorem~\ref{T:coc.divisible}(f) to find $\mathcal C\in\Cocc(S_{\bf p})$ with $p_K\mathcal C=\co(\xi)$. Then, by Theorem~\ref{T:F.and.Cech.td} and (\ref{Eq:Kprp.cntn.unn}),
\begin{equation*}
K^{\perp}\subseteq F(\mathcal C)^{\perp}=\bigcup_{n\in\mathbb N}\frac{1}{d_n}\left(\xi^{\flat}\right)^{-1}(d_n\pi^1(X))\subseteq K^{\perp},
\end{equation*}
and hence $F(\mathcal C)=K$. This verifies condition (1).
\end{proof}

\begin{corollary}\label{T:sbgp.sol.CLAC}
Let $\Flow$ be a minimal flow with $\Gamma\in\mathsf{CLAC}$ simply connected and with $X$ compact. Given an integer $d\geq2$, a solenoid $S_{{\bf p}}$ and a subgroup $K\cong\mathbb Z_d$ of $S_{{\bf p}}$, the following conditions are equivalent:
\begin{enumerate}
\item[(1)] there is an extension $\mathcal C\in\Cocc(S_{{\bf p}})$ with $F(\mathcal C)=K$,
\item[(2)] $\pi^1(X)$ contains a $d$-pure subgroup isomorphic to $K^{\perp}$.
\end{enumerate}
\end{corollary}
\begin{remark}\label{R:sbgp.sol.CLAC}
We wish to mention the following facts.
\begin{itemize}
\item A finite group $K$ is isomorphic to a subgroup of a solenoid $S_{{\bf p}}$ if and only if it is isomorphic to a quotient group of $(S_{{\bf p}})^*$. Since $(S_{{\bf p}})^*$ has rank $1$, it follows that the number of the elementary divisors of $K$ is $1$, and hence $K$ is isomorphic to $\mathbb Z_d$ for some $d\in\mathbb N$. Thus the corollary covers all (non-trivial) finite subgroups of all solenoids.
\item Given a solenoid $S_{{\bf p}}$ and an integer $d\geq2$, the group $\mathbb Z_d$ may or may not be isomorphic to a subgroup of $S_{{\bf p}}$, depending on the sequence ${\bf p}=(p_n)_{n\in\mathbb N}$. As a matter of fact, the following conditions are equivalent:
\begin{enumerate}
\item[(a)] $S_{{\bf p}}$ contains a subgroup isomorphic to $\mathbb Z_d$,
\item[(b)] $\mathbb Z/{\bf p}$ factors onto $\mathbb Z_d$,
\item[(c)] there is $m\in\mathbb N$ such that $\gcd(p_n,d)=1$ for every $n\geq m$.
\end{enumerate}
First, the equivalence of (a) and (b) follows from the isomorphism $(S_{\bf p})^*\cong\mathbb Z/{\bf p}$. Further, the implication (a)$\Rightarrow$(c) follows from the first step of the proof of Theorem~\ref{T:sbgp.sol.inf.CLAC}. Finally, in proving the implication (c)$\Rightarrow$(b) one constructs the re\-qui\-red epimorphism $\mathbb Z/{\bf p}\to \mathbb Z_d$ by successively extending an arbitrary epimorphism $\mathbb Z/p_0\dots p_{m-1}\to\mathbb Z_d$ to the groups $\mathbb Z/p_0\dots p_n$ for $n=m,m+1,m+2,\dots$.
\end{itemize}
\end{remark}
\begin{proof}[Proof of Corollary~\ref{T:sbgp.sol.CLAC}]
Since the group $\mathbb Z_d$ has the singleton $d$ as its generating net of torsions, the corollary follows immediately from Theorem~\ref{T:sbgp.sol.inf.CLAC}.
\end{proof}

\begin{proposition}\label{P:min.ext.prscr.dim}
Let $\Flow$ be a minimal flow with $\Gamma\in\mathsf{CLAC}$ simply connected and with $X$ compact. Given $G\in\mathsf{CAGp}$ connected with $\dim(G)=n\in\mathbb N$, the following conditions are equivalent:
\begin{enumerate}
\item[(1)] $\Cocc(G)$ contains a minimal extension,
\item[(2)] $\Cocc(\mathbb T^n)$ contains a minimal extension.
\end{enumerate}
\end{proposition}
\begin{proof}
In the whole proof we shall use the notation $q\colon G\to\mathbb T(G)$ for the maximal toral quotient of $G$. Recall that the kernel $\ker(q)$ of $q$ is a totally disconnected subgroup of $G$. We divide the proof into two steps.

\emph{1st step.} We show that if $H\sbgp G$ satisfies $q(H)=\mathbb T(G)$ then $H=G$.

Since $\mathbb T(G)=\mathbb T^{\dim(G)}=\mathbb T^n$, we may identify $\mathbb T(G)^*$ with $\mathbb Z^n$. Let $j$ be the inclusion morphism $H\to G$ and set $r=qj\colon H\to\mathbb T(G)$. By our assumptions on $H$, $r$ is an epimorphism. The dual morphism $j^*$ to $j$ is an epimorphism $G^*\to H^*$ with kernel $H^{\perp}$. Therefore, in order to show that $H=G$, it suffices to verify that $j^*$ is a monomorphism.

To this end, fix $\chi\in G^*$ with $j^*(\chi)=0$. Let $e_1,\dots,e_n$ be the standard basis for $\mathbb Z^n$. Then $q^*(e_1),\dots,q^*(e_n)$ are independent in $G^*$ and hence, since $n=\dim(G)=\rank(G^*)$, there exist integers $k_1,\dots,k_n$ and $k\neq0$ such that $k\chi=\sum_{i=1}^nk_iq^*(e_i)$. By applying the morphism $j^*$ to the last equality, we obtain
\begin{equation*}
0=kj^*(\chi)=j^*(k\chi)=j^*\left(\sum_{i=1}^nk_iq^*(e_i)\right)=\sum_{i=1}^nk_ir^*(e_i)=r^*\left(\sum_{i=1}^nk_ie_i\right).
\end{equation*}
Since $r^*$ is a monomorphism, it follows that $k_1=\dots=k_n=0$ and hence $k\chi=0$. By connectedness of $G$, the group $G^*$ is torsion-free, and so $k\chi=0$ and $k\neq0$ imply $\chi=0$. Thus, $j^*$ is a monomorphism, as was to be shown.

\emph{2nd step.} We verify the equivalence of the conditions (1) and (2). First, we recall from the first step of the proof that $\mathbb T(G)=\mathbb T^n$.

The implication (1)$\Rightarrow$(2) is clear, for if $\mathcal C\in\Cocc(G)$ is minimal then so is $q\mathcal C\in\Cocc(\mathbb T(G))$. Conversely, assume that $\mathcal D\in\Cocc(\mathbb T(G))$ is minimal. Since the group $\Cocc$ is divisible and torsion-free by our assumptions on $\Gamma$ and $X$ (see Remark~\ref{R:lift.simply.con} and Theorem~\ref{T:gimel.conn}), and since the kernel $\ker(q)$ of $q$ is totally disconnected, we may apply Theorem~\ref{T:coc.divisible}(f) to find $\mathcal C\in\Cocc(G)$ with $q\mathcal C=\mathcal D$. Now $F(\mathcal C)$ is a closed subgroup of $G$ satisfying $qF(\mathcal C)=F(q\mathcal C)=F(\mathcal D)=\mathbb T(G)$. Thus, by the first step of the proof, $F(\mathcal C)=G$. This means that the extension $\mathcal C\in\Cocc(G)$ is minimal, which verifies condition (1).
\end{proof}

\begin{corollary}\label{C:min.ext.sol.CLAC}
Let $\Flow$ be a minimal flow with $\Gamma\in\mathsf{CLAC}$ simply connected and with $X$ compact. Then the following conditions are equivalent:
\begin{enumerate}
\item[(i)] $\Cocc(S_{\bf p})$ contains a minimal extension for every sequence ${\bf p}$,
\item[(ii)] $\Cocc(S_{\bf p})$ contains a minimal extension for some sequence ${\bf p}$,
\item[(iii)] $\Cocc$ contains a minimal extension.
\end{enumerate}
\end{corollary}
\begin{proof}
Since $S_{\bf p}\in\mathsf{CAGp}$ is connected and $\dim(S_{\bf p})=1$ for every sequence ${\bf p}$, the corollary follows directly from Proposition~\ref{P:min.ext.prscr.dim}.
\end{proof}

\chapter{First Cohomology Groups}\label{S:struct.res}

\section{Groups of minimal extensions}\label{Sub:gp.min.ext}

Let $\Flow$ be a minimal flow and $G\in\mathsf{CAGp}$ be a connected group. In this section we show that under mild assumptions on $\mathcal F$ and $G$, the groupoid $\Cocm(G)$ contains an isomorphic copy of an arbitrary torsion-free abelian group $A$ with $\card(A)\leq\mathfrak{c}$. We concentrate on the following three classes of flows $\mathcal F$:
\begin{itemize}
\item the group $\Gamma\in\mathsf{LCGp}\setminus\mathsf{CGp}$ is amenable, the space $X$ is compact and the flow $\mathcal F$ possesses a free point; see Theorem~\ref{T:grp.min.ext.am},
\item the group $\Gamma\in\mathsf{CLAC}$ is simply connected, the space $X$ is compact and $1\leq\rank(\pi^1(X))<\mathfrak{c}$; see Theorem~\ref{T:grp.min.ext.simp},
\item the group $\Gamma\in\mathsf{LieGp}$ is connected, the space $X$ is a compact connected manifold and the flow $\mathcal F$ possesses a free cycle; see Theorem~\ref{T:gp.min.free.ccl}.
\end{itemize}
One of the main tools used in this section is the following theorem.

\begin{theorem}\label{T:grp.min.ext}
Let $\Flow$ be a minimal flow and $\mathfrak k$ be an infinite cardinal number. Then the following conditions are equivalent.
\begin{enumerate}
\item[(1)] For every torsion-free group $A\in\mathsf{AbGp}$ with $\card(A)\leq\mathfrak k$ and every non-trivial connected group $H\in\mathsf{CAGp}$ with $\w(H)\leq\mathfrak k$, $\Cocm(H)$ contains an isomorphic copy of $A$.
\item[(2)] For every torsion-free group $A\in\mathsf{AbGp}$ with $\card(A)\leq\mathfrak k$, $\Cocm $ contains an isomorphic copy of $A$.
\item[(3)] For every non-trivial connected group $G\in\mathsf{CAGp}$ with $\w(G)\leq\mathfrak k$, $\Cocm(G)\setminus1\neq\emptyset$ (that is, $\Cocc(G)$ contains a minimal extension).
\end{enumerate}
\end{theorem}
\begin{remark}\label{R:grp.min.ext}
Notice the following facts.
\begin{itemize}
\item The assumption of torsion-freeness of the group $A$ is necessary in both (1) and (2). This follows from the fact that for every connected group $G\in\mathsf{CAGp}$, the groupoid $\Cocm(G)$ contains no non-trivial torsion elements. Indeed, if $k\in\mathbb N$ then $\kappa_k\colon G\to G$ is an epimorphism by connectedness of $G$ and hence $F(\mathcal C^k)=F(\kappa_k\mathcal C)=\kappa_kF(\mathcal C)=\kappa_k(G)=G$ for every $\mathcal C\in\Cocm(G)\setminus1$.
\item If $A\in\mathsf{AbGp}$ is torsion-free with $\card(A)\leq\mathfrak{k}$ then the divisible hull $\mathbb Q\otimes A$ of $A$ satisfies $\QLSdim(\mathbb Q\otimes A)=\rank(A)\leq\card(A)\leq\mathfrak{k}$. Hence $\mathbb Q\otimes A$ can be viewed as a subgroup of $\mathbb Q^{(\mathfrak{k})}$, the cardinality of the latter group being $\card(\mathbb Q^{(\mathfrak{k})})=\mathfrak{k}$. Consequently, a groupoid $B$ contains an isomorphic copy of every torsion-free group $A\in\mathsf{AbGp}$ with $\card(A)\leq\mathfrak{k}$ if and only if $B$ contains an isomorphic copy of $\mathbb Q^{(\mathfrak{k})}$. (This applies to the groupoids $B=\Cocm(H)$ and $B=\Cocm$ in conditions (1) and (2), respectively.) In case of $\mathfrak{k}=\mathfrak{c}$ we have $\mathbb Q^{(\mathfrak{k})}=\mathbb Q^{(\mathfrak{c})}\cong\mathbb R$.
\end{itemize}
\end{remark}
\begin{proof}[Proof of Theorem~\ref{T:grp.min.ext}]
We divide the proof into three steps. First, let us fix some notation. Let $A,B$ be abelian groups and let $E$ be a sub-groupoid of $B$. By a (mono)morphism $A\to E$ we mean a (mono)morphism of groups $A\to B$, whose image is contained in $E$. This terminology will be applied twice in the proof. First in the case of $B=\Cocc$ and $E=\Cocm$, and second in the case of $B=\Hom(H^*,\Cocc)$ and $E=\Hom(H^*,\Cocm)$ with $H\in\mathsf{CAGp}$.

\emph{1st step.} We recollect some useful properties of the tensor product in $\mathsf{AbGp}$.

Let $A,C$ be non-trivial torsion-free abelian groups with cardinalities at most $\mathfrak{k}$. Then the following statements hold.
\begin{enumerate}
\item[($\iota$)] If $a\in A\setminus0$ and $c\in C\setminus0$ then $a\otimes c\neq0$; in particular, $A\otimes C$ is a non-trivial abelian group. Indeed, since the groups $A,C$ are torsion-free and the group $\mathbb Q$ is divisible, there are morphisms $\varphi\colon A\to\mathbb Q$ and $\psi\colon C\to\mathbb Q$ with $\varphi(a)=\psi(c)=1$. Let $e\colon A\times C\to A\otimes C$ be the universal bilinear map. Since the map $b\colon A\times C\ni(u,v)\mapsto\varphi(u)\psi(v)\in\mathbb Q$ is also bilinear, there is a morphism $h\colon A\otimes C\to \mathbb Q$ with $he=b$. Then $h(a\otimes c)=b(a,c)=\varphi(a)\psi(c)=1\neq0$, and hence $a\otimes c\neq0$.
\item[($\iota\iota$)] The group $A\otimes C$ is torsion-free. This follows from the following three facts. First, both $A$ and $C$ are direct limits of (non-trivial) finitely generated free abelian groups. Second, the (non-trivial) finitely generated free abelian groups are closed with respect to the tensor product. Finally, the tensor product commutes with direct limits.
\item[($\iota\iota\iota$)] The cardinality of $A\otimes C$ is at most $\mathfrak{k}$. Indeed, the group $A\otimes C$ is generated by the monomial tensors $a\otimes c$ ($a\in A$, $c\in C$) and these form a set with cardinality at most $\mathfrak{k}$.
\end{enumerate}

\emph{2nd step.} We show that (1)$\Rightarrow$(2)$\Rightarrow$(3).

Implication (1)$\Rightarrow$(2) is immediate, since $\w(\mathbb T^1)=\aleph_0\leq\mathfrak k$. To verify implication (2)$\Rightarrow$(3), let $G\in\mathsf{CAGp}$ be non-trivial connected with $\w(G)\leq\mathfrak k$. Then $G^*$ is a non-trivial torsion-free abelian group with $\card(G^*)=\w(G)\leq\mathfrak k$ and so condition (2) yields a monomorphism $\varphi\colon G^*\to\Cocm $. By virtue of Corollary~\ref{C:min.contr.char}, there is $\mathcal C\in\Cocm(G)\setminus1$ with $\varphi=\Phi_G(\mathcal C)=\mathcal C^*$. Thus, $\Cocm(G)\setminus1\neq\emptyset$ and condition (3) hods.

\emph{3rd step.} We finish the proof by  verifying implication (3)$\Rightarrow$(1).

Fix $A$ and $H$ as in (1). We may clearly assume that $A$ is non-trivial. By the first step of the proof, the group $A\otimes H^*$ is non-trivial torsion-free abelian with cardinality at most $\mathfrak k$ and hence its dual group $(A\otimes H^*)^*\in\mathsf{CAGp}$ is non-trivial connected with weight at most $\mathfrak k$. Condition (3) thus yields $\Cocm((A\otimes H^*)^*)\setminus1\neq\emptyset$. By virtue of Corollary~\ref{C:min.contr.char} this means that there is a monomorphism $\varphi\colon A\otimes H^*\to\Cocm $. Now consider the mapping
\begin{equation*}
\psi\colon A\ni a\mapsto(\chi\mapsto\varphi(a\otimes \chi))\in\Hom(H^*,\Cocm ).
\end{equation*}
It is a standard procedure to verify that $\psi$ is well defined (that is, it takes its values in $\Hom(H^*,\Cocm)$) and that it is a morphism from the group $A$ to the groupoid $\Hom(H^*,\Cocm)$. We show that $\psi$ is in fact a monomorphism that takes its values in $\Mon(H^*,\Cocm)$. In view of Corollary~\ref{C:min.contr.char} this will secure that $\Cocm(H)$ contains an isomorphic copy of $A$.

To see that $\psi$ is a monomorphism, let $a\in A\setminus 0$. Fix $\chi\in H^*\setminus 0$. By virtue of ($\iota$), we have $a\otimes\chi\neq0$. Since $\varphi$ is a monomorphism, it follows $\psi(a)(\chi)=\varphi(a\otimes\chi)\neq0$. Thus, $\psi(a)\neq0$ and so $\psi$ is indeed a monomorphism. We finish the proof by showing that $\psi(a)\in\Mon(H^*,\Cocm)$ for every $a\in A$. If $a=0$ then $\psi(a)=0$ and so in this case the statement holds. If $a\neq0$ then, as above, $\psi(a)(\chi)\neq0$ for every $\chi\in H^*\setminus 0$, which means that $\psi(a)$ is a monomorphism.
\end{proof}

\begin{theorem}\label{T:grp.min.ext.am}
Let $\Gamma\in\mathsf{LCGp}\setminus\mathsf{CGp}$ be an amenable group, $X$ be a compact space and $\Flow$ be a minimal flow with a free point. If both $\Gamma$ and $X$ are second countable then the following statements hold:
\begin{enumerate}
\item[(a)] For every torsion-free group $A\in\mathsf{AbGp}$ with $\card(A)\leq\mathfrak{c}$ and every non-trivial connected group $G\in\mathsf{CAGp}$ with $\w(G)\leq\mathfrak{c}$, $\Cocm(G)$ contains an isomorphic copy of $A$.
\item[(b)] For every non-trivial connected second countable group $G\in\mathsf{CAGp}$, the groupoid $\Cocm(G)$ contains an isomorphic copy of the quotient group $\Cocc(G)/\tor(\Cocc(G))$.
\end{enumerate}
\end{theorem}
\begin{remark}\label{R:grp.min.ext.am}
In connection with statement (b) let us mention two facts. First, as observed in Remark~\ref{R:grp.min.ext}, $\Cocm(G)$ contains no non-trivial torsion elements. Therefore, it is necessary in (b) to factor out the torsion subgroup $\tor(\Cocc(G))$ of $\Cocc(G)$. Second, if the space $X$ is connected and the group $\Gamma$ has no non-trivial finite abelian quotient groups then, by Theorem~\ref{T:gimel.conn}, the group $\Cocc(G)$ is torsion-free. In such case the groupoid $\Cocm(G)$ contains an isomorphic copy of $\Cocc(G)$ itself.
\end{remark}
\begin{proof}[Proof of Theorem~\ref{T:grp.min.ext.am}]
We verify statement (a). From the equivalence of conditions (1) and (3) in Theorem~\ref{T:grp.min.ext} it follows that it suffices to check the following statement:
\begin{enumerate}
\item[($\circ$)] for every $G\in\mathsf{CAGp}$ connected with $\w(G)\leq\mathfrak{c}$, there is a minimal extension in $\Cocc(G)$.
\end{enumerate}
So fix $G\in\mathsf{CAGp}$ connected with weight at most $\mathfrak{c}$. The dual group $G^*$ of $G$ is then torsion-free with cardinality at most $\mathfrak{c}$ and so there is a monomorphism $\sigma\colon G^*\to\mathbb R$. Under the usual identifications $\mathbb R^*\cong\mathbb R$ and $G^{**}\cong G$, the dual morphism $\sigma^*\colon\mathbb R\to G$ of $\sigma$ has a dense image.

By virtue of \cite[Theorem~6]{Dir3}, there is an extension $\mathcal C\in\Cocc(\mathbb R)$ such that the flow $\mathcal F_{\mathcal C}\colon\Gamma\curvearrowright X\times\mathbb R$ is topologically transitive. Since the phase space $X\times\mathbb R$ of $\mathcal F_{\mathcal C}$ is a Polish space, transitivity means that $\mathcal F_{\mathcal C}$ possesses a dense orbit. Now, the map $\id_X\times\sigma^*\colon X\times\mathbb R\to X\times G$ is continuous, has a dense image and constitutes a morphism of flows $\mathcal F_{\mathcal C}\to\mathcal F_{\sigma^*\mathcal C}$. It follows that the flow $\mathcal F_{\sigma^*\mathcal C}$ also possesses a dense orbit and hence $F(\sigma^*\mathcal C)=G$. Thus, $\sigma^*\mathcal C\in\Cocc(G)$ is minimal, which verifies condition ($\circ$).

To verify statement (b), fix $G\in\mathsf{CAGp}$ non-trivial connected second countable. Then by our assumptions on $\Gamma$ and $X$,  $\Cocc(G)$ is a Polish group with the topology of uniform convergence on compact sets, and so it has cardinality at most $\mathfrak{c}$. Consequently, the quotient group $A=\Cocc(G)/\tor(\Cocc(G))$ is torsion-free and also has cardinality at most $\mathfrak{c}$. Since $\w(G)=\aleph_0<\mathfrak{c}$, we may apply part (a) of the theorem to conclude that $A$ has an isomorphic copy within $\Cocm(G)$.
\end{proof}

\begin{theorem}\label{T:grp.min.ext.simp}
Let $\Flow$ be a minimal flow with $\Gamma\in\mathsf{CLAC}$ simply connected and with $X$ compact. If $1\leq\rank(\pi^1(X))<\mathfrak{c}$ then the following statement holds:
\begin{enumerate}
\item[(i)] For every torsion-free group $A\in\mathsf{AbGp}$ with $\card(A)\leq\mathfrak{c}$ and every non-trivial connected group $G\in\mathsf{CAGp}$ with $\w(G)\leq\mathfrak{c}$, $\Cocm(G)$ contains an isomorphic copy of $A$.
\end{enumerate}
If, in addition, both $\Gamma$ and $X$ are second countable and $\Gamma$ is locally compact then the following statement also holds.
\begin{enumerate}
\item[(ii)] $\Cocm(G)$ contains an isomorphic copy of $\Cocc(G)$ for every non-trivial connected second countable group $G\in\mathsf{CAGp}$.
\end{enumerate}
\end{theorem}
\begin{proof}
We begin by verifying statement (i). By virtue of Theorem~\ref{T:grp.min.ext}, we need to show that $\Cocm(G)\setminus1$ is nonempty for every non-trivial connected group $G\in\mathsf{CAGp}$ with weight $\w(G)\leq\mathfrak{c}$. So let $G$ be such a group. Denote by $q\colon\Cocc\to\Coch/\tor(\Coch)$ the composition of the canonical quotient morphisms $\pi\colon\Cocc\to\Coch$ and $\Coch\to\Coch/\tor(\Coch)$. By Remark~\ref{R:lift.simply.con} and Theorem~\ref{T:gimel.conn}, the group $\Cocc$ is divisible and torsion-free and hence it is a rational linear space. Consequently, $\Coch/\tor(\Coch)$ is also a rational linear space and $q$ is thus an epimorphism of rational linear spaces. It follows that $q$ possesses a right inverse in the class of morphisms of rational linear spaces, that is, there is a morphism $r\colon\Coch/\tor(\Coch)\to\Cocc$ with $qr=\id$.

Now, by virtue of Theorem~\ref{T:min.exist.simp}, $\Coch/\tor(\Coch)$ has rational linear dimension at least $\mathfrak{c}$. Since the group $G^*$ is abelian torsion-free with rank at most $\mathfrak{c}$, there is a monomorphism $h\colon G^*\to\Coch/\tor(\Coch)$. Set $k=rh$. Since $h$ is a monomorphism and $h=qk$, it follows that $k$ is a monomorphism $G^*\to\Cocc$. Moreover, if $\chi\in G^*\setminus1$ then $1\neq h(\chi)=q(k(\chi))$ and so $\pi(k(\chi))\notin\tor(\Coch)$. Thus $k(\chi)\in\Cocc$ is minimal for every $\chi\in G^*\setminus1$ and so $k\in\Mon(G^*,\Cocm)\setminus 1$. Corollary~\ref{C:min.contr.char} then yields that $\Cocm(G)$ is non-trivial, as was to be shown.

Now assume that both $\Gamma$ and $X$ are second countable and that $\Gamma$ is locally compact. Fix $G\in\mathsf{CAGp}$ as in (ii). With the topology of uniform convergence on compact sets, $\Cocc(G)$ is a Polish group and hence it has cardinality at most $\mathfrak{c}$. Moreover, the group $\Cocc(G)$ is torsion-free by Theorem~\ref{T:gimel.conn}, and so statement (ii) follows from statement (i).
\end{proof}

\begin{theorem}\label{T:gp.min.free.ccl}
Let $\Flow$ be a minimal flow with $\Gamma\in\mathsf{LieGp}$ connected and with $X$ a compact connected manifold. Assume that $\mathcal F$ possesses a free cycle. Then the following statements hold.
\begin{enumerate}
\item[(i)] For every torsion-free group $A\in\mathsf{AbGp}$ with $\card(A)\leq\mathfrak{c}$ and every non-trivial connected group $G\in\mathsf{CAGp}$ with $\w(G)\leq\mathfrak{c}$, $\Cocm(G)$ contains an isomorphic copy of $A$.
\item[(ii)] For every non-trivial connected second countable group $G\in\mathsf{CAGp}$, the groupoid $\Cocm(G)$ contains an isomorphic copy of $\Cocc(G)$.
\end{enumerate}
\end{theorem}
\begin{proof}
We divide the proof into three steps. First, let us recall and fix some notation. Write $n=\rank(H_1^w(\mathcal F))$ and $n+m=\rank(H_1^w(X))$. By the assumptions of the theorem, $m\neq0$. Further, let $\mathcal G$ consist of the extensions $\mathcal C\in\Cocc$, whose induced morphisms $\mathcal C^{\sharp}\colon\pi_1(\Gamma\times X)\to\pi_1(\mathbb T^1)$ vanish. Finally, set $\mathcal G'=\mathcal G\cap\Cob$. Then $\mathcal G$ is clearly a subgroup of $\Cocc$ and $\mathcal G'$ is a subgroup of $\mathcal G$.

\emph{1st step.} We show that there is an exact sequence of abelian groups
\begin{equation*}
0\longrightarrow \mathbb Z^m\longrightarrow\Coch(\mathbb R)\longrightarrow\Coch\longrightarrow\Coch/(\mathcal G/\mathcal G')\longrightarrow0.
\end{equation*}
In order to do this, two facts need a verification. First, $\mathcal G/\mathcal G'$ is isomorphic to a subgroup of $\Coch$. Second, $\Coch(\mathbb R)$ contains an isomorphic copy of $\mathbb Z^m$ in such a way that $\Coch(\mathbb R)/\mathbb Z^m\cong\mathcal G/\mathcal G'$. The first fact follows immediately by definition of $\mathcal G$ and $\mathcal G'$. Now we turn to the proof of the second fact.

Let $p\colon\mathbb R\to\mathbb T^1$ be the covering morphism with kernel $\mathbb Z$. Then, by virtue of Lem\-ma~\ref{P:lifting.cocycle}, $\mathcal G$ consists of those extensions $\mathcal C\in\Cocc$, which lift across $p$ to an element of $\Cocc(\mathbb R)$. Consequently, $\widehat{p}\colon\Cocc(\mathbb R)\ni\mathcal D\to p\mathcal D\in\mathcal G$ is an isomorphism of groups. Set $\mathcal K=\widehat{p}(\Cob(\mathbb R))$. Then $\mathcal K$ is a subgroup of $\mathcal G'$ and we have an isomorphism
\begin{equation*}
\Coch(\mathbb R)=\Cocc(\mathbb R)/\Cob(\mathbb R)\cong\mathcal G/\mathcal K.
\end{equation*} 

Denote by $\psi$ the isomorphism $C_z(X,\mathbb T^1)\to\Cob$ acting by the rule $\psi(\xi)=\co(\xi)$. Given $\xi\in C_z(X,\mathbb T^1)$, the following conditions are equivalent:
\begin{itemize}
\item $\co(\xi)\in\mathcal K$,
\item $\co(\xi)$ lifts across $p$ to an element of $\Cob(\mathbb R)$,
\item $\xi$ lifts across $p$ to an element of $C_z(X,\mathbb R)$.
\end{itemize}
Consequently, under the usual isomorphism $C_z(X,\mathbb T^1)\cong C_z(X,\mathbb R)\oplus\pi^1(X)$, we have 
\begin{equation}\label{Eq:psi.-1.K}
\psi^{-1}(\mathcal K)=C_z(X,\mathbb R).
\end{equation}
Similarly, by virtue of Lemma~\ref{L:tor.Lie.mnfld}, the following conditions are equivalent for every $\xi\in C_z(X,\mathbb T^1)$:
\begin{itemize}
\item $\co(\xi)\in\mathcal G'$,
\item $\co(\xi)^{\sharp}\pi_1(\Gamma\times X)=0$,
\item $\xi^{\sharp}\mathcal F_z^{\sharp}\pi_1(\Gamma)=0$,
\item $\xi^{\sharp}H_1^w(\mathcal F)=0$.
\end{itemize}
Therefore, 
\begin{equation}\label{Eq:psi.-1.G'}
\psi^{-1}(\mathcal G')=C_z(X,\mathbb R)\oplus\{f\in\pi^1(X) : f^{\sharp}H_1^w(\mathcal F)=0\}.
\end{equation}

Now, from the identities (\ref{Eq:psi.-1.K}) and (\ref{Eq:psi.-1.G'}) it follows that the group $\mathcal K$ is a direct summand in $\mathcal G'$ and its complementary summand is isomorphic to the group $A=\{f\in\pi^1(X) : f^{\sharp}H_1^w(\mathcal F)=0\}$. Further, recall from Section~\ref{S:ex.prescr.sect.1} in Chapter~\ref{S:alg.top.asp} that there is an isomorphism of groups $\sigma\colon \pi^1(X)\to\Hom(H_1^w(X),\mathbb Z)$, acting by the rule $\sigma(f)=f^{\sharp}$. Also, under the natural isomorphisms $\Hom(H_1^w(X),\mathbb Z)\cong\mathbb Z^{n+m}\cong\mathbb Z^n\oplus\mathbb Z^m$, the group $\sigma(A)$ corresponds to $\mathbb Z^m$, and so we have a splitting $\mathcal G'=\mathcal K\oplus\mathbb Z^m$. It follows that there are isomorphisms of groups
\begin{equation*}
\mathcal G/\mathcal G'\cong\mathcal G/(\mathcal K\oplus\mathbb Z^m)\cong(\mathcal G/\mathcal K)/\mathbb Z^m\cong\Coch(\mathbb R)/\mathbb Z^m,
\end{equation*}
as was to be shown.

\emph{2nd step.} We verify statement (i) from the theorem. Since every torsion-free abelian group with cardinality at most $\mathfrak{c}$ can be viewed as a subgroup of $\mathbb R$, by virtue of Theorem~\ref{T:grp.min.ext} it suffices to show that $\Cocm$ contains an isomorphic copy of $\mathbb R$.

Let $\pi$ denote the quotient morphism $\mathcal G\to\mathcal G/\mathcal G'$. By the first step of the proof, $\mathbb Z^m$ can be viewed as a subgroup of $\Coch(\mathbb R)$ in such a way that $\mathcal G/\mathcal G'\cong\Coch(\mathbb R)/\mathbb Z^m$. Since $\Coch(\mathbb R)$ is a real (and hence also rational) linear space, we get inclusions $\mathbb Z^m\subseteq\mathbb Q^m\subseteq\Coch(\mathbb R)$. Clearly, $\mathbb Q^m$ is a direct summand in $\Coch(\mathbb R)$; denote its complementary summand by $L$. Since $\Coch(\mathbb R)\neq0$, we have $\LSdim(\Coch(\mathbb R))\geq1$ and hence $\QLSdim(\Coch(\mathbb R))\geq\mathfrak{c}$. Consequently, $L$ is a rational linear space with dimension $\QLSdim(L)\geq\mathfrak{c}$. Also, by definition of $L$, there is an isomorphism $\mathcal G/\mathcal G'\cong L\oplus(\mathbb Q/\mathbb Z)^m$.

Now, the group $\mathcal G$ is torsion free and divisible by its definition, and $(\mathbb Q/\mathbb Z)^m$ is a direct summand in $\mathcal G/\mathcal G'$ with a torsion-free complementary summand $L$. Consequently, $\pi^{-1}((\mathbb Q/\mathbb Z)^m)$ is a divisible subgroup of $\mathcal G$ and we may write $\mathcal G=\mathcal L\oplus\pi^{-1}((\mathbb Q/\mathbb Z)^m)$, where $\mathcal L$ is a rational linear subspace of $\mathcal G$. Given $l\in L$, fix $\mathcal C_l\in\mathcal L$ and $\mathcal D_l\in\pi^{-1}((\mathbb Q/\mathbb Z)^m)$ with $l=\pi(\mathcal C_l+\mathcal D_l)$. An elementary argument shows that for $l_1,l_2\in L$ with $l_1\neq l_2$, one has $\mathcal C_{l_1}\neq\mathcal C_{l_2}$. Therefore, $\card(\mathcal L)\geq\card(L)\geq\mathfrak{c}$, and hence $\mathcal L$ contains an isomorphic copy of $\mathbb R$. Now it remains only to show that $\mathcal L\subseteq\Cocm$. This follows from two facts. First, the non-minimal extensions in $\mathcal G$ project to $\tor(\mathcal G/\mathcal G')=(\mathbb Q/\mathbb Z)^m$ via $\pi$. Second, by definition of $\mathcal L$, we have $\mathcal L\cap \pi^{-1}((\mathbb Q/\mathbb Z)^m)=0$.

\emph{3rd step.} We finish the proof of the theorem by verifying statement (ii).

So let $G\in\mathsf{CAGp}$ be a non-trivial connected second countable group. Then the group $\Cocc(G)$ is torsion-free by connectedness of $\Gamma$ and $X$, see Theorem~\ref{T:gimel.conn}. Moreover, $\Cocc(G)$ is a Polish group with the topology of uniform convergence on compact sets and hence its cardinality is at most $\mathfrak{c}$. Thus, statement (ii) follows immediately by applying statement (i) to $A=\Cocc(G)$.
\end{proof}

Our assumption on the existence of a free cycle for the flow $\mathcal F$ in Theorem~\ref{T:gp.min.free.ccl} can not be dropped, as the following example demonstrates.

\begin{example}\label{E:free.ccl.drpd}
Let $G$ be a compact connected Lie group and $H$ be a closed subgroup of $G$ with the topological dimension $\dim(H)=k$. Consider the usual transitive action $\mathcal F\colon G\curvearrowright G/H$ and choose $z=H$ as the base point for $G/H$. If $G$ acts in a minimal way on a space $X$ then, by compactness of $G$, the action is necessarily transitive. Consequently, the space $X$ is (up to homeomorphism) a homogeneous space of $G$ and hence $\dim(X)\leq\dim(G)$. It follows that for $l>k$ there is no minimal $G$-flow on $G/H\times\mathbb T^l$ and hence $\Cocm(\mathbb T^l)=e$. From the point of view of Theorem~\ref{T:gp.min.free.ccl}, this must be due to the fact that the flow $\mathcal F$ does not possess a free cycle. To see this, let $\mathcal F_z=q\colon G\to G/H$ be the canonical projection. Since $G$ is a locally trivial fibre bundle over $G/H$ with the fibre $H$ and since $G$ and $G/H$ are arc-wise connected by our assumptions, the long exact homotopy sequence of the bundle yields exactness of the upper row of the diagram in Figure~\ref{Fig:no.free.ccl} below (see \cite[Section~9.8]{FomFuch}). It follows that there is an isomorphism of groups
\begin{equation*}
\pi_1(G/H)/\im(q^{\sharp})=\pi_1(G/H)/\ker(\partial)\cong H/H_0,
\end{equation*}
where the group $H/H_0=\pi_0(H)$ of the path-components of $H$ is finite ($H$ being a compact Lie group). Since both $p_G$ and $p_{G/H}$ are epimorphisms and the diagram in Figure~\ref{Fig:no.free.ccl} commutes, it follows that the group
\begin{equation*}
H_1^w(G/H)/H_1^w(\mathcal F)=H_1^w(G/H)/\mathcal F_z^{\sharp}H_1^w(G)=H_1^w(G/H)/\im(q^{\sharp})
\end{equation*}
is also finite. This means that the flow $\mathcal F$ does not possess a free cycle indeed.
\begin{figure}[ht]
\[\minCDarrowwidth20pt\begin{CD}
\pi_1(G) @>q^{\sharp}>> \pi_1(G/H) @>\partial>> H/H_0 @>>> 0\\
@Vp_GVV @Vp_{G/H}VV  \\
H_1^w(G) @>>q^{\sharp}> H_1^w(G/H) 
\end{CD}\]
\caption{The flow $\mathcal F$ does not possess a free cycle}
\label{Fig:no.free.ccl}
\end{figure}
\end{example}

\section{First cohomology groups with integer coefficients}\label{S:circ.case}

Let $\Flow$ be a mi\-ni\-mal flow. Our aim in this section is to express the inclusions $\Cob\subseteq\Cocc$ and $F(\Daleth)\subseteq\Gimel$ in terms of inclusions of familiar (topological) abelian groups, as well as to compute the first cohomology group $\Coch$ of $\mathcal F$ with integer coefficients. We concentrate on the following situations:
\begin{itemize}
\item[(i)] the acting group $\Gamma$ of $\mathcal F$ is a simply connected Lie group and the phase space $X$ of $\mathcal F$ is a compact second countable space with $\pi^1(X)\neq0$, see Theorem~\ref{T:F.in.Gim.CLAC},
\item[(ii)] the acting group $\Gamma$ of $\mathcal F$ is a connected Lie group, the phase space $X$ of $\mathcal F$ is a compact manifold, the flow $\mathcal F$ is topologically free and possesses a free cycle, see Theorem~\ref{T:F.in.Gim.free.ccl}.
\end{itemize}
We also know from Corollary~\ref{C:structure.coc.simply} that, in case (i), the group $\Cocc$ is the additive topological group of a real separable Banach space. In Proposition~\ref{P:Cocc.top.free.ccl} we show that, in case (ii), the structure of a Banach space is carried by the divisible subgroup $\Div(\Cocc)$ of $\Cocc$, which forms an open topological direct summand in $\Cocc$.

Before turning to the mentioned results we prove two propositions, which indicate how the results of this section can be used to determine the possible values of the functor $F$ on the groups $\Cocc(G)$ for every $G\in\mathsf{CAGp}$. First, recall that for a minimal flow $\mathcal F$ we use notation $\Daleth$ for the free extension of $\mathcal F$ and that this extension takes its values in $\Gimel=(\Cocc)_d^*\in\mathsf{CAGp}$.

\begin{proposition}\label{P:pos.val.F.gen}
Let $\Flow$ be a minimal flow, $G\in\mathsf{CAGp}$ and $H\sbgp G$. Then the following conditions are equivalent:
\begin{enumerate}
\item[(a)] there is $\mathcal C\in\Cocc(G)$ with $F(\mathcal C)=H$,
\item[(b)] there is $q\in\Hom(\Gimel,G)$ with $qF(\Daleth)=H$,
\item[(c)] there is $k\in\Hom(G^*,\Cocc)$ with $k^{-1}\Cob=H^{\perp}$.
\end{enumerate}
\end{proposition}
\begin{remark}\label{R:pos.val.F.gen}
Observe that the equality $k^{-1}\Cob=H^{\perp}$ from condition (c) can be reformulated in the following way:
\begin{itemize}
\item $k\left(H^{\perp}\right)\subseteq\Cob$ and $k\left(G^*\setminus H^{\perp}\right)\subseteq\Cocc\setminus\Cob$.
\end{itemize}
\end{remark}
\begin{proof}[Proof of Proposition~\ref{P:pos.val.F.gen}]
The equivalence of (a) and (b) follows from the following facts:
\begin{itemize}
\item every $\mathcal C\in\Cocc(G)$ is of the form $\mathcal C=q\Daleth$ for some $q\in\Hom(\Gimel,G)$; this follows directly from the definition of the free extension of $\mathcal F$,
\item $F(q\Daleth)=qF(\Daleth)$ for every $q\in\Hom(\Gimel,G)$; this follows from Theorem~\ref{T:functor.E.def}(3).
\end{itemize}

The equivalence of (b) and (c) follows from the following facts:
\begin{itemize}
\item the elements of $\Hom(\Gimel,G)$ are in a one-to-one correspondence with the elements of $\Hom(G^*,\Cocc)$; this follows from the isomorphism $(\Gimel)^*\cong(\Cocc)_d$,
\item a morphism $q\in\Hom(\Gimel,G)$ maps $F(\Daleth)$ onto $H$ if and only if its dual morphism $q^*\in\Hom(G^*,(\Gimel)^*)$ satisfies $\left(q^*\right)^{-1}\left(F(\Daleth)^{\perp}\right)=H^{\perp}$; this is clear,
\item under the usual isomorphism $(\Gimel)^*\cong(\Cocc)_d$, $F(\Daleth)^{\perp}$ corresponds to $\Cob$; this follows from Theorem~\ref{T:prop.of.gimel}(2).
\end{itemize}
\end{proof}

\begin{proposition}\label{P:pos.val.F.gen.2}
Let $\Flow$ be a minimal flow, $G\in\mathsf{CAGp}$ and $H\sbgp G$. Then the following conditions are equivalent:
\begin{enumerate}
\item[(i)] there is $\mathcal C\in\Cocc(G)$ with $G_{\mathcal C}=G$ and $F(\mathcal C)=H$,
\item[(ii)] there is an epimorphism $q\colon\Gimel\to G$ with $qF(\Daleth)=H$,
\item[(iii)] there is a monomorphism $k\colon G^*\to\Cocc$ with $k^{-1}\Cob=H^{\perp}$.
\end{enumerate}
\end{proposition}
\begin{remark}\label{R:pos.val.F.gen.2}
We recall from Proposition~\ref{P:further.prop.of.gimel}(d) that for a group $G\in\mathsf{CAGp}$ and an extension $\mathcal C\in\Cocc(G)$, the symbol $G_{\mathcal C}$ is used to denote the closed subgroup of $G$ generated by the values of $\mathcal C$. Also, similarly to Remark~\ref{R:pos.val.F.gen}, condition (iii) can be reformulated in the following way:
\begin{itemize}
\item the group $\Cocc$ contains an isomorphic copy of $G^*$ in such a way that $G^*\cap\Cob=H^{\perp}$.
\end{itemize}
\end{remark}
\begin{proof}[Proof of Proposition~\ref{P:pos.val.F.gen.2}]
The equivalence of (i) and (ii) is proved by the same line of reasoning as the equivalence of (a) and (b) in Proposition~\ref{P:pos.val.F.gen}, with the following additional argument:
\begin{itemize}
\item given $q\in\Hom(\Gimel,G)$ and $\mathcal C=q\Daleth$, the equality $G_{\mathcal C}=G$ holds if and only if $q$ is an epimorphism $\Gimel\to G$; this follows from Proposition~\ref{P:further.prop.of.gimel}($\iota$).
\end{itemize}

Likewise, the equivalence of (ii) and (iii) follows by the same argument as the equivalence of (b) and (c) in Proposition~\ref{P:pos.val.F.gen}; one uses, additionally, that $q\in\Hom(\Gimel,G)$ is an epimorphism if and only if $k=q^*\in\Hom(G^*,\Cocc)$ is a monomorphism.
\end{proof}

\begin{theorem}\label{T:F.in.Gim.CLAC}
Let $\Flow$ be a minimal flow with $\Gamma\in\mathsf{LieGp}$ sim\-ply connected and with $X$ compact. Assume that the space $X$ is second countable and $\pi^1(X)\neq0$. Then there are isomorphisms of groups
\begin{equation}\label{Eq:F.in.Gim.CLAC}
\begin{array}{ccccccccc}
\Cob & \cong & 0 & \oplus & \mathbb R & \oplus & \pi^1(X) & \subseteq & {}\\
{} & \subseteq & \mathbb R & \oplus & \mathbb R & \oplus & \big(\mathbb Q\otimes\pi^1(X)\big) & \cong & \Cocc
\end{array}
\end{equation}
and topological isomorphisms
\begin{equation}\label{Eq:F.in.Gim.CLAC2}
\begin{array}{ccccccccc}
F(\Daleth) & \cong & b\mathbb R & \times & 1 & \times & \pi^1(X)^{\perp} & \subseteq & {}\\
{} & \subseteq & b\mathbb R & \times & b\mathbb R & \times & \big(\mathbb Q\otimes\pi^1(X)\big)^* & \cong & \Gimel.
\end{array}
\end{equation}
Also, there is an isomorphism
\begin{equation}\label{Eq:F.in.Gim.CLAC3}
\Coch\cong\mathbb R\oplus\left((\mathbb Q/\mathbb Z)\otimes\pi^1(X)\right).
\end{equation}
If, in addition, the group $\pi^1(X)$ is finitely generated and $k=\rank\left(\pi^1(X)\right)$ then the three isomorphisms above take the form
\begin{equation}\label{Eq:F.in.Gim.CLAC.CW}
\begin{array}{ccccccccc}
\Cob & \cong & 0 & \oplus & \mathbb R & \oplus & \mathbb Z^k & \subseteq & {}\\
{} & \subseteq & \mathbb R & \oplus & \mathbb R & \oplus & \mathbb Q^k & \cong & \Cocc,
\end{array}
\end{equation}
respectively,
\begin{equation}\label{Eq:F.in.Gim.CLAC2.CW}
\begin{array}{ccccccccc}
F(\Daleth) & \cong & b\mathbb R & \times & 1 & \times & \left(\mathbb Z^{\perp}\right)^k & \subseteq & {}\\
{} & \subseteq & b\mathbb R & \times & b\mathbb R & \times & \left(\mathbb Q^*\right)^k & \cong & \Gimel,
\end{array}
\end{equation}
and, finally,
\begin{equation}\label{Eq:F.in.Gim.CLAC3.CW}
\Coch\cong \mathbb R\oplus (\mathbb Q/\mathbb Z)^k.
\end{equation}
\end{theorem}
\begin{remark}\label{R:F.in.Gim.CLAC}
We wish to mention the following facts.
\begin{itemize}
\item In the topological isomorphisms (\ref{Eq:F.in.Gim.CLAC2}) and (\ref{Eq:F.in.Gim.CLAC2.CW}) the group $\mathbb R$ is equipped with its usual topology and the groups $\mathbb Q\otimes\pi^1(X)$, $\mathbb Q$ are equipped with the discrete topology. All the other isomorphisms are meant in the algebraic sense.
\item The annihilator $\pi^1(X)^{\perp}$ in (\ref{Eq:F.in.Gim.CLAC2}) is taken in $\left(\mathbb Q\otimes\pi^1(X)\right)^*$. The group $\pi^1(X)$ is identified here with the tensor product $\mathbb Z\otimes\pi^1(X)\subseteq\mathbb Q\otimes\pi^1(X)$ under the usual isomorphism $\pi^1(X)\ni f\mapsto 1\otimes f\in\mathbb Z\otimes\pi^1(X)$.
\item Finally, recall that the group $\pi^1(X)$ is finitely generated if $X$ is a compact connected manifold.
\end{itemize}
\end{remark}
\begin{proof}[Proof of Theorem~\ref{T:F.in.Gim.CLAC}]
Since both $\Gamma$ and $X$ are connected, $\Cocc$ is a torsion-free abelian group by Theorem~\ref{T:gimel.conn}. Moreover, the group $\Cocc$ is divisible by virtue of Theorem~\ref{T:lift.simply.con} and Remark~\ref{R:lift.simply.con}. Also, since $\Gamma$ and $X$ are both locally compact and second countable, $\Cocc$ is a Polish group with the topology of uniform convergence on compact sets and it has therefore cardinality at most $\mathfrak{c}$. Finally, it follows from our assumptions on $X$ that $1\leq\rank\left(\pi^1(X)\right)<\mathfrak{c}$. Indeed, the first inequality is secured by the assumption $\pi^1(X)\neq0$ and the second inequality follows since the space $X$ is second countable.

If every element of $\Cob$ is identified with its base point preserving transfer function then there are the usual isomorphisms $\Cob\cong C_z(X,\mathbb T^1)\cong C_z(X,\mathbb R)\oplus\pi^1(X)$. As observed above, the group $\Cocc$ is divisible and so it contains the divisible hull $\mathbb Q\otimes\Cob$ of $\Cob$. Moreover, there are isomorphisms
\begin{equation*}
\mathbb Q\otimes\Cob\cong\mathbb Q\otimes\left(C_z(X,\mathbb R)\oplus\pi^1(X)\right)\cong C_z(X,\mathbb R)\oplus\left(\mathbb Q\otimes\pi^1(X)\right).
\end{equation*}
Since the group $\mathbb Q\otimes\Cob$ is divisible by definition, it is a direct summand in $\Cocc$; its complementary summand $\mathfrak{D}$ intersects $\Cob$ only at $0$ and it is therefore a divisible subgroup of the groupoid $\Cocm$. Thus, we have
\begin{equation*}
\Cocc\cong\mathfrak{D}\oplus\left(\mathbb Q\otimes\Cob\right)\cong\mathfrak{D}\oplus C_z(X,\mathbb R)\oplus\left(\mathbb Q\otimes\pi^1(X)\right).
\end{equation*}

Now $\mathfrak{D}$ is a divisible torsion-free abelian group and it is therefore a rational linear space. Write $\mathfrak{k}=\QLSdim\left(\mathfrak{D}\right)$, we claim that $\mathfrak{k}=\mathfrak{c}$. First, since $\card\left(\Cocc\right)\leq\mathfrak{c}$, it follows that $\mathfrak{k}\leq\mathfrak{c}$. Further, by virtue of Theorem~\ref{T:grp.min.ext.simp}, there is a monomorphism of groups $\sigma\colon\mathbb R\to\Cocc$ with $\im(\sigma)\subseteq\Cocm$. Let $q$ be the projection $\Cocc\to\mathfrak{D}$. Since the complementary summand $\mathbb Q\otimes\Cob$ of $\mathfrak{D}$ in $\Cocc$ intersects $\Cocm$ only at $0$, it follows that $q\sigma$ is a monomorphism $\mathbb R\to\mathfrak{D}$. This shows that $\mathfrak{k}\geq\mathfrak{c}$. Thus, indeed, $\mathfrak{k}=\mathfrak{c}$ and $\mathfrak{D}\cong\mathbb Q^{(\mathfrak{c})}\cong\mathbb R$.

Since $X$ is a non-degenerate second countable continuum, we have $\w(X)=\aleph_0$. Consequently, $C_z(X,\mathbb R)$ is a real Banach space with dimension $\BSdim\left(C_z(X,\mathbb R)\right)=\w(X)=\aleph_0$ and hence its dimension as a real linear space is $\LSdim\left(C_z(X,\mathbb R)\right)=\mathfrak{c}$ (see Subsection~\ref{Sub:spc.real.fnctns}). It follows that as an abelian group, $C_z(X,\mathbb R)$ is isomorphic to $\mathbb R$.

The observations made so far may be summarized in the form of the following isomorphisms
\begin{equation*}
\begin{array}{ccccccccc}
\Cob & \cong & 0 & \oplus & C_z(X,\mathbb R) & \oplus & \pi^1(X) & \cong & {}\\
{} & \cong & 0 & \oplus & \mathbb R & \oplus & \pi^1(X) & \subseteq & {}\\
{} & \subseteq & \mathbb R & \oplus & \mathbb R & \oplus & \big(\mathbb Q\otimes\pi^1(X)\big) & \cong & {}\\
{} & \cong & \mathfrak{D} & \oplus & C_z(X,\mathbb R) & \oplus & \big(\mathbb Q\otimes\pi^1(X)\big) & \cong & \Cocc,
\end{array}
\end{equation*}
which verify (\ref{Eq:F.in.Gim.CLAC}). Further, from (\ref{Eq:F.in.Gim.CLAC}) it follows at once that the following isomorphisms hold:
\begin{equation*}
\begin{array}{ccccccccc}
\left(\Cob\right)^{\perp} & \cong & \left(\mathbb R_d\right)^* & \times & 1 & \times & \pi^1(X)^{\perp} & \subseteq & {}\\
{} & \subseteq & \left(\mathbb R_d\right)^* & \times & \left(\mathbb R_d\right)^* & \times & \big(\mathbb Q\otimes\pi^1(X)\big)^* & \cong & \left(\Cocc\right)_d^*.
\end{array}
\end{equation*}
These isomorphisms lead to (\ref{Eq:F.in.Gim.CLAC2}) with the help of the following facts:
\begin{itemize}
\item $\left(\Cocc\right)_d^*=\Gimel$; this follows from Theorem~\ref{T:exist.free.ext},
\item $\left(\Cob\right)^{\perp}=F(\Daleth)$; this follows from Theorem~\ref{T:prop.of.gimel}(2),
\item $\left(\mathbb R_d\right)^*\cong b\mathbb R$; this is clear.
\end{itemize}
Finally, (\ref{Eq:F.in.Gim.CLAC3}) follows from (\ref{Eq:F.in.Gim.CLAC}) and from the isomorphism $(\mathbb Q\otimes\pi^1(X))/\pi^1(X)\cong(\mathbb Q/\mathbb Z)\otimes\pi^1(X)$.

Now assume that the group $\pi^1(X)$ is finitely generated and set $k=\rank(\pi^1(X))$. Then the inclusion $\pi^1(X)\subseteq\mathbb Q\otimes\pi^1(X)$ takes the form $\mathbb Z^k\subseteq\mathbb Q^k$ and so the isomorphisms (\ref{Eq:F.in.Gim.CLAC.CW})--(\ref{Eq:F.in.Gim.CLAC3.CW}) follow from the isomorphisms (\ref{Eq:F.in.Gim.CLAC})--(\ref{Eq:F.in.Gim.CLAC3}), respectively.
\end{proof}

\begin{theorem}\label{T:F.in.Gim.free.ccl}
Let $\Flow$ be a minimal flow with $\Gamma\in\mathsf{LieGp}$ connected and with $X$ a compact connected manifold. Set $n=\rank(H_1^w(\mathcal F))$, $n+m=\rank(H_1^w(X))$ and denote by $d_1,\dots,d_n$ the elementary divisors of $H_1^w(\mathcal F)$ in $H_1^w(X)$. If $\mathcal F$ is topologically free and possesses a free cycle then there are isomorphisms of groups
\begin{equation}\label{Eq:Cob.in.Cocc}
\begin{array}{ccccccccccc}
\Cob & \cong & 0 & \oplus & \mathbb R & \oplus & \mathbb Z^m & \oplus & \bigoplus_{i=1}^nd_i\mathbb Z & \subseteq & {} \\
{} & \subseteq & \mathbb R & \oplus & \mathbb R & \oplus & \mathbb Q^m & \oplus & \mathbb Z^n & \cong & \Cocc
\end{array}
\end{equation}
and topological isomorphisms
\begin{equation}\label{Eq:Cob.in.Cocc2}
\begin{array}{ccccccccccc}
F(\Daleth) & \cong & b\mathbb R & \times & 1 & \times & \left(\mathbb Z^{\perp}\right)^m & \times & \prod_{i=1}^n\mathbb Z_{d_i} & \subseteq & {} \\
{} & \subseteq & b\mathbb R & \times & b\mathbb R & \times & \left(\mathbb Q^*\right)^m & \times & \mathbb T^n & \cong & \Gimel.
\end{array}
\end{equation}
Also, there is an isomorphism of groups
\begin{equation}\label{Eq:Cob.in.Cocc3}
\Coch\cong \mathbb R\oplus (\mathbb Q/\mathbb Z)^m\oplus\bigoplus_{i=1}^n\mathbb Z_{d_i}.
\end{equation}
\end{theorem}
\begin{remark}\label{R:F.in.Gim.free.ccl}
In the isomorphisms (\ref{Eq:Cob.in.Cocc2}) the groups $\mathbb R$ and $\mathbb T^n$ are equipped with their usual topologies and the group $\mathbb Q$ carries the discrete topology. The isomorphisms (\ref{Eq:Cob.in.Cocc}) and (\ref{Eq:Cob.in.Cocc3}) are meant in the algebraic sense.
\end{remark}
\begin{proof}[Proof of Theorem~\ref{T:F.in.Gim.free.ccl}]
We divide the proof into four steps.

\emph{1st step.} We recall some facts from Theorem~\ref{T:ex.seq.Lie.mfld} and from its proof.

Write $\eta=d_1\dots d_n$. We shall use the usual identifications
\begin{equation*}
H_1^w(\mathcal F)\cong d_1\mathbb Z\oplus\dots\oplus d_n\mathbb Z\subseteq\mathbb Z^n\cong\mathbb Z^n\oplus0\subseteq\mathbb Z^n\oplus\mathbb Z^m\cong H_1^w(X),
\end{equation*}
and denote the corresponding standard basis of $H_1^w(X)$ by $e_1,\dots,e_{n+m}$. By virtue of the isomorphism (\ref{Eq:mrphsm.R}) from the proof of Lemma~\ref{L:(co)hmtp.rel}, there are maps $\xi_i\in\pi^1(X)$ ($i=1,\dots,n$), whose induced morphisms $\left(\xi_i\right)^{\sharp}\colon H_1^w(X)\to\mathbb Z$ satisfy $\left(\xi_i\right)^{\sharp}(e_j)=\delta_{ij}(\eta/d_j)$ for $j=1,\dots,n$ and $\left(\xi_i\right)^{\sharp}\mathbb Z^m=0$.

Given $i\in\{1,\dots,n\}$, the extension $\co(\xi_i)\in\Cob$ lifts across $\kappa_{\eta}\colon\mathbb T^1\to\mathbb T^1$ to an extension $\mathcal C_i\in\Cocc$. The following statements hold:
\begin{itemize}
\item the extensions $\mathcal C_1,\dots,\mathcal C_n$ are independent in $\Cocc$ and so they generate a free abelian group $\mathcal S=\langle\mathcal C_1,\dots,\mathcal C_n\rangle$ with $\rank(\mathcal S)=n$,
\item there is a direct sum $\Cocc=\Div(\Cocc)\oplus\mathcal S$,
\item $\Cob\cap\mathcal S$ is a free abelian group with rank $n$ and the extensions $\mathcal C_1^{d_1},\dots,\mathcal C_n^{d_n}$ form a basis for $\Cob\cap\mathcal S$.
\end{itemize}

\emph{2nd step.} We show that the groupoid $\Cocm$ possesses a divisible subgroup $\mathfrak{D}$ such that there is a direct sum
\begin{equation}\label{Eq:dir.sum.Coc.D}
\Cocc=\mathfrak{D}\oplus C_z(X,\mathbb R)\oplus\mathbb Q^m\oplus\mathbb Z^n.
\end{equation}

For $i=1,\dots,n$ let $\zeta_i$ be the base point preserving transfer function for $\mathcal C_i^{d_i}$. Then 
\begin{equation*}
\co\big(\zeta_i^{\eta}\big)=\co(\zeta_i)^{\eta}=\big(\mathcal C_i^{d_i}\big)^{\eta}=\big(\mathcal C_i^{\eta}\big)^{d_i}=\co(\xi_i)^{d_i}=\co\big(\xi_i^{d_i}\big),
\end{equation*}
and hence $\zeta_i^{\eta}=\xi_i^{d_i}$. It follows that the induced morphisms $(\zeta_i)^{\sharp}\colon H_1^w(X)\to\mathbb Z$ satisfy $(\zeta_i)^{\sharp}(e_j)=\delta_{ij}$ for $j=1,\dots n$ and $(\zeta_i)^{\sharp}\mathbb Z^m=0$. Also, since $\xi_1,\dots,\xi_n\in\pi^1(X)$ and $\pi^1(X)$ is a direct summand in $C_z(X,\mathbb T^1)$ with a torsion-free complementary summand $C_z(X,\mathbb R)$, we get $\zeta_1,\dots,\zeta_n\in\pi^1(X)$. Further, by virtue of the isomorphism (\ref{Eq:mrphsm.R}), there are maps $\vartheta_k\in\pi^1(X)$ ($k=1,\dots,m$) with $(\vartheta_k)^{\sharp}(e_l)=\delta_{kl}$ for $l=1,\dots,n+m$. By applying (\ref{Eq:mrphsm.R}) once more we see that the maps $\vartheta_k$ ($k=1,\dots,m$) and $\zeta_i$ ($i=1,\dots,n$) together form a basis for the group $\pi^1(X)\cong\mathbb Z^{n+m}$. Thus, under the usual identification $\Cob\cong C_z(X,\mathbb T^1)\supseteq C_z(X,\mathbb R)$, there are isomorphisms
\begin{equation}\label{Eq:cob.split.m.n}
\begin{split}
\Cob&=C_z(X,\mathbb R)\oplus\left\langle\co(\vartheta_1),\dots,\co(\vartheta_m)\right\rangle\oplus\left\langle \co(\zeta_1),\dots,\co(\zeta_n)\right\rangle\\
&\cong C_z(X,\mathbb R)\oplus\mathbb Z^m\oplus\mathbb Z^n.
\end{split}
\end{equation}

Write $\mathcal S'=\langle\co(\vartheta_1),\dots,\co(\vartheta_m)\rangle$. Since $\vartheta_k^{\sharp}H_1^w(\mathcal F)=0$ for every $k=1,\dots,m$, it follows from Corollary~\ref{C:lift.cob.Lie.gp} that the elements of $\mathcal S'$ lift across all $\kappa_l\colon\mathbb T^1\to\mathbb T^1$ ($l\in\mathbb N$) to elements of $\Cocc$. In other words, $\mathcal S'\subseteq\Div(\Cocc)$ and hence also $\mathbb Q\otimes\mathcal S'\subseteq\Div(\Cocc)$. Moreover, since $\mathcal S'\cong\mathbb Z^m$, it follows that $\mathbb Q\otimes\mathcal S'\cong\mathbb Q^m$. Further, we have $C_z(X,\mathbb R)\subseteq\Div(\Cocc)$ and, by virtue of (\ref{Eq:cob.split.m.n}), also $C_z(X,\mathbb R)\cap(\mathbb Q\otimes\mathcal S')=0$. This means that $C_z(X,\mathbb R)\oplus(\mathbb Q\otimes\mathcal S')$ is a divisible subgroup of $\Div(\Cocc)$ and hence $\Div(\Cocc)=\mathfrak{D}\oplus C_z(X,\mathbb R)\oplus(\mathbb Q\otimes\mathcal S')$ for an appropriate divisible subgroup $\mathfrak{D}$ of $\Div(\Cocc)$. Thus
\begin{equation}\label{Eq:GpExt.dcmp.tp.free}
\begin{split}
\Cocc&=\Div(\Cocc)\oplus\mathcal S=\mathfrak{D}\oplus C_z(X,\mathbb R)\oplus(\mathbb Q\otimes\mathcal S')\oplus\mathcal S\\
&\cong\mathfrak{D}\oplus C_z(X,\mathbb R)\oplus\mathbb Q^m\oplus\mathbb Z^n,
\end{split}
\end{equation}
and so it remains only to show that $\mathfrak{D}\subseteq\Cocm$.

Given $0\neq\zeta\in\langle\zeta_1,\dots,\zeta_n\rangle$, we have $\zeta^{\sharp}H_1^w(\mathcal F)\neq0$. By virtue of Corollary~\ref{C:lift.cob.Lie.gp} it follows that none of the non-zero elements of $\langle\co(\zeta_1),\dots,\co(\zeta_n)\rangle$ is divisible in $\Cocc$ and hence, by (\ref{Eq:cob.split.m.n}), $\Cob\cap\Div(\Cocc)=C_z(X,\mathbb R)\oplus\mathcal S'$. Now let $\mathcal C\in\mathfrak{D}$ be a non-minimal extension, we show that $\mathcal C$ is trivial. Fix $k\in\mathbb N$ with $\mathcal C^k\in\Cob$. Then
\begin{equation*}
\mathcal C^k\in\Cob\cap\mathfrak{D}=(\Cob\cap\Div(\Cocc))\cap\mathfrak{D}=(C_z(X,\mathbb R)\oplus\mathcal S')\cap\mathfrak{D}=1.
\end{equation*}
This means that $\mathcal C=1$, as was to be shown.

\emph{3rd step.} We determine how the group $\Cob$ fits into the isomorphism (\ref{Eq:dir.sum.Coc.D}).

Recall that for every $i=1,\dots,n$, $\mathcal C_i^{d_i}=\co(\zeta_i)$. Consequently, it follows from (\ref{Eq:cob.split.m.n}) and from the first line in (\ref{Eq:GpExt.dcmp.tp.free}) that the inclusion $\Cob\subseteq\Cocc$ takes the form
\begin{equation*}
\begin{array}{ccccccccccc}
\Cob & \cong & 0 & \oplus & C_z(X,\mathbb R) & \oplus & \mathcal S' & \oplus & \bigoplus_{i=1}^n\langle\mathcal C_i^{d_i}\rangle & \subseteq & {} \\
{} & \subseteq & \mathfrak{D} & \oplus & C_z(X,\mathbb R) & \oplus & (\mathbb Q\otimes\mathcal S') & \oplus & \bigoplus_{i=1}^n\langle\mathcal C_i\rangle & \cong & \Cocc.
\end{array}
\end{equation*}
This leads to the isomorphisms
\begin{equation}\label{Eq:Cob.in.Cocc.iso}
\begin{array}{ccccccccccc}
\Cob & \cong & 0 & \oplus & C_z(X,\mathbb R) & \oplus & \mathbb Z^m & \oplus & \bigoplus_{i=1}^nd_i\mathbb Z & \subseteq & {} \\
{} & \subseteq & \mathfrak{D} & \oplus & C_z(X,\mathbb R) & \oplus & \mathbb Q^m & \oplus & \mathbb Z^n & \cong & \Cocc.
\end{array}
\end{equation}

\emph{4th step.} We finish the proof of the theorem.

First we verify the following statements:
\begin{enumerate}
\item[($\alpha$)] $\mathfrak{D}$ is a rational linear space with dimension $\QLSdim(\mathfrak{D})=\mathfrak{c}$, that is, $\mathfrak{D}\cong\mathbb Q^{(\mathfrak{c})}\cong\mathbb R$,
\item[($\beta$)] $C_z(X,\mathbb R)$ is isomorphic to $\mathbb R$ as an abelian group.
\end{enumerate}

Statement ($\alpha$) follows from the following facts:
\begin{itemize}
\item $\mathfrak{D}$ is a rational linear space; indeed, by definition, it is a divisible subgroup of a torsion-free abelian group $\Cocc$,
\item the rational linear dimension $\QLSdim(\mathfrak{D})$ of $\mathfrak{D}$ is at most $\mathfrak{c}$; this follows from the fact that $\Cocc$ is a Polish group with the topology of uniform convergence on compact sets and hence its cardinality is at most $\mathfrak{c}$,
\item the rational linear dimension of $\mathfrak{D}$ is at least $\mathfrak{c}$; this follows from the following two facts: 
\begin{itemize}
\item first, $\Cocm$ contains a subgroup isomorphic to $\mathbb R\cong\mathbb Q^{(\mathfrak{c})}$; this follows from Theorem~\ref{T:gp.min.free.ccl}, 
\item second, the projection $\Cocc\to\mathfrak{D}$ becomes a monomorphism when restricted to any subgroup of the groupoid $\Cocm$; this follows from the fact that the complementary summand of $\mathfrak{D}$ in $\Cocc$ contains no minimal extensions whatsoever.
\end{itemize}
\end{itemize}

Statement ($\beta$) follows from the following observations:
\begin{itemize}
\item $C_z(X,\mathbb R)$ is a real linear space,
\item the real linear dimension $\LSdim(C_z(X,\mathbb R))$ of $C_z(X,\mathbb R)$ is $\mathfrak{c}$; this follows by the same argument as in the proof of Theorem~\ref{T:F.in.Gim.CLAC},
\item every real linear space $L$ with $\LSdim(L)=\mathfrak{c}$ is isomorphic to $\mathbb R$ as an abelian group.
\end{itemize}

Now the isomorphisms in (\ref{Eq:Cob.in.Cocc}) follow immediately from (\ref{Eq:Cob.in.Cocc.iso}), ($\alpha$) and ($\beta$). The isomorphism (\ref{Eq:Cob.in.Cocc3}) is a direct consequence of (\ref{Eq:Cob.in.Cocc}) and, finally, (\ref{Eq:Cob.in.Cocc2}) follows from (\ref{Eq:Cob.in.Cocc}) similarly to Theorem~\ref{T:F.in.Gim.CLAC}.
\end{proof}

For the purpose of the following corollary we fix some notation and terminology. Let
\begin{equation*}
E_1\colon 0\longrightarrow A_1\stackrel{\mu_1}{\longrightarrow}B_1\stackrel{\nu_1}{\longrightarrow}C_1\longrightarrow0\hspace{4mm}\text{and}\hspace{4mm}
E_2\colon 0\longrightarrow A_2\stackrel{\mu_2}{\longrightarrow}B_2\stackrel{\nu_2}{\longrightarrow}C_2\longrightarrow0
\end{equation*}
be extensions of abelian groups. In Subsection~\ref{Sub:ext.ab.gps} we defined the equivalence of $E_1$ and $E_2$ in a strict sense which required, in particular, that $A_1=A_2$ and $C_1=C_2$. We shall now relax this condition by regarding isomorphic groups as the same. Thus, we call $E_1,E_2$ equivalent if there exist isomorphisms of groups $h\colon A_1\to A_2$, $k\colon B_1\to B_2$ and $l\colon C_1\to C_2$ with $k\mu_1=\mu_2 h$ and $l\nu_1=\nu_2 k$. Further, we call $2$-cocycles $\varphi_1\in Z^2(C_1,A_1)$ and $\varphi_2\in Z^2(C_2,A_2)$ equivalent if they give rise to equivalent extensions $E(\varphi_1)$, $E(\varphi_2)$. (Clearly, if two $2$-cocycles $\varphi_1,\varphi_2\in Z^2(C,A)$ are equivalent in a strict sense, that is, if they differ by a $2$-coboundary, then they are equivalent.) We use the same symbol $\cong$ for the equivalence of extensions and $2$-cocycles. If $A,C$ are abelian groups and one of them is trivial, then $\varphi_0(C,A)$ stands for the unique element of $Z^2(C,A)$ and $E_0(C,A)=E(\varphi_0(C,A))$ stands for the extension underlying to $\varphi_0(C,A)$. Finally, given a minimal flow $\Flow$, we denote by $E_{\mathcal F}$ the usual extension
\begin{equation*}
E_{\mathcal F}\colon 0\longrightarrow\Cob\stackrel{\mu}{\longrightarrow}\Cocc\stackrel{\pi}{\longrightarrow}\Coch\longrightarrow0
\end{equation*}
and by $\varphi_{\mathcal F}$ a fixed element of $Z^2(\Coch,\Cob)$ with $E(\varphi_{\mathcal F})\cong E_{\mathcal F}$.

\begin{corollary}\label{C:coc.for.GpExt}
Under the assumptions of Theorem~\ref{T:F.in.Gim.free.ccl} there are an equivalence of extensions
\begin{equation}\label{Eq:coc.for.GpExt1}
E_{\mathcal F}\cong E_0(\mathbb R,0)\oplus E_0(0,\mathbb R)\oplus\bigoplus_{j=1}^m E(\mathbb Q/\mathbb Z)\oplus\bigoplus_{i=1}^nE(\mathbb Z_{d_i})
\end{equation}
and an equivalence of $2$-cocycles
\begin{equation}\label{Eq:coc.for.GpExt2}
\varphi_{\mathcal F}\cong \varphi_0(\mathbb R,0)\oplus \varphi_0(0,\mathbb R)\oplus\bigoplus_{j=1}^m\psi
\oplus\bigoplus_{i=1}^n\varphi_{d_i},
\end{equation}
where $E(\mathbb Z_{d_i})$, $\varphi_{d_i}$, $E(\mathbb Q/\mathbb Z)$ and $\psi$ are given by \emph{(\ref{Eq:ext.for.Zk})}, \emph{(\ref{Eq:2-coc.for.Zk})}, \emph{(\ref{Eq:ext.for.Q/Z})} and \emph{(\ref{Eq:2-coc.for.Q/Z})}, respectively.
\end{corollary}
\begin{remark}\label{R:coc.for.GpExt}
Notice that under the assumptions of Theorem~\ref{T:F.in.Gim.free.ccl} we have $m\neq0$. This means that the extension $E_{\mathcal F}$ does not split or, in other words, that $\Cob$ is not a direct summand in $\Cocc$. It follows also that $\varphi_{\mathcal F}$ is not a $2$-coboundary.
\end{remark}
\begin{proof}[Proof of Corollary~\ref{C:coc.for.GpExt}]
We recall that the notion of a direct sum of extensions and $2$-cocycles has been defined in Subsection~\ref{Sub:dir.sum.seq.coc}. Now, by virtue of the isomorphisms (\ref{Eq:Cob.in.Cocc}), there are equivalences of extensions
\begin{equation*}
\begin{split}
E(\varphi_{\mathcal F})\cong E_{\mathcal F}&\cong E_0(\mathbb R,0)\oplus E_0(0,\mathbb R)\oplus\bigoplus_{j=1}^m E(\mathbb Q/\mathbb Z)\oplus\bigoplus_{i=1}^nE(\mathbb Z_{d_i})\\
&\cong E(\varphi_0(\mathbb R,0))\oplus E(\varphi_0(0,\mathbb R))\oplus\bigoplus_{j=1}^mE(\psi)\oplus\bigoplus_{i=1}^nE(\varphi_{d_i})\\
&\cong E\left(\varphi_0(\mathbb R,0)\oplus\varphi_0(0,\mathbb R)\oplus\bigoplus_{j=1}^m\psi\oplus\bigoplus_{i=1}^n\varphi_{d_i}\right).
\end{split}
\end{equation*}
These equivalences yield both (\ref{Eq:coc.for.GpExt1}) and (\ref{Eq:coc.for.GpExt2}).
\end{proof}

\begin{proposition}\label{P:Cocc.top.free.ccl}
Let $\Flow$ be a minimal topologically free flow with $\Gamma\in\mathsf{LieGp}$ connected and with $X$ a compact connected manifold. Set $n=\rank(H_1^w(\mathcal F))$. If $\Cocc$ is equipped with the topology of uniform convergence on compact sets then there is a topological direct sum
\begin{equation}\label{Eq:Div.in.Cocc.top}
\Cocc=\Div(\Cocc)\oplus\mathbb Z^n,
\end{equation}
where $\Div(\Cocc)$ is the additive topological group of a real separable Banach space and the group $\mathbb Z^n$ is discrete.
\end{proposition}
\begin{proof}
By virtue of (\ref{Eq:Cocc.div.sum}) from Theorem~\ref{T:ex.seq.Lie.mfld}, the direct sum (\ref{Eq:Div.in.Cocc.top}) holds in the algebraic sense. Let the group $\Cocc(\mathbb R)$ carry the topology of uniform convergence on compact sets. Since the group $\Gamma$ is compactly generated and both $\Gamma$, $X$ are locally compact second countable, $\Cocc(\mathbb R)$ is the additive topological group of a real separable Banach space with the norm given by (\ref{Eq:norm.ext}), see Lemma~\ref{L:Cocc.Ban.sp.Gmm} and Remark~\ref{R:Cocc.Ban.sp.Gmm} in Subsection~\ref{Sub:Ban.sp.of.ext}. Let $p$ be the standard covering morphism $p\colon\mathbb R\to\mathbb T^1$ and consider the induced morphism
\begin{equation*}
\widehat{p}\colon\Cocc(\mathbb R)\ni\mathcal D\mapsto p\mathcal D\in\Cocc.
\end{equation*}
Then $\widehat{p}$ is a morphism of topological groups. Moreover, by connectedness of both $\Gamma$ and $X$, $\widehat{p}$ is a monomorphism. We claim that $\widehat{p}$ is in fact a topological isomorphism onto its image. So let $(\mathcal D_n)_{n\in\mathbb N}$ be a sequence in $\Cocc(\mathbb R)$ with $p\mathcal D_n\to 1$ in $\Cocc$; we show that $\mathcal D_n\to0$ in $\Cocc(\mathbb R)$. To this end, fix a compact set $C\subseteq\Gamma\times X$ and $0<\varepsilon<1/2$. Let $1\in K\subseteq\Gamma$ be a compact connected set with $C\subseteq K\times X$. Then for all but finitely many $n\in\mathbb N$, $p\mathcal D_n(K\times X)\subseteq p(-\varepsilon,\varepsilon)$. Since the sets $\mathcal D_n(K\times X)\subseteq\mathbb R$ are all connected and contain $0$, it follows that the inclusions $\mathcal D_n(C)\subseteq\mathcal D_n(K\times X)\subseteq(-\varepsilon,\varepsilon)$ also hold for all but finitely many $n\in\mathbb N$. Thus, $\mathcal D_n\to0$ in $\Cocc(\mathbb R)$, as was to be shown.

Now, given $\mathcal C\in\Cocc$, the following conditions are equivalent by Lemma~\ref{P:lifting.cocycle}:
\begin{itemize}
\item $\mathcal C\in\Div(\Cocc)$,
\item for every $k\in\mathbb N$, $\mathcal C$ lifts across $\kappa_k\colon\mathbb T^1\to\mathbb T^1$ to an element of $\Cocc$,
\item the induced morphism $\mathcal C^{\sharp}\colon\pi_1(\Gamma\times X)\to\pi_1(\mathbb T^1)$ vanishes,
\item $\mathcal C$ lifts across $p$ to and element of $\Cocc(\mathbb R)$.
\end{itemize}
This means that $\im(\widehat{p})=\Div(\Cocc)$ and, since $\Cocc(\mathbb R)$ is the additive topological group of a real separable Banach space, it follows that so is the group $\Div(\Cocc)$.

In order to finish the proof, we must show that the complementary summand $\mathbb Z^n$ of $\Div(\Cocc)$ in $\Cocc$ is a discrete subgroup of $\Cocc$ and that the splitting (\ref{Eq:Div.in.Cocc.top}) is topological. We verify both these statements at once by showing that $\Div(\Cocc)$ is an open subgroup of $\Cocc$. Clearly, it is sufficient to find a neighbourhood $\mathcal U$ of $1$ in $\Cocc$ with $\mathcal U\subseteq\Div(\Cocc)$. Since $X$ is a compact connected manifold and $\Gamma$ is a connected Lie group, both $\pi_1(X)$ and $\pi_1(\Gamma)$ are finitely generated. Consequently, the group $\pi_1(\Gamma\times X)$ is also finitely generated, say, with generators $f_1,\dots,f_l$. The union of the images of $[0,1]$ under the paths $f_1,\dots,f_l$ is a compact subset $K$ of $\Gamma\times X$ and hence $\mathcal U=[K;\mathbb T^1\setminus-1]=\{\mathcal C\in\Cocc : \mathcal C(K)\subseteq \mathbb T^1\setminus-1\}$ is a neighbourhood of $1$ in $\Cocc$. If $\mathcal C\in\Cocc$ is an element of $\mathcal U$ then $\mathcal C^{\sharp}(f_i)=0$ for every $i=1,\dots,l$, since $\mathbb T^1\setminus-1$ is a contractible space, and hence $\mathcal C^{\sharp}=0$. Thus, $\mathcal U\subseteq\im(\widehat{p})=\Div(\Cocc)$, as was to be shown.
\end{proof}

\section{First cohomology groups with arbitrary coefficients}\label{S:gen.fib}

Let $\Flow$ be a minimal flow and $G\in\mathsf{CAGp}$ be a connected group. In this section we continue our investigation from Section~\ref{S:circ.case}. Our aim is to express the inclusions $\Cob(G)\subseteq\Cocc(G)$ and $\Coch(G)\subseteq\Hom(G^*,\Coch)$ in terms of inclusions of familiar (topological) abelian groups and to compute the first cohomology group $\Coch(G)$ of $\mathcal F$ with coefficients in $G^*$. We concentrate on the two situations (i) and (ii), which have been described at the beginning of the preceding section, see Theorems~\ref{T:Co.in.Cocc.G} and~\ref{C:GpExt.frccl.gen.G}, respectively. We also recall Corollary~\ref{C:structure.coc.simply}, which states that, in case (i), the group $\Cocc(G)$ is the additive topological group of a real separable Fr\'echet or Banach space. In Proposition~\ref{P:Cocc.top.free.ccl.G} we show that, in case (ii), the structure of a Fr\'echet/Banach space is carried by the divisible subgroup $\Div(\Cocc(G))$ of $\Cocc(G)$, which forms a topological direct summand in $\Cocc(G)$.

Before proceeding to the next theorem let us recall from Theorem~\ref{T:intro.embed.G(F)} that for every $G\in\mathsf{CAGp}$, the morphism $\Psi_G\colon\Coch(G)\to\Hom(G^*,\Coch)$ defined by (\ref{Eq:ext.hom.G*}) is a topological isomorphism onto its image, provided $\Coch(G)$, $\Coch$ are equipped with the ext-topology and $\Hom(G^*,\Coch)$ carries the topology of point-wise convergence.

\begin{theorem}\label{C:GpExt.frccl.gen.G}
Let $G\in\mathsf{CAGp}$ be second countable and connected. Write
\begin{equation*}
\mathfrak{r}=\rank(G^*)\,,\hspace{3mm}\mathfrak{f}=\rank(\mathfrak{f}(G^*))\hspace{3mm}
\text{and}\hspace{3mm}\mathfrak{t}=\rank(\mathfrak{t}(G^*)).
\end{equation*}
Then, under the assumptions and notation from Theorem~\ref{T:F.in.Gim.free.ccl}, the inclusion $\Cob(G)\subseteq\Cocc(G)$ takes the form
\begin{equation}\label{Eq:CobG.in.CoccG}
\begin{array}{cccccccccccc}
{} & 0 & \oplus & \mathbb R^{\mathfrak{r}} & \oplus & \mathbb Z^{m\mathfrak{f}} & \oplus & 0 & \oplus & (d_1\mathbb Z\oplus\dots\oplus d_n\mathbb Z)^{\mathfrak{f}} & \subseteq \\
\subseteq & \mathbb R^\mathfrak{r} & \oplus & \mathbb R^{\mathfrak{r}} & \oplus & \mathbb Q^{m\mathfrak{f}} & \oplus & \mathbb Q^{m\mathfrak{t}} & \oplus & \left(\mathbb Z^n\right)^{\mathfrak{f}}, & {}
\end{array}
\end{equation}
and there is an isomorphism of groups
\begin{eqnarray}\label{Eq:Cocch.gen.G.frccl}
\Coch(G)\cong\mathbb R^{\mathfrak{r}}\oplus(\mathbb Q/\mathbb Z)^{m\mathfrak{f}}\oplus \mathbb Q^{m\mathfrak{t}}\oplus\left(\mathbb Z_{d_1}\oplus\dots\oplus\mathbb Z_{d_n}\right)^{\mathfrak{f}}.
\end{eqnarray}
Moreover, the topological morphism $\Psi_G\colon\Coch(G)\to\Hom(G^*,\Coch)$ takes the form of the inclusion of topological groups
\begin{equation*}
\begin{array}{cccccc}
\mathbb R^{\mathfrak{r}} & \oplus & (\mathbb Q/\mathbb Z)^{m\mathfrak{f}} & \oplus & \Hom\left(\mathfrak{t}(G^*),\mathbb Q^m\right) & \oplus \\ 
\left(\bigoplus_{i=1}^n\mathbb Z_{d_i}\right)^{\mathfrak{f}} & \oplus & 0 & \subseteq & {} & \\
\mathbb R^{\mathfrak{r}} &  \oplus & (\mathbb Q/\mathbb Z)^{m\mathfrak{f}} & \oplus & \Hom\left(\mathfrak{t}(G^*),(\mathbb Q/\mathbb Z)^m\right) & \oplus \\
 \left(\bigoplus_{i=1}^n\mathbb Z_{d_i}\right)^{\mathfrak{f}} & \oplus & \bigoplus_{i=1}^n\tor_{d_i}\left(\mathfrak{t}(G)\right), & {} & {} & 
\end{array}
\end{equation*}
where the direct sums and the powers are topological, the groups $\mathbb R$, $\mathbb Q$ and $\mathbb Q/\mathbb Z$ are equipped with the discrete topology, the group $\Hom(\mathfrak{t}(G^*),(\mathbb Q/\mathbb Z)^m)$ is equipped with the topology of point-wise convergence and the groups $\tor_{d_i}(\mathfrak{t}(G))$ ($i=1,\dots,n$) carry the topology inherited from $\mathfrak{t}(G)$.
\end{theorem}
\begin{remark}\label{R:GpExt.frccl.gen.G}
Notice the following facts.
\begin{itemize}
\item Let us clarify the last statement of the theorem. While $\Hom(\mathfrak{t}(G^*),(\mathbb Q/\mathbb Z)^m)$ carries the topology of point-wise convergence, the group $\Hom(\mathfrak{t}(G^*),\mathbb Q^m)$ is equipped with the subgroup topology inherited from $\Hom(\mathfrak{t}(G^*),(\mathbb Q/\mathbb Z)^m)$ via the monomorphism 
\begin{equation*}
\Hom(\mathfrak{t}(G^*),q)\colon\Hom(\mathfrak{t}(G^*),\mathbb Q^m)\ni h\mapsto qh\in\Hom(\mathfrak{t}(G^*),(\mathbb Q/\mathbb Z)^m),
\end{equation*}
where $q\colon\mathbb Q^m\to(\mathbb Q/\mathbb Z)^m$ is the canonical quotient morphism. (This is indeed a monomorphism, since the group $\ker(q)=\mathbb Z^m$ is free abelian and the group $\mathfrak{t}(G^*)$ is torsion-less.) We emphasize that this subgroup topology on $\Hom(\mathfrak{t}(G^*),\mathbb Q^m)$ is in general different from the topology of point-wise convergence, see our remark following the sixth step of the proof. We also wish to mention the isomorphism (\ref{Eq:Hom(T,Q/Zresp.Q)}) from the eighth step of the proof, which states that $\Hom(\mathfrak{t}(G^*),\mathbb Q^m)$ is algebraically a direct summand in $\Hom(\mathfrak{t}(G^*),(\mathbb Q/\mathbb Z)^m)$ with the complementary summand $\Ext(\mathfrak{t}(G^*),\mathbb Z^m)$.
\item Assume that the group $G$ is non-trivial. Then $1\leq\mathfrak{r}\leq\aleph_0$ and we get isomorphisms of abelian groups
\begin{equation*}
\mathbb R^{\mathfrak{r}}\cong\left(\mathbb Q^{(\mathfrak{c})}\right)^{\mathfrak{r}}\cong\mathbb Q^{\left(\mathfrak{c}^{\mathfrak{r}}\right)}\cong\mathbb Q^{(\mathfrak{c})}\cong\mathbb R,
\end{equation*}
which apply to the isomorphisms (\ref{Eq:CobG.in.CoccG}) and (\ref{Eq:Cocch.gen.G.frccl}). (They do not apply to the last isomorphism of the theorem, since they are not topological in general.)
\item It follows from the last statement of the theorem that $\Psi_G$ is an isomorphism if $\mathfrak{t}(G^*)=0$, that is, if $G$ is a torus. The converse implication is also true. If $\Psi_G$ is an isomorphism then, by the first part of the remark, $\Ext(\mathfrak{t}(G^*),\mathbb Z^m)=0$. Thus $\mathfrak{t}(G^*)$ is a Whitehead group and, being countable, it is free abelian by the Pontryagin theorem. Since $\mathfrak{t}(G^*)$ is also torsion-less by definition, it follows that, indeed, $\mathfrak{t}(G^*)=0$.
\end{itemize}
\end{remark}
\begin{proof}[Proof of Theorem~\ref{C:GpExt.frccl.gen.G}]
While the isomorphisms (\ref{Eq:CobG.in.CoccG}) and (\ref{Eq:Cocch.gen.G.frccl}) are meant only in algebraic sense, the last isomorphism describes an inclusion of topological groups. (The groups $\Coch$, $\Coch(G)$ carry the ext-topology and the group $\Hom(G^*,\Coch)$ is equipped with the topology of point-wise convergence.) When dealing with topological isomorphisms, in particular, when proving this last isomorphism of the theorem, we shall follow the following rules:
\begin{itemize}
\item all $\Hom$ groups carry the topology of point-wise convergence,
\item all direct sums are topological,
\item all powers are topological,
\item the groups $\mathbb R$, $\mathbb Q$, $\mathbb Q/\mathbb Z$ and $\mathbb Z$ carry the discrete topology.
\end{itemize}

We divide the proof into several short steps.

\emph{1st step.} We recall some properties of the functor $\Hom(G^*,-)$, which we shall use later.

For $k=1,\dots,n$ let $A_k\in\mathsf{DAGp}$ and $B_k\sbgp A_k$. Set $A=\bigoplus_{k=1}^nA_k$ and $B=\bigoplus_{k=1}^nB_k$. Then $B$ is a subgroup of $A$, $\Hom(G^*,B)$ is a closed subgroup of $\Hom(G^*,A)$ and there are topological isomorphisms
\begin{equation*}
\Hom(G^*,B)\cong\bigoplus_{k=1}^n\Hom(G^*,B_k)\subseteq\bigoplus_{k=1}^n\Hom(G^*,A_k)
\cong\Hom(G^*,A).
\end{equation*}
Consequently, there is a topological isomorphism
\begin{equation*}
\Hom(G^*,A)/\Hom(G^*,B)\cong\bigoplus_{k=1}^n\Big(\Hom(G^*,A_k)/
\Hom(G^*,B_k)\Big).
\end{equation*}

\emph{2nd step.} We show that there is an isomorphism of groups $\Hom(G^*,\mathbb R)\cong\mathbb R^{\mathfrak{r}}$.

Since the group $G^*$ is torsion-free and the group $\mathbb R$ is divisible, every element $h$ of $\Hom(G^*,\mathbb R)$ extends uniquely to an element $h'$ of $\Hom(\mathbb Q\otimes G^*,\mathbb R)$. Moreover, $\mathbb Q\otimes G^*$ is a rational linear space with dimension $\mathfrak{r}$ and so it is isomorphic to $\mathbb Q^{(\mathfrak{r})}$. This leads to isomorphisms of groups
\begin{equation}\label{Eq:QLS.Hom.iso}
\Hom\left(G^*,\mathbb R\right)\cong\Hom\left(\mathbb Q\otimes G^*,\mathbb R\right)\cong\Hom\left(\mathbb Q^{(\mathfrak{r})},\mathbb R\right)\cong\Hom(\mathbb Q,\mathbb R)^{\mathfrak{r}}\cong\mathbb R^{\mathfrak{r}}.
\end{equation}
One verifies that all the isomorphisms in (\ref{Eq:QLS.Hom.iso}) are in fact topological and so the isomorphism $\Hom(G^*,\mathbb R)\cong\mathbb R^{\mathfrak{r}}$ is also topological.

\emph{3rd step.} We describe the inclusion $\Hom\left(G^*,\mathbb Z^m\right)\subseteq\Hom\left(G^*,\mathbb Q^m\right)$ up to isomorphism.

Recall the direct sum $G^*=\mathfrak{f}(G^*)\oplus\mathfrak{t}(G^*)$. The group $\mathfrak{f}(G^*)$ is free abelian with rank $\mathfrak{f}$ and so it is isomorphic to $\mathbb Z^{(\mathfrak{f})}$. The group $\mathfrak{t}(G^*)$ is torsion-less and hence it satisfies $\Hom(\mathfrak{t}(G^*),\mathbb Z^m)=0$. Consequently, there are isomorphisms of groups
\begin{equation}\label{Eq:Hom.A.Zm}
\begin{split}
\Hom(G^*,\mathbb Z^m)&=\Hom(\mathfrak{f}(G^*)\oplus \mathfrak{t}(G^*),\mathbb Z^m)\\
&\cong\Hom(\mathfrak{f}(G^*),\mathbb Z^m)\oplus\Hom(\mathfrak{t}(G^*),\mathbb Z^m)\\
&\cong\Hom\left(\mathbb Z^{(\mathfrak{f})},\mathbb Z^m\right)\cong\left(\mathbb Z^m\right)^{\mathfrak{f}}\cong\mathbb Z^{m\mathfrak{f}}.
\end{split}
\end{equation}

In a similar way (and with the help of the ideas used to verify (\ref{Eq:QLS.Hom.iso})) we obtain isomorphisms of groups
\begin{equation}\label{Eq:Hom.A.Zm2}
\begin{split}
\Hom\left(G^*,\mathbb Q^m\right)&=\Hom\left(\mathfrak{f}(G^*)\oplus \mathfrak{t}(G^*),\mathbb Q^m\right)\\
&\cong\Hom\left(\mathfrak{f}(G^*),\mathbb Q^m\right)\oplus\Hom\left(\mathfrak{t}(G^*),\mathbb Q^m\right)\\
&\cong\Hom\left(\mathfrak{f}(G^*),\mathbb Q^m\right)\oplus\Hom\left(\mathbb Q\otimes \mathfrak{t}(G^*),\mathbb Q^m\right)\\
&\cong\Hom\left(\mathbb Z^{(\mathfrak{f})},\mathbb Q^m\right)\oplus\Hom\left(\mathbb Q^{(\mathfrak{t})},\mathbb Q^m\right)\\
&\cong\left(\mathbb Q^m\right)^{\mathfrak{f}}\oplus\left(\mathbb Q^m\right)^{\mathfrak{t}}\cong\mathbb Q^{m\mathfrak{f}}\oplus\mathbb Q^{m\mathfrak{t}}.
\end{split}
\end{equation}

The isomorphisms (\ref{Eq:Hom.A.Zm}) and (\ref{Eq:Hom.A.Zm2}) now yield the following form of the inclusion $\Hom(G^*,\mathbb Z^m)\subseteq\Hom(G^*,\mathbb Q^m)$:
\begin{equation*}
\Hom\left(G^*,\mathbb Z^m\right)\cong\mathbb Z^{m\mathfrak{f}}\oplus 0\subseteq\mathbb Q^{m\mathfrak{f}}\oplus\mathbb Q^{m\mathfrak{t}}\cong\Hom\left(G^*,\mathbb Q^m\right).
\end{equation*}

\emph{4th step.} We describe the inclusion $\Hom\left(G^*,\bigoplus_{i=1}^nd_i\mathbb Z\right)\subseteq\Hom(G^*,\mathbb Z^n)$ up to isomorphism.

By using the same ideas as in (\ref{Eq:Hom.A.Zm}), one puts the considered inclusion at once into the following form:
\begin{equation*}
\Hom\left(G^*,\bigoplus_{i=1}^nd_i\mathbb Z\right)\cong\left(\bigoplus_{i=1}^nd_i\mathbb Z\right)^{\mathfrak{f}}\subseteq\left(\mathbb Z^n\right)^{\mathfrak{f}}\cong\Hom\left(G^*,Z^n\right).
\end{equation*}

\emph{5th step.} We verify the isomorphisms (\ref{Eq:CobG.in.CoccG}) and (\ref{Eq:Cocch.gen.G.frccl}).

By virtue of Corollary~\ref{C:min.contr.char}, the group $\Cob(G)$ corresponds to $\Hom(G^*,\Cob)$ under the isomorphism $\Phi_G\colon\Cocc(G)\to\Hom(G^*,\Cocc)$. Consequently, the inclusion (\ref{Eq:CobG.in.CoccG}) follows at once from the first four steps of the proof and from the inclusion (\ref{Eq:Cob.in.Cocc}) in Theorem~\ref{T:F.in.Gim.free.ccl}. The isomorphism (\ref{Eq:Cocch.gen.G.frccl}) is a direct consequence of (\ref{Eq:CobG.in.CoccG}).

\emph{6th step.} We describe our strategy for proving the last statement of the theorem.

Under the usual identifications 
\begin{equation*}
\Cocc(G)\cong\Hom(G^*,\Cocc)\hspace{3mm}\text{and}\hspace{3mm}\Cob(G)\cong\Hom(G^*,\Cob),
\end{equation*}
the morphism $\Psi_G\colon\Coch(G)\to\Hom(G^*,\Coch)$ acts by the rule
\begin{equation}\label{Eq:the.form.PsiG}
\Psi_G\colon \frac{\Hom(G^*,\Cocc)}{\Hom(G^*,\Cob)}\ni\varphi+\Hom(G^*,\Cob)\mapsto\pi\varphi\in\Hom(G^*,\Coch),
\end{equation}
where $\pi$ is the quotient morphism $\Cocc\to\Coch$ and $\varphi\in\Hom(G^*,\Cocc)$. Recall that the group $\Coch$ carries the discrete topology, the group $\Hom(G^*,\Coch)$ is equipped with the topology of point-wise convergence and $\Psi_G$ is a topological isomorphism onto its image by virtue of Theorem~\ref{T:intro.embed.G(F)}.

Our aim now is to determine what form the inclusion (\ref{Eq:the.form.PsiG}) takes, when the inclusion $\Cob\subseteq\Cocc$ is expressed in the form of (\ref{Eq:Cob.in.Cocc}) and, correspondingly, the group $\Coch$ takes the form expressed in (\ref{Eq:Cob.in.Cocc3}). Write $A_1=A_2=\mathbb R$, $A_3=\mathbb Q^m$, $A_4=\mathbb Z^n$, $B_1=0$, $B_2=\mathbb R$, $B_3=\mathbb Z^m$ and $B_4=\bigoplus_{i=1}^nd_i\mathbb Z$. For $1\leq i\leq 4$ let $\pi_i$ be the quotient morphism $A_i\to A_i/B_i$. Then by virtue of Theorem~\ref{T:F.in.Gim.free.ccl}, the morphism $\pi$ takes the form of the direct sum of the morphisms $\pi_i$:
\begin{equation*}
\pi=\bigoplus_{i=1}^4\pi_i\colon\bigoplus_{i=1}^4A_i\to\bigoplus_{i=1}^4A_i/B_i.
\end{equation*}
Consequently, the expression of the morphism $\Psi_G$ from (\ref{Eq:the.form.PsiG}) changes to
\begin{equation*}
\Psi_G\colon\frac{\bigoplus_{i=1}^4\Hom(G^*,A_i)}{\bigoplus_{i=1}^4\Hom(G^*,B_i)}\cong\bigoplus_{i=1}^4\frac{\Hom(G^*,A_i)}{\Hom(G^*,B_i)}\to\bigoplus_{i=1}^4\Hom(G^*,A_i/B_i)
\end{equation*}
and $\Psi_G$ now acts by the rule
\begin{equation*}
\Psi_G\colon\bigoplus_{i=1}^4\left(h_i+\Hom(G^*,B_i)\right)\mapsto
\bigoplus_{i=1}^4\pi_ih_i
\end{equation*}
for $h_i\in\Hom(G^*,A_i)$ ($1\leq i\leq 4$). Thus $\Psi_G$ is a direct sum of the morphisms
\begin{equation*}
q_i\colon\frac{\Hom(G^*,A_i)}{\Hom(G^*,B_i)}\ni h_i+\Hom(G^*,B_i)\mapsto\pi_ih_i\in\Hom(G^*,A_i/B_i)
\end{equation*}
for $1\leq i\leq 4$ and we need to determine them all. Observe that all the morphisms $q_i$ ($1\leq i\leq 4$) are monomorphisms.

\emph{Remark.} Before proceeding further we wish to make some observations. First we recall that all the groups $A_i$ and $B_i$ carry the discrete topology and the quotient groups $A_i/B_i$ are thus also discrete. Further, the groups $\Hom(G^*,A_i/B_i)$ are equipped with the to\-po\-lo\-gy of point-wise convergence. Finally, since $\Psi_G$ is a topological isomorphism onto its image, we must consider on the quotient group $\Hom(G^*,A_i)/\Hom(G^*,B_i)$ the subgroup topology inherited from $\Hom(G^*,A_i/B_i)$ via the monomorphism $q_i$. We would like to emphasize that this inherited topology need not coincide with the quotient topology on $\Hom(G^*,A_i)/\Hom(G^*,B_i)$. Indeed, as an example consider the morphism
\begin{equation*}
q\colon\Hom(\mathbb Q,\mathbb Q)/\Hom(\mathbb Q,\mathbb Z)=\Hom(\mathbb Q,\mathbb Q)\ni h\mapsto ph\in\Hom(\mathbb Q,\mathbb Q/\mathbb Z),
\end{equation*}
where $p$ is the quotient morphism $\mathbb Q\to\mathbb Q/\mathbb Z$. Then $q$ is not a topological isomorphism onto its image if $\Hom(\mathbb Q,\mathbb Q)$ is equipped with the topology of point-wise convergence. For if $h_n\in\Hom(\mathbb Q,\mathbb Q)$ is given by the rule $h_n(x)=(n!)x$ for every $n\in\mathbb N$, then the sequence $(q(h_n))_{n\in\mathbb N}$ converges to zero in $\Hom(\mathbb Q,\mathbb Q/\mathbb Z)$ but the sequence $(h_n)_{n\in\mathbb N}$ does not converge to zero in $\Hom(\mathbb Q,\mathbb Q)$.

\emph{7th step.} We determine the morphisms $q_1$ and $q_2$.

Since $A_1=\mathbb R$ and $B_1=0$, the morphism $q_1$ turns at once into the identity on $\Hom\left(G^*,\mathbb R\right)$. Moreover, in the second step of the proof we have shown that the group $\Hom(G^*,\mathbb R)$ is topologically isomorphic to $\mathbb R^{\mathfrak{r}}$. Thus, we have an equivalence
\begin{equation}\label{Eq:q1}
q_1\cong \Id\colon\mathbb R^{\mathfrak{r}}\to\mathbb R^{\mathfrak{r}}.
\end{equation}
Also, since $A_2=B_2=\mathbb R$, $q_2$ is equivalent to the unique endomorphism of the trivial group $0$.

\emph{8th step.} We determine the morphism $q_3$.

Here $A_3=\mathbb Q^m$ and $B_3=\mathbb Z^m$. From the direct sum $G^*=\mathfrak{f}(G^*)\oplus \mathfrak{t}(G^*)$ it follows that the group $\Hom(G^*,A_3/B_3)$ is topologically isomorphic to the topological direct sum $\Hom(\mathfrak{f}(G^*),A_3/B_3)\oplus\Hom(\mathfrak{t}(G^*),A_3/B_3)$ and the morphism $q_3$ takes the form of a direct sum $q_3\cong q_3^{(f)}\oplus q_3^{(t)}$ of the morphisms
\begin{equation*}
q_3^{(f)}\colon\frac{\Hom\left(\mathfrak{f}(G^*),\mathbb Q^m\right)}{\Hom\left(\mathfrak{f}(G^*),\mathbb Z^m\right)}\ni h+\Hom(\mathfrak{f}(G^*),\mathbb Z^m)\mapsto\pi_3h\in\Hom\left(\mathfrak{f}(G^*),(\mathbb Q/\mathbb Z)^m\right)
\end{equation*}
and
\begin{equation*}
q_3^{(t)}\colon\frac{\Hom\left(\mathfrak{t}(G^*),\mathbb Q^m\right)}{\Hom\left(\mathfrak{t}(G^*),\mathbb Z^m\right)}\ni k+\Hom(\mathfrak{t}(G^*),\mathbb Z^m)\mapsto\pi_3k\in\Hom\left(\mathfrak{t}(G^*),(\mathbb Q/\mathbb Z)^m\right).
\end{equation*}

First we determine the morphism $q_3^{(f)}$. Similarly to the third step of the proof one verifies that the group $\Hom(\mathfrak{f}(G^*),(\mathbb Q/\mathbb Z)^m)$ is topologically isomorphic to $\left((\mathbb Q/\mathbb Z)^m\right)^{\mathfrak{f}}\cong(\mathbb Q/\mathbb Z)^{m\mathfrak{f}}$. Moreover, since $\mathfrak{f}(G^*)$ is a free abelian group, it satisfies $\Ext(\mathfrak{f}(G^*),\mathbb Z^m)=0$. Consequently, the covariant Hom-Ext sequence derived from the sequence
\begin{equation}\label{Eq:short.seq.Zm.Qm}
0\longrightarrow\mathbb Z^m\longrightarrow\mathbb Q^m\longrightarrow(\mathbb Q/\mathbb Z)^m\longrightarrow0
\end{equation}
and associated to the group $\mathfrak{f}(G^*)$ takes the form of the short exact sequence
\begin{equation*}
0\longrightarrow\Hom(\mathfrak{f}(G^*),\mathbb Z^m)\longrightarrow\Hom(\mathfrak{f}(G^*),\mathbb Q^m)\longrightarrow\Hom(\mathfrak{f}(G^*),(\mathbb Q/\mathbb Z)^m)\longrightarrow0.
\end{equation*}
This shows that $\Hom(\mathfrak{f}(G^*),\pi_3)$ is an epimorphism with the kernel $\Hom(\mathfrak{f}(G^*),\mathbb Z^m)$ and so $q_3^{(f)}$ is an isomorphism. Thus, as a topological morphism, $q_3^{(f)}$ is equivalent to the identity on $(\mathbb Q/\mathbb Z)^{m\mathfrak{f}}$.

Now we determine the morphism $q_3^{(t)}$. First, since the group $\mathfrak{t}(G^*)$ is torsion-less, we have $\Hom\left(\mathfrak{t}(G^*),\mathbb Z^m\right)=0$, which means that $q_3^{(t)}\cong\Hom(\mathfrak{t}(G^*),\pi_3)$ represents the natural inclusion $\Hom(\mathfrak{t}(G^*),\mathbb Q^m)\subseteq\Hom(\mathfrak{t}(G^*),(\mathbb Q/\mathbb Z)^m)$. Further, since $\mathbb Q^m$ is a divisible group, the group of extensions $\Ext(\mathfrak{t}(G^*),\mathbb Q^m)$ vanishes. Consequently, the covariant Hom-Ext sequence derived from (\ref{Eq:short.seq.Zm.Qm}) and associated to the group $\mathfrak{t}(G^*)$ takes the form of the short exact sequence
\begin{equation*}
0\longrightarrow\Hom(\mathfrak{t}(G^*),\mathbb Q^m)\stackrel{q_3^{(t)}}{\longrightarrow}\Hom(\mathfrak{t}(G^*),(\mathbb Q/\mathbb Z)^m)\longrightarrow\Ext(\mathfrak{t}(G^*),\mathbb Z^m)\longrightarrow0.
\end{equation*}
It follows that the group $\Hom(\mathfrak{t}(G^*),(\mathbb Q/\mathbb Z)^m)$ is an extension of $\Ext(\mathfrak{t}(G^*),\mathbb Z^m)$ by $\im(q_3^{(t)})=\Hom(\mathfrak{t}(G^*),\mathbb Q^m)$. Since the group $\Hom(\mathfrak{t}(G^*),\mathbb Q^m)\cong\mathbb Q^{m\mathfrak{t}}$ is divisible and hence injective, it is a direct summand in $\Hom(\mathfrak{t}(G^*),(\mathbb Q/\mathbb Z)^m)$ and we thus obtain the following isomorphisms of abelian groups:
\begin{equation}\label{Eq:Hom(T,Q/Zresp.Q)}
\begin{split}
\Hom(\mathfrak{t}(G^*),(\mathbb Q/\mathbb Z)^m)&=\Hom(\mathfrak{t}(G^*),\mathbb Q^m)\oplus\Ext(\mathfrak{t}(G^*),\mathbb Z^m)\\
&\cong\mathbb Q^{m\mathfrak{t}}\oplus\Ext(\mathfrak{t}(G^*),\mathbb Z^m).
\end{split}
\end{equation}

\emph{9th step.} We determine the morphism $q_4$.

Here we have $A_4=\mathbb Z^n$ and $B_4=\bigoplus_{i=1}^nd_i\mathbb Z$. Similarly to the eighth step of the proof one verifies that the group $\Hom(G^*,A_4/B_4)$ is topologically isomorphic to the topological direct sum $\Hom(\mathfrak{f}(G^*),A_4/B_4)\oplus\Hom(\mathfrak{t}(G^*),A_4/B_4)$ and that the morphism $q_4$ takes the form of a direct sum $q_4\cong q_4^{(f)}\oplus q_4^{(t)}$ of the natural morphisms
\begin{equation*}
q_4^{(f)}\colon\Hom(\mathfrak{f}(G^*),\mathbb Z^n)/\Hom\left(\mathfrak{f}(G^*),\bigoplus_{i=1}^nd_i\mathbb Z\right)\to\Hom\left(\mathfrak{f}(G^*),\bigoplus_{i=1}^n\mathbb Z_{d_i}\right)
\end{equation*}
and
\begin{equation*}
q_4^{(t)}\colon\Hom(\mathfrak{t}(G^*),\mathbb Z^n)/\Hom\left(\mathfrak{t}(G^*),\bigoplus_{i=1}^nd_i\mathbb Z\right)\to\Hom\left(\mathfrak{t}(G^*),\bigoplus_{i=1}^n\mathbb Z_{d_i}\right).
\end{equation*}

In order to determine the morphism $q_4^{(f)}$, one can use the same ideas as those used to determine $q_3^{(f)}$ in the eighth step of the proof. Indeed, since the group $\mathfrak{f}(G^*)$ is free abelian, we have $\Ext(\mathfrak{f}(G^*),\bigoplus_{i=1}^nd_i\mathbb Z)=0$. It follows that $\Hom(\mathfrak{f}(G^*),\pi_4)$ is an epimorphism onto $\Hom(\mathfrak{f}(G^*),\bigoplus_{i=1}^n\mathbb Z_{d_i})$ with the kernel $\Hom(\mathfrak{f}(G^*),\bigoplus_{i=1}^nd_i\mathbb Z)$ and hence $q_4^{(f)}$ is an isomorphism. Thus, as a morphism of topological groups, $q_4^{(f)}$ is equivalent to the identity on
\begin{equation*}
\Hom\left(\mathfrak{f}(G^*),\bigoplus_{i=1}^n\mathbb Z_{d_i}\right)\cong\Hom\left(\mathbb Z^{(\mathfrak{f})},\bigoplus_{i=1}^n\mathbb Z_{d_i}\right)\cong\left(\bigoplus_{i=1}^n\mathbb Z_{d_i}\right)^{\mathfrak{f}}.
\end{equation*}

Now we determine the morphism $q_4^{(t)}$. The domain of $q_4^{(t)}$ is the trivial group $0$, since $\mathfrak{t}(G^*)$ is a torsion-less group. The co-domain $\Hom(\mathfrak{t}(G^*),\bigoplus_{i=1}^n\mathbb Z_{d_i})$ of $q_4^{(t)}$ is obtained from the following topological isomorphisms:
\begin{equation*}
\begin{split}
\Hom\left(\mathfrak{t}(G^*),\bigoplus_{i=1}^n\mathbb Z_{d_i}\right)&\cong\bigoplus_{i=1}^n\Hom\left(\mathfrak{t}(G^*),\mathbb Z_{d_i}\right)\cong\bigoplus_{i=1}^n\Hom\left(\mathbb Z_{d_i},\mathfrak{t}(G^*)^*\right)\\
&\cong\bigoplus_{i=1}^n\text{$\tor$}_{d_i}(\mathfrak{t}(G^*)^*)\cong\bigoplus_{i=1}^n\text{tor}_{d_i}(\mathfrak{t}(G)).
\end{split}
\end{equation*}

\emph{10th step.} We finish the proof of the theorem by verifying its last statement.

The last statement of the theorem follows at once from the previous steps $6-9$.
\end{proof}

\begin{proposition}\label{P:Cocc.top.free.ccl.G}
Let $\Flow$ be a minimal topologically free flow with $\Gamma\in\mathsf{LieGp}$ connected and with $X$ a compact manifold, and let $G\in\mathsf{CAGp}$ be connected second countable. Set $n=\rank(H_1^w(\mathcal F))$ and $\mathfrak{f}=\rank(\mathfrak{f}(G^*))$. If $\Cocc(G)$ is equipped with the topology of uniform convergence on compact sets then there is a topological direct sum
\begin{equation}\label{Eq:Div.in.Cocc.top.G}
\Cocc(G)=\Div(\Cocc(G))\oplus\mathbb Z^{n\mathfrak{f}},
\end{equation}
where $\Div(\Cocc(G))$ is the additive topological group of a real separable Fr\'echet space, $\mathbb Z$ is a discrete group and $\mathbb Z^{n\mathfrak{f}}$ carries the product topology. If the group $G$ is additionally assumed finite-dimensional then $\Div(\Cocc(G))$ is in fact the additive topological group of a real (separable) Banach space.
\end{proposition}
\begin{remark}\label{R:Cocc.top.free.ccl.G}
Notice the following facts.
\begin{itemize}
\item The group $\Cocc(G)$ is connected if and only if either $n=0$ or $\mathfrak{f}=0$. The first condition means that $H_1^w(\mathcal F)=0$. The second condition means that $G$ is a torus-free group, that is, the maximal torus of $G$ is trivial.
\item Let the group $G$ be finite-dimensional and torus-free (say, a solenoid). While it is not at all apparent from the definition of $\Cocc(G)$ that it is a real separable Banach space, the proposition states that this is indeed the case. In particular, the group $\Cocc(G)$ is arc-wise connected also when the group $G$ is not.
\end{itemize}
\end{remark}
\begin{proof}[Proof of Proposition~\ref{P:Cocc.top.free.ccl.G}]
Before turning to the proof we wish to mention that, similarly to the proof of Theorem~\ref{C:GpExt.frccl.gen.G}, all $\Hom$ groups are assumed to carry the topology of point-wise convergence. Now, in our proof of the isomorphism (\ref{Eq:Div.in.Cocc.top.G}) we shall need the following facts:
\begin{itemize}
\item[($\alpha$)] there is a topological direct sum
\begin{equation*}
\Cocc=\Div(\Cocc)\oplus\mathbb Z^n,
\end{equation*}
where $\Div(\Cocc)$ is the additive topological group of a real separable Banach space and the group $\mathbb Z^n$ is discrete; this follows from Proposition~\ref{P:Cocc.top.free.ccl},
\item[($\beta$)] the topological isomorphism $\Phi_G$ from Theorem~\ref{T:structure.coc} restricts to an isomorphism
\begin{equation*}
\Phi_G\colon\Div(\Cocc(G))\to\Hom(G^*,\Div(\Cocc));
\end{equation*}
indeed, since the group $\Cocc$ is torsion-free, we have
\begin{equation*}
\Div(\Hom(G^*,\Cocc))=\Hom(G^*,\Div(\Cocc)),
\end{equation*}
\item[($\gamma$)] there is a topological isomorphism $\Hom(G^*,\mathbb Z^n)\cong\mathbb Z^{n\mathfrak{f}}$; this follows from the topological isomorphisms
\begin{equation*}
\Hom(G^*,\mathbb Z^n)\cong\Hom(\mathfrak{f}(G^*),\mathbb Z^n)\cong\Hom(\mathbb Z^{(\mathfrak{f})},\mathbb Z^n)\cong\left(\mathbb Z^n\right)^{\mathfrak{f}}\cong\mathbb Z^{n\mathfrak{f}}.
\end{equation*}
\end{itemize}
From statement ($\alpha$) it follows that $\Hom(G^*,\Div(\Cocc))$ is a topological direct summand in $\Hom(G^*,\Cocc)$ and that $\Hom(G^*,\mathbb Z^n)$ is its complementary topological direct summand. By virtue of ($\beta$) and ($\gamma$), this translates into the topological direct sum (\ref{Eq:Div.in.Cocc.top.G}).

Further, since $\Div(\Cocc)$ is the additive topological group of a real separable Banach space and the group $G^*$ is countable, it follows from the first part of Lemma~\ref{L:hom.Ban.Fre} that $\Div(\Cocc(G))\cong\Hom(G^*,\Div(\Cocc))$ is the additive topological group of a real separable Fr\'echet space. If the group $G$ is finite-dimensional then $\rank(G^*)=\dim(G)$ is finite and so $\Div(\Cocc(G))$ is a real separable Banach space by the second part of Lemma~\ref{L:hom.Ban.Fre}.
\end{proof}

\begin{theorem}\label{T:Co.in.Cocc.G}
Let $\Flow$ be a minimal flow with $\Gamma\in\mathsf{LieGp}$ simply connected and with $X$ compact. Assume that the space $X$ is second countable and $\pi^1(X)\neq0$. Let $G\in\mathsf{CAGp}$ be a connected group and write $\mathfrak{r}=\rank(G^*)$. Then there is an isomorphism of groups
\begin{equation*}\label{Eq:Co.in.Cocc.G}
\begin{array}{ccccccccc}
\Cob(G) & \cong & 0 & \oplus & \mathbb R^{\mathfrak{r}} & \oplus & \Hom\left(G^*,\pi^1(X)\right) & \subseteq & {}\\
{} & \subseteq & \mathbb R^{\mathfrak{r}} & \oplus & \mathbb R^{\mathfrak{r}} & \oplus & \Hom\left(G^*,\mathbb Q\otimes\pi^1(X)\right) & \cong & \Cocc(G).
\end{array}
\end{equation*}
Further, there is a direct sum
\begin{equation}\label{Eq:dir.sm.HE}
\Hom\left(G^*,(\mathbb Q/\mathbb Z)\otimes\pi^1(X)\right)=\frac{\Hom(G^*,\mathbb Q\otimes\pi^1(X))}{\Hom(G^*,\pi^1(X))}\oplus\Ext\left(G^*,\pi^1(X)\right).
\end{equation}
Finally, there are isomorphisms of groups
\begin{equation*}
\Coch(G)\cong\mathbb R^{\mathfrak{r}}\oplus\frac{\Hom\left(G^*,\mathbb Q\otimes\pi^1(X)\right)}{\Hom\left(G^*,\pi^1(X)\right)}\cong\mathbb R^{\mathfrak{r}}\oplus\frac{\Hom\left(G^*,(\mathbb Q/\mathbb Z)\otimes\pi^1(X)\right)}{\Ext\left(G^*,\pi^1(X)\right)}.
\end{equation*}
\end{theorem}
\begin{remark}\label{R:Co.in.Cocc.G}
We wish to mention the following facts.
\begin{itemize}
\item If the group $G$ is non-trivial second countable then $1\leq\mathfrak{r}\leq\aleph_0$ and hence, similarly to Remark~\ref{R:GpExt.frccl.gen.G}, we have an isomorphism of abelian groups $\mathbb R^{\mathfrak{r}}\cong\mathbb R$.
\item Assume that the group $\pi^1(X)$ is finitely generated with rank $m$ and that the group $G^*$ is torsion-less. (This is the case, for instance, if $X$ is a manifold and $G$ is a solenoid.) Then $\Hom(G^*,\pi^1(X))=0$ and
\begin{equation*}
\begin{split}
\Hom(G^*,\mathbb Q\otimes\pi^1(X))&\cong\Hom(G^*,\mathbb Q^m)\cong\Hom(\mathbb Q\otimes G^*,\mathbb Q^m)\\
&\cong\Hom(\mathbb Q^{(\mathfrak{r})},\mathbb Q^m)\cong(\mathbb Q^m)^{\mathfrak{r}}\cong\mathbb Q^{m\mathfrak{r}}.
\end{split}
\end{equation*}
Thus, by virtue of the last statement of the theorem, there is an isomorphism $\Coch(G)\cong\mathbb R^{\mathfrak{r}}\oplus\mathbb Q^{m\mathfrak{r}}\cong\mathbb R^{\mathfrak{r}}$.
\item Let the group $\pi^1(X)$ be divisible (say, let $X=\mathbb Q^*$ be the universal metrizable solenoid). Then $\mathbb Q\otimes\pi^1(X)=\pi^1(X)$ and hence, by the last statement of the theorem, $\Coch(G)\cong\mathbb R^{\mathfrak{r}}$.
\item Assume that $G=\mathbb T^{\mathfrak{r}}$ is a torus. Then $G^*$ is free abelian and so $\Ext(G^*,\pi^1(X))=0$. Moreover,
\begin{equation*}
\begin{split}
\Hom(G^*,(\mathbb Q/\mathbb Z)\otimes\pi^1(X))&\cong\Hom(\mathbb Z^{(\mathfrak{r})},(\mathbb Q/\mathbb Z)\otimes\pi^1(X))\\
&\cong\left((\mathbb Q/\mathbb Z)\otimes\pi^1(X)\right)^{\mathfrak{r}}.
\end{split}
\end{equation*}
Consequently, by the last statement of the theorem,
\begin{equation*}
\Coch(G)\cong\mathbb R^{\mathfrak{r}}\oplus\left((\mathbb Q/\mathbb Z)\otimes\pi^1(X)\right)^{\mathfrak{r}}.
\end{equation*}
\end{itemize}
\end{remark}
\begin{proof}[Proof of Theorem~\ref{T:Co.in.Cocc.G}]
The first isomorphism of the theorem follows from (\ref{Eq:F.in.Gim.CLAC}) in the same way as (\ref{Eq:CobG.in.CoccG}) has been obtained from (\ref{Eq:Cob.in.Cocc}) in Theorem~\ref{C:GpExt.frccl.gen.G}.

To verify the second statement of the theorem, consider the injective resolution of $\pi^1(X)$
\begin{equation*}
0\longrightarrow\pi^1(X)\longrightarrow\mathbb Q\otimes\pi^1(X)\longrightarrow(\mathbb Q/\mathbb Z)\otimes\pi^1(X)\longrightarrow0.
\end{equation*}
Since the group $\mathbb Q\otimes\pi^1(X)$ is divisible, we have $\Ext(G^*,\mathbb Q\otimes\pi^1(X))=0$ and so the derived covariant Hom-Ext sequence associated to $G^*$ becomes
\begin{equation}\label{Eq:HE.seq.CW.G}
\begin{split}
0&\longrightarrow\Hom\left(G^*,\pi^1(X)\right)\longrightarrow\Hom\left(G^*,\mathbb Q\otimes\pi^1(X)\right)\longrightarrow\\
&\longrightarrow\Hom\left(G^*,(\mathbb Q/\mathbb Z)\otimes\pi^1(X)\right)\longrightarrow\Ext\left(G^*,\pi^1(X)\right)\longrightarrow0.
\end{split}
\end{equation}
To simplify notation, write (\ref{Eq:HE.seq.CW.G}) as
\begin{equation*}
0\longrightarrow A\stackrel{\alpha}{\longrightarrow}B\stackrel{\beta}{\longrightarrow}C\stackrel{\gamma}{\longrightarrow}D\longrightarrow0.
\end{equation*}
By exactness of (\ref{Eq:HE.seq.CW.G}), we have
\begin{itemize}
\item[(a)] $\im(\beta)\cong B/\ker(\beta)=B/\im(\alpha)=B/A$, and
\item[(b)] $D\cong C/\ker(\gamma)=C/\im(\beta)$.
\end{itemize}
Since the group $\mathbb Q\otimes\pi^1(X)$ is divisible torsion-free, it follows that the group $B=\Hom(G^*,\mathbb Q\otimes\pi^1(X))$ is divisible. Consequently, the group $\im(\beta)$ is also divisible and hence it is a direct summand in $C$. By virtue of (b), $D$ is a complementary summand of $\im(\beta)$ in $C$ and hence, by virtue of (a), there is a direct sum
\begin{equation*}
C\cong\im(\beta)\oplus D\cong(B/A)\oplus D.
\end{equation*}
This is precisely the desired isomorphism (\ref{Eq:dir.sm.HE}). The last statement of the theorem now follows from the first one and from (\ref{Eq:dir.sm.HE}).
\end{proof}

\section{More on first cohomology groups}\label{S:fur.struct.res}

Let $\Flow$ be a minimal flow. In this section we continue our study of the groups $\Cocc(G)$ and $\Coch(G)$ under the assumptions of Theorem~\ref{T:F.in.Gim.free.ccl}; that is, we assume that $\mathcal F$ is a topologically free flow possessing a free cycle. Being inspired by our results from Theorem~\ref{P:tor.gen.G}, which are valid for flows $\mathcal F$ with simply connected acting groups $\Gamma$, we focus our concentration on how the groups $\Coctd(G)$ and $\Cochtd(G)$ are situated in $\Cocc(G)$ and $\Coch(G)$, respectively. Unlike in Theorem~\ref{P:tor.gen.G}, neither $\Cocc$ nor $\Coctd$ is divisible under the assumptions of Theorem~\ref{T:F.in.Gim.free.ccl} in general, but the groups $\Coctd(G)$ and $\Cochtd(G)$ still turn out to form direct summands in $\Cocc(G)$ and $\Coch(G)$, respectively; this is proved in Theorem~\ref{T:cn.td.LG.mf}. Moreover, just like in Theorem~\ref{P:tor.gen.G}, the complementary summands to $\Coctd(G)$ and $\Cochtd(G)$ in $\Cocc(G)$ and $\Coch(G)$ can be chosen in such a way that all their elements $\mathcal C$ have connected sections $F(\mathcal C)$. Finally, we close this section with Theorem~\ref{P:CochG.closed.Y.N}, where we show that $\Coch(G)$ is a direct summand in $\Hom(G^*,\Coch)$ and that its closure $\overline{\Coch(G)}$ is a topological direct summand in $\Hom(G^*,\Coch)$; we also determine the corresponding complementary (topological) direct summands.

Before turning to the mentioned results, we describe one necessary construction. Let $\Flow$ be a minimal flow and $G\in\mathsf{CAGp}$. The key starting point of our construction is the following assumption:
\begin{enumerate}
\item[(a)] $\Coctd$ is a direct summand in $\Cocc$.
\end{enumerate}
Having (a) at our disposal, we proceed in the following steps.
\begin{enumerate}
\item[(i)] Fix a complementary summand $\Coccn$ of $\Coctd$ in $\Cocc$. Then $F(\mathcal D)=\mathbb T^1$ for every $1\neq\mathcal D\in\Coccn$. Thus, every element $\mathcal D$ of $\Coccn$ has a connected section $F(\mathcal D)$.
\item[(ii)] Recall that $\Phi_G\colon\Cocc(G)\to\Hom(G^*,\Cocc)$ is an isomorphism of groups. Moreover, given $\mathcal C\in\Cocc(G)$, the group $F(\mathcal C)$ is totally disconnected if and only if all the groups $\chi F(\mathcal C)=F(\chi\mathcal C)$ ($\chi\in G^*$) are totally disconnected. Consequently, $\Coctd(G)=\Phi_G^{-1}(\Hom(G^*,\Coctd))$.
\item[(iii)] Set $\Coccn(G)=\Phi_G^{-1}(\Hom(G^*,\Coccn))$. Given $\mathcal C\in\Cocc(G)$, the group $F(\mathcal C)$ is connected if and only if all the groups $\chi F(\mathcal C)=F(\chi\mathcal C)$ ($\chi\in G^*$) are connected. Consequently, every extension $\mathcal C\in\Coccn(G)$ has a connected section $F(\mathcal C)$.
\item[(iv)] The first three steps of our construction now lead to the direct sum
\begin{equation*}
\Cocc(G)=\Coccn(G)\oplus\Coctd(G).
\end{equation*}
\item[(v)] Let $\pi_G\colon\Cocc(G)\to\Coch(G)$ be the canonical quotient morphism. Set $\Cochcn(G)=\pi_G(\Coccn(G))$ and $\Cochtd(G)=\pi_G(\Coctd(G))$. Since $\ker(\pi_G)=\Cob(G)\subseteq\Coctd(G)$, it follows that the restriction $\pi_G\colon\Coccn(G)\to\Cochcn(G)$ is an isomorphism and we have a direct sum
\begin{equation*}
\Coch(G)=\Cochcn(G)\oplus\Cochtd(G).
\end{equation*}
\end{enumerate}
We now apply this construction to the situation of Theorem~\ref{T:F.in.Gim.free.ccl}.

\begin{theorem}\label{T:cn.td.LG.mf}
Let $G\in\mathsf{CAGp}$ be second countable and connected. Write
\begin{equation*}
\mathfrak{r}=\rank(G^*)\,,\hspace{3mm}\mathfrak{f}=\rank(\mathfrak{f}(G^*))\hspace{3mm}
\text{and}\hspace{3mm}\mathfrak{t}=\rank(\mathfrak{t}(G^*)).
\end{equation*}
Then, under the assumptions and notation from Theorem~\ref{T:F.in.Gim.free.ccl}, the following statements hold:
\begin{enumerate}
\item[(a)] $\Coctd$ is a direct summand in $\Cocc$,
\item[(b)] there is a direct sum
\begin{equation}\label{Eq:cn.td.LG.mf1}
\begin{array}{ccccc}
\Cocc(G) & = & \Coccn(G) & \oplus & \Coctd(G) \\
{} & \cong & \mathbb R^{\mathfrak{r}} & \oplus & \mathbb R^{\mathfrak{r}}\oplus\mathbb Z^{n\mathfrak{f}},
\end{array}
\end{equation}
\item[(c)] there is a direct sum
\begin{equation}\label{Eq:cn.td.LG.mf2}
\begin{array}{ccccc}
\Coch(G) & = & \Cochcn(G) & \oplus & \Cochtd(G) \\
{} & \cong & \mathbb R^{\mathfrak{r}} & \oplus & \mathbb Q^{m\mathfrak{t}}\oplus(\mathbb Q/\mathbb Z)^{m\mathfrak{f}}\oplus\left(\bigoplus_{i=1}^n\mathbb Z_{d_i}\right)^{\mathfrak{f}}.
\end{array}
\end{equation}
\end{enumerate}
\end{theorem}
\begin{remark}\label{R:cn.td.LG.mf}
If the acting group $\Gamma$ of a flow $\mathcal F$ is a simply connected Lie group and the phase space $X$ of $\mathcal F$ is compact then the group $\Cocc$ is divisible. By Lemma~\ref{L:tor.gen.G}, the group $\Coctd$ is also divisible and hence it is a direct summand in $\Cocc$. Consequently, the assumption (a) above is satisfied and the procedure described in steps (i)-(v) applies (we have seen this already in Theorem~\ref{P:tor.gen.G}). In the situation of Theorem~\ref{T:cn.td.LG.mf} one can not apply the mentioned results from Chapter~\ref{S:E.and.alg-top}, for neither $\Cocc$ nor $\Coctd$ is a divisible group in general. Nevertheless, the group $\Coctd$ is still a direct summand in $\Cocc$.
\end{remark}
\begin{proof}[Proof of Theorem~\ref{T:cn.td.LG.mf}]
By virtue of (\ref{Eq:Cob.in.Cocc}), there is an isomorphism
\begin{equation*}
\Coctd\cong0\oplus\mathbb R\oplus\mathbb Q^m\oplus\mathbb Z^n\cong\mathbb R\oplus\mathbb Q^m\oplus\mathbb Z^n.
\end{equation*}
Thus, $\Coctd$ is a direct summand in $\Cocc$ with a complementary summand
\begin{equation*}
\Coccn\cong\mathbb R\oplus0\oplus0\oplus0\cong\mathbb R,
\end{equation*}
which verifies statement (a).

By using step (ii) above and by applying the isomorphism $\Phi_G$, the inclusion $\Cob(G)\subseteq\Coctd(G)$ translates into the inclusion $\Hom(G^*,\Cob)\subseteq\Hom(G^*,\Coctd)$. The isomorphism (\ref{Eq:Cob.in.Cocc}) and the results from the proof of Theorem~\ref{C:GpExt.frccl.gen.G} thus lead to the isomorphisms
\begin{equation*}
\begin{array}{ccccccccccc}
\Cob(G) & \cong & \mathbb R^{\mathfrak{r}} & \oplus & 0 & \oplus & \mathbb Z^{m\mathfrak{f}} & \oplus & \left(\bigoplus_{i=1}^nd_i\mathbb Z\right)^{\mathfrak{f}} & \subseteq & {} \\
{} & \subseteq & \mathbb R^{\mathfrak{r}} & \oplus & \mathbb Q^{m\mathfrak{t}} & \oplus & \mathbb Q^{m\mathfrak{f}} & \oplus & \left(\mathbb Z^n\right)^{\mathfrak{f}} & \cong & \Coctd(G).
\end{array}
\end{equation*}
This yields, on one hand,
\begin{equation}\label{Eq:cn.td.mf.LG1}
\Coctd(G)\cong\mathbb R^{\mathfrak{r}} \oplus \mathbb Q^{m\mathfrak{t}} \oplus \mathbb Q^{m\mathfrak{f}} \oplus \left(\mathbb Z^n\right)^{\mathfrak{f}}\cong\mathbb R^{\mathfrak{r}}\oplus\mathbb Q^{m\mathfrak{r}}\oplus\mathbb Z^{n\mathfrak{f}}\cong\mathbb R^{\mathfrak{r}}\oplus\mathbb Z^{n\mathfrak{f}}
\end{equation}
and, on the other hand,
\begin{equation}\label{Eq:cn.td.mf.LG2}
\Cochtd(G)=\Coctd(G)/\Cob(G)\cong\mathbb Q^{m\mathfrak{t}}\oplus(\mathbb Q/\mathbb Z)^{m\mathfrak{f}}\oplus\left(\bigoplus_{i=1}^n\mathbb Z_{d_i}\right)^{\mathfrak{f}}.
\end{equation}
Also, by virtue of steps (iii) and (v) above,
\begin{equation}\label{Eq:cn.td.mf.LG3}
\Cochcn(G)\cong\Coccn(G)\cong\Hom(G^*,\Coccn)\cong\Hom(G^*,\mathbb R)\cong\mathbb R^{\mathfrak{r}}.
\end{equation}

Now, statement (b) of the theorem follows from step (iv) and from the isomorphisms (\ref{Eq:cn.td.mf.LG1}) and (\ref{Eq:cn.td.mf.LG3}). Similarly, statement (c) follows from step (v) and from the isomorphisms (\ref{Eq:cn.td.mf.LG2}) and (\ref{Eq:cn.td.mf.LG3}).
\end{proof}

Before turning to the next auxiliary result, we wish to make some observations. Consider the injective resolution of the group $\mathbb Z$:
\begin{equation}\label{Eq:Z.IJ}
0\longrightarrow\mathbb Z\stackrel{j}{\longrightarrow}\mathbb Q\stackrel{p}{\longrightarrow}\mathbb Q/\mathbb Z\longrightarrow0.
\end{equation}
Given an abelian group $D$, we have $\Ext(D,\mathbb Q)=0$ by divisibility of $\mathbb Q$, and so the covariant Hom-Ext sequence derived from (\ref{Eq:Z.IJ}) and associated to the group $D$ becomes
\begin{equation}\label{Eq:Z.IJ.t.HE}
0\longrightarrow\Hom(D,\mathbb Z)\stackrel{\widehat{j}}{\longrightarrow}\Hom(D,\mathbb Q)\stackrel{\widehat{p}}{\longrightarrow}\Hom(D,\mathbb Q/\mathbb Z)\stackrel{\epsilon}{\longrightarrow}\Ext(D,\mathbb Z)\longrightarrow0,
\end{equation}
where $\widehat{j}=\Hom(D,j)$, $\widehat{p}=\Hom(D,p)$ and $\epsilon$ is the connecting morphism. It follows from the exactness of (\ref{Eq:Z.IJ.t.HE}) that there are isomorphisms
\begin{equation*}
\begin{split}
\Hom(D,\mathbb Q/\mathbb Z)&\supseteq\im(\widehat{p})\cong\Hom(D,\mathbb Q)/\ker(\widehat{p})=\Hom(D,\mathbb Q)/\im(\widehat{j})\\
&=\Hom(D,\mathbb Q)/\Hom(D,\mathbb Z).
\end{split}
\end{equation*}
Thus the quotient group $\Hom(D,\mathbb Q)/\Hom(D,\mathbb Z)$ is naturally identified with the subgroup $\im(\widehat{p})$ of $\Hom(D,\mathbb Q/\mathbb Z)$ and we shall use this identification in the following lemma.

\begin{lemma}\label{L:hom.rat.not.clsd}
Let $D$ be a countable torsion-free abelian group. Equip the groups $\mathbb Q$, $\mathbb Q/\mathbb Z$ with the discrete topology and the group $\Hom(D,\mathbb Q/\mathbb Z)$ with the topology of point-wise convergence. Then $\Hom(D,\mathbb Q)/\Hom(D,\mathbb Z)$ is a dense subgroup of $\Hom(D,\mathbb Q/\mathbb Z)$. Moreover, the following conditions are equivalent:
\begin{enumerate}
\item[(a)] $\Hom(D,\mathbb Q)/\Hom(D,\mathbb Z)$ is a closed subgroup of $\Hom(D,\mathbb Q/\mathbb Z)$,
\item[(b)] $\Hom(D,\mathbb Q)/\Hom(D,\mathbb Z)=\Hom(D,\mathbb Q/\mathbb Z)$,
\item[(c)] $\Ext(D,\mathbb Z)=0$,
\item[(d)] $D$ is free abelian, that is, $\mathfrak{t}(D)=0$.
\end{enumerate}
\end{lemma}
\begin{proof}
We begin by verifying the density of $\im(\widehat{p})$ in $\Hom(D,\mathbb Q/\mathbb Z)$. So let $0\neq F\subseteq D$ be a finite set and fix $h\in\Hom(D,\mathbb Q/\mathbb Z)$; we need to show that there is $k\in\Hom(D,\mathbb Q)$ such that $pk$ and $h$ agree on $F$. Since $D$ is a torsion-free group, the subgroup $\langle F\rangle$ of $D$ generated by $F$ is free abelian with a finite rank. Let $d_1,\dots,d_n$ be a basis for $\langle F\rangle$. For $i=1,\dots,n$ let $q_i\in\mathbb Q$ be such that $p(q_i)=h(d_i)$ and denote by $k_0$ the morphism $\langle F\rangle\to\mathbb Q$ with $k_0(d_i)=q_i$ ($i=1,\dots,n$). By divisibility of $\mathbb Q$, $k_0$ extends to a morphism $k\colon D\to\mathbb Q$. Clearly, $pk$ coincides with $h$ on $\langle F\rangle$ and hence also on $F$.

We verify the equivalence of the conditions (a)--(d). First, conditions (a), (b) are equivalent by virtue of the first statement of the theorem. Further, the equivalence of (b) and (c) follows from the exactness of the sequence (\ref{Eq:Z.IJ.t.HE}). Finally, the equivalence of (c) and (d) follows from the Pontryagin theorem, which states that in the class of countable groups $D$, the Whitehead groups coincide with the free abelian groups.
\end{proof}

\begin{theorem}\label{P:CochG.closed.Y.N}
Let the assumptions of Theorem~\ref{T:F.in.Gim.free.ccl} be fulfilled and let $G\in\mathsf{CAGp}$ be connected second countable. Then, under the identification $\Coch(G)\cong\im(\Psi_G)\subseteq\Hom(G^*,\Coch)$, $\Coch(G)$ is a direct summand in $\Hom(G^*,\Coch)$ in the algebraic sense and there is a direct sum
\begin{equation}\label{Eq:Cch(G).dir.sum.hw1}
\Hom\left(G^*,\Coch\right)\cong\Coch(G)\oplus\Ext\left(\mathfrak{t}(G^*),\mathbb Z^m\right)\oplus\left(\bigoplus_{i=1}^n\text{$\tor$}_{d_i}(\mathfrak{t}(G))\right).
\end{equation}
Moreover, the following conditions are equivalent:
\begin{enumerate}
\item[(1)] $\Coch(G)$ is a topological direct summand in $\Hom(G^*,\Coch)$,
\item[(2)] $\Coch(G)$ is a closed subgroup of $\Hom(G^*,\Coch)$,
\item[(3)] $\Coch(G)$ coincides with $\Hom(G^*,\Coch)$,
\item[(4)] $G^*$ is a free abelian group, that is, $\mathfrak{t}(G^*)=0$,
\item[(5)] $G$ is a torus.
\end{enumerate}
Finally, the closure $\overline{\Coch(G)}$ of $\Coch(G)$ in $\Hom(G^*,\Coch)$ is a topological direct summand in $\Hom(G^*,\Coch)$ and there is a topological direct sum
\begin{equation}\label{Eq:Cch(G).dir.sum.hw2}
\Hom\left(G^*,\Coch\right)\cong\overline{\Coch(G)}\oplus\left(\bigoplus_{i=1}^n\text{$\tor$}_{d_i}(\mathfrak{t}(G))\right).
\end{equation}
\end{theorem}
\begin{remark}\label{R:CochG.closed.Y.N}
Notice the following facts.
\begin{itemize}
\item It follows from our proof of the theorem that the isomorphisms (\ref{Eq:Cch(G).dir.sum.hw1}) and (\ref{Eq:Cch(G).dir.sum.hw2}) are compatible in the sense that $\Coch(G)$ is a direct summand in $\overline{\Coch(G)}$ in the algebraic sense and its complementary summand is $\Ext(\mathfrak{t}(G^*),\mathbb Z^m)$.
\item By virtue of (\ref{Eq:Cch(G).dir.sum.hw2}), $\overline{\Coch(G)}$ is a co-compact subgroup of $\Hom(G^*,\Coch)$. It also follows from (\ref{Eq:Cch(G).dir.sum.hw2}) that $\Coch(G)$ is dense in $\Hom(G^*,\Coch)$ if and only if the group $\mathfrak{t}(G)$ contains no torsion elements of orders dividing $d=d_n$.
\end{itemize}
\end{remark}
\begin{proof}[Proof of Theorem~\ref{P:CochG.closed.Y.N}]
The first statement of the theorem and the algebraic splitting (\ref{Eq:Cch(G).dir.sum.hw1}) follow at once from the last statement of Theorem~\ref{C:GpExt.frccl.gen.G} and from the first part of Remark~\ref{R:GpExt.frccl.gen.G}.

We verify the equivalence of conditions (1)--(5). First, the relations (3)$\Rightarrow$(1) and (4)$\Leftrightarrow$(5) are clear. Further, (2) follows from (1) at once, for a topological direct summand is necessarily closed. The implication (4)$\Rightarrow$(3) is a direct consequence of the last statement of Theorem~\ref{C:GpExt.frccl.gen.G}. Now it suffices to show that (4) follows from (2). So assume that $\Coch(G)$ is closed in $\Hom(G^*,\Coch)$. Then, by virtue of Theorem~\ref{C:GpExt.frccl.gen.G}, $\Hom(\mathfrak{t}(G^*),\mathbb Q^m)$ is closed in $\Hom(\mathfrak{t}(G^*),(\mathbb Q/\mathbb Z)^m)$ and hence $\Hom(\mathfrak{t}(G^*),\mathbb Q)/\Hom(\mathfrak{t}(G^*),\mathbb Z)=\Hom(\mathfrak{t}(G^*),\mathbb Q)$ is closed in $\Hom(\mathfrak{t}(G^*),\mathbb Q/\mathbb Z)$. It follows from Lemma~\ref{L:hom.rat.not.clsd} that $\mathfrak{t}(G^*)$ is a free abelian group, that is, $0=\mathfrak{t}(\mathfrak{t}(G^*))=\mathfrak{t}(G^*)$. This verifies condition (4).

We finish the proof of the theorem by verifying the topological direct sum (\ref{Eq:Cch(G).dir.sum.hw2}). By virtue of Lemma~\ref{L:hom.rat.not.clsd}, $\Hom(\mathfrak{t}(G^*),\mathbb Q^m)$ is dense in $\Hom(\mathfrak{t}(G^*),(\mathbb Q/\mathbb Z)^m)$. Thus, (\ref{Eq:Cch(G).dir.sum.hw2}) follows at once from the last statement of Theorem~\ref{C:GpExt.frccl.gen.G}.
\end{proof}

\section{Dense groups of minimal extensions}\label{S:top.min.ext}

Let $\Flow$ be a minimal flow and $G\in\mathsf{CAGp}$ be a connected second countable group. In this section we study density properties of the groupoid $\Cocm(G)$ in $\Cocc(G)$. First, in Theorem~\ref{T:min.ext.res} we show that under mild assumptions on $\mathcal F$, the minimal extensions from $\Cocc(G)$ form a dense $G_{\delta}$ subset of $\Cocc(G)$. Then, in Theorem~\ref{T:dns.sbgp.min.div} and Remark~\ref{R:dns.sbgp.min.div} we present necessary and sufficient conditions, under which a given abelian group $R$ can be embedded into $\Cocm(G)$ in such a way that it forms a dense subgroup of $\Cocc(G)$ or $\Div(\Cocc(G))$.

\begin{theorem}\label{T:min.ext.res}
Let $\Flow$ be a minimal flow and $G\in\mathsf{CAGp}$ be a connected second countable group. Assume that one of the following conditions holds:
\begin{enumerate}
\item[(1)] $\Gamma$ is a simply connected Lie group and $X$ is a compact second countable space with $\pi^1(X)\neq0$,
\item[(2)] $\Gamma$ is a connected Lie group, $X$ is a compact manifold, the flow $\mathcal F$ is topologically free and possesses a free cycle.
\end{enumerate}
Then the minimal elements of $\Cocc(G)$ form a dense $G_{\delta}$ subset of $\Cocc(G)$.
\end{theorem}
\begin{remark}\label{R:min.ext.res}
We wish to add the following remarks.
\begin{itemize}
\item As usual, the group $\Cocc(G)$ is equipped with the topology of u.c.s. convergence. By our topological assumptions on $\Gamma$, $X$ and $G$, $\Cocc(G)$ is a Polish group and hence a Baire space. Thus, it follows from the theorem that a generic element of $\Cocc(G)$ is minimal.
\item It follows from \cite{Ell1} and \cite{JonPar} that the minimal extensions are generic if $\Gamma=\mathbb Z$ or $\Gamma=\mathbb R$. In case when the group $\Gamma$ is amenable locally compact the minimal extensions are also generic, this follows from \cite[Theorem~6]{Dir3}.
\item When proving results of this type one usually proceeds by expressing the set of minimal extensions as a countable intersection of sets, each of which one shows to be open and dense. Our results proved so far allow us to use a different approach - we proceed by verifying both the density and the $G_{\delta}$ property separately. That is, first we check that the minimal extensions form a set of type $G_{\delta}$. Then, arbitrarily close to any prescribed extension in $\Cocc(G)$, we construct (without using any limit argument) one concrete minimal element of $\Cocc(G)$.
\end{itemize}
\end{remark}
\begin{proof}[Proof of Theorem~\ref{T:min.ext.res}]
We divide the proof into three steps.

\emph{1st step.} We show that the set $\Cocm(G)\setminus e$ of the minimal elements of $\Cocc(G)$ is of type $G_{\delta}$ in $\Cocc(G)$.

Fix $z\in X$, a countable local base $(U_m)_{m\in\mathbb N}$ at $z$ in $X$ and a countable basis $(V_n)_{n\in\mathbb N}$ for the topology of $G$. Given $m,n\in\mathbb N$, denote by $\mathcal M_{m,n}$ the set of all $\mathcal C\in\Cocc(G)$ such that there is $\gamma\in\Gamma$ with $T_{\gamma}z\in U_m$ and $\mathcal C(\gamma,z)\in V_n$. Then all the sets $\mathcal M_{m,n}$ are open in $\Cocc(G)$. Furthermore, given $\mathcal C\in\Cocc(G)$, we have $F(\mathcal C)=G$ if and only if $\mathcal C\in\mathcal M_{m,n}$ for all $m,n\in\mathbb N$. Thus, indeed, $\Cocm(G)\setminus e=\bigcap_{m,n\in\mathbb N}\mathcal M_{m,n}$ is of type $G_{\delta}$ in $\Cocc(G)$.

\emph{2nd step.} We recall some useful facts that we shall make use of in the final step of the proof.

First, by virtue Theorem~\ref{T:structure.coc} and Corollary~\ref{C:min.contr.char}, there is a topological isomorphism between $\Cocc(G)$ and $\Hom(G^*,\Cocc)$, under which the minimal elements of $\Cocc(G)$ correspond to the monomorphisms $G^*\to\Cocm$. As usual, the group $\Hom(G^*,\Cocc)$ is assumed to carry the topology of point-wise convergence.

Second, by virtue of (\ref{Eq:F.in.Gim.CLAC}) from Theorem~\ref{T:F.in.Gim.CLAC} and (\ref{Eq:Cob.in.Cocc}) from Theorem~\ref{T:F.in.Gim.free.ccl}, in both cases (1) and (2) of the theorem there is a direct sum $\Cocc=A\oplus B$ such that $A$ is isomorphic to $\mathbb R$ as an abelian group and $\Cob\subseteq B$. We let $q$ denote the projection morphism $\Cocc\to A$.

Third, since both $\Gamma$ and $X$ are connected, the groups $\Cocc$ and $\Cocc(G)$ are torsion-free by virtue of Theorem~\ref{T:gimel.conn}.

\emph{3rd step.} We show that the monomorphisms $G^*\to\Cocm$ form a dense subset of $\Hom(G^*,\Cocc)$. According to the second step, this will finish the proof of the theorem.

Fix $h\in\Hom(G^*,\Cocc)$, a finite set $F\subseteq G^*$ and, for every $\chi\in F$, a neighbourhood $\mathcal U_{\chi}$ of $h(\chi)$ in $\Cocc$; we need to construct $k\in\Hom(G^*,\Cocc)$ with the following properties:
\begin{enumerate}
\item[(a)] $h+k$ is a monomorphism,
\item[(b)] $h+k$ takes its values in $\Cocm$,
\item[(c)] $(h+k)(\chi)\in\mathcal U_{\chi}$ for every $\chi\in F$.
\end{enumerate}

Consider the divisible hull $D=\mathbb Q\otimes qh(G^*)$ of $qh(G^*)$. Since the group $A\supseteq qh(G^*)$ is divisible, it follows that $D\subseteq A$ and hence $D$ is a direct summand in $A$. That is, $A=A'\oplus D$ for an appropriate subgroup $A'$ of $A$. Since both $A$ and $D$ are divisible torsion-free groups, it follows that so is $A'$. Moreover, we have $\rank(A)=\mathfrak{c}$ and $\rank(D)\leq\rank(G^*)\leq\card(G^*)\leq\aleph_0$, and hence $\rank(A')=\mathfrak{c}$. Thus, $A'$ is isomorphic to $\mathbb R$ as an abelian group. We let $q'$ stand for the projection morphism $A\to A'$.

Let $(\Upsilon_i)_{i\in I}$ be a maximal independent family in $G^*$. Then there exist a finite set $J\subseteq I$ and integers $l_{\chi}$, $l_{\chi,j}$ ($j\in J$, $\chi\in F$) such that for every $\chi\in F$, $l_{\chi}\chi=\sum_{j\in J}l_{\chi,j}\Upsilon_j$. By applying the morphism $h$, we get $\sum_{j\in J}l_{\chi,j}h(\Upsilon_j)=l_{\chi}h(\chi)\in l_{\chi}\mathcal U_{\chi}$ for every $\chi\in F$. Consequently, there is an identity neighbourhood $\mathcal U$ in $\Cocc$ such that for every $\chi\in F$,
\begin{equation}\label{Eq:Ups.chi.h.cont}
\sum_{j\in J}l_{\chi,j}\left(h(\Upsilon_j)+\mathcal U\right)\subseteq l_{\chi}\mathcal U_{\chi}.
\end{equation}

Now, since $\card(I)=\rank(G^*)\leq\aleph_0$ and $A'\cong\mathbb R\cong\mathbb Q^{(\mathfrak{c})}$, there is a direct sum $A'=\mathbb Q^{(J)}\oplus\mathbb Q^{(I\setminus J)}\oplus C$ for an appropriate subgroup $C$ of $A'$. Regarding $\mathbb Q^{(J)}$ and $\mathbb Q^{(I\setminus J)}$ as rational linear spaces, let $(s_j)_{j\in J}$ and $(t_i)_{i\in I\setminus J}$ be their respective bases. We have clearly $\mathbb Q^{(J)}\subseteq\Div(\Cocc)$, where the latter group is a separable Banach space by virtue of Corollary~\ref{C:structure.coc.simply} and Proposition~\ref{P:Cocc.top.free.ccl}. Thus, by finiteness of $J$, there is $n\in\mathbb N$ with $(1/n)s_j\in\mathcal U$ for every $j\in J$. Denote by $k$ the morphism 
\begin{equation*}
k\colon G^*\to\mathbb Q^{(J)}\oplus\mathbb Q^{(I\setminus J)}\subseteq A'\subseteq A\subseteq\Cocc
\end{equation*}
with $k(\Upsilon_j)=(1/n)s_j$ for $j\in J$ and $k(\Upsilon_i)=t_i$ for $i\in I\setminus J$. By our choice of $(\Upsilon_i)_{i\in I}$, such a morphism indeed exists and is unique. In fact, $k$ is a monomorphism.

\emph{3rd step continued.} We finish the proof of the theorem by showing that $k$ satisfies the properties (a), (b), (c) formulated above.

First, to see that $h+k$ is a monomorphism, fix $\lambda\in G^*$ with $(h+k)(\lambda)=0$. Since $qh(\lambda)\in D$, we have $q'qh(\lambda)=0$. Since $k(\lambda)\in A'$, we have $q'qk(\lambda)=k(\lambda)$. Consequently, $0=q'q(h+k)(\lambda)=q'qh(\lambda)+q'qk(\lambda)=k(\lambda)$. Since $k$ is a monomorphism, $\lambda=0$, as was to be shown.

Further, to see that $h+k$ takes its values in $\Cocm$, fix $\lambda\in G^*$ and assume that $(h+k)(\lambda)$ is a non-minimal element of $\Cocc$; we show that $\lambda=0$. Then $d(h+k)(\lambda)\in\Cob$ for an appropriate positive integer $d$. By applying the morphism $q'q$ we obtain, in the same way as above,
\begin{equation*}
dk(\lambda)=q'q\left(d(h+k)(\lambda)\right)\in q'q\left(\Cob\right)\subseteq q'q(B)=q'(0)=0.
\end{equation*}
Since the group $\Cocc$ is torsion-free by the second step of the proof, it follows that $k(\lambda)=0$. This implies $\lambda=0$, for $k$ is a monomorphism.

Finally, we need to check that $(h+k)(\chi)\in\mathcal U_{\chi}$ for every $\chi\in F$. So let $\chi\in F$. Then, by virtue of (\ref{Eq:Ups.chi.h.cont}),
\begin{equation*}
\begin{split}
l_{\chi}(h+k)(\chi)&=(h+k)(l_{\chi}\chi)=(h+k)\left(\sum_{j\in J}l_{\chi,j}\Upsilon_j\right)=\sum_{j\in J}l_{\chi,j}\left(h(\Upsilon_j)+k(\Upsilon_j)\right)\\
&=\sum_{j\in J}l_{\chi_j}\left(h(\Upsilon_j)+\frac{1}{n}s_j\right)\in\sum_{j\in J}l_{\chi,j}\left(h(\Upsilon_j)+\mathcal U\right)\subseteq l_{\chi}\mathcal U_{\chi}.
\end{split}
\end{equation*}
Since the group $\Cocc$ is torsion-free, this verifies condition (c).
\end{proof}

\begin{theorem}\label{T:dns.sbgp.min.div}
Let $\Flow$ be a minimal flow and $G\in\mathsf{CAGp}$ be a non-trivial connected second countable group. Assume that one of the following conditions holds:
\begin{enumerate}
\item[(1)] $\Gamma$ is a simply connected Lie group and $X$ is a compact second countable space with $\pi^1(X)\neq0$,
\item[(2)] $\Gamma$ is a connected Lie group, $X$ is a compact manifold, the flow $\mathcal F$ is topologically free and possesses a free cycle.
\end{enumerate}
Given a torsion-free abelian group $R$ with $\aleph_0\leq\rank(R)\leq\mathfrak{c}$, the following statements hold:
\begin{enumerate}
\item[(i)] the groupoid $\Cocm(G)$ contains an isomorphic copy of the group $R$, which forms a dense subset of $\Div(\Cocc(G))$,
\item[(ii)] the groupoid $\Cocm(G)$ contains an isomorphic copy of the group $R\oplus\mathbb Z^{(\aleph_0)}$, which forms a dense subset of $\Cocc(G)$.
\end{enumerate}
\end{theorem}
\begin{remark}\label{R:dns.sbgp.min.div}
We wish to mention the following facts.
\begin{itemize}
\item Under the assumption (1), the group $\Cocc(G)$ is divisible by virtue of Remark~\ref{R:lift.simply.con} and Theorem~\ref{T:coc.divisible}, and so every group $R$ as above has an isomorphic copy within $\Cocm(G)$ in such a way that it forms a dense set in the whole group $\Cocc(G)$. Under the assumption (2), the group $\Cocc(G)$ is in general not divisible. Thus, if the group $R$ is divisible then each of its isomorphic copies within $\Cocc(G)$ can be dense only in $\Div(\Cocc(G))$. (Recall from Proposition~\ref{P:Cocc.top.free.ccl.G} that $\Div(\Cocc(G))$ is a topological direct summand in $\Cocc(G)$ and hence it is a closed subgroup of $\Cocc(G)$.)
\item By virtue of statement (ii), the groupoid $\Cocm(G)$ contains an isomorphic copy of $\mathbb Z^{(\aleph_0)}$, which is dense in the whole group $\Cocc(G)$. This is no longer true if the group $\mathbb Z^{(\aleph_0)}$ is replaced by a finitely generated free abelian group. To demonstrate this, take $G=\mathbb T^1$ and consider case (1) of the theorem. Then $\Cocc$ is an infinite-dimensional Banach space by virtue of Corollary~\ref{C:structure.coc.simply} and Remark~\ref{R:structure.coc.simply}, and hence it does not contain any dense finitely generated group. This shows, in particular, that the assumption $\rank(R)\geq\aleph_0$ on the group $R$ is necessary. The other assumption $\rank(R)\leq\mathfrak{c}$ is also necessary, for $\Cocc$ is a Polish group and hence its cardinality is at most $\mathfrak{c}$.
\end{itemize}
\end{remark}
\begin{proof}[Proof of Theorem~\ref{T:dns.sbgp.min.div}]
In proving the theorem we proceed in four steps. To avoid a repetition of arguments, at certain points we shall invoke some of the ideas from the proof of Theorem~\ref{T:min.ext.res}.

\emph{1st step.} We fix notation and collect some useful facts.

As observed in the second step of the proof of Theorem~\ref{T:min.ext.res}, there is a direct sum $\Cocc=A\oplus B$, where $A$ is isomorphic to $\mathbb R$ as an abelian group and $\Cob\subseteq B$. Clearly, $A$ is a subgroup of $\Cocm$. The projection morphism $\Cocc\to A$ is denoted by $q$.

As in the proof of Theorem~\ref{T:min.ext.res}, we will work in the group of morphisms $\Hom(G^*,\Cocc)$ rather than in the group of extensions $\Cocc(G)$. This is possible due to the topological isomorphism $\Phi_G$ from Theorem~\ref{T:structure.coc}. Recall that $\Phi_G(\Cocm(G))=\Mon(G^*,\Cocm)$. Also, by virtue of Remark~\ref{R:lift.simply.con} and Proposition~\ref{P:Cocc.top.free.ccl}, in both cases (1) and (2) the group $\Div(\Cocc)$ is open in $\Cocc$. In fact, in case (1) we have $\Div(\Cocc)=\Cocc$ and in case (2) there is a topological direct sum $\Cocc=\Div(\Cocc)\oplus\mathbb Z^n$ with $\mathbb Z^n\cap\Cocm=0$. We shall denote by $r$ the projection morphism $\Cocc\to\Div(\Cocc)$. Finally, since the group $\Cocc$ is torsion-free,
\begin{equation*}
\Phi_G\left(\Div(\Cocc(G))\right)=\Div\left(\Hom(G^*,\Cocc)\right)=\Hom\left(G^*,\Div(\Cocc)\right).
\end{equation*}

Write $\mathfrak{r}=\rank(R)$ and denote by $\xi$ be the smallest ordinal with $\card(\xi)=\mathfrak{r}$. Since $\mathfrak{r}\geq\aleph_0$, we have $\xi\geq\omega$. Also, for every $\vartheta<\xi$, the ordinals $\eta$ with $\eta<\vartheta$ form a set, whose cardinality is less than $\card(\xi)=\mathfrak{r}$.

Fix finite sets $F_n\subseteq G^*$ ($n\in\mathbb N$) and open sets $\mathcal U_{\chi,n}\subseteq\Div(\Cocc)$ ($\chi\in F_n$, $n\in\mathbb N$) in such a way that the sets $\mathcal V_n=\bigcap_{\chi\in F_n}[\chi,\mathcal U_{\chi,n}]$ ($n\in\mathbb N$) form a basis for $\Hom(G^*,\Div(\Cocc))$. Since $\Div(\Cocc)$ is open in $\Cocc$, all the sets $\mathcal U_{\chi,n}$ are open in $\Cocc$. Consequently, by virtue of Theorem~\ref{T:min.ext.res}, for every $n\in\mathbb N$ there is a monomorphism $k_n\colon G^*\to\Cocm$ such that $k_n(\chi)\in\mathcal U_{\chi,n}$ for every $\chi\in F_n$. 

Notice from Proposition~\ref{P:Cocc.top.free.ccl.G} that $\Hom(G^*,\Div(\Cocc))$ is in general not an open subgroup of $\Hom(G^*,\Cocc)$ and so the morphisms $k_n$ ($n\in\mathbb N$) need not take their values in $\Div(\Cocc)$. To overcome this difficulty, set $h_n=rk_n$ ($n\in\mathbb N$). Then $h_n\in\mathcal V_n\subseteq\Hom(G^*,\Div(\Cocc))$ for every $n\in\mathbb N$ and we show that all $h_n$ are monomorphisms $G^*\to\Cocm$. Fix $n\in\mathbb N$.
\begin{itemize}
\item We claim that $h_n$ is a monomorphism. To this end, fix $\lambda\in G^*$ with $h_n(\lambda)=0$. Then $rk_n(\lambda)=0$. Since $k_n$ takes its values in the groupoid $\Cocm$, which intersects the complementary summand $C$ of $\Div(\Cocc)$ in $\Cocc$ only at the identity $0$, it follows that $k_n(\lambda)=0$. Thus, $\lambda=0$, for $k_n$ is a monomorphism.
\item We claim that $h_n$ takes its values in $\Cocm$. So assume that $h_n\neq0$ and fix $\lambda\in G^*$ with $h_n(\lambda)$ non-minimal. Then $d(rk_n(\lambda))=dh_n(\lambda)\in\Cob$ for some $d\in\mathbb N$. As observed above, we have $C\cap\Cocm=0$. Thus, without loss of generality, we may assume that $d$ satisfies $d(k_n(\lambda)-rk_n(\lambda))\in\Cob$. Then $dk_n(\lambda)\in\Cob$. Since $k_n\in\Mon(G^*,\Cocm)\setminus 0$, it follows that $\lambda=0$. Thus, indeed, $\im(h_n)\subseteq\Cocm$.
\end{itemize}
Finally, for $\omega\leq\eta<\xi$ set $F_{\eta}=0$, $\mathcal U_{0,\eta}=\Div(\Cocc)$, $\mathcal V_{\eta}=[0,\mathcal U_{0,\eta}]=\Hom(G^*,\Div(\Cocc))$ and $h_{\eta}=0\in\mathcal V_{\eta}$.

\emph{2nd step.} We show that there are morphisms $g_{\eta}\colon G^*\to\Div(\Cocc)$ ($\eta<\xi$), such that the following conditions are satisfied for every $\vartheta<\xi$:
\begin{enumerate}
\item[($a_{\vartheta}$)] $g_{\eta}\in\mathcal V_{\eta}$ for every $\eta\leq\vartheta$,
\item[($b_{\vartheta}$)] $(g_{\eta})_{\eta\leq\vartheta}$ is an independent family in $\Hom(G^*,\Div(\Cocc))$,
\item[($c_{\vartheta}$)] the morphisms $g_{\eta}$ ($\eta\leq\vartheta$) generate a subgroup of $\Mon(G^*,\Cocm)$.
\end{enumerate}
In constructing such a family of morphisms we shall proceed by transfinite induction on $\vartheta$.

Fix $\theta<\xi$ and assume that a family of morphisms $g_{\eta}\colon G^*\to\Div(\Cocc)$ ($\eta<\theta$) is given, satisfying conditions ($a_{\vartheta}$)--($c_{\vartheta}$) for every $\vartheta<\theta$. We need to construct a morphism $g_{\theta}\colon G^*\to\Div(\Cocc)$ so that conditions ($a_{\theta}$)--($c_{\theta}$) be satisfied. Notice that by our choice of $\xi$, $\card(\theta)<\card(\xi)=\mathfrak{r}\leq\mathfrak{c}$.

Denote by $D$ the divisible hull $D=\mathbb Q\otimes E$ of the group
\begin{equation}\label{Eq:def.gp.E}
E=qh_{\theta}(G^*)+\sum_{\eta<\theta}qg_{\eta}(G^*).
\end{equation}
Since $E\subseteq A$ and $A$ is divisible, we have $D\subseteq A$. Moreover, by divisibility of both $D$ and $A$, there is a direct sum $A=A'\oplus D$ for an appropriate divisible subgroup $A'$ of $A$. Further, since $\rank(D)=\rank(E)\leq\card(E)\leq\aleph_0\card(\theta)<\mathfrak{c}$ and $\rank(A)=\rank(\mathbb R)=\mathfrak{c}$, it follows that $\rank(A')=\mathfrak{c}$. Thus, $A'\cong\mathbb Q^{(\mathfrak{c})}\cong\mathbb R$. We let $q'$ stand for the projection morphism $A\to A'$.

Now, since $\rank(G^*)\leq\aleph_0<\mathfrak{c}=\rank(A')$, we may repeat the arguments from the third step of the proof of Theorem~\ref{T:min.ext.res} to obtain a monomorphism $k\colon G^*\to A'$ with $(h_{\theta}+k)(\chi)\in\mathcal U_{\chi,\theta}$ for every $\chi\in F_{\theta}$. Set $g_{\theta}=h_{\theta}+k$. Since both $h_{\theta}$ and $k$ take their values in $\Div(\Cocc)$, $g_{\theta}$ is a morphism $G^*\to\Div(\Cocc)$. We verify the validity of ($a_{\theta}$)--($c_{\theta}$).

First, $g_{\theta}\in\mathcal V_{\theta}$ by definition of $g_{\theta}$ and $g_{\eta}\in\mathcal V_{\eta}$ for $\eta<\theta$ by the induction hypothesis. This verifies condition ($a_{\theta}$). To verify conditions ($b_{\theta}$) and ($c_{\theta}$), the following statements must be checked:
\begin{itemize}
\item the morphisms $g_{\eta}$ ($\eta\leq\theta$) are independent in $\Hom(G^*,\Div(\Cocc))$,
\item every non-zero element of $\langle (g_{\eta})_{\eta\leq\theta}\rangle$ is a monomorphism,
\item every element of $\langle (g_{\eta})_{\eta\leq\theta}\rangle$ takes its values in $\Cocm$.
\end{itemize}
We verify all these three statements at once by showing that the following condition holds.
\begin{enumerate}
\item[(s)] Let $d$ be a positive integer, $l_{\eta}$ ($\eta\leq\theta$) be integers (all but finitely many of them zero) and $1\neq\lambda\in G^*$ be a character with $d\sum_{\eta\leq\theta}l_{\eta}g_{\eta}(\lambda)\in\Cob$. Then $l_{\eta}=0$ for every $\eta\leq\theta$.
\end{enumerate}
So fix $d$, $l_{\eta}$ ($\eta\leq\theta$) and $\lambda$ as above. (Observe that $G^*\neq1$ by our assumptions on $G$.) By virtue of (\ref{Eq:def.gp.E}) and by definition of $D$, $A'$ and $q'$, $q'qh_{\theta}=q'qg_{\eta}=0$ for every $\eta<\theta$. By definition of $A'$, $q'$ and $k$, $q'qk=k$. Consequently,
\begin{equation*}
\begin{split}
0&=q'(0)=q'q(B)\supseteq q'q(\Cob)\ni q'q\left(d\sum_{\eta\leq\theta}l_{\eta}g_{\eta}(\lambda)\right)\\
&=q'q\left(dl_{\theta}h_{\theta}(\lambda)+dl_{\theta}k(\lambda)+d\sum_{\eta<\theta}l_{\eta}g_{\eta}(\lambda)\right)=dl_{\theta}k(\lambda).
\end{split}
\end{equation*}
Since $k$ is a monomorphism and the group $\Cocc$ is torsion-free, it follows that $l_{\theta}=0$ and hence $d\sum_{\eta<\theta}l_{\eta}g_{\eta}(\lambda)\in\Cob$. Fix $\vartheta<\theta$ with $l_{\eta}=0$ for every $\eta>\vartheta$. Then
\begin{equation*}
d\sum_{\eta\leq\vartheta}l_{\eta}g_{\eta}(\lambda)=d\sum_{\eta<\theta}l_{\eta}g_{\eta}(\lambda)\in\Cob.
\end{equation*}
Since $\lambda\neq1$, condition ($c_{\vartheta}$) yields $\sum_{\eta\leq\vartheta}l_{\eta}g_{\eta}(\lambda)=0$. By applying ($c_{\vartheta}$) once more we obtain $\sum_{\eta\leq\vartheta}l_{\eta}g_{\eta}=0$. By virtue of ($b_{\vartheta}$), this implies $l_{\eta}=0$ for every $\eta\leq\vartheta$. Thus, $l_{\eta}=0$ for every $\eta\leq\theta$, which verifies statement (s).

\emph{3rd step.} We verify statement (i) of the theorem.

Let a group $R$ satisfy the assumptions of the theorem. We need to find an isomorphic copy $\mathcal R$ of $R$ within $\Mon(G^*,\Cocm)$ in such a way that $\mathcal R$ form a dense subset of $\Hom(G^*,\Div(\Cocc))$.

First, since $\Div(\Cocc)$ is a divisible torsion-free abelian group, it follows that so is $\Hom(G^*,\Div(\Cocc))$. Further, by the second step of the proof, the morphisms $g_{\eta}$ ($\eta<\xi$) form an independent family in $\Hom(G^*,\Div(\Cocc))$ with cardinality $\card(\xi)$ and the set $\{g_n : n\in\mathbb N\}$ is dense in $\Hom(G^*,\Div(\Cocc))$. Finally, since $\card(\xi)=\mathfrak{r}=\rank(R)$, the group $\Hom(G^*,\Div(\Cocc))$ contains an isomorphic copy $\mathcal R$ of $R$, for which $(g_{\eta})_{\eta<\xi}$ is a maximal independent set. In particular, $\mathcal R\supseteq\{g_n : n\in\mathbb N\}$ is dense in $\Hom(G^*,\Div(\Cocc))$. Also, by the second step of the proof, the subgroup $\mathcal S=\langle(g_{\eta})_{\eta<\xi}\rangle$ of $\Hom(G^*,\Div(\Cocc))$ generated by the morphisms $g_{\eta}$ ($\eta<\xi$) is contained in $\Mon(G^*,\Cocm)$.

To finish the proof of statement (i), we need to show that the group $\mathcal R$ is contained in $\Mon(G^*,\Cocm)$. So fix $h\in\mathcal R$. By definition of $\mathcal S$ and $\mathcal R$, there is a positive integer $d\in\mathbb N$ with $dh\in\mathcal S$. Then $dh\in\Mon(G^*,\Cocm)$. If $dh=0$ then also $h=0$, for the group $\Cocc$ is torsion-free. If $dh$ is a monomorphism $G^*\to\Cocm$ then $h$ is clearly a monomorphism $G^*\to\Cocc$. Moreover, $h$ takes its values in $\Cocm$, for $\kappa_d^{-1}(\Cocm\setminus1)\subseteq\Cocm\setminus1$. This verifies the desired inclusion $\mathcal R\subseteq\Mon(G^*,\Cocm)$.

\emph{4th step.} We verify statement (ii) of the theorem. First observe that under the assumption (1), the group $\Cocc(G)$ is divisible by virtue of Remark~\ref{R:lift.simply.con} and Theorem~\ref{T:coc.divisible}. Thus, in this case, statement (ii) follows at once from statement (i). For this reason, we shall restrict ourselves to the case (2).

Let a group $R$ satisfy the assumptions of the theorem. We need to find an isomorphic copy of $R\oplus\mathbb Z^{(\aleph_0)}$ within $\Mon(G^*,\Cocm)$ that would form a dense subset of $\Hom(G^*,\Cocc)$.

By the third step of the proof, $\Hom(G^*,\Div(\Cocc))$ contains a subgroup of the form $\mathcal R\oplus\mathcal Q$ such that the following conditions hold:
\begin{itemize}
\item $\mathcal R\oplus\mathcal Q$ is contained in $\Mon(G^*,\Cocm)$,
\item the group $\mathcal R$ (in fact, an appropriate countable independent subset of $\mathcal R$) is dense in $\Hom(G^*,\Div(\Cocc))$,
\item there are isomorphisms of groups $\mathcal R\cong R$ and $\mathcal Q\cong\mathbb Q^{(\aleph_0)}$.
\end{itemize}
Let $(f_m)_{m\in\mathbb N}$ be a rational linear basis for $\mathcal Q$. Since the group $G^*$ is countable by the assumptions of the theorem, we may assume that $f_m\to0$ point-wise on $G^*$. (This is accomplished by replacing $(f_m)_{m\in\mathbb N}$ by $((1/k_m)f_m)_{m\in\mathbb N}$ with an appropriately chosen sequence of positive integers $(k_m)_{m\in\mathbb N}$.)

Set $n=\rank(H_1^w(\mathcal F))$. By virtue of Proposition~\ref{P:Cocc.top.free.ccl},
there is a topological direct sum $\Cocc=\Div(\Cocc)\oplus\mathbb Z^n$. Consequently, there is a topological direct sum
\begin{equation}\label{Eq:splt.Hom.G*.G}
\Hom(G^*,\Cocc)=\Hom(G^*,\Div(\Cocc))\oplus\Hom(G^*,\mathbb Z^n).
\end{equation}
By virtue of (\ref{Eq:Cob.in.Cocc}) from Theorem~\ref{T:F.in.Gim.free.ccl}, $\mathbb Z^n\cap\Cocm=0$. Since the group $G^*$ is countable, the group $\Hom(G^*,\mathbb Z^n)$ is separable. Let $(h_m)_{m\in\mathbb N}$ be a sequence in $\Hom(G^*,\mathbb Z^n)$ such that the set $\{h_m : m\geq k\}$ is dense in $\Hom(G^*,\mathbb Z^n)$ for every $k\in\mathbb N$.

Now let $\mathcal Z$ be the subgroup of $\Hom(G^*,\Cocc)$ generated by the morphisms $f_m+h_m$ ($m\in\mathbb N$). We shall verify the following statements:
\begin{enumerate}
\item[($\alpha$)] the morphisms $f_m+h_m$ ($m\in\mathbb N$) are independent in $\Hom(G^*,\Cocc)$ and hence the group $\mathcal Z$ is isomorphic to $\mathbb Z^{(\aleph_0)}$,
\item[($\beta$)] $\mathcal R\cap\mathcal Z=0$,
\item[($\gamma$)] $\mathcal R\oplus\mathcal Z$ is a subgroup of $\Mon(G^*,\Cocm)$,
\item[($\delta$)] $\mathcal R\oplus\mathcal Z$ is dense in $\Hom(G^*,\Cocc)$.
\end{enumerate}
From statements ($\alpha$)--($\delta$) the claim of this step of the proof will follow. First, let us mention that for $m\in\mathbb N$ and $l\in\mathbb Z^m$ we write $l=(l_1,\dots,l_m)$, $\varphi_l=\sum_{i=1}^ml_if_i$ and $\psi_l=\sum_{i=1}^ml_ih_i$. We shall use repeatedly that $\im(\varphi_l)\subseteq\Div(\Cocm)$ and $\im(\psi_l)\subseteq\mathbb Z^n$ for all $l$.

To verify ($\alpha$), fix $m\in\mathbb N$ and $l\in\mathbb Z^m$ with $\varphi_l+\psi_l=0$; we show that $l=0$. From the identity $\Div(\Cocc)\cap\mathbb Z^n=0$ it follows that $\varphi_l=\psi_l=0$. Thus, since $f_m$ ($m\in\mathbb N$) are independent in $\Hom(G^*,\Cocc)$, $l=0$.

To verify ($\beta$), fix $m\in\mathbb N$ and $l\in\mathbb Z^m$ with $\varphi_l+\psi_l\in\mathcal R$; we show that $l=0$ and hence $\varphi_l+\psi_l=0$. Since $\mathcal R\subseteq\Hom(G^*,\Div(\Cocc))$ and $\varphi_l$ takes its values in $\Div(\Cocc)$, it follows that $\psi_l$ also takes its values in $\Div(\Cocc)$. Thus $\im(\psi_l)\subseteq\Div(\Cocc)\cap\mathbb Z^n=0$ and hence $\psi_l=0$. It follows that $\varphi_l\in\mathcal R\cap\mathcal Q=0$ and so, as above, $l=0$.

We verify ($\gamma$). To this end, fix $g\in\mathcal R$, $m\in\mathbb N$ and $l\in\mathbb Z^m$; we check that $g+(\varphi_l+\psi_l)\in\Mon(G^*,\Cocm)$. We proceed in two steps.
\begin{itemize}
\item We show that $g+(\varphi_l+\psi_l)\in\Mon(G^*,\Cocc)$. So assume that $1\neq\lambda\in G^*$ is such that $(g+(\varphi_l+\psi_l))(\lambda)=0$. Since $g(\lambda)+\varphi_l(\lambda)\in\Div(\Cocc)$ and $\psi_l(\lambda)\in\mathbb Z^n$, we obtain $g(\lambda)+\varphi_l(\lambda)=\psi_l(\lambda)=0$. Consequently, it follows from the inclusions $g+\varphi_l\in\mathcal R\oplus\mathcal Q\subseteq\Mon(G^*,\Cocm)$ that $g+\varphi_l=0$. Since $\mathcal R\cap\mathcal Q=0$, we have $g=\varphi_l=0$ and hence $l=0$. Thus, $g+(\varphi_l+\psi_l)=0$, as was to be shown.
\item We show that $g+(\varphi_l+\psi_l)$ takes its values in $\Cocm$. So assume that $1\neq\lambda\in G^*$ is such that $(g+(\varphi_l+\psi_l))(\lambda)$ is not minimal. Fix $d\in\mathbb N$ with $d(g+(\varphi_l+\psi_l))(\lambda)\in\Cob$. Since $\psi_l(\lambda)\in\mathbb Z^n$ and $\mathbb Z^n\cap\Cocm=0$, we may assume that $d\psi_l(\lambda)\in\Cob$. Then $d(g(\lambda)+\varphi_l(\lambda))\in\Cob$. Further, $g+\varphi_l\in\mathcal R\oplus\mathcal Q\subseteq\Mon(G^*,\Cocm)$, whence $g+\varphi_l=0$. Thus, as above, $g+(\varphi_l+\psi_l)=0$.
\end{itemize}

We finish the proof of the theorem by verifying ($\delta$). First, the group $\mathcal R$ is dense in $\Hom(G^*,\Div(\Cocc))$ by definition. Second, the group $\Hom(G^*,\mathbb Z^n)$ is contained in the closure of $\mathcal Z$ by our choice of $(f_m)_{m\in\mathbb N}$ and $(h_m)_{m\in\mathbb N}$. Thus, by virtue of (\ref{Eq:splt.Hom.G*.G}), the direct sum $\mathcal R\oplus\mathcal Z$ is dense in $\Hom(G^*,\Cocc)$, as was to be shown.
\end{proof}

\backmatter

\bibliographystyle{amsalpha}


\printindex[symbol]
\printindex

\end{document}